%% file: habil.tex
\documentclass[12pt]{amsbook}

\usepackage{amssymb}
\usepackage[dvips]{graphics}
\usepackage{rotating}
\usepackage{xypic}
\xyoption{curve}

\DeclareMathAlphabet{\eusm}{U}{}{}{}  
\SetMathAlphabet\eusm{normal}{U}{eus}{m}{n}
\SetMathAlphabet\eusm{bold}{U}{eus}{b}{n}

\DeclareMathAlphabet{\eufrak}{U}{}{}{}  
\SetMathAlphabet\eufrak{normal}{U}{euf}{m}{n}
\SetMathAlphabet\eufrak{bold}{U}{euf}{b}{n}

\swapnumbers

\newtheorem{theorem}{Theorem}[section]
\newtheorem{proposition}[theorem]{Proposition}
\newtheorem{lemma}[theorem]{Lemma}
\newtheorem{corollary}[theorem]{Corollary}

\theoremstyle{definition}
\newtheorem{definition}[theorem]{Definition}
\newtheorem{example}[theorem]{Example}

\theoremstyle{remark}
\newtheorem{remark}[theorem]{Remark}

\numberwithin{equation}{section}

\begin{document}
\pagenumbering{roman}

\include{titlepage}

\tableofcontents

\include{introduction}
\include{results}

\part{The Theory of L\'evy Processes on Involutive Bialgebras}

\include{levy-intro}

\part{Independences and Quantum L\'evy Processes}

\include{independence}

\include{quantum-levy}

\part{Examples and Applications of Quantum L\'evy Processes}

\include{levy-lie}

\include{braided}

\include{mall1}

\include{quasi}

\include{dilations}



\end{document}

%% file: titlepage.tex
\title{The Theory of Quantum L\'evy Processes \\[10mm] Habilitationsschrift}
\author{zur Erlangung des akademischen Grades \\ doctor rerum naturalium habilitatus (Dr.~rer.~nat.~habil.) \\ an der \\ Mathematisch-Naturwissenschaftlichen Fakult\"at \\ der Ernst-Moritz-Arndt-Unversit\"at Greifswald \\[25mm] vorgelegt von Uwe Franz \\
geboren am 19.\ Februar 1966 \\
in Wildeshausen \\[25mm] Greifswald, 2003}
\maketitle

%% file: introduction.tex
\chapter*{Introduction}
\pagenumbering{arabic}

Random variables and stochastic processes are used to describe non-determi\-nistic fluctuations in statistics, finance, actuarial mathematics, computer science, engineering, biology, physics, etc. Often these fluctuations can not be predicted with certainty due to a lack of information on the initial state which can not be avoided for practical reasons, e.g., in meteorology or gambling.  Or the unpredictable behavior of the system may be due to more fundamental reasons as in quantum mechanics. Heisenberg's uncertainty relations limit the accuracy of predictions of measurements of so-called complementary observables, e.g., momentum and position, for any given state of the particle. Even though - at least in principle - there exist states in which one of the observables can be predicted with arbitrary precision, it is never possible to predict both simultaneously.

But often it is reasonable to assume that the experiment can be repeated a large number of times under the same conditions and that the relative frequencies of the possible outcomes can be predicted. It is useful to keep in mind that it is this assumption which justifies the applicability of probabilistic models.

If the random fluctuations do not depend on time or position, then they should be described by time- and space-homogeneous stochastic processes. This leads, in Euclidean space, to the important class of stochastic processes with independent and stationary increments, i.e.\ L\'evy processes. These processes have recently received increasing interest, see \cite{levy65,skorohod91,bertoin96,sato99,barndorff-nielsen+mikosch+resnick01,applebaum} and the references therein.

Since in quantum mechanics even complete knowledge of the state is not sufficient to predict with certainty the outcomes of all possible measurements, the statistical interpretation has to be an essential part of the theory. Quantum probability starts from von Neumann's formulation of standard quantum mechanics \cite{vonneumann} and studies quantum theory from a probabilistic point of view, cf.\ \cite{accardi+frigerio+lewis82,hudson+parthasarathy84}. For an introduction to quantum probability and quantum stochastic calculus, see also \cite{biane90,parthasarathy92,meyer95} or the lecture notes \cite{attal+lindsay03a,attal+lindsay03b}.

A typical situation where quantum noise enters is the description of a small quantum system which is interacting with its environment, cf.\ \cite{gardiner+zoller00}. The state of the environment, also called heat bath, can not be measured or controlled completely, but it is natural to assume that it is homogeneous in time and space and that the influence of the small system on the heat bath can be neglected.

Structures similar to infinite divisibility and L\'evy processes also appear in quantum mechanics in the theory of repeated and continuous measurements, cf.\ \cite{holevo01}.

Another line of research leading to so-called quantum L\'evy processes is von Waldenfels' work on the emission and absorption of light, cf.\ \cite{waldenfels73,waldenfels84}.

Quantum random variables are defined as homomorphisms $j:\mathcal{B}\to\mathcal{A}$ taking values in an algebra $\mathcal{A}$ where a fixed state $\Phi$ has been chosen, and stochastic processes are indexed families of quantum random variables. In order to define quantum L\'evy processes, i.e.\ quantum stochastic processes with independent and stationary increments, we have to explain what we mean by ``independent'' and ``increment.''

One possibility to define independence is to call two elements $a,b\in\mathcal{A}$ independent, if they commute, i.e.
\[
ab=ba
\]
and if the state factorizes on them, i.e.
\[
\Phi(ab)=\Phi(a)\Phi(b).
\]
This definition generalizes the one used in classical probability theory and corresponds to what physicists would consider as independent observables. Two quantum random variables $j_1,j_2$ are called independent, if any pair $a\in {\rm range}\,j_1$, $b\in{\rm range}\,j_2$ of elements of their respective ranges is independent. This is also the notion of independence that underlies the quantum stochastic calculus on the symmetric Fock space developed by Hudson and Parthasarathy, \cite{hudson+parthasarathy84,parthasarathy92}.

In order to define what an increment is, we need a composition for quantum random variables. This is possible, e.g., if the algebra on which the quantum random variables are defined, is a bialgebra. Then it is equipped with a homomorphism $\Delta:\mathcal{B}\to\mathcal{B}\otimes\mathcal{B}$ and the convolution product $j_1\star j_2$ of two quantum random variables $j_1,j_2:\mathcal{B}\to\mathcal{A}$ can be defined as
\[
j_1\star j_1 = m_\mathcal{A}\circ(j_1\otimes j_2)\circ\Delta,
\]
where $m_\mathcal{A}:\mathcal{A}\otimes\mathcal{A}\to\mathcal{A}$ denotes the multiplication in $\mathcal{A}$, $m_\mathcal{A}(a\otimes b)=ab$. If $j_1$ and $j_2$ are independent, then $j_1\star j_2$ is again a homomorphism.

These two choices lead to the theory of L\'evy processes on involutive bialgebras as it was introduced and studied by Sch\"urmann et al., cf.\ \cite{accardi+schuermann+waldenfels88,schuermann93}. For an introduction to this theory, see also \cite[Chapter VII]{meyer95}, \cite{franz+schott99}, or Chapter \ref{levy-intro}.

However, in quantum probability there exist also other notions of independence like, e.g., freeness \cite{voiculescu+dykema+nica92}, see Chapter \ref{chapter-III}. In order to formulate a general theory of L\'evy processes for these independences, bialgebras have to be replaced by the dual groups introduced in \cite{voiculescu87,voiculescu90}. The definition of these dual groups, also called H-algebras or cogroups, is very similar to that of bialgebras and Hopf algebras, but the tensor product is replaced by the free product of algebras, see also \cite{zhang91,bergman+hausknecht96}. There exists a homomorphism $\Delta:\mathcal{B}\to\mathcal{B}\coprod\mathcal{B}$ from $\mathcal{B}$ to the free product of $\mathcal{B}$ with itself, which allows to define the convolution product $j_1\star j_2=m_\mathcal{A}\circ(j_1\coprod j_2)\circ\Delta$ of two random variables $j_1,j_2:\mathcal{B}\to\mathcal{A}$. The theory of L\'evy processes on these algebras has been developed in \cite{schuermann95b,benghorbal+schuermann99,franz01,franz03b}, see also Chapter \ref{chapter-IV}. The special case of L\'evy processes with additive free increments was first studied in \cite{glockner+schuermann+speicher92}, and more recently in \cite{biane98,anshelevich99,anshelevich01b,anshelevich01c,barndorff-nielsen+thorbjornsen01a,barndorff-nielsen+thorbjornsen01b,barndorff-nielsen+thorbjornsen01c}.

Classical L\'evy processes are not only used as models for random events, they are also related to many other areas of mathematics, e.g., potential theory, harmonic analysis, or representation theory. Similarly, quantum L\'evy processes appear in various situations. Many prominent examples of quantum stochastic processes are quantum L\'evy processes or can be constructed from these. Classical L\'evy processes and factorizable representations of current groups and current algebras are special cases of quantum L\'evy processes, see Example \ref{levy-intro}.\ref{I-exa-class-levy} and Subsection \ref{levy-intro}.\ref{I-subsec-lie}. Quantum L\'evy processes also arise in the construction of dilations of quantum dynamical semigroups, cf.\ \cite{schuermann90b}.

This thesis presents several recent contributions to the theory of quantum L\'evy processes. These works were motivated by three central ideas.

The first is the desire to understand better the role of independence in quantum probability. If one requires that the joint law of independent quantum random variables is uniquely determined by their marginal distributions and imposes certain natural conditions on the construction of joint law, then there exist exactly five possibilities, see \cite{speicher96,benghorbal+schuermann99,benghorbal+schuermann02,muraki02a,muraki02,muraki02b} or Chapter \ref{chapter-III}. It has been shown that three of them can be reduced to tensor independence, see \cite{franz03b} or Section \ref{chapter-III}.\ref{II-reduction}. This reduction can also be applied to reduce the corresponding quantum L\'evy processes to L\'evy processes on involutive bialgebras. It is not known if such a reduction is also possible for free independence.

Another central topic of this thesis is the relation between classical and quantum probability. Many interesting classical stochastic processes arise as components of quantum stochastic processes, see, e.g., \cite{parthasarathy90,schuermann91b,franz99} or also Sections \ref{levy-lie}.\ref{class} and \ref{ch-quasi}.\ref{QUASI-examples}. Some properties of these classical processes can then be proved using the whole quantum stochastic process, see, e.g., the proof of the chaos completeness of the Az\'ema martingale in  \cite{parthasarathy99} or the calculation of quasi-invariance formulas for certain increment processes in Chapter \ref{ch-quasi}. On the other hand, ideas from classical probability theory can often be applied successfully for the study of quantum stochastic processes, see Chapter \ref{ch-malliavin}.

The value of a theory is of course also determined by the richness of its examples. Therefore the third central theme is the detailed study of examples. In Chapter \ref{levy-lie}, we have looked at the special case of L\'evy processes on real Lie algebras and classified the Sch\"urmann triples belonging to unitary irreducible representations for several Lie algebras. More general Sch\"urmann triples can then be constructed as direct sums. In Section \ref{braided}.\ref{example}, the Brownian motions, i.e., L\'evy processes with quadratic generators, were classified on several standard examples of so-called braided spaces. Finally, in Chapter \ref{chapter-II}, we determined all generators and all quadratic generators on the non-commutative analogue of the algebra of coefficients of the unitary group. This bialgebra is related to the construction of dilations of quantum dynamical semigroups on the matrix algebra $\mathcal{M}_n$.

This thesis is written in cumulative form. Most of the chapters are taken from separate publications and can also be read individually. In the next section we give a detailed summary of this thesis and its main results.

\subsection*{Acknowledgements}

I would like to thank all my colleagues at the Mathematics Department of the University Greifswald, where most of this research was carried out. It has been a pleasure to share Micheal Sch\"urmann's rich expertise and I am looking forward to further ambitious collaborations. The dialogue with my friend and office mate Rolf Gohm has been very enriching for me and I am looking forward to finally having enough time for pursuing our joint projects.

I have benefited from discussions with many of those working in the area of quantum probability (taken in the large sense), it would be impossible to name all of them.

I am indebted to Luigi Accardi, B.V.R.\ Bhat, Marek Bo\.{z}ejko, Philip Feinsilver, Claus K\"ostler, Romuald Lenczewski, Nobuaki Obata, Habib Ouerdiane, Nicolas Privault, Ren\'e Schott, K.B.\ Sinha, Michael Skeide, and Roland Speicher for their kind hospitality. The visits at their institutes were of great value for my research. 

Last, but not least, I would like to my wife and my family for their support and patience. Without Aline Mio this thesis might have been completed a few months earlier, but my life would not be the same.

%% file: results.tex
\chapter*{Summary of results}

In this thesis various recent results on quantum L\'evy processes are presented. Most chapters are taken from separate publications and can be read by themselves.  The first part provides an introduction to the theory of L\'evy processes on involutive bialgebras. The notion of independence used for these processes is tensor independence, see Definition \ref{levy-intro}.\ref{ch1-def-indep}.

In quantum probability there exist other notions of independence and L\'evy processes can also be defined for the five so-called universal independences. This is the topic of the second part. In particular, we show that boolean, monotone, and anti-monotone independence can be reduced to tensor independence.

Finally, in the third part, we consider several classes of quantum L\'evy processes of special interest, e.g., L\'evy processes on real Lie algebras or Brownian motions on braided spaces. We also present several applications of these processes.

\section{L\'evy processes on involutive bialgebras}

In the first chapter, we give an introduction to the basic theory of L\'evy processes on involutive bialgebras. We recall their definition and the one-to-one correspondence between L\'evy processes, convolution semigroups of states, generators, and Sch\"urmann triples. We also recall without proof Sch\"urmann's representation theorem, see Theorem \ref{levy-intro}.\ref{I-rep-thm}. Section \ref{levy-intro}.\ref{I-cyclic} contains the recent result by Franz and Skeide that the vacuum vector is cyclic for the Fock space realization of a L\'evy process, if the cocycle is surjective, see Theorem \ref{levy-intro}.\ref{I-cyclic-thm}.

This chapter is taken from a series of lectures on quantum L\'evy processes held at the school ``Quantum Independent Increment Processes: Structure and Applications to Physics'' in Greifswald in March 2003, see \cite{franz+schuermann}.

\section{The five universal independences}

In Chapter \ref{chapter-III}, we recall the classification of universal non-commutative independences due to Ben Ghorbal and Sch\"urmann \cite{benghorbal+schuermann99,benghorbal+schuermann02}, and Muraki \cite{muraki02a,muraki02,muraki02b}, see also \cite{schuermann95,speicher96}.

The starting point of their work is the assumption that the joint law of independent random variables should be uniquely determined by their marginals. If one imposes certain natural conditions like associativity on the construction of the joint law, it becomes possible to classify all notions of independence on a given category of probability spaces. For classical probability spaces it turns out that the only possible construction is the tensor product. Therefore the usual notion of independence used in classical probability is the only universal notion of independence available for that case.

But for categories of non-commutative probability spaces several universal notions can exist. If the underlying algebra is an arbitrary unital algebra, then there exist exactly two independences: one based on the tensor product and one based on the free product of states. For non-unital algebras, five universal independence exist. In addition to tensor and free independence, we also have boolean, monotone, and anti-monotone independence, cf.\ \cite{muraki02a,muraki02,muraki02b}.

In \cite{franz02}, it was shown that the axioms imposed in \cite{benghorbal+schuermann99,benghorbal+schuermann02,muraki02} on the construction of the joint law are equivalent to saying that it is a tensor product (in the sense of category theory, see Definition \ref{chapter-III}.\ref{def-tensor}) and that there exist natural transformations embedding the probability spaces into their tensor product, see Section \ref{chapter-III}.\ref{III-class-indep}.

Once the notion of independence is reformulated in this way, it is straight forward that a reduction from one universal notion to another should be a tensor functor, or rather a cotensor functor $F$, see Definition \ref{chapter-III}.\ref{def-reduction}. In addition, we need a natural transformation from the first category to the image of $F$ in the second category. This then allows to embed any products taken in the first category nicely into products taken in the second category. Therefore all calculations requiring the first product can be reduced to calculations using the second product. In particular, it is shown that three of the universal independences, namely boolean, monotone, and anti-monotone independence, can be reduced to tensor independence, see Subsection \ref{chapter-III}.\ref{section-reduction-boolean-etc}.

This chapter is also taken from the series of lectures on quantum L\'evy processes held at the school ``Quantum Independent Increment Processes: Structure and Applications to Physics'' in Greifswald in March 2003, see \cite{franz+schuermann}, and contains a fairly detailed section with preliminaries on category theory, see Section \ref{chapter-III}.\ref{III-prelim}.

\section{L\'evy processes on dual groups}

In Chapter \ref{chapter-IV}, the reduction of boolean, monotone, and anti-monotone independence is applied to quantum L\'evy processes. Bijections between the classes of L\'evy processes on dual groups with boolean, monotonically, or anti-monotonically independent increments and certain classes of L\'evy processes on involutive bialgebras are constructed, see Theorem \ref{chapter-IV}.\ref{theo-one-to-one}. These bijections are used to classify and construct L\'evy processes on dual groups with boolean, monotonically, or anti-monotonically independent increments. These processes can always be realized on the symmetric Fock space. Therefore we can also reduce the corresponding quantum stochastic differential calculi to Hudson-Parthasarathy quantum stochastic differential calculus on the symmetric Fock space, see Subsection \ref{chapter-IV}.\ref{subsec-fock}.

As another application of the reduction we show that monotone Markov processes have a natural Markov structure, see Corollary \ref{chapter-IV}.\ref{cor-mon-markov}.

This chapter is based on the results published in \cite{franz03b}.

\section{Renormalized squares of white noise and other non-Gaussian noises as L\'evy processes on real Lie algebras}
\markboth{SUMMARY OF RESULTS}{RENORMALIZED SQUARES OF WHITE NOISE}

Let $a_t$, $a_s^+$ be the white noise creation and annihilation operators on the Fock space, cf.\ \cite{hida+kuo+potthoff+streit93,obata94}. Formally, they satisfy the relation
\begin{equation}\label{ccr}
[a_t,a^+_s]=\delta(t-s),
\end{equation}
for $s,t\in\mathbb{R}$. This relation can be given meaning by considering integrals
\[
A(f)=\int_{\mathbb{R}}\overline{f(t)}{\rm d}a_t, \qquad A^+(g)=\int_{\mathbb{R}}g(t){\rm d}a^+_t
\]
against test functions $f,g$ from some appropriate function space $\mathcal{S}(\mathbb{R})$. Then Equation \eqref{ccr} becomes
\[
\big[A(f),A^+(g)\big]=\int_{\mathbb{R}}\overline{f(t)}g(s)\delta(t-s){\rm d}t{\rm d}s= \int_{\mathbb{R}}\overline{f(t)}g(t){\rm d}t
\]
for $f,g\in\mathcal{S}(\mathbb{R})$.

Taking squares
\[
b_t=\frac{1}{2}(a_t)^2 \quad\mbox{ and } \quad b_s^+=\frac{1}{2}(a_s^+)^2,
\]
one obtains formally the relation
\begin{eqnarray*}
[b_t,b_s^+] &=& \frac{1}{4}[a_t^2,(a_s^+)^2]= \frac{1}{4}\left(2a_ta_s^+\delta(t-s)+2a_s^+a_t\delta(t-s)\right) \\
&=& a_ts_s^+\delta(t-s)+ \frac{1}{2}\big(\delta(t-s)\big)^2.
\end{eqnarray*}
But this relation cannot be given a rigorous meaning simply by integrating against test functions. To overcome this difficulty, Accardi, Lu, and Volovich \cite{accardi+lu+volovich99} proposed the ``renormalization rule''
\[
\big(\delta(t-s)\big)^2 \to 2\gamma \delta(t-s)
\]
where $\gamma$ is some ``renormalization constant.'' Using this rule, they obtained the following three algebraic relations for $b_t=(a_t)^2/2$, $b^+_s=(a_s^+)^2/2$, and $n_t=a_ta_t^+$,
\begin{subequations}\label{wn-swn-rel}
\begin{equation}
[b_t,b^+_s] = n_t\delta(t-s)+\gamma\delta(t-s),
\end{equation}
\begin{equation}
[n_t,b_s] = - 2 b_t\delta(t-s),
\end{equation}
\begin{equation}
[n_t,b^+_s] = 2 b_t^+\delta(t-s).
\end{equation}
\end{subequations}
These relations can be made rigorous by integrating them against test functions, see Equations (\ref{levy-lie}.\ref{swn-rel}). For $\gamma>0$, Accardi, Lu, and Volovich \cite{accardi+lu+volovich99} have constructed a realization of these relations by operators acting on a Hilbert space with a cyclic vector $\Omega$ which is annihilated by $b_t$ and $n_t$. They have also shown that such a realization does not exist for $\gamma<0$.

Chapter \ref{levy-lie} starts from the relations \eqref{wn-swn-rel} and shows that any factorizable current representation of the current algebra $\eufrak{sl}_2^{\mathbb{R}}$ over the Lie algebra $\eufrak{sl}_s=sl(2,\mathbb{R})$ gives a realization. Furthermore, any realization of \eqref{wn-swn-rel} that satisfies $(b_t)^*=b_t$ and $(n_t)^*=n_t$ and an independence property can be obtained in this way. This includes also the realization constructed in \cite{accardi+lu+volovich99}.

In Chapter \ref{levy-lie}, we also recall how these factorizable current representations can be classified and constructed using the representation theory and cohomological properties of the underlying Lie algebra. For several Lie algebras we compute the relevant cohomological information for their unitary irreducible representations and thereby classify their factorizable current representations.

In the last section of Chapter \ref{levy-lie}, we show that one can associate a classical L\'evy process to these factorizable current representations by fixing an element of the Lie algebra, see Theorem \ref{levy-lie}.\ref{levy-lie-class}. For the realization of \cite{accardi+lu+volovich99}, one finds that the marginal distributions of these L\'evy processes are the measures of orthogonality of three of the five Meixner classes \cite{meixner34} of orthogonal polynomials, namely the Meixner-Pollaczek, Laguerre, and Meixner polynomials. In probability theory these processes are called Meixner, Gamma, and Pascal process.

The results of this chapter have been published in \cite{accardi+franz+skeide02}.

\section{L\'evy processes and Brownian motion on braided spaces}

To get more than the five universal independences, one has to consider categories of probability spaces with additional structure. In Example \ref{chapter-III}.\ref{fermi-indep}, Fermi independence was defined for $\mathbb{Z}_2$-graded quantum probability spaces. Recently, the theory of the so-called $\mathbb{Z}_2$-graded superspaces and supergeometry has been generalized to braided spaces and braided geometry, see \cite{majid95} and the references therein. In Chapter \ref{braided}, the notion of braided independence is developed. Braided quantum probability spaces are $*$-algebras equipped with an action and coaction of some $*$-bialgebra and a state that satisfies an invariance condition.

In Chapter \ref{braided}, we show how a braided category can be constructed from an involutive Hopf algebra or from an involutive bialgebra equipped with an $r$-form. Algebras and bialgebras in such a category are called braided algebras and braided bialgebras, if their structure maps are morphisms of the category. Braided spaces are braided bialgebras which are generated by a finite number of primitive elements. We show how such braided spaces can be constructed for any bi-invertible, real type I R-matrix, see Section \ref{braided}.\ref{construction}. The construction is similar to the one in \cite{majid95}, but the axioms for the involution are new.

In Section \ref{braided}.\ref{symmetrization}, we show that for a braided $*$-bialgebras $\mathcal{B}$ it is possible to construct a $*$-bialgebra $H$ that contains $\mathcal{B}$ as a subalgebra and as a one-sided coideal. Furthermore, it is possible to extend the generator $L:\mathcal{B}\to\mathbb{C}$ in a canonical way to a generator $L^H:H\to\mathbb{C}$ such that the L\'evy process $(j_{st})_{0\le s\le t}$ associated to $L$ can be recovered from the L\'evy process $(j^H_{st})_{0\le s\le t}$ associated to $L^H$. This construction can actually be formulated as a reduction in the sense of Definition \ref{chapter-III}.\ref{def-reduction}.

In Section \ref{braided}.\ref{example}, we classify all quadratic generators on several standard examples of braided spaces. This gives a classification of all Brownian motions on these spaces, i.e., of the analogue of time- and space-homogeneous diffusions.

\section{Malliavin calculus and Skorohod integration for quantum stochastic processes}

In Chapter \ref{ch-malliavin}, we develop a non-commutative analogue of infinite-dimen\-sional analysis on the Wiener space. Our goal is to find sufficient conditions for the regularity of joint densities of non-commuting operators.

Let $\eufrak{h}$ be a real Hilbert space. The Wiener space $W(\eufrak{h})$ is replaced by the algebra of bounded operators on its complexification $\Gamma(\eufrak{h}_\mathbb{C})=W(\eufrak{h})_\mathbb{C}$. On $\Gamma(\eufrak{h}_\mathbb{C})$ we have two non-commuting quantum stochastic processes $Q$ and $P$, indexed by elements of $\eufrak{h}$, which satisfy
\[
[P(h),Q(k)] = 2i\langle h,k\rangle, \qquad \mbox{ for }h,k\in\eufrak{h}.
\]
Operators which can be obtained via Weyl's functional calculus play the role of smooth functionals. For these operators we obtain a Girsanov type formula in Proposition \ref{ch-malliavin}.\ref{prop girsanov} and as its infinitesimal version an integration by parts formula, see Proposition \ref{ch-malliavin}.\ref{prop-int by parts}. The left-hand-side of the integration by parts formula can be interpreted as a directional or Fr\'echet derivative,
\[
D_k B=\frac{i}{2}[Q(k_1)-P(k_2),B]
\]
for $k=(k_1,k_2)$, $k_1,k_2\in\eufrak{h}$, and $B$ a ``smooth'' operator. As a next step a derivation operator $D$ is defined, which is related to the Fr\'echet derivative by
\[
D_kB = \langle k,DB\rangle
\]
for $k=(k_1,k_2)$, $k_1,k_2\in\eufrak{h}$, and $B$ a ``smooth'' operator, see Proposition \ref{ch-malliavin}.\ref{prop frechet}.

Even though $D$ is no longer a map between Hilbert spaces, we can define also an analogue of the divergence operator $\delta$. Both the derivation operator and the divergence operator are shown to be closable.

In the white noise case, i.e., if $\eufrak{h}$ is the space of $L^2$-functions on some measurable space, e.g., $\eufrak{h}=L^2(\mathbb{R}_+)$, the classical divergence operator is the Hitsuda-Skorohod integral, which extends the Wiener integral to not necessarily adapted processes. Here it turns out that the divergence operators coincides with a non-causal extension of the Hudson-Parthasarathy quantum stochastic integral, which had already been studied by Belavkin \cite{belavkin91a,belavkin91b} and Lindsay \cite{lindsay93}.

We also show that the derivation operator can also be used to obtain sufficient conditions for the existence of smooth densities, see Subsection \ref{ch-malliavin}.\ref{mall-smooth}.

This chapter is taken from \cite{franz+leandre+schott00}, see also \cite{franz+leandre+schott99}.

\section{Quasi-invariance formulas for components of quantum L\'evy processes}
\markboth{SUMMARY OF RESULTS}{QUASI-INVARIANCE FORMULAS}

The Girsanov type formula in Proposition \ref{ch-malliavin}.\ref{prop girsanov} actually contains the classical Girsanov formula for Brownian motion with a deterministic drift as a special case. In Chapter \ref{ch-quasi}, we study this situation in more detail. We start from a quantum stochastic process whose restriction to a commutative subalgebra gives a classical stochastic process. Then we choose unitary operators and let them act on the range of the quantum stochastic process. This action can be computed in two ways. We can let the unitary act on the algebra or on the state vector. Since the expectations computed either way have to agree, this leads to a kind of quasi-invariance formula. Under certain conditions this quasi-invariance formula carries over to the classical process.

The basic idea of this construction is described in Section \ref{ch-quasi}.\ref{QUASI-idea}. In the following section we give several examples, including the classical Girsanov formula and a quasi-invariance formula for the Gamma process recently obtained by Tsilevich, Vershik, and Yor \cite{tsilevich+vershik+yor00,tsilevich+vershik+yor01}. We also derive a new invariance formula for the Meixner process, see Subsection \ref{ch-quasi}.\ref{QUASI-exa-meixner}.

\section{L\'evy processes and dilations of completely positive semigroups}

Let $\mathcal{U}_d$ be the free $*$-algebra generated by indeterminates $u_{ij}$, $u^*_{ij}$, $i,j=1,\ldots,d$, with the relations
\begin{eqnarray*}
\sum_{j=1}^d u_{kj}u^*_{\ell j} &=& \delta_{k\ell}, \\
\sum_{j=1}^d u^*_{jk}u_{j\ell} &=& \delta_{k\ell},
\end{eqnarray*}
The $*$-algebra $\mathcal{U}_d$ is turned into a $*$-bialgebra, if we define $\Delta:\mathcal{U}_d\to\mathcal{U}_d\otimes\mathcal{U}_d$ and $\varepsilon:\mathcal{U}_d\to\mathbb{C}$ on the generators by
\begin{eqnarray*}
\Delta(u_{k\ell}) &=& \sum_{j=1}^d u_{kj}\otimes u_{j\ell}, \\
\varepsilon(u_{k\ell}) &=& \delta_{k\ell},
\end{eqnarray*}
and extend them as $*$-algebra homomorphisms. Sch\"urmann \cite{schuermann90b} has shown that the Hudson-Parthasarathy cocycles which are used to construct dilations of quantum dynamical semigroups on $\mathcal{M}_d$ are equivalent to L\'evy processes on $\mathcal{U}_d$.

In Chapter \ref{chapter-II}, we classify the Gaussian generators and all generators on $\mathcal{U}_d$, see Theorems \ref{chapter-II}.\ref{class un gauss} and \ref{chapter-II}.\ref{class un}. We also recall how dilations of quantum dynamical semigroups on $\mathcal{M}_d$ can be constructed from L\'evy processes on $\mathcal{U}_d$. Furthermore, we give the relation between the generator of the quantum dynamical semigroup and the generator of the L\'evy process.

%% file: levy-intro.tex
\chapter{L\'evy Processes on Involutive Bialgebras}\label{levy-intro}

L\'evy processes on involutive bialgebras first appeared in models of the laser, cf.\ \cite{waldenfels73,waldenfels84}. Their algebraic framework was formulated in \cite{accardi+schuermann+waldenfels88}. They are a generalization of both classical stochastic processes with independent and stationary increments, called L\'evy processes, and factorizable current representations of groups and Lie algebras.

In this Chapter we will give the definition of these processes, cf.\ Section \ref{I-def}.

In Section \ref{I-gen} be will begin to develop their basic theory. We will see that the marginal distributions of a L\'evy process form a convolution semigroup of states and that we can associate a generator with a L\'evy process on an involutive bialgebra, that characterizes uniquely its distribution, like in classical probability.  By a GNS-type construction we can get a so-called Sch\"urmann triple from the generator.

This Sch\"urmann triple can be used to obtain a realization of the process on a symmetric Fock space, see Section \ref{I-rep}. This realization can be found as the (unique) solution of a quantum stochastic differential equation. It establishes the one-to-one correspondence between L\'evy processes, convolution semigroups of states, generators, and Sch\"urmann triples. We will not present the proof of the representation theorem here, but refer to \cite[Chapter 2]{schuermann93}.

In Section \ref{I-cyclic}, we present a recent unpublished result by Franz and Skeide. If the cocycle of the Sch\"urmann triple is surjective, then the vacuum vector is cyclic for the L\'evy process constructed on the symmetric Fock space via the representation theorem.

Finally, in Section \ref{I-exa}, we look at several examples.

For more details, see also \cite{schuermann93}\cite[Chapter VII]{meyer95}\cite{franz+schott99}.

\section{Definition of L\'evy processes on involutive bialgebras}\label{I-def}

A {\em quantum probability space} in the purely algebraic sense is a pair $({\mathcal A},\Phi)$ consisting of a unital $*$-algebra ${\mathcal A}$ and a state (i.e.\ a normalized positive linear functional) $\Phi$ on ${\mathcal A}$. A {\em quantum random variable} $j$ over a quantum probability space $(\mathcal{A},\Phi)$ on a $*$-algebra ${\mathcal B}$ is simply a $*$-algebra homomorphism $j:{\mathcal B}\to{\mathcal A}$. A {\em quantum stochastic process} is an indexed family of random variables $(j_t)_{t\in I}$. For a quantum random variable $j:{\mathcal B}\to{\mathcal A}$ we will call $\varphi_j=\Phi\circ j$ its {\em distribution} in the state $\Phi$. For a quantum stochastic process $(j_t)_{t\in I}$ the functionals $\varphi_t=\Phi\circ j_t:\mathcal{B}\to\mathbb{C}$ are called {\em marginal distributions}. The {\rm joint distribution} $\Phi\circ\left(\coprod_{t\in I} j_t\right)$ of a quantum stochastic process is a functional on the free product $\coprod_{t\in I} \mathcal{B}$, see Chapter \ref{chapter-III}.

 Two stochastic processes $\left(j^{(1)}_{t}:\mathcal{B}\to\mathcal{A}_1\right)_{t\in I}$ and $\left(j^{(2)}_{t}:\mathcal{B}\to\mathcal{A}_2\right)_{t\in I}$ on $\mathcal{B}$ over $(\mathcal{A}_1,\Phi_1)$ and $(\mathcal{A}_2,\Phi_2)$ are called {\em equivalent}, if there joint distributions coincide. This is the case, if and only if all their moments agree, i.e.\ if
\[
\Phi_1\left( j^{(1)}_{t_1}(b_1) \cdots j^{(1)}_{t_n}(b_n)\right) = \Phi_2\left( j^{(2)}_{t_1}(b_1) \cdots j^{(2)}_{t_n}(b_n)\right)
\]
holds for all $n\in\mathbb{N}$, $t_1,\ldots,t_n\in I$ and all $b_1,\ldots,b_n\in\mathcal{B}$.

The term `quantum stochastic process' is sometimes also used for an indexed
family $(X_t)_{t\in I}$ of operators on a Hilbert space or more
generally of elements of a quantum probability space. We will reserve
the name {\em operator process} for this. An operator process
$(X_t)_{t\in I}\subseteq{\mathcal A}$ always defines a quantum stochastic process
$\left(j_t:\mathbb{C}\langle a, a^*\rangle\to{\mathcal A}\right)_{t\in I}$ on the free $*$-algebra with
one generator, if we set $j_t(a)=X_t$ and extend $j_t$ as a $*$-algebra
homomorphism. On the other hand operator processes can be obtained
from quantum stochastic processes $(j_t:{\mathcal B}\to{\mathcal A})_{t\in I}$
by choosing an element $x$ of the
algebra ${\mathcal B}$ and setting $X_t=j_t(x)$.

The notion of independence we use for L\'evy processes on involutive bialgebras is the so-called tensor or boson independence. In Chapter \ref{chapter-III} we will see that other interesting notions of independence exist.

\begin{definition}\label{ch1-def-indep}
Let $({\mathcal A},\Phi)$ be a quantum probability space and ${\mathcal B}$ a
*-algebra. The quantum random variables $j_1,\ldots,j_n:{\mathcal B}\to {\mathcal A}$ are called {\em tensor or boson independent} (w.r.t.\ the state $\Phi$), if
\begin{itemize}
\item[(i)]
$\Phi\big(j_{1}(b_1) \cdots j_{n}(b_n)\big) =
\Phi\big(j_{1}(b_1)\big) \cdots \Phi\big(j_{n}(b_n)\big)$ for all
$b_1,\ldots,b_n\in{\mathcal B}$, and
\item[(ii)]
$[j_l(b_1),j_k(b_2)]=0$ for all $k \not= l$ and all $b_1,b_2\in{\mathcal B}$.
\end{itemize}
\end{definition}

Recall that an {\em involutive bialgebra} $(\mathcal{B},\Delta,\varepsilon)$ is a unital $*$-algebra $\mathcal{B}$ with two $*$-homomorphisms $\Delta:\mathcal{B}\to\mathcal{B}\otimes\mathcal{B}$, $\varepsilon:\mathcal{B}\to\mathbb{C}$ called {\em coproduct} and {\em counit}, satisfying
\begin{gather*}
({\rm id}\otimes \Delta)\circ \Delta = (\Delta\otimes{\rm id})\circ \Delta \qquad \text{(coassociativity)} \\
({\rm id}\otimes \varepsilon)\circ \Delta = {\rm id} = (\varepsilon\otimes {\rm id})\circ \Delta \qquad \text{(counit property)}. 
\end{gather*}
Let $j_1,j_2:\mathcal{B}\to\mathcal{A}$ be two linear maps with values in some algebra $\mathcal{A}$, then we define their {\em convolution} $j_1\star j_2$ by
\[
j_1\star j_2 = m_\mathcal{A} \circ (j_1\otimes j_2) \circ \Delta.
\]
If $j_1$ and $j_2$ are two independent quantum random variables, then $j_1\star j_2$ is again a quantum random variable, i.e.\ a $*$-homomorphism. The fact that we can compose quantum random variables allows us to define L\'evy processes, i.e.\ processes with independent and stationary increments.

\begin{definition}\label{I-levy-def}
Let ${\mathcal B}$ be an involutive bialgebra. A quantum stochastic process
$(j_{st})_{0\le s\le t}$ on
${\mathcal B}$ over some quantum probability space $({\mathcal A},\Phi)$ is
called a {\em L{\'e}vy process}, if the
following four conditions are satisfied.
\begin{enumerate}
\item
(Increment property) We have
\begin{eqnarray*}
j_{rs}\star j_{st} &=& j_{rt} \quad \mbox{ for all } 0\le r\le s\le
t, \\
j_{tt} &=& \mathbf{1}\circ \varepsilon \quad \mbox{ for all } 0\le t.
\end{eqnarray*}
\item
(Independence of increments)
The family $(j_{st})_{0\le s\le t}$ is independent, i.e.\ the quantum random variables $j_{s_1t_1},\ldots,j_{s_nt_n}$ are independent for all $n\in\mathbb{N}$ and all $0\le s_1\le t_1\le s_2\le \cdots \le t_n$.
\item
(Stationarity of increments)
The distribution $\varphi_{st}=\Phi\circ j_{st}$ of $j_{st}$ depends
only on the difference $t-s$.
\item
(Weak continuity) The quantum random variables $j_{st}$ converge to
$j_{ss}$ in distribution for $t\searrow s$.
\end{enumerate}
\end{definition}

Recall that an (involutive) Hopf algebra $(\mathcal{B},\Delta,\varepsilon,S)$ is an (involutive) bialgebra $(\mathcal{B},\Delta,\varepsilon)$ equipped with a linear map called {\em antipode} $S:\mathcal{B}\to\mathcal{B}$ satisfying
\begin{equation}\label{I-antipode}
S\star {\rm id} = \mathbf{1}\circ \varepsilon = {\rm id}\star S.
\end{equation}
The antipode is unique, if it exists. Furthermore, it is an algebra and coalgebra anti-homomorphism, i.e.\ it satisfies $S(ab)=S(b)S(a)$ for all $a,b\in\mathcal{B}$ and $(S\otimes S)\circ\Delta= \tau\circ\Delta\circ S$, where $\tau:\mathcal{B}\otimes\mathcal{B}\to\mathcal{B}\otimes\mathcal{B}$ is the flip $\tau(a\otimes b)=b\otimes a$.
If $(\mathcal{B},\Delta,\varepsilon)$ is an involutive bialgebra and $S:\mathcal{B}\to\mathcal{B}$ a linear map satisfying \eqref{I-antipode}, then $S$ satisfies also the relation
\[
S\circ *\circ S\circ * = {\rm id},
\]
in particular, it follows that $S$ is invertible.

One can show that if $(j_t)_{t\ge 0}$ is any quantum stochastic process on an involutive Hopf algebra, then the quantum stochastic process defined by
\[
j_{st}=m_{\mathcal{A}}\circ\big((j_s\circ S)\otimes j_t\big)\circ\Delta,
\]
for $0\le s\le t$, satisfies the increment property (1) in Definition \ref{I-levy-def}. A one-parameter stochastic process $(j_t)_{t\ge 0}$ on a Hopf $*$-algebra $H$ is called a {\em L\'evy process on $H$}, if its {\em increment process} $(j_{st})_{0\le s\le t}$ with $j_{st}=\big((j_s\circ S)\otimes j_t)\circ\Delta$ is a L\'evy process on $H$ in the sense of Definition \ref{I-levy-def}.

Let $(j_{st})_{0\le s\le t}$ be a L\'evy process on some involutive bialgebra. We will denote the marginal distributions of $(j_{st})_{0\le s\le t}$ by $\varphi_{t-s}=\Phi\circ j_{st}$. Due to the stationarity of the increments this is well defined.

\begin{lemma}
The marginal distributions $(\varphi_t)_{t\ge 0}$ of a L\'evy process on an involutive bialgebra $\mathcal{B}$ form a convolution semigroup of states on $\mathcal{B}$, i.e.\ they satisfy
\begin{enumerate}
\item
$\varphi_0=\varepsilon$, $\varphi_s\star\varphi_t=\varphi_{s+t}$ for all $s,t\ge 0$, and $\lim_{t\searrow 0}\varphi_t(b)=\varepsilon(B)$ for all $b\in\mathcal{B}$, and 
\item
$\varphi_t(\mathbf{1})=1$, and $\varphi_t(b^*b)\ge 0 $ for all $t\ge 0$ and all $b\in\mathcal{B}$.
\end{enumerate}
\end{lemma}
\begin{proof}
$\varphi_t=\Phi\circ j_{0t}$ is clearly a state, since $j_{0t}$ is a $*$-homomorphism and $\Phi$ a state.

From the first condition in Definition \ref{I-levy-def} we get
\[
\varphi_{0}=\Phi\circ j_{00} = \Phi(1)\varepsilon = \varepsilon,
\]
and
\[
\varphi_{s+t}(b) = \Phi\big(j_{0,s+t}(b)\big)=\Phi\left(\sum j_{0s}(b_{(1)})j_{s,s+t}(b_{(2)})\right),
\]
for $b\in\mathcal{B}$, $\Delta(b)=\sum b_{(1)}\otimes b_{(2)}$. Using the independence of increments, we can factorize this and get
\begin{eqnarray*}
\varphi_{s+t}(b)&=&\sum \Phi\big(j_{0s}(b_{(1)})\big)\Phi\big(j_{s,s+t}(b_{(2)})\big) = \sum \varphi_s(b_{(1)})\varphi_t(b_{(2)}) \\
&=& \varphi_s\otimes\varphi_t\big(\Delta(b)\big) = \varphi_s\star\varphi_t(b)
\end{eqnarray*}
for all $\in\mathcal{B}$. 

The continuity is an immediate consequence of the last condition in Definition \ref{I-levy-def}.
\end{proof}

\begin{lemma}\label{I-lemma-semigroups}
The convolution semigroup of states characterizes a L\'evy process on an involutive bialgebra up to equivalence.
\end{lemma}
\begin{proof}
This follows from the fact that the increment property and the independence of increments allow to express all joint moments in terms of the marginals. E.g., for $0\le s\le u\le v$ and $a,b,c\in\mathcal{B}$, the moment $\Phi\big( j_{su}(a)j_{st}(b)j_{sv}(c)\big)$ becomes
\begin{eqnarray*}
\Phi\big( j_{su}(a)j_{st}(b)j_{sv}(c)\big)&=&\Phi\big( (j_{st}\star j_{tu})(a) j_{st}(b) (j_{st}\star j_{tu}\star j_{uv})(c)\big) \\
&=& \Phi\big( j_{st}(a_{(1)}) j_{tu}(a_{(2)}) j_{st}(b) j_{st}(c_{(1)} j_{tu}(c_{(2)} j_{uv}(c_{(3)})\big) \\
&=&  \Phi\big( j_{st}(a_{(1)}bc_{(1)}) j_{tu}(a_{(2)}c_{(2)}) j_{uv}(c_{(3)})\big) \\
&=& \varphi_{t-s}(a_{(1)}bc_{(1)}) \varphi_{u-t}(a_{(2)}c_{(2)}) \varphi_{v-u}(c_{(3)}.
\end{eqnarray*}
\end{proof}
It is possible to reconstruct process $(j_{st})_{0\le s\le t}$ from its convolution semigroups, see \cite[Section 1.9]{schuermann93} or \cite[Section 4.5]{franz+schott99}. Therefore, we even have a one-to-one correspondence between L\'evy processes on $\mathcal{B}$ and convolution semigroups of states on $\mathcal{B}$.

\section{The generator and the Sch\"urmann triple of a L\'evy process}\label{I-gen}

In this section we will meet two more objects that classify L\'evy processes, namely their generator and their triple (called Sch\"urmann triple by P.-A.\ Meyer, see \cite[Section VII.1.6]{meyer95}).

We begin with a technical lemma.
\begin{lemma}\label{I-gen-lem}
\begin{description}
\item[(a)]
Let $\psi:\mathcal{C}\to\mathbb{C}$ be a linear functional on  some coalgebra $\mathcal{C}$. Then the series
\[
\exp \psi (b) \stackrel{{\rm def}}{=} \sum_{n=0} \frac{\psi^{\star n}}{n!}(b) = \varepsilon (b)+ \psi(b) + \frac{1}{2} \psi\star\psi(b) + \cdots
\]
converges for all $b\in\mathcal{C}$.
\item[(b)]
Let $(\varphi_t)_{t\ge 0}$ be a convolution semigroup on some coalgebra $\mathcal{C}$. Then the limit
\[
L(b) = \lim_{t\searrow 0} \frac{1}{t}\big(\varphi_t(b)-\varepsilon(b)\big)
\]
exists for all $b\in\mathcal{C}$. Furthermore we have $\varphi_t = \exp tL$ for all $t\ge 0$.
\end{description}
\end{lemma}

The proof of this Lemma uses the fundamental theorem of coalgebras \cite{sweedler69}, which states that for any element $c$ of a coalgebra $\mathcal{C}$ there exists a finite-dimensional subcoalgebra $\mathcal{C}_b$ of $\mathcal{C}$ containing $b$.

\begin{proposition}\label{I-schoenberg}
{\bf (Schoenberg correspondence)}
Let $\mathcal{B}$ be an involutive bialgebra, $(\varphi_t)_{t\ge 0}$ a convolution semigroup of linear functionals on $\mathcal{B}$ and 
\[
L=\lim_{t\searrow 0} \frac{1}{t}\big(\varphi_t-\varepsilon\big).
\]
Then the following are equivalent.
\begin{description}
\item[(i)]
$(\varphi_t)_{t\ge 0}$ is a convolution semigroup of states.
\item[(ii)]
$L:\mathcal{B}\to \mathbb{C}$ satisfies $L(\mathbf{1})=0$, and it is hermitian and conditionally positive, i.e.
\[
L(b^*)=\overline{L(b)}
\]
for all $b\in\mathcal{B}$, and
\[
L(b^*b)\ge 0
\]
for all $b\in\mathcal{B}$ with $\varepsilon(b)=0$.
\end{description}
\end{proposition}
\begin{proof}
We prove only the (easy) direction (i)$\Rightarrow$(ii), the converse follows from the representation theorem \ref{I-rep-thm}, whose proof can be found in  \cite[Chapter 2]{schuermann93}.

The first property follows by differentiating $\varphi_t(\mathbf{1})=1$ w.r.t.\ $t$.

Let $b\in\mathcal{B}$, $\varepsilon(b)=0$. If all $\varphi_t$ are states, then we have $\varphi_t(b^*b)\ge 0$ for all $t\ge 0$ and therefore
\[
L(b^*b)=\lim_{t\searrow 0} \frac{1}{t}\big(\varphi_t(b^*b)-\varepsilon(b^*b)\big)= \lim_{t\searrow 0} \frac{\varphi_t(b^*b)}{t}\ge 0.
\]
Similarly, $L$ is hermitian, since all $\varphi_t$ are hermitian.
\end{proof}

We will call a linear functional satisfying condition (ii) of the preceding Proposition a {\em generator}. Lemma \ref{I-gen-lem} and Proposition \ref{I-schoenberg} show that L\'evy processes can also be characterized by their generator $L=\left.\frac{{\rm d}}{{\rm d}t}\right|_{t=0} \varphi_t$.

Let $D$ be a pre-Hilbert space. Then we denote by $\mathcal{L}(D)$ the set of all linear operators on $D$ that have an adjoint defined everywhere on $D$, i.e.\ \[
\mathcal{L}(D)=\left\{X:D\to D \mbox{ linear}\left|\begin{array}{c} \mbox{ there exists }X^*: D\to D \mbox{ linear s.t.} \\ \langle u,Xv\rangle = \langle X^*u,v\rangle\mbox{ for all }u,v\in D\end{array}\right.\right\}.
\]
$\mathcal{L}(D)$ is clearly a unital $*$-algebra.

\begin{definition}
Let $\mathcal{B}$ be a unital $*$-algebra equipped with a unital hermitian character $\varepsilon:\mathcal{B}\to\mathbb{C}$ (i.e.\ $\varepsilon(\mathbf{1})=1$, $\varepsilon(b^*)=\overline{\varepsilon(b)}$, and $\varepsilon(ab)=\varepsilon(a)\varepsilon(b)$ for all $a,b\in\mathcal{B}$). A {\em Sch\"urmann triple on $(\mathcal{B},\varepsilon)$} is a triple $(\rho,\eta,L)$ consisting of
\begin{itemize}
\item
a unital $*$-representation $\rho:\mathcal{B}\to \mathcal{L}(D)$ of $\mathcal{B}$ on some pre-Hilbert space $D$,
\item
a $\rho$-$\varepsilon$-1-cocycle $\eta:\mathcal{B}\to D$, i.e.\ a linear map $\eta:\mathcal{B}\to D$ such that
\begin{equation}\label{I-eta-cocycle}
\eta(ab)=\rho(a)\eta(b)+\eta(a)\varepsilon(b)
\end{equation}
for all $a,b\in\mathcal{B}$,
\item
and a hermitian linear functional $L:\mathcal{B}\to\mathbb{C}$ that has the linear map $\mathcal{B}\times\mathcal{B}\ni(a,b)\mapsto - \langle\eta(a^*),\eta(b)\rangle$ as a $\varepsilon$-$\varepsilon$-2-coboundary, i.e.\ that satisfies
\begin{equation}\label{I-L-coboundary}
 - \langle\eta(a^*),\eta(b)\rangle = \partial L(a,b) = \varepsilon(a)L(b)-L(ab)+L(a)\varepsilon(b)
\end{equation}
for all $a,b\in\mathcal{B}$.
\end{itemize}
We will call a Sch\"urmann triple {\em surjective}, if the cocycle $\eta:\mathcal{B}\to D$ is surjective.
\end{definition}

\begin{theorem}\label{I-1-1-corresp}
Let $\mathcal{B}$ be an involutive bialgebra. We have one-to-one correspondences between L\'evy processes on $\mathcal{B}$ (modulo equivalence), convolution semigroups of states on $\mathcal{B}$, generators on  $\mathcal{B}$, and surjective Sch\"urmann triples on $\mathcal{B}$ (modulo unitary equivalence).
\end{theorem}
\begin{proof}
It only remains to establish the one-to-one correspondence between generators and Sch\"urmann triples.

Let $(\rho,\eta,L)$ be a Sch\"urmann triple, then we can show $L$ is a generator, i.e.\ a hermitian, conditionally positive linear functional with $L(\mathbf{1})=0$.

The cocycle has to vanish on the unit element $\mathbf{1}$, since
\[
\eta(\mathbf{1})=\eta(\mathbf{1}\cdot\mathbf{1})=\rho(\mathbf{1})\eta(\mathbf{1})+\eta(\mathbf{1})\varepsilon(\mathbf{1}) = 2\eta(\mathbf{1}).
\]
This implies
\[
L(\mathbf{1})=L(\mathbf{1}\cdot\mathbf{1}) = \varepsilon(\mathbf{1})L(\mathbf{1}) + \langle\eta(\mathbf{1}),\eta(\mathbf{1})\rangle+L(\mathbf{1})\varepsilon(\mathbf{1})=2L(\mathbf{1})=0.
\]
Furthermore, $L$ is hermitian by definition and conditionally positive, since by \eqref{I-L-coboundary} we get
\[
L(b^*b)=\langle\eta(b),\eta(b)\rangle =||\eta(b)||^2\ge 0
\]
for $b\in{\rm ker}\,\varepsilon$.

Let now $L$ be a generator. The sesqui-linear form $\langle\cdot,\cdot\rangle_L:\mathcal{B}\times\mathcal{B}\to\mathbb{c}$ defined by
\[
\langle a,b\rangle_L = L\Big(\big(a-\varepsilon(a)\mathbf{1}\big)^*\big(b-\varepsilon(b)\mathbf{1}\big)\Big)
\]
for $a,b\in\mathcal{B}$ is positive semi-definite, since $L$ is conditionally positive. Dividing $\mathcal{B}$ by the null-space
\[
\mathcal{N}_L=\{ a\in\mathcal{B}|\langle a,a\rangle_L=0\}
\]
we obtain a pre-Hilbert space $D=\mathcal{B}/\mathcal{N}_L$ with a positive definite inner product $\langle\cdot,\cdot\rangle$ induced by $\langle\cdot,\cdot\rangle_L$. For the cocycle $\eta:\mathcal{B}\to D$ we take the canonical projection, this is clearly surjective and satisfies Equation \eqref{I-L-coboundary}.

The $*$-representation $\rho$ is induced from the left multiplication on $\mathcal{B}$ on ${\rm ker}\,\varepsilon$, i.e.\
\[
\rho(a)\eta\big(b-\varepsilon(b)\mathbf{1}\big)=\eta\Big(a\big(b-\varepsilon(b)\mathbf{1}\big)\Big)
\qquad\mbox{ or }\qquad
\rho(a)\eta(b) = \eta(ab)-\eta(a)\varepsilon(b)
\]
for $a,b\in\mathcal{B}$. To show that this is well-defined, we have to verify that left multiplication by elements of $\mathcal{B}$ leaves the null-space invariant. Let therefore $a,b\in\mathcal{B}$, $b\in\mathcal{N}_L$, then we have
\begin{eqnarray*}
\left|\left|\Big(a\big(b-\varepsilon(b)\mathbf{1}\big)\Big)\right|\right|^2 &=& L\Big(\big(ab-a\varepsilon(b)\mathbf{1}\big)^*\big(ab-a\varepsilon(b)\mathbf{1}\big)\Big) \\
&=& L\Big(\big(b-\varepsilon(b)\mathbf{1}\big)^*a^*\big(ab-a\varepsilon(b)\mathbf{1}\big)\Big) = \big\langle b-\varepsilon(b)\mathbf{1}, a^*a\big(b-\varepsilon(b)\mathbf{1}\big)\big\rangle_L \\
&\le & ||b-\varepsilon(b)\mathbf{1}||^2 \,\left|\left| a^*a\big(b-\varepsilon(b)\mathbf{1}\right|\right|^2 =0,
\end{eqnarray*}
with Schwarz' inequality.

That the Sch\"urmann triple $(\rho,\eta,L$) obtained in this way is unique up to unitary equivalence follows similarly as for the usual GNS construction.
\end{proof}

\begin{example}\label{I-exa-class-levy}
Let $(X_t)_{t\ge 0}$ be a classical real-valued L\'evy process on a probability space $(\Omega,\mathcal{F},P)$), whose characteristic function
\[
\mathbb{E}(e^{iuX_t})=\exp t\left(idu-\frac{\sigma^2}{2}u^2+\int_{\mathbb{R}\backslash\{0\}} (e^{iux}-1-iux\mathbf{1}_{|x|\le 1}){\rm d}\nu(x)\right)
\]
is analytic in a neighborhood of $0$. We will now define a L\'evy process on the free algebra $\mathbb{C}[x]$ generated by one symmetric
element $x=x^*$ with the coproduct and counit determined by
$\Delta(x)=x\otimes \mathbf{1} + \mathbf{1}\otimes x$ and
$\varepsilon(x)=0$, whose moments agree with those of $(X_t)_{t\ge
  0}$, i.e.\ such that
\[
\Phi\big(j_{st}(x^k)\big) = \mathbb{E}\left((X_t-X_s)^k\right)
\]
holds for all $k\in\mathbb{N}$ and all $0\le s\le t$. Furthermore, we will construct its Sch\"urmann triple.

We take $\mathcal{A}=Pol(X)$, i.e.\ the algebra generated by $\{X_t|t\ge 0\}$, and for the the state $\Phi$ we take the expectation functional. The L\'evy process is defined by $j_{st}(x^k)=(X_t-X_s)^k$, for $0\le s\le t$. It is straightforward to verify that this defines indeed a L\'evy process on $\mathbb{C}[x]$.

Our assumptions allow us compute the generator $L$ as
\begin{eqnarray*}
L(x^k) &=& \left.\frac{{\rm d}}{{\rm d}t}\right|_{t=0} \Phi\left(j_{0t}(x^k)\right) =\left.\frac{{\rm d}}{{\rm d}t}\right|_{t=0} \mathbb{E}(X_t^k)
=\left.\frac{{\rm d}}{{\rm d}t}\right|_{t=0} \left(\frac{1}{i^k}\left.\frac{{\rm d}^k}{{\rm d}^ku}\right|_{u=0} \mathbb{E}(e^{iuX_t})\right) \\
&=&
\left\{
\begin{array}{lcl}
0 & & k=0, \\
d+\int_{|x|>1}x{\rm d}\nu(x) & & k=1, \\
\sigma^2+\int_{\mathbb{R}\backslash\{0\}}x^2{\rm d}\nu(x) & & k=2, \\
\int_{\mathbb{R}\backslash\{0\}}x^k{\rm d}\nu(x) & & k\ge 3.
\end{array}
\right.
\end{eqnarray*}
We can define the Sch\"urmann triple on the Hilbert space $H=\mathbb{C}\oplus L^2(\nu)$. The representation $\rho$ acts as
\[
\rho(p)(\lambda,f)=\big(p(0)\lambda,pf\big)
\]
for a polynomial $p\in\mathbb{C}[x]$, and $\lambda\in\mathbb{C}$, $f\in L^2(\nu)$. The cocycle $\eta:\mathbb{C}[x]\to H$ can be written as
\[
\eta(p) =\big(\sigma p'(0),p-p(0)\big)
\]
where $p'$ denotes the first derivative of $p$. It is not difficult to check that this defines indeed a Sch\"urmann triple.
\end{example}

For the classification of Gaussian and drift generators on an involutive bialgebra $\mathcal{B}$ with counit $\varepsilon$, we need the ideals
\begin{eqnarray*}
K &=& {\rm ker}\,\varepsilon, \\
K^2 &=& {\rm span}\,\{ab|a,b\in K\}, \\
K^3 &=& {\rm span}\,\{abc|a,b,c\in K\}.
\end{eqnarray*}

\begin{proposition}
Let $L$ be a conditionally positive, hermitian linear functional on $\mathcal{B}$. Then the following are equivalent.
\begin{description}
\item[(i)]
$\eta=0$,
\item[(ii)]
$L|_{K^2}=0$,
\item[(iii)]
$L$ is an $\varepsilon$-derivation, i.e.\ $L(ab)=\varepsilon(a)L(b)+L(a)\varepsilon(b)$ for all $a,b\in\mathcal{B}$,
\item[(iv)]
The states $\varphi_t$ are homomorphisms, i.e.\ $\varphi_t(ab)=\varphi_t(a)\varphi_t(b)$ for all $a,b\in\mathcal{B}$ and $t\ge 0$.
\end{description}
If a  conditionally positive, hermitian linear functional $L$ satisfies one of these conditions, then we call it and the associated L\'evy process a {\em drift}.
\end{proposition}

\begin{proposition}\label{I-prop-def-gaussian}
Let $L$ be a conditionally positive, hermitian linear functional on $\mathcal{B}$.

Then the following are equivalent.
\begin{description}
\item[(i)]
$L|_{K^3}=0$,
\item[(ii)]
$L(b^*b)=0$ for all $B\in K^2$,
\item[(iii)]
$L(abc)=L(ab)\varepsilon(c)+L(ac)\varepsilon(b)+L(bc)\varepsilon(a) - \varepsilon(ab)L(c)-\varepsilon(ac)L(b)-\varepsilon(bc)L(a)$ for all $a,b,c\in\mathcal{B}$,
\item[(iv)]
$\rho|_{K}=0$ for the representation $\rho$ in the surjective Sch\"urmann triple $(\rho,\eta,L)$ associated to $L$ by the GNS-type construction presented in the proof of Theorem \ref{I-1-1-corresp},
\item[(v)]
$\rho=\varepsilon\mathbf{1}$, for the representation $\rho$ in the surjective  Sch\"urmann triple $(\rho,\eta,L)$ associated to $L$ by the GNS-type construction presented in the proof of Theorem \ref{I-1-1-corresp},
\item[(vi)]
$\eta|_{K^2}=0$ for the cocycle $\eta$ in any  Sch\"urmann triple $(\rho,\eta,L)$ containing $L$,
\item[(vii)]
$\eta(ab)=\varepsilon(a)\eta(b)+\eta(a)\varepsilon(b)$ for all $a,b\in\mathcal{B}$ and the cocycle $\eta$ in any  Sch\"urmann triple $(\rho,\eta,L)$ containing $L$.
\end{description}
If a  conditionally positive, hermitian linear functional $L$ satisfies one of these conditions, then we call it and also the associated L\'evy process {\em quadratic} or {\em Gaussian}.
\end{proposition}

\begin{proposition}
Let $L$ be a conditionally positive, hermitian linear functional on $\mathcal{B}$. The the following are equivalent.
\begin{description}
\item[(i)]
There exists a state $\psi:\mathcal{B}\to\mathbb{C}$ and a real number $\lambda>0$ such that
\[
L(b)=\lambda\big(\varphi(b)-\varepsilon(b)\big)
\]
for all $b\in\mathcal{B}$.
\item[(ii)]
There exists a  Sch\"urmann triple $(\rho,\eta,L)$ containing $L$, in which the cocycle $\eta$ is trivial, i.e.\ of the form
\[
\eta(b) = \big(\rho(b)-\varepsilon(b)\big)\omega, \qquad \mbox{ for all }b\in\mathcal{B},
\]
for some non-zero vector $\omega\in D$.
\end{description}
If a  conditionally positive, hermitian linear functional $L$ satisfies one of these conditions, then we call it a {\rm Poisson generator} and the associated L\'evy process a {\em compound Poisson process}.
\end{proposition}
\begin{proof}
To show that (ii) implies (i), set $\varphi(b)=\langle\omega,\rho(b)\omega\rangle$ and $\lambda=||\omega||^2$.

For the converse, let $(D,\rho,\omega)$ be the GNS triple for $(\mathcal{B},\varphi)$ and check that $(\rho,\eta,L)$ with $\eta(b) = \big(\rho(b)-\varepsilon(b)\big)\omega$, $b\in\mathcal{B}$ defined a triple.
\end{proof}

\begin{remark}
The  Sch\"urmann triple for a Poisson generator $L=\lambda(\varphi-\varepsilon)$ obtained by the GNS construction for $\varphi$ is not necessarily surjective. Consider, e.g., a classical additive $\mathbb{R}$-valued compound Poisson process, whose L\'evy measure $\mu$ is not supported on a finite set. Then the construction of a surjective  Sch\"urmann triple in the proof of Theorem \ref{I-1-1-corresp} gives the pre-Hilbert space $D_0={\rm span}\,\{x^k|k=1,2,\ldots\}\subseteq L^2(\mathbb{R},\mu)$. On the other hand, the GNS-construction for $\varphi$ leads to the pre-Hilbert space $D={\rm span}\,\{x^k|k=0,1,2,\ldots\}\subseteq L^2(\mathbb{R},\mu)$. The cocycle $\eta$ is the coboundary of the constant function, which is not contained in $D_0$.
\end{remark}

\section{The representation theorem}\label{I-rep}

The representation theorem gives a direct way to construct a L\'evy process from the  Sch\"urmann triple, using quantum stochastic calculus.

\begin{theorem}\label{I-rep-thm}
{\bf (Representation theorem)}
Let $\mathcal{B}$ be an involutive bialgebra and $(\rho,\eta,L)$ a  Sch\"urmann triple on $\mathcal{B}$. Then the quantum stochastic differential equations
\begin{equation}\label{I-rep-thm-qsde}
{\rm d}j_{st} = j_{st}\star \big({\rm d}A^*_t\circ\eta + {\rm d}\Lambda_t\circ(\rho-\varepsilon) + {\rm d}A_t\circ\eta\circ* + L{\rm d}t\big)
\end{equation}
with the initial conditions
\[
j_{ss}=\varepsilon \mathbf{1}
\]
have a solution $(j_{st})_{0\le s\le t}$.  Moreover, in the vacuum state $\Phi(\cdot)=\langle\Omega,\,\cdot\,\Omega\rangle$, $(j_{st})_{0\le s\le t}$ is a L\'evy process with generator $L$.

Conversely, every L\'evy process with generator $L$ is equivalent to $(j_{st})_{0\le s\le t}$.
\end{theorem}
For the proof of the representation theorem we refer to \cite[Chapter 2]{schuermann93}.

Written in integral form and applied to an element $b\in\mathcal{B}$ with $\Delta(b)=\sum b_{(1)}\otimes b_{(2)}$, Equation \eqref{I-rep-thm-qsde} takes the form
\[
j_{st}(b)=\varepsilon(b)\mathbf{1} + \int_s^t j_{s\tau}(b_{(1)}) \Big({\rm d}A^*_\tau\big(\eta(b_{(2)})\big) + {\rm d}\Lambda_\tau\big(\rho(b_{(2)}-\varepsilon(b_{(2)})\mathbf{1}\big) + {\rm d}A_\tau\big(\eta(b_{(2)}^*)\big) + L(b_{(2)}){\rm d}\tau\Big).
\]

One can show that
\[
{\rm d}M_t = {\rm d}A^*_t\circ\eta + {\rm d}\Lambda_t \circ (\rho - \varepsilon) + {\rm d}A_t \circ \tilde\eta + L {\rm d}t
\]
formally defines a $*$-homomorphism on ker $\varepsilon=\mathcal{B}_0$, if we define the algebra of quantum stochastic differentials (or It\^o algebra, cf.\ \cite{belavkin98} and the references therein) over some pre-Hilbert space $D$ as follows.

The algebra of quantum stochastic differentials $\mathcal{I}(D)$ over $D$ is the $*$-algebra generated by
\[
\{{\rm d}\Lambda(F)| F\in \mathcal{L}(D)\} \cup  \{{\rm d}A^*(u)| u \in D\}\cup  \{{\rm d}A(u)| u\in D\} \cup  \{{\rm d}t\},
\]
if we identify
\begin{eqnarray*}
{\rm d}\Lambda(\lambda F+\mu G) &\equiv& \lambda {\rm d}\Lambda(F)+\mu{\rm d}\Lambda(G), \\
{\rm d}A^*(\lambda u+\mu v) &\equiv& \lambda {\rm d}A^*(u)+\mu{\rm d}A^*(v), \\
{\rm d}A(\lambda u+\mu v) &\equiv &\overline{\lambda} {\rm d}A(u)+\overline{\mu}{\rm d}A(v),
\end{eqnarray*}
for all $F,G\in \mathcal{L}(D)$, $u,v\in D$, $\lambda,\mu\in\mathbb{C}$. The involution of $\mathcal{I}(D)$ is defined by
\begin{eqnarray*}
{\rm d}\Lambda(F)^* &=& {\rm d}\Lambda(F^*), \\
\Big({\rm d}A^*(u)\Big)^* &=&{\rm d}A(u), \\
{\rm d}A(u)^* &=&{\rm d}A^*(u),
\end{eqnarray*}
for $F\in \mathcal{L}(D)$, $u\in D$, and the multiplication by the It\^o table
\[
\begin{array}{|c||c|c|c|c|}
\hline
\bullet & {\rm d}A^*(u) & {\rm d}\Lambda(F) & {\rm d}A(u) & {\rm d}t \\
\hline\hline 
{\rm d}A^*(v) & 0 & 0 & 0 & 0 \\ \hline
{\rm d}\Lambda(G) & {\rm d}A^*(Gu) &  {\rm d}\Lambda(GF) & 0 & 0 \\ \hline
{\rm d}A(v) & \langle v,u\rangle {\rm d}t & {\rm d}A(F^*v) & 0 & 0 \\ \hline
{\rm d}t & 0 & 0 & 0 & 0 \\ \hline
\end{array}
\]
for all $F,G\in \mathcal{L}(D)$, $u,v\in D$, i.e.\ we have, for example,
\[
{\rm d}A(v)\bullet {\rm d}A^*(u)= \langle v,u\rangle {\rm d}t, \quad\mbox{ and }\quad   {\rm d}A^*(u)\bullet {\rm d}A(v)=0.
\]

\begin{proposition}
Let $(j_{st})_{0\le s\le t}$ be a L\'evy process on a $*$-bialgebra $\mathcal{B}$ with  Sch\"urmann triple $(\rho,\eta,L)$, realized on the Fock space $\Gamma\big(L^2(\mathbb{R}_+,D)\big)$ over the pre-Hilbert space $D$. Let furthermore $u$ be a unitary operator on $D$ and $\omega\in D$. Then the quantum stochastic differential equation
\[
{\rm d}U_t =U_t\left({\rm d}A_t(\omega) - {\rm d}A^*_t(u\omega) + {\rm d}\Lambda_t(u-\mathbf{1}) - \frac{||\omega||^2}{2}{\rm d}t\right)
\]
with the initial condition $U_0=\mathbf{1}$ has a unique solution $(U_t)_{t\ge 0}$ with $U_t$ a unitary for all $t\ge 0$.

Furthermore, the quantum stochastic process $(\tilde{\jmath}_{st})_{0\le s\le t}$ defined by
\[
\tilde{\jmath}_{st}(b)=U^*_tj_{st}(b)U_t, \qquad\mbox{ for }b\in\mathcal{B},
\]
is again a L\'evy process with respect to the vacuum state. The  Sch\"urmann triple $(\tilde{\rho},\tilde{\eta},\tilde{L})$ of $(\tilde{\jmath}_{st})_{0\le s\le t}$ is given by
\begin{eqnarray*}
\tilde{\rho}(b) &=& u^*\rho(b) u, \\
\tilde{\eta}(b) &=& u^*\eta(b) - u^*\rho(b)u\omega, \\
\tilde{L}(b) &=& L(b)-\langle u\omega,\eta(b)\rangle+\langle\eta(b),u\omega\rangle + \langle u\omega,\rho(b)u\omega\rangle \\
&=& L(b)-\langle \omega,\tilde{\eta}(b)\rangle+\langle\tilde{\eta}(b),\omega\rangle - \langle \omega,\tilde{\rho}(b)\omega\rangle
\end{eqnarray*}
\end{proposition}

On exponential vectors the operator process $(U_t)_{t\ge 0}$ can be given by
\[
U_t=e^{-A^*_t(u\omega)}\Gamma_t(u)e^{A_t(\omega)}e^{-t||\omega||^2/2},
\]
where $\Gamma_t(u)$ denotes the second quantization of $u$. $(U_t)_{t\ge0}$ is a unitary local cocycle, cf.\ \cite{lindsay03,bhat03}.

Setting
\[
k_t(x)=U_t
\]
and extending this as a $*$-homomorphism, we get a L\'evy process on the group algebra $\mathcal{A}=\mathbb{CZ}$. The algebra $\mathcal{A}$ can be regarded as the $*$-algebra generated by a unitary $x$, i.e.\ $\mathbb{CZ}\cong\mathbb{C}\langle x,x^*\rangle/\langle xx^*-\mathbf{1},x^*x-\mathbf{1}\rangle$. Its Hopf algebra structure is given by
\[
\varepsilon(x)=1,\qquad \Delta(x)=x\otimes x, \qquad S(x)=x^*.
\]

Now it is straightforward to verify that $(\tilde{\jmath}_{st})_{0\le s\le t}$ is a L\'evy process, using the information on $(U_t)_{t\ge 0}$ we have due to the fact that it is a local unitary cocycle or a L\'evy process.

Using the quantum It\^o formula, one can then show that $(\tilde{\jmath}_{st})_{0\le s\le t}$ satisfies the quantum stochastic differential equation
\[
{\rm d}\tilde{\jmath}_{st} = j_{st}\star \big({\rm d}A^*_t\circ\tilde{\eta} + {\rm d}\Lambda_t\circ(\tilde{\rho}-\varepsilon) + {\rm d}A_t\circ\tilde{\eta}\circ L + \tilde{L}{\rm d}t\big)
\]
with initial condition $\tilde{\jmath}_{ss}=\varepsilon\mathbf{1}$, and deduce that $(\tilde{\rho},\tilde{\eta},\tilde{L})$ is a  Sch\"urmann triple for $(\tilde{\jmath}_{st})_{0\le s\le t}$.

\begin{corollary}
If the cocycle $\eta$ is trivial, then $(j_{st})_{0\le s\le t}$ is cocycle conjugate to the second quantization $\big(\Gamma_{st}(\rho)\big)_{0\le s\le t}$ of $\rho$.
\end{corollary}

\section{Cyclicity of the vacuum vector}\label{I-cyclic}

Recently, Franz and Skeide \cite{franz+schuermann+skeide03} have shown that the vacuum vector is cyclic for the realization of a L\'evy process over the Fock space given by Theorem \ref{I-rep-thm}, if the cocycle is surjective.

\begin{theorem}\label{I-cyclic-thm}
Let $(\rho,\eta,L)$ be a surjective Sch\"urmann triple on an involutive bialgebra $\mathcal{B}$ over $D$ and let $(j_{st})_{0\le s\le t}$ be the solution of Equation \eqref{I-rep-thm-qsde} on the Fock space $\Gamma\big(L^2(\mathbb{R}_+,D)\big)$. Then the vacuum vector $\Omega$ is cyclic for $(j_{st})_{0\le s\le t}$, i.e.\
\[
{\rm span}\{j_{s_1t_1}(b_1)\cdots j_{s_nt_n}(b_n)\Omega|n\in\mathbb{N},0\le s_1\le t_1\le s_2\le\cdots\le t_n,b_1,\ldots,b_n\in\mathcal{B}\}
\]
is dense in $\Gamma\big(L^2(\mathbb{R}_+,D)\big)$.
\end{theorem}

The proof which we will present here is due to Skeide. It uses the fact that the exponential vectors of indicator functions form a total subset of the Fock space.

\begin{theorem}\label{I-total-thm}
 \cite{parthasarathy+sunder98,skeide00}
Let $\eufrak{h}$ be a Hilbert space and $B\in\eufrak{h}$ a total subset of $\eufrak{h}$. Let furthermore $\mathcal{R}$ denote the ring generated by bounded intervals in $\mathbb{R}_+$. Then
\[
\{\mathcal{E}(v\mathbf{1}_I)|v\in B,I\in\mathcal{R}\}
\]
is total in $\Gamma\big(L^2(\mathbb{R}_+,\eufrak{h})\big)$.
\end{theorem}

We first show how exponential vectors of indicator functions of intervals can be generated from the vacuum vector.

\begin{lemma}
Let $0\le s\le t$ and $b\in {\rm ker}\,\varepsilon$. For $n\in\mathbb{N}$, we define
\[
\Pi^n_{[s,t]}(b)= j_{s,s+\delta}(\mathbf{1}+b)j_{s+\delta,s+2\delta}(\mathbf{1}+b)\cdots j_{t-\delta,t}(\mathbf{1}+b)e^{-(t-s)L(b)},
\]
where $\delta=(t-s)/n$. Then $\Pi^n_{[s,t]}(b)\Omega$ converges to the exponential vector $\mathcal{E}\big(\eta(b)\mathbf{1}_{[s,t]}\big)$
\end{lemma}

\begin{proof}
Let $b\in\mathcal{B}$ and $k\in D$. Then the fundamental lemma of quantum stochastic calculus implies
\[
\langle\mathcal{E}(k\mathbf{1}_{[0,T]}),j_{st}(b)\Omega\rangle=\varepsilon(B)+\int_s^t\langle\mathcal{E}(k\mathbf{1}_{[0,T]}),j_{s\tau}(b_{(1)})\Omega\rangle\big(\langle k,\eta(b_{(2)})\rangle + L(b_{(2)})\big){\rm d}\tau
\]
for $0\le s\le t\le T$. This is an integral equation for a linear functional on $\mathcal{B}$, it has a unique solution given by the convolution exponential
\[
\langle\mathcal{E}(k\mathbf{1}_{[0,T]}),j_{st}(b)\Omega\rangle=\exp_\star(t-s)\big(\langle k,\eta(b)\rangle + L(b)\big).
\]
Let $b\in{\rm ker}\,\varepsilon$, then we have
\[
\langle\mathcal{E}(k\mathbf{1}_{[0,T]}),j_{st}(\mathbf{1}+b)e^{-(t-s)L(b)}\Omega\rangle = 1+(t-s)\langle k,\eta(b)\rangle + O\left(t-s)^2\right).
\]
for all $0\le s\le t\le T$.

Furthermore, we have
\begin{gather*}
\langle j_{st}(\mathbf{1}+b)e^{-(t-s)L(b)}\Omega,j_{st}(\mathbf{1}+b)e^{-(t-s)L(b)}\Omega\rangle \\
= \langle \Omega,j_{st}\big((\mathbf{1}+b)^*(\mathbf{1}+b)\big)e^{-(t-s)\left(L(b)+L(b^*)\right)}\Omega\rangle \\
= \big(1 +\varphi_{t-s}(b^*)+\varphi_{t-s}(b)+\varphi_{t-s}(b^*b)\big)e^{-(t-s)\left(L(b)+L(b^*)\right)}
\end{gather*}
for $b\in {\rm ker}\,\varepsilon$, and therefore
\[
\langle j_{st}(\mathbf{1}+b)e^{-(t-s)L(b)}\Omega,j_{st}(\mathbf{1}+b)e^{-(t-s)L(b)}\Omega\rangle = 1 + (t-s)\langle\eta(b),\eta(b)\rangle + O\left((t-s)^2\right).
\]
These calculations show that $\Pi^n_{[s,t]}(b)\Omega$ converges in norm to the exponential vectors $\mathcal{E}\big(\eta(b)\mathbf{1}_{[s,t]}\big)$, since using the independence of increments of $(j_{st})_{0\le s\le t}$, we get
\begin{eqnarray*}
\left|\left|\Pi^n_{[s,t]}(b)\Omega-\mathcal{E}\big(\eta(b)\mathbf{1}_{[s,t]}\big)\right|\right|^2 &=& \langle \Pi^n_{[s,t]}(b)\Omega,\Pi^n_{[s,t]}(b)\Omega\rangle - \langle \Pi^n_{[s,t]}(b)\Omega,\mathcal{E}\big(\eta(b)\mathbf{1}_{[s,t]}\big)\rangle \\
&&- \langle\mathcal{E}\big(\eta(b)\mathbf{1}_{[s,t]}\big),\Pi^n_{[s,t]}(b)\Omega\rangle + \langle\mathcal{E}\big(\eta(b)\mathbf{1}_{[s,t]}\big),\mathcal{E}\big(\eta(b)\mathbf{1}_{[s,t]}\big)\rangle \\[1mm]
&=& \big(1+\delta||\eta(b)||^2+O(\delta^2)\big)^n-e^{(t-s)||\eta(n)||^2} \\[1mm]
&\stackrel{n\to\infty}{\rightarrow}& 0.
\end{eqnarray*}
\end{proof}

\begin{proof} (of Theorem \ref{I-cyclic-thm})
We can generate exponential vectors of the form $\mathcal{E}(v\mathbf{1}_I)$, with $I=I_1\cup\cdots\cup I_k\in\mathcal{R}$ a union of disjoint intervals by taking products
\[
\Pi^n_I(b)=\Pi_{I_1}^n(b)\cdots \Pi_{I_k}^n(b)
\]
with an element $b\in{\rm ker}\,\varepsilon$, $\eta(b)=v$. If $\eta$ is surjective, then it follows from Theorem \ref{I-total-thm} that we can generate a total subset from the vacuum vector.
\end{proof}

If the L\'evy process is defined on a Hopf algebra, then it is sufficient to consider time-ordered products of increments corresponding to intervals starting at $0$.

\begin{corollary}\label{I-cor-cyclic-hopf}
Let $H$ be a Hopf algebra with antipode $S$. Let furthermore $(\rho,\eta,L)$ be a surjective Sch\"urmann triple on $H$ over $D$ and $(j_{st})_{0\le s\le t}$ the solution of Equation \eqref{I-rep-thm-qsde} on the Fock space $\Gamma\big(L^2(\mathbb{R}_+,D)\big)$. Then the subspaces
\begin{eqnarray*}
\mathcal{H}_\uparrow &=&{\rm span}\{j_{0t_1}(b_1)\cdots j_{0t_n}(b_n)\Omega|n\in\mathbb{N}, 0\le t_1\le t_2\le\cdots\le t_n,b_1,\ldots,b_n\in H\}, \\
\mathcal{H}_\downarrow &=&{\rm span}\{j_{0t_n}(b_1)\cdots j_{0t_1}(b_n)\Omega|n\in\mathbb{N}, 0\le t_1\le t_2\le\cdots\le t_n,b_1,\ldots,b_n\in H\},
\end{eqnarray*}
are dense in $\Gamma\big(L^2(\mathbb{R}_+,D)\big)$.
\end{corollary}

\begin{remark}
Let $(\rho,\eta,L)$ be an arbitrary Sch\"urmann triple on some involutive bialgebra $\mathcal{B}$ and let $(j_{st})_{0\le s\le t}$ be the solution of Equation \eqref{I-rep-thm-qsde} on the Fock space $\Gamma\big(L^2(\mathbb{R}_+,D)\big)$.
Then we have $\mathcal{H}_\uparrow\subseteq \mathcal{H}_0$ and $\mathcal{H}_\downarrow\subseteq \mathcal{H}_0$ for the subspaces $\mathcal{H}_\uparrow,\mathcal{H}_\downarrow,\mathcal{H}_0\subseteq\Gamma\big(L^2(\mathbb{R}_+,D)\big)$ defined as in Theorem \ref{I-cyclic-thm} and Corollary \ref{I-cor-cyclic-hopf}. This follows since any product $j_{s_1t_1}(b_1)\cdots j_{s_nt_n}(b_n)$ with arbitrary bounded intervals $[s_1,t_1],\ldots [s_n,t_n]\subseteq\mathbb{R}_+$ can be decomposed in a linear combination of products with disjoint intervals, see the proof of Lemma \ref{I-lemma-semigroups}.

E.g., for $j_{0s}(a)j_{0t}(b)$, $a,b\in\mathcal{B}$, $0\le s\le t$, we get
\[
j_{0s}(a)j_{0t}(b) = j_{0s}(ab_{(1)})j_{st}(b_{(2)})
\]
where $\Delta(b)=b_{(1)}\otimes b_{(2)}$.
\end{remark}

\begin{proof}
The density of $\mathcal{H}_\uparrow$ follows, if we show $\mathcal{H}_\uparrow=\mathcal{H}_0$. This is clear, if we show that the map $T_1:H\otimes H\to H\otimes H$, $T_1=(m\otimes {\rm id})\circ ({\rm id}\otimes \Delta)$, i.e., $T_1(a\otimes b)=ab_{(1)}\otimes b_{(2)}$ is a bijection, since
\begin{gather*}
j_{0t_1}(b_1)\cdots j_{0t_n}(b_n) \\
= m_{\mathcal{A}}^{(n-1)}\circ (j_{0t_1}\otimes j_{t_1t_2}\otimes \cdots \otimes j_{t_{n-1}t_n})\Big(\big(b_1\otimes1^{\otimes(n-1)}\big)\big(\Delta(b_2)\otimes 1^{\otimes(n-2)}\big)\cdots\big(\Delta^{(n-1)}\big)\Big) \\
=m_{\mathcal{A}}^{(n-1)} \circ (j_{0t_1}\otimes j_{t_1t_2}\otimes \cdots \otimes j_{t_{n-1}t_n})\circ T^{(n)}_1(b_1\otimes\cdots\otimes b_n),
\end{gather*}
where
\[
T_1^{(n)}=(T_1\otimes{\rm id}_{H^{\otimes (n-2)}})\circ({\rm id}_H\otimes T_1\otimes{\rm id}_{H^{\otimes (n-3)}})\circ\cdots\circ({\rm id}_{H^{\otimes(n-2)}}\otimes T_1)
\]
see also \cite[Section 4.5]{franz+schott99}. To prove that $T_1$ is bijective, we give an explicit formula for its inverse,
\[
T_1^{-1} = (m\otimes{\rm id})\circ({\rm id}\otimes S\otimes{\rm id})\circ({\rm id}\otimes\Delta).
\]

To show $\mathcal{H}_\downarrow=\mathcal{H}_0$ it is sufficient to show that the map $T_2:H\otimes H\to H\otimes H$, $T_2=(m\otimes{\rm id})\circ({\rm id}\otimes\tau)\circ(\Delta\otimes{\rm id})$, $T_2(a\otimes b)=a_{(1)}b\otimes a_{(2)}$ is bijective. This follows from the first part of the proof, since $T_1=(*\otimes *)\circ T_2\circ (*\otimes*)$.
\end{proof}

The following simple lemma is useful for checking if a Gaussian Sch\"urmann triple is surjective.
\begin{lemma}\label{I-lemma-surj}
Let $(\rho,\eta,L)$ be a Gaussian Sch\"urmann triple on a $*$-bialgebra $\mathcal{B}$ and let $G\subseteq\mathcal{B}$ be a set of algebraic generators, i.e.
\[
{\rm span}\{g_1\cdots g_n|n\in\mathbb{N},g_1,\ldots,g_n\in G\}=\mathcal{B}.
\]
Then we have
\[
{\rm span}\,\eta(G)=\eta(\mathcal{B}).
\]
\end{lemma}
\begin{proof}
For Gaussian Sch\"urmann triples one can show by induction over $n$,
\[
\eta(g_1\cdots g_n) = \sum_{k=1}^n \varepsilon(g_1\cdots g_{k-1}g_{k+1}\cdots g_n)\eta(g_k).
\]
\end{proof}

\section{Examples}\label{I-exa}

\subsection{Additive L\'evy processes}

For a vector space $V$ the tensor algebra $\mathcal{T}(V)$ is the vector space
\[
\mathcal{T}(V)=\bigoplus_{n\in\mathbb{N}} V^{\otimes n},
\]
where $V^{\otimes n}$ denotes the $n$-fold tensor product of $V$ with itself, $V^{\otimes 0}=\mathbb{C}$, with the multiplication given by
\[
(v_1\otimes \cdots \otimes v_n)(w_1\otimes\cdots\otimes w_m)=v_1\otimes \cdots \otimes v_n\otimes w_1\otimes\cdots\otimes w_m,
\]
for $n,m\in\mathbb{N}$, $v_1,\ldots,v_n,w_1,\ldots,w_m\in V$. The elements of $\bigcup_{n\in\mathbb{N}} V^{\otimes n}$ are called homogeneous, and the degree of a homogeneous element $a\not=0$ is $n$ if $a\in V^{\otimes n}$. If $\{v_i|i\in I\}$ is a basis of $V$, then the tensor algebra $\mathcal{T}(V)$ can be viewed as the free algebra generated by $v_i$, $i\in I$. The tensor algebra can be characterized by the following universal property.

There exists an embedding $\imath:V\to\mathcal{T}(V)$ of $V$ into $\mathcal{T}(V)$ such that for any linear mapping $R:V\to\mathcal{A}$ from $V$ into an algebra there exists a unique algebra homomorphism $\mathcal{T}(R):\mathcal{T}(V)\to \mathcal{A}$ such that the following diagram commutes,
\[
\xymatrix{
V\ar[r]^R\ar[d]_\imath & \mathcal{A} \\
\mathcal{T}(V)\ar[ur]_{\mathcal{T}(R)} &
}
\]
i.e.\ $\mathcal{T}(R)\circ \imath = R$.

Conversely, any homomorphism $Q:\mathcal{T}(V)\to \mathcal{A}$ is uniquely determined by its restriction to $V$.

In a similar way, an involution on $V$ gives rise to a unique extension as an involution on $\mathcal{T}(V)$. Thus for a $*$-vector space $V$ we can form the tensor $*$-algebra $\mathcal{T}(V)$. The tensor $*$-algebra $\mathcal{T}(V)$ becomes a $*$-bialgebra, if we extend the linear $*$-maps
\[
\begin{array}{rclcrcl}
\varepsilon&:&V\to \mathbb{C},&& \varepsilon(v) &=& 0, \\
\Delta&:& V\to \mathcal{T}(V)\otimes\mathcal{T}(V), && \Delta(v)&=&v\otimes\mathbf{1}+\mathbf{1}\otimes v, 
\end{array}
\]
in the same way. We will denote the coproduct $\mathcal{T}(\Delta)$ and the counit $\mathcal{T}(\varepsilon)$ again by $\Delta$ and $\varepsilon$. The tensor $*$-algebra is even a Hopf $*$-algebra with the antipode defined by $S(v)=-v$ on the generators and extended as an anti-homomorphism.

We will now study L\'evy processes on $\mathcal{T}(V)$. Let $D$ be a pre-Hilbert space and suppose we are given
\begin{enumerate}
\item
a linear $*$-map $R:V\to\mathcal{L}(D)$,
\item
a linear map $N:V\to D$, and
\item
a linear $*$-map $\psi:V\to \mathbb{C}$ (i.e.\ a hermitian linear functional),
\end{enumerate}
then
\begin{equation}\label{I-add-levy}
J_t(v) = \Lambda_t\big(R(v)\big) + A^*_t(N(v)\big) + A_t\big(N(v^*)\big) + t\psi(v)
\end{equation}
for $v\in V$ extends to a L\'evy process $(j_t)_{t\ge 0}$, $j_t=\mathcal{T}(J_t)$, on $\mathcal{T}(V)$ (w.r.t.\ the vacuum state).

In fact, all L\'evy processes on $\mathcal{T}(V)$ are of this form, cf.\ \cite{schuermann91c}.

The maps $(R,N,\psi)$ can be extended to a  Sch\"urmann triple on $\mathcal{T}(V)$ as follows
\begin{enumerate}
\item
Set $\rho=\mathcal{T}(R)$.
\item
Define $\eta:\mathcal{T}(V)\to D$ by $\eta(\mathbf{1})=0$, $\eta(v)=N(v)$ for $v\in V$, and 
\[
\eta(v_1\otimes \cdots \otimes v_n)=R(v_1)\cdots R(v_{n-1})N(v_n)
\]
for homogeneous elements $v_1\otimes \cdots \otimes v_n\in V^{\otimes n}$, $n\ge 2$.
\item
Finally, define $L:\mathcal{T}(V)\to \mathbb{C}$ by $L(\mathbf{1})=0$, $L(v)=\psi(v)$ for $v\in V$, and 
\[
\eta(v_1\otimes \cdots \otimes v_n)=
\left\{
\begin{array}{lcl}
\big\langle N(v_1^*),N(v_2)\big\rangle &\mbox{ if } n=2, \\
\big\langle N(v_1),R(v_2)\cdots R(v_{n-1})N(v_n)\big\rangle &\mbox{ if } n\ge 3, \\
\end{array}\right.
\]
for homogeneous elements $v_1\otimes \cdots \otimes v_n\in V^{\otimes n}$, $n\ge 2$.
\end{enumerate}

One can prove furthermore that all Sch\"urmann triples of $\mathcal{T}(V)$ are of this form.

\subsection{L\'evy processes on real Lie algebras}\label{I-subsec-lie}

The theory of factorizable representations was developed in the early seventies by Araki, Streater, Parthasarathy, Schmidt, Guichardet, $\cdots$, see, e.g. \cite{guichardet72,parthasarathy+schmidt72} and the references therein, or Section 5 of the historical survey by Streater \cite{streater00}. In this Section we shall see that in a sense this theory is a special case of the theory of L\'evy processes on involutive bialgebras.

\begin{definition}
A {\em Lie algebra} $\eufrak{g}$ over a field $\mathbb{K}$ is a $\mathbb{K}$-vector space with a linear map $[\cdot,\cdot]:\eufrak{g}\times\eufrak{g}\to\eufrak{g}$ called {\em Lie bracket} that satisfies the following two properties.
\begin{enumerate}
\item
Anti-symmetry: for all $X,Y\in\eufrak{g}$, we have
\[
[X,Y]=[Y,X].
\]
\item
Jacobi identity: for all $X,Y,Z\in\eufrak{g}$, we have
\[
\big[X,[Y,Z]\big]+\big[Y,[Z,X]\big]+\big[Z,[X,Y]\big]=0.
\]
\end{enumerate}
\end{definition}
For $\mathbb{K}=\mathbb{R}$, we call $\eufrak{g}$ a {\em real Lie algebra}, for $\mathbb{K}=\mathbb{C}$ a {\em complex Lie algebra}.

If $\mathcal{A}$ is an algebra, then $[a,b]=ab-ba$ defines a Lie bracket on $\mathcal{A}$.

We will see below that we can associate a Hopf $*$-algebra to a real Lie algebra, namely its universal enveloping algebra. But it is possible to define L\'evy processes on real Lie algebras without explicit reference to any coalgebra structure.

\begin{definition}\label{def levy}
Let $\eufrak{g}$ be a Lie algebra over $\mathbb{R}$, $D$ be a pre-Hilbert space, and $\Omega\in D$ a unit vector. We call a family $\big(j_{st}:\eufrak{g}\to\mathcal{L}(D)\big)_{0\le s\le t}$ of representations of $\eufrak{g}$ by anti-hermitian operators (i.e.\ satisfying $j_{st}(X)^*=- j_{st}(X)$ for all $X\in \eufrak{g}$, $0\le s\le t$) a {\em L\'evy process on $\eufrak{g}$} over $D$ (with respect to $\Omega$), if the following conditions are satisfied.
\begin{enumerate}
\item
(Increment property) We have
\[
j_{st}(X)+j_{tu}(X)=j_{su}(X)
\]
for all $0\le s\le t \le u$ and all $X\in\eufrak{g}$.
\item
(Independence)
We have $[j_{st}(X),j_{s't'}(Y)]=0$ for all $X,Y\in\eufrak{g}$, $0\le s\le t\le s'\le t'$ and
\[
\langle\Omega, j_{s_1t_1}(X_1)^{k_1}\cdots j_{s_nt_n}(X_n)^{k_n}\Omega\rangle = \langle\Omega, j_{s_1t_1}(X_1)^{k_1}\Omega\rangle\cdots\langle\Omega, j_{s_nt_n}(X_n)^{k_n}\Omega\rangle
\]
for all $n,k_1,\ldots,k_n\in\mathbb{N}$, $0\le s_1\le t_1\le s_2\le \cdots \le t_n$, $X_1,\ldots,X_n\in\eufrak{g}$.
\item
(Stationarity) For all $n\in\mathbb{N}$ and all $X\in\eufrak{g}$, the moments
\[
m_n(X;s,t)=\langle\Omega,j_{st}(X)^n\Omega\rangle
\]
depend only on the difference $t-s$.
\item
(Weak continuity)
We have $\lim_{t\searrow s} \langle\Omega,j_{st}(X)^n\Omega\rangle =0$ for all $n\in\mathbb{N}$ and all $X\in\eufrak{g}$. 
\end{enumerate}
\end{definition}

For a construction of L\'evy processes on several Lie algebras of interest to physics and for several examples see also \cite{accardi+franz+skeide02,franz03}.

Let $\eufrak{g}$ be a real Lie algebra. Then the complex vector space $\eufrak{g}_\mathbb{C}=\mathbb{C}\otimes_\mathbb{R}\eufrak{g}=\eufrak{g}\oplus i\eufrak{g}$ is a complex Lie algebra with the Lie bracket
\[
[X+iY,X'+iY']=[X,X']-[Y,Y']+i\big([X,Y']+[Y,X']\big)
\]
for $X,X',Y,Y'\in\eufrak{g}$.

We summarize the relation between real Lie algebras and complex involutive Lie algebras in the following proposition.
\begin{proposition}
\begin{enumerate}
\item
The map $*:\eufrak{g}_\mathbb{C}\to\eufrak{g}_\mathbb{C}$, $Z=X+iY\mapsto Z^*=-X+iY$, $X,Y\in\eufrak{g}$, is an involution on $\eufrak{g}_\mathbb{C}$, i.e.\ it satisfies
\[
(Z^*)^*=Z \qquad \mbox{ and }\qquad [Z_1,Z_2]^* =[Z_2^*,Z_1^*]
\]
for all $Z,Z_1,Z_2\in\eufrak{g}_\mathbb{C}$
\item
The functor $\eufrak{g}\mapsto(\eufrak{g}_\mathbb{C},*)$ is an isomorphism between the category of real Lie algebras and the category of involutive complex Lie algebras.
\end{enumerate}
\end{proposition}

The {\em universal enveloping algebra} $\mathcal{U}(\eufrak{g})$ of a Lie algebra $\eufrak{g}$ can be constructed as the quotient $\mathcal{T}(\eufrak{g})/\mathcal{J}$ of the tensor algebra $\mathcal{T}(\eufrak{g})$ over $\eufrak{g}$ by the ideal $\mathcal{J}$ generated by
\[
\big\{X\otimes Y-Y\otimes X-[X,Y]|X,Y\in\eufrak{g}\big\}.
\]
The universal enveloping algebra is characterized by a universal property. Composing the embedding $\imath:\eufrak{g}\to\mathcal{T}(\eufrak{g})$ with the canonical projection $p:\mathcal{T}(\eufrak{g})\to\mathcal{T}(\eufrak{g})$ we get an embedding $\imath'=p\circ\imath:\eufrak{g}\to\mathcal{U}(\eufrak{g})$ of $\eufrak{g}$ into its enveloping algebra. For every algebra $\mathcal{A}$ and every Lie algebra homomorphism $R:\eufrak{g}\to\mathcal{A}$ there exists a unique algebra homomorphism $\mathcal{U}(R):\mathcal{U}(\eufrak{g})\to\mathcal{A}$ such that the following diagram commutes,
\[
\xymatrix{
\eufrak{g}\ar[r]^R\ar[d]_{\imath'} & \mathcal{A} \\
\mathcal{U}(\eufrak{g})\ar[ur]_{\mathcal{\mathcal{U}}(R)} &
}
\]
i.e.\ $\mathcal{U}(R)\circ \imath' = R$. If $\eufrak{g}$ has an involution, then it can be extended to an involution of $\mathcal{U}(g)$.

The enveloping algebra $\mathcal{U}(\eufrak{g})$ becomes a bialgebra, if we extend the Lie algebra homomorphisms
\[
\begin{array}{rclcrcl}
\varepsilon&:&\eufrak{g}\to \mathbb{C},&& \varepsilon(X) &=& 0, \\
\Delta&:& \eufrak{g}\to \mathcal{U}(\eufrak{g})\otimes\mathcal{U}(\eufrak{g}), && \Delta(X)&=&X\otimes\mathbf{1}+\mathbf{1}\otimes X, 
\end{array}
\]
to $\mathcal{U}(\eufrak{g})$. We will denote the coproduct $\mathcal{U}(\Delta)$ and the counit $\mathcal{U}(\varepsilon)$ again by $\Delta$ and $\varepsilon$. It is even a Hopf algebra with the antipode $S:\mathcal{U}(\eufrak{g})\to\mathcal{U}(\eufrak{g})$ given by $S(X)=-X$ on $\eufrak{g}$ and extended as an anti-homomorphism.

\begin{proposition}
Let $\eufrak{g}$ be a real Lie algebra and $\mathcal{U}=\mathcal{U}(\eufrak{g}_\mathbb{C})$ the enveloping algebra of its complexification.
\begin{enumerate}
\item
Let $(j_{st})_{0\le s\le t}$ be a L\'evy process on $\mathcal{U}$. Then its restriction to $\eufrak{g}$ is a  L\'evy process on $\eufrak{g}$.
\item
Let $(k_{st})_{0\le s\le t}$ now be a L\'evy process on $\eufrak{g}$. Then its extension to $\mathcal{U}$ given by the universal property is a  L\'evy process on $\mathcal{U}$.
\end{enumerate}
This establishes a one-to-one correspondence between L\'evy processes on a real Lie algebra and L\'evy processes on the universal enveloping algebra of its complexification.
\end{proposition}

We will now show that L\'evy processes on real Lie algebras are the same as factorizable representation of current algebras.

Let $\eufrak{g}$ be a real Lie algebra and $(\mathbb{T},\mathcal{T},\mu)$ a measure space (e.g.\ the real line $\mathbb{R}$ with the Lebesgue measure $\lambda$). Then the set of $\eufrak{g}$-valued step functions
\[
\eufrak{g}^I=
\left\{X=\sum_{i=1}^nX_i\mathbf 1_{M_i};X_i\in\eufrak g, M_i\in\mathcal{T},\mu(M_i)<\infty, M_i\subseteq I, n\in\mathbb N\right\}.
\]
on an $I\subseteq\mathbb{T}$ is again a real Lie algebra with the pointwise Lie bracket. For $I_1\subseteq I_2$ we have an inclusion $i_{I_1,I_2}:\eufrak{g}^{I_1}\to\eufrak{g}^{I_2}$, simply extending the functions as zero outside $I_1$. Furthermore, for disjoint subsets $I_1,I_2\in\mathcal{T}$, $\eufrak{g}^{I_1\cup I_2}$ is equal to the direct sum $\eufrak{g}^{I_1}\oplus\eufrak{g}^{I_2}$. If $\pi$ is a representation of $\eufrak{g}^\mathbb{T}$ and $I\in\mathcal{T}$, then have also a representation $\pi^I=\pi\circ i_{I,\mathbb{T}}$ of $\eufrak{g}^I$

Recall that for two representations $\rho_1,\rho_2$ of two Lie algebras $\eufrak{g}_1$ and $\eufrak{g_2}$, acting on (pre-) Hilbert spaces $H_1$ and $H_2$, we can define a representation $\rho=(\rho_1\otimes\rho_2)$ of $\eufrak{g}_1\oplus\eufrak{g}_1$ acting on $H_1\otimes H_2$  by
\[
(\rho_1\otimes\rho_2)(X_1+X_2) = \rho_1(X_1)\otimes\mathbf{1}+\mathbf{1}\otimes\rho_2(X_2),
\]
for $X_1\in\eufrak{g}_1$, $X_2\in\eufrak{g}_2$.

\begin{definition}
A triple $(\pi,D,\Omega)$ consisting of a representation $\pi$ of $\eufrak{g}^\mathbb{T}$ by anti-hermitian operators and a unit vector $\Omega\in D$ is called a {\em factorizable representation} of the simple current algebra $\eufrak{g}^\mathbb{T}$, if the following conditions are satisfied.
\begin{enumerate}
\item
(Factorization property)
For all $I_1,I_2\in\mathcal{T}$, $I_1\cap I_2=\emptyset$, we have
\[
(\pi^{I_1\cup I_2},D,\Omega)\cong (\pi^{I_1}\otimes\pi^{I_2},D\otimes D,\Omega\otimes\Omega).
\]
\item
(Invariance)
The linear functional $\varphi_I:\mathcal{U}(\eufrak{g})\to\mathbb{}$ determined by
\[
\varphi_I(X^n) = \langle\Omega, \pi(X\mathbf{1}_{I})^n\Omega\rangle
\]
for $X\in\eufrak{g}$, $I\in\mathcal{T}$ depends only on $\mu(I)$.
\item
(Weak continuity)
For any sequence $(I_k)_{k\in \mathbb{N}}$ with $\lim_{k\to\infty}\mu(I_k)=0$ we have $\lim_{k\to\infty}\varphi_{I_k}(u)=\varepsilon(u)$ for all $u\in\mathcal{U}(\eufrak{g})$.
\end{enumerate}
\end{definition}

\begin{proposition}
Let $\eufrak{g}$ be a real Lie algebra and take $(\mathbb{T},\mathcal{T},\mu)=(\mathbb{R}_+,\mathcal{B}(\mathbb{R}_+),\lambda)$. Then we have a one-to-one correspondence between factorizable representations of $\eufrak{g}^{\mathbb{R}_+}$ and L\'evy processes on $\eufrak{g}$.
\end{proposition}

The relation which is used to switch from one to the other is
\[
\pi(X\mathbf{1}_{[s,t[})=j_{st}(X)
\]
for $0\le s\le t$ and $X\in\eufrak{g}$.

\begin{proposition}
Let $\eufrak{g}$ be a real Lie algebra and $(\mathbb{T},\mathcal{T},\mu)$ a measure space without atoms. Then all factorizable representations of $\eufrak{g}^\mathbb{T}$ are characterized by generators or equivalently by Sch\"urmann triples on $\mathcal{U}(\eufrak{g}_\mathbb{C})$. They have a realization on the symmetric Fock space $\Gamma\big(L^2(\mathbb{T},\mathcal{T},\mu)\big)$ determined by
\[
\pi(X\mathbf{1}_I) = A^*\big(\mathbf{1}_I\times\rho(X)\big) + \Lambda\big(\mathbf{1}_I\otimes\rho(X)\big)+A\big(\mathbf{1}_I\otimes\eta(X^*)\big) + \mu(I)L(X)
\]
for $I\in\mathcal{T}$ with $\mu(I)<\infty$ and $X\in\eufrak{g}$.
\end{proposition}

\subsection{The quantum Az\'ema martingale}\label{I-sub-azema}

Let $q\in\mathbb{C}$ and $\mathcal{B}_q$ the involutive bialgebra with generators $x,x^*,y,y^*$ and relations
\begin{gather*}
yx=qxy,\qquad x^*y=qyx^*,\\
\Delta(x)=x\otimes y + \mathbf{1}\otimes x, \qquad \Delta(y)=y\otimes y, \\
\varepsilon(x)=0, \qquad \varepsilon(y)=1.
\end{gather*}

\begin{proposition}
There exists a unique Sch\"urmann triple on $\mathcal{B}_q$ acting on $D=\mathbb{C}$ with
\[
\begin{array}{lclclcl}
\rho(y)&=&q, && \rho(x)=0, \\
\eta(y)&=&0, && \eta(x)=1, \\
L(y)&=&0, && L(x)=0.
\end{array}
\]
\end{proposition}

Let $(j_{st})_{0\le s\le t}$ be the associated L\'evy process on $\mathcal{B}_q$ and set $Y_t=j_{0t}(y)$, $X_t=j_{0t}(x)$, and $X^*_t=j_{0t}(x^*)$. These operator processes are determined by the quantum stochastic differential equations
\begin{eqnarray}
{\rm d}Y_t &=& (q-1) Y_t{\rm d}\Lambda_t, \label{I-azema-y} \\
{\rm d}X_t &=& {\rm d}A^*_t + (q-1) X_t{\rm d}\Lambda_t, \\
{\rm d}X^*_t &=& {\rm d}A_t + (\overline{q}-1) X_t{\rm d}\Lambda_t,
\end{eqnarray}
with initial conditions $Y_0=\mathbf{1}$, $X_0=X^*_0=0$. This process is the quantum Az\'ema martingale introduced by Parthasarathy \cite{parthasarathy90}, see also \cite{schuermann91b}. The first Equation \eqref{I-azema-y} can be solved explicitely, the operator process $(Y_t)_{t\ge 0}$ is the second quantization of multiplication by $q$, i.e.,
\[
Y_t = \Gamma_t(q), \qquad \mbox{ for }t\ge 0
\]
Its action on exponential vectors is given by
\[
Y_t\mathcal{E}(f)=\mathcal{E}\big(qf\mathbf{1}_{[0,t[}+f\mathbf{1}_{[t,+\infty[}\big). 
\]
The hermitian operator process $(Z_t)_{t\ge 0}$ defined by $Z_t=X_t+X_t^*$ has the same joint moments as the classical Az\'ema martingale $(M_t)_{t\ge0}$ introduced by Az\'ema and Emery, cf.\ \cite{emery89}, i.e.\ is has the same joint moments,
\[
\langle\Omega,Z^{n_1}_{t_1}\cdots Z^{n_k}_{t_k}\Omega\rangle=
E\left(M^{n_1}_{t_1}\cdots M^{n_k}_{t_k}\right)
\]
for all $n_1,\ldots,n_k\in\mathbb{N}$, $t_1,\ldots,t_k\in\mathbb{R}_+$. This was the first example of a classical normal martingale having the so-called {\em chaotic representation property}, which is not a classical L\'evy process.

%% file: independence.tex
\chapter{The Five Universal Independences}\label{chapter-III}

In classical probability theory there exists only one canonical notion of independence. But in quantum probability many different notions of independence have been used, e.g., to obtain central limit theorems or to develop a quantum stochastic calculus. If one requires that the joint law of two independent random variables should be determined by their marginals, then an independence gives rise to a product. Imposing certain natural condition, e.g., that functions of independent random variables should again be independent or an associativity property, it becomes possible to classify all possible notions of independence. This program has been carried out in recent years by Sch\"urmann \cite{schuermann95}, Speicher \cite{speicher96}, Ben Ghorbal and Sch\"urmann \cite{benghorbal+schuermann99,benghorbal+schuermann02}, and Muraki \cite{muraki02a,muraki02}. In this chapter we will present the results of these classifications. Furthermore we will formulate a category theoretical approach to the notion of independence and show that boolean, monotone, and anti-monotone independence can be reduced to tensor independence in a similar way as the bosonization of fermi independence \cite{hudson+parthasarathy86} or the symmetrization of \cite[Section 3]{schuermann93}.

\section{Preliminaries on category theory}\label{III-prelim}

We give the basic definitions and properties from category theory that we shall use. For a thorough introduction, see \cite{maclane98}.

\begin{definition}
A {\em category}\index{category} $\mathcal C$ consists of 
\begin{itemize}
\item[(a)] a class ${\rm Ob}\,{\mathcal C}$ of {\em objects}\index{object} denoted by $A,B,C,\ldots$,
\item[(b)] a class ${\rm Mor}\,{\mathcal C}$ of {\em morphism}\index{morphism} (or {\em arrows}\index{arrow!|see{morphism}}) denoted by $g,f,h,\ldots$,
\item[(c)] mappings ${\rm tar}, {\rm src}:{\rm Mor}\,{\mathcal C}\to{\rm Ob}\,{\mathcal C}$ assigning to each morphism $f$ its {\em source}\index{source} (or {\em domain}\index{domain|see{source}}) ${\rm src}(f)$ and its {\em target}\index{target} (or {\em codomain}\index{codomain|see{target}}) ${\rm tar}(f)$. We will say that $f$ is a morphism in ${\mathcal C}$ from $A$ to $B$ or write ``$f:A\to B$ is a morphism in ${\mathcal C}$'' if $f$ is a morphism in ${\mathcal C}$ with source  ${\rm src}(f)=A$ and target ${\rm tar}(f)=B$,
\item[(d)]
a composition $(f,g)\mapsto g\circ f$ for pairs of morphisms $f,g$ that satisfy ${\rm src}(g)={\rm tar}(f)$, 
\item[(e)] and a map ${\rm id}:{\rm Ob}\,{\mathcal C}\to{\rm Mor}\,{\mathcal C}$ assigning to an object $A$ of $\mathcal C$ the {\em identity morphism}\index{identity morphism} ${\rm id}_A:A\to A$,
\end{itemize}
such that the
\begin{itemize}
\item[(1)] {\em associativity property}\index{associativity property}: for all morphisms $f:A\to B$, $g:B\to C$, and $h:C\to D$ of  $\mathcal C$, we have
\[
(h\circ g)\circ f = h\circ (g\circ f),
\]
and the
\item[(2)] {\em identity property}\index{identity property}:  ${\rm id}_{{\rm tar}(f)}\circ f = f$ and $f\circ {\rm id}_{{\rm src}(f)}=f$ holds for all morphisms $f$ of $\mathcal C$,
\end{itemize}
are satisfied.
\end{definition}

Let us emphasize that it is not so much the objects, but the morphisms that contain the essence of a category (even though categories are usually named after their objects). Indeed, it is possible to define categories without referring to the objects at all, see the definition of ``arrows-only metacategories'' in \cite[Page 9]{maclane71,maclane98}. The objects are in one-to-one correspondence with the identity morphisms, in this way ${\rm Ob}\,\mathcal{C}$ can always be recovered from ${\rm Mor}\,\mathcal{C}$.

We give an example.

\begin{example}\label{exa-set}
Let ${\rm Ob}\,\eufrak{Set}$ be the class of all sets (of a fixed universe) and ${\rm Mor}\,\eufrak{Set}$ the class of total functions between them. Recall that a {\em total function}\index{total function}\index{function!total} (or simply {\em function}\index{function}) is a triple $(A,f,B)$, where $A$ and $B$ are sets, and $f\subseteq A\times B$ is a subset of the cartesian product of $A$ and $B$ such that for a given $x\in A$ there exists a unique $y\in B$ with $(x,y)\in f$. Usually one denotes this unique element by $f(x)$, and writes $x\mapsto f(x)$ to indicate $\big(x,f(x)\big)\in f$. The triple $(A,f,B)$ can also be given in the form $f:A\to B$. We define
\[
{\rm src}\big((A,f,B)\big)= A , \quad \text{ and }\quad {\rm tar}\big((A,f,B)\big)= B.
\]
The composition of two morphisms $(A,f,B)$ and $(B,g,C)$ is defined as
\[
(B,g,C)\circ(A,f,B) = (A,g\circ f,C),
\]
where $g\circ f$ is the usual composition of the functions $f$ and $g$, i.e.
\[
g\circ f= \{(x,z)\in A\times C; \text{ there exists a }y\in B\text{ s.t. }(x,y)\in f \text{ and }(y,z)\in g\}.
\]
The identity morphism assigned to an object $A$ is given by $(A,{\rm id}_A,A)$, where ${\rm id}_A\subseteq A\times A$ is the identity function, ${\rm id}_A=\{(x,x); x\in A\}$. It is now easy to check that these definitions satisfy the associativity property and the identity property, and therefore define a category. We shall denote this category by $\eufrak{Set}$\index{$\eufrak{Set}$}.
\end{example}

\begin{definition}
Let  $\mathcal{C}$ be a category. A morphism $f:A\to B$ in $\mathcal C$ is called an {\em isomorphism}\index{isomorphism} (or {\em invertible}\index{invertible morphism|see{isomorphism}}\index{morphism!invertible|see{isomorphism}}), if there exists a morphism $g:B\to A$ in $\mathcal C$ such that $g\circ f = {\rm id}_A$ and $f\circ g={\rm id}_B$. Such a morphism $g$ is uniquely determined, if it exists, it is called the {\em inverse}\index{inverse morphism}\index{morphism!inverse} of $f$ and denoted by $f^{-1}$. Objects $A$ and $B$ are called {\em isomorphic}\index{isomorphic}, if there exists an isomorphism $f:A\to B$.

Morphisms $f$ with ${\rm tar}(f)={\rm src}(f)=A$ are called {\em endomorphisms}\index{endomorphism} of $A$. Isomorphic endomorphism are called {\em automorphisms}\index{automorphism}.
\end{definition}

For an arbitrary pair of objects $A,B\in {\rm Ob}\,\mathcal{C}$ we define ${\rm Mor}_{\mathcal{C}}(A,B)$ to be the collection of morphisms from $A$ to $B$, i.e.
\[
{\rm Mor}_{\mathcal{C}}(A,B)= \{ f\in {\rm Mor}\, \mathcal{C}; {\rm src}(f)=A \text{ and } {\rm tar}(f)=B\}.
\]
Often the collections ${\rm Mor}_{\mathcal{C}}(A,B)$ are also denoted by ${\rm hom}_\mathcal{C}(A,B)$ and called the {\em hom-sets}\index{hom-set} of $\mathcal{C}$. In particular, ${\rm Mor}_{\mathcal{C}}(A,A)$ contains exactly the endomorphisms of $A$, they form a semigroup with identity element with respect to the composition of $\mathcal{C}$ (if ${\rm Mor}_{\mathcal{C}}(A,A)$ is a set).

Compositions and inverses of isomorphisms are again isomorphisms. The automorphisms of an object form a group (if they form a set).

\begin{example}\label{exa-part}
We define now the category of sets and partial functions, $\eufrak{Part}$\index{$\eufrak{Part}$}. As objects we have again the class of all sets  (of a fixed universe) and as morphisms we take the class of partial functions between them. A {\em partial function}\index{partial function}\index{function!partial} is a triple $(A,f,B)$ where $A$ and $B$ are sets, and $f\subseteq A\times B$ is a subset of the cartesian product of $A$ and $B$ such that for a given $x\in A$ there exists at most one $y\in B$ with $(x,y)\in f$, i.e.\ if $(x,y)\in f$ and $(x,y')\in f$ then $y=y'$. If for a given $x\in A$ there exists no $y\in B$ with $(x,y)\in f$, then we say that $f$ is undefined at $x$.
\end{example}

\begin{example}
In $\eufrak{Set}$ and $\eufrak{Part}$ the morphisms are maps. We will now define a category in which the arrows are not given by maps. The objects of the category $\eufrak{Rel}$\index{$\eufrak{Rel}$} are all sets (or the elements of some set of sets $M$, if we want to define the small category $\eufrak{Rel}_M$), and the morphisms between two objects $A$ and $B$ are all binary relations $R\subseteq A\times B$. The identity morphism is given by the identity relation ${\rm id}_A=\{(a,a); a\in A\}$, and the composition $R\circ S:A\to C$ of two morphisms $S:A\to B$ and $R:B\to C$ is defined by the relative product
\[
R\circ S = \{ (a,c)\subseteq A\times C; \text{ there exists a } b\in B \text{ s.t. } (a,b)\in S\text{ and } (b,c)\in R\}.
\]
This category has an additional structure: for each morphism $R:A\to B$ there is a converse relation $R^{\#}:B\to A$ consisting of all pairs $(b,a)$ with $(a,b)\in R$.
\end{example}

\begin{example}\label{exa-semigroup}
Let $(G,\circ,e)$ be a semigroup with identity element $e$. Then $(G,\circ,e)$ can be viewed as a category. The only object of this category is $G$ itself, and the morphisms are the elements of $G$. The identity morphism is $e$ and the composition is given by the composition of $G$.
\end{example}

\begin{definition}\label{def-basic-opp}
For every category $\mathcal{C}$ we can define its {\em dual or opposite category}\index{dual category}\index{opposite category}\index{category!dual}\index{category!opposite} $\mathcal{C}^{\rm op}$. It has the same objects and morphisms, but target and source are interchanged, i.e.\
\[
{\rm tar}_{\mathcal{C}^{\rm op}}(f) = {\rm src}_{\mathcal{C}}(f) \text{ and } {\rm src}_{\mathcal{C}^{\rm op}}(f) ={\rm tar}_{\mathcal{C}}(f)
\]
and the composition is defined by $f\circ_{\rm op} g=g \circ f$. We obviously have $\mathcal{C}^{\rm op\,op}=\mathcal{C}$.
\end{definition}

Dualizing, i.e.\ passing to the opposite category, is a very useful concept in category theory. Whenever we define something in a category, like an epimorphism, a terminal object, a product, etc., we get a definition of a ``cosomething'', if we take the corresponding definition in the opposite category. For example, an {\em epimorphism}\index{epimorphism} or {\em epi}\index{epi|see{epimorphism}} in $\mathcal{C}$ is a morphism in $\mathcal{C}$ which is right cancellable, i.e.\ $h\in{\rm Mor}\,\mathcal{C}$ is called an epimorphism, if for any morphisms $g_1,g_2\in{\rm Mor}\,\mathcal{C}$ the equality $g_1\circ h=g_2\circ h$ implies $g_1=g_2$. The dual notion of a epimorphism is a morphism, which is an epimorphism in the category $\mathcal{C}^{\rm op}$, i.e.\ a morphism that is left cancellable. It could therefore be called a ``coepimorphism'', but the generally accepted name is {\em monomorphism}\index{monomorphism} or {\em monic}\index{monic|see{monomorphism}}. The same technique of dualizing applies not only to definitions, but also to theorems. A morphism $r:B\to A$ in $\mathcal{C}$ is called a {\em right inverse}\index{right inverse}\index{inverse!right}\index{morphism!inverse!right} of $h:A\to B$ in $\mathcal{C}$, if $h\circ r = {\rm id}_B$. If a morphism has a right inverse, then it is necessarily an epimorphism, since $g_1\circ g=g_2\circ h$ implies $g_1=g_1\circ g\circ r=g_2\circ h \circ r=g_2$, if we compose both sides of the equality with a right inverse $r$ of $h$. Dualizing this result we see immediately that a morphism $f:A\to B$ that has a {\em left inverse}\index{left inverse}\index{inverse!left}\index{morphism!inverse!left} (i.e.\ a morphism $l:B\to A$ such that $l\circ f={\rm id}_A$) is necessarily a monomorphism. Left inverses are also called {\em retractions}\index{retraction|see{left inverse}} and right inverses are also called {\em sections}\index{section|see{right inverse}}. Note that one-sided inverses are usually not unique.

\begin{definition}
A category $\mathcal{D}$ is called a {\em subcategory}\index{subcategory}\index{category!sub-} of the category $\mathcal{C}$, if
\begin{itemize}
\item[(1)]
the objects of $\mathcal{D}$ form a subclass of ${\rm Ob}\,\mathcal{C}$, and the morphisms of  $\mathcal{D}$ form a subclass of ${\rm Mor}\,\mathcal{C}$,
\item[(2)]
for any morphism $f$ of $\mathcal{D}$, the source and target of $f$ in $\mathcal{C}$ are objects of  $\mathcal{D}$ and agree with the source and target taken in  $\mathcal{D}$,
\item[(3)]
for every object $D$ of $\mathcal{D}$, the identity morphism ${\rm id}_D$ of $\mathcal{C}$ is a morphism of $\mathcal{D}$, and
\item[(4)]
for any pair $f:A\to B$ and $g:B\to C$ in $\mathcal{D}$, the composition $g\circ f$ in $\mathcal{C}$ is a morphism of $\mathcal{D}$ and agrees with the composition of $f$ and $g$ in $\mathcal{D}$.
\end{itemize}
A subcategory $\mathcal{D}$ of $\mathcal{C}$ is called {\em full}\index{full subcategory}\index{category!full sub-}\index{subcategory!full}, if for any two objects $A,B\in{\rm Ob}\,\mathcal{D}$ all $\mathcal{C}$-morphisms from $A$ to $B$ belong also to $\mathcal{D}$, i.e.\ if
\[
{\rm Mor}_{\mathcal{D}}(A,B)={\rm Mor}_{\mathcal{C}}(A,B).
\]
\end{definition}

\begin{remark}
If $D$ is an object of $\mathcal{D}$, then the identity morphism of $D$ in $\mathcal{D}$ is the same as that in $\mathcal{C}$, since the identity element of a semigroup is unique, if it exists.
\end{remark}

\begin{definition}
Let $\mathcal{C}$ and $\mathcal{D}$ be two categories. A {\em covariant functor}\index{covariant functor}\index{functor!covariant} (or simply {\em functor}\index{functor}) $T:\mathcal{C}\to\mathcal{D}$ is a map for objects and morphisms, every object $A\in{\rm Ob}\,\mathcal{C}$ is mapped to an object $T(A)\in{\rm Ob}\,\mathcal{D}$, and every morphism $f:A\to B$ in $\mathcal{C}$ is mapped to a morphism $T(f):T(A)\to T(B)$ in $\mathcal{D}$, such that the identities and the composition are respected, i.e.\ such that
\begin{gather*}
T({\rm id}_A) = {\rm id}_{T(A)}, \quad\text{ for all } A\in {\rm Ob}\,\mathcal{C} \\
T(g\circ f) = T(g)\circ T(f), \quad\text{ whenever } g\circ f \text{ is defined in } \mathcal{C}.
\end{gather*}
We will denote the collection of all functors between two categories $\mathcal{C}$ and $\mathcal{D}$ by ${\rm Funct}(\mathcal{C},\mathcal{D})$.

A {\em contravariant functor}\index{contravariant functor}\index{functor!contravariant} $T:\mathcal{C}\to\mathcal{D}$ maps an object $A\in{\rm Ob}\,\mathcal{C}$ to  an object $T(A)\in{\rm Ob}\,\mathcal{D}$, and a morphism $f:A\to B$ in $\mathcal{C}$ to a morphism $T(f):T(B)\to T(A)$ in $\mathcal{D}$, such such that
\begin{gather*}
T({\rm id}_A) = {\rm id}_{T(A)}, \quad\text{ for all } A\in {\rm Ob}\,\mathcal{C} \\
T(g\circ f) = T(f)\circ T(g), \quad\text{ whenever } g\circ f \text{ is defined in } \mathcal{C}.
\end{gather*}
\end{definition}

\begin{example}
Let $\mathcal{C}$ be a category. The {\em identity functor}\index{identity functor}\index{functor!identity} ${\rm id}_{\mathcal{C}}:\mathcal{C}\to\mathcal{C}$ is defined by ${\rm id}_{\mathcal{C}}(A)=A$ and ${\rm id}_{\mathcal{C}}(f)=f$.
\end{example}

\begin{example}
The {\em inclusion}\index{inclusion} of a subcategory $\mathcal{D}$ of $\mathcal{C}$ into $\mathcal{C}$ also defines a functor, we can denote it by $\subseteq:\mathcal{D}\to\mathcal{C}$ or by  $\mathcal{D}\subseteq\mathcal{C}$.
\end{example}

\begin{example}
The functor ${\rm op}:\mathcal{C}\to\mathcal{C}^{\rm op}$ that is defined as the identity map on the objects and morphisms is a contravariant functor. This functor allows to obtain covariant functors from contravariant ones. Let $T:\mathcal{C}\to\mathcal{D}$ be a contravariant functor, then $T\circ {\rm op}:\mathcal{C}^{\rm op}\to\mathcal{D}$ and ${\rm op}\circ T:\mathcal{C}\to\mathcal{D}^{\rm op}$ are covariant.
\end{example}

\begin{example}\label{exa-functor-semigroup}
Let $G$ and $H$ be unital semigroups, then the functors $T:G\to H$ are precisely the identity preserving semigroup homomorphisms from $G$ to $H$.
\end{example}

Functors can be composed, if we are given two functors $S:\mathcal{A}\to\mathcal{B}$ and $T:\mathcal{B}\to\mathcal{C}$, then the composition $T\circ S:\mathcal{A}\to\mathcal{C}$,
\begin{gather*}
(T\circ S)(A)=T(S(A)), \quad \text{ for } A\in {\rm Ob}\, \mathcal{A}, \\
(T\circ S)(f)=T(S(f)), \quad \text{ for } f\in {\rm Mor}\, \mathcal{A},
\end{gather*}
is again a functor. The composite of two covariant or two contravariant functors is covariant, whereas the composite of a covariant and a contravariant functor is contravariant. The identity functor obviously is an identity w.r.t.\ to this composition. Therefore we can define categories of categories, i.e.\ categories whose objects are categories and whose morphisms are the functors between them.

\begin{definition}
Let $\mathcal{C}$ and $\mathcal{D}$ be two categories and let $S,T:\mathcal{C}\to\mathcal{D}$ be two functors between them. A {\em natural transformation}\index{natural transformation}\index{transformation!natural} (or {\em morphism of functors}\index{morphism of functors|see{natural transformation}}) $\eta:S\to T$ assigns to every object $A\in{\rm Ob}\,\mathcal{C}$ of $\mathcal{C}$ a morphism $\eta_A:S(A)\to T(A)$ such that the diagram
\[
\xymatrix{
S(A) \ar[r]^{\eta_A} \ar[d]_{S(f)} &  T(A)\ar[d]^{T(f)}  \\
S(B) \ar[r]^{\eta_B} &  T(B)
}
\]
is commutative for every morphisms $f:A\to B$ in $\mathcal{C}$. The morphisms $\eta_A$, $A\in{\rm Ob}\,\mathcal{C}$ are called the {\em components}\index{components of a natural transformation} of $\eta$. If every component $\eta_A$ of $\eta:S\to T$ is an isomorphism, then $\eta:S\to T$ is called a {\em natural isomorphism}\index{natural isomorphism}\index{isomorphism!natural} (or a {\em natural equivalence}\index{natural equivalence|see{natural isomomorphism}}\index{equivalence!natural|see{natural isomorphism}}), in symbols this is expressed as $\eta:S \cong T$.

We will denote the collection of all natural transformations between two functors $S,T:\mathcal{C}\to\mathcal{D}$ by ${\rm Nat}(S,T)$.
\end{definition}

\begin{definition}
Natural transformations can also be composed. Let $S,T,U:\mathcal{B}\to\mathcal{C}$ and let $\eta:S\to T$ and $\vartheta:T\to U$ be two natural transformations. Then we can define a natural transformation $\vartheta\cdot\eta:S\to U$, its components are simply $(\vartheta\cdot\eta)_A = \vartheta_A\circ\eta_A$. To show that this defines indeed a natural transformation, take a morphism $f:A\to B$ of $\mathcal{B}$. Then the following diagram is commutative, because the two trapezia are.
\[
\xymatrix{
S(A) \ar[rd]^{\eta_A}\ar[ddd]_{S(f)} \ar[rr]^{(\vartheta\cdot\eta)_A=\vartheta_A\circ\eta_A} &  &  U(A) \ar[ddd]^{U(f)} \\
& T(A) \ar[ur]^{\vartheta_A} \ar[d]^{T(f)} & \\
&  T(B) \ar[dr]^{\vartheta_B} &\\
S(B) \ar[ur]^{\eta_B} \ar[rr]_{(\vartheta\cdot\eta)_B=\vartheta_B\circ\eta_B}  &  &  U(B)
}
\]
For a given functor $S:\mathcal{B}\to\mathcal{C}$ there exists also the {\em identical natural transformation}\index{identical natural transformation}\index{natural transformation!identical}\index{transformation!identical natural} ${\rm id}_S:S\to S$ that maps $A\in{\rm Ob}\,\mathcal{B}$ to ${\rm id}_{S(A)}\in {\rm Mor}\,\mathcal{C}$, it is easy to check that it behaves as a unit for the composition defined above.

Therefore we can define the {\em functor category}\index{functor category}\index{category!functor} $\mathcal{C}^\mathcal{B}$ that has the functors from $\mathcal{B}$ to $\mathcal{C}$ as objects and the natural transformations between them as morphisms.
\end{definition}

\begin{remark}
Note that a natural transformation $\eta:S\to T$ has to be defined as the triple $\left(S, (\eta_A)_A, T\right)$ consisting of its the source $S$, its components $(\eta_A)_A$ and its target $T$. The components $(\eta_A)_A$ do not uniquely determine the functors $S$ and $T$, they can also belong to a natural transformation between another pair of functors $(S',T')$.
\end{remark}

\begin{definition}
Two categories $\mathcal{B}$ and $\mathcal{C}$ can be called {\em isomorphic}\index{isomorphic categories}\index{categories!isomorphic}, if there exists an invertible functor $T:\mathcal{B}\to\mathcal{C}$. A useful weaker notion is that of {\em equivalence}\index{equivalence} or {\em categorical equivalence}\index{categorical equivalence}\index{equivalence!categorical}. Two categories $\mathcal{B}$ and $\mathcal{C}$ are equivalent, if there exist functors $F:\mathcal{B}\to\mathcal{C}$ and $G:\mathcal{C}\to\mathcal{B}$ and natural isomorphisms $G\circ F\cong {\rm id}_\mathcal{B}$ and $F\circ G\cong {\rm id}_\mathcal{C}$.
\end{definition}

We will look at products and coproducts of objects in a category. The idea of the product of two objects is an abstraction of the Cartesian product of two sets. For any two sets $M_1$ and $M_2$ their Cartesian product $M_1\times M_2$ has the property that for any pair of maps $(f_1,f_2)$, $f_1:N\to M_1$, $f_2:N\to M_2$, there exists a unique map $h:N\to M_1\times M_2$ such that $f_i=\pi_i\circ h$ for $i=1,2$, where $p_i:M_1\times M_2\to M_i$ are the canonical projections $p_i(m_1,m_2)=m_i$. Actually, the Cartesian product $M_1\times M_2$ is characterized by this property up to isomorphism (of the category $\eufrak{Set}$, i.e.\ set-theoretical bijection).

\begin{definition}\label{III-def-product}
A tuple $(A\,\Pi\, B,\pi_A,\pi_B)$ is called a {\em product}\index{product} (or {\em binary product}\index{binary product}\index{product!binary}) of the objects $A$ and $B$ in the category $\mathcal{C}$, if for any object $C\in{\rm Ob}\,\mathcal{C}$ and any morphisms $f:C\to A$ and $g:C\to B$ there exists a unique morphism $h$ such that the following diagram commutes,
\[
\xymatrix{
& C \ar[dl]_f \ar@{.>}[d]|-h \ar[dr]^g & \\
A & A\,\Pi\, B \ar[l]^{\pi_A} \ar[r]_{\pi_B} & B
}
\]
We will also denote the mediating morphism $h:C\to A\,\Pi\, B$ by $[f,g]$.
\end{definition}
Often one omits the morphisms $\pi_A$ and $\pi_B$ and simply calls $A\,\Pi\,B$ the product of $A$ and $B$. The product of two objects is sometimes also denoted by $A\times B$.

\begin{proposition}\label{prop-prod-bifunctor}
\begin{itemize}
\item[(a)]
The product of two objects is unique up to isomorphism, if it exists.
\item[(b)]
Let $f_1:A_1\to B_1$ and $f_2:A_2\to B_2$ be two morphisms in a category $\mathcal{C}$ and assume that the products $A_1\,\Pi\,A_2$ and $B_1\,\Pi\,B_2$ exist in $\mathcal{C}$. Then there exists a unique morphism $f_1\,\Pi\,f_2:A_1\,\Pi\,A_2\to B_1\,\Pi\,B_2$ such that the following diagram commutes,
\[
\xymatrix{
& A_1 \ar[r]^{f_1}& B_1 & \\
A_1\,\Pi\,A_2 \ar[ur]^{\pi_{A_1}} \ar[dr]_{\pi_{A_2}} \ar@{.>}[rrr]|-{f_1\,\Pi\,f_2} &&& B_1\,\Pi\,B_2 \ar[ul]_{\pi_{B_1}} \ar[dl]^{\pi_{B_2}}\\
& A_2 \ar[r]_{f_2}& B_2 &
}
\]
\item[(c)]
Let $A_1,A_2,B_1,B_2,C_1,C_2$ be objects of a category  $\mathcal{C}$ and suppose that the products $A_1\,\Pi\,A_2$, $B_1\,\Pi\,B_2$ and $C_1\,\Pi\,C_2$ exist in $\mathcal{C}$. Then we have
\[
{\rm id}_{A_1}\,\Pi\,{\rm id}_{A_2} = {\rm id}_{A_1\,\Pi\,A_2} \text{ and } (g_1\,\Pi\,g_2)\circ(f_1\,\Pi\,f_2) = (g_1\circ f_1)\,\Pi\,(g_2\circ f_2)
\]
for all morphisms $f_i:A_i\to B_i$, $g_i:B_i\to C_i$, $i=1,2$.
\end{itemize}
\end{proposition}
\begin{proof}
\begin{itemize}
\item[(a)]
Suppose we have two candidates $(P,\pi_A,\pi_B)$ and $(P',\pi'_A,\pi'_B)$ for the product of $A$ and $B$, we have to show that $P$ and $P'$ are isomorphic. Applying the defining property of the product to  $(P,\pi_A,\pi_B)$ with $C=P'$ and to $(P',\pi'_A,\pi'_B)$ with $C=P'$, we get the following two commuting diagrams,
\[
\xymatrix{
& P' \ar[dl]_{\pi'_A} \ar@{.>}[d]|-h \ar[dr]^{\pi'_B} & \\
A & P \ar[l]^{\pi_A} \ar[r]_{\pi_B} & B
}
\qquad
\xymatrix{
& P \ar[dl]_{\pi_A} \ar@{.>}[d]|-{h'} \ar[dr]^{\pi_B} & \\
A & P' \ar[l]^{\pi'_A} \ar[r]_{\pi'_B} & B
}
\]
We get $\pi_A\circ h \circ h'= \pi'_A\circ h'=\pi_A$ and $\pi_B\circ h\circ h'= \pi'_B\circ h'=\pi_B$, i.e.\ the diagram
\[
\xymatrix{
& P \ar[dl]_{\pi_A} \ar@{.>}[d]|-{h\circ h'} \ar[dr]^{\pi_B} & \\
A & P \ar[l]^{\pi_A} \ar[r]_{\pi_B} & B
}
\]
is commutative. It is clear that this diagram also commutes, if we replace $h\circ h'$ by ${\rm id}_P$, so the uniqueness implies $h\circ h'={\rm id}_P$. Similarly one proves $h'\circ h={\rm id}_{P'}$, so that $h:P'\to P$ is the desired isomorphism.
\item[(b)]
The unique morphism $f_1\,\Pi\,f_2$ exists by the defining property of the product of $B_1$ and $B_2$, as we can see from the diagram
\[
\xymatrix{
& A_1\,\Pi\,A_2 \ar[dl]_{f_1\circ\pi_{A_1}} \ar@{.>}[d]|-{f_1\,\Pi\,f_2} \ar[dr]^{f_2\circ\pi_{A_2}} & \\
B_1 & B_1\,\Pi\,B_2 \ar[l]^{\pi_{B_1}} \ar[r]_-{\pi_{B_2}} & B_2
}
\]
\item[(c)]
Both properties follow from the uniqueness of the mediating morphism in the defining property of the product. To prove ${\rm id}_{A_1}\,\Pi\,{\rm id}_{A_2} = {\rm id}_{A_1\,\Pi\,A_2}$ one has to show that both expressions make the diagram
\[
\xymatrix{
 & A_1\,\Pi\,A_2 \ar[dl]_{{\rm id}_{A_1}} \ar@{.>}[d] \ar[dr]^{{\rm id}_{A_2}} & \\
A_1 & A_1\,\Pi\,A_2 \ar[l]^{\pi_{A_1}} \ar[r]_-{\pi_{A_2}} & A_2
}
\]
commutative, for the the second equality one checks that $(g_1\,\Pi\,g_2)\circ(f_1\,\Pi\,f_2)$ and $(g_1\circ f_1)\,\Pi\,(g_2\circ f_2)$ both make the diagram
\[
\xymatrix{
 & A_1\,\Pi\,A_2 \ar[dl]_{g_1\circ f_1} \ar@{.>}[d] \ar[dr]^{g_2\circ f_2} & \\
C_1 & C_1\,\Pi\,C_2 \ar[l]^{\pi_{C_1}} \ar[r]_-{\pi_{C_2}} & C_2
}
\]
commutative.
\end{itemize}
\end{proof}

The notion of product extends also to more then two objects.

\begin{definition}
Let $(A_i)_{i\in I}$ be a family of objects of a category $\mathcal{C}$, indexed by some set $I$. The pair $\left(\prod_{i\in I} A_i, \left(\pi_j:\prod_{i\in I} A_i\to A_j\right)_{j\in I}\right)$ consisting of an object $\prod_{i\in I} A_i$ of $\mathcal{C}$ and a family of morphisms $\left(\pi_j:\prod_{i\in I} A_i\to A_j\right)_{j\in I}$ of $\mathcal{C}$ is a {\em product}\index{product} of the family $(A_i)_{i\in I}$ if for any object $C$ and any family of morphisms $(f_i:C\to A_i)_{i\in I}$ there exists a unique morphism $h:C\to \prod_{i\in I} A_i$ such that
\[
\pi_j \circ h = f_j, \quad \text{ for all }j\in I
\]
holds. The morphism $\pi_j:\prod_{i\in I} A_i\to A_j$ for $j\in I$ is called the $j$th {\em product projection}\index{product projection}\index{projection!product}. We will also write $[f_i]_{i\in I}$ for the morphism $h:C\to \prod_{i\in I}A_i$.
\end{definition}

An object $T$ of a category $\mathcal{C}$ is called {\em terminal}\index{terminal object}\index{object!terminal}, if for any object $C$ of $\mathcal{C}$ there exists a unique morphism from $C$ to $T$. A terminal object is unique up to isomorphism, if it exists. A product of the empty family is a terminal object.

\begin{remark}
Let $\mathcal{C}$ be a category that has finite products. Then the product is associative and commutative. More precisely, there exist natural isomorphisms $\alpha_{A,B,C}:A\,\Pi\, (B\,\Pi\, C)\to (A\,\Pi\,B)\Pi\,C$ and $\gamma_{A,B}:B\,\Pi\,A\to A\,\Pi\,B$ for all objects $A,B,C\in{\rm Ob}\,\mathcal{C}$.
\end{remark}

The notion {\em coproduct}\index{coproduct} is the dual of the product, i.e.\ $\left(\coprod_{i\in I} A_i, \left(\imath_j: A_j\to \coprod_{i\in I} A_i\right)_{j\in I}\right)$ is called a coproduct of the family $(A_i)_{i\in I}$ of objects in $\mathcal{C}$, if it is a product of the same family in the category  $\mathcal{C}^{\rm op}$. Formulated in terms of objects and morphisms of $\mathcal{C}$ only, this amounts to the following.

\begin{definition}
Let $(A_i)_{i\in I}$ be a family of objects of a category $\mathcal{C}$, indexed by some set $I$. The pair $\left(\coprod_{i\in I} A_i, \left(\imath_j:A_k\to\coprod_{i\in I} A_i\right)_{j\in I}\right)$ consisting of an object $\coprod_{i\in I} A_i$ of $\mathcal{C}$ and a family of morphisms $\left(\imath_j:A_j\to \coprod_{i\in I} A_i\right)_{j\in I}$ of $\mathcal{C}$ is a coproduct of the family $(A_i)_{i\in I}$ if for any object $C$ and any family of morphisms $(f_i: A_i\to C)_{i\in I}$ there exists a unique morphism $h: \coprod_{i\in I} A_i\to C$ such that
\[
h\circ \imath_j = f_j, \quad \text{ for all }j\in I
\]
holds. The morphism $\imath_j: A_j\to\prod_{i\in I} A_i$ for $j\in I$ is called the $j$th {\em coproduct injection}\index{coproduct injection}\index{injection!coproduct}. We will write $[f_i]_{i\in I}$ for the morphism $h:\prod_{i\in I}A_i\to C$.
\end{definition}

A coproduct of the empty family in $\mathcal{C}$ is an {\em initial object}\index{initial object}\index{object!initial}, i.e.\ an object $I$ such that for any object $A$ of $\mathcal{C}$ there exists exactly one morphism from $I$ to $A$.

It is straightforward to translate Proposition \ref{prop-prod-bifunctor} to its counterpart for the coproduct.

\begin{example}
In the trivial unital semigroup $(G=\{e\},\circ,e)$, viewed as a category (note that is is isomorphic to the discrete category over a set with one element) its only object $G$ is a terminal and initial object, and also a product and coproduct for any family of objects. The product projections and coproduct injections are given by the unique morphism $e$ of this category.

In any other unital semigroup there exist no initial or terminal objects and no binary or higher products or coproducts. 
\end{example}

\begin{example}
In the category $\eufrak{Set}$ a binary product of two sets $A$ and $B$ is given by their Cartesian product $A\times B$ (together with the obvious projections) and any set with one element is terminal. A coproduct of $A$ and $B$ is defined by their disjoint union $A\dot{\cup}B$ (together with the obvious injections) and the empty set is an initial object. Recall that we can define the disjoint union as $A\dot{\cup}B=(A\times\{A\})\cup(B\times\{B\})$.
\end{example}

The following example shall be used throughout this chapter and the following.

\begin{example}\label{III-exa-free}
The coproduct in the category of unital algebras $\eufrak{Alg}$ is the free product of $*$-algebras {\em with} identification of the units. Let us recall its defining universal property. Let $\{\mathcal{A}_k\}_{k\in I}$ be a family of unital $*$-algebras and $\coprod_{k\in I}\mathcal{A}_k$ their free product, with canonical inclusions $\{i_k:\mathcal{A}_k\to\coprod_{k\in I}\mathcal{A}_k\}_{k\in I}$. If $\mathcal{B}$ is any unital $*$-algebra, equipped with unital $*$-algebra homomorphisms $\{i'_k:\mathcal{A}_k\to\mathcal{B}\}_{k\in I}$, then there exists a unique unital $*$-algebra homomorphism $h:\coprod_{k\in I}\mathcal{A}_k\to\mathcal{B}$ such that
\[
h\circ i_k =i'_k, \qquad \mbox{ for all }\quad k\in I.
\]
It follows from the universal property that for any pair of unital $*$-algebra homomorphisms $j_1:\mathcal{A}_1\to\mathcal{B}_1$,  $j_2:\mathcal{A}_2\to\mathcal{B}_2$
there exists a unique unital $*$-algebra homomorphism $j_1\coprod j_2:\mathcal{A}_1\coprod\mathcal{A}_2\to\mathcal{B}_1\coprod\mathcal{B}_2 $ such that the diagram
\[
\xymatrix{
& \mathcal{A}_1\ar[dl]_{i_{\mathcal{A}_1}} \ar[r]^{j_1}& \mathcal{B}_1 \ar[dr]^{i_{\mathcal{B}_1}} & \\
\mathcal{A}_1\coprod\mathcal{A}_2 \ar@{.>}[rrr]|-{j_1\coprod j_2} &&& \mathcal{B}_1\coprod\mathcal{B}_2\\
& \mathcal{A}_2 \ar[ul]^{i_{\mathcal{A}_2}}\ar[r]_{j_2}& \mathcal{B}_2\ar[ur]_{i_{\mathcal{B}_2}} &
}
\]
commutes.

The free product $\coprod_{k\in I}\mathcal{A}_k$ can be constructed as a sum of tensor products of the $\mathcal{A}_k$, where neighboring elements in the product belong to different algebras. For simplicity, we illustrate this only for the case of the free product of two algebras. Let
\[
\mathbb{A}=\bigcup_{n\in\mathbb{N}}\{ \epsilon\in \{1,2\}^n| \epsilon_1\not=\epsilon_2\not=\cdots\not=\epsilon_n\}
\]
and decompose $\mathcal{A}_i=\mathbb{C}\mathbf{1}\oplus\mathcal{A}_i^0$, $i=1,2$, into a direct sum of vector spaces.
Then $\mathcal{A}_1\coprod \mathcal{A}_2$ can be constructed as
\[
\mathcal{A}_1\coprod \mathcal{A}_2 = \bigoplus_{\epsilon\in\mathbb{A}} \mathcal{A}^\epsilon,
\]
where $\mathcal{A}^\emptyset=\mathbb{C}$, $\mathcal{A}^\epsilon=\mathcal{A}^0_{\epsilon_1}\otimes \cdots \otimes\mathcal{A}^0_{\epsilon_n}$ for $\epsilon=(\epsilon_1,\ldots,\epsilon_n)$. The multiplication in $\mathcal{A}_1\coprod \mathcal{A}_2$ is inductively defined by
\[
(a_1\otimes \cdots\otimes  a_n)\cdot (b_1\otimes \cdots\otimes  b_m) =
\left\{
\begin{array}{lcl}
a_1\otimes \cdots\otimes  (a_n\cdot b_1)\otimes \cdots\otimes  b_m & \mbox{ if } \epsilon_n=\delta_1, \\
a_1\otimes \cdots\otimes  a_n\otimes b_1\otimes \cdots\otimes  b_m & \mbox{ if } \epsilon_n\not=\delta_1,
\end{array}
\right.
\]
for $a_1\otimes \cdots\otimes a_n\in \mathcal{A}^\epsilon$, $b_1\otimes \cdots\otimes b_m\in \mathcal{A}^\delta$. Note that in the case $\epsilon_n=\delta_1$ the product $a_n\cdot b_1$ is not necessarily in $\mathcal{A}_{\epsilon_n}^0$, but is in general a sum of a multiple of the unit of $\mathcal{A}_{\epsilon_n}$ and an element of $\mathcal{A}_{\epsilon_n}^0$. We have to identify $a_1\otimes \cdots a_{n-1}\otimes 1\otimes b_2\otimes \cdots b_m$ with $a_1\otimes \cdots \otimes a_{n-1}\cdot b_2\otimes\cdots b_m$.

Since $\coprod$ is the coproduct of a category, it is commutative and associative in the sense that there exist natural isomorphisms
\begin{eqnarray}
\gamma_{\mathcal{A}_1,\mathcal{A}_2}&:& \mathcal{A}_1\coprod\mathcal{A}_2 \stackrel{\cong}{\to} \mathcal{A}_2\coprod\mathcal{A}_1, \label{eq-comm}\\
\alpha_{\mathcal{A}_1,\mathcal{A}_2,\mathcal{A}_3}&:& \mathcal{A}_1\coprod\left(\mathcal{A}_2\coprod\mathcal{A}_3\right) \stackrel{\cong}{\to} \left(\mathcal{A}_1\coprod\mathcal{A}_2\right)\coprod\mathcal{A}_3 \nonumber
\end{eqnarray}
for all unital $*$-algebras $\mathcal{A}_1,\mathcal{A}_2,\mathcal{A}_3$. Let $i_\ell:\mathcal{A}_\ell\to \mathcal{A}_1\coprod\mathcal{A}_2$ and  $i'_\ell:\mathcal{A}_\ell\to \mathcal{A}_2\coprod\mathcal{A}_1$, $\ell=1,2$ be the canonical inclusions. The commutativity constraint $\gamma_{\mathcal{A}_1,\mathcal{A}_2}: \mathcal{A}_1\coprod\mathcal{A}_2 \to \mathcal{A}_2\coprod\mathcal{A}_1$ maps an element of $\mathcal{A}_1\coprod\mathcal{A}_2$ of the form $i_1(a_1)i_2(b_1)\cdots i_2(b_n)$ with $a_1,\ldots,a_n\in\mathcal{A}_1$, $b_1,\ldots,b_n\in\mathcal{A}_2$ to
\[
\gamma_{\mathcal{A}_1,\mathcal{A}_2}\big(i_1(a_1)i_2(b_1)\cdots i_2(b_n)\big)=i'_1(a_1)i'_2(b_1)\cdots i'_2(b_n)\in\mathcal{A}_2\coprod\mathcal{A}_1.
\]
\end{example}

We also consider non-unital algebras. The free product of algebras {\em without} identification of units is a coproduct in the category $\eufrak{nuAlg}$ of non-unital (or rather not necessarily unital) algebras.

The following defines a a functor from the category of non-unital algebras $\eufrak{nuAlg}$ to the category of unital algebras $\eufrak{Alg}$. For an algebra $\mathcal{A}\in{\rm Ob}\,\eufrak{nuAlg}$, $\tilde{\mathcal{A}}$ is equal to $\tilde{\mathcal{A}}=\mathbb{C}\mathbf{1}\oplus\mathcal{A}$ as a vector space and the multiplication is defined by
\[
(\lambda\mathbf{1}+a)(\lambda'\mathbf{1}+a')=\lambda\lambda'\mathbf{1}+\lambda' a+\lambda a'+aa'
\]
for $\lambda,\lambda'\in\mathbb{C}$, $a,a'\in\mathcal{A}$. We will call $\tilde{\mathcal{A}}$ the {\em unitization} of $\mathcal{A}$. Note that $\mathcal{A}\cong 0\mathbf{1}+\mathcal{A}\subseteq\tilde{\mathcal{A}}$ is not only a subalgebra, but even an ideal in $\tilde{\mathcal{A}}$. 

The following relation holds between the free product {\em with} identification of units $\coprod_{\eufrak{Alg}}$ and the free product {\em without} identification of units $\coprod_{\eufrak{nuAlg}}$,
\[
\widetilde{\mathcal{A}_1\coprod_{\eufrak{nuAlg}}\mathcal{A}_2} \cong \tilde{\mathcal{A}}_1\coprod_{\eufrak{Alg}}\tilde{\mathcal{A}_2}
\]
for all $\mathcal{A}_1,\mathcal{A}_2\in{\rm Ob}\,\eufrak{nuAlg}$.

Note furthermore that the range of this functor consists of all algebras that admit a decomposition of the form $\mathcal{A}=\mathbb{C}\mathbf{1}\oplus\mathcal{A}_0$, where $\mathcal{A}_0$ is a subalgebra. This is equivalent to having a one-dimensional representation. The functor is not surjective, e.g., the algebra $\mathcal{M}_2$ of $2\times 2$-matrices can not be obtained as a unitization of some other algebra.

Let us now recall the definition of a tensor category.
\begin{definition}\label{def-tensor}
A category $(\mathcal{C},\square)$ equipped with a bifunctor $\square:\mathcal{C}\times\mathcal{C}\to\mathcal{C}$, called {\em tensor product}, that is associative up to a natural isomorphism
\[
\alpha_{A,B,C}:A\square(B\square C) \stackrel{\cong}{\to} (A\square B)\square C,\qquad \text{ for all } A,B,C\in{\rm Ob}\,\mathcal{C},
\]
 and an element $E$ that is, up to natural isomorphisms
\[
\lambda_A:E\square A\stackrel{\cong}{\to}A, \quad\text{ and }\quad \rho_A:A\square E\stackrel{\cong}{\to} A, \quad \text{ for all } A\in{\rm Ob}\,\mathcal{C},
\]
a unit for $\square$, is called a {\em tensor category} or {\em monoidal category}, if the {\em pentagon axiom}
\[
\xymatrix{
 & (A\square B)\square(C\square D)\ar[dr]^{\alpha_{A\square B,C,C}} & \\
A\square \big(B \square (C \square D)\big) \quad\ar[ur]^{\alpha_{A,B,C\square D}}\ar[d]_{{\rm id}_A\square\alpha_{A,B,C}} & & \quad\big((A\square B) \square C\big) \square D \\
  \qquad A\square \big((B \square C) \square D\big)\ar[rr]_{\alpha_{A,B\square C.D}} & & \big(A\square (B \square C)\big) \square D \ar[u]_{\alpha_{A,B,C}\square{\rm id}_D} \qquad
}
\]
and the {\rm triangle axiom}
\[
\xymatrix{
A\square (E\square C)\ar[rr]^{\alpha_{A,E,C}}\ar[dr]_{{\rm id}_A\square \lambda_C} & & (A\square E)\square C\ar[dl]^{\rho_A\square{\rm id}_C} \\
& A\square C & \\
}
\]
are satisfied for all objects $A,B,C,D$ of $\mathcal{C}$.
\end{definition}

If a category has products or coproducts for all finite sets of objects, then the universal property guarantees the existence of the isomorphisms $\alpha$, $\lambda$, and $\rho$ that turn it into a tensor category.

A functor between tensor categories that behaves ``nicely'' with respect to the tensor products, is called a tensor functor or monoidal functor, see, e.g., Section XI.2 in MacLane\cite{maclane98}.

\begin{definition}
Let $(\mathcal{C},\square)$ and $(\mathcal{C}',\square')$ be two tensor categories.
A {\em cotensor functor} or {\em comonoidal functor} $F:(\mathcal{C},\square)\to(\mathcal{C}',\square')$ is an ordinary functor $F:\mathcal{C}\to\mathcal{C}'$ equipped with a morphism $F_0: F(E_\mathcal{C}) \to E_{\mathcal{C}'}$ and a natural transformation $F_2:F(\,\cdot\, \square \,\cdot\,) \to  F(\,\cdot\,) \square' F(\,\cdot\,)$, i.e.\ morphisms $F_2(A,B):F(A \square B) \to  F(A) \square' F(B)$ for all $A,B\in{\rm Ob}\,\mathcal{C}$ that are natural in $A$ and $B$, such that the diagrams
\begin{equation}\label{III-cotensor-alpha}
\xymatrix{
F\big(A\square (B\square C)\big)\ar[rr]^{F(\alpha_{A,B,C})} \ar[d]_{F_2(A,B\square C)} && F\big((A\square B)\square C\big)\ar[d]^{F_2(A\square B,C)} \\
F(A)\square' F(B\square C) \ar[d]_{{\rm id}_{F(A)}\square'F_2(B,C)} && F(A\square B)\square' F(C)\ar[d]^{F_2(A,B)\square'{\rm id}_{F(C)}} \\
F(A)\square'\big(F(B)\square'F(C)\big) \ar[rr]_{\alpha'_{F(A),F(B),F(C)}}&& \big(F(A)\square' F(B)\big)\square' F(C)
}
\end{equation}
\begin{equation}\label{III-cotensor-rho}
\xymatrix{
F(B\square E_\mathcal{C}) \ar[rr]^{F_2(B,E_\mathcal{C})} \ar[d]_{F(\rho_B)} && F(B)\square' F(E_\mathcal{C})\ar[d]^{{\rm id}_B\square' F_0} \\
F(B) && F(B)\square'E_{\mathcal{C}'} \ar[ll]^{\rho'_{F(B)}}
}
\end{equation}
\begin{equation}
\xymatrix{
F(E_\mathcal{C}\square B) \ar[rr]^{F_2(E_\mathcal{C},B)} \ar[d]_{F(\lambda_B)} && F(E_\mathcal{C})\square'F(B) \ar[d]^{F_0\square'{\rm id}_B} \\
F(B) && E_{\mathcal{C}'} \square'F(B) \ar[ll]^{\lambda'_{F(B)}}
}
\end{equation}
commute for all $A,B,C\in{\rm Ob}\,\mathcal{C}$.
\end{definition}
We have reversed the direction of $F_0$ and $F_2$ in our definition. In the case of a strong tensor functor, i.e.\ when all the morphisms are isomorphisms, our definition of a cotensor functor is equivalent to the usual definition of a tensor functor as, e.g., in MacLane\cite{maclane98}.

The conditions are exactly what we need to get morphisms
\[
F_n(A_1,\ldots,A_n):F(A_1\square\cdots\square A_n)\to F(A_1)\square'\cdots\square' F(A_n)
\]
for all finite sets $\{A_1,\ldots,A_n\}$ of objects of $\mathcal{C}$ such that, up to these morphisms, the functor $F:(\mathcal{C},\square)\to(\mathcal{C}',\square')$ is a homomorphism.

\section{Classical stochastic independence and the product of
  probability spaces}

Two random variables $X_1:(\Omega,\mathcal{F},P)\to(E_1,\mathcal{E}_1)$ and
$X_2:(\Omega,\mathcal{F},P)\to(E_2,\mathcal{E}_2)$, defined on the same
probability space $(\Omega,\mathcal{F},P)$ and with values in two possibly
distinct measurable spaces $(E_1,\mathcal{E}_1)$ and $(E_2,\mathcal{E}_2)$,
are called {\em stochastically independent}\index{independent!stochastically} (or simply {\em independent}) w.r.t.\ $P$, if the $\sigma$-algebras $X_1^{-1}(\mathcal{E}_1)$ and $X_2^{-1}(\mathcal{E}_2)$ are independent w.r.t.\ $P$, i.e.\ if
\[
P\big((X^{-1}_1(M_1)\cap X^{-1}_2(M_2)\big) =P\big((X^{-1}_1(M_1)\big)P\big( X^{-1}_2(M_2)\big) 
\]
holds for all $M_1\in\mathcal{E}_1$, $M_2\in\mathcal{E}_2$. If there is no
danger of confusion, then the reference to the measure $P$ is often omitted.

This definition can easily be extended to arbitrary families of random
variables. A family
$\big(X_j:(\Omega,\mathcal{F},P)\to(E_j,\mathcal{E}_j))_{j\in J}$, indexed by
some set $J$, is called independent, if
\[
P\left(\bigcap_{k=1}^n(X^{-1}_{j_k}(M_{j_k})\right) = \prod_{k=1}^n
P\big(X^{-1}_{j_k}(M_{j_k})\big)
\]
holds for all $n\in\mathbb{N}$ and all choices of indices $k_1,\ldots,k_n\in
J$ with $j_k\not=j_\ell$ for $j\not=\ell$, and all choices of measurable sets
$M_{j_k}\in\mathcal{E}_{j_k}$.

There are many equivalent formulations for independence, consider, e.g., the
following proposition.

\begin{proposition}
Let $X_1$ and $X_2$ be two real-valued random variables. The following are
equivalent.
\begin{itemize}
\item[(i)]
$X_1$ and $X_2$ are independent.
\item[(ii)] For all bounded measurable functions $f_1,f_2$ on $\mathbb{R}$ we
  have
\[
\mathbb{E}\big(f_1(X_1)f_2(X_2)\big)=
\mathbb{E}\big(f_1(X_1)\big)\mathbb{E}\big(f_2(X_2)\big).
\]
\item[(iii)]
The probability space $(\mathbb{R}^2,\mathcal{B}(\mathbb{R}^2),P_{(X_1,X_2)})$ is the
product of the probability spaces
$(\mathbb{R},\mathcal{B}(\mathbb{R}),P_{X_1})$ and
$(\mathbb{R},\mathcal{B}(\mathbb{R}),P_{X_2})$, i.e.
\[
P_{(X_1,X_2)} = P_{X_1}\otimes P_{X_2}.
\]
\end{itemize}
\end{proposition}

We see that stochastic independence can be reinterpreted as a rule to compute
the joint distribution of two random variables from their marginal
distribution. More precisely, their joint distribution can be computed as a
product of their marginal distributions. This product is associative and can
also be iterated to compute the joint distribution of more than two independent
random variables.

The classifications of independence for non-commutative probability
\cite{speicher96,benghorbal+schuermann99,benghorbal01,muraki02a,muraki02}
that we are interested in are based on redefining independence as a
product satisfying certain natural axioms.

\section{Definition of independence in the language of category theory}

We will now define the notion of independence in the language of category theory. The usual notion of independence for classical probability theory and the independences classified in \cite{speicher96,benghorbal+schuermann99,benghorbal01,muraki02a,muraki02} will then be instances of this general notion obtained by considering the category of classical probability spaces or categories of algebraic probability spaces.

In order to define a notion of independence we need less than a (co-) product, but a more than a tensor product. What we need are inclusions or projections that allow us to view the objects $A$, $B$ as subsystems of their product $A\square B$.

\begin{definition}\label{def-tensor-cat-with}
A {\em tensor category with projections} $(\mathcal{C},\square,\pi)$ is a tensor category $(\mathcal{C},\square)$ equipped with two natural transformations $\pi_1:\square \to P_1$ and $\pi_2:\square \to P_2$, where the bifunctors $P_1,P_2:\mathcal{C}\times\mathcal{C}\to\mathcal{C}$ are defined by $P_1(B_1,B_2)=B_1$, $P_2(B_1,B_2)=B_2$, on pairs of objects $B_1,B_2$ of $\mathcal{C}$, and similarly on pairs of morphisms. In other words, for any pair of objects $B_1,B_2$ there exist two morphisms $\pi_{B_1}:B_1\square B_2\to B_1$, $\pi_{B_2}:B_1\square B_2\to B_2$, such that for any pair of morphisms $f_1:A_1\to B_1$, $f_2:A_2\to B_2$, the following diagram commutes,
\[
\xymatrix{
A_1 \ar[d]_{f_1}& A_1\square A_2 \ar[d]|-{f_1\square f_2} \ar[l]_{\pi_{A_1}} \ar[r]^{\pi_{A_2}} & A_2\ar[d]^{f_2} \\
B_1 & B_1\square B_2\ar[l]^{\pi_{B_1}} \ar[r]_{\pi_{B_2}} & B_2 .
}
\]

Similarly, a {\em tensor product with inclusions}
$(\mathcal{C},\square,i)$ is a tensor category $(\mathcal{C},\square)$
equipped with two natural transformations $i_1:P_1\to\square$ and
$i_2:P_2\to\square$,  i.e.\ for any pair of objects $B_1,B_2$ there exist two morphisms $i_{B_1}:B_1\to B_1\square B_2$, $i_{B_2}:B_2\to B_1\square B_2$, such that for any pair of morphisms $f_1:A_1\to B_1$, $f_2:A_2\to B_2$, the following diagram commutes,
\[
\xymatrix{
A_1 \ar[r]_{i_{A_1}} \ar[d]_{f_1}& A_1\square A_2 \ar[d]|-{f_1\square f_2} & A_2\ar[d]^{f_2} \ar[l]^{i_{A_2}} \\
B_1 \ar[r]^{i_{B_1}}& B_1\square B_2 & B_2\ar[l]_{i_{B_2}} .
}
\]
\end{definition}

In  a tensor category with projections or with inclusions we can define a notion of independence for morphisms.

\begin{definition}\label{def-independence}
Let $(\mathcal{C},\square,\pi)$ be a tensor category with
projections. Two morphism $f_1:A\to B_1$ and $f_2:A\to B_2$ with the same
source $A$ are called
{\em independent}\index{independent} (with respect to $\square$), if there exists a morphism $h:A\to B_1\square B_2$ such that the
diagram
\begin{equation}\label{diag-indep1}
\xymatrix{
& A \ar[dl]_{f_1} \ar@{.>}[d]|-h\ar[dr]^{f_2} & \\
B_1 & B_1\square B_2 \ar[l]^{\pi_{B_1}} \ar[r]_{\pi_{B_2}} & B_2
}
\end{equation}
commutes.

In a tensor category with inclusions $(\mathcal{C},\square,i)$, two morphisms $f_1:B_1\to A$ and $f_2:B_2\to A$ with the same
target $B$ are called independent, if there exists a morphism $h:B_1\square
B_2\to A$ such that the diagram
\begin{equation}\label{diag-indep2}
\xymatrix{
& A  & \\
B_1 \ar[ur]^{f_1}\ar[r]_{i_{B_1}}& B_1\square B_2 \ar@{.>}[u]|-h & B_2\ar[l]^{i_{B_2}}\ar[ul]_{f_2}
}
\end{equation}
commutes.
\end{definition}
This definition can be extended in the obvious way to arbitrary sets of morphisms.

If $\square$ is actually a product (or coproduct, resp.), then the universal property in Definition \ref{III-def-product} implies that for all pairs of morphisms with the same source (or target,
resp.) there exists even a unique morphism
that makes diagram \eqref{diag-indep1} (or \eqref{diag-indep2}, resp.)
commuting. Therefore in that case
all pairs of morphism with the same source (or target, resp.) are
independent.

We will now consider several examples. We will show that for the category of classical probability spaces we recover usual stochastic independence, if we take the product of probability spaces, cf.\ Proposition \ref{prop-stoch-independence}.

\subsection{Example: Independence in the Category of Classical Probability Spaces}\label{example-class-indep}

The category $\eufrak{Meas}$ of measurable spaces consists of pairs $(\Omega,\mathcal{F})$, where $\Omega$ is a set and $\mathcal{F}\subseteq\mathcal{P}(\Omega)$ a $\sigma$-algebra. The morphisms are the measurable maps. This category has a product,
\[
(\Omega_1,\mathcal{F}_1)\,\Pi\,(\Omega_2,\mathcal{F}_2)=(\Omega_1\times\Omega_2,\mathcal{F}_1\otimes\mathcal{F}_2)
\]
where $\Omega_1\times\Omega_2$ is the Cartesian product of $\Omega_1$ and $\Omega_2$, and $\mathcal{F}_1\otimes\mathcal{F}_2$ is the smallest $\sigma$-algebra on $\Omega_1\times\Omega_2$ such that the canonical projections $p_1:\Omega_1\times\Omega_2\to\Omega_1$ and $p_2:\Omega_1\times\Omega_2\to\Omega_2$ are measurable.

The category of probability spaces $\eufrak{Prob}$ has as objects triples $(\Omega,\mathcal{F},P)$ where $(\Omega,\mathcal{F})$ is a measurable space and $P$ a probability measure on $(\Omega,\mathcal{F})$. A morphism $X:(\Omega_1,\mathcal{F}_1,P_1)\to (\Omega_1,\mathcal{F}_2,P_2)$ is a measurable map $X:(\Omega_1,\mathcal{F}_1)\to (\Omega_1,\mathcal{F}_2)$ such that
\[
P_1\circ X^{-1} = P_2.
\]
This means that a random variable $X:(\Omega,\mathcal{F})\to(E,\mathcal{E})$ automatically becomes a morphism, if we equip $(E,\mathcal{E})$ with the measure
\[
P_X=P\circ X^{-1}
\]
induced by $X$.

This category does not have universal products. But one can check that the product of measures turns $\eufrak{Prob}$ into a tensor category,
\[
(\Omega_1,\mathcal{F}_1,P_1)\otimes(\Omega_2,\mathcal{F}_2,P_2)=(\Omega_1\times\Omega_2,\mathcal{F}_1\otimes\mathcal{F}_2,P_1\otimes P_2),
\]
where $P_1\otimes P_2$ is determined by
\[
(P_1\otimes P_2)(M_1\times M_2) =P_1(M_1)P_2(M_2),
\]
for all $M_1\in\mathcal{F}_1$, $M_2\in\mathcal{F}_2$. It is even a tensor category with projections in the sense of Definition \ref{def-tensor-cat-with} with the canonical projections $p_1:(\Omega_1\times\Omega_2,\mathcal{F}_1\otimes\mathcal{F}_2,P_1\otimes P_2)\to (\Omega_1,\mathcal{F}_1,P_1)$, $p_2:(\Omega_1\times\Omega_2,\mathcal{F}_1\otimes\mathcal{F}_2,P_1\otimes P_2)\to (\Omega_2,\mathcal{F}_2,P_2)$ given by $p_1\big((\omega_1,\omega_2)\big)=\omega_1$, $p_2\big((\omega_1,\omega_2)\big)=\omega_2$ for $\omega_1\in\Omega_1$, $\omega_2\in\Omega_2$.

The notion of independence associated to this tensor product with projections is exactly the one used in probability.

\begin{proposition}\label{prop-stoch-independence}
Two random variables $X_1:(\Omega,\mathcal{F},P)\to(E_1,\mathcal{E}_1)$ and $X_2:(\Omega,\mathcal{F},P)\to(E_2,\mathcal{E}_2)$, defined on the same probability space $(\Omega,\mathcal{F},P)$ and with values in measurable spaces $(E_1,\mathcal{E}_1)$ and $(E_2,\mathcal{E}_2)$, are stochastically independent, if and only if they are independent in the sense of Definition \ref{def-independence} as morphisms $X_1:(\Omega,\mathcal{F},P)\to(E_1,\mathcal{E}_1,P_{X_1})$ and $X_2:(\Omega,\mathcal{F},P)\to(E_2,\mathcal{E}_2,P_{X_2})$ of the tensor category with projections $(\eufrak{Prob},\otimes,p)$.
\end{proposition}
\begin{proof}
Assume that $X_1$ and $X_2$ are stochastically independent. We have to find a morphism $h:(\Omega,\mathcal{F},P)\to(E_1\times E_2,\mathcal{E}_1\otimes\mathcal{E}_2,P_{X_1}\otimes P_{X_2})$ such that the diagram
\[
\xymatrix{
& (\Omega,\mathcal{F},P) \ar[dl]_{X_1} \ar@{.>}[d]|-h\ar[dr]^{X_2} & \\
(E_1,\mathcal{E}_1,P_{X_1}) &  (E_1\times E_2,\mathcal{E}_1\otimes\mathcal{E}_2,P_{X_1}\otimes P_{X_2})\ar[l]^(.65){p_{E_1}} \ar[r]_(.65){p_{E_2}} & (E_2,\mathcal{E}_2,P_{X_2})
}
\]
commutes. The only possible candidate is $h(\omega)=\big(X_1(\omega),X_2(\omega)\big)$ for all $\omega\in\Omega$, the unique map that completes this diagram in the category of measurable spaces and that exists due to the universal property of the product of measurable spaces. This is a morphism in $\eufrak{Prob}$, because we have
\begin{eqnarray*}
P\big(h^{-1}(M_1\times M_2)\big)&=&P\big(X_1^{-1}(M_1)\cap X_2^{-1}(M_2)\big) = P\big(X_1^{-1}(M_1)\big)P\big(X_2^{-1}(M_2)\big)\\
&=& P_{X_1}(M_1)P_{X_2}(M_2)=(P_{X_1}\otimes P_{X_2})(M_1\times M_2)
\end{eqnarray*}
for all $M_1\in\mathcal{E}_1$, $M_2\in\mathcal{E}_2$, and therefore
\[
P\circ h^{-1}=P_{X_1}\otimes P_{X_2}.
\]
Conversely, if $X_1$ and $X_2$ are independent in the sense of Definition \ref{def-independence}, then the morphism that makes the diagram commuting has to be again $h:\omega\mapsto \big(X_1(\omega),X_2(\omega)\big)$. This implies
\[
P_{(X_1,X_2)}=P\circ h^{-1} = P_{X_1}\otimes P_{X_2}
\]
and therefore
\[
P\big(X_1^{-1}(M_1)\cap X_2^{-1}(M_2)\big) =  P\big(X_1^{-1}(M_1)\big)P\big(X_2^{-1}(M_2)\big)
\]
for all $M_1\in\mathcal{E}_1$, $M_2\in\mathcal{E}_2$.
\end{proof}

\subsection{Example: Tensor Independence in the Category of Algebraic Probability Spaces}\label{tensor indep}

By the category of algebraic probability spaces $\eufrak{AlgProb}$ we denote the category of associative unital algebras over $\mathbb{C}$ equipped with a unital linear functional. A morphism $j:(\mathcal{A}_1,\varphi_1)\to(\mathcal{A}_2,\varphi_2)$ is a quantum random variable, i.e.\ an algebra homomorphism $j:\mathcal{A}_1\to\mathcal{A}_2$ that preserves the unit and the functional, i.e.\ $j(\mathbf{1}_{\mathcal{A}_1})=\mathbf{1}_{\mathcal{A}_2}$ and $\varphi_2\circ j=\varphi_1$.

The tensor product we will consider on this category is just the usual tensor product $(\mathcal{A}_1\otimes\mathcal{A}_2,\varphi_1\otimes\varphi_2)$, i.e.\ the algebra structure of $\mathcal{A}_1\otimes\mathcal{A}_2$ is defined by
\begin{eqnarray*}
\mathbf{1}_{\mathcal{A}_1\otimes\mathcal{A}_2} &=& \mathbf{1}_{\mathcal{A}_1}\otimes\mathbf{1}_{\mathcal{A}_2}, \\
(a_1\otimes a_2)(b_1\otimes b_1) &=& a_1b_2 \otimes a_2b_2,
\end{eqnarray*}
and the new functional is defined by
\[
(\varphi_1\otimes\varphi_2)(a_1\otimes a_2) = \varphi_1(a_1)\varphi_2(a_2),
\]
for all $a_1,b_1\in\mathcal{A}_1$, $a_2,b_2\in\mathcal{A}_2$.

This becomes a tensor category with inclusions with the inclusions defined by
\begin{eqnarray*}
i_{\mathcal{A}_1}(a_1) &=& a_1\otimes \mathbf{1}_{\mathcal{A}_2}, \\
i_{\mathcal{A}_2}(a_2) &=& \mathbf{1}_{\mathcal{A}_1}\otimes a_2,
\end{eqnarray*}
for $a_1\in\mathcal{A}_1$, $a_2\in\mathcal{A}_2$.

One gets the category of $*$-algebraic probability spaces, if one assumes that the underlying algebras have an involution and the functional are states, i.e.\ also positive. Then an involution is defined on $\mathcal{A}_1\otimes\mathcal{A}_2$ by $(a_1\otimes a_2)^*= a_1^*\otimes a_2^*$ and $\varphi_1\otimes\varphi_2$ is again a state.

The notion of independence associated to this tensor product with inclusions by Definition \ref{def-independence} is the usual notion of {\em Bose} or {\em tensor independence} used in quantum probability, e.g., by Hudson and Parthasarathy. See also definition \ref{ch1-def-indep} in Chapter \ref{levy-intro}.

\begin{proposition}
Two quantum random variables $j_1:(\mathcal{B}_1,\psi_1)\to(\mathcal{A},\varphi)$ and $j_2:(\mathcal{B}_2,\psi_2)\to(\mathcal{A},\varphi)$, defined on algebraic probability spaces $(\mathcal{B}_1,\psi_1),(\mathcal{B}_2,\psi_2)$ and with values in the same algebraic probability space $(\mathcal{A},\varphi)$ are independent if and only if the following two conditions are satisfied.
\begin{description}
\item[(i)]
The images of $j_1$ and $j_2$ commute, i.e.
\[
\big[j_1(a_1),j_2(a_2)\big] = 0,
\]
for all $a_1\in\mathcal{A}_1$, $a_2\in\mathcal{A}_2$.
\item[(ii)]
$\varphi$ satisfies the factorization property
\[
\varphi\big(j_1(a_1)j_2(a_2)\big) = \varphi\big(j_1(a_1)\big)\varphi\big(j_2(a_2)\big),
\]
for all $a_1\in\mathcal{A}_1$, $a_2\in\mathcal{A}_2$.
\end{description}
\end{proposition}
We will not prove this Proposition since it can be obtained as a special case of Proposition \ref{prop-fermi-indep}, if we equip the algebras with the trivial $\mathbb{Z}_2$-grading $\mathcal{A}^{(0)}=\mathcal{A}$, $\mathcal{A}^{(1)}=\{0\}$.

\subsection{Example: Fermi Independence}\label{fermi-indep}

Let us now consider the category of $\mathbb{Z}_2$-graded algebraic probability spaces $\mathbb{Z}_2$-$\eufrak{AlgProb}$. The objects are pairs $(\mathcal{A},\varphi)$ consisting of a $\mathbb{Z}_2$-graded unital algebra $\mathcal{A}=\mathcal{A}^{(0)}\oplus\mathcal{A}^{(1)}$ and an even unital functional $\varphi$, i.e.\ $\varphi|_{\mathcal{A}^{(1)}}=0$. The morphisms are random variables that don't change the degree, i.e., for $j:(\mathcal{A}_1,\varphi_1)\to(\mathcal{A}_2,\varphi_2)$, we have
\[
j(\mathcal{A}_1^{(0)})\subseteq\mathcal{A}_2^{(0)} \quad \mbox{ and }\quad j(\mathcal{A}_1^{(1)})\subseteq\mathcal{A}_2^{(1)}.
\]
The tensor product $(\mathcal{A}_1\otimes_{\mathbb{Z}_2}\mathcal{A}_2,\varphi_1\otimes\varphi_2)= (\mathcal{A}_1,\varphi_1)\otimes_{\mathbb{Z}_2}(\mathcal{A}_2,\varphi_2)$ is defined as follows. The algebra $\mathcal{A}_1\otimes_{\mathbb{Z}_2}\mathcal{A}_2$ is the graded tensor product of $\mathcal{A}_1$ and $\mathcal{A}_2$, i.e.\ $(\mathcal{A}_1\otimes_{\mathbb{Z}_2}\mathcal{A}_2)^{(0)}=\mathcal{A}^{(0)}_1\otimes\mathcal{A}^{(0)}_2\oplus \mathcal{A}^{(1)}_1\otimes\mathcal{A}^{(1)}_2$,  $(\mathcal{A}_1\otimes_{\mathbb{Z}_2}\mathcal{A}_2)^{(1)}=\mathcal{A}^{(1)}_1\otimes\mathcal{A}^{(0)}_2\oplus \mathcal{A}^{(0)}_1\otimes\mathcal{A}^{(1)}_2$, with the algebra structure given by
\begin{eqnarray*}
\mathbf{1}_{\mathcal{A}_1\otimes_{\mathbb{Z}_2}\mathcal{A}_2} &=& \mathbf{1}_{\mathcal{A}_1} \otimes\mathbf{1}_{\mathcal{A}_2}, \\
(a_1\otimes a_2) \cdot (b_1\otimes b_2) &=& (-1)^{\deg a_2\deg b_1} a_1b_1\otimes a_2b_2, 
\end{eqnarray*}
for all homogeneous elements $a_1,b_1\in\mathcal{A}_1$, $a_2,b_2\in\mathcal{A}_2$. The functional $\varphi_1\otimes\varphi_2$ is simply the tensor product, i.e.\ $(\varphi_1\otimes\varphi_2)(a_1\otimes a_2)=\varphi_1(a_1)\otimes \varphi_2(a_2)$ for all $a_1\in\mathcal{A}_1$, $a_2\in\mathcal{A}_2$. It is easy to see that $\varphi_1\otimes\varphi_2$ is again even, if $\varphi_1$ and $\varphi_2$ are even. The inclusions $i_1:(\mathcal{A}_1,\varphi_1)\to(\mathcal{A}_1\otimes_{\mathbb{Z}_2}\mathcal{A}_2,\varphi_1\otimes\varphi_2)$ and  $i_2:(\mathcal{A}_2,\varphi_2)\to(\mathcal{A}_1\otimes_{\mathbb{Z}_2}\mathcal{A}_2,\varphi_1\otimes\varphi_2)$ are defined by
\[
i_1(a_1) = a_1\otimes \mathbf{1}_{\mathcal{A}_2}\quad\mbox{ and } i_2(a_2)= \mathbf{1}_{\mathcal{A}_1}\otimes a_2,
\]
for $a_1\in\mathcal{A}_1$, $a_2\in\mathcal{A}_2$.

If the underlying algebras are assumed to have an involution and the functionals to be states, then the involution on the $\mathbb{Z}_2$-graded tensor product is defined by $(a_1\otimes a_2)^* = (-1)^{\deg a_1 \deg a_2} a_1^*\otimes a_2^*$, this gives the category of $\mathbb{Z}_2$-graded $*$-algebraic probability spaces.

The notion of independence associated to this tensor category with inclusions is called {\em Fermi independence} or {\em anti-symmetric independence}.

\begin{proposition}\label{prop-fermi-indep}
Two random variables
$j_1:(\mathcal{B}_1,\psi_1)\to(\mathcal{A},\varphi)$ and
$j_2:(\mathcal{B}_2,\psi_2)\to(\mathcal{A},\varphi)$, defined on two
$\mathbb{Z}_2$-graded algebraic probability spaces
$(\mathcal{B}_1,\psi_1)$, $(\mathcal{B}_2,\psi_2)$ and with values in the same $\mathbb{Z}_2$-algebraic probability space $(\mathcal{A},\varphi)$ are independent if and only if the following two conditions are satisfied.
\begin{description}
\item[(i)]
The images of $j_1$ and $j_2$ satisfy the commutation relations
\[
j_2(a_2)j_1(a_1) = (-1)^{\deg a_1 \deg a_2}j_1(a_1)j_2(a_2)
\]
for all homogeneous elements $a_1\in\mathcal{B}_1$, $a_2\in\mathcal{B}_2$.
\item[(ii)]
$\varphi$ satisfies the factorization property
\[
\varphi\big(j_1(a_1)j_2(a_2)\big) = \varphi\big(j_1(a_1)\big)\varphi\big(j_2(a_2)\big),
\]
for all $a_1\in\mathcal{B}_1$, $a_2\in\mathcal{B}_2$.
\end{description}
\end{proposition}
\begin{proof}
The proof is similar to that of Proposition \ref{prop-stoch-independence}, we will only outline it. It is clear that the morphism $h:(\mathcal{B}_1,\psi_1)\otimes_{\mathbb{Z}_2}(\mathcal{B}_2,\psi_2)\to(\mathcal{A},\varphi)$ that makes the diagram in Definition \ref{def-independence} commuting, has to act on elements of $\mathcal{B}_1\otimes \mathbf{1}_{\mathcal{B}_2}$ and $\mathbf{1}_{\mathcal{B}_1}\otimes\mathcal{B}_2$ as
\[
h(b_1\otimes \mathbf{1}_{\mathcal{B}_2}) = j_1(b_1) \quad\mbox{ and } \quad h(\mathbf{1}_{\mathcal{B}_1}\otimes b_2)=j_2(b_2).
\]
This extends to a homomorphism from $(\mathcal{B}_1,\psi_1)\otimes_{\mathbb{Z}_2}(\mathcal{B}_2,\psi_2)$ to $(\mathcal{A},\varphi)$, if and only if the commutation relations are satisfied. And the resulting homomorphism is a quantum random variable, i.e.\ satisfies $\varphi\circ h=\psi_1\otimes\psi_2$, if and only if the factorization property is satisfied.
\end{proof}

\subsection{Example: Free Independence}

We will now introduce another tensor product with inclusions for the category of algebraic probability spaces $\eufrak{AlgProb}$. On the algebras we take simply the free product of algebras with identifications of units introduced in Example \ref{III-exa-free}. This is the coproduct in the category of algebras, therefore we also have natural inclusions. It only remains to define a unital linear functional on the free product of the algebras.

Voiculescu's\cite{voiculescu+dykema+nica92} free product $\varphi_1*\varphi_2$ of two unital linear functionals $\varphi_1:\mathcal{A}_1\to\mathbb{C}$ and  $\varphi_2:\mathcal{A}_2\to\mathbb{C}$ can be defined recursively by
\[
(\varphi_1*\varphi_2)(a_1a_2\cdots a_m)= \sum_{I\subsetneqq \{1,\ldots,m\}} (-1)^{m-\sharp I+1} (\varphi_1*\varphi_2)\left(\prod_{k\in I}^{\rightarrow}a_k\right) \prod_{k\not\in I} \varphi_{\epsilon_k}(a_k)
\]
for a typical element $a_1a_2\cdots a_m\in\mathcal{A}_1\coprod\mathcal{A}_2$, with $a_k\in\mathcal{A}_{\epsilon_k}$, $\epsilon_1\not=\epsilon_2\not=\cdots\not=\epsilon_m$, i.e.\ neighboring $a$'s don't belong to the same algebra. $\sharp I$ denotes the number of elements of $I$ and $\prod_{k\in I}^{\rightarrow}a_k$ means that the $a$'s are to be multiplied in the same order in which they appear on the left-hand-side. We use the convention $(\varphi_1*\varphi_2)\left(\prod_{k\in \emptyset}^{\rightarrow}a_k\right)=1$.

It turns out that this product has many interesting properties, e.g., if $\varphi_1$ and $\varphi_2$ are states, then their free product is a again a state. For more details, see \cite{barndorff-nielsen+thorbjoensen03} and the references given there.

\subsection{Examples: Boolean, Monotone, and Anti-monotone Independence}

Ben Ghorbal and Sch\"urmann\cite{benghorbal01,benghorbal+schuermann99} and Muraki\cite{muraki02} have also considered the category of non-unital algebraic probability $\eufrak{nuAlgProb}$ consisting of pairs $(\mathcal{A},\varphi)$ of a not necessarily unital algebra $\mathcal{A}$ and a linear functional $\varphi$. The morphisms in this category are algebra homomorphisms that leave the functional invariant. On this category we can define three more tensor products with inclusions corresponding to the boolean product $\diamond$, the monotone product $\triangleright$ and the anti-monotone product $\triangleleft$ of states. They can be defined by
\begin{eqnarray*}
\varphi_1\diamond\varphi_2(a_1a_2\cdots a_m) &=& \prod_{k=1}^m \varphi_{\epsilon_k}(a_k), \\
\varphi_1\triangleright\varphi_2(a_1a_2\cdots a_m) &=& \varphi_1\left(\prod_{k:\epsilon_k=1}^{\rightarrow} a_k\right)\prod_{k:\epsilon_k=2} \varphi_2(a_k), \\
\varphi_1\triangleleft\varphi_2(a_1a_2\cdots a_m) &=& \prod_{k:\epsilon_k=1} \varphi_1(a_k)\,\,\varphi_2\left(\prod_{k:\epsilon_k=2}^{\rightarrow} a_k\right), 
\end{eqnarray*}
for $\varphi_1:\mathcal{A}_1\to\mathbb{C}$ and  $\varphi_2:\mathcal{A}_2\to\mathbb{C}$ and a typical element $a_1a_2\cdots a_m\in\mathcal{A}_1\coprod\mathcal{A}_2$, $a_k\in\mathcal{A}_{\epsilon_k}$, $\epsilon_1\not=\epsilon_2\not=\cdots\not=\epsilon_m$, i.e.\ neighboring $a$'s don't belong to the same algebra. Note that for the algebras and the inclusions we here use the free product without units, the coproduct in the category of not necessarily unital algebras.

The monotone and anti-monotone product are not commutative, but related by
\[
\varphi_1\triangleright\varphi_2 = (\varphi_2\triangleleft\varphi_1)\circ\gamma_{\mathcal{A}_1,\mathcal{A}_2}, \quad \mbox{ for all linear functionals }\varphi_1:\mathcal{A}_1\to\mathbb{C}, \varphi_2:\mathcal{A}_2\to\mathbb{C},
\]
where $\gamma_{\mathcal{A}_1,\mathcal{A}_2}:\mathcal{A}_1\coprod\mathcal{A}_2\to\mathcal{B}_2\coprod\mathcal{B}_1$ is the commutativity constraint (for the commutativity constraint for the free product of unital algebras see Equation \eqref{eq-comm}). The boolean product is commutative, i.e.\ it satisfies
\[
\varphi_1\diamond\varphi_2 = (\varphi_2\diamond\varphi_1)\circ\gamma_{\mathcal{B}_1,\mathcal{B}_2},
\]
for all linear functionals $\varphi_1:\mathcal{A}_1\to\mathbb{C}$, $\varphi_2:\mathcal{A}_2\to\mathbb{C}$.

\begin{remark}\label{III-xca-boolean}
The boolean, the monotone and the anti-monotone product can also be defined for unital algebras, if they are in the range of the unitization functor.

Let $\varphi_1:\mathcal{A}_1\to\mathbb{C}$ and $\varphi_2:\mathcal{A}_2\to\mathbb{C}$ be two unital functionals on algebras $\mathcal{A}_1$, $\mathcal{A}_2$, which can be decomposed as $\mathcal{A}_1=\mathbb{C}\mathbf{1}\oplus\mathcal{A}_1^0$, $\mathcal{A}_2=\mathbb{C}\mathbf{1}\oplus\mathcal{A}_2^0$. Then we define the boolean, monotone, or anti-monotone product of $\varphi_1$ and $\varphi_2$ as the unital extension of the boolean, monotone, or anti-monotone product of their restrictions $\varphi_1|_{\mathcal{A}_1^0}$ and $\varphi_2|_{\mathcal{A}_2^0}$.

This leads to the following formulas.
\begin{eqnarray*}
\varphi_1\diamond\varphi_2(a_1a_2\cdots a_n) &=&\prod_{i=1}^n \varphi_{\epsilon_i}(a_i), \\
\varphi_1\triangleright\varphi_2(a_1a_2\cdots a_n) &=&\varphi_1\left(\prod_{i:\epsilon_i=1} a_i\right)\prod_{i:\epsilon_i=2}\varphi_2(a_i), \\
\varphi_1\triangleleft\varphi_2(a_1a_2\cdots a_n) &=&\prod_{i:\epsilon_i=1}\varphi_1(a_i)\varphi_2\left(\prod_{i:\epsilon_i=2} a_i\right),
\end{eqnarray*}
for $a_1a_2\cdots a_n\in\mathcal{A}_1\coprod\mathcal{A}_2$, $a_i\in\mathcal{A}_{\epsilon_i}^0$, $\epsilon_1\not=\epsilon_2\not=\cdots\not=\epsilon_n$. We use the convention that the empty product is equal to the unit element.
\end{remark}

These products can be defined in the same way for $*$-algebraic probability spaces, where the algebras are unital $*$-algebras having such a decomposition $\mathcal{A}=\mathbb{C}\mathbf{1}\oplus\mathcal{A}_0$ and the functionals are states.
To check that $\varphi_1\diamond\varphi_2,\varphi_1\triangleright\varphi_2,\varphi_1\triangleleft\varphi_2$
are again states, if $\varphi_1$ and $\varphi_2$ are states, one can verify
that the following constructions give their GNS representations. Let
$(\pi_1,H_1,\xi_1)$ and $(\pi_2,H_2,\xi_2)$ denote the GNS representations of
$(\mathcal{A}_1,\varphi_1)$ and $(\mathcal{A}_2,\varphi_2)$. The GNS
representations of $(\mathcal{A}_1\coprod\mathcal{A}_2,\varphi_1\diamond\varphi_2)$, $(\mathcal{A}_1\coprod\mathcal{A}_2,\varphi_1\triangleright\varphi_2)$, and $(\mathcal{A}_1\coprod\mathcal{A}_2,\varphi_1\triangleleft\varphi_2)$
can all be defined on the Hilbert space $H=H_1\otimes H_2$ with the state
vector $\xi=\xi_1\otimes \xi_2$. The representations are defined by
$\pi(\mathbf{1})={\rm id}$ and 
\[
\begin{array}{rccrccll}
\pi|_{\mathcal{A}^0_1} &=& \pi_1\otimes P_2, &\pi|_{\mathcal{A}^0_2} &=& P_1\otimes \pi_2, & \mbox{ for } & \varphi_1\diamond\varphi_2, \\
\pi|_{\mathcal{A}^0_1} &=& \pi_1\otimes P_2, &\pi|_{\mathcal{A}^0_2} &=& {\rm id}_{H_2}\otimes\pi_2, & \mbox{ for } & \varphi_1\triangleright\varphi_2, \\
\pi|_{\mathcal{A}^0_1} &=& \pi_1\otimes{\rm id}_{H_2}, &\pi|_{\mathcal{A}^0_2} &=& P_1\otimes\pi_2, & \mbox{ for } & \varphi_1\triangleleft\varphi_2, \end{array}
\]
where $P_1,P_2$ denote the orthogonal projections $P_1:H_1\to \mathbb{C}\xi_1$, $P_2:H_2\to \mathbb{C}\xi_2$. For the boolean case, $\xi=\xi_1\otimes\xi_2\in H_1\otimes H_2$ is not cyclic for $\pi$, only the subspace $\mathbb{C}\xi\oplus H_1^0\oplus H_2^0$ can be generated from $\xi$.

\section{Reduction of an independence}\label{II-reduction}

For a reduction of independences we need a little bit more than a cotensor functor.

\begin{definition}\label{def-reduction}
Let $(\mathcal{C},\square,i)$ and $(\mathcal{C}',\square',i')$ be two tensor categories with inclusions and assume that we are given functors $I:\mathcal{C}\to\mathcal{D}$ and $I':\mathcal{C}'\to\mathcal{D}$ to some category $\mathcal{D}$. A {\em reduction} $(F,J)$ of the tensor product $\square$ to the tensor product $\square'$ (w.r.t.\ $(\mathcal{D},I,I')$) is a cotensor functor $F:(\mathcal{C},\square)\to(\mathcal{C}',\square')$ and a natural transformation $J:I\to I'\circ F$, i.e.\ morphisms $J_A:I(A)\to I'\circ F(A)$ in $\mathcal{D}$ for all objects $A\in{\rm Ob}\,\mathcal{C}$ such that the diagram
\[
\xymatrix{
I(A) \ar[d]_{I(f)} \ar[r]^{J_A} & I'\circ F(A)\ar[d]^{I'\circ F(f)} \\
I(B) \ar[r]_{J_B} & I'\circ F(B)
}
\]
commutes for all morphisms $f:A\to B$ in $\mathcal{C}$.
\end{definition}

In the simplest case, $\mathcal{C}$ will be a subcategory of  $\mathcal{C}'$, $I$ will be the inclusion functor from  $\mathcal{C}$ into $\mathcal{C}'$, and $I'$ the identity functor on $\mathcal{C}'$. Then such a reduction provides us with a system of inclusions
\[
J_n(A_1,\ldots,A_n)=F_n(A_1,\ldots,A_n)\circ J_{A_1\square\cdots
  \square A_n}:A_1\square\cdots \square A_n\to
F(A_1)\square'\cdots\square'F(A_n)
\]
with $J_1(A)=J_A$ that satisfies,
e.g., $J_{n+m}(A_1,\ldots,A_{n+m}) = F_2\big(F(A_1)\square' \cdots
\square'F(A_n),$ $F(A_{n+1})\square' \cdots \square'F(A_{n+m})\big) \circ \big(J_n(A_1, \ldots ,A_n)\square J_m(A_{n+1}, \ldots ,A_{n+m})\big)$ for all $n,m\in\mathbb{N}$ and $A_1,\ldots,A_{n+m}\in{\rm Ob}\,\mathcal{C}$.

In our applications we will often encounter the case where $\mathcal{C}$ is not a subcategory of $\mathcal{C}'$, but we have, e.g., a forgetful functor $U$ from $\mathcal{C}$ to $\mathcal{C}'$ that ``forgets'' an additional structure that $\mathcal{C}$ has. An example for this situation is the reduction of Fermi independence to tensor independence in following subsection. Here we have to forget the $\mathbb{Z}_2$-grading of the objects of $\mathbb{Z}_2$-$\eufrak{AlgProb}$ to get objects of $\eufrak{AlgProb}$. In this situation a reduction of the tensor product with inclusions $\square$ to the tensor product with inclusions $\square'$ is a tensor function $F$ from $(\mathcal{C},\square)$ to $(\mathcal{C}',\square')$ and a natural transformation $J:U\to F$.

\begin{example}
The identity functor can be turned into a reduction from the category $(\eufrak{Alg},\coprod)$ of unital associative algebras with the free product to the category $(\eufrak{Alg},\otimes)$ of unital associative algebras with the tensor product(with the obvious inclusions). The morphism $F_0:\mathbb{K}\to\mathbb{K}$ is the identity map and $F_2=[i_1,i_2]:\mathcal{A}_1\coprod\mathcal{A}_2\to \mathcal{A}_1\otimes\mathcal{A}_2$ is the unique morphism that makes the diagram
\[
\xymatrix{
\mathcal{A}_1 \ar[r]\ar[rd] & \mathcal{A}_1\coprod\mathcal{A}_2 \ar@{.>}[d]|-{F_2} & \mathcal{A}_2\ar[l]\ar[dl] \\
 & \mathcal{A}_1 \otimes\mathcal{A}_2 &
}
\]
commuting.
\end{example}

\subsection{The symmetric Fock space as a tensor functor}

The category $\eufrak{Vec}$ with the direct product $\oplus$ is of course a tensor category with inclusions and with projections, since the direct sum of vector spaces is both a product and a coproduct.

Not surprisingly, the usual tensor product of vector spaces is also a tensor product in the sense of category theory, but there are no canonical inclusions or projections. We can fix this by passing to the category $\eufrak{Vek}_*$ of pointed vector spaces, whose objects are pairs $(V,v)$ consisting of a vector space $V$ and a non-zero vector $v\in V$. The morphisms $h:(V_1,v_1)\to(V_2,v_2)$ in this category are the linear maps $h:V_1\to V_2$ with $h(v_1)=v_2$. In this category (equipped with the obvious tensor product $(V_1,v_1)\otimes(V_2,v_2)=(V_1\otimes V_2,v_1\otimes v_2)$ inclusions can be defined by $I_1:V_1\ni u \mapsto u\otimes v_2\in V_1\otimes V_2$ and $I_2:V_1\ni u \mapsto v_1\otimes u\in V_1\otimes V_2$.

In $(\eufrak{Vek}_*,\otimes,I)$ all pairs of morphisms are independent, even though the tensor product is not a coproduct.

\begin{proposition}
Take $\mathcal{D}=\eufrak{Vek}$, $I={\rm id}_\eufrak{Vek}$, and $I':\eufrak{Vek}_*\to\eufrak{Vek}$ the functor that forgets the fixed vector.

The symmetric Fock space $\Gamma$ is a reduction from $(\eufrak{Vek},\oplus,i)$ to $(\eufrak{Vek}_*,\otimes,I)$ (w.r.t.\ $(\eufrak{Vek},{\rm id}_\eufrak{Vek},I')$).
\end{proposition}

We will not prove this proposition, we will only define all the natural transformations.

On the objects, $\Gamma$ maps a vector space $V$ to the pair $\big(\Gamma(V),\Omega\big)$ consisting of the algebraic symmetric Fock space
\[
\Gamma(V)=\bigoplus_{n\in\mathbb{N}} V^{\otimes n}
\]
and the vacuum vector $\Omega$. The trivial vector space $\{0\}$ gets mapped to the field $\Gamma(\{0\})=\mathbb{K}$ with the unit $1$ as fixed vector. Linear maps $h:V_1\to V_2$ get mapped to their second quantization $\Gamma(h):\Gamma(V_1)\to \Gamma(V_2)$. $F_0:\Gamma(\{0\})=(\mathcal{K},\mathbf{1})\to (\mathcal{K},\mathbf{1})$ is just the identity and $F_2$ is the natural isomorphism from $\Gamma(V_1\oplus V_2)$ to $\Gamma(V_1)\otimes \Gamma(V_2)$ which acts on exponential vectors as
\[
F_2:\mathcal{E}(u_1+u_2)\mapsto \mathcal{E}(u_1)\otimes\mathcal{E}(u_2)
\]
for $u_1\in V_1$, $u_2\in V_2$.

The natural transformation $J:{\rm id}_\eufrak{Vec}\to \Gamma$ finally is the embedding of $V$ into $\Gamma(V)$ as one-particle space.

\subsection{Example: Bosonization of Fermi Independence}\label{reduction-fermi-bose}

We will now define the bosonization of Fermi independence as a reduction from $(\mathbb{Z}_2$-$\eufrak{AlgProb},\otimes_{\mathbb{Z}_2},i)$ to $(\eufrak{AlgProb},\otimes,i)$. We will need the group algebra $\mathbb{C}\mathbb{Z}_2$ of $\mathbb{Z}_2$ and the linear functional $\varepsilon:\mathbb{C}\mathbb{Z}_2\to\mathbb{C}$ that arises as the linear extension of the trivial representation of $\mathbb{Z}_2$, i.e.
\[
\varepsilon(\mathbf{1})=\varepsilon(g)=1,
\]
if we denote the even element of $\mathbb{Z}_2$ by $\mathbf{1}$ and the odd element by $g$.

The underlying functor $F:\mathbb{Z}_2$-$\eufrak{AlgProb}\to\eufrak{AlgProb}$ is given by
\[
F:
\begin{array}{lcl}
{\rm Ob}\,\mathbb{Z}_2\mbox{-}\eufrak{AlgProb}\ni(\mathcal{A},\varphi) & \mapsto & (\mathcal{A}\otimes_{\mathbb{Z}_2}\mathbb{C}\mathbb{Z}_2, \varphi\otimes\varepsilon)\in{\rm Ob}\,\eufrak{AlgProb}, \\[1mm]
{\rm Mor}\,\mathbb{Z}_2\mbox{-}\eufrak{AlgProb}\ni f & \mapsto &  f\otimes {\rm id}_{\mathbb{C}\mathbb{Z}_2} \in{\rm Mor}\,\eufrak{AlgProb}.
\end{array}
\]

The unit element in both tensor categories is the one-dimensional unital algebra $\mathbb{C}\mathbf{1}$ with the unique unital functional on it. Therefore $F_0$ has to be a morphism from $F(\mathbb{C}\mathbf{1})\cong \mathbf{C}\mathbb{Z}_2$ to $\mathbb{C}\mathbf{1}$. It is defined by $F_0(\mathbf{1})=F_0(g)=\mathbf{1}$.

The morphism $F_2(\mathcal{A}_1,\mathcal{A}_2)$ has to go from $F(\mathcal{A}\otimes_{\mathbb{Z}_2}\mathcal{B})=(\mathcal{A}\otimes_{\mathbb{Z}_2}\mathcal{B})\otimes \mathbb{C}\mathbb{Z}_2$ to $F(\mathcal{A})\otimes F(\mathcal{B})=(\mathcal{A}\otimes_{\mathbb{Z}_2}\mathbf{C}\mathbb{Z}_2)\otimes(\mathcal{B}\otimes_{\mathbb{Z}_2}\mathbf{C}\mathbb{Z}_2)$. It is defined by
\[
a\otimes b \otimes \mathbf{1} \mapsto
\left\{
\begin{array}{lcl}
(a\otimes \mathbf{1}) \otimes (b\otimes \mathbf{1}) & \mbox{ if } & b\mbox{ is even}, \\
(a\otimes g) \otimes (b\otimes \mathbf{1}) & \mbox{ if } & b\mbox{ is odd},
\end{array}
\right.
\]
and
\[
a\otimes b \otimes g \mapsto
\left\{
\begin{array}{lcl}
(a\otimes g) \otimes (b\otimes g) & \mbox{ if } & b\mbox{ is even}, \\
(a\otimes \mathbf{1}) \otimes (b\otimes g) & \mbox{ if } & b\mbox{ is odd},
\end{array}
\right.
\]
for $a\in\mathcal{A}$ and homogeneous $b\in\mathcal{B}$.

Finally, the inclusion $J_\mathcal{A}:\mathcal{A}\to \mathcal{A}\otimes_{\mathbb{Z}_2}\mathbf{C}\mathbb{Z}_2$ is defined by
\[
J_\mathcal{A}(a) = a\otimes\mathbf{1}
\]
for all $a\in\mathcal{A}$.

In this way we get inclusions $J_n=J_n(\mathcal{A}_1,\ldots,\mathcal{A}_n)=F_n(\mathcal{A}_1,\ldots,\mathcal{A}_n)\circ J_{\mathcal{A}_1\otimes_{\mathbb{Z}_2}\ldots\otimes_{\mathbb{Z}_2}\mathcal{A}_n}$ of the graded tensor product $\mathcal{A}_1\otimes_{\mathbb{Z}_2}\cdots\otimes_{\mathbb{Z}_2}\mathcal{A}_n$ into the usual tensor product $(\mathcal{A}_1\otimes_{\mathbb{Z}_2}\mathbb{C}\mathbb{Z}_2)\otimes\cdots\otimes (\mathcal{A}_n\otimes_{\mathbb{Z}_2}\mathbb{C}\mathbb{Z}_2)$ which respect the states and allow to reduce all calculations involving the graded tensor product to calculations involving the usual tensor product on the bigger algebras $F(\mathcal{A}_1)=\mathcal{A}_1\otimes_{\mathbb{Z}_2}\mathbb{C}\mathbb{Z}_2,\ldots,F(\mathcal{A}_n)=\mathcal{A}_n\otimes_{\mathbb{Z}_2}\mathbb{C}\mathbb{Z}_2$. These inclusions are determined by
\begin{eqnarray*}
J_n(\underbrace{\mathbf{1}\otimes\cdots\otimes\mathbf{1}}\otimes a \otimes\underbrace{\mathbf{1}\otimes\cdots\otimes\mathbf{1}}) &=& \underbrace{\tilde{g}\otimes\cdots\otimes\tilde{g}}\otimes \tilde{a}\otimes\underbrace{\tilde{\mathbf{1}}\otimes\cdots\otimes\tilde{\mathbf{1}}}, \\
\hspace*{11mm} k-1 \mbox{ times} \hspace{8mm} n-k \mbox{ times} && k-1 \mbox{ times} \hspace{8mm} n-k \mbox{ times}
\end{eqnarray*}
for $a\in \mathcal{A}_k$, $1\le k\le n$, where we used the abbreviations
\[
\tilde{g}=\mathbf{1}\otimes g, \qquad \tilde{a}=a\otimes\mathbf{1}, \qquad \tilde{\mathbf{1}}=\mathbf{1}\otimes\mathbf{1}.
\]

\subsection{The Reduction of Boolean, Monotone, and Anti-Monotone Independence to Tensor Independence}\label{section-reduction-boolean-etc}

We will now present the unification of tensor, monotone, anti-monotone, and boolean independence of Franz\cite{franz03b} in our category theoretical framework. It resembles closely the bosonization of Fermi independence in Subsection \ref{reduction-fermi-bose}, but the group $\mathbb{Z}_2$ has to be replaced by the semigroup $M=\{\mathbf{1},p\}$ with two elements, $\mathbf{1}\cdot \mathbf{1}=\mathbf{1}$, $\mathbf{1}\cdot p=p\cdot \mathbf{1}=p\cdot p=p$. We will need the linear functional $\varepsilon:\mathbb{C}M\to\mathbb{C}$ with $\varepsilon(\mathbf{1})=\varepsilon(p)=1$.

The underlying functor and the inclusions are the same for the reduction of the boolean, the monotone and the anti-monotone product. They map the algebra $\mathcal{A}$ of $(\mathcal{A},\varphi)$ to the free product $F(\mathcal{A})=\tilde{\mathcal{A}}\coprod\mathbb{C}M$ of the unitization $\tilde{\mathcal{A}}$ of $\mathcal{A}$ and the group algebra $\mathbb{C}M$ of $M$. For the unital functional $F(\varphi)$ we take the boolean product $\tilde{\varphi}\diamond\varepsilon$ of the unital extension $\tilde{\varphi}$ of $\varphi$ with $\varepsilon$. The elements of $F(\mathcal{A})$ can be written as linear combinations of terms of the form
\[
p^\alpha a_1p\cdots p a_mp^\omega
\]
with $m\in\mathbb{N}$, $\alpha,\omega\in\{0,1\}$, $a_1,\ldots.a_m\in\mathcal{A}$, and $F(\varphi)$ acts on them as
\[
F(\varphi)(p^\alpha a_1p\cdots p a_mp^\omega) = \prod_{k=1}^m \varphi(a_k).
\]
The inclusion is simply
\[
J_{\mathcal{A}}: \mathcal{A}\ni a\mapsto a\in F(\mathcal{A}).
\]
The morphism $F_0:F(\mathbb{C}\mathbf{1})=\mathbb{C}M\to\mathbb{C}\mathbf{1}$ is given by the trivial representation of $M$, $F_0(\mathbf{1})=F_0(p)=\mathbf{1}$. 

The only part of the reduction that is different for the three cases are the morphisms
\[
F_2(\mathcal{A}_1,\mathcal{A}_2):F\left(\mathcal{A}_1\coprod\mathcal{A}_2\right)\to F(\mathcal{A}_1)\otimes F(\mathcal{A}_2)=(\tilde{\mathcal{A}}\coprod\mathbb{C}M)\otimes(\tilde{\mathcal{A}}\coprod\mathbb{C}M).
\]
We set
\[
F^{\rm B}_2(\mathcal{A}_1,\mathcal{A}_2)(a) =
\left\{
\begin{array}{lcl}
a\otimes p & \mbox{ if } a\in\mathcal{A}_1, \\
p\otimes a & \mbox{ if } a\in\mathcal{A}_2,
\end{array}
\right.
\]
for the boolean case,
\[
F^{\rm M}_2(\mathcal{A}_1,\mathcal{A}_2)(a) =
\left\{
\begin{array}{lcl}
a\otimes p & \mbox{ if } a\in\mathcal{A}_1, \\
\mathbf{1}\otimes a & \mbox{ if } a\in\mathcal{A}_2,
\end{array}
\right.
\]
for the monotone case, and
\[
F^{\rm AM}_2(\mathcal{A}_1,\mathcal{A}_2)(a) =
\left\{
\begin{array}{lcl}
a\otimes \mathbf{1} & \mbox{ if } a\in\mathcal{A}_1, \\
p\otimes a & \mbox{ if } a\in\mathcal{A}_2,
\end{array}
\right.
\]
for the anti-monotone case. Furthermore, we have $F^\bullet_2(\mathcal{A}_1,\mathcal{A}_2)(p)=p\otimes p$ in all three cases.

For the higher order inclusions $J^\bullet_n= F^\bullet_n(\mathcal{A}_1,\ldots,\mathcal{A}_n)\circ J_{\mathcal{A}_1\coprod\cdots\coprod\mathcal{A}_n}$, $\bullet\in\{{\rm B},{\rm M},{\rm AM}\}$, one gets
\begin{eqnarray*}
J_n^{\rm B}(a) &=& p^{\otimes (k-1)}\otimes a \otimes p^{\otimes (n-k)}, \\
J_n^{\rm M}(a) &=& \mathbf{1}^{\otimes (k-1)}\otimes a \otimes p^{\otimes (n-k)}, \\
J_n^{\rm AM}(a) &=& p^{\otimes (k-1)}\otimes a \otimes \mathbf{1}^{\otimes (n-k)},
\end{eqnarray*}
if $a\in\mathcal{A}_k$.

One can verify that this indeed defines reductions $(F^{\rm B},J)$,  $(F^{\rm M},J)$, and $(F^{\rm AM},J)$ from the categories $(\eufrak{nuAlgProb},\diamond,i)$, $(\eufrak{nuAlgProb},\triangleright,i)$, and $(\eufrak{nuAlgProb},\triangleleft,i)$ to the category $(\eufrak{AlgProb},\otimes,i)$ of algebraic probability spaces with the usual tensor product. The functor $U:\eufrak{nuAlgProb}\to\eufrak{AlgProb}$ is the unitization of the algebra and the unital extension of the functional and the morphisms.

This reduces all calculations involving the boolean, monotone or anti-monotone product to the tensor product. These constructions can also be applied to reduce the quantum stochastic calculus on the boolean, monotone, and anti-monotone Fock space to the boson Fock space. Furthermore, they allow to reduce the theories of boolean, monotone, and anti-monotone L\'evy processes to Sch\"urmann's\cite{schuermann93} theory of L\'evy processes on involutive bialgebras, see Franz\cite{franz03b} or Section \ref{IV-reduction} in Chapter \ref{chapter-IV}.

A similar reduction exists for the category of unital algebras $\mathcal{A}$ having a decomposition $\mathcal{A}=\mathbb{C}\mathbf{1}\oplus\mathcal{A}_0$ and the boolean, monotone, or anti-monotone product defined for these algebras in Remark \ref{III-xca-boolean}

\section{Classification of the universal independences}\label{III-class-indep}

In the previous Section we have seen how a notion of independence can be defined in the language of category theory and we have also encountered several examples. 

We are mainly interested in different categories of algebraic
probability spaces. There objects are pairs consisting of an algebra
$\mathcal{A}$ and a linear functional $\varphi$ on
$\mathcal{A}$. Typically, the algebra has some additional structure,
e.g., an involution, a unit, a grading, or a topology (it can be, e.g., a von Neumann algebra or a $C^*$-algebra), and the functional behaves nicely with respect to this additional structure, i.e., it is positive, unital, respects the grading, continuous, or normal. The morphisms are algebra homomorphisms, which leave the linear functional invariant, i.e., $j:(\mathcal{A},\varphi)\to (\mathcal{B},\psi)$ satisfies
\[
\varphi=\psi\circ j
\]
and behave also nicely w.r.t.\ additional structure, i.e., they can be
required to be $*$-algebra homomorphisms, map the unit of
$\mathcal{A}$ to the unit of $\mathcal{B}$, respect the grading, etc. We have already seen one example in Subsection \ref{fermi-indep}.

The tensor product then has to specify a new algebra with a linear functional
and inclusions for every pair of of algebraic probability spaces. If the category of algebras obtained from our algebraic
probability space by forgetting the linear functional has a coproduct, then it
is sufficient to consider the case where the new algebra is the coproduct of
the two algebras. 

\begin{proposition}\label{prop-forget}
Let $(\mathcal{C},\square,i)$ be a tensor category with inclusions and $F:\mathcal{C}\to\mathcal{D}$ a functor from $\mathcal{C}$ into another category $\mathcal{D}$ which has a coproduct $\coprod$ and an initial object $E_\mathcal{D}$. Then $F$ is a tensor functor. The morphisms $F_2(A,B):F(A)\coprod F(B)\to F(A\square B)$ and $F_0:E_\mathcal{D}\to F(E)$ are those guaranteed by the universal property of the coproduct and the initial object, i.e.\ $F_0:E_\mathcal{D}\to F(E)$ is the unique morphism from $E_\mathcal{D}$ to $F(E)$ and $F_2(A,B)$ is the unique morphism that makes the diagram
\[
\xymatrix{
F(A) \ar[r]^(0.4){F(i_A)}\ar[dr]_{i_{F(A)}} & F(A\square B) & \ar[l]_(0.4){F(i_B)} \ar[dl]^{i_{F(B)}} F(B) \\
 & F(A)\coprod F(B) \ar@{.>}[u]|-{F_2(A,B)} &
}
\]
commuting.
\end{proposition}
\begin{proof}
Using the universal property of the coproduct and the definition of $F_2$, one shows that the triangles containing the $F(A)$ in the center of the diagram
\[
\xymatrix{
F(A)\coprod \big(F(B)\coprod F(C)\big) \ar[rr]^{\alpha_{F(A),F(B),F(C)}} \ar[d]_{{\rm id}_{F(A)}\coprod F_2(B,C)} && \big(F(A)\coprod F(B)\big)\coprod F(C)\ar[d]^{F_2(A,B)\coprod{\rm id}_{F(C)}} \\
F(A)\coprod F(B\square C)\ar[d]_{F_2(A,B\square C)} & F(A) \ar[ul]|-{i_{F(A)}} \ar[l]|-{i_{F(A)}} \ar[dl]|-{F(i_A)} \ar[ur]|-{i_{F(A)}} \ar[r] \ar[rd]|-{F(i_A)} & F(A\square B)\coprod F(C)\ar[d]^{F_2(A\square B,C)} \\
F\big(A\square(B\square C)\big)\ar[rr]_{F(\alpha_{A,B,C})} & & F\big((A\square B)\square C\big)
}
\]
commute (where the morphism from $F(A)$ to $ F(A\square B)\coprod F(C)$ is $F(i_A)\coprod{\rm id}_{F(C)}$), and therefore that the morphisms corresponding to all the different paths form $F(A)$ to $F\big((A\square B)\square C\big)$ coincide. Since we can get similar diagrams with $F(B)$ and $F(C)$, it follows from the universal property of the triple coproduct $F(A)\coprod \big(F(B)\coprod F(C)\big)$ that there exists only a unique morphism from $F(A)\coprod \big(F(B)\coprod F(C)\big)$ to $F\big((A\square B)\square C\big)$ and therefore that the whole diagram commutes.

The commutativity of the two diagrams involving the unit elements can be shown similarly.
\end{proof}

Let $\mathcal{C}$ now be a category of algebraic probability spaces and $F$
the functor that maps a pair $(\mathcal{A},\varphi)$ to the algebra $\mathcal{A}$, i.e., that ``forgets'' the linear functional $\varphi$. Suppose that $\mathcal{C}$ is equipped with a tensor product $\square$ with inclusions and that $F(\mathcal{C})$ has a coproduct $\coprod$. Let $(\mathcal{A},\varphi)$, $(\mathcal{B},\psi)$ be two algebraic probability spaces in $\mathcal{C}$, we will denote the pair $(\mathcal{A},\varphi)\square(\mathcal{B},\psi)$ also by $(\mathcal{A}\square\mathcal{B},\varphi\square\psi)$. By Proposition \ref{prop-forget} we have morphisms $F_2(\mathcal{A},\mathcal{B}):\mathcal{A}\coprod\mathcal{B} \to\mathcal{A}\square\mathcal{B}$ that define a natural transformation from the bifunctor $\coprod$ to the bifunctor $\square$. With these morphisms we can define a new tensor product $\widetilde{\square}$ with inclusions by
\[
(\mathcal{A},\varphi)\widetilde{\square}(\mathcal{B},\psi)= \left(\mathcal{A}\coprod \mathcal{B},(\varphi\square\psi)\circ F_2(\mathcal{A},\mathcal{B})\right).
\]
The inclusions are those defined by the coproduct.

\begin{proposition}
If two random variables $f_1:(\mathcal{A}_1,\varphi_1)\to(\mathcal{B},\psi)$ and
$f_1:(\mathcal{A}_1,\varphi_1)\to(\mathcal{B},\psi)$ are independent with respect to $\square$, then they are also independent with respect to $\widetilde{\square}$.
\end{proposition}
\begin{proof}
If $f_1$ and $f_2$ are independent with respect to $\square$, then there exists a random variable $h:(\mathcal{A}_1\square \mathcal{A}_2,\varphi_1\square\varphi_2)\to (\mathcal{B},\psi)$ that makes diagram (\ref{diag-indep2}) in Definition \ref{def-independence} commuting. Then $h\circ F_2(\mathcal{A}_1,\mathcal{A}_2):(\mathcal{A}_1\coprod \mathcal{A}_2,\varphi_1\tilde{\square}\varphi_2)\to (\mathcal{B},\psi)$ makes the corresponding diagram for $\tilde{\square}$ commuting.
\end{proof}
The converse is not true. Consider the category of algebraic probability spaces with the tensor product, see Subsection \ref{tensor indep}, and take $B=\mathcal{A}_1\coprod\mathcal{A}_2$ and $\psi=(\varphi_1\otimes\varphi_2)\circ F_2(\mathcal{A}_1,\mathcal{A}_2)$. The canonical inclusions $i_{\mathcal{A}_1}:(\mathcal{A}_1,\varphi_1)\to (\mathcal{B},\psi)$ and $i_{\mathcal{A}_2}:(\mathcal{A}_2,\varphi_2)\to (\mathcal{B},\psi)$ are independent w.r.t.\ $\tilde{\otimes}$, but not with respect to the tensor product itself, because their images do not commute in $\mathcal{B}=\mathcal{A}_1\coprod\mathcal{A}_2$.

We will call a tensor product with inclusions in a category of quantum probability spaces {\em universal}, if it is equal to the coproduct of the corresponding category of algebras on the algebras. The preceding discussion shows that every tensor product on the category of algebraic quantum probability spaces $\eufrak{AlgProb}$ has a universal version. E.g., for the tensor independence defined in the category of algebraic probability spaces in Subsection \ref{tensor indep}, the universal version is defined by
\[
\varphi_1\tilde{\otimes}\varphi_2(a_1a_2\cdots a_m)=\varphi_1\left(\prod^{\rightarrow}_{i:\epsilon_i=1}a_i\right)\varphi_2\left(\prod^{\rightarrow}_{i:\epsilon_i=2}a_i\right)
\]
for two unital functionals $\varphi_1:\mathcal{A}_1\to\mathbb{C}$ and  $\varphi_2:\mathcal{A}_2\to\mathbb{C}$ and  a typical element $a_1a_2\cdots a_m\in\mathcal{A}_1\coprod\mathcal{A}_2$, with $a_k\in\mathcal{A}_{\epsilon_k}$, $\epsilon_1\not=\epsilon_2\not=\cdots\not=\epsilon_m$, i.e.\ neighboring $a$'s don't belong to the same algebra.

We will now reformulate the classification by Muraki\cite{muraki02} and by Ben Ghorbal and Sch\"urmann\cite{benghorbal01,benghorbal+schuermann99} in terms of universal tensor products with inclusions for the category of algebraic probability spaces $\eufrak{AlgProb}$.

In order to define a universal tensor product with inclusions on $\eufrak{AlgProb}$ one needs a map that associates to a pair of unital functionals $(\varphi_1,\varphi_2)$ on two algebras $\mathcal{A}_1$ and $\mathcal{A}_2$ a unital functional $\varphi_1\cdot\varphi_2$ on the free product $\mathcal{A}_1\coprod\mathcal{A}_2$ (with identification of the units) of $\mathcal{A}_1$ and $\mathcal{A}_2$ in such a way that the bifunctor
\[
\square:(\mathcal{A}_1,\varphi_1)\times(\mathcal{A}_2,\varphi_1)\mapsto (\mathcal{A}_1\coprod\mathcal{A}_2,\varphi_1\cdot\varphi_2)
\]
satisfies all the necessary axioms. Since $\square$ is equal to the coproduct $\coprod$ on the algebras, we don't have a choice for the isomorphisms $\alpha,\lambda,\rho$ implementing the associativity and the left and right unit property. We have to take the ones following from the universal property of the coproduct. The inclusions and the action of $\square$ on the morphisms also have to be the ones given by the coproduct.

The associativity gives us the condition
\begin{equation}\label{cond-associativity}
\big((\varphi_1\cdot\varphi_2)\cdot\varphi_3\big)\circ \alpha_{\mathcal{A}_1,\mathcal{A}_2,\mathcal{A}_3} = \varphi_1\cdot(\varphi_2\cdot\varphi_3),
\end{equation}
for all $(\mathcal{A}_1,\varphi_1),(\mathcal{A}_2,\varphi_2),(\mathcal{A}_3,\varphi_3)$ in $\eufrak{AlgProb}$. Denote the unique unital functional on $\mathbb{C}\mathbf{1}$ by $\delta$, then the unit properties are equivalent to
\[
(\varphi\cdot \delta)\circ\rho_\mathcal{A} = \varphi\quad\mbox{ and }\quad (\delta\cdot\varphi)\circ\lambda_\mathcal{A} = \varphi,
\]
for all $(\mathcal{A},\varphi)$ in $\eufrak{AlgProb}$.
The inclusions are random variables, if and only if
\begin{equation}\label{cond-inclusion}
(\varphi_1\cdot\varphi_2)\circ i_{\mathcal{A}_1}=\varphi_1\quad\mbox{ and }\quad(\varphi_1\cdot\varphi_2)\circ i_{\mathcal{A}_2}=\varphi_2
\end{equation}
for all $(\mathcal{A}_1,\varphi_1),(\mathcal{A}_2,\varphi_2)$ in $\eufrak{AlgProb}$. Finally, from the functoriality of $\square$ we get the condition
\begin{equation}\label{cond-functoriality}
(\varphi_1\cdot\varphi_2)\circ (j_1\coprod j_2) = (\varphi_1\circ j_1)\cdot(\varphi_2\circ j_2)
\end{equation}
for all pairs of morphisms $j_1:(\mathcal{B}_1,\psi_1)\to(\mathcal{A}_1,\varphi_1)$, $j_2:(\mathcal{B}_2,\psi_2)\to(\mathcal{A}_2,\varphi_2)$ in $\eufrak{AlgProb}$.

Our Conditions (\ref{cond-associativity}), (\ref{cond-inclusion}), and (\ref{cond-functoriality}) are exactly the axioms (P2), (P3), and (P4) in Ben Ghorbal and Sch\"urmann\cite{benghorbal+schuermann99}, or the axioms (U2), the first part of (U4), and (U3) in Muraki\cite{muraki02}.
\begin{theorem}
(Muraki\cite{muraki02}, Ben Ghorbal and Sch\"urmann\cite{benghorbal01,benghorbal+schuermann99}) There exist exactly two universal tensor products with inclusions on the category of algebraic probability spaces $\eufrak{AlgProb}$, namely the universal version $\tilde{\otimes}$ of the tensor product defined in Section \ref{tensor indep} and the one associated to the free product $*$ of states.
\end{theorem}

For the classification in the non-unital case, Muraki imposes the additional condition
\begin{equation}\label{cond-factorisation}
(\varphi_1\cdot\varphi_2)(a_1a_2) = \varphi_{\epsilon_1}(a_1)\varphi_{\epsilon_2}(a_2)
\end{equation}
for all $(\epsilon_1,\epsilon_2)\in\big\{(1,2),(2,1)\big\}$, $a_1\in\mathcal{A}_{\epsilon_1}$, $a_2\in\mathcal{A}_{\epsilon_2}$.

\begin{theorem}
(Muraki\cite{muraki02})
There exist exactly five universal tensor products with inclusions satisfying (\ref{cond-factorisation}) on the category of non-unital algebraic probability spaces $\eufrak{nuAlgProb}$, namely the universal version $\tilde{\otimes}$ of the tensor product defined in Section \ref{tensor indep} and the ones associated to the free product $*$, the boolean product $\diamond$, the monotone product $\triangleright$ and the anti-monotone product $\triangleleft$.
\end{theorem}
The monotone and the anti-monotone are not symmetric, i.e.\ $(\mathcal{A}_1\coprod\mathcal{A}_2,\varphi_1\triangleright\varphi_2)$ and $(\mathcal{A}_2\coprod\mathcal{A}_2,\varphi_2\triangleright\varphi_1)$ are not isomorphic in general. Actually, the anti-monotone product is simply the mirror image of the monotone product,
\[
(\mathcal{A}_1\coprod\mathcal{A}_2,\varphi_1\triangleright\varphi_2)\cong(\mathcal{A}_2\coprod\mathcal{A}_1,\varphi_2\triangleleft\varphi_1)
\]
for all $(\mathcal{A}_1,\varphi_1),(\mathcal{A}_2,\varphi_2)$ in the category of non-unital algebraic probability spaces. The other three products are symmetric.

In the symmetric setting of Ben Ghorbal and Sch\"urmann, Condition (\ref{cond-factorisation}) is not essential. If one drops it and adds symmetry, one finds in addition the degenerate product
\[
(\varphi_1\bullet_0 \varphi_2) (a_1a_2\cdots a_m) =
\left\{
\begin{array}{lcl}
\varphi_{\epsilon_1}(a_1) & \mbox{ if } & m=1, \\
0 & \mbox{ if } & m>1.
\end{array}
\right.
\]
and families
\[
\varphi_1\bullet_q \varphi_2 = q\big((q^{-1}\varphi_1)\cdot(q^{-1}\varphi_2)\big),
\]
parametrized by a complex number $q\in\mathbb{C}\backslash\{0\}$, for each of the three symmetric products, $\bullet\in\{\tilde{\otimes},*,\diamond\}$.

If one adds the condition that products of states are again states, then one can also show that the constant has to be equal to one.

\begin{remark}
Consider the category of non-unital $*$-algebraic probability spaces, whose objects are pairs $(\mathcal{A},\varphi)$ consisting of a not necessarily unital $*$-algebra $\mathcal{A}$ and a state $\varphi:\mathcal{A}\to\mathbb{C}$. Here a state is a linear functional $\varphi:\mathcal{A}\to\mathbb{C}$ whose unital extension $\tilde{\varphi}:\tilde{\mathcal{A}}\cong\mathbb{C}\mathbf{1}\oplus\mathcal{A}\to\mathbb{C}$, $\lambda\mathbf{1}+a\mapsto \tilde{\varphi}(\lambda\mathbf{1}+a)=\lambda+\varphi(a)$, to the unitization of $\mathcal{A}$ is a state.

Assume we have a product $\cdot:S(\mathcal{A}_1)\times S(\mathcal{A}_2)\to S(\mathcal{A}_1\coprod\mathcal{A}_2)$ of linear functionals on non-unital algebras $\mathcal{A}_1,\mathcal{A}_2$ that satisfies
\begin{eqnarray*}
(\varphi_1\cdot\varphi_2)(a_1a_2) &=& c_1 \varphi_1(a_1)\varphi_2(a_2), \\
(\varphi_1\cdot\varphi_2)(a_2a_1) &=& c_2 \varphi_1(a_1)\varphi_2(a_2),
\end{eqnarray*}
for all linear functionals $\varphi_1:\mathcal{A}_1\to\mathbb{C}$, $\varphi_2:\mathcal{A}_2\to\mathbb{C}$, and elements $a_1\in\mathcal{A}_1$, $a_2\in\mathcal{A}_2$ with ``universal'' constants $c_1,c_2\in\mathbb{C}$, i.e. constants that do not depend on the algebras, the functionals, or the algebra elements. That for every universal independence such constants have to exist is part of the proof of the classifications in \cite{benghorbal01,benghorbal+schuermann99,muraki02}.

Take now $\mathcal{A}_1=\mathbb{C}[x_1]$ and $\mathcal{A}_2=\mathbb{C}[x_2]$, with the states $\varphi_1,\varphi_2$ determined by
\[
\varphi_i(x_i^n)=1
\]
for $n\in\mathbb{N}$, $i=1,2$. If $\varphi_1\cdot\varphi_2$ is a state, then
\[
\varphi_1\cdot\varphi_2\left((\lambda\mathbf{1}+\mu x_1 + \nu x_2)^*(\lambda\mathbf{1}+\mu x_1 + \nu x_2)\right)\ge 0
\]
for all $\lambda,\mu,\nu\in\mathbb{C}$. This is equivalent to the matrix
\[
A=\left(\begin{array}{ccc}
1 & 1 & 1 \\
1 & 1 & c_1 \\
1 & c_2 & 1
\end{array}
\right)
\]
being positive semi-definite. This is only possible, if $c_1=\overline{c_2}$ and $\det A=-(1-c_1)(1-c_2)\ge 0$, i.e.\ if $c_1=c_2=1$.
\end{remark}

The proof of the classification of universal independences can be split into three steps.

Using the ``universality'' or functoriality of the product, one can show that there exist some ``universal constants'' - not depending on the algebras - and a formula for evaluating
\[
(\varphi_1\otimes\varphi_2)(a_1a_2\cdots a_m)
\]
for $a_1a_2\cdots a_m\in\mathcal{A}_1\coprod\mathcal{A}_2$, with $a_k\in\mathcal{A}_{\epsilon_k}$, $\epsilon_1\not=\epsilon_2\not=\cdots\not=\epsilon_m$, as a linear combination of products $\varphi_1(M_1)$, $\varphi_2(M_2)$, where $M_1$, $M_2$ are ``sub-monomials'' of $a_1a_2\cdots a_m$. Then in a second step it is shown by associativity that only products with {\em ordered} monomials $M_1$, $M_2$ contribute. This is the content of \cite[Theorem 5]{benghorbal+schuermann02} in the commutative case and of \cite[Theorem 2.1]{muraki02} in the general case.

The third step, which was actually completed first in both cases, see \cite{speicher96} and \cite{muraki02a}, is to find the conditions that the universal constants have to satisfy, if the resulting product is associative. It turns out that the universal coefficients for $m>5$ are already uniquely determined by the coefficients for $1\le m\le 5$. Detailed analysis of the non-linear equations obtained for the coefficients of order up to five then leads to the classifications stated above.

%% file: quantum-levy.tex
\chapter{L\'evy Processes on Dual Groups}\label{chapter-IV}

We now want to study quantum stochastic processes whose increments are free or independent in the sense of boolean, monotone, or anti-monotone independence. The approach based on bialgebras that we followed in the first chapter is based on the tensor product and fails in the other cases because the corresponding products are not defined on the tensor product, but on the free product of the algebra. The algebraic structure which has to replace bialgebras was first introduced by Voiculescu \cite{voiculescu87,voiculescu90}, who named them dual groups. In this chapter we will introduce these algebras and develop the theory of their L\'evy processes. It turns out that L\'evy processes on dual groups with boolean, monotonically, or anti-monotonically independent increments can be reduced to L\'evy processes on involutive bialgebra. We do not know if this is also possible for L\'evy processes on dual groups with free increments.

In the literature additive free L\'evy processes have been studied most intensively, see, e.g., \cite{glockner+schuermann+speicher92,biane98,anshelevich01b,anshelevich01c,barndorff-nielsen+thorbjornsen01a,barndorff-nielsen+thorbjornsen01b,barndorff-nielsen+thorbjornsen01c}.

\section{Preliminaries on dual groups}

Denote by $\eufrak{ComAlg}$ the category of commutative unital algebras and let $\mathcal{B}\in {\rm Ob}\,\eufrak{ComAlg}$ be a commutative bialgebra. Then the mapping
\[
{\rm Ob}\,\eufrak{ComAlg}\ni \mathcal{A}\mapsto {\rm Mor}_\eufrak{ComAlg}(\mathcal{B},\mathcal{A})
\]
can be understood as a functor from $\eufrak{ComAlg}$ to the category of unital semigroups. The multiplication in ${\rm Mor}_\eufrak{Alg}(\mathcal{B},\mathcal{A})$ is given by the convolution, i.e.
\[
f\star g = m_{\mathcal{A}}\circ (f\otimes g)\circ \Delta_\mathcal{B}
\]
and the unit element is $\varepsilon_\mathcal{B}\mathbf{1}_\mathcal{A}$. A unit-preserving algebra homomorphism $h:\mathcal{A}_1\to\mathcal{A}_2$ gets mapped to the unit-preserving semigroup homomorphism ${\rm Mor}_\eufrak{ComAlg}(\mathcal{B},\mathcal{A}_1)\ni f\to h\circ f\in{\rm Mor}_\eufrak{ComAlg}(\mathcal{B},\mathcal{A}_2)$, since
\[
h\circ(f\star g) = (h\circ f)\star (h\circ g)
\]
for all $\mathcal{A}_1,\mathcal{A}_2\in{\rm Ob}\,\eufrak{ComAlg}$, $h\in{\rm Mor}_\eufrak{ComAlg}(\mathcal{A}_1,\mathcal{A}_2)$, $f,g\in{\rm Mor}_\eufrak{ComAlg}(\mathcal{B},\mathcal{A}_1)$.

If $\mathcal{B}$ is even a Hopf algebra with antipode $S$, then ${\rm Mor}_\eufrak{ComAlg}(\mathcal{B},\mathcal{A})$ is a group with respect to the convolution product. The inverse of a homomorphism $f:\mathcal{B}\to\mathcal{A}$ with respect to the convolution product is given by $f\circ S$.

The calculation
\begin{eqnarray*}
(f\star g)(ab) &=& m_{\mathcal{A}}\circ (f\otimes g)\circ \Delta_\mathcal{B}(ab) \\
&=& f(a_{(1)}b_{(1)})g(a_{(2)}b_{(2)})= f(a_{(1)})f(b_{(1)})g(a_{(2)})g(b_{(2)}) \\
&=& f(a_{(1)})g(a_{(2)})g(b_{(1)})g(b_{(2)}) = (f\star g)(a)(f\star g)(b)
\end{eqnarray*}
shows that the convolution product $f\star g$ of two homomorphisms $f,g:\mathcal{B}\to\mathcal{A}$ is again a homomorphism. It also gives an indication why non-commutative bialgebras or Hopf algebras do not give rise to a similar functor on the category of non-commutative algebras, since we had to commute $f(b_{(1)})$ with $g(a_{(2)})$.

Zhang \cite{zhang91}, Berman and Hausknecht \cite{bergman+hausknecht96} showed that if one replaces the tensor product in the definition of bialgebras and Hopf algebras by the free product, then one arrives at a class of algebras that do give rise to a functor from the category of non-commutative algebras to the category of semigroups or groups.

A dual group \cite{voiculescu87,voiculescu90} (called  H-algebra or cogroup in the category of unital associative $*$-algebras in \cite{zhang91} and \cite{bergman+hausknecht96}, resp.) is a unital  $*$-algebra $\mathcal{B}$ equipped with three unital $*$-algebra homomorphisms $\Delta:\mathcal{B}\to\mathcal{B}\coprod\mathcal{B}$, $S:\mathcal{B}\to\mathcal{B}$ and $\varepsilon:\mathcal{B}\to\mathbb{C}$ (also called comultiplication, antipode, and counit) such that
\begin{eqnarray}
\left(\Delta\coprod{\rm id}\right)\circ \Delta &=& \left({\rm id}\coprod\Delta\right)\circ \Delta, \label{coassociative}\\
\left(\varepsilon\coprod{\rm id}\right)\circ \Delta = &{\rm id}& = \left({\rm id}\coprod\varepsilon\right)\circ \Delta, \label{counit}\\
m_{\mathcal{B}}\circ\left(S\coprod{\rm id}\right)\circ \Delta = &{\rm id}& = m_{\mathcal{B}}\circ\left({\rm id}\coprod S\right)\circ \Delta,
\end{eqnarray}
where $m_{\mathcal{B}}:\mathcal{B}\coprod\mathcal{B}\to\mathcal{B}$, $m_{\mathcal{B}}(a_1\otimes a_2\otimes\cdots\otimes a_n)=a_1\cdot a_2\cdot\,\cdot\cdot\cdot\,\cdot a_n$, is the multiplication of $\mathcal{B}$. Besides the formal similarity, there are many relations between dual groups on the one side and Hopf algebras and bialgebras on the other side, cf.\ \cite{zhang91}. For example, let $\mathcal{B}$ be a dual group with comultiplication $\Delta$, and let $R:\mathcal{B}\coprod\mathcal{B}\to\mathcal{B}\otimes\mathcal{B}$ be the unique unital $*$-algebra homomorphism with
\[
R_{\mathcal{B},\mathcal{B}}\circ i_1(b)=b\otimes \mathbf{1}, \qquad
R_{\mathcal{B},\mathcal{B}}\circ i_2(b)=\mathbf{1}\otimes b,
\]
for all $b\in\mathcal{B}$. Here $i_1,i_2:\mathcal{B}\to\mathcal{B}\coprod\mathcal{B}$ denote the canonical inclusions of $\mathcal{B}$ into the first and the second factor of the free product $\mathcal{B}\coprod\mathcal{B}$. Then $\mathcal{B}$ is a bialgebra with the comultiplication $\overline{\Delta}=R_{\mathcal{B},\mathcal{B}}\circ\Delta$, see \cite[Theorem 4.2]{zhang91}, but in general it is not a Hopf algebra.

We will not really work with dual groups, but the following weaker notion. A dual semigroup is a unital $*$-algebra $\mathcal{B}$ equipped with two unital $*$-algebra homomorphisms $\Delta:\mathcal{B}\to\mathcal{B}\coprod\mathcal{B}$ and $\varepsilon:\mathcal{B}\to\mathbb{C}$ such that Equations \eqref{coassociative} and \eqref{counit} are satisfied. The antipode is not used in the proof of \cite[Theorem 4.2]{zhang91}, and therefore we also get an involutive bialgebra $(\mathcal{B},\overline{\Delta},\varepsilon)$ for every dual semigroup $(\mathcal{B},\Delta,\varepsilon)$.

Note that we can always write a dual semigroup $\mathcal{B}$ as a direct sum $\mathcal{B}=\mathbb{C}\mathbf{1}\oplus\mathcal{B}^0$, where $\mathcal{B}^0=\mbox{ker}\,\varepsilon$ is even a $*$-ideal. Therefore it is in the range of the unitization functor and the boolean, monotone, and anti-monotone product can be defined for unital linear functionals on $\mathcal{B}$, cf.\ Remark \ref{III-xca-boolean} in Chapter \ref{chapter-III}.

The comultiplication of a dual semigroup can also be used to define a convolution product. The convolution $j_1\star j_2$ of two unital $*$-algebra homomorphisms $j_1,j_2:\mathcal{B}\to\mathcal{A}$ is defined as
\[
j_1\star j_2=m_{\mathcal{A}}\circ\left(j_1\coprod j_2\right)\circ \Delta.
\]
As the composition of the three unital $*$-algebra homomorphisms $\Delta:\mathcal{B}\to\mathcal{B}\coprod\mathcal{B}$, $j_1\coprod j_2:\mathcal{B}\coprod\mathcal{B}\to\mathcal{A}\coprod\mathcal{A}$, and $m_{\mathcal{A}}:\mathcal{A}\coprod\mathcal{A}\to\mathcal{A}$, this is obviously again a unital $*$-algebra homomorphism. Note that this convolution can not be defined for arbitrary linear maps on $\mathcal{B}$ with values in some algebra, as for bialgebras, but only for unital $*$-algebra homomorphisms.

\section{Definition of L\'evy processes on dual groups}

\begin{definition}\label{def-indep}
Let $j_1:\mathcal{B}_1\to(\mathcal{A},\Phi),\ldots, j_n:\mathcal{B}_n\to(\mathcal{A},\Phi)$ be quantum random variables over the same quantum probability space $(\mathcal{A},\Phi)$ and denote their marginal distributions by $\varphi_i=\Phi\circ j_i$, $i=1,\ldots,n$. The quantum random variables $(j_1,\ldots,j_n)$ are called tensor independent (respectively boolean independent, monotonically independent, anti-monotonically independent or free), if the state $\Phi\circ m_\mathcal{A}\circ(j_1\coprod \cdots \coprod j_n)$ on the free product $\coprod_{i=1}^n \mathcal{B}_i$ is equal to the tensor product (boolean, monotone, anti-monotone, or free product, respectively) of $\varphi_1,\ldots,\varphi_n$.
\end{definition}

Note that tensor, boolean, and free independence do not depend on the order, but monotone and anti-monotone independence do. An $n$-tuple $(j_1,\ldots,j_n)$ of quantum random variables is monotonically independent, if and only if $(j_n,\ldots,j_1)$ is anti-monotonically independent.

We are now ready to define  tensor, boolean, monotone, anti-monotone, and free L{\'e}vy processes on dual semigroups.

\begin{definition}\label{def-levy-dual} \cite{schuermann95b}
Let $({\mathcal B},\Delta,\varepsilon)$ be a dual semigroup. A quantum stochastic process $\{j_{st}\}_{0\le s\le t\le T}$ on ${\mathcal B}$ over some quantum probability space $({\mathcal A},\Phi)$ is called a {\em tensor (resp.\ boolean, monotone, anti-monotone, or free) L{\'e}vy process on the dual semigroup ${\mathcal B}$}, if the following four conditions are satisfied.
\begin{enumerate}
\item
(Increment property) We have
\begin{eqnarray*}
j_{rs}\star j_{st} &=& j_{rt} \quad \mbox{ for all } 0\le r\le s\le
t\le T, \\
j_{tt} &=& \mathbf{1}_\mathcal{A}\circ \varepsilon \quad \mbox{ for all } 0\le t\le T.
\end{eqnarray*}
\item
(Independence of increments)
The family $\{j_{st}\}_{0\le s\le t\le T}$ is tensor independent (resp.\ boolean, monotonically, anti-monotonically independent, or free) w.r.t.\ $\Phi$, i.e.\ the $n$-tuple $(j_{s_1t_2}, \ldots, j_{s_nt_n})$ is tensor independent (resp.\ boolean, monotonically, anti-monotonically independent, or free) for all $n\in\mathbb{N}$ and all $0\le s_1\le t_1\le s_2\le \cdots\le t_n\le T$.
\item
(Stationarity of increments)
The distribution $\varphi_{st}=\Phi\circ j_{st}$ of $j_{st}$ depends
only on the difference $t-s$.
\item
(Weak continuity) The quantum random variables $j_{st}$ converge to
$j_{ss}$ in distribution for $t\searrow s$.
\end{enumerate}
\end{definition}

\begin{remark}
The independence property depends on the products and therefore for boolean,
monotone and anti-monotone L\'evy processes on the choice of a decomposition
$\mathcal{B}=\mathbb{C}\mathbf{1}\oplus\mathcal{B}^0$. In order to show that
the convolutions defined by $(\varphi_1\diamond\varphi_2)\circ \Delta$,
$(\varphi_1\triangleright\varphi_2)\circ \Delta$, and
$(\varphi_1\triangleleft\varphi_2)\circ \Delta$ are associative and that the
counit $\varepsilon$ acts as unit element w.r.t.\ these convolutions, one
has to use the functoriality property \cite[Condition
(P4)]{benghorbal+schuermann99}, see also \eqref{cond-functoriality} in Chapter \ref{chapter-III}. In our setting it is only satisfied for
morphisms that respect the decomposition. Therefore we are forced to choose
the decomposition given by $\mathcal{B}^0={\rm ker}\,\varepsilon$.
\end{remark}

The marginal distributions $\varphi_{t-s}:=\varphi_{st}=\Phi\circ j_{st}$ form
again a convolution semigroup $\{\varphi_t\}_{t\in\mathbb{R}_+}$, with respect to the tensor (boolean, monotone, anti-monotone, or free respectively) convolution defined by $(\varphi_1\tilde{\otimes}\varphi_2)\circ \Delta$ ($(\varphi_1\diamond\varphi_2)\circ \Delta$, $(\varphi_1\triangleright\varphi_2)\circ \Delta$, $(\varphi_1\triangleleft\varphi_2)\circ \Delta$, or $(\varphi_1*\varphi_2)\circ \Delta$, respectively). It has been shown that the generator $\psi:\mathcal{B}\to\mathbb{C}$,
\[
\psi(b)=\lim_{t\searrow 0} \frac{1}{t} \big(\varphi_t(b)-\varepsilon(b)\big)
\]
is well-defined for all $b\in\mathcal{B}$ and uniquely characterizes the semigroup $\{\varphi_t\}_{t\in\mathbb{R}_+}$, cf.\ \cite{schuermann95b,benghorbal+schuermann99,franz01}.

Denote by $S$ be the flip map $S:\mathcal{B}\coprod\mathcal{B}\to\mathcal{B}\coprod\mathcal{B}$, $S=m_{\mathcal{B}\coprod\mathcal{B}}\circ \left(i_2\coprod i_1\right)$, where $i_1,i_2:\mathcal{B}\to\mathcal{B}\coprod\mathcal{B}$ are the inclusions of $\mathcal{B}$ into the first and the second factor of the free product $\mathcal{B}\coprod\mathcal{B}$. The flip map $S$ acts on $i_1(a_1)i_2(b_1)\cdots i_2(b_n)\in \mathcal{B}\coprod\mathcal{B}$ with $a_1,\ldots,a_n,b_1,\ldots,b_n\in\mathcal{B}$ as $S\big(i_1(a_1)i_2(b_1)\cdots i_2(b_n)\Big)=i_2(a_1)i_1(b_1)\cdots i_1(b_n)$. If $j_1:\mathcal{B}\to\mathcal{A}_1$ and $j_2:\mathcal{B}\to\mathcal{A}_2$ are two unital $*$-algebra homomorphisms, then we have $\left(j_2\coprod j_1\right)\circ S=\gamma_{\mathcal{A}_1,\mathcal{A}_2}\circ \left(j_1\coprod j_2\right)$. Like for bialgebras, the opposite comultiplication $\Delta^{\rm op}=S\circ\Delta$ of a dual semigroup $(\mathcal{B},\Delta,\varepsilon)$ defines a new dual semigroup  $(\mathcal{B},\Delta^{\rm op},\varepsilon)$.

\begin{lemma}\label{lem-mon-anti-mon}
Let $\{j_{st}:\mathcal{B}\to(\mathcal{A},\Phi)\}_{0\le s\le t\le T}$ be a quantum stochastic process on a dual semigroup $(\mathcal{B},\Delta,\varepsilon)$ and define its time-reversed process $\{j^{\rm op}_{st}\}_{0\le s\le t\le T}$ by
\[
j^{\rm op}_{st} = j_{T-t,T-s}
\]
for $0\le s\le t\le T$.
\begin{itemize}
\item[(i)]
The process $\{j_{st}\}_{0\le s\le t\le T}$ is a tensor (boolean, free, respectively) L\'evy process on the dual semigroup $(\mathcal{B},\Delta,\varepsilon)$ if and only if the time-reversed process $\{j^{\rm op}_{st}\}_{0\le s\le t\le T}$ is a tensor (boolean, free, respectively) L\'evy process on the dual semigroup $(\mathcal{B},\Delta^{\rm op},\varepsilon)$.
\item[(ii)]
The process $\{j_{st}\}_{0\le s\le t\le T}$ is a monotone L\'evy process on the dual semigroup $(\mathcal{B},\Delta,\varepsilon)$ if and only if the time-reversed process $\{j^{\rm op}_{st}\}_{0\le s\le t\le T}$ is an anti-monotone L\'evy process on the dual semigroup $(\mathcal{B},\Delta^{\rm op},\varepsilon)$.
\end{itemize}
\end{lemma}
\begin{proof}
The equivalence of the stationarity and continuity property for $\{j_{st}\}_{0\le s\le t\le T}$ and $\{j^{\rm op}_{st}\}_{0\le s\le t\le T}$ is clear.

The increment property for $\{j_{st}\}_{0\le s\le t\le T}$ with respect to $\Delta$ is equivalent to the increment property of $\{j^{\rm op}_{st}\}_{0\le s\le t\le T}$ with respect to $\Delta^{\rm op}$, since
\begin{eqnarray*}
m_\mathcal{A}\circ \left(j^{\rm op}_{st}\coprod j^{\rm op}_{tu}\right)\circ \Delta^{\rm op} &=& m_\mathcal{A}\circ \left(j_{T-t,T-s}\coprod j_{T-u,T-t}\right)\circ S\circ \Delta \\
&=& m_\mathcal{A}\circ \gamma_{\mathcal{A},\mathcal{A}}\circ\left(j_{T-u,T-t}\coprod j_{T-t,T-s}\right)\circ\Delta \\
&=& m_\mathcal{A}\circ\left(j_{T-u,T-t}\coprod j_{T-t,T-s}\right)\circ\Delta
\end{eqnarray*}
for all $0\le s\le t\le u\le T$.

If $\{j_{st}\}_{0\le s\le t\le T}$ has monotonically independent increments, i.e.\ if  the $n$-tuples $(j_{s_1t_2},$ $ \ldots, j_{s_nt_n})$ are monotonically independent for all $n\in\mathbb{N}$ and all $0\le s_1\le t_1\le s_2\le \cdots\le t_n$, then the $n$-tuples $(j_{s_nt_n}, \ldots, j_{s_1t_1})=(j^{\rm op}_{T-t_n,T-s_n}, \ldots, j^{\rm op}_{T-t_1,T-s_1})$ are anti-monotonically independent and $\{j^{\rm op}_{st}\}_{0\le s\le t\le T}$ has anti-monotonically independent increments, and vice versa.

Since tensor and boolean independence and freeness do not depend on the order, $\{j_{st}\}_{0\le s\le t\le T}$ has tensor (boolean, free, respectively) independent increments, if and only $\{j^{\rm op}_{st}\}_{0\le s\le t\le T}$ has tensor (boolean, free, respectively) independent increments.
\end{proof}

Before we study boolean, monotone, and anti-monotone L\'evy processes in more detail, we will show how the theory of tensor L\'evy processes on dual semigroups reduces to the theory of L\'evy processes on involutive bialgebras, see also \cite{schuermann95b}. If quantum random variables $j_1,\ldots, j_n$ are independent in the sense of Condition 2 in Definition I.\ref{I-levy-def}, then they are also tensor independent in the sense of Definition \ref{def-indep}. Therefore every L\'evy process on the bialgebra $(\mathcal{B},\overline{\Delta},\varepsilon)$ associated to a dual semigroup $(\mathcal{B},\Delta,\varepsilon)$ is automatically also a tensor L\'evy process on the dual semigroup $(\mathcal{B},\Delta,\varepsilon)$. To verify this, it is sufficient to note that the increment property in Definition I.\ref{I-levy-def} with respect to $\overline{\Delta}$ and the commutativity of the increments imply the increment property in Definition \ref{def-levy-dual} with respect to $\Delta$.

But tensor independence in general does not imply independence in the sense of Condition 2 in Definition I.\ref{I-levy-def}, because the commutation relations are not necessarily satisfied. Therefore, in general, a tensor L\'evy process on a dual semigroup $(\mathcal{B},\Delta,\varepsilon)$ will {\em not} be a L\'evy process on the involutive bialgebra $(\mathcal{B},\overline{\Delta},\varepsilon)$. But we can still associate an equivalent L\'evy process on the involutive bialgebra $(\mathcal{B},\overline{\Delta},\varepsilon)$ to it. To do this, note that the convolutions of two unital functionals $\varphi_1,\varphi_2:\mathcal{B}\to\mathbb{C}$ with respect to the dual semigroup structure and the tensor product and with respect to the bialgebra structure coincide, i.e.\
\[
(\varphi_1\tilde{\otimes}\varphi_2)\circ \Delta = (\varphi_1\otimes\varphi_2)\circ \overline{\Delta}.
\]
for all unital functionals $\varphi_1,\varphi_2:\mathcal{B}\to\mathbb{C}$. Therefore the semigroup of marginal distributions of a tensor L\'evy process on the dual semigroup $(\mathcal{B},\Delta,\varepsilon)$ is also a convolution semigroup of states on the involutive bialgebra $(\mathcal{B},\overline{\Delta},\varepsilon)$. It follows that there exists a unique (up to equivalence) L\'evy process on the involutive bialgebra $(\mathcal{B},\overline{\Delta},\varepsilon)$ that has this semigroup as marginal distributions. It is easy to check that this process is equivalent to the given tensor L\'evy process on the dual semigroup  $(\mathcal{B},\Delta,\varepsilon)$. We summarize our result in the following theorem.

\begin{theorem}\label{theo-tensor}
Let $(\mathcal{B},\Delta,\varepsilon)$ be a dual semigroup, and $(\mathcal{B},\overline{\Delta},\varepsilon)$ with $\overline{\Delta}=R_{\mathcal{B},\mathcal{B}}\circ\Delta$ the associated involutive bialgebra. The tensor L\'evy processes on the dual semigroup $(\mathcal{B},\Delta,\varepsilon)$ are in one-to-one correspondence (up to equivalence) with the L\'evy processes on the involutive bialgebra $(\mathcal{B},\overline{\Delta},\varepsilon)$.

Furthermore, every L\'evy process on the involutive bialgebra $(\mathcal{B},\overline{\Delta},\varepsilon)$ is also a tensor L\'evy process on the dual semigroup $(\mathcal{B},\Delta,\varepsilon)$.
\end{theorem}

\section{Reduction of boolean, monotone, and  anti-monotone L\'evy
  processes to L\'evy processes on involutive bialgebras}\label{IV-reduction}
\markright{{\sc REDUCTION OF BOOLEAN, MONOTONE, $\dots$}}

In this section we will construct three involutive bialgebras for every dual semigroup $(\mathcal{B},\Delta,\varepsilon)$ and establish a one-to-one correspondence between boolean, monotone, and anti-monotone L\'evy processes on the dual semigroup $(\mathcal{B},\Delta,\varepsilon)$ and a certain class of L\'evy processes on one of those involutive bialgebras.

We start with some general remarks.

Let $(\mathcal{C},\square)$ be a tensor category. Then we call an object $\mathcal{D}$ in $\mathcal{S}$ equipped with morphisms
\[
\varepsilon:\mathcal{D}\to E,\qquad
\Delta:\mathcal{D}\to\mathcal{D}\square\mathcal{D}
\]
a {\em dual semigroup in $(\mathcal{C},\square)$}, if the following diagrams commute.
\[
\xymatrix{
& \mathcal{D} \ar[dr]^{\Delta}\ar[dl]_{\Delta} & \\
\mathcal{D}\square\mathcal{D}\ar[d]_{{\rm id}_\mathcal{D}\square\Delta} & & \mathcal{D}\square\mathcal{D}\ar[d]^{\Delta\square{\rm id}_\mathcal{D}} \\
\mathcal{D}\square(\mathcal{D}\square\mathcal{D})\ar[rr]_{\alpha_{\mathcal{D},\mathcal{D},\mathcal{D}}} & & (\mathcal{D}\square\mathcal{D})\square\mathcal{D}
}
\qquad
\xymatrix{
E\square\mathcal{D} \ar[ddr]_{\lambda_\mathcal{D}} & \mathcal{D}\square\mathcal{D} \ar[l]_{\varepsilon\square{\rm id}_\mathcal{D}} \ar[r]^{{\rm id}_\mathcal{D}\square\varepsilon} & \mathcal{D}\square E \ar[ddl]_{\rho_\mathcal{D}} \\
 & \mathcal{D}\ar[u]^{\Delta} \ar[d]^{{\rm id}}& \\
 & \mathcal{D} &
}
\]

\begin{proposition}\label{IV-prop-bialgebra}
Let $\mathcal{D}$ be a dual semigroup in a tensor category and let $F:\mathcal{C}\to\eufrak{Alg}$ be a cotensor functor with values in the category of unital algebras (equipped with the usual tensor product). Then $F(\mathcal{D})$ is a bialgebra with the counit $F_0\circ F(\varepsilon)$ and the coproduct $F_2(\mathcal{D},\mathcal{D})\circ F(\Delta)$.
\end{proposition}
\begin{proof}
We only prove one right half of the counit property. Applying $F$ to $\lambda_\mathcal{D}\circ (\varepsilon\square{\rm id}_\mathcal{D})\circ \Delta = {\rm id}_\mathcal{D}$, we get $F(\lambda_\mathcal{D})\circ F(\varepsilon\square{\rm id}_\mathcal{D})\circ F\Delta = {\rm id}_{F(\mathcal{D})}$. Using the naturality of $F_2$ and Diagram (III.\ref{III-cotensor-rho}), we can extend this to the following commutative diagram,
\[
\xymatrix{
F(\mathcal{D})\otimes F(\mathcal{D})\ar[rrr]^{{\rm id}_{F(\mathcal{D})}\otimes F(\varepsilon)} & & & F(\mathcal{D})\otimes F(E) \ar[dddll]^{{\rm id}_{F(\mathcal{D})}\otimes F_0} \\
F(\mathcal{D}\square\mathcal{D})\ar[u]^{F_2(\mathcal{D},\mathcal{D})} \ar[r]^{F({\rm id}_\mathcal{D}\square\varepsilon)} & F(\mathcal{D}\square E) \ar[ddl]^{F(\rho_\mathcal{D})} \ar[urr]^{F_2(\mathcal{D},E)} & & \\
F(\mathcal{D})\ar[u]^{F(\Delta)} \ar[d]_{{\rm id}_{F(\mathcal{D})}} & & & \\
F(\mathcal{D}) & F(\mathcal{D})\otimes \mathbb{C} \ar[l]_{\cong} & & &
}
\]
which proves the right counit property of $F(\mathcal{D})$. The proof of the left counit property is of course done by taking the mirror image of this diagram and replacing $\rho$ by $\lambda$. The proof of the coassociativity requires a bigger diagram which makes use of (\ref{III-cotensor-alpha}) in Chapter \ref{chapter-III}.
\end{proof}

Assume now that we have a family $(\mathcal{D}_t)_{t\ge 0}$ of objects in $\mathcal{C}$ equipped with morphisms $\varepsilon:\mathcal{D}_0\to E$ and $\delta_{st}:\mathcal{D}_{s+t}\to:\mathcal{D}_{s}\square\mathcal{D}_{t}$ for $s,t\ge 0$ such that the diagrams
\[
\xymatrix{
& \mathcal{D}_{s+t+u} \ar[dr]^{\delta_{s,t+u}}\ar[dl]_{\delta_{s+t,u}} & \\
\mathcal{D}_s\square\mathcal{D}_{t+u}\ar[d]_{{\rm id}\square\delta_{tu}} & & \mathcal{D}_{s+t}\square\mathcal{D}_u\ar[d]^{\delta_{st}\square{\rm id}} \\
\mathcal{D}\square(\mathcal{D}\square\mathcal{D})\ar[rr]_{\alpha_{\mathcal{D}_s,\mathcal{D}_t,\mathcal{D}_u}} & & (\mathcal{D}_s\square\mathcal{D}_t)\square\mathcal{D}_u
}
\qquad
\xymatrix{
\mathcal{D}_0\square\mathcal{D}_t \ar[dd]_{\varepsilon\square{\rm id}} & \mathcal{D}_t \ar[dd]_{{\rm id}} \ar[l]_{\delta_{0t}}\ar[r]^{\delta_{t0}} & \mathcal{D}_t\square\mathcal{D}_0 \ar[dd]^{{\rm id}\square\varepsilon} \\
 & & \\
E\square\mathcal{D}_t \ar[r]_{\lambda_{\mathcal{D}_t}}& \mathcal{D}_t & \mathcal{D}_t\square E \ar[l]^{\rho_{\mathcal{D}_t}}
}
\]
In the application we have in mind the objects $\mathcal{D}_t$ will be pairs consisting of a fixed dual semigroup $\mathcal{B}$ and a state $\varphi_t$ on $\mathcal{B}$ that belongs to a convolution semigroup $(\Phi_t)_{t\ge 0}$ on $\mathcal{B}$. The morphisms $\delta_{st}$ and $\varepsilon$ will be the coproduct and the counit of $\mathcal{B}$.

If there exists a cotensor functor $F:\mathcal{C}\to\eufrak{AlgProb}$, $F(\mathcal{D}_t)=(\mathcal{A}_t,\varphi_t)$ such that the algebras ${\rm Alg}\big(F(\mathcal{D}_t)\big)=\mathcal{A}_t$ and the morphisms $F_2(\mathcal{D}_s,\mathcal{D}_t)\circ F(\delta_{st})$ are do not depend on $s$ and $t$, then $\mathcal{A}={\rm Alg}\big(F(\mathcal{D}_t)\big)$ is again a bialgebra with coproduct $\tilde{\Delta}=F_2(\mathcal{D}_s,\mathcal{D}_t)\circ F(\delta_{st})$ and the counit $\tilde{\varepsilon}=F_0\circ F(\varepsilon)$, as in Proposition \ref{IV-prop-bialgebra}.

Since morphisms in $\eufrak{AlgProb}$ leave the states invariant, we have $\Phi_s\otimes\Phi_t\circ \tilde{\Delta}=\Phi_{s+t}$ and $\Phi_0=\tilde{\varepsilon}$, i.e.\ $(\Phi_t)_{t\ge 0}$ is a convolution semigroup on $\mathcal{A}$ (up to the continuity property). From this convolution semigroup we can then reconstruct a L\'evy process on $\mathcal{A}$.

\subsection{Construction of a L\'evy process on an involutive bialgebra}

After the category theoretical considerations of the previous subsection we shall now explicitly construct one-to-one correspondences between boolean, monotone, and anti-monotone L\'evy processes on dual groups and certain classes of L\'evy processes on involutive bialgebras.

Let $M=\{\mathbf{1},p\}$ be the unital semigroup with two elements and the
multiplication $p^2=\mathbf{1}p=p\mathbf{1}=p$, $\mathbf{1}^2=\mathbf{1}$. Its
`group algebra' $\mathbb{C}M=\mbox{span}\,\{\mathbf{1},p\}$ is an involutive
bialgebra with comultiplication
$\Delta(\mathbf{1})=\mathbf{1}\otimes\mathbf{1}$, $\Delta(p)=p\otimes p$,
counit $\varepsilon(\mathbf{1})=\varepsilon(p)=1$, and involution
$\mathbf{1}^*=\mathbf{1}$, $p^*=p$. The involutive bialgebra $\mathbb{C}M$ was
already used by Lenczweski \cite{lenczewski98b,lenczewski01} to give a tensor
product construction for a large family of products of quantum probability
spaces including the boolean and the free product and to define and study the
additive convolutions associated to these products. As a unital $*$-algebra it
is also used in Skeide's approach to boolean calculus, cf.\ \cite{skeide01},
where it is introduced as the unitization of $\mathbb{C}$. It also plays an
important role in \cite{schuermann00,franz+schuermann00}.

Let $\mathcal{B}$ be a unital $*$-algebra, then we define its $p$-extension $\tilde{\mathcal{B}}$ as the free product $\tilde{\mathcal{B}}=\mathcal{B}\coprod\mathbb{C}M$. Due to the identification of the units of $\mathcal{B}$ and $\mathbb{C}M$, any element of $\tilde{\mathcal{B}}$ can be written as sums of products of the form $p^\alpha b_1 p b_2 p\cdots pb_np^\omega$ with $n\in\mathbb{N}$, $b_1,\ldots,b_n\in\mathcal{B}$ and $\alpha,\omega=0,1$. This representation can be made unique, if we choose a decomposition of $\mathcal{B}$ into a direct sum of vector spaces $\mathcal{B}=\mathbb{C}\mathbf{1}\oplus\mathcal{V}^0$ and require $b_1,\ldots,b_n\in\mathcal{V}^0$. We define the $p$-extension $\tilde{\varphi}:\tilde{\mathcal{B}}\to\mathbb{C}$ of a unital functional $\varphi:\mathcal{B}\to\mathbb{C}$ by
\begin{equation}\label{bool-ext}
\tilde{\varphi}(p^\alpha b_1 p b_2 p\cdots pb_np^\omega)=\varphi(b_1)\varphi(b_2)\cdots\varphi(b_n)
\end{equation}
and $\tilde{\varphi}(p)=1$. The $p$-extension does not depend on the
decomposition $\mathcal{B}=\mathbb{C}\mathbf{1}\oplus\mathcal{V}^0$, since
Equation \eqref{bool-ext} actually holds not only for $b_1,\ldots,b_n\in\mathcal{V}^0$, but also for $b_1,\ldots,b_n\in\mathcal{B}$.

If $\mathcal{B}_1,\ldots,\mathcal{B}_n$ are unital $*$-algebras that can be written as direct sums $\mathcal{B}_i=\mathbb{C}\mathbf{1}\oplus\mathcal{B}_i^0$ of $*$-algebras, then we can define unital $*$-algebra homomorphisms $I_{k,\mathcal{B}_1,\ldots,\mathcal{B}_n}^{\rm B}, I_{k,\mathcal{B}_1,\ldots,\mathcal{B}_n}^{\rm M},$ $I_{k,\mathcal{B}_1,\ldots,\mathcal{B}_n}^{\rm AM}:\mathcal{B}_k\to \tilde{\mathcal{B}}_1\otimes \cdots\otimes \tilde{\mathcal{B}}_n$ for $k=1,\ldots,n$ by
\begin{eqnarray*}
I_{k,\mathcal{B}_1,\ldots,\mathcal{B}_n}^{\rm B}(b) &=& \underbrace{p\otimes \cdots \otimes p} \otimes b \otimes \underbrace{p\otimes \cdots \otimes p}, \\
&& \hspace*{5mm} k-1 \mbox{ times } \hspace*{9mm} n-k \mbox{ times } \\
I_{k,\mathcal{B}_1,\ldots,\mathcal{B}_n}^{\rm M}(b) &=& \underbrace{ \mathbf{1}\otimes \cdots \otimes\mathbf{1}} \otimes b \otimes \underbrace{p\otimes \cdots \otimes p}, \\
&& \hspace*{5mm} k-1 \mbox{ times } \hspace*{9mm} n-k \mbox{ times } \\
I_{k,\mathcal{B}_1,\ldots,\mathcal{B}_n}^{\rm AM}(b) &=& \underbrace{p\otimes \cdots \otimes p} \otimes b \otimes \underbrace{\mathbf{1}\otimes \cdots \otimes\mathbf{1}}, \\
&& \hspace*{5mm} k-1 \mbox{ times } \hspace*{9mm} n-k \mbox{ times }
\end{eqnarray*}
for $b\in\mathcal{B}_k^0$.

Let $n\in\mathbb{N}$, $1\le k\le n$, and denote the canonical inclusions of $\mathcal{B}_k$ into the $k^{\rm th}$ factor of the free product $\coprod_{j=1}^n\mathcal{B}_j$ by $i_k$. Then, by the universal property, there exist unique unital $*$-algebra homomorphisms $R_{\mathcal{B}_1,\ldots,\mathcal{B}_n}^{\bullet}:\coprod_{k=1}^n\mathcal{B}_k\to \otimes_{k=1}^n\tilde{\mathcal{B}}_k$ such that
\[
R_{\mathcal{B}_1,\ldots,\mathcal{B}_n}^{\bullet}\circ i_k = I_{k,\mathcal{B}_1,\ldots,\mathcal{B}_n}^{\bullet},
\]
for $\bullet\in\{{\rm B},{\rm M},{\rm AM}\}$.

\begin{proposition}
Let $(\mathcal{B},\Delta,\varepsilon)$ be a dual semigroup. Then we have the following three involutive bialgebras $(\tilde{\mathcal{B}},\overline{\Delta}_{\rm B},\tilde{\varepsilon})$, $(\tilde{\mathcal{B}},\overline{\Delta}_{\rm M},\tilde{\varepsilon})$, and $(\tilde{\mathcal{B}},\overline{\Delta}_{\rm AM},\tilde{\varepsilon})$, where the comultiplications are defined by
\begin{eqnarray*}
\overline{\Delta}_{\rm B} &=& R^{\rm B}_{\mathcal{B},\mathcal{B}}\circ\Delta, \\
\overline{\Delta}_{\rm M} &=& R^{\rm M}_{\mathcal{B},\mathcal{B}}\circ\Delta,\\
\overline{\Delta}_{\rm AM} &=& R^{\rm AM}_{\mathcal{B},\mathcal{B}}\circ\Delta,
\end{eqnarray*}
on $\mathcal{B}$ and by
\[
\overline{\Delta}_{\rm B} (p) = \overline{\Delta}_{\rm M} (p) = \overline{\Delta}_{\rm AM}  (p) = p\otimes p
\]
on $\mathbb{C}M$.
\end{proposition}
\begin{remark}
This is actually an application of Proposition \ref{IV-prop-bialgebra}. Below we give an direct proof.
\end{remark}
\begin{proof}
We will prove that $(\tilde{\mathcal{B}},\overline{\Delta}_{\rm B},\tilde{\varepsilon})$ is an involutive bialgebra, the proofs for $(\tilde{\mathcal{B}},\overline{\Delta}_{\rm M},\tilde{\varepsilon})$ and $(\tilde{\mathcal{B}},\overline{\Delta}_{\rm AM},\tilde{\varepsilon})$ are similar.

It is clear that $\overline{\Delta}_{\rm B}:\tilde{\mathcal{B}}\to \tilde{\mathcal{B}}\otimes\tilde{\mathcal{B}}$ and $\tilde{\varepsilon}:\tilde{\mathcal{B}}\to\mathbb{C}$ are unital $*$-algebra homomorphisms, so we only have to check the coassociativity and the counit property. That they are satisfied for $p$ is also immediately clear. The proof for elements of $\mathcal{B}$ is similar to the proof of \cite[Theorem 4.2]{zhang91}. We get
\begin{eqnarray*}
\left.\left(\overline{\Delta}_{\rm B}\otimes{\rm id}_{\tilde{\mathcal{B}}}\right)\circ \overline{\Delta}_{\rm B}\right|_{\mathcal{B}} &=& R^{\rm B}_{\mathcal{B},\mathcal{B},\mathcal{B}}\circ \left(\Delta\coprod{\rm id}_{\mathcal{B}}\right)\circ \Delta  \\
&=& R^{\rm B}_{\mathcal{B},\mathcal{B},\mathcal{B}}\circ \left({\rm id}_{\mathcal{B}}\coprod\Delta\right)\circ \Delta  \\
&=& \left.\left({\rm id}_{\tilde{\mathcal{B}}}\otimes\overline{\Delta}_{\rm B}\right)\circ \overline{\Delta}_{\rm B}\right|_{\mathcal{B}} 
\end{eqnarray*}
and
\begin{eqnarray*}
\left.\left(\tilde\varepsilon\otimes {\rm id}_{\tilde{\mathcal{B}}}\right)\circ \overline{\Delta}_{\rm B}\right|_{\mathcal{B}}  &=& \left(\tilde\varepsilon\otimes {\rm id}_{\tilde{\mathcal{B}}}\right)\circ R^{\rm B}_{\mathcal{B},\mathcal{B}}\circ \Delta \\
&=& \left(\varepsilon\coprod {\rm id}_{\mathcal{B}}\right)\circ \Delta
\,\,=\,\, {\rm id}_\mathcal{B} \\
&=& \left({\rm id}_{\mathcal{B}}\coprod\varepsilon\right)\circ \Delta \\
&=&  \left({\rm id}_{\tilde{\mathcal{B}}}\otimes\tilde\varepsilon\right)\circ R^{\rm B}_{\mathcal{B},\mathcal{B}}\circ \Delta \\
&=& \left.\left({\rm id}_{\tilde{\mathcal{B}}}\otimes\tilde\varepsilon\right)\circ \overline{\Delta}_{\rm B}\right|_{\mathcal{B}}.
\end{eqnarray*}
\end{proof}

These three involutive bialgebras are important for us, because the boolean convolution (monotone convolution, anti-monotone convolution, respectively) of unital functionals on a dual semigroup $(\mathcal{B},\Delta,\varepsilon)$ becomes the convolution with respect to the comultiplication $\overline{\Delta}_{\rm B}$ ($\overline{\Delta}_{\rm M}$, $\overline{\Delta}_{\rm AM}$, respectively) of their $p$-extension on $\tilde{\mathcal{B}}$.

\begin{proposition}\label{prop-conv}
Let $(\mathcal{B},\Delta,\varepsilon)$ be a dual semigroup and $\varphi_1,\varphi_2:\mathcal{B}\to\mathbb{C}$ two unital functionals on $\mathcal{B}$. Then we have
\begin{eqnarray*}
\widetilde{(\varphi_1\diamond\varphi_2)\circ \Delta} &=& (\tilde\varphi_1\otimes\tilde\varphi_2)\circ \overline{\Delta}_{\rm B}, \\
\widetilde{(\varphi_1\triangleright\varphi_2)\circ \Delta} &=& (\tilde\varphi_1\otimes\tilde\varphi_2)\circ \overline{\Delta}_{\rm M}, \\
\widetilde{(\varphi_1\triangleleft\varphi_2)\circ \Delta} &=& (\tilde\varphi_1\otimes\tilde\varphi_2)\circ \overline{\Delta}_{\rm AM}.
\end{eqnarray*}
\end{proposition}
\begin{proof}
Let $b\in\mathcal{B}^0$. As an element of $\mathcal{B}\coprod\mathcal{B}$, $\Delta(b)$ can be written in the form $\Delta(b)=\sum_{\epsilon\in\mathbb{A}} b^\epsilon\in\bigoplus_{\epsilon\in\mathbb{A}}\mathcal{B}_\epsilon$. Only finitely many terms of this sum are non-zero. The individual summands are tensor products $b^\epsilon =b^\epsilon_1\otimes\cdots \otimes b^\epsilon_{|\epsilon|}$ and due to the counit property we have $b^\emptyset=0$. Therefore we have
\[
(\varphi_1\diamond\varphi_2)\circ \Delta(b) = \sum_{{\epsilon\in\mathbb{A}}\atop{\epsilon\not=\emptyset}} \prod_{k=1}^{|\epsilon|} \varphi_{\epsilon_k}(b^\epsilon_k).
\]
For the right-hand-side, we get the same expression on $\mathcal{B}$,
\begin{eqnarray*}
(\tilde\varphi_1\otimes\tilde\varphi_2)\circ \overline{\Delta}_{\rm B}(b) &=&(\tilde\varphi_1\otimes\tilde\varphi_2)\circ R^{\rm B}_{\mathcal{B},\mathcal{B}}\circ\Delta(b) \\
&=& (\tilde\varphi_1\otimes\tilde\varphi_2)\circ R^{\rm B}_{\mathcal{B},\mathcal{B}} \sum_{\epsilon\in\mathbb{A}}b^\epsilon_1\otimes\cdots \otimes b^\epsilon_{|\epsilon|} \\
&=& \sum_{{\epsilon\in\mathbb{A}}\atop{\epsilon_1=1}}\tilde\varphi_1(b^\epsilon_1pb^\epsilon_3 \cdots )\tilde\varphi_2 (pb^\epsilon_2p\cdots) \\
&& + \sum_{{\epsilon\in\mathbb{A}}\atop{\epsilon_1=2}}\tilde\varphi_1(pb^\epsilon_2p \cdots )\tilde\varphi_2 (b^\epsilon_1pb^\epsilon_3\cdots) \\
&=& \sum_{{\epsilon\in\mathbb{A}}\atop{\epsilon\not=\emptyset}} \prod_{k=1}^{|\epsilon|} \varphi_{\epsilon_k}(b^\epsilon_k).
\end{eqnarray*}
To conclude, observe
\[
\overline{\Delta}_{\rm B}(p^\alpha b_1 p \cdots pb_np^\omega)= (p^\alpha\otimes p^\alpha) \overline{\Delta}_{\rm B}(b_1) (p\otimes p)\cdots (p\otimes p)\overline{\Delta}_{\rm B}(b_n)(p^\omega\otimes p^\omega)
\]
for all $b_1,\ldots,b_n\in\mathcal{B}$, $\alpha,\omega\in\{0,1\}$, and therefore
\[
(\tilde\varphi_1\otimes\tilde\varphi_2)\circ \overline{\Delta}_{\rm B} = \widetilde{\left.(\tilde\varphi_1\otimes\tilde\varphi_2)\circ \overline{\Delta}_{\rm B}\right|_{\mathcal{B}}} = \widetilde{(\varphi_1\diamond\varphi_2)\circ \Delta}.
\]
The proof for the monotone and anti-monotone product is similar.
\end{proof}

We can now state our first main result.

\begin{theorem}\label{theo-one-to-one}
Let $(\mathcal{B},\Delta,\varepsilon)$ be a dual semigroup. We have a one-to-one correspondence between boolean (monotone, anti-monotone, respectively) L\'evy processes on the dual semigroup $(\mathcal{B},\Delta,\varepsilon)$ and L\'evy processes on the involutive bialgebra $(\tilde{\mathcal{B}},\overline{\Delta}_{\rm B},\tilde{\varepsilon})$ ($(\tilde{\mathcal{B}},\overline{\Delta}_{\rm M},\tilde{\varepsilon})$, $(\tilde{\mathcal{B}},\overline{\Delta}_{\rm AM},\tilde{\varepsilon})$, respectively), whose marginal distributions satisfy
\begin{equation}\label{eq-p-invariance}
\varphi_t(p^\alpha b_1 p \cdots pb_np^\omega) = \varphi_t(b_1)\cdots\varphi_t(b_n) 
\end{equation}
for all $t\ge 0$, $b_1,\ldots,b_n\in\mathcal{B}$, $\alpha,\omega\in\{0,1\}$.
\end{theorem}
\begin{proof}
Condition (\ref{eq-p-invariance}) says that the functionals $\varphi_t$ on $\tilde{\mathcal{B}}$ are equal to the $p$-extension of their restriction to $\mathcal{B}$.

Let $\{j_{st}\}_{0\le s\le t\le T}$ be a boolean (monotone, anti-monotone, respectively) L\'evy process on the dual semigroup $(\mathcal{B},\Delta,\varepsilon)$ with convolution semigroup $\varphi_{t-s}=\Phi\circ j_{st}$. Then, by Proposition \ref{prop-conv}, their $p$-extensions $\{\tilde\varphi_t\}_{t\ge 0}$ form a convolution semigroup on the involutive bialgebra $(\tilde{\mathcal{B}},\overline{\Delta}_{\rm B},\tilde\varepsilon)$ ($(\tilde{\mathcal{B}},\overline{\Delta}_{\rm M},\tilde\varepsilon)$, $(\tilde{\mathcal{B}},\overline{\Delta}_{\rm AM},\tilde\varepsilon)$, respectively). Thus there exists a unique (up to equivalence) L\'evy process $\{\tilde\jmath_{st}\}_{0\le s\le t\le T}$ on the involutive bialgebra $(\tilde{\mathcal{B}},\overline{\Delta}_{\rm B},\tilde\varepsilon)$ ($(\tilde{\mathcal{B}},\overline{\Delta}_{\rm M},\tilde\varepsilon)$, $(\tilde{\mathcal{B}},\overline{\Delta}_{\rm AM},\tilde\varepsilon)$, respectively) with these marginal distribution.

Conversely, let $\{j_{st}\}_{0\le s\le t\le T}$ be a L\'evy process on the involutive bialgebra $(\tilde{\mathcal{B}},\overline{\Delta}_{\rm B},\tilde\varepsilon)$ ($(\tilde{\mathcal{B}},\overline{\Delta}_{\rm M},\tilde\varepsilon)$, $(\tilde{\mathcal{B}},\overline{\Delta}_{\rm AM},\tilde\varepsilon)$, respectively) with marginal distributions $\{\varphi_t\}_{t\ge 0}$ and suppose that the functionals $\varphi_t$ satisfy Equation (\ref{eq-p-invariance}). Then, by  Proposition \ref{prop-conv}, their restrictions to $\mathcal{B}$ form a convolution semigroup on the dual semigroup $(\mathcal{B},\Delta,\varepsilon)$ with respect to the boolean (monotone, anti-monotone, respectively) convolution and therefore there exists a unique (up to equivalence) boolean (monotone, anti-monotone, respectively) L\'evy process on the dual semigroup $(\mathcal{B},\Delta,\varepsilon)$ that has these marginal distributions.

The correspondence is one-to-one, because the $p$-extension establishes a bijection between unital functionals on $\mathcal{B}$ and unital functionals on $\tilde{\mathcal{B}}$ that satisfy Condition (\ref{eq-p-invariance}). Furthermore, a unital functional on  $\mathcal{B}$ is positive if and only if its $p$-extension is positive on $\tilde{\mathcal{B}}$.
\end{proof}

We will now reformulate Equation (\ref{eq-p-invariance}) in terms of the
generator of the process. Let $n\ge 1$,   $b_1,\ldots,b_n\in\mathcal{B}^0={\rm
  ker}\,\varepsilon$, $\alpha,\omega\in\{0,1\}$, then we have
\begin{eqnarray*}
\psi(p^\alpha b_1 p \cdots pb_np^\omega) &=& \lim_{t\searrow 0} \frac{1}{t} \big(\varphi_t(p^\alpha b_1 p \cdots pb_np^\omega)-\tilde\varepsilon(p^\alpha b_1 p \cdots pb_np^\omega)\big) \\
&=& \lim_{t\searrow 0} \frac{1}{t}\big(\varphi_t(b_1) \cdots \varphi_t(b_n)-\varepsilon(b_1)\cdots\varepsilon(b_n)\big) \\
&=& \sum_{k=1}^n \varepsilon(b_1)\cdots\varepsilon(b_{k-1})\psi(b_k)\varepsilon(b_{k+1})\cdots\varepsilon(b_n) \\
&=&
\left\{\begin{array}{lcl}
\psi(b_1) & \mbox{ if }& n=1, \\
0 & \mbox{ if }& n>1.
\end{array}
\right.
\end{eqnarray*}
Conversely, let $\{\varphi_t:\tilde{\mathcal{B}}\to\mathbb{C}\}_{t\ge 0}$ be a convolution semigroup on $(\tilde{\mathcal{B}},\overline{\Delta}_\bullet,\tilde\varepsilon)$, $\bullet\in\{{\rm B},{\rm M},{\rm AM}\}$, whose generator $\psi:\tilde{\mathcal{B}}\to\mathbb{C}$ satisfies $\psi(\mathbf{1})=\psi(p)=0$ and
\begin{equation}\label{eq-psi-p-inv}
\psi(p^\alpha b_1 p \cdots pb_np^\omega) =
\left\{\begin{array}{lcl}
\psi(b_1) & \mbox{ if }& n=1, \\
0 & \mbox{ if }& n>1,
\end{array}
\right.
\end{equation}
for all $n\ge 1$, $b_1,\ldots,b_n\in\mathcal{B}^0={\rm ker}\,\varepsilon$,
$\alpha,\omega\in\{0,1\}$. For $b_1,\ldots,b_n\in\mathcal{B}^0$,
$\overline{\Delta}_\bullet(b_i)$ is of the form
$\overline{\Delta}_\bullet(b_i)= b_i\otimes \mathbf{1}+\mathbf{1}\otimes b_i +
\sum_{k=1}^{n_i} b^{(1)}_{i,k}\otimes b^{(2)}_{i,k}$, with $b^{(1)}_{i,k},
b^{(2)}_{i,k}\in\mbox{ker}\,\tilde\varepsilon$. By the fundamental theorem of
coalgebras \cite{sweedler69} there exists a finite-dimensional
subcoalgebra $\mathcal{C}\subseteq\tilde{\mathcal{B}}$ of
$\tilde{\mathcal{B}}$ that contains all possible products of
$\mathbf{1},b_i,b^{(1)}_{i,k_i}, b^{(2)}_{i,k_i}$, $i=1,\ldots,n$, $k_i=1,\ldots,n_i$.

Then we have
\begin{eqnarray*}
& \left.\varphi_{s+t}\right|_\mathcal{C}(p^\alpha b_1 p \cdots pb_np^\omega) = & \\
& (\left.\varphi_s\right|_\mathcal{C}\otimes\left.\varphi_t\right|_\mathcal{C})\big((p\alpha\otimes p^\alpha)\overline{\Delta}_\bullet(b_1)(p\otimes p) \cdots (p\otimes p) \overline{\Delta}_\bullet(b_n)(p^\omega\otimes p^\omega)\big) &
\end{eqnarray*}
and, using (\ref{eq-psi-p-inv}), we find the differential equation
\begin{eqnarray}
\left.\dot\varphi_{s}\right|_\mathcal{C}(p^\alpha b_1 p \cdots pb_np^\omega) &=& \sum_{i=1}^n \left.\varphi_s\right|_\mathcal{C}(p^\alpha b_1p\cdots b_{i-1}p\mathbf{1}pb_{i+1}p\cdots b_np^\omega)\psi(b_i) \nonumber\\
&& + \sum_{i=1}^n\sum_{k_i=1}^{n_i} \left.\varphi\right|_\mathcal{C}(p^\alpha b_1p\cdots b_{i-1}pb^{(1)}_{i,k_i}pb_{i+1}p\cdots b_np^\omega)\psi(b^{(2)}_{i,k_i}) \label{eq-diff}
\end{eqnarray}
for $\{\left.\varphi_t\right|_\mathcal{C}\}_{t\ge 0}$. This a linear inhomogeneous differential equation for a function with values in the finite-dimensional complex vector space $\mathcal{C}^*$ and it has a unique global solution for every initial value $\left.\varphi_0\right|_{\mathcal{C}}$. Since we have
\begin{eqnarray*}
\dot\varphi_s(b_i) &=& (\varphi_s\otimes\psi)\left(b_i\otimes \mathbf{1}+\mathbf{1}\otimes b_i + \sum_{k=1}^{n_i} b^{(1)}_{i,k}\otimes b^{(2)}_{i,k}\right) \\
&=&\psi(b_i) +  \sum_{k_i=1}^{n_i}\varphi_s(b^{(1)}_{i,k_i})\psi(b^{(2)}_{i,k_i}),
\end{eqnarray*}
we see that $\left\{\left.\widetilde{(\left.\varphi_t\right|_\mathcal{B})}\right|_\mathcal{C}\right\}_{t\ge 0}$ satisfies the differential equation (\ref{eq-diff}). The initial values also agree,
\[
\varphi_0(p^\alpha b_1 p \cdots pb_np^\omega) =\tilde\varepsilon(p^\alpha b_1 p \cdots pb_np^\omega) = \varepsilon(b_1)\cdots \varepsilon(b_n)=\varphi_0(b_1)\cdots \varphi_0(b_n)
\]
and therefore it follows that $\{\varphi_t\}_{t\ge 0}$ satisfies Condition (\ref{eq-p-invariance}).

We have shown the following.

\begin{lemma}
Let $\{\varphi_t:\tilde{\mathcal{B}}\to\mathbb{C}\}_{t\ge 0}$ be a convolution semigroup of unital functionals on the involutive bialgebra $(\tilde{\mathcal{B}},\overline{\Delta}_\bullet,\tilde\varepsilon)$, $\bullet\in\{{\rm B},{\rm M},{\rm AM}\}$, and let $\psi:\tilde{\mathcal{B}}\to\mathbb{C}$ be its infinitesimal generator.

Then the functionals of the convolution semigroup $\{\varphi_t\}_{t\ge 0}$ satisfy (\ref{eq-p-invariance}) for all $t\ge 0$, if and only if its generator $\psi$ satisfies (\ref{eq-psi-p-inv}).
\end{lemma}

For every linear functional $\psi:\mathcal{B}\to\mathbb{C}$ on $\mathcal{B}$ there exists only one unique functional $\hat\psi:\tilde{\mathcal{B}}\to\mathbb{C}$ with $\hat\psi|_\mathcal{B}=\psi$ that satisfies Condition (\ref{eq-psi-p-inv}). And since this functional $\hat\psi$ is hermitian and conditionally positive, if and only if $\psi$ is hermitian and conditionally positive, we have shown the following.

\begin{corollary}
We have a one-to-one correspondence between boolean L\'evy processes, monotone L\'evy processes, and anti-monotone L\'evy processes on a dual semigroup $(\mathcal{B},\Delta,\varepsilon)$ and generators, i.e.\ hermitian, conditionally positive, linear functionals $\psi:\mathcal{B}\to\mathbb{C}$ on $\mathcal{B}$ with $\psi(\mathbf{1})=0$.
\end{corollary}

Another corollary of Theorem \ref{theo-one-to-one} is the Schoenberg correspondence for the boolean, monotone, and anti-monotone convolution.

\begin{corollary}{\bf (Schoenberg correspondence)}\label{cor-schoenberg} Let $\{\varphi_t\}_{t\ge 0}$ be a convolution semigroup of unital functionals with respect to the tensor, boolean, monotone, or anti-monotone convolution on a dual semigroup $(\mathcal{B},\Delta,\varepsilon)$ and let $\psi:\mathcal{B}\to\mathbb{C}$ be defined by
\[
\psi(b)=\lim_{t\searrow 0} \frac{1}{t} \big(\varphi_t(b)-\varepsilon(b)\big)
\]
for $b\in\mathcal{B}$. Then the following statements are equivalent.
\begin{itemize}
\item[(i)]
$\varphi_t$ is positive for all $t\ge 0$.
\item[(ii)]
$\psi$ is hermitian and conditionally positive.
\end{itemize}
\end{corollary}

We have now obtained a classification of boolean, monotone, and anti-monotone L\'evy processes on a given dual semigroup in terms of a class of L\'evy processes on a certain involutive bialgebra and in terms of their generators. In the next subsection we will see how to construct realizations.

\subsection{Construction of boolean, monotone, and anti-monotone L\'evy processes}\label{sub-real}

The following theorem gives us a way to construct realizations of boolean, monotone, and anti-monotone L\'evy processes.

\begin{theorem}\label{theo-real}
Let $\{k_{st}^{\rm B}\}_{0\le s\le t\le T}$ ($\{k_{st}^{\rm M}\}_{0\le s\le t\le T}$, $\{k_{st}^{\rm AM}\}_{0\le s\le t\le T}$, respectively) be a boolean (monotone, anti-monotone, respectively) L\'evy process with generator $\psi$ on some dual semigroup $(\mathcal{B},\Delta,\varepsilon)$. Denote the unique extension of $\psi:\mathcal{B}\to\mathbb{C}$ determined by Equation (\ref{eq-psi-p-inv}) by $\hat\psi:\tilde{\mathcal{B}}\to\mathbb{C}$.

If $\{\tilde\jmath_{st}^{\rm B}\}_{0\le s\le t\le T}$ ($\{\tilde\jmath_{st}^{\rm M}\}_{0\le s\le t\le T}$, $\{\tilde\jmath_{st}^{\rm AM}\}_{0\le s\le t\le T}$, respectively) is a L\'evy process on the involutive bialgebra $(\tilde{\mathcal{B}},\overline{\Delta}_{\rm B},\tilde\varepsilon)$ ($(\tilde{\mathcal{B}},\overline{\Delta}_{\rm M},\tilde\varepsilon)$, 
$(\tilde{\mathcal{B}},\overline{\Delta}_{\rm AM},\tilde\varepsilon)$, respectively), then the quantum stochastic process $\{j_{st}^{\rm B}\}_{0\le s\le t\le T}$ ($\{j_{st}^{\rm M}\}_{0\le s\le t\le T}$, $\{j_{st}^{\rm AM}\}_{0\le s\le t\le T}$, respectively) on $\mathcal{B}$ defined by
\[
\begin{array}{rclcrclll}
j^{\rm B}_{st}(\mathbf{1}) &=& {\rm id}, && j^{\rm B}_{st}(b) &=& \tilde\jmath^{\rm B}_{0s}(p)\tilde\jmath^{\rm B}_{st}(b)\tilde\jmath^{\rm B}_{tT}(p) & \mbox{ for } &  b\in \mathcal{B}^0={\rm ker}\, \varepsilon, \\
j^{\rm M}_{st}(\mathbf{1}) &=& {\rm id}, && j^{\rm B}_{st}(b) &=&
\tilde\jmath^{\rm M}_{st}(b)\tilde\jmath^{\rm M}_{tT}(p) & \mbox{ for } &
b\in \mathcal{B}^0={\rm ker}\, \varepsilon, \\
j^{\rm AM}_{st}(\mathbf{1}) &=& {\rm id}, && j^{\rm AM}_{st}(b) &=& \tilde\jmath^{\rm AM}_{0s}(p)\tilde\jmath^{\rm AM}_{st}(b) & \mbox{ for } &  b\in \mathcal{B}^0={\rm ker}\, \varepsilon, \\
\end{array}
\]
for $0\le s\le t\le T$, is a boolean (monotone, anti-monotone, respectively) L\'evy process on the dual semigroup $(\mathcal{B},\Delta,\varepsilon)$. Furthermore, if $\{\tilde\jmath_{st}^{\rm B}\}_{0\le s\le t\le T}$ ($\{\tilde\jmath_{st}^{\rm M}\}_{0\le s\le t\le T}$, $\{\tilde\jmath_{st}^{\rm AM}\}_{0\le s\le t\le T}$, respectively) has generator $\hat\psi$, then $\{j_{st}^{\rm B}\}_{0\le s\le t\le T}$ ($\{j_{st}^{\rm M}\}_{0\le s\le t\le T}$, $\{j_{st}^{\rm AM}\}_{0\le s\le t\le T}$, respectively) is equivalent to $\{k_{st}^{\rm B}\}_{0\le s\le t\le T}$ ($\{k_{st}^{\rm M}\}_{0\le s\le t\le T}$, $\{k_{st}^{\rm AM}\}_{0\le s\le t\le T}$, respectively).
\end{theorem}
\begin{remark}
Every L\'evy process on an involutive bialgebra can be realized on boson Fock space as solution of quantum stochastic differential equations, see Theorem I.\ref{I-rep-thm} or \cite[Theorem 2.5.3]{schuermann93}. Therefore Theorem \ref{theo-real} implies that boolean, monotone, and anti-monotone L\'evy processes can also always be realized on a boson Fock space. We will refer to the realizations obtained in this way as standard Fock realization.

It is natural to conjecture that monotone and anti-monotone L\'evy processes can also be realized on their respective Fock spaces (see Subsection \ref{subsec-fock}) as solutions of monotone or anti-monotone quantum stochastic differential equations, like this has been proved for the tensor case in \cite[Theorem 2.5.3]{schuermann93} and discussed for free and boolean case in \cite{schuermann95b,benghorbal01}. We will show in Subsection \ref{sub-real-fock} that this is really possible.
\end{remark}
\begin{proof}
$\{\tilde\jmath^{\rm \bullet}_{st}\}_{0\le s\le t\le T}$ is a L\'evy process on the involutive bialgebra $(\tilde{\mathcal{B}},\overline{\Delta}_{\rm B},\tilde\varepsilon)$, $\bullet\in\{{\rm B},{\rm M},{\rm AM}\}$, and therefore, by the independence property of its increments, we have
\[
\big[\tilde\jmath^{\rm \bullet}_{st}(b_1),\tilde\jmath^{\rm \bullet}_{s't'}(b_2)\big] =0
\]
for all $0\le s\le t\le T$, $0\le s'\le t'\le T$ with $]s,t[\cap]s',t'[=\emptyset$ and all $b_1,b_2\in\tilde{\mathcal{B}}$. Using this property one immediately sees that the $j_{st}^\bullet$ are unital $*$-algebra homomorphisms. Using again the independence of the increments of $\{\tilde\jmath^{\rm \bullet}_{st}\}_{0\le s\le t\le T}$ and the fact that its marginal distributions $\varphi^\bullet_{st}=\Phi\circ\tilde\jmath^\bullet_{0s}$, $0\le s\le t\le T$, satisfy Equation (\ref{eq-p-invariance}), we get
\[
\Phi\big(j^{\rm B}_{st}(b)\big) = \Phi\big(\tilde\jmath^{\rm B}_{0s}(p)\tilde\jmath^{\rm B}_{st}(b)\tilde\jmath^{\rm B}_{tT}(p)\big)
=\Phi\big(\tilde\jmath^{\rm B}_{0s}(p)\big)\Phi\big(\tilde\jmath^{\rm B}_{st}(b)\big)\Phi\big(\tilde\jmath^{\rm B}_{tT}(p)\big)
= \varphi^{\rm B}_{st}(b)
\]
and similarly
\begin{eqnarray*}
\Phi\big(j^{\rm M}_{st}(b)\big) &=& \varphi^{\rm M}_{st}(b), \\
\Phi\big(j^{\rm AM}_{st}(b)\big) &=& \varphi^{\rm AM}_{st}(b),
\end{eqnarray*}
for all $b\in\mathcal{B}^0$. Thus marginal distributions of $\{j^\bullet_{st}\}_{0\le s\le t\le T}$ are simply the restrictions of the marginal distributions of $\{\tilde\jmath^{\rm \bullet}_{st}\}_{0\le s\le t\le T}$. This proves the stationarity and the weak continuity of $\{j^\bullet_{st}\}_{0\le s\le t\le T}$, it only remains to show the increment property and the independence of the increments. We check these for the boolean case, the other two cases are similar. Let $b\in\mathcal{B}^0$ with $\Delta(b)=\sum_{\epsilon\in\mathbb{A}} b^\epsilon$, where $b^\epsilon=b^\epsilon_1\otimes \cdots b^\epsilon_{\epsilon_{|\epsilon|}} \in \mathcal{B}_\epsilon=(\mathcal{B}^0)^{\otimes |\epsilon|}$, then we have
\begin{equation}\label{eq-comult-bool}
\overline{\Delta}_{\rm B}(b) =
\sum_{{\epsilon\in\mathbb{A}}\atop{\epsilon_1=1}}b^\epsilon_1pb^\epsilon_3 \cdots \otimes pb^\epsilon_2p\cdots  + \sum_{{\epsilon\in\mathbb{A}}\atop{\epsilon_1=2}}pb^\epsilon_2p \cdots \otimes b^\epsilon_1p b^\epsilon_3\cdots
\end{equation}
We set $j^{\rm B}_{st}=j_1$, $j^{\rm B}_{tu}=j_2$, and get
\begin{eqnarray*}
&&m_{\mathcal{A}}\circ\left(j^{\rm B}_{st}\coprod j^{\rm B}_{tu}\right)\circ \Delta(b)\\
 &=& \sum_{{\epsilon\in\mathbb{A}}\atop{\epsilon\not=\emptyset}} j_{\epsilon_1}(b^\epsilon_1)j_{\epsilon_2}(b^\epsilon_2)\cdots j_{\epsilon_{|\epsilon|}}(b^\epsilon_{|\epsilon|}) \\
&=& \sum_{{\epsilon\in\mathbb{A}}\atop{\epsilon_1=1}}\tilde\jmath^{\rm B}_{0s}(p)\tilde\jmath^{\rm B}_{st}(b^\epsilon_1)\tilde\jmath^{\rm B}_{tT}(p)\tilde\jmath^{\rm B}_{0t}(p)\tilde\jmath^{\rm B}_{tu}(b^\epsilon_2)\tilde\jmath^{\rm B}_{uT}(p)\cdots\tilde\jmath^{\rm B}_{0s}(p)\tilde\jmath^{\rm B}_{st}(b^\epsilon_{|\epsilon|})\tilde\jmath^{\rm B}_{tT}(p) \\
&&+ \sum_{{\epsilon\in\mathbb{A}}\atop{\epsilon_1=2}}\tilde\jmath^{\rm B}_{0t}(p)\tilde\jmath^{\rm B}_{tu}(b^\epsilon_1)\tilde\jmath^{\rm B}_{uT}(p)\tilde\jmath^{\rm B}_{0s}(p)\tilde\jmath^{\rm B}_{st}(b^\epsilon_2)\tilde\jmath^{\rm B}_{tT}(p)\cdots\tilde\jmath^{\rm B}_{0t}(p)\tilde\jmath^{\rm B}_{tu}(b^\epsilon_{|\epsilon|})\tilde\jmath^{\rm B}_{uT}(p)\\
&=& \tilde\jmath^{\rm B}_{0s}(p)\left(\sum_{{\epsilon\in\mathbb{A}}\atop{\epsilon_1=1}}\tilde\jmath^{\rm B}_{st}(b^\epsilon_1)\tilde\jmath^{\rm B}_{st}(p)\tilde\jmath^{\rm B}_{st}(b^\epsilon_3)\cdots\tilde\jmath^{\rm B}_{tu}(p)\tilde\jmath^{\rm B}_{tu}(b^\epsilon_2)\tilde\jmath^{\rm B}_{tu}(p)\cdots\right)\tilde\jmath^{\rm B}_{uT}(p) \\
&&+ \tilde\jmath^{\rm B}_{0s}(p)\left(\sum_{{\epsilon\in\mathbb{A}}\atop{\epsilon_1=2}}\tilde\jmath^{\rm B}_{st}(p)\tilde\jmath^{\rm B}_{st}(b^\epsilon_2)\tilde\jmath^{\rm B}_{st}(p)\cdots\tilde\jmath^{\rm B}_{tu}(b^\epsilon_1)\tilde\jmath^{\rm B}_{tu}(p)\tilde\jmath^{\rm B}_{tu}(b^\epsilon_3)\cdots\right)\tilde\jmath^{\rm B}_{uT}(p) \\
&=& \tilde\jmath^{\rm B}_{0s}(p)\left(m_{\mathcal{A}}\circ(\tilde\jmath^{\rm B}_{st}\otimes\tilde\jmath^{\rm B}_{tu})\circ\overline{\Delta}_{\rm B}(b)\right)\tilde\jmath^{\rm B}_{uT}(p) \\
&=&\tilde\jmath^{\rm B}_{0s}(p)\tilde\jmath^{\rm B}_{su}(b)\tilde\jmath^{\rm B}_{uT}(p) = j^{\rm B}_{su}(b).
\end{eqnarray*}
For the boolean independence of the increments of $\{j^{\rm B}_{st}\}_{0\le s\le t\le T}$, we have to check
\[
\Phi\circ m_{\mathcal{A}}\circ\left(j^{\rm B}_{s_1t_1}\coprod\cdots\coprod j^{\rm B}_{s_nt_n}\right) = \varphi^{\rm B}_{s_1t_1}|_\mathcal{B}\diamond\cdots \diamond\varphi^{\rm B}_{s_nt_n}|_{\mathcal{B}}
\]
for all $n\in\mathbb{N}$ and $0\le s_1\le t_1\le s_2\le \cdots\le t_n\le T$. Let, e.g., $n=2$, and take an element of $\mathcal{B}\coprod\mathcal{B}$ of the form $i_1(a_1)i_2(b_1)\cdots i_n(b_n)$, with $a_1,\ldots,a_n,b_1,\ldots,b_n\in\mathcal{B}^0$. Then we have
\begin{eqnarray*}
&&\Phi\circ m_{\mathcal{A}}\circ\left(j^{\rm B}_{s_1t_1}\coprod j^{\rm B}_{s_2t_2}\right)\big(i_1(a_1)i_2(b_1)\cdots i_n(b_n)\big) \\
 &=& \Phi\left( \tilde\jmath^{\rm B}_{0s_1}(p)\tilde\jmath^{\rm B}_{s_1t_1}(a_1)\tilde\jmath^{\rm B}_{t_1T}(p)\tilde\jmath^{\rm B}_{0s_2}(p)\tilde\jmath^{\rm B}_{s_2t_2}(b_1)\tilde\jmath^{\rm B}_{t_2T}(p) \cdots \tilde\jmath^{\rm B}_{0s_2}(p)\tilde\jmath^{\rm B}_{s_2t_2}(b_n)\tilde\jmath^{\rm B}_{t_2T}(p)\right) \\
&=& \Phi\left( \tilde\jmath^{\rm B}_{0s_1}(p)\tilde\jmath^{\rm B}_{s_1t_1}(a_1)\tilde\jmath^{\rm B}_{s_1t_1}(p)\cdots \tilde\jmath^{\rm B}_{s_1t_1}(a_n)\jmath^{\rm B}_{s_1t_2}(p)\tilde\jmath^{\rm B}_{s_2t_2}(b_1) \tilde\jmath^{\rm B}_{s_2t_2}(p)\cdots\tilde\jmath^{\rm B}_{s_2t_2}(b_n)\tilde\jmath^{\rm B}_{t_2T}(p)\right) \\
&=&\varphi^{\rm B}_{s_1t_1}(a_1pa_2p\cdots pa_n)\varphi^{\rm B}_{s_2t_2}(pb_1p\cdots pb_n) = \prod_{j=1}^n \varphi^{\rm B}_{s_1t_1}(a_j)\prod_{j=1}^n \varphi^{\rm B}_{s_2t_2}(b_j) \\
&=& \left(\varphi^{\rm B}_{s_1t_1}\diamond\varphi^{\rm B}_{s_2t_2}\right)\big(i_1(a_1)i_2(b_1)\cdots i_n(b_n)\big).
\end{eqnarray*}
The calculations for the other cases and general $n$ are similar.
\end{proof}

For the actual construction of $\{\tilde\jmath_{st}^{\rm B}\}_{0\le s\le t\le T}$ ($\{\tilde\jmath_{st}^{\rm M}\}_{0\le s\le t\le T}$, $\{\tilde\jmath_{st}^{\rm AM}\}_{0\le s\le t\le T}$, respectively) via quantum stochastic calculus, we need to know the Sch\"urmann triple of $\hat\psi$.

\begin{proposition}\label{prop-triple}
Let $\mathcal{B}$ be a unital $*$-algebra, $\psi:\mathcal{B}\to\mathbb{C}$ a generator, i.e.\ a hermitian, conditionally positive linear functional with $\psi(\mathbf{1})=0$, and $\hat\psi:\tilde{\mathcal{B}}\to\mathbb{C}$ the extension of $\psi$ to $\tilde{\mathcal{B}}$ given by Equation (\ref{eq-psi-p-inv}). If $(\rho,\eta,\psi)$ is a Sch\"urmann triple of $\psi$, then a Sch\"urmann triple $(\hat\rho,\hat\eta,\hat\psi)$ for $\hat\psi$ is given by
\[
\begin{array}{rclcrcl}
\hat\rho|_\mathcal{B} &=& \rho, && \hat\rho(p)=0, \\
\hat\eta|_\mathcal{B} &=& \eta, && \hat\eta(p)=0, \\
\hat\psi|_\mathcal{B} &=& \psi, && \hat\psi(p)=0,  
\end{array}
\]
in particular, it can be defined on the same pre-Hilbert space as  $(\rho,\eta,\psi)$.

\end{proposition}
\begin{proof}
The restrictions of $\hat\rho$ and $\hat\eta$ to $\mathcal{B}$ have to be unitarily equivalent to $\rho$ and $\eta$, respectively, since $\hat\psi|_\mathcal{B}=\psi$. We can calculate the norm of $\hat\eta(p)$ with Equation (I.\ref{I-L-coboundary}), we get
\[
\hat\psi(p)=\hat\psi(p^2)=\tilde\varepsilon(p)\hat\psi(p) + \langle\hat\eta(p^*),\hat\eta(p)\rangle + \hat\psi(p)\tilde\varepsilon(p)
\]
and therefore $||\hat\eta(p)||^2=-\hat\psi(p)=0$. From Equation (I.\ref{I-eta-cocycle}) follows
\[
\hat\eta(p^\alpha b_1 p b_2 p\cdots pb_np^\omega) =\left\{\begin{array}{ccl}
\eta(b_1) &\mbox{ if } & n=1, \alpha=0, \omega\in\{0,1\}, \\
0 &\mbox{ if } & n>1 \mbox{ or } \alpha=1.
\end{array}\right.
\]
For the representation $\hat\rho$ we get
\[
\hat\rho(p)\eta(b)=\hat\eta(pb)-\hat\eta(p)\varepsilon(b)=0
\]
for all $b\in\mathcal{B}$.
\end{proof}

The L\'evy processes $\{\tilde\jmath^\bullet_{st}\}_{0\le s\le t\le T}$ on the involutive bialgebras $(\tilde{\mathcal{B}},\overline{\Delta}_\bullet,\tilde\epsilon)$, $\bullet\in\{{\rm B},{\rm M},{\rm AM}\}$, with the generator $\hat\psi$ can now be constructed as solutions of the quantum stochastic differential equations
\[
\tilde\jmath^\bullet_{st}(b)=\tilde\varepsilon(b){\rm id}+\left( \int_s^t \tilde\jmath^\bullet_{s\tau}\otimes {\rm d}I_\tau\right)\overline{\Delta}_\bullet(b), \qquad\mbox{ for all } b\in\tilde{\mathcal{B}},
\]
where the integrator ${\rm d}I$ is given by
\[
{\rm d}I_t(b) = {\rm d}\Lambda_t(\hat\rho(b)-\tilde\varepsilon(b){\rm id}) + {\rm d}A_t^+(\hat\eta(b))+{\rm d}A_t(\hat\eta(b^*))+\hat\psi(b){\rm d}t.
\]

The element $p\in\tilde{\mathcal{B}}$ is group-like, i.e.\ $\overline{\Delta}_\bullet(p)=p\otimes p$, and mapped to zero by any Sch\"urmann triple $(\hat\rho,\hat\eta,\hat\psi)$ on $\tilde{\mathcal{B}}$ that is obtained by extending a Sch\"urmann triple $(\rho,\eta,\psi)$ on $\mathcal{B}$ as in Proposition \ref{prop-triple}. Therefore we can compute $\{\tilde\jmath^\bullet_{st}(p)\}_{0\le s\le t\le T}$ without specifying $\bullet\in\{{\rm B},{\rm M},{\rm AM}\}$ or knowing the Sch\"urmann triple $(\rho,\eta,\psi)$.

\begin{proposition}
Let $\{\tilde\jmath^\bullet_{st}\}_{0\le s\le t\le T}$ be a L\'evy process on $(\tilde{\mathcal{B}},\overline{\Delta}_\bullet,\tilde\epsilon)$, $\bullet\in\{{\rm B},{\rm M},{\rm AM}\}$, whose Sch\"urmann triple $(\hat\rho,\hat\eta,\hat\psi)$ is of the form given in Proposition \ref{prop-triple}. Denote by $0_{st}$ the projection from $L^2([0,T[,D)$ to $L^2([0,s[,D)\oplus L^2([t,T[,D)\subseteq L^2([0,T[,D)$,
\[
0_{st} f(\tau) =
\left\{\begin{array}{ll}
f(\tau) &\mbox{ if } \tau\not\in[s,t[, \\
0 & \mbox{ if } \tau\in[s,t[,
\end{array}\right.
\]
Then
\[
\tilde\jmath_{st}^\bullet(p)=\Gamma(0_{st}) \quad\mbox{ for all } 0\le s\le t\le T,
\]
i.e.\ $\tilde\jmath_{st}^\bullet(p)$ is equal to the second quantization of $0_{st}$ for all $0\le s\le t\le T$ and $\bullet\in\{{\rm B},{\rm M},{\rm AM}\}$.
\end{proposition}
\begin{proof}
This follows immediately from the quantum stochastic differential equation
\[
\tilde\jmath^\bullet_{st}(p) = {\rm id} - \int_s^t \tilde\jmath^\bullet_{s\tau}(p){\rm d}\Lambda_\tau({\rm id}).
\]
\end{proof}

\subsection{Boson Fock space realization of boolean, monotone, and anti-mono\-tone quantum stochastic calculus}\label{subsec-fock}

For each of the independences trea\-ted in this article, we can define a Fock space with a creation, annihilation and conservation process, and develop a quantum stochastic calculus. For the monotone case, this was done in \cite{muraki97a,lu97b}, for the boolean calculus see, e.g., \cite{benghorbal+dogan+schuermann01} and the references therein.

Since the integrator processes of these calculi have independent and stationary increments, we can use our previous results to realize them on a boson Fock space. Furthermore, we can embed the corresponding Fock spaces into a boson Fock space and thus reduce the boolean, monotone, and anti-monotone quantum stochastic calculus to the quantum stochastic calculus on boson Fock space defined in \cite{hudson+parthasarathy84} (but the integrands one obtains in the boolean or monotone case turn out to be not adapted in general). For the anti-monotone creation and annihilation process with one degree of freedom, this was already done in \cite{parthasarathy99} (see also \cite{liebscher99}).

Let $H$ be a Hilbert space. Its conjugate or dual is, as a set, equal to $\overline{H}=\{\overline{u}|u\in H\}$. The addition and scalar multiplication are defined by
\[
\overline{u}+\overline{v}=\overline{u+v},\qquad, z \overline{u}=\overline{\overline{z}u},\qquad \mbox{ for } u,v\in H,\quad z\in\mathbb{C}.
\]
Then $V(H)= H\otimes\overline{H}\oplus \overline{H}\oplus H$ (algebraic tensor product and direct sum, no completion) is an involutive complex vector space with the involution
\[
\left(v\otimes \overline{u}+ \overline{x}+y\right)^* = u\otimes\overline{v} + \overline{y}+x, \qquad \mbox{ for } u,v,x,y \in H.
\]
We will also write $|u\rangle\langle v|$ for $u\otimes \overline{v}$. Let now $\mathcal{B}_H$ be the free unital $*$-algebra over $V(H)$. This algebra can be made into a dual semigroup, if we define the comultiplication and counit by
\[
\Delta v = i_1(v)+ i_2(v),
\]
and $\varepsilon(v)=0$ for $v\in V(H)$ and extend them as unital $*$-algebra homomorphisms. On this dual semigroup we can define the fundamental noises for all our independences. For the Sch\"urmann triple we take the Hilbert space $H$, the representation $\rho$ of $\mathcal{B}_H$ on $H$ defined by
\[
\rho(u)=\rho(\overline{u})=0, \quad \rho\big(|u\rangle\langle v|\big): H\ni x\mapsto \langle v,x\rangle u\in H,
\]
the cocycle $\eta:\mathcal{B}_H\to H$ with
\[
\eta(u)=u,\quad \eta(\overline{u})=\eta\big(|u\rangle\langle v|\big) = 0,
\]
and the generator $\psi:\mathcal{B}_H\to\mathbb{C}$ with
\[
\psi(\mathbf{1})=\psi(u)=\psi(\overline{u})=\psi\big(|u\rangle\langle v|\big)=0,
\]
for all $u,v\in H$.

A realization of the tensor L\'evy process  $\{j_{st}\}_{0\le s\le t}$ on the dual semigroup $(\mathcal{B}_H,\Delta,\varepsilon)$ with this Sch\"urmann triple on the boson Fock space $\Gamma\big(L^2(\mathbb{R}_+,H)\big)$ is given by
\[
j_{st} (u) = A^+_{st}(u), \quad
j_{st} (\overline{u}) = A_{st}(u), \quad
j_{st} (|u\rangle\langle v|) = \Lambda_{st}\big(|u\rangle\langle v|\big),
\]
for all $0\le s\le t\le T$, $u,v\in H$.

\subsubsection{Boolean calculus}

Let $H$ be a Hilbert space. The boolean Fock space over $L^2([0,T[;H)\cong L^2([0,T])\otimes H$ is defined as $\Gamma_{\rm B}\big(L^2([0,T[,H)\big)=\mathbb{C}\oplus L^2([0,T[,H)$. We will write the elements of $\Gamma_{\rm B}\big(L^2([0,T[,H)\big)$ as vectors
\[
\left(\begin{array}{c} \lambda \\ f \end{array}\right)
\]
with $\lambda\in\mathbb{C}$ and $f\in L^2([0,T[,H)$. The boolean creation, annihilation, and conservation processes are defined as
\begin{eqnarray*}
A^{\rm B+}_{st}(u)\left(\begin{array}{c} \lambda \\ f \end{array}\right) &=&\left(\begin{array}{c} 0 \\ \lambda u\mathbf{1}_{[s,t[} \end{array}\right), \\
A^{\rm B}_{st}(u)\left(\begin{array}{c} \lambda \\ f \end{array}\right) &=& \left(\begin{array}{c} \int_s^t \langle u,f(\tau)\rangle{\rm d}\tau \\ 0 \end{array}\right), \\
\Lambda^{\rm B}_{st}\big(|u\rangle\langle v|\big)\left(\begin{array}{c} \lambda \\ f \end{array}\right) &=&\left(\begin{array}{c} 0 \\  \mathbf{1}_{[s,t[}(\cdot)\langle v,f(\cdot)\rangle u \end{array}\right),
\end{eqnarray*}
for $\lambda\in\mathbb{C}$, $f\in L^2([0,T[,H)$, $u,v\in H$. These operators define a boolean L\'evy process $\{k^{\rm B}_{st}\}_{0\le s\le t \le T}$ on the dual semigroup $(\mathcal{B}_H,\Delta,\varepsilon)$ with respect to the vacuum expectation, if we set
\[
k_{st}^{\rm B}(u) = A^{\rm B+}_{st}(u), \quad
k_{st}^{\rm B}(\overline{u}) = A^{\rm B}_{st}(u), \quad
k_{st}^{\rm B}\big(|u\rangle\langle v|\big) = \Lambda^{\rm B}_{st}\big(|u\rangle\langle v|\big),
\]
for all $0\le s\le t\le T$, $u,v\in H$, and extend the $k_{st}^{\rm B}$ as unital $*$-algebra homomorphisms to $\mathcal{B}_H$.

On the other hand, using Theorem \ref{theo-real} and Proposition \ref{prop-triple}, we can define a realization of the same L\'evy process on a boson Fock space. Since the comultiplication $\overline{\Delta}_{\rm B}$ acts on elements of the involutive bialgebra $(\tilde{\mathcal{B}}_H,\overline{\Delta}_{\rm B},\tilde\varepsilon)$ as
\[
\overline{\Delta}_{\rm B}(v)=v\otimes p+p\otimes v, \qquad\mbox{ for } v\in V(H),
\]
we have to solve the quantum stochastic differential equations
\begin{eqnarray*}
\tilde\jmath^{\rm B}_{st}(u) &=& \int_s^t \Gamma(0_{s\tau}){\rm d}A^+_\tau(u) - \int_s^t \tilde\jmath^{\rm B}_{s\tau}(u){\rm d}\Lambda_\tau({\rm id}_H), \\
\tilde\jmath^{\rm B}_{st}(\overline{u}) &=& \int_s^t \Gamma(0_{s\tau}){\rm d}A_\tau(u) - \int_s^t \tilde\jmath^{\rm B}_{s\tau}(\overline{u}){\rm d}\Lambda_\tau({\rm id}_H), \\
\tilde\jmath^{\rm B}_{st}\big(|u\rangle\langle v|\big) &=& \int_s^t \Gamma(0_{s\tau}){\rm d}\Lambda_\tau\big(|u\rangle\langle v|\big) - \int_s^t \tilde\jmath^{\rm B}_{s\tau}\big(|u\rangle\langle v|\big){\rm d}\Lambda_\tau({\rm id}_H),
\end{eqnarray*}
and set
\begin{eqnarray*}
j^{\rm B}_{st}(u) &=& \Gamma(0_{0s})\tilde\jmath^{\rm B}_{st}(u)\Gamma(0_{tT}), \\
j^{\rm B}_{st}(\overline{u}) &=& \Gamma(0_{0s})\tilde\jmath^{\rm B}_{st}(\overline{u})\Gamma(0_{tT}), \\
j^{\rm B}_{st}\big(|u\rangle\langle v|\big) &=& \Gamma(0_{0s})\tilde\jmath^{\rm B}_{st}\big(|u\rangle\langle v|\big)\Gamma(0_{tT}),
\end{eqnarray*}
These operators act on exponential vectors as
\begin{eqnarray*}
j^{\rm B}_{st}(u)\mathcal{E}(f) &=& u\mathbf{1}_{[s,t[}, \\
j^{\rm B}_{st}(\overline{u})\mathcal{E}(f) &=& \int_s^t \langle u,f(\tau)\rangle{\rm d}\tau \Omega, \\
j^{\rm B}_{st}\big(|u\rangle\langle v|\big)\mathcal{E}(f) &=& \mathbf{1}_{[s,t[}\langle v,f(\cdot)\rangle u,
\end{eqnarray*}
for $0\le s\le t\le T$, $f\in L^2([0,T[)$, $u,v\in H$.

Since $\{k^{\rm B}_{st}\}_{0\le s\le t \le T}$ and $\{j^{\rm B}_{st}\}_{0\le s\le t \le T}$ are boolean L\'evy processes on the dual semigroup $(\mathcal{B}_H,\Delta,\varepsilon)$ with the same generator, they are equivalent.

If we isometrically embed the boolean Fock space $\Gamma_{\rm B}\big(L^2([0,T[,H)\big)$ into the boson Fock space $\Gamma\big(L^2([0,T[,H)\big)$ in the natural way, $\theta_{\rm B}:\Gamma_{\rm B}\big(L^2([0,T[,H)\big)\to\Gamma\big(L^2([0,T[,H)\big)$,
\[
\theta_{\rm B}\left(\begin{array}{c} \lambda \\ f \end{array}\right) =\lambda\Omega+f
\]
for $\lambda\in\mathbb{C}$, $f\in L^2([0,T[,H)$, then we have
\[
k^{\rm B}_{st}(b)=\theta_{\rm B}^*j^{\rm B}_{st}(b)\theta_{\rm B}
\]
for all $b\in\mathcal{B}$.

\subsubsection{Anti-monotone calculus}\label{subsub-anti-mon}

We will treat the anti-monotone calculus first, because it leads to simpler quantum stochastic differential equations. The monotone calculus can then be constructed using time-reversal, cf.\ Lemma \ref{lem-mon-anti-mon}.

We can construct the monotone and the anti-monotone calculus on the same Fock space. Let
\[
\mathbb{T}_n=\{(t_1,\ldots,t_n)| 0\le t_1\le t_2\le\cdots\le t_n\le T\}\subseteq [0,T[^n\subseteq \mathbb{R}^n,
\]
then the monotone and anti-monotone Fock space $\Gamma_{\rm M}\big(L^2([0,T[,H)\big)$ over $L^2([0,T[,H)$  can be defined as
\[
\Gamma_{\rm M}\big(L^2([0,T[,H)\big) = \mathbb{C}\Omega\oplus \bigoplus_{n=1}^\infty L^2(\mathbb{T}_n, H^{\overline\otimes n}),
\]
where where $H^{\overline\otimes n}$ denotes the $n$-fold Hilbert space tensor product of $H$ and the measure on $\mathbb{T}_n$ is the restriction of the Lebesgue measure on $\mathbb{R}^n$ to $\mathbb{T}_n$. Since $\mathbb{T}_n\subseteq [0,T[^n$, we can interprete $f_1\otimes\cdots\otimes f_n\in L^2([0,T[,H)^{\overline\otimes n} \cong L^2([0,T[^n,H^{\overline\otimes n})$ also as an element of $L^2(\mathbb{T}_n,H^{\overline\otimes n})$ (by restriction).

The anti-monotone creation, annihilation, and conservation operator are defined by
\begin{eqnarray*}
A^{\rm AM+}_{st}(u)f_1\otimes\cdots\otimes f_n(t_1,\ldots,t_{n+1})&=&\mathbf{1}_{[s,t[}(t_1)u\otimes f_1\otimes\cdots\otimes f_n(t_2,\ldots,t_{n+1}) \\
A^{\rm AM}_{st}(u)f_1\otimes\cdots\otimes f_n(t_1,\ldots,t_{n-1})&=&\int^{\min(t,t_1)}_s \langle u, f_1(\tau)\rangle{\rm d}\tau f_2\otimes \cdots\otimes f_n (t_1,\ldots,t_{n-1}) \\
\Lambda^{\rm AM}_{st}\big(|u\rangle\langle v|\big)f_1\otimes\cdots\otimes f_n(t_1,\ldots,t_n)&=& \mathbf{1}_{[s,t[}(t_1)\langle v,f_1(t_1)\rangle u\otimes f_2\otimes\cdots\otimes f_n(t_2,\ldots,t_n),
\end{eqnarray*}
for $0\le s\le t\le T$, $0\le t_1\le t_2 \le \cdots \le t_n\le t_{n+1}\le T$, $u,v\in H$.

These operators define an anti-monotone L\'evy process $\{k^{\rm AM}_{st}\}_{0\le s\le t \le T}$ on the dual semigroup $\mathcal{B}$ with respect to the vacuum expectation, if we set
\[
k_{st}^{\rm AM}(u) = A^{\rm AM+}_{st}(u), \quad
k_{st}^{\rm AM}(\overline{u}) = A^{\rm AM}_{st}(u), \quad
k_{st}^{\rm AM}\big(|u\rangle\langle v|\big) = \Lambda^{\rm AM}_{st}\big(|u\rangle\langle v|\big),
\]
for all $0\le s\le t\le T$, $u,v\in H$, and extend the $k_{st}^{\rm AM}$ as unital $*$-algebra homomorphisms to $\mathcal{B}$.

We can define a realization of the same L\'evy process on a boson Fock space with Theorem \ref{theo-real}. The anti-monotone annihilation operators $j^{\rm AM}_{st}(\overline{u})$, $u\in H$, obtained this way act on exponential vectors as
\[
j^{\rm AM}_{st}(u)\mathcal{E}(f) = u\mathbf{1}_{[s,t[}(\cdot)\otimes_s\mathcal{E}(0_{0\cdot}f), \qquad f\in L^2([0,T[,H),
\]
and the anti-monotone creation operators are given by $j^{\rm AM}_{st}(u)=j^{\rm AM}_{st}(\overline{u})^*$, $u\in H$. On symmetric simple tensors $f_1\otimes\cdots\otimes f_n\in L^2([0,T[,H^{\overline{\otimes}n})$ they act as
\begin{eqnarray*}
& j^{\rm AM}_{st}(u)f_1\otimes\cdots\otimes f_n(t_1,\ldots,t_{n+1}) = & \\
& f(t_1)\otimes\cdots\otimes f_{k-1}(t_{k-1})\otimes u\mathbf{1}_{[s,t[}(t_k)\otimes f_{k+1}(t_{k+1})\otimes \cdots \otimes f_n(t_n) &
\end{eqnarray*}
where $k$ has to be chosen such that $t_k=\min\{t_1,\ldots,t_{n+1}\}$.

Since $\{k^{\rm AM}_{st}\}_{0\le s\le t \le T}$ and $\{j^{\rm AM}_{st}\}_{0\le s\le t \le T}$ are boolean L\'evy processes on the dual semigroup $\mathcal{B}$ with the same generator, they are equivalent.

A unitary map $\theta_{\rm M}:\Gamma_{\rm M}\big(L^2([0,T[,H)\big)\to \Gamma\big(L^2([0,T[,H)\big)$ can be defined by extending functions on $\mathbb{T}_n$ to symmetric functions on $[0,T[^n$ and dividing them by $\sqrt{n!}$. The adjoint $\theta_{\rm M}^*:\Gamma\big(L^2([0,T[,H)\big)\to \Gamma_{\rm M}\big(L^2([0,T[,H)\big)$ of $\theta_{\rm M}$ acts on simple tensors $f_1\otimes\cdots\otimes f_n\in L^2([0,T[,H)^{\overline\otimes n}\cong L^2([0,T[^n,H^{\overline\otimes n})$ as restriction to $\mathbb{T}_n$ and multiplication by $\sqrt{n!}$, i.e.\
\[
\theta_{\rm M}^*f_1\otimes\cdots\otimes f_n(t_1,\ldots,t_n)=\sqrt{n!} f_1(t_1)\otimes \cdots \otimes f_n(t_n),
\]
for all $f_1,\ldots,f_n\in L^2([0,T[,H)$, $(t_1,\ldots,t_n)\in\mathbb{T}_n$.

This isomorphism intertwines between $\{k^{\rm AM}_{st}\}_{0\le s\le t \le T}$ and $\{j^{\rm AM}_{st}\}_{0\le s\le t \le T}$, we have
\[
k^{\rm AM}_{st}(b) = \theta_{\rm M}^*j^{\rm AM}_{st}(b)\theta_{\rm M}
\]
for all $0\le s\le t\le T$ and $b\in\mathcal{B}_H$.

\subsubsection{Monotone calculus}\label{subsub-mon}

The monotone creation, annihilation, and conservation operator on the monotone Fock space $\Gamma_{\rm M}\big(L^2([0,T[,H)\big)$ can be defined by
\begin{eqnarray*}
A^{\rm M+}_{st}(u)f_1\otimes\cdots\otimes f_n(t_1,\ldots,t_{n+1})&=&f_1\otimes\cdots\otimes f_n(t_1,\ldots,t_n)\otimes \mathbf{1}_{[s,t[}(t_{n+1})u \\
A^{\rm AM}_{st}(u)f_1\otimes\cdots\otimes f_n(t_1,\ldots,t_{n-1})&=&\int_{\max(s,t_{n-1})}^t \langle u, f_n(\tau)\rangle{\rm d}\tau f_1\otimes \cdots\otimes f_{n-1} (t_1,\ldots,t_{n-1}) \\
\Lambda^{\rm AM}_{st}\big(|u\rangle\langle v|\big)f_1\otimes\cdots\otimes f_n(t_1,\ldots,t_n)&=& f_1\otimes\cdots\otimes f_{n-1}(t_1,\ldots,t_{n-1})\mathbf{1}_{[s,t[}(t_n)\langle v,f_n(t_n)\rangle u,
\end{eqnarray*}
for $0\le s\le t\le T$, $u,v\in H$. These operators define a monotone L\'evy process $\{k^{\rm M}_{st}\}_{0\le s\le t \le T}$ on the dual semigroup $\mathcal{B}$ with respect to the vacuum expectation, if we set
\[
k_{st}^{\rm M}(u) = A^{\rm M+}_{st}(u), \quad
k_{st}^{\rm M}(\overline{u}) = A^{\rm M}_{st}(u), \quad
k_{st}^{\rm M}\big(|u\rangle\langle v|\big) = \Lambda^{\rm M}_{st}\big(|u\rangle\langle v|\big),
\]
for all $0\le s\le t\le T$, $u,v\in H$, and extend the $k_{st}^{\rm M}$ as unital $*$-algebra homomorphisms to $\mathcal{B}$.

Define a time-reversal $R:\Gamma_M\big(L^2([0,T[,H)\big)\to\Gamma_M\big(L^2([0,T[,H)\big)$ for the monotone Fock space by $R\Omega=\Omega$ and
\[
Rf_1\otimes\cdots\otimes f_n(t_1,\ldots,t_n)=f_n(T-t_n)\otimes\cdots\otimes f_1(T-t_1),
\]
for $(t_1,\ldots,t_n)\in\mathbb{T}_n$, $f_,\ldots,f_n\in L^2(\mathbb{T}_n)$. The time-reversal $R$ is unitary and satisfies $R^2={\rm id}_{\Gamma_{\rm M}(L^2([0,T[;H))}$. It intertwines between the monotone and anti-monotone noise on the monotone Fock space, i.e.\ we have
\[
k^{\rm AM}_{st}(b)=Rk^{\rm M}_{T-t,T-s}(b)R
\]
for all $0\le s\le t\le T$, $b\in\mathcal{B}_H$.
On the boson Fock space we have to consider $R_{\rm M}=\theta_{\rm M}R\theta^*_{\rm M}:\Gamma\big(L^2([0,T[,H)\big)\to \Gamma\big(L^2([0,T[,H)\big)$. This map is again unitary and satisfies also $R_{\rm M}^2={\rm id}$. It follows that the realization $\{j^{\rm M}_{st}\}_{0\le s\le t\le T}$ of $\{k^{\rm M}_{st}\}_{0\le s\le t\le T}$ on boson Fock space can be defined via
\begin{eqnarray*}
j^{\rm M}_{st}(u) &=& \int_s^t {\rm d}\tilde{A}^+_\tau(u) \Gamma(0_{\tau T}),\\
j^{\rm M}_{st}(\overline{u}) &=& \int_s^t {\rm d}\tilde{A}_\tau(u) \Gamma(0_{\tau T}), \\
j^{\rm M}_{st}\big(|u\rangle\langle v|\big) &=&  \int_s^t {\rm d}\tilde{\Lambda}_\tau\big(|u\rangle\langle v|\big) \Gamma(0_{\tau T}),
\end{eqnarray*}
where the integrals are {\em backward} quantum stochastic integrals.

\begin{remark}\label{rem-azema}
Taking $H=\mathbb{C}$ and comparing these equations with \cite[Section 4.3]{schuermann93}, one recognizes that our realization of the monotone creation and annihilation process on the boson Fock space can be written as
\begin{eqnarray*}
\theta_{\rm M}A^{\rm M+}_{st}(1)\theta^*_{\rm M}&=& j^{\rm M}_{st}(1) = X^*_{st}\Gamma(0_{tT}), \\
\theta_{\rm M}A^{\rm M}_{st}(1)\theta^*_{\rm M} &=& j^{\rm M}_{st}(\overline{1}) = X_{st}\Gamma(0_{tT}),
\end{eqnarray*}
where $\{(X^*_{st},X_{st})\}_{0\le s\le t\le T}$ is the quantum Az\'ema martingale \cite{parthasarathy90,schuermann91b} with parameter $q=0$, cf.\ Subsection I.\ref{I-sub-azema}. Note that here $1$ denotes the unit of $H=\mathbb{C}$, not the unit of $\mathcal{B}_\mathbb{C}$.
\end{remark}

\subsection{Markov structure of boolean, monotone, and anti-monotone L\'evy processes}
\label{subsec-markov}
By a conditional expectation $\mathbb{E}$ on a quantum probability space $(\mathcal{A},\Phi)$ we will mean a completely positive map $\mathbb{E}:\mathcal{A}\to\mathcal{A}$ that satisfies
\begin{eqnarray*}
\Phi\circ\mathbb{E}&=&\Phi, \\
\mathbb{E}\circ\mathbb{E} &=& \mathbb{E}, \\
\mathbb{E}(xyz)&=&x\mathbb{E}(y)z, \qquad \mbox{ for all }x,z\in\mathbb{E}(\mathcal{A}),y\in\mathcal{A}.
\end{eqnarray*}
Note that we do {\em not} require that $\mathbb{E}$ maps the unit of $\mathcal{A}$ to itself. This will only be the case for the family of conditional expectations that we will use for L\'evy processes on involutive bialgebras and anti-monotone L\'evy processes, but not for boolean or monotone L\'evy processes on dual semi-groups.

We call a process $\{j_t:\mathcal{B}\to(\mathcal{A},\Phi)\}_{0\le t\le T}$ markovian (w.r.t.\ $\{\mathbb{E}_t\}_{0\le t\le T}$), if there exists a family of conditional expectations $\{\mathbb{E}_t:\mathcal{A}\to\mathcal{A}\}_{0\le t\le T}$ such that
\[
\mathbb{E}_s\big(j_t(\mathcal{B})\big) \subseteq j_s(\mathcal{B}),
\]
for all $0\le s\le t\le T$. A semi-group $\{p_t:\mathcal{B}\to\mathcal{B}\}_{t\ge 0}$ is called a markovian semi-group of $\{j_t:\mathcal{B}\to(\mathcal{A},\Phi)\}_{0\le t\le T}$ (w.r.t.\ $\{\mathbb{E}_t\}_{0\le t\le T}$), if
\[
\mathbb{E}_s\big(j_t(b)\big) = j_s\big(p_{t-s}(b)\big),
\]
for all $0\le s\le t\le T$ and $b\in\mathcal{B}$.

A L\'evy process on an involutive bialgebra $(\mathcal{B},\Delta,\varepsilon)$ has a natural Markov structure, cf.\ \cite{franz99}. There exists a semi-group of completely positive maps $\{P_t:\mathcal{B}\to\mathcal{B}\}_{t\ge 0}$ such that
\[
\mathbb{E}_s\big(j_{0t}(b)\big) = j_{0s}\big(P_{t-s}(b)\big)
\]
for all $0\le s\le t\le T$ and $b\in\mathcal{B}$. In \cite{franz99}, we took the inductive limit realization and $\{\mathbb{E}_t\}_{0\le t\le T}$ was a family of conditional expectations constructed on that realization. But the same holds for the realization on the boson Fock space obtained by solving the quantum stochastic differential equation (\ref{levy-intro}.\ref{I-rep-thm-qsde}) with respect to the standard conditional expectations $\{\mathbb{E}^1_t\}_{0\le t\le T}$ of the boson Fock space $\Gamma\big(L^2([0,T[,H)\big)$. The conditional expectation $\mathbb{E}^1_t(X)$, $0\le t\le T$, of a (possibly unbounded) operator on $\Gamma\big(L^2([0,T[,H)\big)$ is determined by the condition
\begin{equation}\label{eq-cond-exp-1}
\langle u,\mathbb{E}^1_t(X) v\rangle = \langle u_{t]}\otimes \Omega_{[t},X v_{t]}\otimes \Omega_{[t}\rangle\langle u_{[t},v_{[t}\rangle
\end{equation}
for vectors $u=u_{t]}\otimes u_{[t},v=v_{t]}\otimes v_{[t}\in \Gamma\big(L^2([0,t[,H)\big)\otimes \Gamma\big(L^2([t,T[,H)\big)\cong  \Gamma\big(L^2([0,T[,H)\big)$ such that $v_{t]}\otimes \Omega_{[t}$ is in the domain of $X$. It satisfies $\mathbb{E}^1_t({\rm id})={\rm id}$ for all $0\le t\le T$. The markovian semi-group $\{P_t:\mathcal{B}\to\mathcal{B}\}_{t\ge 0}$ is defined by
\[
P_t=({\rm id}_\mathcal{B}\otimes \varphi_t)\circ \Delta.
\]
Since we have $\varphi_0=\varepsilon$ and $\varphi_s\star\varphi_t=\varphi_{s+t}$ for all $s,t\ge 0$, is follows that $\{P_t\}_{t\ge 0}$ is a semi-group, i.e.\
\[
P_0={\rm id}_{\mathcal{B}}, \quad \mbox{ and } \quad P_s\circ P_t=P_{s+t}, \quad \mbox{ for all } s,t\ge 0.
\]
Due to the fact that $\Delta$ is a unital $*$-algebra homomorphism and that the $\varphi_t$, $t\ge 0$ are states, the $P_t$, $t\ge 0$, are completely positive and preserve the unit of $\mathcal{B}$.

Let $\{k_{st}\}_{0\le s\le t\le T}$ be a boolean (monotone, anti-monotone, respectively) L\'evy process on a dual semi-group $(\mathcal{B},\Delta,\varepsilon)$. The L\'evy process $\{\tilde\jmath^{\rm B}_{st}\}_{0\le s\le t\le T}$ ($\{\tilde\jmath^{\rm M}_{st}\}_{0\le s\le t\le T}$, $\{\tilde\jmath^{\rm AM}_{st}\}_{0\le s\le t\le T}$, respectively) on the involutive bialgebra $(\tilde{\mathcal{B}},\overline{\Delta}_{\rm B},\tilde\varepsilon)$  ($(\tilde{\mathcal{B}},\overline{\Delta}_{\rm M},\tilde\varepsilon)$, $(\tilde{\mathcal{B}},\overline{\Delta}_{\rm AM},\tilde\varepsilon)$ respectively) that we associate to it, is markovian. Its markovian semi-group is given by $\tilde{P}^{\bullet}_t=({\rm id}_{\tilde{\mathcal{B}}}\otimes \tilde\varphi_t)\circ \overline{\Delta}_{\bullet}$, $\bullet\in\{{\rm B},{\rm M},{\rm AM}\}$. We will now check in what sense this gives a markovian semi-group for the boolean (monotone, anti-monotone, respectively) L\'evy process on the  dual semi-group $(\mathcal{B},\Delta,\varepsilon)$.

In the following discussion of the Markov property of boolean, monotone, and anti-monotone L\'evy processes on dual semi-groups, we will assume that all our processes are realized on a boson Fock space as solutions of quantum stochastic differential equations in the way given by Sch\"urmann's representation theorem \cite[Theorem 2.5.3]{schuermann93} and Theorem \ref{theo-real}, i.e.\ as their standard Fock realization.

\subsubsection{Markov structure of boolean L\'evy processes}

Here we take the family of (non-unital!) conditional expectations $\{\mathbb{E}^0_t\}_{0\le t\le T}$ defined by
\[
\mathbb{E}^0_t(X) = \Gamma(0_{tT})X\Gamma(0_{tT}), \qquad \mbox{ for } 0\le t\le T,
\]
for a (possible unbounded) operator $X$ on $\Gamma\big(L^2([0,T[,H)\big)$, where
\[
\Gamma(0_{tT}):\Gamma\big(L^2([0,T[,H)\big)\to \Gamma\big(L^2([0,T[,H)\big)
\]
is the second quantization of the projection $0_{tT}$ of $L^2([0,T[,H)$ onto $L^2([0,t[,H)$.

Note that these conditional expectations do not preserve the identity, instead we have
\[
\mathbb{E}^0_t({\rm id}) = \Gamma(0_{tT}), \qquad \mbox{ for } 0\le t\le T.
\]
Therefore we have to modify the definition of our monotone L\'evy processes on the unit element.

\begin{theorem}\label{theo-markov-bool}
Let $\{j^{\rm B}_{st}\}_{0\le s\le t\le T}$ be the standard Fock realization of a boolean L\'evy process on a dual semi-group $(\mathcal{B},\Delta,\varepsilon)$. Define $\{\tilde{k}^{\rm B}_{t}\}_{0\le t\le T}$ on $\tilde{\mathcal{B}}$ by
\begin{eqnarray*}
\tilde{k}^{\rm B}_{t}(\mathbf{1})&=&\Gamma(0_{tT}), \\
\tilde{k}^{\rm B}_{t}(p) &=& \Gamma(0_{0T}), \\
\tilde{k}^{\rm B}_{t}(b) &=& j_{0t}^{\rm B}(b),\quad \mbox{ for
  }b\in\mathcal{B}^0={\rm ker}\,\varepsilon,
\end{eqnarray*}
on $\mathcal{B}$ and $\mathbb{C}M$ and extend as (non-unital!) $*$-algebra homomorphisms.

Then $\{\tilde{k}^{\rm B}_{t}\}_{0\le t\le T}$ is a Markov process on $\tilde{\mathcal{B}}$ with respect to the family of conditional expectations $\{\mathbb{E}^0_t\}_{0\le t\le T}$. Furthermore, $\{\tilde{P}^{\rm B}_t=({\rm id}_{\tilde{\mathcal{B}}}\otimes \tilde\varphi_t)\circ \overline{\Delta}_{\rm B}\}_{0\le t\le T}$ is a markovian semi-group for it, i.e.\ we have
\[
\mathbb{E}^0_s\big(\tilde{k}^{\rm B}_{t}(b)\big) = \tilde{k}^{\rm B}_{s}\big(\tilde{P}^{\rm B}_{t-s}(b)\big)
\]
for all $0\le s\le t\le T$ and $b\in \tilde{\mathcal{B}}$.
\end{theorem}
\begin{proof}
Let $\{\tilde\jmath^{\rm B}_{st}\}_{0\le s\le t\le T}$ be the L\'evy process on the involutive bialgebra $(\tilde{\mathcal{B}},\overline{\Delta}_{\rm B},\tilde\varepsilon)$, from which $\{j^{\rm B}_{st}\}_{0\le s\le t\le T}$ is constructed as in Theorem \ref{theo-real}.

For general elements $b\in\tilde{\mathcal{B}}$ we can define $\{\tilde{k}^{\rm B}_{t}\}_{0\le t\le T}$ also by
\[
\tilde{k}^{\rm B}_{t}(b)=\tilde\jmath^{\rm B}_{0t}(b)\Gamma(0_{tT}), \quad 0\le t\le T.
\]

Let $b\in\mathcal{B}^0$, then $\overline{\Delta}_{\rm B}(b)$ has the form given in Equation (\ref{eq-comult-bool}), and, using the increment property of $\{\tilde\jmath^{\rm B}_{st}\}_{0\le s\le t\le T}$, we get
\begin{eqnarray*}
\mathbb{E}^0_s\big(\tilde{k}^{\rm B}_{t}(b)\big) &=& \mathbb{E}^0_s\big(\tilde\jmath^{\rm B}_{0t}(b)\Gamma(0_{tT})\big)
=
\mathbb{E}^0_s\Big(\big(m_\mathcal{A}\circ\left(\tilde\jmath^{\rm B}_{0s}\otimes\tilde\jmath^{\rm B}_{st}\right)\circ\overline{\Delta}_{\rm B}(b)\big)\Gamma(0_{tT})\Big) \\
&=& \Gamma(0_{sT})\Big(\sum_{{\epsilon\in\mathbb{A}}\atop{\epsilon_1=1}}\tilde\jmath^{\rm B}_{0s}(b^\epsilon_1pb^\epsilon_3\cdots)\tilde\jmath^{\rm B}_{st}(pb^\epsilon_2p\cdots) \Gamma(0_{tT}) \\
&&+ \sum_{{\epsilon\in\mathbb{A}}\atop{\epsilon_1=2}}\tilde\jmath^{\rm B}_{0s}(pb^\epsilon_2p\cdots)\tilde\jmath^{\rm B}_{st}(b^\epsilon_1pb^\epsilon_3\cdots)\Gamma(0_{tT})\Big)\Gamma(0_{sT})\\
&=& \sum_{{\epsilon\in\mathbb{A}}\atop{\epsilon_1=1}}\tilde\jmath^{\rm B}_{0s}(b^\epsilon_1pb^\epsilon_3\cdots)\Gamma(0_{st})\tilde\jmath^{\rm B}_{st}(pb^\epsilon_2p\cdots) \Gamma(0_{st})\Gamma(0_{tT}) \\
&&+\sum_{{\epsilon\in\mathbb{A}}\atop{\epsilon_1=2}}\tilde\jmath^{\rm B}_{0s}(pb^\epsilon_2p\cdots)\Gamma(0_{st})\tilde\jmath^{\rm B}_{st}(b^\epsilon_1pb^\epsilon_3\cdots)\Gamma(0_{st})\Gamma(0_{tT})\\
&=& \sum_{{\epsilon\in\mathbb{A}}\atop{\epsilon_1=1}}\tilde\jmath^{\rm B}_{0s}(b^\epsilon_1pb^\epsilon_3\cdots)\tilde\varphi_{st}(pb^\epsilon_2p\cdots)\Gamma(0_{st})\Gamma(0_{tT}) \\
&&+\sum_{{\epsilon\in\mathbb{A}}\atop{\epsilon_1=2}}\tilde\jmath^{\rm B}_{0s}(pb^\epsilon_2p\cdots)\tilde\varphi_{st}(b^\epsilon_1pb^\epsilon_2\cdots)\Gamma(0_{st})\Gamma(0_{tT}),
\end{eqnarray*}
since
\[
\Gamma(0_{st})\tilde\jmath^{\rm B}_{st}(p^\alpha b_1p\cdots pb_np^\omega) \Gamma(0_{st}) = \tilde\jmath^{\rm B}_{st}(p b_1p\cdots pb_np) = \tilde\varphi_{st}(p^\alpha b_1p\cdots pb_np^\omega)\Gamma(0_{st})
\]
for all $n\in\mathbb{N}$, $b_1,\ldots,b_n\in\mathcal{B}^0$, $\alpha,\omega\in\{0,1\}$. But this is nothing else than
\begin{eqnarray*}
\mathbb{E}^0_s\big(\tilde{k}^{\rm B}_{0t}(b)\big) &=& \tilde\jmath^{\rm B}_{0s}\Big(\sum_{{\epsilon\in\mathbb{A}}\atop{\epsilon_1=1}}b^\epsilon_1pb^\epsilon_3\cdots\tilde\varphi_{st}(pb^\epsilon_2p\cdots)+\sum_{{\epsilon\in\mathbb{A}}\atop{\epsilon_1=2}}pb^\epsilon_2p\cdots\tilde\varphi_{st}(b^\epsilon_1pb^\epsilon_2\cdots)\Big)\Gamma(0_{sT}) \\
&=&\tilde\jmath^{\rm B}_{0s}\big(({\rm id}_{\tilde{\mathcal{B}}}\otimes\tilde\varphi_{st})\circ\overline{\Delta}_{\rm B}(b)\big)\Gamma(0_{sT}) = j^{\rm B}_{0s}\big(\tilde{P}^{\rm B}_{t-s}(b)\big)\Gamma(0_{sT}) \\
&=&\tilde{k}^{\rm B}_s\big(\tilde{P}^{\rm B}_{t-s}(b)\big)
\end{eqnarray*}
This implies that the Markov property also holds for all elements of $\tilde{\mathcal{B}}$.
\end{proof}

\subsubsection{Markov structure of monotone L\'evy processes}

For monotone L\'evy processes on dual semi-groups we have to use the same family of non-unital conditional expectations as for the boolean case.

\begin{theorem}
Let $\{j^{\rm M}_{st}\}_{0\le s\le t\le T}$ be the standard Fock realization of a  monotone L\'evy process on a dual semi-group $(\mathcal{B},\Delta,\varepsilon)$. Define $\{\tilde{k}^{\rm B}_{t}\}_{0\le t\le T}$ on $\tilde{\mathcal{B}}$ by
\begin{eqnarray*}
\tilde{k}^{\rm M}_{t}(\mathbf{1})&=&\Gamma(0_{tT}), \\
\tilde{k}^{\rm M}_{t}(p) &=& \Gamma(0_{0T}), \\
\tilde{k}^{\rm M}_{t}(b) &=& j^{\rm B}(b),\quad \mbox{ for
  }b\in\mathcal{B}^0={\rm ker}\,\varepsilon,
\end{eqnarray*}
on $\mathcal{B}$ and $\mathbb{C}M$ and extend as (non-unital!) $*$-algebra homomorphisms.

Then $\{\tilde{k}^{\rm M}_{t}\}_{0\le t\le T}$ is a Markov process on $\tilde{\mathcal{B}}$ w.r.t.\ $\{\mathbb{E}^0_t\}_{0\le t\le T}$, with markovian semi-group $\tilde{P}^{\rm M}_t=({\rm id}\otimes \tilde\varphi^{\rm M}_t)\circ\overline{\Delta}_{\rm M}$, i.e.\ we have
\[
\mathbb{E}^{\rm 0}_s\big(k^{\rm M}_{t}(b)\big) = k^{\rm M}_{s}\big(\tilde{P}^{\rm M}_{t-s}(b)\big)
\]
for all $0\le s\le t\le T$ and $b\in \tilde{\mathcal{B}}$.
\end{theorem}
\begin{proof}
Similar to the proof of Theorem \ref{theo-markov-bool}.
\end{proof}

The comultiplication of the involutive bialgebra $(\tilde{\mathcal{B}},\overline{\Delta}_{\rm M},\tilde\varepsilon)$ is given by
\[
\overline{\Delta}_{\rm M} = \sum_{{\epsilon\in\mathbb{A}}\atop{\epsilon_1=1}} b^\epsilon_1b_3^\epsilon\cdots\otimes p b^\epsilon_2p\cdots + \sum_{{\epsilon\in\mathbb{A}}\atop{\epsilon_1=2}} b^\epsilon_2b_4^\epsilon\cdots\otimes b^\epsilon_1pb^\epsilon_3p\cdots
\]
for an element $b\in\mathcal{B}^0$ with $\Delta(b)=\sum_{\epsilon\in\mathbb{A}} b^\epsilon$, $b^\epsilon\in \mathcal{B}_\epsilon=\left(\mathcal{B}^0\right)^{\otimes |\epsilon|}$, for the comultiplication of the dual semi-group $(\mathcal{B},\Delta,\varepsilon)$. The fact that there are no $p$'s on the left-hand-side implies that $\mathcal{B}$ is a right coideal in $(\tilde{\mathcal{B}},\overline{\Delta}_{\rm M},\tilde\varepsilon)$, i.e.\
\[
\overline{\Delta}_{\rm M}(\mathcal{B})\subseteq \mathcal{B}\otimes\tilde{\mathcal{B}}.
\]
This implies that the restriction of $\{\tilde\jmath^{\rm M}_{st}\}_{0\le s\le t\le T}$ to $\mathcal{B}$ is again a Markov process, cf.\ \cite[Theorem 3.2]{franz99}. The same also holds for $\{\tilde{k}_t^{\rm M}\}_{0\le s\le t\le T}$ restricted to $\mathcal{B}$, i.e.\ the non-unital version $\{\tilde{k}^{\rm M}_t|_\mathcal{B}\}_{0\le s\le t\le T}$ of a monotone L\'evy process on a dual semi-group defined in the previous theorem is a Markov process on $\mathcal{B}$ itself, not only on its $p$-extension.

\begin{corollary}\label{cor-mon-markov}
Define $\{k^{\rm M}_t\}_{0\le t\le T}$ by $k^{\rm M}_t=\tilde{k}^{\rm M}_t|_\mathcal{B}$ and $\{P^{\rm M}_t\}_{0\le t\le T}$ by $P^{\rm M}_t=\tilde{P}^{\rm M}_t|_\mathcal{B}$ for $0\le t\le T$. Then $\{k^{\rm M}_t\}_{0\le t\le T}$ is a Markov process on $\mathcal{B}$ w.r.t. $\{\mathbb{E}^0_t\}_{0\le t\le T}$, with markovian semi-group $\{P^{\rm M}_t\}_{0\le t\le T}$, i.e.\
\[
\mathbb{E}^0_s\big(k^{\rm M}_t(b)\big) = k^{\rm M}_s\big(P^{\rm M}_{t-s}(b)\big)
\]
for all $0\le s\le t\le T$, $b\in\mathcal{B}$.
\end{corollary}

\subsubsection{Markov structure of anti-monotone L\'evy processes}

For anti-monotone L\'evy processes on dual semi-groups we use the same family of unital conditional expectations $\{\mathbb{E}^1_t\}_{0\le t\le T}$ as for L\'evy processes on involutive bialgebras.

\begin{theorem}
Let $\{j^{\rm AM}_{st}\}_{0\le s\le t\le T}$ be the standard Fock realization of an anti-monotone L\'evy process on a dual semi-group $(\mathcal{B},\Delta,\varepsilon)$. Denote by $\{\tilde{k}^{\rm AM}_{t}\}_{0\le t\le T}$ the extension of $\{j^{\rm AM}_{0t}\}_{0\le t\le T}$ to $\tilde{\mathcal{B}}$ as unital $*$-homomorphisms with
\[
\tilde{k}^{\rm AM}_t(p)=\Gamma(0_{0t}), \qquad \mbox{ for all } 0\le t\le T,
\]
i.e.\ $\tilde{k}^{\rm AM}_t=\tilde\jmath^{\rm AM}_{0t}$.

Then $\{\tilde{k}^{\rm AM}_t\}_{0\le t\le T}$ is a Markov process on $\tilde{\mathcal{B}}$ w.r.t.\ $\{\mathbb{E}^1_t\}_{0\le t\le T}$, with markovian semi-group $\tilde{P}^{\rm AM}_t$, i.e.\ we have
\[
\mathbb{E}^1_s\big(\tilde{k}^{\rm B}_t(b)\big) = \tilde{k}^{\rm AM}_s\big(\tilde{P}^{\rm AM}_{t-s}(b)\big)
\]
for all $0\le s\le t\le T$ and $b\in \tilde{\mathcal{B}}$.
\end{theorem}
\begin{proof}
This is simply \cite[Theorem 3.1]{franz99} and \cite[Corollary 3.1]{franz99}, because $\{\tilde{k}^{\rm AM}_{t}\}_{0\le t\le T}$ is obtained from the L\'evy process $\{\tilde\jmath^{\rm AM}_{st}\}_{0\le s\le t\le T}$ on the involutive bialgebra $(\tilde{\mathcal{B}},\overline{\Delta}_{\rm AM},\tilde\varepsilon)$ by fixing $s=0$.
\end{proof}

\subsubsection{Realization of boolean, monotone, and anti-monotone L\'evy process on boolean, monotone, and anti-monotone Fock spaces} \label{sub-real-fock}

Free and boolean L\'evy processes on dual semigroups can be realized as solutions of free or boolean quantum stochastic equations on the free or boolean Fock space, see e.g.\ \cite{schuermann95b}. A full proof of this fact is still missing, because it would require a generalization of their calculi to unbounded coefficients, but for a large class of examples this has been shown in \cite[Section 6.5]{benghorbal01} for the boolean case. For dual semigroups that are generated by primitive elements (i.e.\ $\Delta(v)=i_1(v)+1_2(v)$) it is sufficient to determine the operators $j_{0t}(v)$, which have additive free or boolean increments. It turns out that they can always be represented as a linear combination of the corresponding creators, annihilators, conservation operators and time (which contains the projection $\Gamma(0_{0T})$ to the vacuum in the boolean case), cf.\ \cite{glockner+schuermann+speicher92,benghorbal01}.

We will sketch, how one can show the same for monotone and anti-monotone L\'evy processes on dual semigroups.

We can write the fundamental integrators of the anti-monotone calculus on the monotone Fock space $\Gamma_{\rm M}\big(L^2([0,t[,H)\big)$ as
\begin{eqnarray*}
{\rm d}A^{\rm AM+}_t(u) &=& \theta_{\rm M}^*\Gamma(0_{0t}){\rm d}A^+_t(u)\theta_{\rm M}, \\
{\rm d}A^{\rm AM}_t(u) &=& \theta_{\rm M}^*\Gamma(0_{0t}){\rm d}A_t(u)\theta_{\rm M}, \\
{\rm d}\Lambda^{\rm AM}_t\big(|u\rangle\langle v|\big) &=& \theta_{\rm M}^*\Gamma(0_{0t}){\rm d}\Lambda_t\big(|u\rangle\langle v|\big)\theta_{\rm M},
\end{eqnarray*}
where $\theta_{\rm M}:\Gamma_{\rm M}\big(L^2([0,t[,H)\big)\to\Gamma\big(L^2([0,t[,H)\big)$ is the unitary isomorphism introduced in \ref{subsub-anti-mon}. Anti-monotone stochastic integrals can be defined using this isomorphism. We call an operator process $\{X_t\}_{0\le t\le T}$ on the monotone Fock space anti-monotonically adapted, if $\{\theta^*_{\rm M}X_t\theta_{\rm M}\}_{0\le t\le T}$ is adapted on the boson Fock space $\Gamma\big(L^2([0,t[,H)\big)$ and define the integral by
\[
\int_0^T X_t{\rm d}I_t^{\rm AM} := \theta_{\rm M}\left(\int_0^T \theta^*_{\rm M}X_t\theta_{\rm M}{\rm d}I_t\right)\theta^*_{\rm M}
\]
for
\begin{eqnarray*}
{\rm d}I_t^{\rm AM} &=& {\rm d}\Lambda^{\rm AM}_t\big(|x\rangle\langle y|\big) + {\rm d}A^{\rm AM+}_t(u)+ {\rm d}A^{\rm AM}_t(v), \\ 
{\rm d}I_t &=& \Gamma(0_{0t})\Big({\rm d}\Lambda_t\big(|x\rangle\langle y|\big) + {\rm d}A^{\rm AM+}_t(u)+ {\rm d}A^{\rm AM}_t(v)\Big),
\end{eqnarray*}
for $x,y,u,v\in H$. In this way all the domains, kernels, etc., defined in \cite[Chapter 2]{schuermann93} can be translated to the monotone Fock space.

Using the form of the comultiplication of $(\tilde{\mathcal{B}},\overline{\Delta}_{\rm AM},\tilde\varepsilon)$, the quantum stochastic equation for the L\'evy process on the involutive bialgebra $(\tilde{\mathcal{B}},\overline{\Delta}_{\rm AM},\tilde\varepsilon)$ that we associated to an anti-monotone L\'evy process on the dual semigroup $(\mathcal{B},\Delta,\varepsilon)$ in Theorem \ref{theo-one-to-one}, and Theorem \ref{theo-real}, one can now derive a representation theorem for anti-monotone L\'evy processes on dual semigroups.

To state our result we need the free product $\coprod^0$ without unification of units. This is the coproduct in the category of all $*$-algebras (not necessarily unital). The two free products $\coprod$ and $\coprod^0$ are related by
\[
(\mathbb{C}\mathbf{1}\oplus \mathcal{A})\coprod(\mathbb{C}\mathbf{1}\oplus \mathcal{B})\cong \mathbb{C}\mathbf{1}\oplus \left(\mathcal{A}{\coprod}^0\mathcal{B}\right).
\]

We will use the notation $\Gamma_{\rm M}(0_{st})=\theta_{\rm M}^*\Gamma(0_{st})\theta_{\rm M}$, $0\le s\le t \le T$.

\begin{theorem}
Let $(\mathcal{B},\Delta,\varepsilon)$ be a dual semigroup and let $(\rho,\eta,\psi)$ be a Sch\"urmann triple on $\mathcal{B}$ over some pre-Hilbert space $D$. Then the anti-monotone stochastic differential equations
\begin{equation}\label{eq-qsde-anti-mon}
j_{st}(b)=\int_s^t \Big(j_{s\tau}{\coprod}^0 {\rm d}I_\tau^{\rm AM}\Big)\circ
\Delta(b), \quad\mbox{ for }b\in\mathcal{B}^0={\rm ker}\,\varepsilon,
\end{equation}
with
\[
{\rm d}I^{\rm AM}_\tau(b)={\rm d}\Lambda^{\rm AM}_t\big(\rho(b)\big) + {\rm d}A^{\rm AM+}_t\big(\eta(b)\big)+ {\rm d}A^{\rm AM}_t\big(\eta(b^*)\big) + \psi(b)\Gamma_{\rm M}(0_{0\tau}){\rm d}t,
\]
have solutions (unique in $\theta_{\rm M}^*\underline{\mathcal{A}}_D\theta_{\rm M}$). If we set $j_{st}(\mathbf{1}_\mathcal{B})={\rm id}$, then $\{j_{st}\}_{0\le s\le t\le T}$ is an anti-monotone L\'evy process on the dual semigroup $(\mathcal{B},\Delta,\varepsilon)$ with respect to the vacuum state. Furthermore, any anti-monotone L\'evy process on the dual semigroup $(\mathcal{B},\Delta,\varepsilon)$ with generator $\psi$ is equivalent to  $\{j_{st}\}_{0\le s\le t\le T}$.
\end{theorem}
\begin{remark}
Let $b\in\mathcal{B}^0$, $\Delta(b)=\sum_{\epsilon\in\mathbb{A}} b^\epsilon$, $b^\epsilon\in \mathcal{B}_\epsilon$, then Equation (\ref{eq-qsde-anti-mon}) has to be interpreted as
\[
j_{st}(b) = \sum_{{\epsilon\in\mathbb{A}}\atop{\epsilon_1=1,\epsilon\not=(1)}} \int_s^t j_{s\tau}(b^\epsilon_1){\rm d}I_\tau^{\rm AM}(b^\epsilon_2)j_{s\tau}(b^\epsilon_3)\cdots + \sum_{{\epsilon\in\mathbb{A}}\atop{\epsilon_1=2}} \int_s^t {\rm d}I_\tau^{\rm AM}(b^\epsilon_1)j_{s\tau}(b^\epsilon_2){\rm d}I_\tau^{\rm AM}(b^\epsilon_3)\cdots,
\]
see also \cite{schuermann95b}.
This equation can be simplified using the relation
\[
{\rm d}I_t^{\rm AM}(b_1)X_t{\rm d}I_t^{\rm AM}(b_2) = \langle\Omega,X_t\Omega\rangle \left({\rm d}I_t^{\rm AM}(b_1)\bullet{\rm d}I_t^{\rm AM}(b_2)\right)
\]
for $b_1,b_2\in\mathcal{B}^0$ and anti-monotonically adapted operator processes $\{X_t\}_{0\le t\le T}$, where the product `$\bullet$' is defined by the anti-monotone It\^o table
\[
\begin{array}{|c||c|c|c|c|}
\hline
\bullet & {\rm d}A^{\rm AM+}(u_1) & {\rm d}\Lambda^{\rm AM}\big(|x_1\rangle\langle y_1|\big) & {\rm d}A^{\rm AM}(v_1) & {\rm d}t \\
\hline\hline 
{\rm d}A^{\rm AM+}(u_2) & 0 & 0 & 0 & 0 \\ \hline
{\rm d}\Lambda^{\rm AM}\big(|x_2\rangle\langle y_2|\big) & \langle y_2,u_1\rangle {\rm d}A^{\rm AM+}(x_2) & \langle y_2,x_1\rangle {\rm d}\Lambda^{\rm AM}\big(|x_2\rangle\langle y_1|\big) & 0 & 0 \\ \hline
{\rm d}A^{\rm AM}(v_2) & \langle v_2,u_1\rangle \Gamma_{\rm M}(0_{0t}){\rm d}t & \langle v_2,x_1\rangle{\rm d}A^{\rm AM}(y_1) & 0 & 0 \\ \hline
{\rm d}t & 0 & 0 & 0 & 0 \\ \hline
\end{array}
\]
for $u_i,v_i,x_i,y_i\in D$, $i=1,2$.

One can check that ${\rm d}I^{\rm AM}_t$ is actually a homomorphism on $\mathcal{B}^0$ for the It\^o product, i.e.\
\[
{\rm d}I_t^{\rm AM}(b_1)\bullet{\rm d}I_t^{\rm AM}(b_2)={\rm d}I_t^{\rm AM}(b_1b_2),
\]
for all $b_1,b_2\in\mathcal{B}^0$.
\end{remark}

Using the time-reversal $R$ defined in \ref{subsub-mon}, we also get a realization of monotone L\'evy processes on the monotone Fock space as solutions of backward monotone stochastic differential equations.

It follows also that operator processes with monotonically or anti-monotonically independent additive increments can be written as linear combination of the four fundamental noises, where the time process has to be taken as $T^{\rm AM}_{st}=\int_s^t\Gamma_{\rm M}(0_{0\tau}){\rm d}\tau$, $0\le s\le t\le T$, for the anti-monotone case and $T^{\rm M}_{st}=\int_s^t\Gamma_{\rm M}(0_{\tau T}){\rm d}\tau$ for the monotone case.

%% file: levy-lie.tex
\chapter{Renormalized Squares of White Noise and other Non-Gaussian Noises as L\'evy Processes on Real Lie Algebras}\label{levy-lie}

\markboth{L\'EVY PROCESSES ON REAL LIE ALGEBRAS}{L\'EVY ON REAL LIE ALGEBRAS}

\vspace*{2cm}

\begin{quote}
It is shown how the relations of the renormalized squared white noise defined by Accardi, Lu, and Volovich \cite{accardi+lu+volovich99} can be realized as factorizable current representations or L\'evy processes on the real Lie algebra $\eufrak{sl}_2$. This allows to obtain its It\^o table, which turns out to be infinite-dimensional. The linear white noise without or with number operator is shown to be a L\'evy process on the Heisenberg-Weyl Lie algebra or the oscillator Lie algebra. Furthermore, a joint realization of the linear and quadratic white noise relations is constructed, but it is proved that no such realizations exist with a vacuum that is an eigenvector of the central element and the annihilator. Classical L\'evy processes are shown to arise as components of L\'evy processes on real Lie algebras and their distributions are characterized. In particular the square of white noise analogue of the quantum Poisson
process is shown to have a $\chi^2$ probability density and the analogue
of the field operators to have a density proportional to
$|\Gamma ( {m_0 + ix \over 2})|^2$ where $\Gamma $ is the usual
$\Gamma $--function and $m_0$ a real parameter.
\end{quote}

\vspace*{2cm}

Joint work with Luigi Accardi and Michael Skeide. Published in Communications in Mathematical Physics Vol.~220, No.~1, pp.~123-150, 2002.

\newpage

\section{Introduction}

The stochastic limit of quantum theory \cite{accardi+lu+volovich02} shows that stochastic
equations (both classical and quantum) are equivalent to white noise
Hamiltonian equations. This suggests a natural extension of stochastic
calculus to higher powers of white noise. The program to develop such an
extension was formulated in \cite{accardi+lu+volovich98} where it was also shown that
it requires some kind of renormalization. As a first step towards the
realization of this program a new type of renormalization was introduced
in \cite{accardi+lu+volovich99} which led to a closed set of algebraic
relations for the renormalized square of white noise (SWN) and to the
construction of a Hilbert space representation for these relations.
This construction was extended by \'Sniady \cite{sniady00} to a family of
processes including non Boson noises and simplified in
\cite{accardi+skeide99a} who also showed that the interacting Fock space
constructed in \cite{accardi+lu+volovich99} was in fact canonically
isomorphic to the Boson Fock space of the finite difference algebra,
introduced by Feinsilver \cite{feinsilver89} and Boukas
\cite{boukas88,boukas91a}. Commenting upon this result U.Franz, and independently a few months later
K.R.Parthasarathy, (private communications) pointed out that the commutation
relations of the SWN define a L\'evy process on the Lie algebra of
$SL(2, \mathbb{R})$ or, equivalently, a representation of a current algebra
over this Lie algebra, and suggested that the theory of representations
of current algebras, developed in the early seventies by Araki, Streater,
Parthasarathy, Schmidt, Guichardet, ...\ (see \cite{parthasarathy+schmidt72,guichardet72} and the references therein) might be used to produce a more
direct construction of the Fock representation of the SWN as well as
different ones.
In the present paper we prove that this is indeed the case.
As a by--product we reduce the stochastic integration with respect to
the SWN to the usual stochastic integration in the sense of Hudson and
Parthasarathy \cite{parthasarathy92} and this also allows to write down
their corresponding It\^o tables (see Equation \eqref{eq def dL}).

After the renormalization procedure (which we shall not discuss here,
simply taking its output as our starting point) the algebraic relations,
defining the SWN are:
\begin{subequations}\label{swn-rel}
\begin{equation}
b_\phi b^+_\psi - b^+_\psi b_\phi = \gamma \langle \phi,\psi\rangle +
n_{\overline{\phi}\psi},
\end{equation}
\begin{equation}
n_\phi b_\psi - b_\psi n_\phi = -2b_{\overline{\phi}\psi},
\end{equation}
\begin{equation}
n_\phi b^+_\psi - b^+_\psi n_\phi = 2b^+_{\phi\psi},
\end{equation}
\begin{equation}
(b_\phi)^* = b^+_\phi, \qquad (n_\phi)^* = n_{\overline{\phi}},
\end{equation}
\end{subequations}
where $\gamma$ is a fixed strictly positive real parameter (coming from the
renormalization) and
$$\phi,\psi\in
\Sigma(\mathbb{R}_+)=\{\phi=\sum_{i=1}^n\phi_i
\mathbf 1_{[s_i,t_i[};\phi_i\in\mathbb C,
s_i<t_i\in\mathbb R_+,n\in\mathbb N\}$$
the algebra of step functions on $\mathbb R_+$ with bounded support and
finitely many values. Furthermore $b^+$ and $n$ are linear and $b$ is
anti-linear in the test functions.

We want to find a Hilbert space representation of these relations, i.e.
we want to construct an Hilbert space $\mathcal{H}$, a dense subspace
$D\subseteq \mathcal{H}$ and three maps $b,b^+,n$ from $\Sigma(\mathbb{R}_+)$
to $\mathcal{L}(D)$, the algebra of adjointable linear operators on $D$,
such that the above relations are satisfied.

The simple current algebra $\eufrak g^{\mathbb{T}}$ of a real Lie algebra
$\eufrak g$ over a measure space $(\mathbb{T},\mathcal{T},\mu)$ is defined
as the space of simple functions on $\mathbb{T}$ with values in $\eufrak g$,
\[
\eufrak g^{\mathbb{T}} =
\left\{X=\sum_{i=1}^nX_i\mathbf 1_{M_i};X_i\in\eufrak g, M_i\in\mathcal{T},
n\in\mathbb N\right\}.
\]
This is a real Lie algebra with the Lie bracket and the involution defined
pointwise. The SWN relations (\ref{swn-rel}) imply that any realization of
SWN on a pre-Hilbert space $D$ defines a representation $\pi$ of the current
algebra $\eufrak{sl}_2^{\mathbb{R}_+}$ of the real Lie algebra
$\eufrak{sl}_2$
over $\mathbb{R}_+$ (with the Borel $\sigma$-algebra and the Lebesgue
measure)
on $D$ by
\[
B^-\mathbf 1_{[s,t[} \mapsto b_{1_{[s,t[}}, \quad B^+\mathbf 1_{[s,t[}
\mapsto
b^+_{1_{[s,t[}}, \quad M \mathbf 1_{[s,t[} \mapsto \gamma(t-s)+
n_{1_{[s,t[}},
\]
where $\eufrak{sl}_2$ is the three-dimensional real Lie algebra spanned
by $\{B^+,B^-,M\}$, with the commutation relations
\[
[B^-,B^+]=M, \quad [M,B^\pm]=\pm 2 B^\pm,
\]
and the involution $(B^-)^*=B^+$, $M^*=M$. The converse is obviously also
true, every representation of the current algebra
$\eufrak{sl}_2^{\mathbb{R}_+}$ defines a realization of the SWN
relations (\ref{swn-rel}). Looking only at indicator functions of intervals
we get a family of $*$-representations $(j_{st})_{0\le s\le t}$ on $D$ of
the Lie algebra $\eufrak{sl}_2$,
\[
j_{st}(X) = \pi(X\mathbf 1_{[s,t[}), \qquad
\text{ for all } X\in \eufrak{sl}_2.
\]
By the universal property these $*$-representations extend to
$*$-representations of the universal enveloping algebra
$\mathcal U(\eufrak{sl}_2)$ of $\eufrak{sl}_2$. If there exists a
vector $\Omega$ in $\mathcal L(D)$ such that the representations
corresponding to disjoint intervals are independent (in the sense of
Definition \ref{def levy-lie}, Condition 2), i.e.\ if they commute and their
expectations in the state $\Phi(\cdot)=\langle\Omega,\cdot\,\Omega\rangle$
factorize, then $(j_{st})_{0\le s\le t}$ is a L\'evy process on $\eufrak{sl}_2$ (in the sense of Definition \ref{def
levy}).
This condition is satisfied in the constructions in
\cite{accardi+lu+volovich99,accardi+skeide99a,sniady00}.
They are of `Fock type' and have a fixed special vector, the so-called
vacuum, and the corresponding vector state has the desired factorization
property.

On the other hand, given a L\'evy process on $\eufrak{sl}_2$ on a
pre-Hilbert space $D$, we can construct a realization of the SWN
relations (\ref{swn-rel}) on $D$. Simply set
\[
b_\phi=
\sum_{i=1}^n\overline{\phi_i} j_{s_i,t_i}(B^-),\,
b^+_\phi=\sum_{i=1}^n\phi_i j_{s_i,t_i}(B^+),\,
n_\phi= \sum_{i=1}^n\phi_i\big(j_{s_i,t_i}(M) - \gamma(t_i-s_i){\rm
id}_D\big),
\]
for $\phi=\sum_{i=1}^n\phi_i\mathbf 1_{[s_i,t_i[}\in\Sigma(\mathbb R_+)$.

We see that in order to construct realizations of the SWN relations we can
construct L\'evy processes on $\eufrak{sl}_2$. Furthermore, all realizations
that have a vacuum vector in which the expectations factorize, will arise
in this way.

In this paper we show how to classify the L\'evy processes on
$\eufrak{sl}_2$
and how to construct realizations of these L\'evy processes acting on
(a dense subspace of) the symmetric Fock space over $L^2(\mathbb{R}_+,H)$ for
some Hilbert space $H$. Given the {\itshape generator} $L$ of a L\'evy
process, we immediately can write down a realization of the process;
see Equation \eqref{eq def j}.
The theory of L\'evy processes has been developed for arbitrary involutive
bialgebras, cf.\ \cite{accardi+schuermann+waldenfels88,schuermann93}, but
here it will be sufficient to consider enveloping algebras of Lie algebras.
This allows some simplification, in particular we do not need to
make explicit use of the coproduct.
The construction of this sub--class of L\'evy process is based on the
theory of ``factorizable unitary representation of current algebras" and
the abelian subprocesses of these processes are the stationary
independent increment processes of classical probability (cf.\ Section 4
below).

As already specified, the SWN naturally leads to the real Lie algebra
$\eufrak{sl}_2$, but we shall also consider
several other real Lie algebras, including the Heisenberg-Weyl Lie
algebra $\eufrak{hw}$, the oscillator Lie algebra $\eufrak{osc}$, and
the finite-difference Lie algebra $\eufrak{fd}$.

This paper is organized as follows.
In Section \ref{prelim}, we recall the definition of L\'evy processes on
real Lie algebras and present their fundamental properties. We also outline
how the L\'evy processes on a given real Lie algebra can be characterized
and
constructed as a linear combination of the four basic processes of
Hudson-Parthasarathy quantum stochastic calculus: number, creation,
annihilation and time.

In Section \ref{examples}, we list all Gaussian L\'evy processes or L\'evy
processes associated to integrable unitary irreducible representations for
several real Lie algebras in terms of their generators or Sch{\"u}rmann
triples (see Definition \ref{def schuermann}). We also give explicit
realizations on a boson Fock space for several examples. These examples
include the processes on the finite-difference Lie algebra defined by Boukas
\cite{boukas88,boukas91a} and by Parthasarathy and Sinha
\cite{parthasarathy+sinha91} as well as a process on $\eufrak{sl}_2$
that has been considered previously by Feinsilver and Schott
\cite[Section 5.IV]{feinsilver+schott93}.
See also \cite{gelfand+graev+vershik73} for factorizable current
representations of current groups over $SL(2,\mathbb{R})$.

Finally, in Section \ref{class}, we show that the restriction of a L\'evy
process to one single hermitian element of the real Lie algebra always gives rise
to a classical L\'evy process. We give a characterization of this process in
terms of its Fourier transform. For several examples we also explicitly
compute its L\'evy measure or its marginal distribution. It turns out
that the densities of self-adjoint linear combinations of the SWN
operators $b_{\mathbf{1}_{[s,t[}}$, $b^+_{\mathbf{1}_{[s,t[}}$,
$n_{\mathbf{1}_{[s,t[}}$ in the realization considered
in \cite{accardi+lu+volovich99,accardi+skeide99a,sniady00} are the
measures of orthogonality of the Laguerre, Meixner, and Meixner-Pollaczek
polynomials.

\section{L\'evy processes on real Lie algebras}\label{prelim}

In this section we give the basic definitions and properties of L\'evy processes on real Lie algebras. This is a special case of the theory of L\'evy processes on involutive bialgebras, for more detailed accounts on these processes see \cite{schuermann93},\cite[Chapter VII]{meyer95},\cite{franz+schott99}. For a list of references on factorizable representations of current groups and algebras and a historical survey, we refer to \cite[Section 5]{streater00}.

By a real Lie algebra we will mean a pair $\eufrak{g}_\mathbb{R}=(\eufrak{g},*)$ consisting of a  Lie algebra $\eufrak{g}$ over the field of complex numbers $\mathbb{C}$ and an involution $*:\eufrak{g}\to\eufrak{g}$. These pairs are in one-to-one correspondence with the Lie algebras over the field of real numbers $\mathbb{R}$. To recover a Lie algebra $\eufrak{g}_0$ over $\mathbb{R}$ from a pair $(\eufrak{g},*)$, simply take the anti-hermitian elements, i.e.\ set $\eufrak{g}_0=\{X\in \eufrak{g}| X^*=-X\}$. Note that it is not possible to take the hermitian elements, because the commutator of two hermitian elements in not again hermitian. Given  a Lie algebra $\eufrak{g}_0$ over $\mathbb{R}$, the involution on its complexification $\eufrak{g}=\eufrak{g}_0\oplus i\eufrak{g}_0$ is defined by $(X+iY)^*=-X+iY$ for $X,Y\in\eufrak{g}_0$.

We denote by $\mathcal{U}(\eufrak{g})$ the universal enveloping algebra of $\eufrak{g}$ and by $\mathcal{U}_0(\eufrak{g})$ the non-unital subalgebra of $\mathcal{U}$ generated by $\eufrak{g}$. If $X_1,\ldots,X_d$ is a basis of $\eufrak{g}$, then
\[
\{X_1^{n_1}\cdots X_d^{n_d} | n_1,\ldots,n_d\in\mathbb{N}, n_1+\cdots+n_d\ge 1\}
\]
is a basis of $\mathcal{U}_0(\eufrak{g})$. Furthermore, we extend the involution on $\eufrak{g}$ as an anti-linear anti-homomorphism to $\mathcal{U}(\eufrak{g})$ and $\mathcal{U}_0(\eufrak{g})$.

\begin{definition}\label{def levy-lie}
Let $D$ be a pre-Hilbert space and $\Omega\in D$ a unit vector. We call a family $\big(j_{st}:\mathcal{U}(\eufrak{g})\to\mathcal{L}(D)\big)_{0\le s\le t}$ of unital $*$-representations of $\mathcal{U}(\eufrak{g})$ a {\em L\'evy process} on $\eufrak{g}_\mathbb{R}$ over $D$ (with respect to $\Omega$), if the following conditions are satisfied.
\begin{enumerate}
\item
(Increment property) We have
\[
j_{st}(X)+j_{tu}(X)=j_{su}(X)
\]
for all $0\le s\le t \le u$ and all $X\in\eufrak{g}$.
\item
(Boson independence)
We have $[j_{st}(X),j_{s't'}(Y)]=0$ for all $X,Y\in\eufrak{g}$, $0\le s\le t\le s'\le t'$ and
\[
\langle\Omega, j_{s_1t_1}(u_1)\cdots j_{s_nt_n}(u_n)\Omega\rangle = \langle\Omega, j_{s_1t_1}(u_1)\Omega\rangle\cdots\langle\Omega, j_{s_nt_n}(u_n)\Omega\rangle
\]
for all $n\in\mathbb{N}$, $0\le s_1\le t_1\le s_2\le \cdots \le t_n$, $u_1,\ldots,u_n\in\mathcal{U}(\eufrak{g})$.
\item
(Stationarity) The functional $\varphi_{st}:\mathcal{U}(\eufrak{g})\to\mathbb{C}$ defined by
\[
\varphi_{st}(u)=\langle\Omega,j_{st}(u)\Omega\rangle, \qquad u\in\mathcal{U}_0(\eufrak{g}),
\]
depends only on the difference $t-s$.
\item
(Weak continuity)
We have $\lim_{t\searrow s} \langle\Omega,j_{st}(u)\Omega\rangle =0$ for all $u\in\mathcal{U}_0(\eufrak{g})$.
\end{enumerate}
\end{definition}

If $(j_{st})_{0\le s\le t}$ is a L\'evy process on $\eufrak{g}_\mathbb{R}$, then the functionals $\varphi_{t}=\langle\Omega,j_{0t}(\cdot)\Omega\rangle:\mathcal{U}(\eufrak{g})\to\mathbb{C}$ are actually states. Furthermore, they are differentiable w.r.t.\ $t$ and
\[
L(u)=\lim_{t\searrow 0}\frac{1}{t}\varphi_{t}(u), \qquad u\in\mathcal{U}_0(\eufrak{g}),
\]
defines a positive hermitian linear functional on $\mathcal{U}_0(\eufrak{g})$. In fact one can prove that the family $(\varphi_{t})$ is a convolution
semigroup on $\eufrak{g}_\mathbb{R}$ whose generator is $L$. The functional
$L$ is also called the {\em generator} of the process.

Let $\big(j^{(1)}_{st}:\mathcal{U}(\eufrak{g})\to\mathcal{L}(D^{(1)})\big)_{0\le s\le t}$ and $\big(j^{(2)}:\mathcal{U}(\eufrak{g})\to\mathcal{L}(D^{(2)})\big)_{0\le s\le t}$ be two L\'evy processes on $\eufrak{g}_\mathbb{R}$ with respect to the state vectors $\Omega^{(1)}$ and $\Omega^{(2)}$, resp. We call them {\em equivalent}, if all their moments agree, i.e.\ if
\[
\langle\Omega^{(1)}, j^{(1)}_{s_1t_1}(u_1)\cdots j^{(1)}_{s_nt_n}(u_n)\Omega^{(1)}\rangle = \langle\Omega^{(2)}, j^{(2)}_{s_1t_1}(u_1)\cdots j^{(2)}_{s_nt_n}(u_n)\Omega^{(2)}\rangle,
\]
for all $n\in\mathbb{N}$, $0\le s_1\le t_1\le s_2\le \cdots \le t_n$, $u_1,\ldots,u_n\in\mathcal{U}(\eufrak{g})$.

By a GNS-type construction, one can associate to every generator a Sch{\"u}rmann triple.

\begin{definition}\label{def schuermann}
A {\em Sch\"urmann triple} on $\eufrak{g}_\mathbb{R}$ is a triple $(\rho,\eta,L)$, where $\rho$ is a $*$-representation of $\mathcal{U}_0(\eufrak{g})$ on some pre-Hilbert space $D$, $\eta:\mathcal{U}_0(\eufrak{g})\to D$ is a surjective $\rho$-1-cocycle, i.e.\ it satisfies
\[
\eta(uv)=\rho(u)\eta(v),
\]
for all $u,v\in\mathcal{U}_0(\eufrak{g})$, and $L:\mathcal{U}_0(\eufrak{g})\to\mathbb{C}$ is a hermitian linear functional such that the bilinear map $(u,v)\mapsto -\langle\eta(u^*),\eta(v)\rangle$ is the 2-coboundary of $L$ (w.r.t.\ the trivial representation), i.e.\
\[
L(uv)=\langle\eta(u^*),\eta(v)\rangle
\]
for all $u,v\in\mathcal{U}_0(\eufrak{g})$.
\end{definition}

Let $(\rho,\eta,L)$ be a Sch{\"u}rmann triple on $\eufrak{g}_\mathbb{R}$, acting on a pre-Hilbert space $D$. We can define a L\'evy process on the symmetric Fock space $\Gamma\big(L^2(\mathbb{R}_+,D)\big)=\bigoplus_{n=0}^\infty L^2(\mathbb{R}_+,D)^{\odot n}$ by setting
\begin{equation}\label{eq def j}
j_{st}(X)=\Lambda_{st}\big(\rho(X)\big)+A_{st}^*\big(\eta(X)\big)+A_{st}\big(\eta(X^*)\big)+L(X)(t-s){\rm id},
\end{equation}
for $X\in\eufrak{g}$, where $\Lambda_{st}$, $A^*_{st}$, $A_{st}$ denote the conservation, creation, and annihilation processes on $\Gamma\big(L^2(\mathbb{R}_+,D)\big)$, cf.\ \cite{parthasarathy92,meyer95}. It is straightforward to check that we have
\[
\big[j_{st}(X),j_{st}(Y)\big]=j_{st}\big([X,Y]\big), \quad\mbox{ and }\quad j_{st}(X)^*=j_{st}(X^*)
\]
for all $0\le s\le t$, $X,Y\in\eufrak{g}$. By the universal property, the family
\[
\left(j_{st}:\eufrak{g}\to\mathcal{L}\Big(\Gamma\big(L^2(\mathbb{R}_+,D)\big)\Big)\right)_{0\le s\le t}
\]
extends to a unique family $(j_{st})_{0\le s\le t}$ of unital $*$-representations of $\mathcal{U}(\eufrak{g})$, and it is not difficult to verify that this family is a L\'evy process with generator $L$ on $\eufrak{g}_\mathbb{R}$ over $\Gamma\big(L^2(\mathbb{R}_+,D)\big)$ with respect to the Fock vacuum $\Omega$.

The following theorem shows that the correspondence between (equivalence
classes of) L\'evy processes and Sch{\"u}rmann triples is one--to--one
and that the representation (\ref{eq def j}) is universal.

\begin{theorem}\label{theo class}\cite{schuermann93}
Two L\'evy processes on $\eufrak{g}_\mathbb{R}$ are equivalent if and only if their Sch{\"u}rmann triples are unitarily equivalent.
A L\'evy process $(k_{st})_{0\le s\le t}$ with generator $L$ and
Sch{\"u}rmann
triple $(\rho,\eta,L)$ is equivalent to the L\'evy process
$(j_{st})_{0\le s\le t}$ associated to $(\rho,\eta,L)$
defined in Equation (\ref{eq def j}).
\end{theorem}

\begin{remark}
Since we know the It\^o table for the four H-P integrators,
\[
\begin{array}{|c||c|c|c|c|}
\hline
\bullet & {\rm d}A^*(u) & {\rm d}\Lambda(F) & {\rm d}A(u) & {\rm d}t \\
\hline\hline
{\rm d}A^*(v) & 0 & 0 & 0 & 0 \\ \hline
{\rm d}\Lambda(G) & {\rm d}A^*(Gu) &  {\rm d}\Lambda(GF) & 0 & 0 \\ \hline
{\rm d}A(v) & \langle v,u\rangle {\rm d}t & {\rm d}A(F^*v) & 0 & 0 \\ \hline
{\rm d}t & 0 & 0 & 0 & 0 \\ \hline
\end{array}
\]
for all $F,G\in \mathcal L(D)$, $u,v\in D$, we can deduce the It\^o tables for the L\'evy processes on $\eufrak{g}_\mathbb{R}$. The map ${\rm d}_L$ associating elements $u$ of the universal enveloping algebra to the corresponding quantum stochastic differentials ${\rm d}_Lu$ defined by
\begin{equation}\label{eq def dL}
{\rm d}_L u={\rm d}\Lambda\big(\rho(u)\big)+ {\rm d}A^*\big(\eta(u)\big) + {\rm d}A\big(\eta(u^*)\big) + L(u){\rm d}t,
\end{equation}
is a $*$-homomorphism from $\mathcal U_0(\eufrak g)$ to the It\^o algebra over $D$, see \cite[Proposition 4.4.2]{franz+schott99}. It follows that the dimension of the It\^o algebra generated by $\{{\rm d}_LX;X\in\eufrak g\}$ is at least the dimension of $D$ (since $\eta$ is supposed surjective) and not bigger than $(\dim D +1)^2$. If $D$ is infinite-dimensional, then its dimension is also infinite. Note that it depends on the choice of the L\'evy process.
\end{remark}

Due to Theorem \ref{theo class}, the problem of characterizing and constructing all L\'evy processes on a given real Lie algebra can be decomposed into the following steps. First, classify all $*$-representations of $\mathcal{U}(\eufrak{g})$ (modulo unitary equivalence), this will give the possible choices for the representation $\rho$ in the Sch\"urmann triple. Next determine all surjective $\rho$-1-cocycles. We distinguish between trivial cocycles, i.e.\ cocycles which are of the form
\[
\eta(u)=\rho(u)\omega, \qquad u\in\mathcal{U}_0(\eufrak{g})
\]
for some vector $\omega\in D$ in the representation space of $\rho$, and non-trivial cocycles, i.e.\ cocycles, which can not be written in this form. We will denote the space of all cocycles of a given $*$-representation $\rho$ on some pre-Hilbert space $D$ by $Z^1(\mathcal{U}_0(\eufrak{g}),\rho,D)$, that of trivial ones by $B^1(\mathcal{U}_0(\eufrak{g}),\rho,D)$. The quotient $H^1(\mathcal{U}_0(\eufrak{g}),\rho,D)=Z^1(\mathcal{U}_0(\eufrak{g}),\rho,D)/B^1(\mathcal{U}_0(\eufrak{g}),\rho,D)$ is called the first cohomology group of $\rho$. In the last step we determine all generators $L$ that turn a pair $(\rho,\eta)$ into a Sch\"urmann triple $(\rho,\eta,L)$. This can again also be viewed as a cohomological problem. If $\eta$ is a $\rho$-1-cocycle, then the bilinear map $(u,v)\mapsto -\langle\eta(u^*),\eta(v)\rangle$ is a 2-cocycle for the trivial representation, i.e.\ it satisfies $-\big\langle\eta\big((uv)^*\big),\eta(w)\big\rangle+\big\langle\eta(u^*),\eta(vw)\big\rangle=0$ for all $u,v,w\in\mathcal{U}_0(\eufrak{g})$. For $L$ we can take any hermitian functional that has the map $(u,v)\mapsto -\langle\eta(u^*),\eta(v)\rangle$ as coboundary, i.e.\ $L$ has to satisfy $L(u^*)=\overline{L(u)}$ and $L(uv)=\langle\eta(u^*),\eta(v)\rangle$ for all $u,v\in\mathcal{U}_0(\eufrak{g})$. If $\eta$ is trivial, then such a functional always exists, we can take $L(u)=\langle\omega,\rho(u)\omega\rangle$. For a given pair $(\rho,\eta)$, $L$ is determined only up to a hermitian $0$-$1$-cocycle, i.e.\ a hermitian functional $\ell$ that satisfies $\ell(uv)=0$ for all $u,v\in\mathcal{U}_0(\eufrak{g})$.

\begin{remark}
A linear $*$-map $\pi:\eufrak{g}\to\mathcal{L}(D)$ is called a projective $*$-representation of $\eufrak{g}$, if there exists a bilinear map $\alpha:\eufrak{g}\times\eufrak{g}\to\mathbb{C}$, such that
\[
\big[\pi(X),\pi(Y)\big]=\pi\big([X,Y]\big)+\alpha(X,Y){\rm id},
\]
for all $X,Y\in\eufrak{g}$. Every projective $*$-representation defines a $*$-representation of a central extension $\tilde{\eufrak{g}}$ of $\eufrak{g}$. As a vector space $\tilde{\eufrak{g}}$ is defined as $\tilde{\eufrak{g}}=\eufrak{g}\oplus \mathbb{C}$. The Lie bracket and the involution are defined by
\[
\big[(X,\lambda),(Y,\mu)\big] = \big([X,Y],\alpha(X,Y)\big), \qquad (X,\lambda)^*=(X^*,\overline{\lambda})
\]
for $(X,\lambda),(Y,\mu)\in\tilde{\eufrak{g}}$. It is not hard to check that
\[
\tilde{\pi}\big((X,\lambda)\big)=\pi(X)+\lambda\,{\rm id}
\]
defines a $*$-representation of $\tilde{\eufrak{g}}$. If the cocycle $\alpha$ is trivial, i.e.\ if there exists a (hermitian) linear functional $\beta$ such that $\alpha(X,Y)=\beta([X,Y])$ for all $X,Y\in\eufrak{g}$, then the central extension is trivial, i.e.\ $\tilde{\eufrak{g}}$ is isomorphic to the direct sum of $\eufrak{g}$ with the (abelian) one-dimensional Lie algebra $\mathbb{C}$. Such an isomorphism is given by $\eufrak{g}\oplus\mathbb{C}\ni(X,\mu)\mapsto (X,\beta(X)+\mu)\in\tilde{\eufrak{g}}$. This implies that in this case
\[
\pi_\beta(X)=\tilde\pi\big((X,\beta(X))\big)=\pi(X)+\beta(X){\rm id}
\]
defines a $*$-representation of $\eufrak{g}$.

For a pair $(\rho,\eta)$ consisting of a $*$-representation $\rho$ and a $\rho$-1-cocycle $\eta$ we can always define a family of projective $*$-representations $(k_{st})_{0\le s\le t}$ of $\eufrak{g}$ by setting
\[
k_{st}(X) = \Lambda_{st}\big(\rho(X)\big)+A_{st}^*\big(\eta(X)\big)+A_{st}\big(\eta(X^*)\big),
\]
for $X\in\eufrak{g}$, $0\le s\le t$. Using the commutation relations of the creation, annihilation, and conservation operators, one finds that the 2-cocycle $\alpha$ is given by $(X,Y)\mapsto \alpha(X,Y)=\langle\eta(X^*),\eta(Y)\rangle-\langle\eta(Y^*),\eta(X)\rangle$. If it is trivial, then $(k_{st})_{0\le s\le t}$ can be used to define a L\'evy process on $\eufrak{g}$. More precisely, if there exists a hermitian functional $\psi$ on $\mathcal{U}_0(\eufrak{g})$ such that $\psi(uv)=\langle\eta(u^*),\eta(v)\rangle$ holds for all $u,v\in\mathcal{U}_0(\eufrak{g})$, then $(\rho,\eta,\psi)$ is a Sch\"urmann triple on $\eufrak{g}$ and therefore defines a L\'evy process on $\eufrak{g}$. But even if such a hermitian functional $\psi$ does not exist, we can define a L\'evy process on $\tilde{\eufrak{g}}$ by setting
\[
\tilde{k}_{st}\big((X,\lambda)\big) = \Lambda_{st}\big(\rho(X)\big)+A^*\big(\eta(X)\big)+A\big(\eta(X^*)\big)+(t-s)\lambda\,{\rm id},
\]
for $(X,\lambda)\in\tilde{\eufrak{g}}$, $0\le s\le t$.
\end{remark}

We close this section with two useful lemmata on cohomology groups.

Sch\"urmann triples $(\rho,\eta,L)$, where the $*$-representation $\rho$ is equal to the trivial representation defined by $0:\mathcal{U}_0(\eufrak{g})\ni u\mapsto 0\in \mathcal{L}(D)$ are called {\em Gaussian}, as well as the corresponding processes, cocycles, and generators
(cf. Corollary \ref{class cor 1} for a justification of this definition). The following lemma completely classifies all Gaussian cocycles of a given Lie algebra.

\begin{lemma}\label{lem H1(0)}
Let $D$ be an arbitrary complex vector space, and $0$ the trivial representation of $\eufrak g$ on $D$. We have
\[
Z^1(\mathcal U_0(\eufrak g),0,D) \cong \big(\eufrak g/[\eufrak g,\eufrak g]\big)^*, \qquad B^1(\mathcal U_0(\eufrak g),0,D)=\{0\},
\]
and therefore $\dim H^1(\mathcal U_0(\eufrak g),0,D)=\dim \eufrak g/[\eufrak g,\eufrak g]$.
\end{lemma}
\begin{proof}
Let $\phi$ be a linear functional on $\eufrak g/[\eufrak g,\eufrak g]$, then we can extend it to a unique $0$-$1$-cocycle on the algebra $\mathcal U_0(\eufrak g/[\eufrak g,\eufrak g])$ (this is the free abelian algebra over $\eufrak g/[\eufrak g,\eufrak g]$), which we denote by $\tilde\phi$. Denote by $\pi$ the canonical projection from $\eufrak g$ to $\eufrak g/[\eufrak g,\eufrak g]$, by the universal property of the enveloping algebra it has a unique extension $\tilde\pi:\mathcal U_0(\eufrak g)\to\mathcal U_0(\eufrak g/[\eufrak g,\eufrak g])$. We can define a cocycle $\eta_\phi$ on $\mathcal U_0(\eufrak g)$ by $\eta_\phi=\tilde\phi\circ\tilde\pi$. Furthermore, since any $0$-$1$-cocycle on $\mathcal U_0(\eufrak g)$ has to vanish on $[\eufrak g,\eufrak g]$ (because $Y=[X_1,X_2]$ implies $\eta(Y)=0\eta(X_2)-0\eta(X_1)=0$), the map $\phi\mapsto\eta_\phi$ is bijective. 
\end{proof}

The following lemma shows that a representation of $\mathcal U(\eufrak g)$ can only have non-trivial cocycles, if the center of $\mathcal U_0(\eufrak g)$ acts trivially.

\begin{lemma}\label{lem Hn(rho)}
Let $\rho$ be a representation of $\eufrak g$ on some vector space $D$ and let $C\in\mathcal U_0(\eufrak g)$ be central. If $\rho(C)$ is invertible, then
\[
H^1(\mathcal U_0(\eufrak g),\rho,D)=\{0\}.
\]
\end{lemma}
\begin{proof}
Let $\eta$ be a $\rho$-cocycle on $\mathcal{U}_0(\eufrak{g})$ and $C\in\mathcal{U}_0(\eufrak{g})$ such that $\rho(C)$ is invertible. Then we get
\[
\rho(C)\eta(u)=\eta(Cu)=\eta(uC)=\rho(u)\eta(C)
\]
and therefore $\eta(u)=\rho(u)\rho(C)^{-1}\eta(C)$ for all $u\in\mathcal{U}_0(\eufrak{g})$, i.e.\ $\eta(u)=\rho(u)\omega$, where $\omega=\rho(C)^{-1}\eta(C)$. This shows that all $\rho$-cocycles are trivial.
\end{proof}

\section{Examples}\label{examples}

In this section we completely classify the Gaussian generators for several real Lie algebras and determine the non-trivial cocycles for some or all of their integrable unitary irreducible representations, i.e.\ those representations that arise by differentiating unitary irreducible representations of the corresponding Lie group. These are $*$-representations of the enveloping algebra $\mathcal U(\eufrak{g})$ on some pre-Hilbert space $D$ for which the Lie algebra elements are mapped to essentially self-adjoint operators. For some of the processes we give explicit realizations on the boson Fock space.

\subsection{White noise or L\'evy processes on $\eufrak{hw}$ and $\eufrak{osc}$}\label{exa hw-osc}

The Heisenberg-Weyl Lie algebra $\eufrak{hw}$ is the three-dimensional Lie algebra with basis $\{A^+,A^-,E\}$, commutation relations
\[
[A^-,A^+]=E, \quad [A^\pm,E]=0,
\]
and involution $(A^-)^*=A^+$, $E^*=E$. Adding a hermitian element $N$ with commutation relations
\[
[N,A^\pm] = \pm A^\pm, \quad [N,E]=0,
\]
we obtain the four-dimensional oscillator Lie algebra $\eufrak{osc}$.

We begin with the classification of all Gaussian generators on these two Lie algebras.

\begin{proposition}\label{prop hw and osc}
\begin{description}
\item[a)]
Let $v_1,v_2\in \mathbb C^2$ be two vectors and $z\in \mathbb C$ an arbitrary complex number. Then
\begin{gather*}
\rho(A^+)=\rho(A^-)=\rho(E)=0, \\
\eta(A^+)=v_1, \quad \eta(A^-)=v_2,\quad \eta(E)=0, \\
L(A^+)=z,\quad L(A^-)=\overline{z}, \quad L(E)=||v_1||^2-||v_2||^2,
\end{gather*}
defines the Sch\"urmann triple on $D={\rm span}\,\{v_1,v_2\}$ of a Gaussian generator on $\mathcal{U}_0(\eufrak{hw})$. Furthermore, all Gaussian generators on $\mathcal{U}_0(\eufrak{hw})$ arise in this way. 
\item[b)]
The Sch\"urmann triples of Gaussian generators on $\mathcal{U}_0(\eufrak{osc})$ are all of the form
\begin{gather*}
\rho(N)=\rho(A^+)=\rho(A^-)=\rho(E)=0, \\
\eta(N)=v,\quad \eta(A^+)= \eta(A^-)=\eta(E)=0, \\
L(N)=b,\quad L(A^+)=L(A^-)=L(E)=0,
\end{gather*}
with $v\in\mathbb C$, $b\in\mathbb R$.
\end{description}
\end{proposition}
\begin{proof}
The form of the Gaussian cocycles on $\mathcal{U}_0(\eufrak{hw})$ and $\mathcal{U}_0(\eufrak{osc})$ follows from Lemma \ref{lem H1(0)}. Then one checks that for all these cocycles there do indeed exist generators and computes their general form.
\end{proof}

Therefore from (\ref{eq def dL})  we get, for an arbitrary Gaussian L\'evy process
on $\eufrak{hw}$:
\begin{eqnarray*}
{\rm d}_LA^+ &=& {\rm d}A^*(v_1)+{\rm d}A(v_2)+z{\rm d}t, \\
{\rm d}_LA^- &=&{\rm d}A^*(v_2)+{\rm d}A(v_1)+\overline{z}{\rm d}t, \\
{\rm d}_LE &=& \big( ||v_1||^2 - ||v_2||^2\big){\rm d}t,
\end{eqnarray*} 
and the It\^o table
\[
\begin{array}{|c||c|c|c|}
\hline
\bullet & {\rm d}_LA^+ & {\rm d}_LA^- & {\rm d}_LE \\
\hline\hline
{\rm d}_LA^+ & \langle v_2,v_1\rangle {\rm d}t &\langle v_2,v_2\rangle {\rm d}t & 0 \\ \hline
{\rm d}_LA^- & \langle v_1,v_1\rangle {\rm d}t & \langle v_1,v_2\rangle {\rm d}t  & 0 \\ \hline
{\rm d}_LE & 0 & 0 & 0 \\ \hline
\end{array}.
\]
For $||v_1||^2=1$ and $v_2=0$, this is the usual It\^o table for the creation and annihilation process in Hudson-Parthasarathy calculus.

Any integrable unitary irreducible representation of $\eufrak{hw}$ is equivalent either to one of the one-dimensional representations defined by
\[
\pi_z(A^+)=z,\quad \pi_z(A^-)=\overline{z},\quad \pi_z(E)=0,
\]
for some $z\in\mathbb C$, or to one of the infinite-dimensional representations defined by
\begin{equation}\label{fock}
\rho_h(A^+)e_n=\sqrt{(n+1)h}\,e_{n+1}, \quad \rho_h(A^-)e_n=\sqrt{nh}\,e_{n-1}, \quad \rho_h(E)e_n=he_n,
\end{equation}
and
\[
\rho_{-h}(A^-)e_n=\sqrt{(n+1)h}\,e_{n+1}, \quad \rho_{-h}(A^+)e_n=\sqrt{nh}\,e_{n-1}, \quad \rho_{-h}(E)e_n=-he_n,
\]
where $h>0$, and $\{e_0,e_1,\ldots\}$ is a orthonormal basis of $\ell^2$. By Lemma \ref{lem Hn(rho)}, the representations $\rho_h$ have no non-trivial cocycles. But by a simple computation using the defining relations of $\eufrak{hw}$ we see that, for $z\not=0$, the representations of the form $\pi_z\,{\rm id}_D$ also have only one trivial cocycle. From $A^+E=EA^+$ we get
\[
z\eta(E)=\eta(A^+E)=\eta(EA^+)=\pi_z(E)\eta(A^+)=0
\]
and therefore $\eta(E)=0$. But $E=A^-A^+-A^+A^-$ implies
\[
0=\eta(E)=\pi_z(A^-)\eta(A^+)-\pi_z(A^+)\eta(A^-)=\overline{z}\eta(A^+)-z\eta(A^-),
\]
and we see that $\eta$ is the coboundary of $\omega=z^{-1}\eta(A^+)$. Thus the integrable unitary irreducible representations (except the trivial one) of $\eufrak{hw}$ have no non-trivial cocycles.

Let us now consider the oscillator Lie algebra $\eufrak{osc}$. The elements $E$ and $NE-A^+A^-$ generate the center of $\mathcal U_0(\eufrak{osc})$. If we want an irreducible representation of $\mathcal U(\eufrak{osc})$, which has non-trivial cocycles, they have to be represented by zero. But this implies that we have also $\rho(A^+)=\rho(A^-)=0$ (since we are only interested in $*$-representations). Thus we are lead to study the representations $\rho_\nu$ defined by
\[
\rho_\nu(N)=\nu\, {\rm id_D}, \quad\rho_\nu(A^+)=\rho_\nu(A^-)=\rho_\nu(E)=0,
\]
with $\nu\in\mathbb R\backslash\{0\}$. It is straightforward to determine all their cocycles and generators.

\begin{proposition}
For $\nu\in\mathbb R$, $\nu\not\in\{-1,0,1\}$, all cocycles of $\rho_\nu$ are of the form
\[
\eta(N)=v,\quad \eta(A^+)=\eta(A^-)=\eta(E)=0,
\]
for some $v\in D$ and thus trivial (coboundaries of $\omega=\nu^{-1}v$).

For $\nu=1$ they are of the form
\[
\eta(N)=v_1,\quad\eta(A^+)=v_2,\quad\eta(A^-)=\eta(E)=0,
\]
and for $\nu=-1$ of the  form
\[
\eta(N)=v_1,\quad\eta(A^-)=v_2,\quad\eta(A^+)=\eta(E)=0,
\]
with some vectors $v_1,v_2\in D$. Therefore we get
\[
\dim H^1(\mathcal U_0(\eufrak{osc}),\rho_{\pm 1},D)=1, \quad \dim B^1(\mathcal U_0(\eufrak{osc}),\rho_{\pm 1},D)=1
\]
and
\[
\dim H^1(\mathcal U_0(\eufrak{osc}),\rho_{\nu},D)=0, \quad \dim B^1(\mathcal U_0(\eufrak{osc}),\rho_{\nu},D)=1
\]
for $\nu\in\mathbb R\backslash\{-1,0,1\}$.
\end{proposition}

Let now $\nu=1$, the case $\nu=-1$ is similar, since $\rho_1$ and $\rho_{-1}$ are related by the automorphism $N\mapsto -N$, $A^+\mapsto A^-$, $A^-\mapsto A^+$, $E\mapsto -E$. It turns out that for all the cocycles given in the preceding proposition there exists a generator, and we obtain the following result.

\begin{proposition}
Let $v_1,v_2\in \mathbb C^2$ and $b\in \mathbb R$. Then $\rho=\rho_1$,
\begin{gather*}
\eta(N)=v_1,\quad\eta(A^+)=v_2,\quad\eta(A^-)=\eta(E)=0, \\
L(N)=b, \quad L(E)=||v_2||^2, \quad L(A^+)=\overline{L(A^-)}=\langle v_1,v_2\rangle,
\end{gather*}
defines a Sch\"urmann triple on $\eufrak{osc}$ acting on $D={\rm span}\,\{v_1,v_2\}$. The corresponding quantum stochastic differentials are
\begin{eqnarray*}
{\rm d}_LN &=& {\rm d}\Lambda({\rm id})+{\rm d}A^*(v_1)+{\rm d}A(v_1)+b {\rm d}t  , \\
{\rm d}_LA^+ &=& {\rm d}A^*(v_2) +\langle v_1,v_2\rangle{\rm d}t , \\
{\rm d}_LA^- &=& {\rm d}A(v_2)+\langle v_2,v_1\rangle {\rm d}t, \\
{\rm d}_LE &=&   ||v_2||^2{\rm d}t,
\end{eqnarray*}
and they satisfy the following It\^o table
\[
\begin{array}{|c||c|c|c|c|}
\hline
\bullet & {\rm d}_LA^+ & {\rm d}_LN & {\rm d}_LA^- & {\rm d}_LE \\ \hline\hline
{\rm d}_LA^+ & 0 & 0 & 0 & 0 \\ \hline
{\rm d}_LN & {\rm d}_LA^+ & {\rm d}_LN+\big(||v_1||^2-b\big){\rm d}t & 0 & 0 \\ \hline
{\rm d}_LA^- & {\rm d}_LE  & {\rm d}_LA^- & 0 & 0 \\ \hline
{\rm d}_LE & 0 & 0 & 0 & 0 \\ \hline
\end{array}.
\]  
\end{proposition}

Note that for $||v_1||^2=b$, this is the usual It\^o table of the four fundamental noises of Hudson-Parthasarathy calculus.

\subsection{SWN or L\'evy processes on $\eufrak{sl}_2$}\label{exa sl2}

The Lie algebra $\eufrak{sl}_2$ is the three-dimensional simple Lie algebra with basis $\{B^+,B^-,M\}$, commutation relations 
\[
[B^-,B^+]=M, \quad [M,B^\pm]=\pm 2 B^\pm,
\]
and involution $(B^-)^*=B^+$, $M^*=M$.
Its center is generated by the Casimir element
\[
C=M(M-2)-4B^+B^-=M(M+2)-4B^-B^+.
\]

We have $[\eufrak{sl}_2,\eufrak{sl}_2]=\eufrak{sl}_2$, and so $\mathcal U_0(\eufrak{sl}_2)$ has no Gaussian cocycles, cf.\ Lemma \ref{lem H1(0)}, and therefore no Gaussian generators either. Let us now determine all the non-trivial cocycles for the integrable unitary irreducible representations of $\eufrak{sl}_2$.

It is known that, beyond the trivial representation $\rho_0$ there are
three families of equivalence classes of integrable unitary irreducible representation
of $\eufrak{sl}_2$ (given in Equations (\ref{lowest weight}),
(\ref{heighest weight}), (\ref{other sl2 reps}) below), see, e.g., \cite{gruber+lenczewski+lorente90} and the references therein. We will
consider them separately. We begin to consider the lowest and highest weight
representations. These families of representations are parametrized by a
real number $m_0$ and are
induced by $\rho(M)\Omega=m_0\Omega$, $\rho(B^-)\Omega=0$, and $\rho(M)\Omega=-m_0\Omega$, $\rho(B^+)\Omega=0$, respectively. The lowest weight representations are spanned by the vectors $v_n=\rho(B^+)^n\Omega$, with $n\in\mathbb N$. We get
\begin{eqnarray*}
\rho(B^+)v_n &=& v_{n+1}, \\
\rho(B^-)v_n &=& \rho\left(B^-(B^+)^n\right)\Omega = \rho\left(\frac{1}{4}\big(M(M+2)-C\big)(B^+)^{n-1}\right)\Omega \\
&=& n(n+m_0-1)\rho(B^+)^{n-1}\Omega=n(n+m_0-1)v_{n-1}, \\
\rho(M)v_n &=&(2n+m_0)v_n.
\end{eqnarray*}
If we want to define an inner product on ${\rm span}\,\{v_n;n\in\mathbb N\}$ such that $\rho(M)^*=\rho(M)$ and $\rho(B^-)^*=\rho(B^+)$, then the $v_n$ have to be orthogonal and their norms have to satisfy the recurrence relation
\begin{equation}\label{eq norms}
||v_{n+1}||^2 = \langle \rho(B^+)v_n,v_{n+1}\rangle = \langle v_n,\rho(B^-)v_{n+1}\rangle = (n+1)(n+m_0) ||v_n||^2.
\end{equation}
It follows there exists an inner product on  ${\rm span}\,\{v_n;n\in\mathbb N\}$ such that the lowest weight representation with $\rho(M)\Omega=m_0\Omega$, $\rho(B^-)\Omega=0$ is a $*$-representation, if and only if the coefficients $(n+1)(n+m_0)$ in Equation (\ref{eq norms}) are non-negative for all $n=0,1,\ldots$, i.e.\ if and only if $m_0\ge 0$. For $m_0=0$ we get the trivial one-dimensional representation $\rho_0(B^+)\Omega=\rho_0(B^-)\Omega=\rho_0(M)\Omega=0$ (since $||v_1||^2=0$), for $m_0>0$ we get
\begin{subequations}\label{lowest weight}
\begin{equation}
\rho^+_{m_0}(B^+)e_n= \sqrt{(n+1)(n+m_0)}\,e_{n+1},
\end{equation}
\begin{equation}
\rho^+_{m_0}(M)e_n=(2n+m_0)e_n,
\end{equation}
\begin{equation}
\rho^+_{m_0}(B^-)e_n= \sqrt{n(n+m_0-1)}\,e_{n-1},
\end{equation}
\end{subequations}
where $\{e_0,e_1,\ldots\}$ is an orthonormal basis of $\ell^2$. Note that the Casimir element acts as $\rho^+_{m_0}(C)e_n=m_0(m_0-2)e_n$.
Similarly we see that there exists a $*$-representation $\rho$ containing a vector $\Omega$ such that $\rho(B^+)\Omega=0$, $\rho(M)\Omega=-m_0\Omega$, if and only if $m_0\ge 0$. For $m_0=0$ this is the trivial representation, for $m_0>0$ it is of the form
\begin{subequations}\label{heighest weight}
\begin{equation}
\rho^-_{m_0}(B^-)e_n= \sqrt{(n+1)(n+m_0)}\,e_{n+1},
\end{equation}
\begin{equation}
\rho^-_{m_0}(M)e_n=-(2n+m_0)e_n,
\end{equation}
\begin{equation}
\rho^-_{m_0}(B^+)e_n= \sqrt{n(n+m_0-1)}\,e_{n-1},
\end{equation}
\end{subequations}
and $\rho^-_{m_0}(C)e_n=m_0(m_0-2)e_n$.
The integrable unitary irreducible representations of
$\eufrak{sl}_2$, belonging to the third class, have no highest
or lowest weight vector. They are parametrized by two real numbers $m_0,c$
and are  induced by $\rho(M)\Omega=m_0\Omega$, $\rho(C)\Omega= c \Omega$. Note
that since $C$ is central, the second relation implies actually
$\rho(C)=c\,{\rm id}$. The vectors $\{v_{\pm n}=\rho(B^\pm)^n\Omega;n\in\mathbb N\}$ form a basis for the induced representation,
\begin{eqnarray*}
\rho(M) v_n &=& (2n +m_0)v_n \\
\rho(B^+)v_n &=&
\left\{
\begin{array}{lcl}
v_{n+1} & \mbox{ if} & n\ge 0, \\
\frac{(m_0+2n+2)(m_0+2n)-c}{4}v_{n+1} & \mbox{ if} & n< 0,
\end{array}
\right. \\
\rho(B^-)v_n &=&
\left\{
\begin{array}{lcl}
\frac{(m_0+2n-2)(m_0+2n)-c}{4}v_{n-1} & \mbox{ if} & n>0, \\
v_{n-1} & \mbox{ if} & n\le 0, \\
\end{array}
\right.
\end{eqnarray*}
We look again for an inner product that turns this representation into a $*$-representation. The $v_n$ have to be orthogonal for such an inner product and their norms have to satisfy the recurrence relations
\begin{gather*}
||v_{n+1}||^2 = \frac{(m_0+2n+2)(m_0+2n)-c}{4}||v_n||^2, \quad \mbox{ for } n\ge 0, \\
||v_{n-1}||^2 = \frac{(m_0+2n-2)(m_0+2n)-c}{4}||v_n||^2, \quad \mbox{ for } n\le 0.
\end{gather*}
Therefore we can define a positive definite inner product on ${\rm span}\,\{v_n;n\in\mathbb Z\}$, if and only if $\lambda(\lambda+2)>c$ for all $\lambda\in m_0+2\mathbb Z$. We can restrict ourselves to $m_0\in[0,2[$, because the representations induced by $(c,m_0)$ and $(c,m_0+2k)$, $k\in\mathbb Z$ turn out to be unitarily equivalent. We get the following family of integrable unitary irreducible representations of $\mathcal U(\eufrak{sl}_2)$,
\begin{subequations}\label{other sl2 reps}
\begin{equation}
\rho_{cm_0}(B^+)e_n = \frac{1}{2}\sqrt{(m_0+2n+2)(m_0+2n)-c\,}\, e_{n+1},
\end{equation}
\begin{equation}
\rho_{cm_0}(M)e_n = (2n+m_0)e_n,
\end{equation}
\begin{equation}
\rho_{cm_0}(B^-)e_n = \frac{1}{2}\sqrt{(m_0+2n-2)(m_0+2n)-c\,}\, e_{n-1},
\end{equation}
\end{subequations}
where $\{e_n;n\in\mathbb Z\}$ is an orthonormal basis of $\ell^2(\mathbb Z)$, $m_0\in[0,2[$, $c<m_0(m_0-2)$.

Due to Lemma \ref{lem Hn(rho)}, we are interested in representations in which $C$ is mapped to zero. There are, up to unitary equivalence, only three such representations, the trivial or zero representation (which has no non-zero cocycles at all, by Lemma \ref{lem H1(0)}), and the two representations $\rho^\pm=\rho^\pm_2$ on $\ell^2$ defined by
\begin{eqnarray*}
\rho^\pm(M)e_n &=& \pm (2n+2) e_n, \\
\rho^+(B^+)e_n &=& \sqrt{(n+1)(n+2)}\,e_{n+1}, \\
\rho^+(B^-)e_n &=& \sqrt{n(n+1)}\,e_{n-1}, \\
\rho^-(B^+)e_n &=& \sqrt{n(n+1)}\,e_{n-1}, \\
\rho^-(B^-)e_n &=& \sqrt{(n+1)(n+2)}\,e_{n+1},
\end{eqnarray*}
for $n\in\mathbb{N}$, where $\{e_0,e_1,\ldots\}$ is an orthonormal basis of $\ell^2$. The representations $\rho^+$ and $\rho^-$ are not unitarily equivalent, but they are related by the automorphism $M\mapsto -M$, $B^+\mapsto B^-$, $B^-\mapsto B^+$. Therefore it is sufficient to study $\rho^+$. Let $\eta$ be a $\rho^+$-$1$-cocycle. Since $\rho^+(B^+)$ is injective, we see that $\eta$ is already uniquely determined by $\eta(B^+)$, since the relations $[M,B^+]=2B^+$ and $[B^-,B^+]=M$ imply
\begin{eqnarray*}
\eta(M) &=& \rho^+(B^+)^{-1}\big(\rho^+(M)-2\big)\eta(B^+), \\
\eta(B^-) &=& \rho^+(B^+)^{-1}\big(\rho^+(B^-)\eta(B^+) - \eta(M)\big).
\end{eqnarray*}
In fact, we can choose any vector for $\eta(B^+)$, the definitions above and the formula $\eta(uv)=\rho^+(u)\eta(v)$ for $u,v\in\mathcal U_0(\eufrak{sl}_2)$ will extend it to a unique $\rho^+$-$1$-cocycle. This cocycle is a coboundary, if and only if the coefficient $v_0$ in the expansion $\eta(B^+)=\sum_{n=0}^\infty v_n e_n$ of $\eta(B^+)$ vanishes, and an arbitrary $\rho^+$-$1$-cocycle is a linear combination of the non-trivial cocyle $\eta_1$ defined by
\begin{equation}\label{def eta1}
\eta_1\left((B^+)^nM^m(B^-)^r\right) = \left\{
\begin{array}{lcl}
0 & \mbox{ if } & n=0, \\
\delta_{r,0}\delta_{m,0} \rho(B^+)^{n-1}e_0 &\mbox{ if } & n\ge 1,
\end{array}
\right.
\end{equation}
and a coboundary. In particular, for $\eta$ with $\eta(B^+)=\sum_{n=0}^\infty v_n e_n$, we get $\eta= v_0\eta_1 + \partial\omega$ with $\omega=\sum_{n=0}^\infty \frac{v_{n+1}}{\sqrt{(n+1)(n+2)}} e_n$. Thus we have shown the following.

\begin{proposition}
We have
\[
\dim H^1(\mathcal U_0(\eufrak{sl}_2),\rho^\pm,\ell^2)= 1
\]
and $\dim H^1(\mathcal U_0(\eufrak{sl}_2),\rho,\ell^2)=0$ for all other integrable unitary irreducible representations of $\eufrak{sl}_2$.
\end{proposition}

Since $[\eufrak{sl}_2,\eufrak{sl}_2]=\eufrak{sl}_2$, all elements of $\mathcal U_0(\eufrak{sl}_2)$ can be expressed as linear combinations of products of elements of $\mathcal U_0(\eufrak{sl}_2)$. Furthermore one checks that
\[
L(u)=\langle\eta(u^*_1),\eta(u_2)\rangle, \qquad \mbox{ for } u=u_1u_2, \quad u_1,u_2\in \mathcal U_0(\eufrak{sl}_2)
\]
is independent of the decomposition of $u$ into a product and defines a hermitian linear functional. Thus there exists a unique generator for every cocycle on $\eufrak{sl}_2$.

\begin{example}\label{exa boukas}
We will now construct the L\'evy process for the cocycle $\eta_1$ defined in Equation (\ref{def eta1}) and the corresponding generator. We get
\begin{eqnarray*}
L(M) &=& \langle\eta_1(B^+),\eta_1(B^+)\rangle - \langle\eta_1(B^-),\eta_1(B^-)\rangle = 1, \\
L(B^+) &=& L(B^-) = 0,
\end{eqnarray*}
and therefore
\begin{equation}\label{eq swn real}
\begin{array}{rcl}
{\rm d}_LM &=& {\rm d}\Lambda\big(\rho^+(M)\big) + {\rm d}t, \\
{\rm d}_LB^+ &=& {\rm d}\Lambda\big(\rho^+(B^+)\big) + {\rm d}A^*(e_0), \\
{\rm d}_LB^- &=& {\rm d}\Lambda\big(\rho^+(B^-)\big) + {\rm d}A(e_0). 
\end{array}
\end{equation}
The It\^o table is infinite-dimensional. This is the process that leads to the realization of SWN that was constructed in the previous works \cite{accardi+lu+volovich99,accardi+skeide99a,sniady00}.

For the Casimir element we get
\[
{\rm d}_LC = -2{\rm d}t.
\]

For this process we have $j_{st}(B^-)\Omega=0$ and $j_{st}(M)\Omega=(t-s)\Omega$ for all $0\le s\le t$. From our previous considerations about the lowest weight representation of $\eufrak{sl}_2$ we can now deduce that for fixed $s$ and $t$ the representation $j_{st}$ of $\eufrak{sl}_2$ restricted to the subspace $j_{st}\big(\mathcal{U}(\eufrak{sl}_2)\big)\Omega$ is equivalent to the representation $\rho^+_{t-s}$ defined in Equation (\ref{lowest weight}).
\end{example}

\begin{example}\label{exa PaSi}
Let now $\rho$ be one of the lowest weight representations defined in (\ref{lowest weight}) with $m_0>0$, and let $\eta$ be the trivial cocycle defined by
\[
\eta(u)=\rho^+_{m_0}(u)e_0,
\]
for $u\in\mathcal U_0(\eufrak{sl}_2)$. There exists a unique generator for this cocycle, and the corresponding L\'evy process is defined by
\begin{equation}\label{eq swn real PaSi}
\begin{array}{rcl}
{\rm d}_LM &=& {\rm d}\Lambda\big(\rho^+_{m_0}(M)\big)+ m_0{\rm d}A^*(e_0)+ m_0{\rm d}A(e_0)+ m_0 {\rm d}t, \\
{\rm d}_LB^+ &=& {\rm d}\Lambda\big(\rho^+_{m_0}(B^+)\big) + \sqrt{m_0}{\rm d}A^*(e_1), \\
{\rm d}_LB^- &=& {\rm d}\Lambda\big(\rho^+_{m_0}(B^-)\big) + \sqrt{m_0}{\rm d}A(e_1).
\end{array}
\end{equation}
For the Casimir element we get
\[
{\rm d}_LC = m_0(m_0-2)\big({\rm d}\Lambda({\rm id}) + {\rm d}A^*(e_0)+ {\rm d}A(e_0) + {\rm d}t\big).
\]
\end{example}

\subsection{White noise and its square or L\'evy processes on $\eufrak{sl}_2\oplus_\alpha \eufrak{hw}$}\label{exa sl2+hw}

We can define an action $\alpha$ of the Lie algebra $\eufrak{sl}_2$ on $\eufrak{hw}$ by
\[
\alpha(M) :
\begin{array}{lll}
A^+ &  \mapsto & A^+, \\
E & \mapsto & 0, \\
A^- & \mapsto & -A^-,
\end{array}
\quad
\alpha(B^+) :
\begin{array}{lll}
A^+ &  \mapsto & 0, \\
E & \mapsto & 0, \\
A^- & \mapsto & -A^+,
\end{array}
\quad
\alpha(B^-):
\begin{array}{lll}
A^+ &  \mapsto & A^-, \\
E & \mapsto & 0, \\
A^- & \mapsto & 0.
\end{array}
\]
The $\alpha(X)$ are derivations and satisfy $\big(\alpha(X)Y\big)^*= -\alpha(X^*)Y^*$ for all $X\in \eufrak{sl}_2$, $Y\in \eufrak{hw}$. Therefore we can define a new Lie algebra $\eufrak{sl}_2\oplus_\alpha \eufrak{hw}$ as the semi-direct sum of $\eufrak{sl}_2$ and $\eufrak{hw}$, it has the commutation relations $\big[(X_1,Y_1),(X_2,Y_2)\big]=\big([X_1,X_2],[Y_1,Y_2]+\alpha(X_1)Y_2 - \alpha(X_2)Y_1\big)$ and the involution $(X,Y)^*=(X^*,Y^*)$. In terms of the basis $\{B^\pm,M,A^\pm,E\}$ the commutation relations are
\begin{gather*}
[B^-,B^+]=M \quad [M,B^\pm]=\pm 2 B^\pm, \\{}
[A^-,A^+]=E, \quad [E,A^\pm]=0, \\{}
[B^\pm,A^\mp]=\mp A^\pm, \quad [B^\pm,A^\pm]=0, \\{}
[M,A^\pm]=\pm A^\pm, \quad [E,B^\pm]=0, \quad [M,E]=0. 
\end{gather*}
The action $\alpha$ has been chosen in order to obtain these relations, which also follow from the renormalization rule introduced in \cite{accardi+lu+volovich02}.

In the following we identify $\mathcal U(\eufrak{hw})$ and $\mathcal U(\eufrak{sl}_2)$ with the corresponding subalgebras in $\mathcal U(\eufrak{sl}_2\oplus_\alpha \eufrak{hw})$.

Note that for any $c\in \mathbb R$, ${\rm span}\{N=M+cE,A^+,A^-,E\}$ forms a Lie subalgebra of $\eufrak{sl}_2\oplus_\alpha \eufrak{hw}$ that is isomorphic to $\eufrak{osc}$.

There exist no Gaussian L\'evy processes on $\eufrak{sl}_2\oplus_\alpha \eufrak{hw}$, since $[\eufrak{sl}_2\oplus_\alpha \eufrak{hw},\eufrak{sl}_2\oplus_\alpha \eufrak{hw}]=\eufrak{sl}_2\oplus_\alpha \eufrak{hw}$. But, like for every real Lie algebra, there exist non-trivial $*$-representations of $\eufrak{sl}_2\oplus_\alpha \eufrak{hw}$, and thus also L\'evy processes, it is sufficient to take, e.g., a trivial cocycle.

The following result shows that the usual creation and annihilation calculus cannot be extended to a joint calculus of creation and annihilation and their squares.
\begin{proposition}\label{prop no ext of hw}
Let $(\rho,\eta,L)$ be the Sch\"urmann triple on $\eufrak{hw}$ defined in Proposition \ref{prop hw and osc} a), and denote the corresponding L\'evy process by $(j_{st})_{0\le s\le t}$. There exists no L\'evy process $(\tilde \jmath_{st})_{0\le s\le t}$ on $\eufrak{sl}_2\oplus_\alpha \eufrak{hw}$ such that
\[
\left(\tilde \jmath_{st}|_{\mathcal U(\eufrak{hw})}\right) \cong (j_{st}),
\]
unless $(j_{st})_{0\le s\le t}$ is trivial, i.e.\ $j_{st}(u)=0$ for all $u\in\mathcal{U}_0(\eufrak{hw})$.
\end{proposition}
\begin{proof}
We will assume that $(\tilde \jmath_{st})$ exists and show that this implies $||v_1||^2=||v_2||^2=|z|^2=0$, i.e.\ $L=0$.

Let $(\tilde\rho,\tilde\eta,\tilde L)$ be the Sch\"urmann triple of $(\tilde \jmath_{st})$. If $\left(\tilde \jmath_{st}|_{\mathcal U(\eufrak{hw})}\right) \cong (j_{st})$, then we have $\tilde L|_{\mathcal U_0(\eufrak{hw})}=L$, and therefore the triple on $\eufrak{hw}$ obtained by restriction of $(\tilde\rho,\tilde\eta,\tilde L)$ is equivalent to $(\rho,\eta,L)$ and there exists an isometry from $D=\eta\big(\mathcal U_0(\eufrak{hw})\big)$ into $\tilde D=\tilde\eta\big(\mathcal U_0(\eufrak{sl}_2\oplus_\alpha \eufrak{hw})\big)$, such that we have
\[
\tilde\rho|_{\mathcal U(\eufrak{hw})\times D}=\rho, \quad \mbox{ and } \tilde\eta|_{\mathcal U_0(\eufrak{hw})}=\eta,
\]
if we identify $D$ with its image in $\tilde D$.

{}From $[B^+,A^-]=- A^+$ and $[B^-,A^+]=A^-$, we get
\begin{gather*}
-\tilde\eta(A^+)=\tilde\rho(B^+)\eta(A^-)-\tilde\rho(A^-)\tilde\eta(B^+), \\
\tilde\eta(A^-)=\tilde\rho(B^-)\eta(A^+)-\tilde\rho(A^+)\tilde\eta(B^-).
\end{gather*}
Taking the inner product with $\tilde\eta(A^+)=\eta(A^+)=v_1$ and $\tilde\eta(A^-)=\eta(A^-)=v_2$, resp., we get
\begin{eqnarray*}
-||v_1||^2 &=& \langle v_1,\rho(B^+)v_2\rangle -\langle v_1, \tilde\rho(A^-)\tilde\eta(B^+)\rangle \\
&=& \langle v_1,\rho(B^+)v_2\rangle -\langle \rho(A^+) v_1,\tilde\eta(B^+)\rangle
= \langle v_1,\rho(B^+)v_2\rangle, \\
||v_2||^2 &=& \langle v_2,\rho(B^-)v_1\rangle,
\end{eqnarray*}
since $\tilde\rho(A^\pm)|_D=\rho(A^\pm)$. Therefore
\[
-||v_1||^2 = \langle v_1,\rho(B^+)v_2\rangle = \overline{\langle v_2,\rho(B^-)v_1\rangle} = ||v_2||^2,
\]
and thus $||v_1||^2=||v_2||^2=0$. But $A^+=-[B^+,A^-]$ and
\begin{eqnarray*}
L(A^+) &=& \tilde L(A^+) = \langle \tilde\eta(A^+),\tilde\eta(B^+)\rangle - \langle\tilde\eta(B^-),\tilde\eta(A^-)\rangle \\
&=& \langle v_1,\tilde\eta(B^+)\rangle - \langle\tilde\eta(B^-),v_2\rangle
\end{eqnarray*}
which now implies that $z=L(A^+)=0$.
\end{proof}

\'Sniady \cite{sniady00} has posed the question, if it is possible to define a joint calculus for the linear white noise and the square of white noise. Formulated in our context, his answer to this question is that there exists no L\'evy process on $\eufrak{sl}_2\oplus_\alpha\eufrak{hw}$ such that
\[
j_{st}(E)=(t-s){\rm id},\quad\mbox{ and }\quad j_{st}(A^-)\Omega=j_{st}(B^-)\Omega=0,
\]
for all $0\le s\le t$. We are now able to show the same under apparently much weaker hypotheses. 

\begin{corollary}
Every L\'evy process on $\eufrak{sl}_2\oplus_\alpha\eufrak{hw}$ such that the state vector $\Omega$ is an eigenvector for $j_{st}(E)$ and $j_{st}(A^-)$ for some pair $s$ and $t$ with $0\le s < t$ is trivial on $\eufrak{hw}$, i.e.\ it has to satisfy $j_{st}|_{\mathcal{U}_0(\eufrak{hw})}=0$ for all $0\le s\le t$. 
\end{corollary}
\begin{proof}
Assume that such a L\'evy process exists. Then it would be equivalent to its realization on a boson Fock space defined by Equation (\ref{eq def j}). Therefore we see that the state vector is an eigenvector of $j_{st}(E)$ and $j_{st}(A^-)$, if and only if the Sch\"urmann triple of $(j_{st})_{0\le s\le t}$ satisfies $\eta(E)=\eta(A^-)=0$. If we show that the only Sch{\"u}rmann triples on $\eufrak{hw}$ satisfying this condition are the Gaussian Sch{\"u}rmann triples, then our result follows from Proposition \ref{prop no ext of hw}.

Let $(\rho,\eta,L)$ be a Sch\"urmann triple on $\eufrak{hw}$ such that $\eta(E)=\eta(A^-)=0$. Then the vector $\eta(A^+)$ has to be cyclic for $\rho$. We get
\[
\rho(E)\eta(A^+)=\rho(A^+)\eta(E)=0,
\]
since $E$ and $A^+$ commute. From $[A^-,A^+]=E$, we get
\[
\rho(A^-)\eta(A^+)=\rho(A^+)\eta(A^-)+\eta(E)=0.
\]
But
\begin{eqnarray*}
||\rho(A^+)\eta(A^+)||^2 &=& \langle \eta(A^+),\rho(A^-)\rho(A^+)\eta(A^+)\rangle \\
&=& \langle \eta(A^-),\rho(A^+)\rho(A^-)\eta(A^+)\rangle+\langle\eta(A^+),\rho(E)\eta(A^+)\rangle \\
&=& 0
\end{eqnarray*}
shows that $\rho(A^+)$ also acts trivially on $\eta(A^+)$ and therefore the restriction of the triple $(\rho,\eta,L)$ to $\mathcal U(\eufrak{hw})$ is Gaussian.
\end{proof}

The SWN calculus defined in Example \ref{exa boukas} can only be extended in the trivial way, i.e.\ by setting it equal to zero on $\eufrak{hw}$, $\tilde\jmath_{st}|_{\eufrak{hw}}=0$.

\begin{proposition}\label{prop no ext of sl2}
Let $(j_{st})_{0\le s\le t}$ be the L\'evy process on $\eufrak{sl}_2$ defined in (\ref{eq swn real}). The only L\'evy process $(\tilde\jmath_{st})_{0\le s\le t}$ on $\eufrak{sl}_2\oplus_\alpha \eufrak{hw}$ such that
\[
(\tilde\jmath_{st}|_{\mathcal U(\eufrak{sl}_2)})\cong (j_{st})
\]
is the process defined by $\tilde\jmath_{st}=j_{st}\circ \pi$ for $0\le s\le t$, where $\pi$ is the canonical homomorphism $\pi:\mathcal U(\eufrak{sl}_2\oplus_\alpha \eufrak{hw})\to\mathcal U\big((\eufrak{sl}_2\oplus_\alpha \eufrak{hw})/\eufrak{hw}\big)\cong\mathcal U(\eufrak{sl}_2)$.
\end{proposition}
\begin{proof}
We proceed as in the proof of Proposition \ref{prop no ext of hw}, we assume that $(\tilde\jmath_{st})_{0\le s\le t}$ is such an extension, and then we show that this necessarily implies $\tilde\rho|_{\mathcal{U}_0(\eufrak{hw})}=0$, $\tilde\eta|_{\mathcal{U}_0(\eufrak{hw})}=0$, and $\tilde{L}|_{\mathcal{U}_0(\eufrak{hw})}=0$ for its Sch{\"u}rmann triple $(\tilde\rho,\tilde\eta,\tilde{L})$. We know that the restriction of the Sch{\"u}rmann triple $(\tilde\rho,\tilde\eta,\tilde{L})$ to the subalgebra $\eufrak{sl}_2$ and the representation space $D=\tilde{\eta}\big(\mathcal{U}_0(\eufrak{sl_2})\big)$ has to be equivalent to the Sch{\"u}rmann triple $(\rho,\eta,L)$ defined in Example \ref{exa boukas}.

Our main tool are the following two facts, which can be deduced from our construction of the irreducible $*$-representations of $\eufrak{sl}_2$ in Subsection \ref{exa sl2}. Let $\pi$ be an arbitrary $*$-representation of $\eufrak{sl}_2$. Then $\pi(B^-)v=0$ and $\pi(M)v=\lambda v$, with $\lambda<0$ implies $v=0$. And if we have a vector $v\not=0$ that satisfies $\pi(B^-)v=0$ and $\pi(M)v=\lambda v$ with $\lambda\ge 0$, then $\pi$ restricted to $\pi\big(\mathcal{U}(\eufrak{sl_2})\big)v$ is equivalent to the lowest weight representation $\rho^+_{m_0}$ with $m_0=\lambda$.

First, we show in several steps that $\eta(B^+)$ is cyclic for $\tilde\rho$ and exhibit several vectors in $\tilde{D}=\tilde\eta\big(\mathcal{U}_0(\eufrak{sl}_2\oplus_\alpha \eufrak{hw})\big)$ which are lowest weight vectors for $\eufrak{sl}_2$. Using this information we can then prove that $\tilde\rho$, $\tilde\eta$, and $\tilde{L}$ vanish on $\eufrak{hw}$ (and therefore also on $\mathcal{U}_0(\eufrak{hw})$).

\begin{description}
\item[Step 1:] $\tilde\eta(A^-)=0$.

The relations $[B^-,A^-]=0$ and $[M,A^-]=-A^-$ imply $\tilde\rho(B^-)\tilde\eta(A^-)=\tilde\rho(A^-)\eta(B^-)=0$ and $-\tilde\eta(A^-)=\tilde\rho(M)\tilde\eta(A^-)-\tilde\rho(A^-)\eta(M)=\tilde\rho(M)\tilde\eta(A^-)$.

\item[Step 2:] If $u_0=\tilde\rho(A^-)\eta(B^+)=\tilde\eta(A^+)\not=0$, then it generates an $\eufrak{sl}_2$-representation that is equivalent to $\rho^+_1$.

Since $\tilde\eta(A^-)=0$, the relation $[A^-,B^+]=A^+$ implies $\tilde\eta(A^+)=\tilde\rho(A^-)\eta(B^+)-\tilde\rho(B^+)\tilde\eta(A^-)=\tilde\rho(A^-)\eta(B^+)$. Furthermore $[B^-,A^+]=A^-$ and $[M,A^+]=A^+$ yield $\tilde\rho(B^-)\tilde\eta(A^+)=\tilde\rho(A^+)\eta(B^-)+\tilde\eta(A^-)=0$ and  $\tilde\rho(M)\tilde\eta(A^+)=\tilde\rho(A^+)\eta(M)+\tilde\eta(A^+)=\tilde\eta(A^+)$.

\item[Step 3:] The  $\eufrak{sl}_2$-representation generated from $v_0=\tilde\rho(A^-)\tilde\eta(A^+)=\tilde\eta(E)$ is equivalent to the trivial one, i.e.\ $\tilde\rho(B^-)\tilde\eta(E)=\tilde\rho(M)\tilde\eta(E)=\tilde\rho(B^+)\tilde\eta(E)=0$.

We get  $\tilde\rho(B^-)\tilde\eta(E)=\tilde\rho(M)\tilde\eta(E)=0$ from the relations $[M,E]=0$ and $[B^-,E]=0$, and $\tilde\rho(B^+)\tilde\eta(E)=0$ follows from our basic facts on $\eufrak{sl}_2$-represen\-tations.

\item[Step 4:]  $\tilde\eta(E)=0$ and $w_0=\tilde\rho(A^+)\tilde\eta(A^+)$ is the lowest weight vector of an $\eufrak{sl}_2$-representation equivalent to $\rho^+_2$ (unless $w_0=0$).

Applying twice the relation $[B^-,A^+]=A^-$ and once $[A^-,A^+]=E$, we get
\begin{eqnarray*}
\tilde\rho(B^-)\tilde\rho(A^+)\tilde\eta(A^+) &=&\tilde\rho(A^+)\tilde\rho(B^-)\tilde\eta(A^+)+ \tilde\rho(A^-)\tilde\eta(A^+) \\
&=& \tilde\rho(A^+)\tilde\rho(A^+)\eta(B^-) + \tilde\rho(A^+)\tilde\eta(A^-)+ \tilde\rho(A^+)\tilde\eta(A^-) +\tilde\eta(E) \\
&=&\tilde\eta(E).
\end{eqnarray*}
We can use this relation to compute the norm of $\tilde\eta(E)$,
\[
||\tilde\eta(E)||^2= \langle\tilde\eta(E),\tilde\rho(B^-)\tilde\rho(A^+)\tilde\eta(A^+)\rangle = \langle\tilde\rho(B^+)\tilde\eta(E),\tilde\rho(A^+)\tilde\eta(A^+)\rangle=0,
\]
since $\tilde\rho(B^+)\tilde\eta(E)=0$.

Using twice the relation $[M,A^+]=A^+$, one also obtains $\tilde\rho(M)w_0=2w_0$.

\item[Step 5:] $\tilde\rho(E)=0$.

The results of steps 1, 2, and 4 and the surjectivity of $\tilde\eta$ imply that $\eta(B^+)$ is cyclic for $\tilde\rho$, i.e.\ any vector $v\in D$ can be written in the from $v=\tilde\rho(u)\eta(B^+)$ for some $u\in\mathcal{U}(\eufrak{sl}_2\oplus_\alpha \eufrak{hw})$. Since $E$ is central, we get
\[
\tilde\rho(E)v=\tilde\rho(E)\tilde\rho(u)\eta(B^+) =\tilde\rho(uB^+)\tilde\eta(E)=0.
\]
for all $v\in D$.

\item[Step 6:] $w_0=0$.

We can compute the norm of $\tilde\rho(B^+)w_0=\tilde\rho(B^+)\tilde\rho(A^+)\tilde\eta(A^+)$ in two different ways. Since $A^+$ and $B^+$ commute, we get
\begin{eqnarray*}
||\tilde\rho(B^+)w_0||^2 &=& ||\tilde\rho(A^+)^2\eta(B^+)||^2
=\langle\eta(B^+),\tilde\rho(A^-)^2\tilde\rho(A^+)^2\eta(B^+)\rangle \\
&=&\langle\tilde\rho(A^-)^2\eta(B^+),\tilde\rho(A^-)^2\eta(B^+)\rangle
=||\tilde\rho(A^-)\tilde\eta(A^+)||^2=||\tilde\eta(E)||^2=0,
\end{eqnarray*}
where we also used $\tilde\rho(E)=0$.

If $w_0\not=0$, then $\tilde\rho$ restricted to
$\tilde\rho(\mathcal{U}(\eufrak{sl}_2)w_0$ is equivalent to $\rho^+_2$,
so in particular the vectors $w_n=\tilde\rho(B^+)^n$, $n\ge 0$, must be an
orthogonal family of non-zero vectors with $||w_1||^2=6||w_0||^2$ by Equation (\ref{eq norms}). But we have
just shown $||w_1||^2=0$.

\item[Step 7:] $u_0=0$ and $\tilde\rho|_{\eufrak{hw}}=0$.

We get
\begin{eqnarray*}
||u_0||^2 &=& \langle\tilde\rho(A^-)\eta(B^+),\tilde\rho(A^-)\eta(B^+)\rangle =\langle\eta(B^+),\tilde\rho(A^+)\tilde\rho(A^-)\eta(B^+)\rangle \\
&=&\langle\eta(B^+),\tilde\rho(A^+)\tilde\eta(A^+)\rangle = \langle\eta(B^+),w_0\rangle =0.
\end{eqnarray*}

Therefore we have $\tilde\eta|_{\eufrak{hw}}=0$ and $\tilde D=D={\rm span}\,\left\{\eta\big((B^+)^k\big)| k=1,2,\ldots\right\}$. From this we can deduce $\tilde\rho(A^+)\eta\big((B^+)^k\big)=\tilde\rho((B^+)^k\big)\tilde\eta(A^+)=0$, i.e.\ $\tilde\rho(A^+)=0$ and therefore also $\tilde\rho(A^-)=\tilde\rho(A^+)^*=0$.

\item[Step 8:] $\tilde{L}|_{\eufrak{hw}}=0$.

Finally, using, e.g., the relations $[M,A^\pm]=\pm A^\pm$ and $E=[A^-,A^+]$, one can show that the generator $\tilde{L}$ also vanishes on $\eufrak{hw}$,
\begin{eqnarray*}
\pm\tilde{L}(A^\pm) &=& \langle\eta(M),\tilde\eta(A^\pm)\rangle -\langle\tilde\eta(A^\mp),\eta(M)\rangle=0, \\
\tilde{L}(E) &=& ||\tilde\eta(A^+)||^2-||\tilde\eta(A^-)||^2=0.
\end{eqnarray*}

\end{description}
\end{proof}

But there do exist non-trivial L\'evy processes such that $j_{st}(A^-)\Omega=j_{st}(B^-)\Omega=0$ for all $0\le s\le t$, as the following example shows.

\begin{example}
Let $h>0$ and let $\rho_h$ be the Fock representation of $\mathcal U(\eufrak{hw})$ defined in (\ref{fock}). This extends to a representation of $\mathcal U(\eufrak{sl}_2\oplus_\alpha \eufrak{hw})$, if we set
\[
\rho_h(B^+)=\frac{\rho_h(A^+)^2}{2h},\quad \rho_h(M)=\frac{\rho_h(A^+A^-+A^-A^+)}{2h}, \quad \rho_h(B^+)=\frac{\rho_h(A^-)^2}{2h}.
\]
The restriction of this representation to $\eufrak{sl}_2$ is a direct sum of the two lowest weight representations $\rho^+_{1/2}$ and $\rho^+_{3/2}$, the respective lowest weight vectors are $e_0$ and $e_1$.
For the cocycle we take the coboundary of the `lowest weight vector' $e_0\in\ell^2$, i.e.\ we set
\[
\eta(u)=\rho_h(u)e_0
\]
for $u\in\mathcal U_0(\eufrak{sl}_2\oplus_\alpha \eufrak{hw})$, and for the generator
\[
L(u)=\langle e_0,\rho_h(u)e_0\rangle
\]
for $u\in\mathcal U_0(\eufrak{sl}_2\oplus_\alpha \eufrak{hw})$. This defines a Sch\"urmann triple on $\eufrak{sl}_2\oplus_\alpha \eufrak{hw}$ over $\ell^2$ and therefore
\begin{eqnarray*}
{\rm d}_L B^+ &=& \frac{1}{2h}{\rm d}\Lambda\left(\rho_h(A^+)^2\right) + \frac{1}{\sqrt 2}{\rm d}A^*(e_2), \\
{\rm d}_L A^+ &=& {\rm d}\Lambda\big(\rho_h(A^+)\big)+\sqrt{h}{\rm d}A^*(e_1), \\
{\rm d}_L M &=& \frac{1}{2h}{\rm d}\Lambda\big(\rho_h(A^+A^-+A^-A^+)\big)+ \frac{1}{2}{\rm d}A^*(e_0) + \frac{1}{2}{\rm d}A(e_0) + \frac{1}{2}{\rm d}t, \\
{\rm d}_L E &=& h {\rm d}\Lambda({\rm id}) + h{\rm d}A^*(e_0)+h{\rm d}(e_0)+h{\rm d}t, \\
{\rm d}_L A^- &=&{\rm d}\Lambda\big(\rho_h(A^-)\big)+\sqrt{h}{\rm d}A(e_1), \\
{\rm d}_L B^- &=& \frac{1}{2h}{\rm d}\Lambda\left(\rho_h(A^-)^2\right) +\frac{1}{\sqrt 2}{\rm d}A(e_2), 
\end{eqnarray*}
defines a L\'evy process $\eufrak{sl}_2\oplus_\alpha \eufrak{hw}$, acting on the Fock space over $L^2(\mathbb R_+,\ell^2)$. The It\^o table of this process is infinite-dimensional. The restriction of this process to $\eufrak{sl}_2$ is equivalent to the process defined in Example \ref{exa PaSi} with $m_0=\frac{1}{2}$.

One can easily verify that $j_{st}(A^-)$ and $j_{st}(B^-)$ annihilate the vacuum vector of $\Gamma\big(L^2(\mathbb R_+,\ell^2)\big)$.

We have $\rho_h(C)=- \frac{3}{4}{\rm id}$, and therefore
\[
{\rm d}_L C = - \frac{3}{4} \big({\rm d}\Lambda({\rm id})+ {\rm d}A^*(e_0)+{\rm d}A(e_0)+{\rm d}t\big).
\]
\end{example}

\subsection{Higher order noises}\label{exa higher}

Let us now consider the infinite-dimensional real Lie algebra $\eufrak{wn}$ that is spanned by $\{B_{n,m};n,m\in\mathbb N\}$ with the commutation relations obtained by the natural extension, to
higher powers of the white noise, of the renormalization rule introduced
in \cite{accardi+lu+volovich99}, i.e.:

\begin{eqnarray*}
[B_{n_1,m_1},B_{n_2,m_2}] &=& \sum_{k=1}^{n_2\wedge m_1} \frac{m_1!n_2!}{(m_1-k)!(n_2-k)!k!}c^k B_{n_1+n_2-k,m_1+m_2-k} \\
&&- \sum_{k=1}^{n_1\wedge m_2} \frac{m_2!n_1!}{(m_2-k)!(n_1-k)!k!}c^k B_{n_1+n_2-k,m_1+m_2-k}
\end{eqnarray*}
for $n_1,n_2,m_1,m_2\in\mathbb N$, and involution $\left(B_{n,m}\right)^*=B_{m,n}$, where $c\ge 0$ is some fixed positive parameter. These relations can be obtained by taking the quotient of the universal enveloping algebra $\mathcal U(\eufrak{hw})$ of $\eufrak{hw}$ with respect to the ideal generated by $E=c\mathbf 1$. The basis elements $B_{n,m}$ are the images of $(A^+)^n(A^-)^m$.

We can embed $\eufrak{hw}$ and $\eufrak{sl}_2\oplus_\alpha \eufrak{hw}$ into $\eufrak{wn}$ by
\begin{gather*}
A^+\mapsto \frac{B_{1,0}}{\sqrt{c}}, \quad A^- \mapsto \frac{B_{0,1}}{\sqrt{c}}, \quad E\mapsto B_{0,0}, \\
B^+\mapsto \frac{1}{2c}B_{2,0},\quad B^-\mapsto \frac{1}{2c}B_{0,2}, \quad M\mapsto \frac{1}{c}B_{1,1}+\frac{1}{2}B_{0,0}.
\end{gather*}

There exist no Gaussian L\'evy processes on $\eufrak{wn}$, since $[\eufrak{wn},\eufrak{wn}]=\eufrak{wn}$.

Let $\rho_c$ be the Fock representation defined in Equation (\ref{fock}). Setting
\[
\rho(B_{n,m})=\rho_c\left((A^+)^n(A^-)^m\right), \qquad n,m\in\mathbb{N},
\]
we get a $*$-representation of $\mathcal{U}(\eufrak{wn})$. If we set $\eta(u)=\rho(u)e_0$ and $L(u)=\langle e_0,\rho(u)e_0\rangle$ for $u\in\mathcal{U}_0(\eufrak{wn})$, then we obtain a Sch\"urmann triple  on $\eufrak{wn}$. For this triple we get
\[
{\rm d}_LB_{n,m}={\rm d}\Lambda\left(\rho_c(A^+)^n\rho_c(A^-)^m\right)+\delta_{m0}\sqrt{c^nn!\,}{\rm d}A^*(e_n)+\delta_{n0}\sqrt{c^mm!\,}{\rm d}A(e_m)+\delta_{n0}\delta_{m0}{\rm d}t,
\]
for the differentials. Note that we have $j_{st}(B_{nm})\Omega=0$ for all $m\ge 1$ and $0\le s\le t$ for the associated L\'evy process.

\subsection{Other examples: L\'evy processes on $\eufrak{fd}$ and $\eufrak{gl}_2$}

The goal of this subsection is to explain the relation of the present paper to previous works by Boukas \cite{boukas88,boukas91a} and Parthasarathy and Sinha \cite{parthasarathy+sinha91}.

We introduce the two real Lie algebras $\eufrak{\eufrak{fd}}$ and $\eufrak{gl}_2$.
The finite-difference Lie algebra $\eufrak{fd}$ is the three-dimensional solvable real Lie algebra with basis $\{P,Q,T\}$, commutation relations
\[
[P,Q]=[T,Q]=[P,T]=T,
\]
and involution $P^*=Q$, $T^*=T$, cf.\ \cite{feinsilver87b}. This Lie algebra is actually the direct sum of the unique non-abelian two-dimensional real Lie algebra and the one-dimensional abelian Lie algebra, its center is spanned by $T-P-Q$.

The Lie algebra $\eufrak{gl}_2$ of the general linear group $GL(2;\mathbb R)$ is the direct sum of $\eufrak{sl}_2$ with the one-dimensional abelian Lie algebra. As a basis of $\eufrak{gl}_2$ we will choose $\{B^+,B^-,M,I\}$, where $B^+$, $B^-$, and $M$ are a basis of the Lie subalgebra $\eufrak{sl}_2$, and $I$ is hermitian and central. Note that $T\mapsto M+B^++B^-$, $P\mapsto (M-I)/2 + B^-$, $Q\mapsto (M-I)/2+B^+$ defines an injective Lie algebra homomorphism from $\eufrak{fd}$ into $\eufrak{gl}_2$, i.e.\ we can regard $\eufrak{fd}$ as a Lie subalgebra of $\eufrak{gl}_2$.

Following ideas by Feinsilver \cite{feinsilver89}, Boukas \cite{boukas88,boukas91a} constructed a calculus for $\eufrak{fd}$, i.e.\ he constructed a L\'evy process on it and defined stochastic integrals with respect to it. He also derived the It\^o formula for these processes and showed that their It\^o table is infinite-dimensional. His realization is not defined on the boson Fock space, but on the so-called finite-difference Fock space especially constructed for his $\eufrak{fd}$ calculus. Parthasarathy and Sinha constructed another L\'evy process on $\eufrak{fd}$, acting on a boson Fock space, in \cite{parthasarathy+sinha91}. They gave an explicit decomposition of the operators into conservation, creation, annihilation, and time, thereby reducing its calculus to Hudson-Parthasarathy calculus.

Accardi and Skeide \cite{accardi+skeide99a,accardi+skeide99b} noted that they were able to recover Boukas' $\eufrak{fd}$ calculus from their SWN calculus. In fact, since $\eufrak{gl}_2$ is a direct sum of $\eufrak{sl}_2$ and the one-dimensional abelian Lie algebra, any L\'evy process $(j_{st})_{0\le s\le t}$ on $\eufrak{sl}_2$ can be extended (in many different ways) to a L\'evy process $(\tilde\jmath_{st})_{0\le s\le t}$ on $\eufrak{gl}_2$. We will only consider the extensions defined by
\[
\tilde\jmath_{st}|_{\eufrak{sl}_2} = j_{st}, \quad\mbox{ and }\quad \tilde\jmath_{st}(I) = \lambda (t-s){\rm id}, \qquad \mbox{ for } 0\le s\le t,
\]
for $\lambda\in\mathbb R$.
Since $\eufrak{fd}$ is a Lie subalgebra of $\eufrak{gl}_2$, we also get a L\'evy process on $\eufrak{fd}$ by restricting $(\tilde\jmath_{st})_{0\le s\le t}$ to $\mathcal U(\eufrak{fd})$.

If we take the L\'evy process on $\eufrak{sl}_2$ defined in Example \ref{exa boukas} and $\lambda=1$, then we get
\begin{eqnarray*}
{\rm d}_LP &=& {\rm d}\Lambda\big(\rho^+(M/2+B^-)\big) + {\rm d}A(e_0), \\
{\rm d}_LQ &=& {\rm d}\Lambda\big(\rho^+(M/2+B^+)\big) + {\rm d}A^*(e_0), \\
{\rm d}_LT &=& {\rm d}\Lambda\big(\rho^+(M+B^++B^-)\big) + {\rm d}A^*(e_0) + {\rm d}A(e_0) + {\rm d}t.
\end{eqnarray*}
It can be checked that this L\'evy process is equivalent to the one defined by Boukas.

If we take instead the L\'evy process on $\eufrak{sl}_2$ defined in Example \ref{exa PaSi}, then we get
\begin{eqnarray*}
{\rm d}_LP &=& {\rm d}\Lambda\big(\rho^+_{m_0}(M/2+B^-)\big) + {\rm d}A^*\left(\frac{m_0}{2}e_0\right) +{\rm d}A\left(\frac{m_0}{2}e_0+\sqrt{m_0}e_1\right) + \frac{m_0-\lambda}{2}{\rm d}t , \\
{\rm d}_LQ &=& {\rm d}\Lambda\big(\rho^+_{m_0}(M/2+B^+)\big) +{\rm d}A^*\left(\frac{m_0}{2}e_0+\sqrt{m_0}e_1\right) + {\rm d}A\left(\frac{m_0}{2}e_0\right) + \frac{m_0-\lambda}{2}{\rm d}t , \\ 
{\rm d}_LT &=& {\rm d}\Lambda\big(\rho^+_{m_0}(M+B^++B^-)\big) + {\rm d}A^*(m_0e_0+\sqrt{m_0}e_1) + {\rm d}A(m_0e_0+\sqrt{m_0}e_1) +m_0 {\rm d}t \\
&=& {\rm d}_LP + {\rm d}_LQ + \lambda {\rm d}t.
\end{eqnarray*}
For $m_0=\lambda=2$, this is exactly the L\'evy process defined in \cite{parthasarathy+sinha91}. Note that in that case the representation $\rho^+_2=\rho^+$ and the Fock space agree with those of Boukas' process, but the cocycle and the generator
are different. Therefore the construction of \cite{parthasarathy+sinha91} leads to the same
algebra as Boukas', but not to the same quantum process -- a fact that
had already been noticed by Accardi and Boukas \cite{accardi+boukas02}.

\section{Classical processes}\label{class}

Let $(j_{st})_{0\le s\le t}$ be a L\'evy process on a real Lie algebra $\eufrak{g}_\mathbb{R}$ over $\Gamma=\Gamma\big(L^2(\mathbb{R}_+,D)\big)$, fix a hermitian element $Y$, $Y^*=Y$, of $\eufrak g_\mathbb R$, and define a map $y:\Sigma(\mathbb R_+)\to \mathcal L(\Gamma)$ by
\[
y_\phi = \sum_{k=1}^n\phi_kj_{s_kt_k}(Y), \qquad \mbox{ for }\phi= \sum_{k=1}^n\phi_k\mathbf{1}_{[s_k,t_k[}\in\Sigma(\mathbb R_+)
\]
It is clear that the operators $\{y_\phi;\phi\in \Sigma(\mathbb R_+)\}$ commute, since $y$ is the restriction of $\pi:\eufrak{g}^{\mathbb{R}_+}\ni\psi=\sum_{k=1}^n\psi_k\mathbf{1}_{[s_k,t_k[}\mapsto \sum_{k=1}^nj_{s_kt_k}(\psi_k)\in\mathcal{L}(\Gamma)$ to the abelian current algebra $\mathbb C Y^{\mathbb R_+}$ over $\mathbb C Y$. Furthermore, if $\phi$ is real-valued, then $y_\phi$ is hermitian, since $Y$ is hermitian. Therefore there exists a classical stochastic process $(\tilde Y_t)_{t\ge 0}$ whose moments are given by
\[
\mathbb{E}(\tilde{Y}_{t_1}\cdots \tilde{Y}_{t_n})=\langle\Omega, y_{\mathbf 1_{[0,t_1[}}\cdots  y_{\mathbf 1_{[0,t_n[}}\Omega\rangle,\qquad \mbox{ for all } t_1,\ldots, t_n\in\mathbb{R}_+.
\]
Since the expectations of $(j_{st})_{0\le s\le t}$ factorize, we can choose $(\tilde{Y}_t)_{t\ge 0}$ to be a L\'evy process. If $j_{st}(Y)$ is even essentially self-adjoint, then the marginal distribution of $(\tilde{Y}_t)_{t\ge 0}$ is uniquely determined.

We will now give a characterization of $(\tilde{Y}_t)_{t\ge 0}$. First, we need two lemmas.

\begin{lemma}\label{lem class 1}
Let $X\in \mathcal{L}(D)$, $u,v\in D$, and suppose furthermore that the series $\sum_{n=0}^\infty \frac{(tX)^n}{n!} w$ and $\sum_{n=0}^\infty \frac{(tX^*)^n}{n!} w$ converge in $D$ for all $w\in D$. Then we have
\begin{eqnarray*}
e^{\Lambda(X)}A(v) &=& A\left(e^{-X^*}v\right)e^{\Lambda(X)} \\
e^{A^*(u)}A(v) &=& \big(A(v)-\langle v,u\rangle\big) e^{A^*(u)} \\
e^{A^*(u)}\Lambda(X) &=& \big(\Lambda(X)-A^*(Xu)\big)e^{A^*(u)}
\end{eqnarray*}
on the algebraic boson Fock space over $D$.
\end{lemma}
\begin{proof}
This can be deduced from the formula for the adjoint actions, ${\rm Ad}e^{X}Y=e^XYe^{-X}=Y+[X,Y]+\frac{1}{2}\big[X,[X,Y]\big]+\cdots=e^{{\rm ad}X}Y$.
\end{proof}

The following formula gives the normally ordered form of the generalized Weyl operators and is a key tool to calculate the characteristic functions of classical subprocesses of L\'evy processes on real Lie algebras.

\begin{lemma}\label{lem class 2}
Let $X\in\mathcal{L}(D)$ and $u,v\in D$ and suppose furthermore that the series $\sum_{n=0}^\infty \frac{(tX)^n}{n!} w$ and $\sum_{n=0}^\infty \frac{(tX^*)^n}{n!} w$ converge in $D$ for all $w\in D$. Then we have
\[
\exp\big(\Lambda(X)+A^*(u)+A(v)+\alpha\big) = \exp\big(A^*(\tilde{u})\big)\exp\big(\Lambda(X)\big)\exp\big(A(\tilde{v})\big)\exp(\tilde\alpha)
\]
on the algebraic boson Fock space over $D$, where
\[
\tilde{u} = \sum_{n=1}^\infty \frac{X^{n-1}}{n!}u,\quad
\tilde{v} = \sum_{n=1}^\infty \frac{(X^*)^{n-1}}{n!}v,\quad
\tilde{\alpha} = \alpha + \sum_{n=2}^\infty \langle v,\frac{X^{n-2}}{n!} u\rangle.
\]
\end{lemma}
\begin{proof}
Let $\omega\in D$ and set $\omega_1(t)=\exp t\big(\Lambda(X)+A^*(u)+A(v)+\alpha\big)\omega$ and
\[
\omega_2(t)= \exp\Big(A^*\big(\tilde{u}(t)\big)\Big)\exp\big(t\Lambda(X)\big)\exp\Big(A\big(\tilde{v}(t)\big)\Big)\exp\big(\tilde\alpha(t)\big)\omega
\]
for $t\in[0,1]$, where
\[
\tilde{u}(t) = \sum_{n=1}^\infty \frac{t^nX^{n-1}}{n!}u,\,
\tilde{v}(t) = \sum_{n=1}^\infty \frac{t^n(X^*)^{n-1}}{n!}v,\,
\tilde{\alpha}(t) = t\alpha +\sum_{n=2}^\infty \langle v,\frac{t^nX^{n-2}}{n!} u\rangle .
\]
Then we have $\omega_1(0)=\omega=\omega_2(0)$. Using Lemma \ref{lem class 1}, we can also check that
\[
\frac{{\rm d}}{{\rm d}t} \omega_1(t) = \big(\Lambda(X)+A^*(u)+A(v)+\alpha\big)\omega\exp t\big(\Lambda(X)+A^*(u)+A(v)+\alpha\big)\omega
\]
and
\begin{eqnarray*}
\frac{{\rm d}}{{\rm d}t} \omega_2(t) &=& A^*\left(\frac{{\rm d}\tilde{u}}{{\rm d}t}(t)\right)\exp\Big(A^*\big(\tilde{u}(t)\big)\Big)\exp\big(t\Lambda(X)\big)\exp\Big(A\big(\tilde{v}(t)\big)\Big)\exp\big(\tilde\alpha(t)\big)\omega \\
&&+\exp\Big(A^*\big(\tilde{u}(t)\big)\Big)\Lambda(X)\exp\big(t\Lambda(X)\big)\exp\Big(A\big(\tilde{v}(t)\big)\Big)\exp\big(\tilde\alpha(t)\big)\omega\\
&&+ \exp\Big(A^*\big(\tilde{u}(t)\big)\Big)\exp\big(t\Lambda(X)\big)A\left(\frac{{\rm d}\tilde{v}}{{\rm d}t}(t)\right)\exp\Big(A\big(\tilde{v}(t)\big)\Big)\exp\big(\tilde\alpha(t)\big)\omega \\
&&+ \exp\Big(A^*\big(\tilde{u}(t)\big)\Big)\exp\big(t\Lambda(X)\big)\exp\Big(A\big(\tilde{v}(t)\big)\Big)\frac{{\rm d}\tilde\alpha}{{\rm d}t}(t)\exp\big(\tilde\alpha(t)\big)\omega
\end{eqnarray*}
coincide for all $t\in[0,1]$. Therefore we have $\omega_1(1)=\omega_2(1)$.
\end{proof}

\begin{theorem}\label{levy-lie-class}
Let $(j_{st})_{0\le s\le t}$ be a L\'evy process on a real Lie algebra $\eufrak g_\mathbb R$ with Sch{\"u}rmann triple $(\rho,\eta,L)$. Then for any hermitian element
$Y$ of $\eufrak g_\mathbb R$ such that $\eta(Y)$ is analytic for $\rho(Y)$,
the associated classical L\'evy process $(\tilde{Y}_t)_{t\ge 0}$ has
characteristic exponent
\[
\Psi(\lambda)= i\lambda L(Y)+\sum_{n=2}^\infty \frac{(i\lambda)^n}{n!}\langle\eta(Y^*),\rho(Y)^{n-2}\eta(Y)\rangle,
\]
($\rho(Y)^0={\rm id}$) for $\lambda$ in some neighborhood of zero.
\end{theorem}
\begin{proof}
The characteristic exponent $\Psi(\lambda)$, $\lambda\in\mathbb{R}$, is defined by $\mathbb{E}(e^{i\lambda\tilde{Y}_t})=e^{t\Psi(\lambda)}$, so we have to compute
\[
\mathbb{E}\left(e^{i\lambda\tilde{Y}_t}\right) = \langle\Omega,e^{i\lambda j_{0t}(Y)}\Omega\rangle
\]
for $j_{0t}(Y)=\Lambda_{0t}\big(\rho(Y)\big)+A^*_{0t}\big(\eta(Y)\big)+A_{0t}\big(\eta(Y)\big) + tL(Y)$. Using Lemma \ref{lem class 2}, we get
\[
\mathbb{E}\left(e^{i\lambda\tilde{Y}_t}\right)= \exp\left( it\lambda L(Y)+t\sum_{n=2}^{\infty}\left\langle\eta(Y^*),\frac{(i\lambda)^n\rho(Y)^{n-2}}{n!}\eta(Y)\right\rangle\right).
\]
\end{proof}
\begin{remark}
Note that $\Psi(\lambda)$ is nothing else than $\sum_{n=1}^\infty \frac{(i\lambda)^n}{n!}L(Y^n)$. It is also possible to give a more direct proof of the theorem, using the convolution of functionals on $\mathcal{U}(\eufrak{g})$ instead of the boson Fock space realization of $(j_{st})_{0\le s\le t}$.
\end{remark}

We give two corollaries of this result, the first justifies our definition of Gaussian generators.

\begin{corollary}\label{class cor 1}
Let $L$ be a Gaussian generator on $\eufrak{g}_\mathbb{R}$ with corresponding L\'evy process $(j_{st})_{0\le s\le t}$. Then for any hermitian element $Y$ the associated classical L\'evy process $(\tilde{Y}_t)_{t\ge 0}$ is Gaussian with mean and variance
\[
\mathbb E(\tilde{Y}_t)=tL(Y), \qquad \mathbb E(\tilde{Y}^2_t)=||\eta(Y)||^2t,\qquad\mbox{ for }t\ge 0.
\]
\end{corollary}

We see that in this case we can take $\big(||\eta(Y)||B_t+L(Y)t\big)_{t\ge 0}$ for $(\tilde{Y}_t)_{t\ge 0}$, where $(B_t)_{t\ge 0}$ is a standard Brownian motion.

The next corollary deals with the case where $L$ is the restriction to $\mathcal{U}_0(\eufrak{g})$ of a positive functional on $\mathcal{U}(\eufrak{g})$.

\begin{corollary}\label{class cor 2}
Let $(\rho,\eta,L)$ be a Sch\"urmann triple on $\eufrak{g}_\mathbb{R}$ whose cocycle is trivial, i.e.\ there exists a vector $\omega\in D$ such that $\eta(u)=\rho(u)\omega$ for all $u\in\mathcal{U}_0(\eufrak{g})$, and whose generator is of the form $L(u)=\langle \omega,\rho(u)\omega\rangle$, for all $u\in\mathcal{U}_0(\eufrak{g})$. Suppose furthermore that the vector $\omega$ is analytical for $\rho(Y)$, i.e.\ that $e^{u\rho(Y)}\omega:=\sum_{n=1}^\infty \frac{u^n\rho(Y)^n}{n!}\omega$ converges for sufficiently small $u$. Then the classical stochastic process $(\tilde{Y}_t)_{t\ge 0}$ associated to $(j_{st})_{0\le s\le t}$ and $Y$ is a compound Poisson process with characteristic exponent
\[
\Psi(u)=\langle\omega,\left(e^{iu\rho(Y)}-1\right)\omega\rangle.
\]
\end{corollary}
\begin{remark}
If the operator $\rho(Y)$ is even (essentially) self-adjoint, then we get the L\'evy measure of $(\tilde{Y}_t)_{t\ge 0}$ by evaluating its spectral measure in the state vector $\omega$,
\[
\mu({\rm d}\lambda)=\langle\omega,{\rm d}\mathbb{P}_\lambda\omega\rangle,
\]
where $\rho(Y)=\int \lambda {\rm d}\mathbb{P}_\lambda$ is the spectral resolution of (the closure of) $\rho(Y)$. 
\end{remark}

Corollary \ref{class cor 2} suggests to call a L\'evy process on $\eufrak{g}$ with trivial cocycle $\eta(u)=\rho(u)\omega$ and generator $L(u)=\langle\omega,\rho(u)\omega\rangle$ for $u\in\mathcal{U}_0(\eufrak{g})$ a {\em Poisson process} on $\eufrak{g}$.

\begin{example}
Let $(j_{st})_{0\le s\le t}$ be the L\'evy process on $\eufrak{sl}_2$ defined in Example \ref{exa PaSi} and let $Y=B^+ + B^-+\beta M$ with $\beta\in\mathbb{R}$. The operator $X=\rho^+_{m_0}(Y)$ is essentially self-adjoint. We now want to characterize the classical L\'evy process $(\tilde{Y}_t)_{t\ge 0}$ associated to $Y$ and $(j_{st})_{0\le s\le t}$ in the manner described above. Corollary \ref{class cor 2} tells us that $(\tilde{Y}_t)_{t\ge 0}$ is a compound Poisson process with characteristic exponent
\[
\Psi(u)=\langle e_0,\left(e^{iuX}-1\right)e_0\rangle.
\]
We want to determine the L\'evy measure of $(\tilde{Y}_t)_{t\ge 0}$, i.e.\ we want to determine the measure $\mu$ on $\mathbb{R}$, for which
\[
\Psi(u)=\int \left( e^{iux}-1\right)\mu({\rm d}x).
\]
This is the spectral measure of $X$ evaluated in the state $\langle e_0,\cdot\,e_0\rangle$. Note that the polynomials $p_n\in\mathbb{R}[x]$ defined by the condition
\[
e_n=p_n(X)e_0,
\]
$n=0,1,\ldots,$ are orthogonal w.r.t.\ $\mu$, since
\[
\int p_n(x)p_m(x)\mu({\rm d}x) = \langle e_0,p_n(X)p_m(X)e_0\rangle = \langle p_n(X)e_0,p_m(X)e_0\rangle =\delta_{nm},
\]
for $n,m\in\mathbb{N}$. Looking at the definition of $X$, we can easily identify the three-term-recurrence relation satisfied by the $p_n$. We get
\[
X e_n = \sqrt{(n+1)(n+m_0)} e_{n+1} + \beta(2n+m_0) e_n + \sqrt{n(n+m_0-1)}e_{n-1},
\]
for $n\in \mathbb{N}$, and therefore
\[
(n+1)P_{n+1} + (2\beta n+\beta m_0-x)P_n + (n+m_0-1) P_{n-1} = 0,
\]
with initial condition $P_{-1}=0$, $P_1=1$, for the rescaled polynomials
\[
P_n=\prod^n_{k=1} \sqrt{\frac{n}{n+m_0}}\, p_n.
\]
According to the value of $\beta$ we have to distinguish three cases.
\begin{enumerate}
\item
$|\beta|=1$: In this case we have, up to rescaling, Laguerre polynomials, i.e.\
\[
P_n(x)=(-\beta)^n L^{(m_0-1)}_n(\beta x)
\]
where the Laguerre polynomials $L^{(\alpha)}_n$ are defined as in \cite[Equation (1.11.1)]{koekoek+swarttouw94}. The measure $\mu$ can be obtained by normalizing the measure of orthogonality of the Laguerre polynomials, it is equal to
\[
\mu({\rm d}x) = \frac{|x|^{m_0-1}}{\Gamma(m_0)} e^{-\beta x}\mathbf{1}_{\beta\mathbb{R}_+}{\rm d}x.
\]
If $\beta=+1$, then this measure is, up to a normalization parameter,
the usual $\chi^2$-distribution (with parameter $m_0$) of probability
theory.
The operator $X$ is then positive and therefore $(\tilde{Y}_t)_{t\ge 0}$ is a subordinator, i.e.\ a L\'evy process with values in $\mathbb{R}_+$, or, equivalently, a L\'evy process with non-decreasing sample paths.
\item
$|\beta|<1$: In this case we find the Meixner-Pollaczek polynomials after rescaling,
\[
P_n(x)= P_n^{(m_0/2)}\left(\frac{x}{2\sqrt{1-\beta^2}};\pi-\arccos \beta\right).
\]
For the definition of these polynomials see, e.g., \cite[Equation (1.7.1)]{koekoek+swarttouw94}. For the measure $\mu$ we get
\[
\mu({\rm d}x) = C \exp\left(\frac{(\pi-2\arccos \beta)x}{2\sqrt{1-\beta^2}}\right)\left|\Gamma\left(\frac{m_0}{2}+\frac{ix}{2\sqrt{1-\beta^2}}\right)\right|^2{\rm d}x,
\]
where $C$ has to be chosen such that $\mu$ is a probability measure.
\item
$|\beta|>1$: In this case we get Meixner polynomials after rescaling,
\[
P_n(x)=
(-c\,{\rm sgn}\beta)^n\prod_{k=1}^n \frac{k+m_0-1}{k}M_n\left(\frac{x}{1/c-c}{\rm sgn}\beta-\frac{m_0}{2};m_0;c^2\right)
\]
where
\[
c=|\beta|-\sqrt{\beta^2-1}
\]
\end{enumerate}
The definition of these polynomials can be found, e.g., in \cite[Equation (1.9.1)]{koekoek+swarttouw94}. The density $\mu$ is again the measure of orthogonality of the polynomials $P_n$ (normalized to a probability measure). We therefore get
\[
\mu= C\sum_{n=0}^\infty \frac{c^{2n}(m_0)_n}{n!} \delta_{x_n}
\]
where
\[
x_n=\left(n+\frac{m_0}{2}\right)\left(\frac{1}{c}-c\right){\rm sgn}\beta, \qquad \mbox{ for } n\in\mathbb{N}
\]
and $C^{-1}=\sum_{n=0}^\infty \frac{c^{2n}(m_0)_n}{n!}=(1-c^2)^{-m_0}$. Here $(m_0)_n$ denotes the Pochhammer symbol, $(m_0)_n=m_0(m_0+1)\cdots(m_0+n-1)$.
\end{example}

\begin{example}
Let now $(j_{st})_{0\le s\le t}$ be the L\'evy process on $\eufrak{sl}_2$ defined in Example \ref{exa boukas} and let again $Y=B^+ +B^-+\beta M$ with $\beta\in\mathbb{R}$. We already noted in Example \ref{exa boukas} that $j_{st}$ is equivalent to $\rho^+_{t-s}$ for fixed $s$ and $t$. Therefore the marginal distributions of the classical L\'evy process $(\tilde{Y}_t)_{t\ge 0}$ are exactly the distributions of the operator $X$ that we computed in the previous example (with $m_0=t$).

For $\beta=1$, we recover \cite[Theorem 2.2]{boukas91a}. The classical L\'evy process associated to $T=B^++B^-+M$ is an exponential or Gamma process with Fourier transform
\[
\mathbb{E}\left(e^{iu \tilde{Y}_t}\right) = (1-iu)^{-t}
\]
and marginal distribution
$\nu_t({\rm d}x)=\frac{x^{t-1}}{\Gamma(t)} e^{- x}\mathbf{1}_{\mathbb{R}_+}{\rm d}x$. This is a subordinator with L\'evy measure is $x^{-1}e^{-x}\mathbf{1}_{\mathbb{R}_+}{\rm d}x$, see, e.g., \cite{bertoin96}.

For $\beta>1$, we can write the Fourier transform of the marginal distributions $\nu_t$ as
\[
\mathbb{E}(e^{iu\tilde{Y}_t})=\exp t\left(\frac{iu(c-1/c)}{2} + \sum_{n=1}^\infty \frac{c^{2n}}{n}\left(e^{iun(c-1/c)}-1\right)\right).
\]
This shows that we can define $(\tilde{Y}_t)_{t\ge 0}$ as sum of Poisson processes with a drift, i.e.\ if $\left(\big(N^{(n)}_t\big)_{t\ge 0}\right)_{n\ge 1}$ are independent Poisson processes (with intensity and jump size equal to one), then we can take
\[
\tilde{Y}_t=(c-1/c)\left(\sum_{n=1}^\infty n N^{(n)}_{c^{2n}t/n} + \frac{t}{2}\right), \qquad\mbox{ for } t\ge 0.
\]

The marginal distributions of these processes for the different values of $\beta$ and their relation to orthogonal polynomials are also discussed in \cite[Chapter 5]{feinsilver+schott93}.
\end{example}

\section{Conclusion}

We have shown that the theories of factorizable current representations of Lie algebras and L\'evy processes on $*$-bialgebras provide an elegant and efficient formalism for defining and studying quantum stochastic calculi with respect to additive operator processes satisfying Lie algebraic relations. The theory of L\'evy processes on $*$-bialgebras can also handle processes whose increments are not simply additive, but are composed by more complicated formulas, the main restriction is that they are independent (in the tensor sense). This allows to answer questions that could not be handled by direct computational methods, such as the computation of the SWN It\^o table, the simultaneous realization of linear and squared white noise on the same Hilbert space, or the characterization of the associated classical processes.

After the completion of the present article, Accardi, Hida, and Kuo \cite{accardi+hida+kuo01} have shown that using white noise calculus it is possible to obtain a closed It\^o table for the quadratic covariations of the three basic square of white noise operators. But the coefficients in their It\^o table contain functions of the Hida derivative and its adjoint.

%% file: braided.tex
\chapter{L{\'e}vy Processes and Brownian Motion on Braided Spaces}\label{braided}

\vspace*{3cm}

\begin{quote}
We define and study L\'evy processes and Brownian motion on braided spaces. First, we introduce the necessary notions and generalize some fundamental results, in particular about the positivity of the convolution of positive functionals and the one-to-one correspondences between generators, convolution semi-groups and L{\'e}vy processes. Next we exhibit a construction of examples of braided $*$-Hopf algebras, and also give their symmetrization. Finally, we study Brownian motion, i.e.\ L{\'e}vy processes with quadratic generators, on these algebras. We show how Brownian motions on these braided $*$-spaces can be classified and we give the quantum stochastic differential equations that define a realization of these processes on Fock space.
\end{quote}

\vspace*{3cm}

Joint work with Ren\'e Schott and Michael Sch\"urmann. 

\newpage

\section{Introduction}

L\'evy processes and especially Brownian motion are important examples of stochastic processes and they have been crucial in the development of the stochastic calculus and the theory of Markov processes. Stochastic integrals were first defined for Brownian motion, before their theory was extended to semi-martingales. Many important properties of Markov processes were first discovered for L\'evy processes.

In quantum probability, L{\'e}vy processes, i.e.\ processes with independent and stationary increments, can be defined whenever we have a notion of independence with an associated convolution. Under some weak conditions (cf.\ \cite{schuermann95,speicher96,benghorbal+schuermann99,muraki02a,muraki02}) this leaves only five possibilities, Voicu\-lescu's \cite{voiculescu+dykema+nica92} free independence, boolean independence (see e.g.\ \cite{speicher+woroudi93}), tensor independence, which is the closest to the notion of independence used in classical probability theory, and the monotone and anti-monotone independence \cite{lu97b,muraki97a}. Works by two of the present authors \cite{schuermann95b,franz03b} and by R.~Lenczewski \cite{lenczewski98b} indicate that many questions about the free, boolean, monotone, or anti-monotone convolution and L{\'e}vy processes can be resolved by reducing them to the tensor case.

In this paper we remove an asymmetry in the present theory of L{\'e}vy processes related to tensor independence, i.e.\ L{\'e}vy processes on bialgebras \cite{accardi+schuermann+waldenfels88,schuermann93,franz+schott99}. The algebras on which the processes are defined were $*$-bialgebras, but the `twisting' of the tensor product was supposed to stem from a group action and coaction. But, in general, a `twisted' (or braided) tensor category does not come from a group but from a quasitriangular bialgebra, see e.g.\ \cite[Chapter 9]{majid95} and the references therein. Here we will define and study L{\'e}vy processes on $*$-bialgebras in the category of representations of another $*$-bialgebra. In particular, for explicit calculations we will consider processes on primitively generated $*$-bialgebras with a braiding defined by an R-matrix. Here the bialgebra `behind' the braiding is that associated to the R-matrix by the Faddeev-Reshetikhin-Takhtajan (FRT) construction. Following Majid we dub these algebras `braided ($*$-) spaces'.

There are several reasons for the generalization presented in the paper. First, as explained above, it is an unnecessary restriction to consider only braidings that can be defined via an action and coaction of a group. Taking a bialgebra (satisfying some additional conditions) instead of the group seems to be the natural framework. But this requires some changes in the theory, e.g.\ the compatibility condition for the action and the coaction used in \cite{schuermann93} only works for cocommutative bialgebras, in the general case it has to be replaced by the one stated in Lemma \ref{lem-cond-comp}.

Second, it leads to an interesting new class of processes. In dimension one there is only one R-matrix, the $1\times 1$-matrix $R=(q)$ with $q\not=0$. In this case the braiding can be defined via a group, and we obtain e.g.\ the Az{\'e}ma process, whose amazing properties have been studied in \cite{azema85,emery89,parthasarathy90,schuermann91b}. But in higher dimension there are many possibilities to choose the R-matrix, and most of them cannot be obtained from groups. We will classify the Brownian motions (i.e.\ L\'evy processes with quadratic generator) on the braided spaces associated to the $sl_2$- and $sl_3$-R-matrix in Section \ref{example}.

Another motivation is the possibility of so-called braid statistics in two- and three-dimensional quantum field theories, see \cite{majid93c,oeckl01}.

This paper is organized as follows.

In Section \ref{theory} we recall the basic definitions and fundamental results and extend them to braided $*$-bialgebras. In particular, this requires a new definition of $*$-structures, as the one proposed previously by S.~Majid \cite{majid94d,majid94c,majid95b} would lead to a definition of positivity of functionals that is not preserved by the convolution, see also \cite{franz+schott96}.

In Section \ref{categories}, we show how the braided categories can be constructed from a Hopf algebra or a coquasi-triangular bialgebra. Our particular setting has been chosen since it will be convenient when we consider L\'evy processes on braided $*$-algebras in these categories. We consider two cases. In the first we have a Hopf $*$-algebra $\mathcal{A}$ and the objects of the category are involutive vector spaces carrying an action and a coaction of $\mathcal{A}$ that satisfy compatibility conditions with respect to each other and with respect to the involution, see Equations \eqref{cond-comp}, \eqref{cond-alpha-*}, and \eqref{cond-gamma-*}. In the second case $\mathcal{A}$ is an involutive coquasi-triangular bialgebra, i.e.\ it is a $*$-bialgebra and it is equipped with a universal R-form. We show that if the R-form satisfies the  compatibility condition \eqref{cond-comp-r-*} with respect to the involution, then we can build a braided category from the comodules whose coaction satisfies \eqref{cond-gamma-*}.

In Section \ref{symmetrization} we show that for every braided $*$-bialgebra $\mathcal{B}$ in either of the kinds of categories we can construct a symmetrization, i.e., we can embed it into a bigger $*$-bialgebra $H$ as an algebra in a way that allows to `lift' L\'evy processes on $\mathcal{B}$ to L\'evy processes on $H$. This reduces the construction and classification of L\'evy processes on braided $*$-bialgebras to usual involutive bialgebras.

In the following section (Section \ref{construction}), we explicitely construct a family of braided $*$-spaces and their bosonization or symmetrization. The construction starts from a bi-invertible $R$-matrix of real type I, it yields a braided $*$-space in the category of comodules of the quasi-triangular $*$-bialgebra associated with this R-matrix by the Faddeev-Reshetikhin-Takhtajan construction \cite{reshetikhin+takhtadzhyan+faddeev90}. Furthermore, these braided $*$-spaces always come with a canonical quadratic generator that gives rise to a standard Brownian motion on these spaces.

In Section \ref{levy-proc}, we prove the existence of these invariant, conditionally positive, quadratic functionals on our braided $*$-spaces and study the associated Brownian motions. These processes can be considered as multi-dimensional analogues of the Az\'ema martingales \cite{azema85,emery89,parthasarathy90,attal+emery94,attal+emery96}.

Section \ref{example} contains the explicit quantum stochastic differential equations defining the standard Brownian motion on the braided line and the braided plane as well as the classification of all quadratic generators for the braided $*$-spaces associated to the $1\times1$-R-matrix $(q)$, and the standard $sl_2$- and $sl_3$-R-matrices. In classical probability so-called structure equations expressing the It\^o products ${\rm d}Z_i\cdot{\rm d}Z_j$ of the differentials of multi-dimensional Az\'ema martingales as linear combinations of the differentials ${\rm d}Z_i$ and ${\rm d}t$ play a fundamental r\^ole. We show that in general their quantum analogues do not satisfy a similar set of nice closed structure equations.

Finally, in the last section (Section \ref{conclusion}), we indicate two directions of possible further generalization.

\section{L{\'e}vy processes on braided $*$-Hopf algebras}\label{theory}

In this section we introduce the notions of braided $*$-Hopf algebras or braided $*$-bialgebras and L{\'e}vy processes on a braided $*$-bialgebra and develop some of their elementary theory, in particular the one-to-one correspondence between these processes and their generators, see Theorem \ref{th-gen-proc}. Most of this is a straightforward generalization of \cite{schuermann93}, therefore we will focus on what is new. The important results are the new axioms for $*$-structures on braided Hopf algebras (Definition \ref{*-hopf}), the notion of $\Psi$-invariance (Definition \ref{def-invariance-Psi} and Lemma \ref{lem-invariant}) and Lemma \ref{pos-conv} which shows that convolutions of ($\Psi$- or $\Psi^{-1}$-invariant) positive functionals are again positive. This allows to extend \cite[Theorem 1.9.2]{schuermann93} and \cite[Theorem 3.2.8]{schuermann93}, i.e.\ the one-to-one correspondence between invariant hermitian, conditionally positive linear functionals and L{\'e}vy processes, see Theorem \ref{th-gen-proc}.

A tensor category or monoidal category is a category $\mathcal{C}$ equipped with a bifunctor $\otimes:\mathcal{C}\times \mathcal{C} \to \mathcal{C}$ satisfying certain conditions, cf. \cite{maclane98,kassel95}. A braiding $\Psi$ in a tensor category is a natural isomorphism between the functors $\otimes:(A,B)\mapsto A\otimes B$ and $\otimes\circ\tau:(A,B)\mapsto B\otimes A$ satisfying the so-called Hexagon Axioms. It is called a symmetry, if it is involutive, i.e.\ if $\Psi_{B,A}\circ\Psi_{A,B}={\rm id}_{A\otimes B}$ for all objects $A,B$ of $\mathcal{C}$. A braided tensor category or braided monoidal category is a pair $(\mathcal{C},\Psi)$ consisting of a tensor category $\mathcal{C}$ and braiding $\Psi$ of $\mathcal{C}$, it is called a symmetric tensor category, if the braiding is a symmetry.

The notions of bialgebras and Hopf algebras can also be defined in braided tensor category, this leads to braided bialgebras and braided Hopf algebras. The product, coproduct, unit, counit and antipode now have to be morphisms of the braided tensor category and satisfy similar axioms as in the usual case, cf.\ \cite{majid95}. Bialgebras and Hopf algebras are special cases of braided bialgebras and braided Hopf algebras, it suffices to take the flip automorphism $\tau:A\otimes B\to B\otimes A$, $\tau(u\otimes v)=v\otimes u$ for the braiding.

We begin by giving the definition of a braided involutive bialgebra.

\begin{definition}\label{*-hopf}
A braided $*$-bialgebra is a braided bialgebra $(\mathcal{A},\Delta,\varepsilon,m,1,\Psi)$ over $\mathbb{C}$ with an anti-linear map $*:\mathcal{A}\to \mathcal{A}$ that satisfies
\begin{itemize}
\item
$(\mathcal{A},m,1,*)$ is a $*$-algebra,
\item
$\mathcal{A}\otimes \mathcal{A}$ admits a $*$-structure such that the canonical inclusions $\mathcal{A}\stackrel{i_1}{\longrightarrow} \mathcal{A}{\otimes}\mathcal{A} \stackrel{i_2}{\longleftarrow}\mathcal{A}$ and the coproduct $\Delta:\mathcal{A}\to \mathcal{A}\otimes \mathcal{A}$ are $*$-algebra homomorphisms for this $*$-structure.
\end{itemize}
If there exists an antipode $S$ (i.e. a linear map $S:\mathcal{A}\to \mathcal{A}$ s.t.\ $(\mathcal{A},\Delta,\varepsilon,m,1,S,\Psi)$ is a braided Hopf algebra), then we will call $(\mathcal{A},\Delta,\varepsilon,m,1,S,\Psi,*)$ a braided $*$-Hopf algebra.
\end{definition}
\begin{remark} 
{}From the condition that the canonical inclusions $\mathcal{A}\stackrel{i_1}{\longrightarrow} \mathcal{A}{\otimes}\mathcal{A} \stackrel{i_2}{\longleftarrow}\mathcal{A}$ are $*$-homomorphisms follows $(a\otimes 1)^* = a^*\otimes 1$, $(1\otimes a)^* = 1\otimes a^*$ for all $a\in \mathcal{A}$. Now the condition that this map is an anti-homomorphism uniquely determines it, $(a\otimes b)^* = \Big((a\otimes1)(1\otimes b)\Big)^* = (1\otimes b^*)(a^*\otimes 1)= \Psi(b^*\otimes a^*)$ for all $a,b\in \mathcal{A}$. This defines a *-structure on $\mathcal{A}\otimes \mathcal{A}$ if and only if it is also an involution, i.e.\ if $(\Psi\circ (*\otimes *)\circ \tau)^2={\rm id}$. To see this, let $\{x_i\}$ be a basis of $\mathcal{A}$ consisting of self-adjoint elements, i.e.\ such that $x^*_i=x_i$. Define the coefficients of the braiding $\Psi$ in this basis by $\Psi(x_i\otimes x_j)= \sum_{kl}\Psi_{ij}^{kl} x_l\otimes x_k$ (note that this sum is finite, even though the index set can be infinite). $\Psi\circ(*\otimes *)\circ \tau$ is an involution, if and only if $\sum_{nm} \overline{\Psi_{ji}^{nm}} \Psi_{nm}^{kl} = \delta^l_i \delta^k_j$. With this relation we can show that $\Psi\circ(*\otimes *)\circ \tau$ is an anti-homomorphism. It is sufficient to check this on the basis elements, we get
\begin{eqnarray*}
&&\Big( (x_i\otimes x_j)\cdot (x_k\otimes x_l)\Big)^* \\
&=&
\sum_{nm}(\Psi_{jk}^{nm} x_ix_m \otimes x_n x_l)^* = \sum_{nm}\overline{\Psi_{jk}^{nm}}\Psi\Big((x_nx_l)^* \otimes (x_i x_m)^*\Big) \\
&=& \sum_{nm}\overline{\Psi_{jk}^{nm}} (m\otimes m) \circ ({\rm id}\otimes \Psi
\otimes {\rm id})\circ (\Psi\otimes \Psi)\circ ({\rm id}\otimes \Psi
\otimes {\rm id})(x_l \otimes x_n \otimes x_m \otimes x_i) \\
&=& \sum_{nmpq} \overline{\Psi_{jk}^{nm}} \Psi_{nm}^{pq} (m\otimes m) \circ ({\rm id}\otimes \Psi \otimes {\rm id})\circ (\Psi\otimes \Psi) (x_l\otimes x_q \otimes x_p \otimes x_i) \\
&=& (m\otimes m) \circ ({\rm id}\otimes \Psi \otimes {\rm id})\circ
(\Psi\otimes \Psi) (x_l \otimes x_k \otimes x_j \otimes x_i) \\
&=& (m\otimes m) \circ ({\rm id}\otimes \Psi \otimes {\rm id}) \Big(\Psi(x_l
\otimes x_k) \otimes \Psi(x_j \otimes x_i)\Big) \\
&=& (x_k \otimes x_l)^* \cdot (x_i\otimes x_j)^*.
\end{eqnarray*}
If $\Psi\circ(*\otimes *)\circ \tau$ is an involution, then we also have $\Psi\circ * \circ \Psi \circ * = {\rm id}_{\mathcal{A}\otimes \mathcal{A}}$.
\end{remark}
\begin{remark}
Note that the axioms for the $*$-structure differ from those proposed by S.~Majid \cite{majid94d,majid94c,majid95b}.
\end{remark}
\begin{remark}
If the coproduct is fixed, then the counit and the antipode (if it exists) are unique. Therefore one expects their axioms to be a consequence of those given in the definition. One can show that $\varepsilon$ is a $*$-algebra homomorphism and that $S$ satisfies $S\circ *\circ S\circ*={\rm id}$, if it exists.
\end{remark}

\begin{definition}
A braided $*$-space is a braided space (cf.\ \cite{majid96}) equipped with an involution that turns it into a braided $*$-Hopf algebra in the sense of the preceding definition.
\end{definition}

Let us now come to the definition of L\'evy processes. A quantum probability space is a pair $(\mathcal{A},\Phi)$ consisting of a $*$-algebra $\mathcal{A}$ and a state (i.e.\ a normalized positive linear functional) $\Phi$ on $\mathcal{A}$. A quantum random variable $j$ over a quantum probability space $(\mathcal{A},\Phi)$ on a $*$-algebra $\mathcal{B}$ is a $*$-algebra homomorphism $j:\mathcal{B}\to\mathcal{A}$. A quantum stochastic process is simply a family of quantum random variables over the same quantum probability space, indexed by some set, and defined on the same algebra.

\begin{definition}
Let $(\mathcal{A},\Phi)$ be a quantum probability space, $\mathcal{B}$ a $*$-algebra, and $\Psi:\mathcal{B}\otimes \mathcal{B} \to \mathcal{B}\otimes\mathcal{B}$ a linear map. An $n$-tuple $(j_1,\ldots,j_n)$ of quantum random variables $j_i:\mathcal{B}\to \mathcal{A}$, $i=1,\ldots,n$, over $(\mathcal{A},\Phi)$ on $\mathcal{B}$ is $\Psi$-independent or braided independent, if
\begin{itemize}
\item[(i)]
$\Phi\Big(j_{\sigma(1)}(b_1) \cdots j_{\sigma(n)}(b_n)\Big) = \Phi(j_{\sigma(1)}(b_1)) \cdots \Phi(j_{\sigma(n)}(b_n))$ for all permutations $\sigma\in\mathcal{S}(n)$ and all $b_1,\ldots,b_n\in\mathcal{B}$, and
\item[(ii)]
$m_{\mathcal{A}} \circ (j_l \otimes j_k) = m_{\mathcal{A}} \circ (j_k\otimes j_l) \circ \Psi$ for all $1\le k < l \le n$.
\end{itemize}
\end{definition}

\begin{definition}
Let $\mathcal{B}$ be a braided $*$-bialgebra. A quantum stochastic process $(j_{st})_{0\le s\le t\le T}$, $T\in\mathbb{R}_+\cup\{\infty\}$ on $\mathcal{B}$ over some quantum probability space $(\mathcal{A},\Phi)$ is called a L{\'e}vy process, if the following conditions are satisfied.
\begin{enumerate}
\item
(Increment property)
\begin{eqnarray*}
j_{rs}\star j_{st} &=& j_{rt} \quad \mbox{ for all } 0\le r\le s\le t\le T, \\
j_{tt} &=& e\circ \varepsilon \quad \mbox{ for all } 0\le t\le T.
\end{eqnarray*}
\item
(Independence of increments)
The family $(j_{st})_{0\le s\le t\le T}$ is independent, i.e.\ $(j_{s_1t_1},$ $\ldots, j_{s_nt_n})$ is independent for all $n\in\mathbb{N}$ and all $0\le s_1\le t_1\le s_2\le \cdots t_n\le T$.
\item
(Stationarity of increments)
The distribution $\varphi_{st}=\Phi\circ j_{st}$ of $j_{st}$ depends only on the difference $t-s$,
\item
(Weak continuity) $j_{st}$ converges to $j_{ss}\,(\,=e\circ\varepsilon)$ in distribution for $t\searrow s$.
\end{enumerate}
\end{definition}

\begin{remark}
If $\mathcal{B}$ is a Hopf $*$-algebra, then we can define the increments of a process $(j_t)_{t\in[0,T]}$ as follows:
\[
j_{st}= m_\mathcal{A}\circ \Big((j_s\circ S) \otimes j_t\Big)\circ \Delta, \qquad 0\le s\le t\le T.
\]
With this definition the increment property is automatically satisfied. We call a process $(j_t)_{t\in[0,T]}$ on a Hopf $*$-algebra $\mathcal{B}$ a L{\'e}vy process, if its increment process $(j_{st})_{0\le s\le t\le T}$ is a L{\'e}vy process on $\mathcal{B}$.
\end{remark}

\begin{definition}\label{def-invariance-Psi}
Let $(\mathcal{C},\Psi)$ be a braided tensor category whose objects are vector spaces. A linear map $L:V\to W$ is called $\Psi$-invariant, if $(L \otimes {\rm id}_X)\circ \Psi_{X,V} = \Psi_{X,W}\circ ({\rm id}_X\otimes L)$ for all $X$.
\end{definition}

\begin{remark}
Note that due to the naturality of $\Psi$, morphisms of the category $(\mathcal{C},\Psi)$ always have to be $\Psi$- and $\Psi^{-1}$-invariant.
\end{remark}

\begin{lemma}\label{pos-conv}
Let $\mathcal{A}$ be a braided $*$-bialgebra and let $\phi,\theta$ be two positive functionals on $\mathcal{A}$. If $\phi$ is $\Psi$-invariant or $\theta$ is $\Psi^{-1}$-invariant, then $\phi\star\theta = (\phi \otimes \theta)\circ {\Delta}$ is also positive.
\end{lemma}
\begin{proof}
Let $a\in\mathcal{A}$, and ${\Delta}(a) = \sum_i a_i^{(1)}\otimes a_i^{(2)}$ (Sweedler's notation). Then ${\Delta}(a^*)= \sum_j \Psi\left(\left(a_j^{(2)}\right)^*\otimes \left(a_j^{(1)}\right)^*\right)$. If $\phi$ is $\Psi$-invariant, then
\begin{eqnarray*}
&& (\phi \otimes \theta)\circ {\Delta}(a^*a)\\
&& = (\phi \otimes \theta) \circ (m\otimes m)\circ ({\rm id}\otimes \Psi \otimes {\rm id})\circ (\Psi\otimes{\rm id}\otimes {\rm id})\left(\sum_{i,j}\left(a_j^{(2)}\right)^*\otimes \left(a_j^{(1)}\right)^*\otimes a_i^{(1)}\otimes a_i^{(2)} \right) \\
&& = \theta\circ (\phi\otimes m)\circ (\Psi\circ {\rm id})\circ({\rm id}\otimes m \otimes {\rm id})\left(\sum_{i,j}\left(a_j^{(2)}\right)^*\otimes \left(a_j^{(1)}\right)^*\otimes a_i^{(1)}\otimes a_i^{(2)} \right) \\
&& = \sum_{i,j}\theta\left(\left(a_j^{(2)}\right)^*a_i^{(2)}\right) \phi\left(\left(a_j^{(1)}\right)^* a_i^{(1)}\right)
\end{eqnarray*}
is positive, since it is the Schur product of two positive definite matrices.

If instead $\theta$ is $\Psi^{-1}$-invariant, then we can show in the same way that
\[
(\phi \otimes \theta)\circ {\Delta}(aa^*) = \sum_{i,j}\phi\left( a_j^{(1)} \left(a_i^{(1)}\right)^*\right) \theta\left(a_j^{(2)} \left(a_i^{(2)}\right)^*\right)
\]
is positive.
\end{proof}

Let $(j_{st})$ be a L{\'e}vy process on some *-bialgebra $\mathcal{B}$. The states on $\mathcal{B}$ defined by $\varphi_t=\Phi\circ j_{0t}$ are called the marginal distributions of $(j_{st})$. They form a convolution semi-group and can be written in the form $\varphi_t = exp_{\star} t L$, where $L:\mathcal{B}\to \mathbb{C}$ is a hermitian, conditionally positive (i.e.\ positive on the kernel of $\varepsilon$) linear functional.

\begin{definition}
Two quantum stochastic processes $(j_t)_{t\in I}$ and $(k_t)_{t\in I}$ over quantum probability spaces $(\mathcal{A},\Phi)$ and $(\mathcal{A}',\Phi')$, resp., on the same *-algebra $\mathcal{B}$ are called equivalent, if their finite joint moments agree, i.e.\ if
\[
\Phi(j_{t_1}(b_1)\cdots j_{t_n}(b_n)) = \Phi'(k_{t_1}(b_1)\cdots k_{t_n}(b_n))
\]
for all $n\in\mathbb{N}$, $t_1,\ldots,t_n\in I$, and $b_1,\ldots, b_n\in \mathcal{B}$.
\end{definition}

\begin{theorem}\label{th-gen-proc}
Let $\mathcal{B}$ be a braided $*$-bialgebra. The set of L{\'e}vy processes $(j_{st})_{0\le s\le t \le T}$ on $\mathcal{B}$ (modulo equivalence), the set of convolution semi-groups $(\varphi_t)_{t\in[0,T]}$ of $\Psi$-invariant states on $\mathcal{B}$, and the set of $\Psi$-invariant, hermitian, conditionally positive linear functionals $L:\mathcal{B}\to \mathbb{C}$ are in one-to-one correspondence.
\end{theorem}
\begin{proof}
This theorem is the straightforward generalization of \cite[Theorem 1.9.2]{schuermann93} and \cite[Theorem 3.2.8]{schuermann93}, the only non-trivial part is the positivity of the convolution semi-group generated by a $\Psi$-invariant, hermitian, conditionally positive linear functional, and this follows immediately from Lemma \ref{pos-conv}.
\end{proof}

\section{A construction of braided categories}\label{categories}

We will now study the special case of tensor categories whose objects are left modules and left comodules of some given bialgebra $\mathcal{A}$. For the rest of this section we will assume that for all objects $V,W$ of our category $\mathcal{C}$ we have actions $\alpha_V:\mathcal{A}\otimes V\to V$, $\alpha_W:\mathcal{A}\otimes W\to W$ and coactions $\gamma_V:V\to\mathcal{A}\otimes V$ and $\gamma_W:W\to\mathcal{A}\otimes W$ and that all morphisms of the category are module and comodule maps with respect to these actions and coactions. Our goal is to determine sufficient conditions on $\alpha$ and $\gamma$ such that $\Psi_{V,W}=(\alpha_W\otimes {\rm id})\circ ({\rm id}\otimes \tau)\circ (\gamma_V\otimes{\rm id})$ defines a braiding for $\mathcal{C}$. We will sometimes omit the subscripts, if it is clear to what objects $\alpha,\gamma,\Psi$ belong.

Analogously, one could study categories consisting of right modules and comodules and define a braiding as $\Psi=({\rm id}\otimes\alpha)\circ(\tau\otimes {\rm id})\circ({\rm id}\otimes \gamma)$.

\begin{lemma}\label{lem-inv}
Assume that $\mathcal{A}$ has an invertible antipode $S$. Then the maps $\Psi_{V,W}=(\alpha_W\otimes {\rm id})\circ ({\rm id}\otimes \tau)\circ (\gamma_V\otimes{\rm id})$ are invertible, their inverses are given by $\Psi^{-1}_{V,W}=\tau\circ(\alpha_W\otimes {\rm id})\circ(S^{-1}\otimes {\rm id}\otimes{\rm id})\circ(\tau\otimes {\rm id})\circ ({\rm id}\otimes\gamma_V)$.
\end{lemma}
\begin{proof} We get
\begin{eqnarray*}
&&\Psi^{-1}_{V,W}\circ\Psi_{V,W} \\
&=&\tau\circ(\alpha\otimes {\rm id})\circ (S^{-1}\otimes {\rm id}\otimes{\rm id})\circ(\tau\otimes {\rm id})\circ ({\rm id}\otimes\gamma)\circ(\alpha\otimes {\rm id})\circ ({\rm id}\otimes \tau)\circ (\gamma\otimes{\rm id})\\
&=&\tau\circ(\alpha\otimes{\rm id})\circ({\rm id}\otimes\alpha\otimes{\rm id})\circ (S^{-1}\otimes{\rm id}\otimes{\rm id}\otimes{\rm id})\circ (\tau\otimes \tau) \circ({\rm id}\otimes\gamma\otimes{\rm id})\circ(\gamma\otimes{\rm id})\\
&=&\tau\circ(\alpha\otimes{\rm id})\circ(m \otimes{\rm id}\otimes{\rm id})\circ (S^{-1}\otimes{\rm id}\otimes{\rm id}\otimes{\rm id})\circ (\tau\otimes \tau) \circ(\Delta\otimes{\rm id}\otimes{\rm id})\circ(\gamma\otimes{\rm id})\\
&=&\tau\circ(\alpha\otimes{\rm id}) \circ ((1\circ\varepsilon)\otimes\tau) \circ(\gamma\otimes{\rm id})\\
&=&\tau\circ\tau={\rm id}_{V\otimes W}
\end{eqnarray*}
where we used first the action and coaction axioms $\alpha\circ({\rm id}\otimes\alpha)= \alpha\circ(m \otimes{\rm id})$ and $({\rm id}\otimes\gamma)\circ\gamma=(\Delta\otimes{\rm id})\circ \gamma$, and then the fact that $S^{-1}$ is an antipode for $\mathcal{A}^{\rm coop}= (\mathcal{A},\Delta^{\rm coop}=\tau\circ\Delta,m)$, i.e.\
\[
m\circ(S^{-1}\otimes{\rm id})\circ\tau\circ\Delta=1\circ \varepsilon.
\]
Similarly, we can show $\Psi_{V,W}\circ\Psi^{-1}_{V,W}={\rm id}_{W\otimes V}$.
\end{proof}

\begin{lemma}\label{lem-cond-comp}
Assume that the actions and coactions satisfy the relation
\begin{equation}\label{cond-comp}
(m\otimes\alpha)\circ({\rm id}\otimes\tau\otimes{\rm id})\circ(\Delta\otimes\gamma) = (m\otimes {\rm id})\circ ({\rm id}\otimes \tau)\circ (\gamma\otimes{\rm id}) \circ(\alpha\otimes {\rm id})\circ ({\rm id}\otimes \tau)\circ (\Delta\otimes{\rm id}).
\end{equation}
Then $\Psi=(\alpha\otimes {\rm id})\circ ({\rm id}\otimes \tau)\circ (\gamma\otimes{\rm id})$ is a module map and a comodule map, i.e. it satisfies
\begin{eqnarray*}
\Psi_{V,W}\circ\alpha_{V\otimes W} &=& \alpha_{W\otimes V}\circ({\rm id}_\mathcal{A}\otimes\Psi_{V,W}), \\
\gamma_{W\otimes V}\circ \Psi_{V,W} &=& ({\rm id}_\mathcal{A}\otimes\Psi_{V,W})\circ \gamma_{V\otimes W}.
\end{eqnarray*}
\end{lemma}
\begin{proof}
The action $\alpha_{V\otimes W}$ is given by
\[
\alpha_{V\otimes W}= (\alpha_V\otimes \alpha_W)\circ({\rm id}\otimes\tau\otimes{\rm id})\circ(\Delta\otimes{\rm id}\otimes{\rm id}).
\]
We get 
\begin{eqnarray*}
&&\Psi_{V,W}\circ\alpha_{V\otimes W} \\
&=&
(\alpha\otimes {\rm id})\circ ({\rm id}\otimes \tau)\circ (\gamma\otimes{\rm id})\circ(\alpha\otimes \alpha)\circ({\rm id}\otimes\tau\otimes{\rm id})\circ(\Delta\otimes{\rm id}_{V\otimes W}) \\
&=&(\alpha\otimes {\rm id})\circ({\rm id}\otimes\alpha\otimes{\rm id})\circ (\gamma\otimes{\rm id}_{\mathcal{A}\otimes W})\circ(\alpha\otimes{\rm id}_{\mathcal{A}\otimes W})\circ({\rm id}\otimes\tau\otimes{\rm id}) \circ(\Delta\otimes{\rm id}_{V\otimes W}) \\
&=&(\alpha\otimes {\rm id})\circ(m \otimes{\rm id}_{W\otimes V})\circ({\rm id}\otimes\tau_{V,\mathcal{A}\otimes W})\circ (\gamma\otimes{\rm id}_{\mathcal{A}\otimes W})\circ(\alpha\otimes{\rm id}_{\mathcal{A}\otimes W})\circ \\
&& \circ({\rm id}\otimes\tau\otimes{\rm id})\circ(\Delta\otimes{\rm id}_{V\otimes W}) \\
&=& (\alpha\otimes {\rm id})\circ({\rm id}\otimes\tau)\circ((m\otimes {\rm id})\circ ({\rm id}\otimes \tau)\circ (\gamma\otimes{\rm id}) \circ(\alpha\otimes {\rm id})\circ ({\rm id}\otimes \tau)\circ (\Delta\otimes{\rm id}))\otimes{\rm id})\\
&=&(\alpha\otimes {\rm id})\circ({\rm id}\otimes\tau)\circ((m\otimes\alpha)\circ({\rm id}\otimes\tau\otimes{\rm id})\circ(\Delta\otimes\gamma))\otimes{\rm id})\\
&=&
(\alpha\otimes\alpha)\circ(m\otimes{\rm id}\otimes{\rm id}\otimes{\rm id})\circ ({\rm id}\otimes{\rm id}\otimes\tau\otimes{\rm id})\circ({\rm id}\otimes\tau\otimes\tau)\circ(\Delta\otimes\gamma\otimes{\rm id})\\
&=&
(\alpha\otimes\alpha)\circ({\rm id}\otimes\alpha\otimes{\rm id}\otimes{\rm id})\circ ({\rm id}\otimes{\rm id}\otimes\tau\otimes{\rm id})\circ({\rm id}\otimes\tau\otimes\tau)\circ(\Delta\otimes\gamma\otimes{\rm id}) \\
&=&
(\alpha\otimes\alpha)\circ({\rm id}\otimes\tau\otimes{\rm id})\circ(\Delta\otimes{\rm id}\otimes{\rm id})\circ({\rm id}\otimes((\alpha\otimes
{\rm id})\circ ({\rm id}\otimes \tau)\circ (\gamma\otimes{\rm id})))\\
&=&\alpha_{W\otimes V}\circ({\rm id}\otimes\Psi_{V,W}),
\end{eqnarray*}
where, besides some reordering, we used first the action axiom $\alpha\circ({\rm id}\otimes\alpha)= \alpha\circ(m \otimes{\rm id})$, then relation (\ref{cond-comp}), and finally the action axiom again in the opposite direction.

Similarly we get the corresponding identity for the coaction
\[
\gamma_{V\otimes W}=(m\otimes{\rm id}\otimes{\rm id})\circ({\rm id}\otimes\tau\otimes{\rm id})\circ(\gamma_V\otimes\gamma_W).
\]
\end{proof}

\begin{lemma}\label{lem-*-inv}
Assume now that $\mathcal{A}$ is a $*$-Hopf algebra and suppose that there also exists a $*$-structure on $V$. If $\alpha$ and $\gamma$ satisfy
\begin{eqnarray}
*\circ \alpha &=& \alpha\circ(*\otimes *)\circ(S\otimes{\rm id}), \label{cond-alpha-*}\\
\gamma \circ * &=& (*\otimes *)\circ \gamma, \label{cond-gamma-*}
\end{eqnarray}
then $\Psi\circ (*\otimes *) \circ \tau: V\otimes V \to V\otimes V$ is an involution.
\end{lemma}
\begin{remark}
The analogous result for $\Psi=({\rm id}\otimes \alpha)\circ(\tau\otimes {\rm id})\circ({\rm id}\otimes \gamma)$ with a right action $\alpha$ and a right coaction $\gamma$ also holds.
\end{remark}
\begin{proof}
\begin{eqnarray*}
&&(\Psi\circ (*\otimes *)\circ \tau)^2 \\
&=& \Psi\circ (*\otimes *)\circ \tau \circ (\alpha\otimes {\rm id})\circ ({\rm id}\otimes \tau)\circ (\gamma\otimes{\rm id}) (*\otimes *)\circ \tau \\
&=& \Psi \circ \tau \circ (\alpha\otimes {\rm id})\circ(*\otimes * \otimes *)\circ(S\otimes{\rm id}\otimes{\rm id})\circ(*\otimes * \otimes *)\circ({\rm id}\otimes \tau) \circ(\gamma\otimes {\rm id})\circ\tau\\
&=& \Psi \circ\tau \circ(\alpha \otimes {\rm id}) \circ(S^{-1}\otimes \tau) \circ (\gamma\circ{\rm id}) \circ \tau \\
&=& (\alpha\otimes {\rm id})\circ ({\rm id}\otimes \tau)\circ (\gamma\otimes{\rm id}) \circ\tau \circ(\alpha \otimes {\rm id}) (S^{-1}\otimes \tau) \circ (\gamma\circ{\rm id}) \circ \tau \\
&=&(\alpha\otimes{\rm id})\circ({\rm id}\otimes\alpha\otimes{\rm id}) \circ(S^{-1}\otimes{\rm id}\otimes{\rm id}\otimes{\rm id})\circ(\tau\otimes\tau) \circ({\rm id}\otimes\gamma\otimes{\rm id})\circ(\gamma\otimes{\rm id})\circ\tau\\
&=&(\alpha\otimes{\rm id})\circ(m \otimes{\rm id}\otimes{\rm id})\circ (S^{-1}\otimes{\rm id}\otimes{\rm id}\otimes{\rm id})\circ (\tau\otimes \tau) \circ(\Delta\otimes{\rm id}\otimes{\rm id})\circ(\gamma\otimes{\rm id})\circ\tau\\
&=&(\alpha\otimes{\rm id}) \circ ((1\circ\varepsilon)\otimes\tau) \circ(\gamma\otimes{\rm id})\circ\tau\\
&=&\tau\circ\tau={\rm id}_{V\otimes W}
\end{eqnarray*}
where we used first the relations $*\circ \alpha = \alpha\circ (*\otimes *)\circ(S\otimes{\rm id})$ and $\gamma \circ * = (*\otimes *)\circ \gamma$ and $S^{-1}=*\circ S\circ *$, and then the same arguments as in the proof of Lemma \ref{lem-inv}.
\end{proof}

\begin{definition}\label{def-invariance-alpha-gamma}
\begin{enumerate}
\item
Let $\alpha_V:\mathcal{A}\otimes V\to V$ and $\alpha_W:\mathcal{A}\otimes W \to W$ be two actions of $\mathcal{A}$. We will call a linear map $L:V\to W$ $\alpha$-invariant, if it is a module map, i.e.\ if $L\circ \alpha_V = \alpha_W \circ ({\rm id}\otimes L)$.
\item
Let $\gamma_V:V\to \mathcal{A}\otimes V$ and $\gamma_W:W\to \mathcal{A}\otimes W$ be two coactions of $\mathcal{A}$. We will call a linear map $L:V\to W$ $\gamma$-coinvariant, if it is a comodule map, i.e.\ if $\gamma_W\circ L = ({\rm id}\otimes L)\circ \gamma_V$.
\end{enumerate}
\end{definition}
   The corresponding definitions for right actions or coactions are
   obvious. In the cases that interest us most, $L$ will be a functional and thus
   $W$ the base field $\mathbb{C}$ (equipped with the trivial action
   $\alpha_\mathbb{C}=\varepsilon\otimes {\rm id}$ and coaction
   $\gamma_\mathbb{C}=e\otimes {\rm id}$, and the trivial braiding $\mathbb{C}\otimes
   X\cong X\cong X\otimes\mathbb{C}$).

   \begin{lemma}\label{lem-invariant}
   Suppose that $\Psi$ is of the form $\Psi=(\alpha\otimes {\rm id})\circ ({\rm
   id}\otimes \tau)\circ (\gamma\otimes{\rm id})$, for some action
   $\alpha$ and coaction $\gamma$, and let $L:V\to W$ be a linear map.
   \begin{enumerate}
   \item
   If $L$ is $\alpha$-invariant, then it is $\Psi$-invariant.
   \item
   If $L$ is $\gamma$-coinvariant, then it is $\Psi^{-1}$-invariant.
   \end{enumerate}
   \end{lemma}
\begin{remark}
If $\Psi$ is of the form $\Psi=({\rm id}\otimes \alpha)\circ(\tau\otimes {\rm id})\circ({\rm id}\otimes \gamma)$ with a right action $\alpha$ and a right coaction $\gamma$, then the preceding proposition holds, if we interchange $\Psi$ and $\Psi^{-1}$.
\end{remark}
\begin{proof}
\begin{enumerate}
\item
Let $L$ be $\alpha$-invariant, then
\begin{eqnarray*}
   (L\otimes {\rm id})\circ\Psi_{X,V}&=&(L\otimes {\rm id})\circ(\alpha_V\otimes {\rm id})\circ ({\rm
   id}\otimes \tau)\circ (\gamma_X\otimes{\rm id}) \\
   &=& (\alpha_W\otimes{\rm id})\circ({\rm id}\otimes L\otimes{\rm id})\circ({\rm id}\otimes \tau)\circ(\gamma_X \otimes {\rm id}) \\
   &=& (\alpha_W\otimes {\rm id})\circ ({\rm
   id}\otimes \tau)\circ (\gamma_X\otimes{\rm id})\circ(L\otimes {\rm id}) \\
   &=&\Psi_{X,W}\circ({\rm id}\otimes L),
   \end{eqnarray*}
   for all $X$ in $\mathcal{C}$.
   \item
   We have $L\otimes {\rm id}=(L\otimes {\rm id})\circ \Psi^{-1}\circ
   \Psi = ({\rm id}\otimes L)\circ \Psi)$, i.e.\ $L$ is $\Psi^{-1}$
   invariant, if and only if $({\rm id}\otimes L)\circ\Psi = L\otimes
   {\rm id}$. Let now $L$ be $\gamma$-coinvariant, then
   \begin{eqnarray*}
   ({\rm id}\otimes L)\circ \Psi_{V,X} &=& ({\rm id}\otimes L)\circ(\alpha_X\otimes {\rm id})\circ ({\rm id}\otimes \tau)\circ (\gamma_V\otimes{\rm id}) \\
   &=&(\alpha_X\otimes {\rm id})\circ ({\rm id}\otimes \tau)\circ (\gamma_W\otimes{\rm id})\circ (L\otimes {\rm id})\\
   &=& \Psi_{W,X}\circ (L\otimes {\rm id}),
   \end{eqnarray*}
   for all $X$ in $\mathcal{C}$.
   \end{enumerate}
   \end{proof}

The preceding lemmas actually show that $\Psi$ defines a braiding. The construction is analog to the construction of the crossed bimodules which uses a left action and a right coaction, see \cite[Section IX.5]{kassel95} and the references cited there. 

\begin{proposition}\label{prop-braided}
Let $\mathcal{A}$ be a Hopf algebra with invertible antipode. Then we can define a braided category $(\mathcal{C}(\mathcal{A}),\Psi)$ as follows. The objects of $\mathcal{C}(\mathcal{A})$ are triples $(V,\alpha_V,\gamma_V)$ consisting of a vector space $V$ equipped with an action $\alpha_V$ and a coaction $\gamma_V$ such that Equation (\ref{cond-comp}) is satisfied.
The morphism between two objects $(V,\alpha_V,\gamma_V)$ and $(W,\alpha_W,\gamma_W)$ are the linear maps from $V$ to $W$ that are $\alpha$-invariant and $\gamma$-coinvariant, i.e. that are module and comodule maps.

The tensor product of two objects $(V,\alpha_V,\gamma_V)$ and $(W,\alpha_W,\gamma_W)$ is given by
\[
(V,\alpha_V,\gamma_V) \otimes (W,\alpha_W,\gamma_W) = (V\otimes W, \alpha_{V\otimes W},\gamma_{V\otimes W}),
\] 
where
\begin{eqnarray*}
\alpha_{V\otimes W}&=& (\alpha_V\otimes \alpha_W)\circ({\rm id}\otimes\tau\otimes{\rm id})\circ(\Delta\otimes{\rm id}\otimes{\rm id}), \\
\gamma_{V\otimes W}&=&(m\otimes{\rm id}\otimes{\rm id})\circ({\rm id}\otimes\tau\otimes{\rm id})\circ(\gamma_V\otimes\gamma_W).
\end{eqnarray*}
The tensor product of two morphisms $f_1:V_1\to W_1$ and $f_2:V_2\to W_2$ is the usual tensor product $f_1\otimes f_2$ of linear maps.

The braiding $\Psi: \otimes \stackrel{\cong}{\to} \otimes\circ\tau$ is defined by
\[
\Psi_{V,W}=(\alpha_W\otimes {\rm id})\circ ({\rm id}\otimes \tau)\circ (\gamma_V\otimes{\rm id})
\]
\end{proposition}
\begin{proof}
It is well-known that the modules and comodules of a bialgebra form tensor categories in the way described in the Proposition. To show that we can combine this into one category whose objects are modules and comodules one just checks that if $\alpha_V,\alpha_W$ and $\gamma_V,\gamma_W$ satisfy (\ref{cond-comp}), then $\alpha_{V\otimes W}$ and $\gamma_{V\otimes W}$ do so, too.

We will now show that $\Psi$ defines a braiding.

That the $\Psi_{V,W}$ are morphisms was shown in Lemma \ref{lem-cond-comp}. Furthermore, by Lemma \ref{lem-inv}, we know that they are invertible.

To show that $\Psi$ is a natural isomorphism, it remains to show that for all morphisms $f_1:V_1\to W_1$ and $f_2:V_2\to W_2$ the diagram
\begin{equation}\label{diag-comm}
\xymatrix{
V_1\otimes V_2 \ar[d]_{f_1\otimes f_2}\ar[r]^{\Psi_{V_1,V_2}} & V_2\otimes V_1 \ar[d]^{f_2\otimes f_1} \\
W_1\otimes W_2 \ar[r]_{\Psi_{V_1,V_2}}  & W_2\otimes W_1
}
\end{equation}
commutes. But this follows immediately from the definition of $\Psi$ and the fact that the morphisms $f_1,f_2$ are module and comodule maps.

That $\Psi$ satisfies the hexagon axioms
\[
\xymatrix{
& U\otimes (V\otimes W) \ar[r]^{\Psi_{U,V\otimes W}} & (V\otimes W)\otimes U\ar[dr]^{\alpha_{V,W,U}} & \\
(U\otimes V)\otimes W\ar[ur]^{\alpha_{U,V,W}} \ar[dr]_{\Psi_{U,V}\otimes{\rm id}} & & & V\otimes (W\otimes U) \\
& (V\otimes U)\otimes W \ar[r]_{\alpha_{V,U,W}} & V\otimes ( U\otimes W) \ar[ur]_{{\rm id}\otimes \Psi_{U,W}} &
}
\]
and
\[
\xymatrix{
& (U\otimes V)\otimes W \ar[r]^{\Psi_{U\otimes V,W}} & W\otimes(U\otimes V) \ar[dr]^{\alpha^{-1}_{W,U,V}} & \\
U\otimes (V\otimes W) \ar[ur]^{\alpha^{-1}_{U,V,W}} \ar[dr]_{{\rm id}\otimes\Psi_{V,W}} & & & (W\otimes U)\otimes V \\
& U\otimes (W\otimes V)\ar[r]_{\alpha^{-1}_{U,W,V}} & (U\otimes W)\otimes V \ar[ur]_{\Psi_{U,W}\otimes{\rm id}} & \\
}
\]
for all objects $U,V,W$
reduces to the conditions
\begin{eqnarray*}
\Psi_{U,V\otimes W} &=& ({\rm id}_V\otimes \Psi_{U,W})\circ(\Psi_{U,V}\otimes{\rm id}_W) \\
\Psi_{U\otimes V,W} &=& (\Psi_{U,W}\otimes{\rm id}_V)\circ({\rm id}_U\otimes\Psi_{V,W}).
\end{eqnarray*}
since the associativity constraint $\alpha$ is the trivial one, i.e., $\alpha:(u\otimes v)\otimes w\mapsto u\otimes (v\otimes w)$.
We get
\begin{eqnarray*}
&& \Psi_{U,V\otimes W} \\
&=& (\alpha_{V\otimes W}\otimes{\rm id}_U)\circ({\rm id}_V\otimes\tau_{U,V\otimes W})\circ(\gamma_U\otimes{\rm id}_{V\otimes W}) \\
&=& (\alpha_V\otimes\alpha_W\otimes{\rm id}_U)\circ({\rm id}_\mathcal{A}\otimes\tau_{\mathcal{A},V}\otimes{\rm id}_{W\otimes U})\circ(\Delta\otimes {\rm id_{V\otimes W\otimes U}})\circ({\rm id}_V\otimes\tau_{U,V\otimes W}) \\
&& \circ(\gamma_U\otimes{\rm id}_{V\otimes W}) \\
&=&(\alpha_V\otimes\alpha_W\otimes{\rm id}_U)\circ({\rm id}_\mathcal{A}\otimes\tau_{\mathcal{A},V}\otimes {\rm id}_{W\otimes U})\circ ({\rm id}_{\mathcal{A}\otimes\mathcal{A}}\otimes\tau_{U,V\otimes W} ) \\
&&\circ({\rm id}_\mathcal{A}\otimes \gamma_U\otimes{\rm id}_{V\otimes W})\circ(\gamma_U\otimes{\rm id}_{V\otimes W}) \\
&=& ({\rm id}_V\otimes \Psi_{U,W})\circ(\Psi_{U,V}\otimes{\rm id}_W)
\end{eqnarray*}
where we used first the definition of $\alpha_{V\otimes W}$ and then the coaction axiom $({\rm id}\otimes\gamma)\circ\gamma=(\Delta\otimes{\rm id})\circ \gamma$. Similarly, using the definition of $\gamma_{U\otimes V}$ and the action axiom, we can show that the second condition is also satisfied. Therefore $(\mathcal{C}(\mathcal{A}),\Psi)$ is indeed a braided category.
\end{proof}

We are interested in the case where $\mathcal{A}$ and the objects in $(\mathcal{C}(\mathcal{A}),\Psi)$ have an involution, because we want to construct braided $*$-bialgebras.

 \begin{theorem}\label{theo-braided}
   Let $\mathcal{A}$ be an involutive Hopf algebra. Then we can define a braided category $(\mathcal{C}_*(\mathcal{A}),\Psi)$ as follows. The
   objects of $\mathcal{C}_*(\mathcal{A})$ are quadruples $(V,\alpha_V,\gamma_V,*_V)$ consisting of a complex vector space $V$ equipped with a conjugation $*_V$, an
   action $\alpha_V$ and a coaction $\gamma_V$ such Equations (\ref{cond-comp}), (\ref{cond-alpha-*}), and (\ref{cond-gamma-*}) are satisfied.
   The morphism between two objects $(V,\alpha_V,\gamma_V)$ and $(W,\alpha_W,\gamma_W)$ are again the linear maps $f:V\to W$ from $V$ to $W$ that are module and comodule maps.

The involution on the tensor product of two objects is given by
   \[
 {}*_{V\otimes W} =\Psi_{W,V}\circ(*_W\otimes*_V)\circ\tau_{V,W}
\]
\end{theorem}
\begin{proof}
That $*_{V\otimes W}$ is an involution was shown in Lemma \ref{lem-*-inv}. In order to prove that we do indeed get a tensor category in this way, it remains to check that $\alpha_{V\otimes W}$, $\gamma_{V\otimes W}$, and $*_{V\otimes W}$ satisfy again conditions (\ref{cond-alpha-*}) and (\ref{cond-gamma-*}), i.e.
\begin{eqnarray*}
{}*_{V\otimes W}\circ\alpha_{V\otimes W} &=& \alpha_{V\otimes W}\circ (*\otimes *_{V\otimes W})\circ (S\otimes {\rm id}) \\
\gamma_{V\otimes W}\circ *_{V\otimes W} &=& (*\otimes *_{V\otimes W})\circ \gamma_{V\otimes W}.
\end{eqnarray*}
Inserting the definitions and using the fact that $\alpha_V$ and $\alpha_W$ satisfy \eqref{cond-alpha-*}, we get
\begin{gather*}
{}*_{V\otimes W}\circ\alpha_{V\otimes W} \\
= 
(\alpha_V\otimes{\rm id})\circ ({\rm id}\otimes\tau)\circ(\gamma_W\otimes{\rm id})\circ(\alpha_V\otimes\alpha_W)\circ\\
\circ (*\otimes*_V\otimes*\otimes*_W)\circ ({\rm id}\otimes\tau\otimes{\rm id})\circ(S\otimes S\otimes{\rm id})\circ(\Delta\otimes{\rm id}).
\end{gather*}
Using the action property of $\alpha_V$ and that the antipode $S$ is a coalgebra anti-homomorphism, this expression can be seen to be equal to
\[
(\alpha_V\otimes{\rm id})\circ(((m\otimes {\rm id})\circ ({\rm id}\otimes \tau)\circ (\gamma_W\otimes{\rm id}) \circ(\alpha_W\otimes {\rm id})\circ ({\rm id}\otimes \tau)\circ (\Delta\otimes{\rm id}) )\otimes{\rm id})\circ (*\otimes*_W\otimes*_V)\circ(S\otimes \tau).
\]
To this we can now apply the compatibility condition \ref{cond-comp}. Using the action property for $\alpha_v$ now in the opposite direction and reordering everything, we get
\begin{gather*}
(\alpha_V\otimes\alpha_W)\circ({\rm id}\otimes\tau\otimes{\rm id})\circ(\Delta\otimes((\alpha_V\otimes{\rm id})\circ({\rm id}\otimes\tau)\circ(\gamma_W\otimes{\rm id})))\circ \\
\circ(*\otimes*_W\otimes*_V)\circ(S\otimes \tau) \\
 = \alpha_{V\otimes W}\circ (*\otimes *_{V\otimes W})\circ (S\otimes {\rm id}),
\end{gather*}
as required.

The proof for $\gamma_{V\otimes W}$ is similar.
\end{proof}

\begin{remark}
We have $\Psi_{W,V}\circ *_{W\otimes V}\circ\Psi_{V,W}\circ *_{V\otimes W}={\rm id}_{V\otimes W}$ for all objects $V,W$ of $(\mathcal{C}(\mathcal{A}),\Psi)$ (compare to the first Remark following Definition \ref{*-hopf}). This shows that the functor acting as the identity on the objects of $(\mathcal{C}_*(\mathcal{A}),\Psi)$ and taking a morphism $f:V\to W$ to $*_W\circ f\circ *_V:V\to W$ defines an isomorphism between the braided categories $(\mathcal{C}_*(\mathcal{A}),\Psi)$ and $(\mathcal{C}_*(\mathcal{A}),\Psi^{-1})$.
\end{remark}

The subcategory of $(\mathcal{C}_*(\mathcal{A}),\Psi)$ defined by taking the same objects, but only the morphisms which are compatible with the involution, i.e., that satisfy $f\circ *=*\circ f$, is again a tensor category since 
the tensor product of two such morphisms has again the same property,
\begin{eqnarray*}
(f_1\otimes f_2) \circ *_{V_1\otimes V_2} &=&  (f_1\otimes f_2) \circ\Psi_{V_2,V_1}\circ(*_{V_2}\otimes*_{V_1})\circ\tau_{V_1,V_2} \\
&=& \Psi_{W_2,W_1}\circ  (f_2\otimes f_1)\circ(*_{V_2}\otimes*_{V_1})\circ\tau_{V_1,V_2} \\
&=& \Psi_{W_2,W_1}\circ  (*_{W_2}\otimes*_{W_1})\circ(f_2\otimes f_1)\circ\tau_{V_1,V_2} \\
&=&\Psi_{W_2,W_1}\circ  (*_{W_2}\otimes*_{W_1})\circ\tau_{W_1,W_2}\circ(f_1\otimes f_2)\\
&=& *_{W_1\otimes W_2}\circ (f_1\otimes f_2)
\end{eqnarray*}
where we used $(f_1\otimes f_2) \circ\Psi_{V_2,V_1}=\Psi_{W_2,W_1}\circ  (f_2\otimes f_1)$, i.e., the commutativity of diagram (\ref{diag-comm}). It is braided if and only if the $\Psi_{V,W}$ are morphisms of the subcategory, i.e., if  $\Psi\circ *=*\circ \Psi$. But since we also have $\Psi\circ *\circ \Psi\circ *={\rm id}$, this is the case, if and only if $\Psi$ is a symmetry.

   One can check that $(\mathcal{A},{\rm ad},\Delta)$ is an object of the braided category $(\mathcal{C}_*(\mathcal{A}),\Psi)$, where
   \[
   {\rm ad}:\mathcal{A}\otimes \mathcal{A}\to\mathcal{A},\quad {\rm ad}=m\circ(m\otimes{\rm id})\circ({\rm id}\otimes\tau)\circ({\rm id}\otimes S\otimes{\rm id})\circ
   (\Delta\otimes{\rm id})
   \]
   is the adjoint action, since Equations (\ref{cond-comp}), (\ref{cond-alpha-*}), and (\ref{cond-gamma-*}) are satisfied for $\alpha_\mathcal{A}={\rm ad}$, $\gamma_\mathcal{A}=\Delta$ and the involution of $\mathcal{A}$.
   Furthermore, the unit $e:(\mathbb{C},\alpha_\mathbb{C},\gamma_\mathbb{C})\to (\mathcal{A},{\rm ad},\Delta)$ and the multiplication $m:(\mathcal{A},{\rm ad},\Delta)\otimes(\mathcal{A},{\rm
   ad},\Delta)\to(\mathcal{A},{\rm ad},\Delta)$ are morphisms of $(\mathcal{C}_*(\mathcal{A}),\Psi)$.

This can be rephrased by saying that $(\mathcal{A},{\rm ad},\Delta)$ is an algebra in the braided category $(\mathcal{C}_*(\mathcal{A}),\Psi)$.

Thus the tensor product $(\mathcal{B}\otimes\mathcal{A})$ of $\mathcal{A}$ with any other algebra $\mathcal{B}$ in $(\mathcal{C}_*(\mathcal{A}),\Psi)$ is also an algebra in the braided category. Furthermore, by the same arguments as in the first Remark following Definition \ref{*-hopf}, the involution $*_{\mathcal{B}\otimes\mathcal{A}}=\Psi\circ(*\otimes *)\circ\tau$ turns $\mathcal{B}\otimes\mathcal{A}$ into a $*$-algebra.

By dualization, $(\mathcal{A},m,{\rm coad})$ with
\[
{\rm coad}=  (m\otimes{\rm id})\circ({\rm id}\otimes S\otimes{\rm id})\circ({\rm id}\otimes \tau)\circ(\Delta\otimes{\rm id})\circ\Delta
\]
is a coalgebra in the braided category $(\mathcal{C}(\mathcal{A}),\Psi)$. But only Equation (\ref{cond-comp}) is satisfied for $\mathcal{A}$ with the action $m$ and
   the coaction ${\rm coad}$, not (\ref{cond-alpha-*}) and (\ref{cond-gamma-*}), and therefore it is not an object of $(\mathcal{C}_*(\mathcal{A}),\Psi)$.

For us the following situation will be of interest. We have an involutive Hopf algebra $\mathcal{A}$, its category $(\mathcal{C}(\mathcal{A}),\Psi)$ as in Proposition \ref{prop-braided}, i.e., without the involutions, and a braided bialgebra $(\mathcal{B},\Delta,\varepsilon,m,1,\Psi)$ in $(\mathcal{C}(\mathcal{A}),\Psi)$. What are the conditions we have to impose on an involution $*_\mathcal{B}$ to turn $\mathcal{B}$ into a braided $*$-bialgebra?

\begin{corollary}\label{cor-braided}
Let $\mathcal{A}$ be an involutive Hopf algebra and let $(\mathcal{C}_*(\mathcal{A}),\Psi)$ and $(\mathcal{C}(\mathcal{A}),\Psi)$ be the categories constructed in Theorem \ref{theo-braided} and Proposition \ref{prop-braided}. Let $(\mathcal{B},\Delta,\varepsilon,m,1,\Psi)$ be a braided bialgebra in the category $\mathcal{C}(\mathcal{A})$. If we have an involution $*_\mathcal{B}$ of $\mathcal{B}$ such that $\mathcal{B}$ is a $*$-algebra, (\ref{cond-alpha-*}) and (\ref{cond-gamma-*}) are satisfied, and $\Delta$ is a $*$-algebra homomorphisms, then $(\mathcal{B},\Delta,\varepsilon,m,1,\Psi,*_\mathcal{B})$ is a braided $*$-bialgebra in  $(\mathcal{C}_*(\mathcal{A}),\Psi)$.
\end{corollary}

\subsection{Coquasi-triangular bialgebras}

We will consider the special case that shall be used later.
Suppose that we have a coquasi-triangular structure on our bialgebra, i.e.\ a universal R-form on $\mathcal{A}$, i.e.\ a linear functional $r:\mathcal{A}\otimes\mathcal{A}\to\mathbb{C}$ that is invertible w.r.t.\ the convolution product (i.e.\ there exists another functional $\bar{r}$ such that $r\star\bar{r}=\bar{r}\star r=\varepsilon\otimes\varepsilon$) that satisfies
\begin{eqnarray}
m^{\rm op} &=& r\star m\star \bar{r} \\
r_{13}\star r_{23} &=& r\circ (m\otimes {\rm id}) \label{r-left}\\
r_{13}\star r_{12} &=& r\circ ({\rm id}\otimes m), \label{r-right}
\end{eqnarray}
where $r_{12},r_{23},r_{13}:\mathcal{A}\otimes \mathcal{A}\otimes \mathcal{A}\to \mathbb{C}$ are defined by $r_{12}=r\otimes \varepsilon$, $r_{23}=\varepsilon\otimes r$, and $r_{13}= (r\otimes \varepsilon)\circ ({\rm id}\otimes \tau)$. This implies that $r$ also satisfies
\[
r_{12}\star r_{13}\star r_{23}=r_{23}\star r_{13}\star r_{12}.
\]
Furthermore we have $r\circ(1\otimes {\rm id})=r\circ({\rm id}\otimes 1)=\varepsilon$ and, if $\mathcal{A}$ has an antipode, then the antipode is invertible and we have $r\circ(S\otimes {\rm id})=\bar{r}=r\circ({\rm id}\otimes S^{-1})$, cf.\ \cite[Section 2.2]{majid95}. The convolution-inverse of the R-form has similar properties,
\begin{eqnarray*}
\bar{r}_{12}\star \bar{r}_{13}\star \bar{r}_{23} &=& \bar{r}_{23}\star \bar{r}_{13}\star\bar{r}_{12} \\
\bar{r}_{23}\star\bar{r}_{13} &=& \bar{r}\circ (m\otimes {\rm id}) \\
\bar{r}_{12}\star\bar{r}_{13} &=& = \bar{r}\circ ({\rm id}\otimes m),
\end{eqnarray*}

\begin{lemma}\label{r-action}
If $\gamma$ is a coaction of $\mathcal{A}$ on $V$, then $\alpha=(\bar{r}\otimes {\rm id})\circ ({\rm id}\otimes \gamma)$ defines an action of $\mathcal{A}$ on $V$. Moreover, $\alpha$ and $\gamma$ satisfy Equation (\ref{cond-comp}).

Furthermore, if a linear map $f:V\to W$ between two comodules $V,W$ is $\gamma$-invariant, i.e.\ a comodule map, then it is also $\alpha$-invariant.
\end{lemma}
\begin{proof} We get
\[
\alpha(1\otimes v) = (\alpha\otimes {\rm id})\circ (1\otimes {\rm id})\circ \gamma(v) = (\varepsilon\otimes{\rm id})\circ \gamma (v)=v
\]
and
\begin{eqnarray*}
\alpha\circ({\rm id}\otimes \alpha) &=&(\bar{r}\otimes {\rm id})\circ ({\rm id}\otimes \gamma)\circ ({\rm id}\otimes \bar{r}\otimes {\rm id})\circ ({\rm id}\otimes{\rm id}\otimes \gamma) \\
&=&
(\bar{r}\otimes {\rm id})\circ ({\rm id}\otimes \bar{r} \otimes {\rm id}\otimes{\rm id})\circ({\rm id}\otimes {\rm id}\otimes \Delta\otimes {\rm id})\circ ({\rm id}\otimes{\rm id}\otimes \gamma) \\
&=&(\bar{r}_{23}\star\bar{r}_{13}\otimes {\rm id})\circ({\rm id}\otimes {\rm id}\otimes \gamma) \\
&=& (\bar{r}\otimes{\rm id})\circ (m\otimes {\rm id}\otimes{\rm id})\circ ({\rm id}\otimes {\rm id}\otimes \gamma) \\
&=& (\bar{r}\otimes{\rm id})\circ ({\rm id}\otimes\gamma)\circ (m\otimes {\rm id}) = \alpha\circ(m\otimes {\rm id}),
\end{eqnarray*}
i.e.\ $\alpha$ is an action. Let us now check Equation (\ref{cond-comp}),
\begin{eqnarray*}
&&(m\otimes\alpha)\circ({\rm id}\otimes\tau\otimes{\rm id})\circ(\Delta\otimes\gamma) \\
&=&({\rm id}\otimes \bar{r}\otimes {\rm id}) \circ(\Delta\otimes{\rm id}\otimes\gamma)\circ({\rm id}\otimes\tau\otimes{\rm id})\circ(\Delta\otimes\gamma) \\
&=&(m\otimes \bar{r}\otimes{\rm id})({\rm id}\otimes \tau\otimes {\rm id}\otimes{\rm id})\circ(\Delta\otimes \Delta\otimes{\rm id}))\circ({\rm id}\otimes \gamma)\\
&=&
((m\star\bar{r})\otimes{\rm id})\circ({\rm id}\otimes \gamma)\\
&=&
((\bar{r}\star m^{\rm op})\otimes{\rm id})\circ({\rm id}\otimes \gamma)\\
&=&(\bar{r}\otimes (m\circ\tau)\otimes{\rm id})\circ({\rm id}\otimes \tau\otimes {\rm id}\otimes{\rm id})\circ(\Delta\otimes \Delta\otimes{\rm id}))\circ({\rm id}\otimes \gamma)\\
&=&
(\bar{r}\otimes (m\circ\tau)\otimes{\rm id})\circ({\rm id}\otimes \tau\otimes {\rm id}\otimes{\rm id})\circ({\rm id}\otimes{\rm id}\otimes {\rm id}\otimes\gamma))\circ(\Delta\otimes \gamma)\\
&=&(m\otimes {\rm id})\circ({\rm id}\otimes\tau)\circ(\gamma\otimes{\rm id})\circ (\bar{r}\otimes{\rm id}\otimes{\rm id})\circ({\rm id}\otimes\gamma\otimes{\rm id})\circ\\
&&\circ({\rm id}\otimes\tau)\circ(\Delta\otimes{\rm id}) \\
&=&(m\otimes {\rm id})\circ ({\rm id}\otimes \tau)\circ (\gamma\otimes{\rm id}) \circ(\alpha\otimes {\rm id})\circ ({\rm id}\otimes \tau)\circ (\Delta\otimes{\rm id}).
\end{eqnarray*}

For the last statement check
\begin{eqnarray*}
\alpha_W\circ ({\rm id}\otimes f) &=& (\bar{r}\otimes {\rm id})\circ ({\rm id}\otimes \gamma_W)\circ ({\rm id}\otimes f) \\
&=& (\bar{r}\otimes {\rm id})\circ({\rm id}\otimes{\rm id}\otimes f)\circ ({\rm id}\otimes \gamma_V) \\
&=&f\circ \alpha_W.
\end{eqnarray*}
\end{proof}

In this case the braiding $\Psi=(\alpha\otimes{\rm id}\circ({\rm id}\otimes\tau)\circ(\gamma\otimes{\rm id})$ has the form
\[
\Psi_{V,W}= \tau\circ(\bar{r}\otimes{\rm id}\otimes{\rm id})\circ({\rm id}\otimes\tau\otimes{\rm id})\circ(\gamma_V\otimes \gamma_W).
\]
The following two lemmas show that in this special case we actually do not need an invertible antipode to construct the braided tensor category $(\mathcal{C}(\mathcal{A}),\Psi)$.

\begin{lemma}
Let $\mathcal{A}$ be a bialgebra and $r$ a universal R-form on $\mathcal{A}$. Then $\Psi_{V,W}= \tau\circ(\bar{r}\otimes{\rm id}\otimes{\rm id})\circ({\rm id}\otimes\tau\otimes{\rm id})\circ(\gamma_V\otimes \gamma_W)$ is invertible and its inverse is given by $\Psi^{-1}_{V,W}= \tau\circ((r\circ\tau)\otimes{\rm id}\otimes{\rm id})\circ({\rm id}\otimes\tau\otimes{\rm id})\circ(\gamma_W\otimes \gamma_V)$.
\end{lemma}
\begin{proof} This follows from the fact that $\bar{r}$ is the convolution inverse of $r$ and $(\varepsilon\otimes{\rm id})\circ\gamma_V={\rm id}_V$.
\end{proof}
\begin{lemma}
Suppose that $\mathcal{A}$ is a $*$-bialgebra and that there exists a $*$-structure on the comodules $V,W$. If $\gamma_V,\gamma_W$ satisfy Equation (\ref{cond-gamma-*}) and $r$ satisfies
\begin{equation}\label{cond-comp-r-*}
r = \overline{\phantom{a}}\circ\bar{r}\circ(*_\mathcal{A}\otimes*_\mathcal{A}),
\end{equation}
where $\overline{\phantom{a}}$ denotes the complex conjugation of $\mathbb{C}$, then $*_{V\otimes W}=\Psi_{W,V}\circ(*_W\otimes*_V)\circ\tau_{V,W}$ with $\Psi= \tau\circ(\bar{r}\otimes{\rm id}\otimes{\rm id})\circ({\rm id}\otimes\tau\otimes{\rm id})\circ(\gamma\otimes \gamma)$ is an involution.
\end{lemma}
\begin{proof}
Due to Equations  (\ref{cond-gamma-*}) and (\ref{cond-comp-r-*}), we get
\begin{eqnarray*}
&&(*\otimes *)\circ \Psi\circ (*\otimes *) \\
&=&(*\otimes *)\circ \tau\circ(\bar{r}\otimes{\rm id}\otimes{\rm id})\circ({\rm id}\otimes\tau\otimes{\rm id})\circ(\gamma\otimes \gamma)\circ (*\otimes *) \\
&=&\tau\circ(*\otimes *)\circ(\bar{r}\otimes{\rm id}\otimes{\rm id})\circ({\rm id}\otimes\tau\otimes{\rm id})\circ(*_\mathcal{A}\otimes*_V\otimes*_\mathcal{A}\otimes*_V)\circ(\gamma\otimes \gamma) \\
&=&\tau\circ(*\otimes *)\circ((\bar{r}\circ(*_\mathcal{A}\otimes*_\mathcal{A}))\otimes*_V\otimes*_V)\circ({\rm id}\otimes\tau\otimes{\rm id})\circ(\gamma\otimes \gamma) \\
&=& \tau\circ\Psi^{-1}\circ\tau,
\end{eqnarray*}
and therefore
\begin{eqnarray*}
*_{V\otimes W}\circ*_{V\otimes W}&=&\Psi\circ(*\otimes*)\circ\tau\circ\Psi\circ(*\otimes*)\circ\tau\\
&=& \Psi\circ\tau\circ(*\otimes*)\circ\Psi\circ(*\otimes*)\circ\tau \\
&=&\Psi\circ\tau\circ\tau\circ\Psi^{-1}\circ\tau\circ\tau = {\rm id}.
\end{eqnarray*}
\end{proof}

As a consequence of the last three lemmas we see that the comodules of a bialgebra $\mathcal{A}$ with a universal R-form $r$ form also a braided tensor category.

\begin{theorem}\label{theo-r-braided}
Let $\mathcal{A}$ be a coquasi-triangular $*$-bialgebra whose universal R-form satisfies Equation (\ref{cond-comp-r-*}). Then we can construct a braided category $(\mathcal{C}_*(\mathcal{A}),\Psi)$ as follows. The objects of $(\mathcal{C}_*(\mathcal{A}),\Psi)$ are triples $(V,\gamma_V,*_V)$ consisting of a complex vector space $V$, a coaction $\gamma_V$, and an involution $*_V$ such that Equation (\ref{cond-gamma-*}) is satisfied. The morphisms between two objects  $(V,\gamma_V,*_V)$ and  $(W,\gamma_W,*_W)$ are the comodule maps from $V$ to $W$. The tensor product of two objects is defined as in Proposition \ref{prop-braided} and Theorem \ref{theo-braided}. The braiding $\Psi$ is given by
\[
\Psi_{V,W}=\tau\circ(\bar{r}\otimes{\rm id}\otimes{\rm id})\circ({\rm id}\otimes\tau\otimes{\rm id})\circ(\gamma_V\otimes \gamma_W).
\]
\end{theorem}

A similar construction is also possible for quasi-triangular bialgebras, i.e.\ bialgebras equipped with a universal R-matrix, but we shall not use it here.

\subsection{Cocommutative bialgebras}

If the bialgebra $\mathcal{A}$ is cocommutative, then $1\otimes 1$ defines an R-matrix. The corresponding braiding is simply the flip $\tau$. But we can also get more interesting examples in the same way as in Theorem \ref{theo-braided}, see, e.g., \cite{schuermann93}, if we have actions and coactions of $\mathcal{A}$ on the objects $V$ that satisfy the compatibility condition 
\begin{equation}\label{cond-comp-comm}
\gamma\circ\alpha =({\rm ad}\otimes\alpha)\circ({\rm id}\otimes\tau\otimes{\rm id})\circ(\Delta\otimes\gamma)
\end{equation}
given on Page 14 in \cite{schuermann93}, where ${\rm ad}=m\circ({\rm id}\otimes m)\circ({\rm id}\otimes{\rm id}\otimes S)\circ ({\rm id}\otimes\tau)\circ(\Delta\otimes{\rm id})$ is the adjoint action of $\mathcal{A}$ on itself. The following lemma shows that our construction is a generalization of the one presented there.

\begin{lemma}
Let $\mathcal{A}$ be a cocommutative bialgebra. Let $\alpha$ and $\gamma$ be an action and a coaction of $\mathcal{A}$ on some vector space $V$. If $\alpha$ and $\gamma$ satisfy Equation (\ref{cond-comp-comm}), then they satisfy also Equation (\ref{cond-comp}).
\end{lemma}
\begin{proof}
Substituting (\ref{cond-comp-comm}) into the expression $ (m\otimes {\rm id})\circ ({\rm id}\otimes \tau)\circ (\gamma\otimes{\rm id}) \circ(\alpha\otimes {\rm id})\circ ({\rm id}\otimes \tau)\circ (\Delta\otimes{\rm id})$, we get after some simplification
\[
(m^{(4)}\otimes m)\circ ({\rm id}\otimes\tau\otimes\tau\otimes{\rm id})\circ({\rm id}\otimes\tau_{\mathcal{A}\otimes\mathcal{A},\mathcal{A}}\otimes{\rm id})\circ({\rm id}_\mathcal{A}\otimes S\otimes{\rm id})\circ(\Delta^{(4)}\otimes\Delta)
\]
where $\Delta^{(4)}$ denotes $(\Delta\otimes{\rm id})\circ(\Delta\otimes{\rm id})\circ\Delta$. Using the cocommutativity, we can produce a term of the form $m\circ(S\otimes{\rm id})\circ\Delta$, to which we can apply the antipode axiom. Using the unit and counit axiom to clean up the resulting expression, we get the desired result.
\end{proof}

\section{Symmetrization of braided $*$-bialgebras and their L\'evy processes}\label{symmetrization}\markright{SYMMETRIZATION}

Now we will show a construction that will allow later on to define a L\'evy process on some symmetric bialgebra for every L\'evy process on a braided bialgebra $\mathcal{B}$. The idea is to construct a bigger symmetric bialgebra $H$ that contains the braided bialgebra $\mathcal{B}$ as a subalgebra and whose coproduct is related to that of $\mathcal{B}$ in a ``nice'' way. For the case where the braiding is defined through the action and coaction of a group, this construction can be found in \cite[Chapter 3]{schuermann93}. For the general case a similar construction, called bosonization, was introduced by Majid, cf.\ \cite[Section 9.4]{majid95} and the references indicated there. But the role of the involution is not studied there. We will show that if we have a braided $*$-bialgebra in the sense of Definition \ref{*-hopf} in a braided tensor category defined as in Theorem \ref{theo-braided}, then the symmetrization is an involutive bialgebra also and the inclusion is a $*$-algebra homomorphism.

\begin{theorem}\label{theo-symm}
Let $(\mathcal{A},\Delta_\mathcal{A},\varepsilon_\mathcal{A},m_\mathcal{A},1_\mathcal{A},*_\mathcal{A})$ be an involutive Hopf algebra or coquasi-triangular bialgebra and let $(\mathcal{C}_*(\mathcal{A}),\Psi)$ be defined through $\mathcal{A}$ as in Theorem \ref{theo-braided} or \ref{theo-r-braided}. If $(\mathcal{B},\Delta_\mathcal{B},\varepsilon_\mathcal{B},m_\mathcal{B},1_\mathcal{B},\Psi,*_\mathcal{B})$ is a braided $*$-bialgebra in $(\mathcal{C}_*(\mathcal{A}),\Psi)$, then $H=\mathcal{B}\otimes\mathcal{A}$ (as a vector space) becomes a $*$-bialgebra with
\begin{eqnarray*}
m_H &=& (m_\mathcal{B}\otimes m_\mathcal{A})\circ({\rm id}\otimes \alpha\otimes{\rm id}\otimes{\rm id})\circ({\rm id}\otimes{\rm id}\otimes\tau\otimes{\rm id})\circ({\rm id}\otimes\Delta\otimes{\rm id}\otimes{\rm id}),\\
1_H &=& 1_\mathcal{B}\otimes1_\mathcal{A}, \\
\Delta_H &=&({\rm id}\otimes m\otimes{\rm id}\otimes{\rm id})\circ ({\rm id}\otimes{\rm id}\otimes\tau\otimes{\rm id})\circ ({\rm id}\otimes \gamma\otimes{\rm id}\otimes{\rm id})\circ(\Delta_\mathcal{B}\otimes\Delta_\mathcal{A}),\\
\varepsilon_H &=& \varepsilon_\mathcal{B}\otimes\varepsilon_\mathcal{A}, \\{}
*_H&=&(\alpha\otimes{\rm id})\circ (\tau\otimes{\rm id})\circ ({\rm id}\otimes \Delta)\circ (*\otimes*)
\end{eqnarray*}
The map ${\rm id}_\mathcal{B}\otimes1_\mathcal{A}$ defines a $*$-algebra isomorphism between $\mathcal{B}$ and $\mathcal{B}\otimes1_\mathcal{A}\subseteq H$.

Furthermore we have the following commutative diagram
\[
\xymatrix{
\mathcal{B} \ar[d]_{\Delta_\mathcal{B}}\ar[rr]^{{\rm id}_\mathcal{B}\otimes1_\mathcal{A}} & & H \ar[d]^{\Delta_H} \\
\mathcal{B}\otimes\mathcal{B} \ar[rr]_{{\rm id}_\mathcal{B}\otimes\gamma\otimes 1_\mathcal{A}}  & & H\otimes H
}
\]
and
\[
\xymatrix{
\mathcal{B}\ar[d]_{\Delta_\mathcal{B}} & & H \ar[d]^{\Delta_H}\ar[ll]^{{\rm id}_\mathcal{B}\otimes\varepsilon_\mathcal{A}} \\
\mathcal{B}\otimes\mathcal{B} & & H\otimes H \ar[ll]^{{\rm id}_\mathcal{B}\otimes\varepsilon_\mathcal{A}\otimes{\rm id}_\mathcal{B}\otimes\varepsilon_\mathcal{A}}
}
\]
\end{theorem}
\begin{proof}
That $(H,m_H,1_H)$ is an algebra, follows from \cite[Proposition 1.6.6]{majid95}. Since $(\mathcal{A},{\rm ad},\Delta)$ is an object of $(\mathcal{C}_*(\mathcal{A}),\Psi)$, it is even an algebra in $(\mathcal{C}_*(\mathcal{A}),\Psi)$ and a $*$-algebra, see the discussion following Theorem \ref{theo-braided}. Furthermore $(H,\Delta_H,\varepsilon_H)$ is a coalgebra (in $(\mathcal{C}(\mathcal{A}),\Psi)$), since $(\mathcal{A},m,{\rm coad})$ is an object of $(\mathcal{C}(\mathcal{A}),\Psi)$.

The proof that $H$ is a symmetric bialgebra, i.e., that $\Delta_H$ and $\varepsilon_H$ are algebra homomorphisms, does not involve the involutions and is similar to the corresponding proof in \cite[Section 9.4]{majid95}.

In order to show that $H$ is a $*$-bialgebra, we still have to verify that $\Delta_H$ is a $*$-homomorphism, i.e.,
\[
\Delta_H\circ *_H = (*_H\otimes *_H)\circ \Delta_H.
\]
The proof of this property is rather long, it consists in applying the conditions satisfied by the actions and coactions and some of the axioms of symmetric and braided bialgebras in the right order to transform the expression on the left-hand-side
\begin{eqnarray*}
&& \Delta_H\circ *_H \\
&=& ({\rm id}\otimes((\gamma\otimes{\rm id})\circ({\rm id}\otimes\tau)\circ(m\otimes{\rm id}))\otimes{\rm id})\circ(\Delta_\mathcal{B}\otimes\Delta_\mathcal{A})\circ(\alpha\otimes{\rm id})\circ \\
&&\circ(\tau\otimes{\rm id})\circ({\rm id}\otimes\Delta_\mathcal{A})\circ(*_\mathcal{B}\otimes *_\mathcal{A})
\end{eqnarray*}
into the expression on the right-hand-side, which we will write in the form
\begin{eqnarray*}
&&(*_H\otimes *_H)\circ \Delta_H \\
&=& (((\alpha\otimes{\rm id})\circ(\tau\otimes{\rm id})\circ({\rm id}\otimes \Delta_\mathcal{A}))\otimes((\alpha\otimes{\rm id})\circ(\tau\otimes{\rm id})\circ({\rm id}\otimes \Delta_\mathcal{A})))\circ\\
&&\circ(*_\mathcal{B}\otimes*_\mathcal{A}\otimes*_\mathcal{B}\otimes*_\mathcal{A})\circ({\rm id}\otimes((m\otimes{\rm id}))\circ({\rm id}\otimes\tau)\circ(\gamma\otimes{\rm id}))\circ (\Delta_\mathcal{B}\otimes\Delta_\mathcal{A}) \\
&=&(*_\mathcal{B}\otimes*_\mathcal{A}\otimes*_\mathcal{B}\otimes*_\mathcal{A})\\
&&\circ((((\alpha\circ S^{-1})\otimes{\rm id})\circ(\tau\otimes{\rm id})\circ({\rm id}\otimes \Delta_\mathcal{A}))\otimes(((\alpha\circ S^{-1})\otimes{\rm id})\circ(\tau\otimes{\rm id})\circ({\rm id}\otimes \Delta_\mathcal{A})))\\
&&\circ({\rm id}\otimes((m\otimes{\rm id}))\circ({\rm id}\otimes\tau)\circ(\gamma\otimes{\rm id}))\circ (\Delta_\mathcal{B}\otimes\Delta_\mathcal{A})
\end{eqnarray*}
We will outline the main steps. First we use the fact that $\Delta_B$ is a morphism of $(\mathcal{C}_*(\mathcal{A}),\Psi)$ and therefore a module map, to move $\alpha$ past $\Delta_\mathcal{B}$. After some reordering one gets a term
\[
(m\otimes{\rm id})\circ({\rm id}\otimes\tau)\circ(\gamma\otimes{\rm id})\circ(\alpha\otimes{\rm id})\circ({\rm id}\otimes\tau)\circ(\Delta_\mathcal{A}\otimes{\rm id})
\]
inside the resulting expression. Since $\alpha$ and $\gamma$ have to satisfy (\ref{cond-comp}), we can replace this term by
\[
(m\otimes \alpha)\circ({\rm id}\otimes\tau\otimes{\rm id})\circ(\Delta_\mathcal{A}\otimes\gamma).
\]
Next we move the involutions past the two coproducts, using the identities
\[
\Delta_\mathcal{A}\circ*_\mathcal{A} = (*_\mathcal{A}\otimes*_\mathcal{A})\circ \Delta_\mathcal{A}
\]
and
\[
\Delta_\mathcal{B}\circ*_\mathcal{B} = *_{\mathcal{B}\otimes\mathcal{B}}\circ \Delta_\mathcal{B} = (\alpha\otimes{\rm id})\circ(\tau\otimes {\rm id})\circ({\rm id}\otimes \gamma)\circ(*_\mathcal{B}\otimes*_\mathcal{B})\circ\Delta_\mathcal{B}.
\]
Moving the involutions now to the beginning of the expression and using (\ref{cond-alpha-*}), one obtains the expression
\begin{eqnarray*}
&&(*_\mathcal{B}\otimes*_\mathcal{A}\otimes*_\mathcal{B}\otimes*_\mathcal{A})\circ (\alpha\otimes m \otimes \alpha\otimes{\rm id})\circ (S^{-1}\otimes{\rm id}\otimes\tau\otimes S^{-1}\otimes{\rm id}_{\mathcal{B}\otimes\mathcal{A}})\circ\\
&&\circ(\tau\otimes\Delta_\mathcal{A}\otimes \gamma\otimes{\rm id})\circ(\alpha\otimes{\rm id}\otimes\tau\otimes{\rm id})\circ(S^{-1}\otimes {\rm id}\otimes \tau\otimes{\rm id}_{\mathcal{A}\otimes\mathcal{A}})\circ\\
&&\circ(\tau\otimes{\rm id})\circ({\rm id}\otimes\gamma\otimes\Delta_\mathcal{A}\otimes{\rm id})\circ(\Delta_\mathcal{B}\otimes\Delta_\mathcal{A})
\end{eqnarray*}
After applying the action and coaction axiom to replace an $\alpha$ by an $m$ and a $\gamma$ by a $\Delta_\mathcal{A}$, one obtains an expression which contains the term
\[
(m\otimes m)\circ({\rm id}\otimes\tau\otimes{\rm id})\circ(\Delta_\mathcal{A}\otimes\Delta_\mathcal{A}).
\]
If we replace this by
\[
\Delta_\mathcal{A}\circ m,
\]
we obtain the desired expression.

That ${\rm id}_\mathcal{B}\otimes1_\mathcal{A}$ is bijective, is clear, and that it is a $*$-algebra homomorphism follows immediately from the definitions of $m_H$ and $*_H$.

The commutativity of the diagrams can easily be verified using the definition of $\Delta_H$.
\end{proof}

The following proposition is important for symmetrizing L\'evy processes, i.e.\ for finding L\'evy processes on $H$ for a given L\'evy process on $\mathcal{B}$ that have the same distribution.

\begin{proposition}\label{sym-functionals}
The map $F:\mathcal{B}'\to H'$, $F\varphi= \varphi\otimes\varepsilon_{\mathcal A}$ is a unital injective algebra homomorphism with respect to the convolution product. Furthermore, it maps positive (or hermitian, conditionally positive) $\Psi$-invariant functionals on $\mathcal{B}$ to positive (or hermitian, conditionally positive, resp.) functionals on $H$.
\end{proposition}
\begin{proof}
The injectivity is clear, since the map from $H'$ to $B'\cong (B\otimes 1_\mathcal{A})'$ defined by $\psi\mapsto\psi\circ{\rm id}_\mathcal{B}\otimes 1_\mathcal{A}$ is a left inverse of $F$.
It preserves the unit of the two convolution algebras, since $F\varepsilon_\mathcal{B}=\varepsilon_\mathcal{B}\otimes\varepsilon_\mathcal{A}=\varepsilon_H$, and the convolution product, since
\begin{eqnarray*}
&&F\varphi_1\star F\varphi_2 \\
&=& (F\varphi_1\otimes F\varphi_2)\circ \Delta_H \\
&=&(\varphi_1\otimes\varepsilon_\mathcal{A}\otimes\varphi_2\otimes\varepsilon_\mathcal{A})\circ ({\rm id}\otimes m\otimes{\rm id}\otimes{\rm id})\circ ({\rm id}\otimes{\rm id}\otimes\tau\otimes{\rm id})\circ \\
&&\circ ({\rm id}\otimes \gamma\otimes{\rm id}\otimes{\rm id})\circ(\Delta_\mathcal{B}\otimes\Delta_\mathcal{A}) \\
&=&(\varphi_1\otimes\varepsilon_\mathcal{A}\otimes\varphi_2\otimes\varepsilon_\mathcal{A}\otimes\varepsilon_\mathcal{A})\circ ({\rm id}\otimes\gamma\otimes{\rm id}_{\mathcal{A}\otimes\mathcal{A}})\circ (\Delta_\mathcal{B}\otimes\Delta_\mathcal{A}) \\
&=& ((\varphi_1\otimes\varphi_2)\circ\Delta_\mathcal{B})\otimes\varepsilon_\mathcal{A} = F(\varphi_1\star\varphi_2)
\end{eqnarray*}
for all $\varphi_1,\varphi_2\in \mathcal{B}'$.

Assume now that $\varphi$ is positive and $\Psi$-invariant. Let $c=\sum_kb_k\otimes a_k\in\mathcal{B}\otimes\mathcal{A}\cong H$. We have to show that
\begin{eqnarray*}
F\varphi(c^*c) &=& \varphi\otimes\varepsilon_\mathcal{A}\left( \sum_{k,\ell} m_H\left((b_k\otimes a_k)^*\otimes b_\ell\otimes a_\ell\right)\right)
\end{eqnarray*}
is positive. Note that
\[
(b\otimes a)^*=\Psi_{\mathcal{A},\mathcal{B}}(a^*\otimes b^*)
\]
for all $b\in\mathcal{B}$, $a\in\mathcal{A}$, and
\[
m_H= (m_\mathcal{B}\otimes m_\mathcal{A})\circ({\rm id}\otimes\Psi_{\mathcal{A},\mathcal{B}}\otimes{\rm id}) 
\]
Therefore we get, as in the proof of Lemma \ref{pos-conv},
\begin{eqnarray*}
&&F\varphi(c^*c) \\
&=&(\varphi\otimes\varepsilon_\mathcal{A})\circ(m_\mathcal{B}\otimes m_\mathcal{A})\circ({\rm id}\otimes\Psi\otimes{\rm id})\circ(\Psi\otimes{\rm id}\otimes{\rm id})\left(\sum_{k,\ell} m_H(a^*_k\otimes b^*_k\otimes b_\ell\otimes a_\ell\right) \\
&=&\sum_{k,\ell} \varphi(b^*_k b_\ell) \varepsilon_\mathcal{A}(a^*_k a_\ell)
=\sum_{k,\ell} \varphi(b^*_k b_\ell) \overline{\varepsilon_\mathcal{A}(a_k)}\varepsilon_\mathcal{A} (a_\ell),
\end{eqnarray*}
which is positive, since it is the Schur product of two positive definite matrices.
\end{proof}

This proposition gives us a one-to-one correspondence between the class of all convolution semigroups of $\Psi$-invariant states on $\mathcal{B}$ and a class of certain convolution semigroups on its symmetrization $H$. By  \cite[Theorem 1.9.2]{schuermann93}, \cite[Theorem 3.2.8]{schuermann93}, and Theorem \ref{th-gen-proc} these convolution semigroups on $*$-bialgebras and braided $*$-bialgebras are in one-to-one correspondence with L\'evy processes. Therefore we get also a one-to-one correspondence between the L\'evy processes on $\mathcal{B}$ and certain L\'evy processes on $H$. It turns out that this can be used to construct realizations of L\'evy processes on braided $*$-bialgebras.
\begin{theorem}\label{theo-symm-levy}
Let $(j_{st})_{0\le s\le t\le T}$ be a L\'evy process on $\mathcal{B}$ with convolution semigroup $(\varphi_t)_{t\ge 0}$, and let $(j^H_{st})_{0\le s\le t\le T}$ be a L\'evy process on $H\cong\mathcal{B}\otimes\mathcal{A}$ with convolution semigroup $(F\varphi_t)_{t\ge 0}$. Then $(\hat\jmath_{st})_{0\le s\le t\le T}$ with
\[
\hat\jmath_{st}=(j^H_{0s}\otimes j^H_{st})\circ(1_\mathcal{B}\otimes\gamma\otimes 1_\mathcal{A})
\]
defines a L\'evy process on $\mathcal{B}$. Furthermore, $(\hat\jmath_{st})_{0\le s\le t\le T}$ is equivalent to $(j_{st})_{0\le s\le t\le T}$.
\end{theorem}
\begin{remark}
This Theorem generalizes \cite[Theorem 3.3.1]{schuermann93}.
\end{remark}

\begin{proof}
We will use Sweedler's notation $\Delta_\mathcal{B}(b)=\sum b^{(1)}\otimes b^{(2)}\in\mathcal{B}\otimes\mathcal{B}$ and $\gamma(b)=\sum b_{(1)}\otimes b_{(2)}\in\mathcal{A}\otimes\mathcal{B}$ for the coproduct and coaction on an element $b\in\mathcal{B}$. In this notation the coassociativity and the coaction axiom become
\[
\sum  b^{(1)}{}^{(1)}\otimes b^{(1)}{}^{(2)}\otimes b^{(2)} = \sum  b^{(1)}\otimes b^{(2)}{}^{(1)}\otimes b^{(2)}{}^{(2)}
\]
and
\[
\sum  b_{(1)}{}^{(1)}\otimes b_{(1)}{}^{(2)}\otimes b_{(2)} = \sum  b_{(1)}\otimes b_{(2)}{}_{(1)}\otimes b_{(2)}{}_{(2)},
\]
and we can write $\hat\jmath_{st}(b)$ as
\[
\hat\jmath_{st}(b)= j^H_{0s}(b_{(1)})j^H_{st}(b_{(2)}).
\]
Let $0\le r\le s\le t\le T$, then we have
\begin{eqnarray*}
\hat\jmath_{rs}\star\hat\jmath_{st}(b) &=& \hat\jmath_{rs}(b^{(1)})\hat\jmath_{st}(b^{(2)}) \\
&=&
j^H_{0r}\left(1_\mathcal{B}\otimes b^{(1)}{}_{(1)}\right)j^H_{rs}\left(b^{(1)}{}_{(2)}\otimes 1_\mathcal{A}\right)j^H_{0s}\left(1_\mathcal{B}\otimes b^{(2)}{}_{(1)}\right)j^H_{st}\left(b^{(2)}{}_{(2)}\otimes 1_\mathcal{A}\right) \\
&=&j^H_{0r}\left(1_\mathcal{B}\otimes b^{(1)}{}_{(1)}b^{(2)}{}_{(1)}{}^{(1)}\right)j^H_{rs}\left(b^{(1)}{}_{(2)}\otimes b^{(2)}{}_{(1)}{}^{(2)}\right)j^H_{st}\left(b^{(2)}{}_{(2)}\otimes 1_\mathcal{A}\right) \\
&=&j^H_{0r}\left(1_\mathcal{B}\otimes b^{(1)}{}_{(1)}b^{(2)}{}_{(1)}\right)j^H_{rs}\left(b^{(1)}{}_{(2)}\otimes b^{(2)}{}_{(2)}{}_{(1)}\right)j^H_{st}\left(b^{(2)}{}_{(2)}{}_{(2)}\otimes 1_\mathcal{A}\right) \\
\end{eqnarray*}
where so far we have used that $(j^H_{st})_{0\le s\le t\le T}$ is a L\'evy process and the coaction axiom. Now we use the fact that $\Delta_\mathcal{B}$ is a comodule map, which can be written as
\[
b^{(1)}{}_{(1)}b^{(2)}{}_{(1)}\otimes b^{(1)}{}_{(2)}\otimes b^{(2)}{}_{(2)}= b_{(1)}\otimes b_{(2)}{}^{(1)}\otimes b_{(2)}{}^{(2)}
\]
in Sweedler's notation. We get
\[
\hat\jmath_{rs}\star\hat\jmath_{st}(b)= j^H_{0r}\left(1_\mathcal{B}\otimes b_{(1)}\right)j^H_{rs}\left(b_{(2)}{}^{(1)}\otimes b_{(2)}{}^{(2)}{}_{(1)}\right)j^H_{st}\left(b_{(2)}{}^{(2)}{}_{(2)}\otimes 1_\mathcal{A}\right).
\]
The arguments of $j^H_{rs}$ and $j^H_{st}$ are nothing but the two factors of $\Delta_H(b_{(2)}\otimes 1_\mathcal{A})$, and therefore the increment property for $(j^H_{st})_{0\le s\le t\le T}$ implies
\[
\hat\jmath_{rs}\star\hat\jmath_{st}(b)= j^H_{0r}\left(1_\mathcal{B}\otimes b_{(1)}\right)j^H_{rt}\left(b_{(2)}\otimes 1_\mathcal{A}\right)=\hat\jmath_{rt}(b),
\]
i.e., we have shown that $(\hat\jmath_{st})_{0\le s\le t\le T}$ also satisfies the increment property.

We have
\begin{eqnarray*}
\Phi\circ\hat\jmath_{st} &=&\Phi\circ(j^H_{0s}\otimes j^H_{st})\circ(1_\mathcal{B}\otimes\gamma\otimes 1_\mathcal{A}) = (F\varphi_s\otimes F\varphi_{t-s})\circ (1_\mathcal{B}\otimes\gamma\otimes 1_\mathcal{A}) \\
&=& (\varepsilon_\mathcal{A}\otimes \varphi_{t-s})\circ\gamma = \varphi_{t-s},
\end{eqnarray*}
i.e., the processes $(j_{st})_{0\le s\le t\le T}$ and  $(\hat\jmath_{st})_{0\le s\le t\le T}$ have the same marginal distributions. This implies the stationarity and the weak continuity of the increments of $(\hat\jmath_{st})_{0\le s\le t\le T}$ and completes the proof that  $(\hat\jmath_{st})_{0\le s\le t\le T}$ is a L\'evy process. Furthermore, it also establishes the equivalence of the two processes and thus completes the proof of the Theorem.
\end{proof}

L\'evy processes on symmetric $*$-bialgebras can be realized on boson Fock spaces using quantum stochastic differential calculus \cite{hudson+parthasarathy84,parthasarathy92}. The necessary input is a triple $(\rho,\eta,L)$, where $\rho$ is a representation of $\mathcal{A}$ on some pre-Hilbert space $D$, $\eta:\mathcal{A}\to D$ a $\rho$-cocycle, i.e.\
\[
\eta(ab)=\rho(a)\eta(b)-\eta(a)\varepsilon(b),
\]
for $a,b\in\mathcal{A}$, and $L:\mathcal{A}\to \mathbb{C}$ is a hermitian
linear functional s.t.\
\[
\langle\eta(a^*),\eta(b)\rangle =
-\varepsilon(a)L(b)+L(ab)-L(a)\varepsilon(b),
\]
for $a,b\in\mathcal{A}$.

This triple can by constructed from the generator by a GNS-type construction, see \cite[Section 2.3]{schuermann93}.

If we know the triple for a generator $L$ on a braided $*$-bialgebra
$\mathcal{B}$, then the following proposition tells us how to extend it
to a triple for $L^H=FL$.

\begin{proposition}\label{sym-triple}
Let $(j_{st})_{0\le s\le t\le T}$ be a L\'evy process on a braided
$*$-bialgebra $\mathcal{B}$ with $\alpha$-invariant generator $L$ and triple $(\rho,\eta,L)$. Then the triple
$(\rho^H,\eta^H,L^H)$ of the symmetrization $(j^H_{st})_{0\le s\le
  t\le T}$ defined by Theorem \ref{theo-symm-levy} lives on the same
pre-Hilbert space as $(\rho,\eta,L)$ and the restrictions of
$(\rho^H,\eta^H,L^H)$ to $\mathcal{B}$ are equal to $(\rho,\eta,L)$.

Furthermore, $L^H=F(L)=L\otimes \varepsilon_\mathcal{A}$, $\eta^H$
vanishes on $1_\mathcal{B}\otimes\mathcal{A}$ since
\[
\eta^H(a\otimes b)=\varepsilon(a)\eta(b),
\]
and $\rho^H$ is determined by
\[
\rho^H(1\otimes a)\eta(b) = \eta(\alpha(a\otimes b)) \quad\mbox{ and
}\quad\rho^H(1\otimes b)=\rho(b)
\]
for for $a\in\mathcal{A}$, $b\in\mathcal{B}$.
\end{proposition}
\begin{proof}
The pre-Hilbert space $D$ on which the triple acts is obtained by defining a sesqui-linear form on ${\rm ker}\,\varepsilon_H$
\[
\langle u,v\rangle_{L^H} = L(u^*v), \qquad u,v\in{\rm ker}\,\varepsilon_H
\]
and taking the quotient $D={\rm ker}\,\varepsilon_H/\mathcal{N}_{L^H}$ with respect to the null space
\[
\mathcal{N}_{L^H} = \{u\in{\rm ker}\,\varepsilon_H; L^H(u^*u)=0\}.
\]
The cocycle $\eta^H$ is then given by
\[
\eta^H(u) = \overline{u-\varepsilon_H(u)}
\]
for $u\in H$, where $\overline{\phantom{aa}}$ denotes the canonical projection from ${\rm ker}\,\varepsilon_H$ to $D$.

Let $b\otimes a,d\otimes c\in{\rm ker}\,\varepsilon_H$, $b,d\in \mathcal{B}$, $a,c\in\mathcal{A}$, then we get
\begin{eqnarray*}
\langle b\otimes a,d\otimes c\rangle_{L^H} &=& L^H\left((b\otimes a)^*d\otimes c\right) \\
&=& (\varepsilon_\mathcal{A}\otimes L)\left(\left(\alpha\left((a^{(1)})^*\otimes b^*\right)\otimes(a^{(2)})^*\right)d\otimes c\right) \\
&=& (\varepsilon_\mathcal{A}\otimes L)\left(\left(\alpha\left((a^{(1)})^*\otimes
b^*\right) \alpha\left((a^{(2)})^*\otimes d\right)\right) \otimes
(a^{(3)})^* d\right) \\
&=&\varepsilon_\mathcal{A}(c) L\left(\alpha(a^*\otimes b^*d)\right)
= \varepsilon_\mathcal{A}(a^*c)L(b^*d)
\end{eqnarray*}
because $m_\mathcal{B}$ is a module map and $L$ is $\alpha$-invariant. Since $b\otimes a$ and $d\otimes c$ are in the kernel of $\varepsilon_H=\varepsilon_\mathcal{B}\otimes\varepsilon_\mathcal{A}$, we get
\[
\langle b\otimes a,d\otimes c\rangle_{L^H} = \varepsilon_\mathcal{A}(a^*c)L(b^*d) =  \varepsilon_\mathcal{A}(a^*c)\langle\eta(b),\eta(d)\rangle,
\]
which proves $D=\eta(\mathcal{A})$ and implies also the formula for $\eta^H$.

The representation $\rho^H$ is the action on $D$ induced from left multiplication on ${\rm ker}\,\varepsilon_H$, this can also be written as
\[
\rho^H(b\otimes a)\eta(d\otimes c) = \eta\big((b\otimes a)(d\otimes c)\big)
=\eta(b\alpha(a^{(1)}\otimes d)\otimes a^{(2)}c)
\]
for $b\otimes a\in H$, $d\otimes c\in{\rm ker}\,\varepsilon_H$, $b,d\in\mathcal{B}$, $a,c\in\mathcal{A}$.

Choosing $d\in{\rm ker}\,\varepsilon_\mathcal{B}$, $c=1$, and $b\otimes a=b\otimes 1$ or $b\otimes a=1\otimes a$, we obtain the formulas given in the proposition. On the other hand these characterize $\rho^H$ completely, since it has to be a homomorphism.
\end{proof}

\section{A construction of braided $*$-spaces}\label{construction}

In this section we will construct a large class of braided $*$-spaces and their symmetrization.

Let $R=(R^i{}_k{}^j{}_l)_{i,j,k,l=1,\ldots,n}\in \mathbb{C}^{n\times n}\otimes \mathbb{C}^{n\times n}$ be a solution of the quantum Yang-Baxter equation
\[
R_{12} R_{13} R_{23} = R_{23} R_{13} R_{12}
\]
where $R_{12}=R\otimes 1$, etc., i.e.\
\[
\sum_{a,b,c=1}^n R^i{}_a{}^k{}_b R^a{}_j{}^d{}_c R^b{}_l{}^c{}_r =
\sum_{a,b,c=1}^n R^k{}_b{}^d{}_c
R^i{}_a{}^c{}_r R^a{}_j{}^b{}_l
\]
for all $d,i,j,k,l,r=1,\ldots,n$. Usually we will not write the summation symbol, but use the summation convention, all indices that appear twice are summed over (from 1 to $n$), and the identities are supposed to hold for all values (running from 1 to $n$) of the indices that appear only once.

We suppose furthermore that $R$ is of real type I, i.e.
\[
\overline{R^i{}_k{}^j{}_l} = R^l{}_j{}^k{}_i,
\]
and that $R$ is bi-invertible, i.e.\ that there exist matrices $R^{-1},\tilde{R}\in \mathbb{C}^{n\times n}\otimes \mathbb{C}^{n\times n}$, such that
\begin{eqnarray*}
(R^{-1})^i{}_k{}^j{}_l R^k{}_p{}^l{}_q &=& R^i{}_k{}^j{}_l
(R^{-1})^k{}_p{}^l{}_q = \delta^i_p \delta^j_q, \\
\tilde{R}^i{}_k{}^j{}_l R^k{}_p{}^q{}_j &=& R^i{}_k{}^j{}_l
R^k{}_p{}^q{}_j = \delta^i_p \delta^q_l.
\end{eqnarray*}

We denote by $\mathcal{A}(R)=(A,m_A, e_A, \Delta_A,\varepsilon_A)$ the *-bialgebra generated by $a^i_j$, $i,j=1,\ldots,n$, and their adjoints $b^i_j=(a^j_i)^*$ with the relations
\begin{eqnarray}
R^i{}_p{}^k{}_q a^p_j a^q_l &=& a^k_q a^i_p R^p{}_j{}^q{}_l \label{aa-rel}\\
R^i{}_k{}^j{}_l b^k_q a^p_j &=& a^u_l b^i_v R^v{}_q{}^p{}_u \label{ab-rel}\\
\Delta(a^i_j) &=& a^i_k \otimes a^k_j \label{a-copr} \\
\varepsilon(a^i_j) &=& \delta^i_j. \label{a-counit}
\end{eqnarray}
see also \cite{reshetikhin+takhtadzhyan+faddeev90}.

In many cases, e.g.\ for the standard R-matrices, this will even be a *-Hopf algebra, cf.\ \cite{majid95}. The defining relations imply, e.g.,
\begin{eqnarray*}
\Delta(b^i_j) &=& b^k_j\otimes b^i_k \\
R^q{}_k{}^p{}_i b^j_p b^l_q &=& b^q_k b^p_i R^l{}_q{}^j{}_p,
\end{eqnarray*}
and $\varepsilon(b^i_j) = \delta^i_j$.

A representation of $\mathcal{A}(R)$ on $\mathbb{C}^n$ is given by
\begin{eqnarray}
\rho(a^j_l) v^i &=& R^i{}_k{}^j{}_l v^k \label{rep-a}\\
\rho(b^j_l) v^i &=& R^j{}_l{}^i{}_k v^k \label{rep-b}
\end{eqnarray}
since $R$ satisfies the quantum Yang-Baxter equation,
\begin{eqnarray*}
R^i{}_p{}^k{}_q \rho(a^p_j) \rho(a^q_l) v^a &=& R^i{}_p{}^k{}_q
R^b{}_c{}^p{}_j R^a{}_b{}^q{}_l v^c = 
R^b{}_c{}^k{}_q R^a{}_b{}^i{}_p R^p{}_j{}^q{}_l v^c \\
&=& \rho(a^k_q) \rho(a^i_p) R^p{}_j{}^q{}_lv^a, \\
R^i{}_k{}^j{}_l \rho(b^k_q) \rho(a^p_j) v^a &=& R^i{}_k{}^j{}_l
R^k{}_q{}^b{}_c R^a{}_b{}^p{}_j v^c = R^b{}_c{}^u{}_l R^i{}_v{}^a{}_b
R^u{}_q{}^p{}_u v^c \\
&=& \rho(a^u_l) \rho(b^i_v) R^u{}_q{}^p{}_u v^a.
\end{eqnarray*}

We can define a universal R-form on $\mathcal{A}(R)$ by
\[
r(a^i_j\otimes a^k_l)= R^i{}_j{}^k{}_l, \qquad r(a^j_l\otimes b^i_k)=R^i{}_k{}^k{}_l,
\]
on the generators. This definition admits a unique extension to all of  $\mathcal{A}(R)$ such that Conditions \eqref{cond-comp-r-*}, \eqref{r-left}, and \eqref{r-right} are satisfied, i.e.\
\begin{eqnarray*}
r(a^*\otimes b^*) &=& \overline{\bar{r}(a\otimes b)}, \\
r(ab\otimes c) &=& r(a\otimes c_{(1)})r(b\otimes c_{(2)}), \\
r(a\otimes bc) &=& r(a_{(1)}\otimes c)r(a_{(2)}\otimes c),
\end{eqnarray*}
hold for all $a,b,c\in\mathcal{A}(R)$.

We will now construct the free braided $*$-space $V(R)$. As an algebra this is the free algebra generated by $x_1,\ldots,x_n$ and their adjoints $v^1=(x_i)^*,\ldots,v^n=(x_n)^*$. The braiding is defined on the generators by
\begin{eqnarray}
\Psi(x_i\otimes x_j) &=& R^k{}_i{}^l{}_j x_l\otimes x_k, \label{psi-xx}\\
\Psi(x_i\otimes v^j) &=& R^j{}_l{}^k{}_i v^l\otimes x_k, \label{psi-xv}\\
\Psi(v^i\otimes x_j) &=& \tilde{R}^i{}_k{}^l{}_j  x_l\otimes v^k, \label{psi-vx}\\
\Psi(v^i\otimes v^j) &=& (R^{-1})^j{}_l{}^i{}_k v^l\otimes v^k, \label{psi-vv}
\end{eqnarray}
and extended to arbitrary elements by
\begin{eqnarray*}
\Psi(1\otimes u) &=& u\otimes 1, \\
\Psi(u\otimes 1) &=& 1\otimes u, \\
\Psi(u\otimes u_1u_2) &=& ({\rm id}\otimes\Psi)(\Psi(u\otimes u_1)\otimes u_2), \\
\Psi(u_1u_2\otimes u) &=& (\Psi\otimes {\rm id})(u_1\otimes \Psi(u_2\otimes u)),
\end{eqnarray*}
for $u,u_1,u_2\in V(R)$. The coproduct is defined by
\begin{eqnarray*}
\Delta(x_i) &=& x_i\otimes 1+ 1\otimes x_i, \\
\Delta(v^i) &=& v^i\otimes 1+1\otimes v^i,
\end{eqnarray*}
on the generators and extended such that it is a homomorphism from  $V(R)$ to  $V(R)\otimes V(R)$, where the latter is equipped with the multiplication $(m\otimes m)\circ({\rm id}\otimes \Psi\otimes {\rm id})$. The counit is the unit algebra homomorphism $\varepsilon:V(R)\to\mathbb{C}$ with $\varepsilon(1)=1$, $\varepsilon(x_i)=\varepsilon(v^i)=0$.

\subsection{The braided $*$-space $V(R)$ as a left comodule of $\mathcal{A}(R')$}

We show how $V(R)$ can be obtained in two ways as a braided $*$-bialgebra in a category of modules and comodules of a coquasi-triangular $*$-bialgebra of FRT-type. This leads to two different symmetrizations.

Let $R'=\tau(R)$. For the first construction we choose the coquasi-triangular $*$-bialgebra $\mathcal{A}(R')$ and the left coaction $\gamma:V(R)\to \mathcal{A}(R')\otimes V(R)$ defined on the generators by
\[
\gamma(1)=1\otimes 1,\qquad \gamma(x_i)=b^j_k\otimes x_j\,\qquad \gamma(v^i)=a^i_j\otimes v^j,
\]
and extended to general elements such that the multiplication in $V(R)$ is a comodule map, i.e.\
\[
\gamma(uv)=u_{(1)}v_{(1)}\otimes u_{(2)}v_{(2)}
\]
for $u,v\in V(R)$, $\gamma(u)=u_{(1)}\otimes u_{(2)}$, $\gamma(v)=v_{(1)}\otimes v_{(2)}$. Let $r':\mathcal{A}(R')\otimes\mathcal{A}(R')\to\mathbb{C}$ denote the universal R-form on $\mathcal{A}(R')$ determined by $r'(a^i_j\otimes a^k_l)= R'{}^i{}_j{}^k{}_l= R{}^k{}_l{}^i{}_j$, $r'(a^j_l\otimes b^i_k)=R'{}^i{}_k{}^j{}_l=R{}^j{}_l{}^i{}_k$.

The left action $\alpha=(\bar{r'}\otimes {\rm id})\circ ({\rm id}\otimes \gamma)$ introduced in Lemma \ref{r-action}, is given by
\begin{eqnarray*}
\alpha(a^i_k\otimes v^j) &=& \overline{r'}(a^i_k\otimes a^j_l)v^l= ({R'}^{-1})^i{}_k{}^j{}_lv^l = (R^{-1})^j{}_l{}^i{}_k v^l, \\
\alpha(a^i_k\otimes x_j) &=& \tilde{R}{}^i{}_k{}^l{}_j x_l, \\
\alpha(b^i_k\otimes v^j) &=& R^j{}_l{}^i{}_k  v^l, \\
\alpha(b^i_k\otimes x_j) &=& R^i{}_k{}^l{}_j x_l.
\end{eqnarray*}
The braiding in the category of left $\mathcal{A}(R')$-comodules is given by
\[
\Psi= (\overline{r}\otimes \tau)\circ(\gamma\otimes\gamma)=(\alpha\otimes{\rm id})\circ({\rm id}\otimes\tau)\circ(\gamma\otimes {\rm id}).
\]
Let us now give the relations of the (left) symmetrization $H_L=V(R)\otimes\mathcal{A}(R')$ of $V(R)$, cf.\ Theorem \ref{theo-symm}. $H_L$ is the $*$-bialgebra generated by $x_i$, $a^i_j$, and their adjoints $v^i=(x_i)^*$ and $b^i_j=(a^j_i)^*$, $i,j=1,\ldots,n$, with the relations
\begin{eqnarray*}
R^k{}_q{}^i{}_p a^p_j a^q_l &=& a^k_q a^i_p R^q{}_l{}^p{}_j,\\
R^j{}_l{}^i{}_k b^k_q a^p_j &=& a^u_l b^i_v R^p{}_u{}^v{}_q, \\
a^j_k x_i &=& \tilde{R}^j{}_l{}^m{}_i x_m a^l_k, \\
b^j_k x_i &=& R^l{}_k{}^m{}_i x_m b^j_l, \\
a^j_k v^i &=& (R^{-1})^i{}_m{}^j{}_l v^m a^l_l, \\
b^j_k v^i &=& R^i{}_m{}^l{}_k v^m b^j_l, \\
\Delta(a^i_j) &=& a^i_k\otimes a^k_j,\\
\Delta(b^i_j) &=& b^k_j \otimes b^i_k,\\
\Delta(x_i) &=& x_j\otimes 1 + b^j_i\otimes x_j, \\
\Delta(v^i) &=& v^j\otimes 1 + a^i_j \otimes v^j \\
\varepsilon(a^i_j) &=& \varepsilon(b^i_j) = \delta^i_j,\\
\varepsilon(x_i) &=& \varepsilon(v^i) = 0.
\end{eqnarray*}

\subsection{The braided $*$-space $V(R)$ as a right comodule of $\mathcal{A}(R)$}

For the second construction we choose the coquasi-triangular $*$-bialgebra $\mathcal{A}(R)$ and the right coaction $\gamma_R:V(R)\to V(R)\otimes\mathcal{A}(R)$ defined on the generators by
\[
\gamma_R(1)=1\otimes 1,\qquad \gamma_R(x_i)=x_j\otimes a^j_i,\qquad \gamma_R(v^i)=v^j\otimes b^i_j,
\]
and again extended to general elements such that the multiplication in $V(R)$ is a comodule map, i.e.\
\[
\gamma_R(uv)=u_{(1)}v_{(1)}\otimes u_{(2)}v_{(2)}
\]
for $u,v\in V(R)$, $\gamma_R(u)=u_{(1)}\otimes u_{(2)}$, $\gamma_R(v)=v_{(1)}\otimes v_{(2)}$.
An analog of Lemma \ref{r-action} holds also for right coactions, i.e.\ $\alpha_R=({\rm id}\otimes r)\circ(\gamma_R\otimes{\rm id})$ defines a right action of $\mathcal{A}(R)$ on $V(R)$. On the generators this right action is given by
\begin{eqnarray*}
\alpha_R(x_i\otimes a^k_j) &=& R^l{}_i{}^k{}_j x_l, \\
\alpha_R(x_i\otimes b^k_j) &=& R^k{}_j{}^l{}_i x_i, \\
\alpha_R(v^i\otimes a^k_k) &=& \tilde{R}^i{}_l{}^k{}_j v^l, \\
\alpha_R(v^l\otimes b^k_j) &=& (R^{-1})^k{}_j{}^i{}_l v^l.
\end{eqnarray*}
The braiding in the category of right $\mathcal{A}(R)$-comodules is defined by
\[
\Psi_R=(\tau\otimes r)\circ(\gamma_R\otimes\gamma_R)=({\rm id}\otimes \alpha_R)\circ(\tau\otimes{\rm id})\circ({\rm id}\otimes \gamma_R).
\]
Just like in Theorem \ref{theo-symm}, we can also define a symmetrization for braided $*$-bialgebras in the category of right $\mathcal{A}(R)$-comodules. Here the underlying vector space is $H_R=\mathcal{A}(R)\otimes V(R)$ and the operations are defined by
\begin{eqnarray*}
m_{H} &=& (m_\mathcal{A}\otimes m_V)\circ({\rm id}\otimes {\rm id}\otimes\alpha_R\otimes{\rm id})\circ({\rm id}\otimes\tau\otimes{\rm id}\otimes{\rm id})\circ({\rm id}\otimes{\rm id}\otimes\Delta_\mathcal{A}\otimes{\rm id}), \\
1_H &=& 1_\mathcal{A}\otimes 1_V, \\
\Delta_H &=& ({\rm id}\otimes {\rm id}\otimes m_\mathcal{A}\otimes{\rm id})\circ({\rm id}\otimes\tau\otimes{\rm id}\otimes{\rm id})\circ({\rm id}\otimes{\rm id}\otimes\gamma_R\otimes{\rm id})\circ(\Delta_\mathcal{A}\otimes\Delta_V), \\
\varepsilon_H &=& \varepsilon_\mathcal{A}\otimes \varepsilon_V, \\
*_H &=& ({\rm id}\otimes\alpha_R)\circ({\rm id}\otimes\tau)\circ(\Delta\otimes{\rm id})\circ(*\otimes *).
\end{eqnarray*}
We will call $H_R$ the right symmetrization of $V(R)$. It is the $*$-bialgebra generated by $x_i$, $a^i_j$, and their adjoints $v^i=(x_i)^*$ and $b^i_j=(a^j_i)^*$, $i,j=1,\ldots,n$, with the relations
\begin{eqnarray*}
R^i{}_p{}^k{}_q a^p_j a^q_l &=& a^k_q a^i_p R^p{}_j{}^q{}_l,\\
R^i{}_k{}^j{}_l b^k_q a^p_j &=& a^u_l b^i_v R^v{}_q{}^p{}_u, \\
x_ia^j_k &=& R^m{}_i{}^l{}_k a^i_l x_m, \\
x_ib^j_k &=& R^j{}_l{}^m{}_i b^l_k x_m, \\
v^i a^j_k &=& \tilde{R}{}^i{}_k{}^j{}_m a^i_l v^m, \\
v^i b^j_k &=& (R^{-1})^j{}_l{}^i{}_m b^l_k v^m, \\
\Delta(a^i_j) &=& a^i_k\otimes a^k_j,\\
\Delta(b^i_j) &=& b^k_j \otimes b^i_k,\\
\Delta(x_i) &=& x_j\otimes a^j_i + 1\otimes x_i,\\
\Delta(v^i) &=& v^j\otimes b^i_j + 1\otimes v^i,\\
\varepsilon(a^i_j) &=& \varepsilon(b^i_j) = \delta^i_j,\\
\varepsilon(x_i) &=& \varepsilon(v^i) = 0.
\end{eqnarray*}

The results of Theorem \ref{theo-symm-levy} and Propositions \ref{sym-functionals} and \ref{sym-triple} for the symmetrization of L\'evy processes on $V(R)$ can be formulated for the right symmetrization, too.

\section{L\'evy processes on braided $*$-spaces}\label{levy-proc}

In this section we will show that there always exists a L\'evy process on the braided $*$-spaces constructed in the previous section that can be considered as a standard Brownian motion on these spaces.

\begin{definition}
A linear functional $\phi:\mathcal{B}\to\mathbb{C}$ on a (braided) bialgebra $(\mathcal{B}, \Delta,\varepsilon,m,1)$ is called quadratic, if it satisfies
\[
\phi(abc)=0
\]
for all $a,b,c\in\mathcal{B}$ with $\varepsilon(a)=\varepsilon(b)=\varepsilon(c)=0$, see also \cite[Section 5.1]{schuermann93}.

A L\'evy process whose generator is quadratic is called a Brownian motion.
\end{definition}

For the rest of this section $R$ will denote a fixed bi-invertible R-matrix of real type I, and $V$ the associated free braided *-space. For explicite calculations we will use the basis $\mathcal{B}$ consisting of the words in the generators $x_1,\ldots,x_n,v^1,\ldots,v_n$. Let $L:V\to \mathbb{C}$ be the functional defined by $L(x_iv^j) = \delta^j_i$ on basis elements of the form $x_i v^j$, and zero on all other basis elements.

\begin{proposition}
The functional $L$ is quadratic, $\Psi^{-1}$-invariant, hermitian, and conditionally positive (i.e.\ positive on ker $\varepsilon_V$).
\end{proposition}
\begin{proof}
Equations (\ref{psi-vx}) and (\ref{psi-xx}) imply
\begin{eqnarray*}
({\rm id}\otimes L)\circ\Psi(x_iv^j\otimes x_k) &=&  \tilde{R}^j{}_l{}^m{}_k R^r{}_i{}^p{}_m x_p L(x_rv^l) = \tilde{R}^j{}_l{}^m{}_k R^l{}_i{}^p{}_m x_p = \delta^j_i \delta^p_k x_p \\
&=& x_k L(x_iv^j).
\end{eqnarray*}
Similarly, using Equations (\ref{psi-vv}) and (\ref{psi-xv}), we get
\[
({\rm id}\otimes L)\circ\Psi(x_iv^j\otimes v^k) = v_k L(x_iv^j) ,
\]
and we see that $L$ is $\Psi^{-1}$ invariant.
\end{proof}

We can now carry out the construction described in \cite[Chapter 2]{schuermann93} to obtain quantum stochastic differential equations for the symmetrization of the process $(j_{st})_{0\le s\le t}$ associated to the generator $L$. Theorem \ref{theo-symm-levy} shows how the process on $H$ with generator $L_H$ can be used to construct a process on $V(R)$ that has generator $L$.

Let
\[
\mathcal{N}=\left\{u\in V; L\Big((u-\varepsilon(u)1)^* (u-\varepsilon(u)1)\Big)=0\right\},
\]
and $K=V(R)/\mathcal{N}$. Then $K$ is a Hilbert space with the inner product induced by $\langle u,v\rangle=L\Big((u-\varepsilon(u)1)^* (v-\varepsilon(v)1)\Big)$. It turns out that $K$ is isomorphic to $\mathbb{C}^n$, and for the canonical projection $\eta_L:H\to K$ we find that $\{\eta_L(v^i)\}$ forms an orthonormal basis of $K$, and $\eta_L(x_i)=0$.

The other ingredient we need is the induced action $\rho_L$ of $H$ on $K$. Here we get $\rho_L(x_i)=0$ and $\rho(v^i)=0$.

One verifies that $L$ is $\alpha_R$-invariant for the right action $\alpha_R=({\rm id}\otimes r)\circ(\gamma_R\otimes{\rm id})$, i.e.\
\[
L(\alpha(u\otimes a))=\varepsilon_\mathcal{A}(a)L(u)
\]
for all $a\in\mathcal{A}(R)$, $u\in V(R)$. Therefore we can get the triple on the right symmetrization $H_R$ as in Proposition \ref{sym-triple}.

Denote by $N$ the number operator on the Fock space $\Gamma (L^2 (\mathbb{R}_+ , K))$ over $L^2 (\mathbb{R}_+ , K)$, and let ${\rm dom} (\alpha ^ N )$ be the domain of the self-adjoint operator $\alpha ^ N$, i.e.
\begin{equation}
{\rm dom} (\alpha^N ) = \{ F \in \Gamma (L^2 (\mathbb{R}_+ , K)) : \sum_{n \in \mathbb{N}} \alpha^{2n} |F^{(n)}| ^ 2 < \infty \} .
\end{equation}
We consider the dense linear subspace
\begin{equation}
\mathcal{E}_K = \bigcap _{\alpha \in \mathbb{R}_+} {\rm dom} (\alpha^N )
\end{equation}
of $\Gamma (L^2 (\mathbb{R}_+ , K ))$. Next let $\underline{\mathcal{A}}_K$ be the space of families $(\underline{k}_{st} )_{0 \leq s\leq t}$ of linear operators on $\mathcal{E}_K$ given by kernels as defined in [Sch93], Section 2.4.

We put $j^H_{st} = (j^H_s \circ S ) \star j^H_t$, i. e. $j^H_{st}$ is the increment process associated with $j^H_t$. Sometimes we call $j^H_{st}$ the L\'evy process rather than $j^H_t$.

\begin{theorem}\label{qsde}
Let $j^H_{st}$ be the L\'evy process on $H_R$ with generator
$L_H=\varepsilon_\mathcal{A}\otimes V(R)$. Then a realization of the
right symmetrization $j^H_{st}$ on the Fock space $\Gamma (L^2
(\mathbb{R}_+ , K ))$ is given by the unique solution in
$\underline{\mathcal{A}}_K$ of the quantum stochastic differential equations
\[
{\rm d} j^H_{st} (b) = (j^H_{st} \star {\rm d} I_t )(b) ; \  j^H_{ss} (b) = \epsilon (b),
\]
where the differential ${\rm d} I_t$ is given by ${\rm d}I_t = {\rm d} \Lambda \circ (\rho-1\cdot\varepsilon) + {\rm d} A^+ \circ \eta + {\rm d} A \circ \tilde{\eta} + L_H {\rm d} t$, and $\tilde\eta=\eta\circ *$.
If we put $X_i (t) = j^H_t (x_i )$, $A_l^k (t) = j^H_t (a_l^k )$, $B_l^k (t) = (A_k^l (t))^* = j_t (b_l^k)$, these equations can be written in the following way:
\begin{eqnarray*}
{\rm d}X_i &=& X_j {\rm d}\Lambda\left((R^k{}_l{}^j{}_i-\delta^k_l \delta^j_i)_{1\le k,l\le n}\right) + {\rm d} A_i \\
{\rm d}X^*_i &=& X^*_j {\rm d}\Lambda\left((R^i{}_j{}^k{}_l- \delta^k_l \delta^j_i)_{1\le k,l\le n}\right) + {\rm d} A^*_i \\
{\rm d}A^i_j &=& A^i_m {\rm d}\Lambda\left( (R^k{}_l{}^m{}_j-\delta^k_l \delta^m_j)_{1\le k,l\le n}\right) \\
{\rm d}B^i_j &=& B^m_j {\rm d}\Lambda\left((R^i{}_m{}^k{}_l-\delta^k_l \delta^i_m)_{1\le k,l\le n}\right)
\end{eqnarray*}
\end{theorem}
\begin{proof} The construction preceding the theorem is exactly the one described in \cite[Chapter 2]{schuermann93} to get the triple $(\rho_L,\eta_L, L_H)$ that is used in \cite[Theorem 2.3.5]{schuermann93} to formulate the quantum stochastic differential equations for $j_{st}$.
\end{proof}

\section{Examples}\label{example}

\subsection{The one-dimensional R-matrix $R=(q)$}

Let us first consider the one-dimensional R-matrix $R=(q)$. For $q\in \mathbb{R}$, $q\not=0$, this is a bi-invertible R-matrix of real type I and defines therefore a braided $*$-space $V(q)$. As an algebra, $V(q)$ is the free algebra generated by $x$ and $x^*=v$. We will use the words in $x$ and $v$ as a basis for $V(q)$. The braiding is given by
\[
\begin{array}{rclrcl}
   \Psi(x\otimes x)&=&q x\otimes x, & \Psi(x\otimes v ) &=& q v\otimes x,\\
   \Psi(v\otimes x)&=&q^{-1} x\otimes v, & \Psi(v\otimes v ) &=& q^{-1} v\otimes x.\\
\end{array}
\]
A functional $L:V(q)\to \mathbb{C}$ has to satisfy
\begin{eqnarray*}
L(x) x &=& \Psi\circ ({\rm id} \otimes L)(x\otimes x) = (L\otimes{\rm id})\circ \Psi(x\otimes x) = q L(x) x \\
   L(v) x &=& \Psi\circ ({\rm id} \otimes L)(v\otimes x) = (L\otimes{\rm id})\circ \Psi(v\otimes x) = q^{-1} L(v) x \\
   L(xx) x &=& \Psi\circ ({\rm id} \otimes L)(xx\otimes x) = (L\otimes{\rm id})\circ \Psi(xx\otimes x) = q^{2} L(xx) x \\
   L(xv) x &=& \Psi\circ ({\rm id} \otimes L)(xv\otimes x) = (L\otimes{\rm id})\circ \Psi(xv\otimes x) = L(xv) x \\
   L(vx) x &=& \Psi\circ ({\rm id} \otimes L)(vx\otimes x) = (L\otimes{\rm id})\circ \Psi(vx\otimes x) = L(vx) x \\
   L(vv) x &=& \Psi\circ ({\rm id} \otimes L)(vv\otimes x) = (L\otimes{\rm id})\circ \Psi(vv\otimes x) = q^{-2} L(vv) x 
   \end{eqnarray*}
   and as well as a similar set of equations for $v$ to be $\Psi$-invariant. Thus, for $q^2\not=1$, a $\Psi$-invariant quadratic functional can have
   non-zero values only on $xv$ and $vx$.

   A quadratic functional on $V(q)$ is conditionally positive if and only if the matrix
   \[
   \left(\begin{array}{cc}
   L(xv) & L(xx) \\
   L(vv) & L(vx)
   \end{array}\right)
   \]
   is positive semi-definite. It is hermitian, if we also have $L(v)=\overline{L(x)}$. Thus we get the following classification for the quadratic
   generators on $V(q)$.

   \begin{theorem}\label{theo-class-gen-1}
   A quadratic functional $L:V(q)\to\mathbb{C}$ on $V(q)$ is a generator of a L\'evy process on $V(q)$ if and only if
   \begin{enumerate}
   \item
   for $q=1$: $\left(\begin{array}{cc}
   L(xv) & L(xx) \\
   L(vv) & L(vx)\end{array}\right)$ is positive semi-definite and $L(v)=\overline{L(x)}$,
   \item
   for $q=-1$: $\left(\begin{array}{cc}
   L(xv) & L(xx) \\
   L(vv) & L(vx)\end{array}\right)$ is positive semi-definite and $L(x)=L(v)=0$.
   \item
   for $q^2\not=1$: $L(xv)\ge 0$, $L(vx)\ge 0$, and $L$ vanishes on all other basis elements. 
   \end{enumerate}
   \end{theorem}

   The symmetrization or bosonization of this braided $*$-space gives for $\mathcal{A}(q)$ the free commutative algebra with group-like generator $a$
   and its adjoint $b$. $H$ is generated by $a$, $x$ and their adjoints $b=a^*$ and $v=x^*$. The algebraic relations are
   \[
   ab=ba, \qquad xa=qax, \quad bx=qxb,
   \]
   and the coalgebraic relations are
   \[
   \Delta(a)=a\otimes a, \qquad \Delta(x)=x\otimes a + 1\otimes x.
   \]
   Let $L$ now be the generator with $L(xv)=1$ and $L(u)=0$ on all other basis elements. The construction of the triplet gives $K=\mathbb{C}$,
   $\eta_L(x)=\eta_L(a)=\eta_L(b)=0$, $\eta_L(v)=1$, and $\rho_L(a)=\rho_L(b)=q$, $\rho_L(x)=\rho_L(v)=0$. Thus we get
   \begin{eqnarray*}
   {\rm d}X &=& X{\rm d}\Lambda(q-1) + {\rm d}A(1) \\
   {\rm d}V &=& V{\rm d}\Lambda(q-1) + {\rm d}A^*(1) \\
   {\rm d}A &=& A {\rm d}\Lambda(q-1) \\
   {\rm d}B &=& B {\rm d}\Lambda(q-1)
   \end{eqnarray*}
   for the processes $X(t)=j_t(x)$, $V(t)=j_t(v)$, $A(t)=j_t(a)$, $B(t)=j_t(b)$. The solution of this system of quantum stochastic differential
   equations is the quantum Az\'ema martingale, cf. \cite{parthasarathy90,schuermann91b}, in particular, $Z(t)=X(t)+V(t)$ has the same distribution as
   the classical Az\'ema martingale with parameter $\beta=q-1$ defined in \cite{emery89}.

\subsection{The $sl_2$-R-matrix}

Let $R_2$ be the R-matrix of the standard two-dimensional quantum plane, i.e.
\[
R_2=\left(\begin{array}{cccc}
   q^2 & 0 & 0 & 0 \\
   0 & q & q^2-1 & 0 \\
   0 & 0 & q & 0 \\
   0 & 0 & 0 & q^2
\end{array}\right)
\]
Then $tR_2$ with $q,t\in\mathbb{R}$, $q,t\not=0$ is bi-invertible and of real type I, we can therefore define a braided $*$-space $V(tR_2)$ for it. As an algebra this is the free algebra generated by $x_1$, $x_2$ and there adjoints $x^*_1=v^1$, $x^*_2=v^2$. We will use the words in these four elements as a basis of $V(tR_2)$. It turns out that the $\Psi$-invariance restricts very much the possible generators.

\begin{proposition}\label{prop-inv}
Let $L$ be a quadratic functional on $V(tR_2)$, $q,t\not=0$. Then $L$ is characterized by $A=(A_{ij})= (L(x_ix_j))$, $B=(B_i{}^j)= (L(x_iv^j))$, $C=(C^i{}_j)=(L(v^ix_j))$, $D=(D^{ij})=(L(v^iv^j))$, $a=(a_i)= (L(x_i))$, and $b=(b^i)=(L(v^i))$. The functional $L$ is $\Psi$-invariant if and only if
   \begin{enumerate}
   \item
   for $q=1$ and
   \begin{itemize}
   \item[a)]
   $t=1$: no restrictions on $L$,
   \item[b)]
   $t=-1$: $a=b=0$,
   \item[c)]
   $t^2\not=1$: $A=D=0$ and $a=b=0$,
   \end{itemize}
   \item
   for $q=-1$ and
   \begin{itemize}
   \item[a)]
   $t=\pm 1$: $A$, $B$, $C$, and $D$ are diagonal, $a=b=0$,
   \item[b)]
   $t^2\not=1$: $B$ and $C$ are diagonal, $A=D=0$, $a=b=0$,
   \end{itemize}
   \item
   for $q^2\not=1$ and
   \begin{itemize}
   \item[a)]
   $t^2q^3=1$: $q A_{12}+ A_{21}=0$, $B_1{}^1=B_2{}^2$, $C^1{}_1=q^2C^2{}_2$, $D^{12}+qD^{21}=0$, and all other coefficients vanish,
   \item[b)]
   $t^2q^3\not=1$: $B_1{}^1=B_2{}^2$, $C^1{}_1=q^2C^2{}_2$, and all other coefficients vanish.
   \end{itemize}
   \end{enumerate}
   \end{proposition}
   \begin{proof}
   Let $L$ be an arbitrary quadratic functional and set $L(x_ix_j)=A_{ij}$, $L(x_iv^j)=B_i{}^j$, $L(v^ix_j)=C^i{}_j$, $L(v^iv^j)=D^{ij}$,
   $L(x_i)=a_i$, and $L(v^i)=b^i$. The functional $L$ is $\Psi$-invariant if and only if the following equations are satisfied
   \[
   \begin{array}{rclrcl}
   A_{ij}\delta^l_k &=&t^2 A_{n_3n_1} R^{n_3}{}_i{}^l{}_{n_2} R^{n_1}{}_j{}^{n_2}{}_k, & A_{ij}\delta^k_l &=& t^2 A_{n_3n_1} R^{n_2}{}_l{}^{n_3}{}_{i}
   R^{k}{}_{n_2}{}^{n_1}{}_j, \\[2mm]
   B_i{}^j\delta^l_k &=& B_{n_3}{}^{n_1} R^{n_3}{}_i{}^l{}_{n_2} \tilde{R}^{j}{}_{n_1}{}^{n_2}{}_k, & B_i{}^j\delta^k_l &=& B_{n_3}{}^{n_1}
   R^{n_2}{}_l{}^{n_3}{}_{i} (R^{-1})^{k}{}_{n_2}{}^{j}{}_{n_1}, \\[2mm]
   C^i{}_j\delta^l_k &=& C^{n_3}{}_{n_1} \tilde{R}^{i}{}_{n_3}{}^l{}_{n_2} R^{n_1}{}_{j}{}^{n_2}{}_k, & C^i{}_j\delta^k_l &=& C^{n_3}{}_{n_1}
   (R^{-1})^{n_2}{}_{l}{}^i{}_{n_3} R^{k}{}_{n_2}{}^{n_1}{}_j, \\[2mm]
   D^{ij}\delta^l_k &=&t^{-2} D^{n_3n_1} \tilde{R}^{i}{}_{n_3}{}^l{}_{n_2} \tilde{R}^{j}{}_{n_1}{}^{n_2}{}_k, & D^{ij}\delta^k_l &=&t^{-2} D^{n_3n_1}
   (R^{-1})^{n_2}{}_{l}{}^i{}_{n_3} (R^{-1})^{k}{}_{n_2}{}^{j}{}_{n_1}, \\[2mm]
   a_i\delta^k_j &=& t R^{n_1}{}_i{}^k{}_j a_{n_1}, & a_i\delta^j_k &=& tR^j{}_k{}^{n_1}{}_i a_{n_1}, \\[2mm]
   b^i\delta^k_j &=& t^{-1} \tilde{R}^i{}_{n_1}{}^k{}_j b^{n_1}, & b^i\delta^j_k &=& t^{-1}(R^{-1})^{j}{}_k{}^i{}_{n_1} b^{n_1} ,
\end{array}
\]
for all $i,j,k,l=1,\ldots,n$. These equations follow directly from the invariance condition, for the first equation, e.g., we apply $\Psi\circ (L\otimes {\rm id}) = ({\rm id}\otimes L)\circ \Psi$ to $x_jx_j\otimes x_k$. Solving this system of linear equations (using, e.g., a computer program for symbolic computations like Maple or Mathematica) one arrives at the results listed in the proposition.
\end{proof}

The functional $L$ is conditionally positive semi-definite, if and only if the matrix $\left(L\left(\begin{array}{cc} x_iv^j & x_ix_j \\ v^iv^j & v^ix_j\end{array}\right)\right) = \left(\begin{array}{cc} B & A \\ D & C\end{array}\right)$ is positive semi-definite. For $L$ to be also hermitian, we need to
   impose furthermore $a_i=L(x_i)=\overline{L(v^i)}= \overline{b^i}$, for $i=1,\ldots,n$. This leads to the following classification of the generators.

\begin{theorem}\label{theo-class-gen-2}
Suppose $q^2\not=1$.
\begin{itemize}
\item[a)]
If $t^2q^3=1$, then all $\Psi$-invariant generators on $V(tR_2)$ are of the form
\[
\left(L\left(\begin{array}{cc} x_iv^j & x_ix_j \\ v^iv^j & v^ix_j \end{array}\right)\right) =
\left(\begin{array}{cccc} 
b & 0 & 0 & -qa \\
0 & b & a & 0 \\
0 & \overline{a} & q^2 c & 0 \\
-q\overline{a} & 0 & 0 & c
\end{array}\right),
\]
and $L(x_i)=L(v^i)=0$ for $i=1,\ldots,n$, where $b\ge 0$, $c\ge 0$, $bc\ge q^2 |a|^2$, and $bc\ge q^{-2}|a|^2$.
\item[b)]
If $t^2q^3\not=1$, then all $\Psi$-invariant generators on $V(tR_2)$ are of the form
\[
\left(L\left(\begin{array}{cc} x_iv^j & x_ix_j \\ v^iv^j & v^ix_j \end{array}\right)\right) =
\left(\begin{array}{cccc} 
b & 0 & 0 & 0 \\
0 & b & 0 & 0 \\
0 & 0 & q^2 c & 0 \\
0 & 0 & 0 & c
\end{array}\right),
\]
and $L(x_i)=L(v^i)=0$ for $i=1,\ldots,n$, where $b\ge 0$, and $c\ge 0$.
\end{itemize}
\end{theorem}
\begin{proof}
\begin{itemize}
\item[a)]
Proposition \ref{prop-inv}.3 a) implies that a $\Psi$-invariant functional $L:V(tR_2)\to\mathbb{C}$ is of the form
\[
\left(L\left(\begin{array}{cc} x_iv^j & x_ix_j \\ v^iv^j & v^ix_j \end{array}\right)\right) =
\left(\begin{array}{cccc} 
b & 0 & 0 & -qa \\
0 & b & a & 0 \\
0 & d & q^2 c & 0 \\
-qd & 0 & 0 & c
\end{array}\right),
\]
and $L(x_i)=L(v^i)=0$, where $a=L(x_2x_1)$, $b=L(x_1v^1)$, $c=L(v^2x_2)$, and $d=v^1x_2$. Such an invariant functional $L$ is a generator, if and only if this matrix is positive semi-definite. This is the case if and only if the matrices $\left(\begin{array}{cc} b & -qa \\ -qd & c \end{array}\right)$ and $\left(\begin{array}{cc} b & a \\ d & q^2c \end{array}\right)$ are positive semi-definite, which leads immediately to the conditions given in the theorem.
\item[b)]
Proposition \ref{prop-inv}.3 b) shows us that if $t^2q^3\not=1$, then $L$ is $\Psi$-invariant if and only if we also have $a=d=0$.
\end{itemize}
\end{proof}

   If we take the generator $L$ defined in part b) of the preceding Theorem with $b=1$ and $c=0$, then we get the quantum stochastic differential
   equations
   \begin{eqnarray*}
   {\rm d} X_1 &=& X_1 {\rm d}\Lambda \left( \begin{array}{cc}tq^2-1 & 0 \\
       0 & tq-1\end{array}\right) + X_2 {\rm d} \Lambda
   \left(\begin{array}{cc} 0 & 0 \\ t(q^2-1) & 0 \end{array}\right) + {\rm
     d}A\left( \begin{array}{c} 1 \\ 0 \end{array}\right), \\
   {\rm d}X_2 &=& X_2 {\rm d}\Lambda \left(\begin{array}{cc} tq-1 & 0 \\ 0 &
       tq^2-1 \end{array}\right) + {\rm d}A\left( \begin{array}{c} 0 \\ 1
     \end{array}\right), \\ 
   {\rm d} X^*_1 &=& X^*_1 {\rm d}\Lambda \left( \begin{array}{cc}tq^2-1 & 0 \\
       0 & tq-1 \end{array}\right) + X^*_2 {\rm d} \Lambda
   \left(\begin{array}{cc} 0 & t(q^2-1) \\ 0 & 0 \end{array}\right) + {\rm
     d}A^*\left( \begin{array}{c} 1 \\ 0 \end{array}\right), \\
   {\rm d}X^*_2 &=& X^*_2 {\rm d}\Lambda \left(\begin{array}{cc} tq-1 & 0 \\ 0 &
       tq^2-1 \end{array}\right) + {\rm d}A^*\left( \begin{array}{c} 0 \\ 1
     \end{array}\right),
   \end{eqnarray*}
and
\begin{eqnarray*}
{\rm d}A^1_1 &=& A^1_1 {\rm d}\Lambda \left( \begin{array}{cc}tq^2-1 & 0 \\
0 & tq-1\end{array}\right) + A^1_2 {\rm d} \Lambda
\left(\begin{array}{cc} 0 & 0 \\ t(q^2-1) & 0 \end{array}\right), \\
{\rm d}A^1_2 &=& A^1_2 {\rm d}\Lambda \left(\begin{array}{cc} tq-1 & 0 \\ 0 &
tq^2-1 \end{array}\right)
\\
{\rm d}A^2_1 &=& A^2_1 {\rm d}\Lambda \left( \begin{array}{cc}tq^2-1 & 0 \\
0 & tq-1\end{array}\right) + A^2_2 {\rm d} \Lambda
\left(\begin{array}{cc} 0 & 0 \\ t(q^2-1) & 0 \end{array}\right), \\
{\rm d}A^2_2 &=& A^2_2 {\rm d}\Lambda \left(\begin{array}{cc} tq-1 & 0 \\ 0 &
tq^2-1 \end{array}\right),
\end{eqnarray*}
since $\rho_L(a^1_1) = t\left( \begin{array}{cc}q^2 & 0 \\ 0 & q\end{array}\right)$, $\rho_L(a^1_2) = t\left(\begin{array}{cc} 0 & 0 \\ 0 & 0 \end{array} \right)$, $\rho_L(a^2_1) = t\left(\begin{array}{cc} 0 & 0 \\ q^2-1 & 0 \end{array}\right)$, and $\rho_L(a^2_2) = t \left(\begin{array}{cc} q & 0 \\ 0 & q^2 \end{array}\right)$.

\begin{remark}
In general, we do not get closed structure equations for these processes. E.g., applying the It\^o formula to compute ${\rm d}Z_2\cdot {\rm d} Z_2$ for ${\rm d} Z_2 = {\rm d}X_2 + {\rm d} X^*_2$, we obtain
\[
{\rm d}Z_2 \cdot {\rm d} Z_2 = {\rm d}t + (tq^2-1) Z_2\left\{ {\rm d}A^*\left( \begin{array}{c} 0 \\ 1 \end{array}\right) + {\rm d}A\left( \begin{array}{c} 0 \\ 1 \end{array}\right)\right\} +
Z_2^2{\rm d} \Lambda \left(\begin{array}{cc} (tq-1)^2 & 0 \\ 0 &
(tq^2-1)^2 \end{array}\right),
\]
which can not be expressed only in terms of the differentials ${\rm d}t$, ${\rm d}Z_1$, and ${\rm d}Z_2$, except for $t=0$, which is excluded since $tR$ has to be invertible, and for $t=\frac{1}{q}$, we get
\begin{eqnarray*}
{\rm d}Z_2 \cdot {\rm d} Z_2 &=& {\rm d}t + (q-1) Z_2\left\{ {\rm d}A^*\left( \begin{array}{c} 0 \\ 1 \end{array}\right) + {\rm d}A\left( \begin{array}{c} 0 \\ 1 \end{array}\right)\right\} + Z_2^2{\rm d} \Lambda \left(\begin{array}{cc} 0 & 0 \\ 0 & (q-1)^2 \end{array}\right)\\
&=&{\rm d}t + (q-1) Z_2{\rm d}Z_2.
\end{eqnarray*}
But even then this is not possible for ${\rm d}Z_1\cdot {\rm d}Z_1$.
\end{remark}

\subsection{The $sl_3$-R-matrix}

Let now
\[
R_3=
\left(\begin{array}{ccccccccc}
q^2 & 0 & 0 & 0 & 0 & 0 & 0 & 0 & 0 \\
0 & q & 0 & q^2-1 & 0 & 0 & 0 & 0 & 0 \\
0 & 0 & q & 0 & 0 & 0 & q^2-1 & 0 & 0 \\
0 & 0 & 0 & q & 0 & 0 & 0 & 0 & 0 \\
0 & 0 & 0 & 0 & q^2 & 0 & 0 & 0 & 0 \\
0 & 0 & 0 & 0 & 0 & q & 0 & q^2-1 & 0 \\
0 & 0 & 0 & 0 & 0 & 0 & q & 0 & 0 \\
0 & 0 & 0 & 0 & 0 & 0 & 0 & q & 0 \\
0 & 0 & 0 & 0 & 0 & 0 & 0 & 0 & q^2
\end{array}\right)
\]
be the $sl_3$-R-matrix. We get a similar classification for the generators on $V(tR_3)$ as in the previous subsection, but there are no additional generators for the special case $t^2q^2=1$ (as in Theorem \ref{theo-class-gen-2} c)).

\begin{theorem}\label{theo-class-gen-3}
Suppose $q^2\not=1$. Then all $\Psi$-invariant generators on $V(tR_3)$ are of the form
\[
\left(L\left(\begin{array}{cc} x_iv^j & x_ix_j \\ v^iv^j & v^ix_j \end{array}\right)\right) =
\left(\begin{array}{cccccc} 
b & 0 & 0 & 0 & 0 & 0 \\
0 & b & 0 & 0 & 0 & 0 \\
0 & 0 & b & 0 & 0 & 0 \\
0 & 0 & 0 & q^4c & 0 & 0 \\
0 & 0 & 0 & 0 & q^2 c & 0 \\
0 & 0 & 0 & 0 & 0 & c
\end{array}\right),
\]
and $L(x_i)=L(v^i)=0$ for $i=1,\ldots,n$, where $b\ge 0$, and $c\ge 0$.
\end{theorem}
\begin{proof}
We first determine all invariant quadratic functionals and then we check positivity, as in Proposition \ref{prop-inv} and Theorem
\ref{theo-class-gen-2} in the previous subsection.
\end{proof}

\section{Conclusion}\label{conclusion}

We have shown how L\'evy processes with quadratic generators, i.e.\ Brownian motions, can be defined for a large class of braidings and braided spaces. These processes are a natural generalization of the quantum Az\'ema processes.

We have shown on three well-known examples how the quadratic invariant generators and thus the Brownian motions on braided $*$-spaces can be classified. It turned out that $\Psi$-invariance is a very strong condition that restricts very much the possible Brownian motions, see Theorems \ref{theo-class-gen-1}, \ref{theo-class-gen-2} and \ref{theo-class-gen-3}.

Let us briefly indicate two directions of further generalization.

Other, more general forms of the compatibility condition of the involution with the coproduct exist and have been studied (see e.g.\ \cite{buffenoir+roche99}). Is it possible to study L\'evy processes on $*$-bialgebras where the coproduct is not a $*$-homomor\-phism?

In the theory of quasi-triangular Hopf algebras the trivial commutativity constraint $\tau:u\otimes v\mapsto v\otimes u$ is replaced by a more general braiding. Similarly, one can allow other associativity constraints than the trivial one: $\Phi: (x\otimes u)\otimes v \mapsto x\otimes (u\otimes v )$. This leads to quasi-Hopf algebras. It should be possible to extend the theory of L\'evy processes to these algebras, but this would require a theory of involutive quasi-Hopf algebras (see also \cite{majid93c}).

%% file: mall1.tex
\chapter{Malliavin Calculus and Skorohod Integration for Quantum Stochastic Processes}\label{ch-malliavin}\markboth{MALLIAVIN CALCULUS AND SKOROHOD INTEGRATON}{MALLIAVIN CALCULUS AND SKOROHOD INTEGRATON}

\vspace*{3cm}

\begin{quote}
A derivation operator and a divergence operator are defined on the algebra of bounded operators on the symmetric Fock space over the complexification of a real Hilbert space $\eufrak{h}$ and it is shown that they satisfy similar properties as the derivation and divergence operator on the Wiener space over $\eufrak{h}$. The derivation operator is then used to give sufficient conditions for the existence of smooth Wigner densities for pairs of operators satisfying the canonical commutation relations. For $\eufrak{h}=L^2(\mathbb{R}_+)$, the divergence operator is shown to coincide with the Hudson-Parthasarathy quantum stochastic integral for adapted integrable processes and  with the non-causal quantum stochastic integrals defined by Lindsay and Belavkin for integrable processes.
\end{quote}

\vspace*{3cm}

Joint work with R\'emi L\'eandre and Ren\'e Schott. Published in Inf.\ Dim.\ Anal., Quantum Prob.\ and Rel.\ Topics Vol.~4, No.~1, pp.~11-38, 2001. 

\newpage

\section{Introduction}

Infinite-dimensional analysis has a long history: it began in the sixties (work of Gross \cite{gross67}, Hida, Elworthy, Kr\'ee, $\ldots$), but it is Malliavin \cite{malliavin78} who has applied it to diffusions in order to give a probabilistic proof of H\"ormander's theorem. Malliavin's approach needs a heavy functional analysis apparatus, as the Ornstein-Uhlenbeck operator and the definition of suitable Sobolev spaces, where the diffusions belong. Bismut \cite{bismut81} has given a simpler approach based upon a suitable choice of the Girsanov formula, which gives quasi-invariance formulas. These are differentiated, in order to get integration by parts formulas for the diffusions, which where got by Malliavin in another way.

Our goal is to generalize the hypoellipticity result of Malliavin for non-com\-mutative quantum processes, by using Bismut's method, see also \cite{franz+leandre+schott99}. For that we consider the case of a non-commutative Gaussian process, which is the couple of the position and momentum Brownian motions on Fock space, and we consider the vacuum state. We get an algebraic Girsanov formula, which allows to get integration by parts formulas for the Wigner densities associated to the non-commutative processes, when we differentiate. This allows us to show that the Wigner functional has a density which belongs to all Sobolev spaces over $\mathbb{R}^2$. Let us remark that in general the density is not positive.

If we consider the deterministic elements of the underlying Hilbert space of the Fock space, the derivation of the Girsanov formula leads to a gradient operator satisfying some integration by parts formulas. This shows it is closable as it is in classical infinite-dimensional analysis. But in classical infinite-dimensional analysis, especially in order to study the Malliavin matrix of a functional, we need to be able to take the derivation along a random element of the Cameron-Martin space. In the commutative set-up, this does not pose any problem. Here, we have some difficulty, which leads to the definition of a right-sided and a left-sided gradient, which can be combined to a two-sided gradient.

We can define a divergence operator as a kind of adjoint of the two-sided gradient for cylindrical (non-commutative) vector fields, but since the vacuum state does not define a Hilbert space, it is more difficult to extend it to general (non-commutative) vector fields.

We show that the non-commutative differential calculus contains in some sense the commutative differential calculus.

In the white noise case, i.e.\ if the underlying Hilbert space is the $L^2$-space of some measure space, the classical divergence operator defines an anticipating stochastic integral, known as the Hitsuda-Skorohod integral. We compute the matrix elements between exponential vectors for our divergence operator and use them to show that the divergence operator coincides with the non-causal creation and annihilation integrals defined by Belavkin \cite{belavkin91a,belavkin91b} and Lindsay \cite{lindsay93} for integrable processes, and therefore with the Hudson-Parthasarathy \cite{hudson+parthasarathy84} integral for adapted processes.

\section{Analysis on Wiener space}\label{comm Wiener}

Let us first briefly recall a few definitions and facts from analysis on Wiener space, for more details see, e.g., \cite{janson97,malliavin97,nualart95,nualart98,ustunel95}. Let $\eufrak{h}$ be a real separable Hilbert space. Then there exists a probability space $(\Omega,\mathcal{F},\mathbb{P})$ and a linear map $W:\eufrak{h}\to L^2(\Omega)$ such that the $W(h)$ are centered Gaussian random variables with covariances given by
\[
\mathbb{E}\big(W(h)W(k)\big) = \langle h,k\rangle, \qquad \text{ for all } h,k\in\eufrak{h}.
\]
Set $\mathcal{H}_1=W(\eufrak{h})$, this is a closed Gaussian subspace of $L^2(\Omega)$ and $W:\eufrak{h}\to\mathcal{H}_1\subseteq L^2(\Omega)$ is an isometry. We will assume that the $\sigma$-algebra $\mathcal{F}$ is generated by the elements of $\mathcal{H}_1$. We introduce the algebra of bounded smooth functionals
\[
\mathcal{S}=\{F=f\big(W(h_1),\ldots,W(h_n)\big)| n\in\mathbb{N}, f\in C^\infty_b(\mathbb{R}^n), h_1,\ldots,h_n\in\eufrak{h}\},
\]
and define the derivation operator $\tilde{D}:\mathcal{S}\to L^2(\Omega)\otimes \eufrak{h}\cong L^2(\Omega;\eufrak{h})$ by
\[
\tilde{D}F=\sum_{i=1}^n \frac{\partial f}{\partial x_i}\big(W(h_1),\ldots,W(h_n)\big)\otimes h_i
\]
for $F=f\big(W(h_1),\ldots,W(h_n)\big)\in\mathcal{S}$. Then one can verify the following properties of $\tilde{D}$.
\begin{enumerate}
\item
$\tilde{D}$ is a derivation (w.r.t.\ the natural $L^\infty(\Omega)$-bimodule structure of $L^2(\Omega;\eufrak{h})$), i.e.\
\[
\tilde{D}(FG)= F(\tilde{D}G)+(\tilde{D}G)F, \qquad \text{ for all } F,G\in\mathcal{S}.
\]
\item
The scalar product $\langle h, \tilde{D}F\rangle $ coincides with the Fr\'echet derivative
\[
\tilde{D}_hF = \left.\frac{{\rm d}}{{\rm d}\varepsilon}\right|_{\varepsilon=0} f\big(W(h_1)+\varepsilon\langle h,h_1\rangle,\ldots,W(h_n)+\varepsilon\langle h,h_n\rangle\big)
\]
for all $F=f\big(W(h_1),\ldots,W(h_n)\big)\in\mathcal{S}$ and all $h\in\eufrak{h}$.
\item
We have the following integration by parts formulas,
\begin{eqnarray}
\mathbb{E}\big(FW(h)\big) &=& \mathbb{E}\big(\langle h, \tilde{D}F\rangle\big) \label{int by parts1}\\
\mathbb{E}\big(FGW(h)\big) &=& \mathbb{E}\big(\langle h, \tilde{D}F\rangle G+F \langle h, \tilde{D}G\rangle \big) \label{int by parts2}
\end{eqnarray}
for all $F,G\in\mathcal{S}$, $h\in\eufrak{h}$.
\item
The derivation operator $\tilde{D}$ is a closable operator from $L^p(\Omega)$ to $L^p(\Omega;\eufrak{h})$ for $1\le p\le \infty$. We will denote its closure again by $\tilde{D}$.
\end{enumerate}
We can also define the gradient $\tilde{D}_u F= \langle u,DF\rangle$ w.r.t.\ $\eufrak{h}$-valued random variables $u\in L^2(\Omega;\eufrak{h})$, this is $L^\infty(\Omega)$-linear in the first argument and a derivation in the second, i.e.\
\begin{eqnarray*}
\tilde{D}_{Fu}G &=& F \tilde{D}_uG, \\
\tilde{D}_u (FG) &=& F (\tilde{D}_uG)+ (\tilde{D}_uF)G.
\end{eqnarray*}

$L^2(\Omega)$ and $L^2(\Omega;\eufrak{h})$ are Hilbert spaces (with the obvious inner products), therefore the closability of $\tilde{D}$ implies that it has an adjoint. We will call the adjoint of $\tilde{D}:L^2(\Omega)\to L^2(\Omega;\eufrak{h})$ the divergence operator and denote it by $\tilde\delta:L^2(\Omega;\eufrak{h})\to L^2(\Omega)$. Denote by \[
\mathcal{S}_\eufrak{h}=\left\{u=\sum_{j=1}^n F_j \otimes h_j \Big| n\in\mathbb{N},F_1,\ldots,F_n\in \mathcal{S}, h_1,\ldots,h_n\in\eufrak{h}\right\} 
\]
the smooth elementary $\eufrak{h}$-valued random variables, then $\tilde\delta(u)$ is given by
\[
\tilde\delta(u)=\sum_{j+1}^n F_j W(h_j) - \sum_{j+1}^n \langle h_j, \tilde{D}F_j\rangle
\]
for $u=\sum_{j=1}^n F_j\otimes h_j\in\mathcal{S}_\eufrak{h}$. If we take, e.g., $\eufrak{h}=L^2(\mathbb{R}_+)$, then $B_t=W(\mathbf{1}_{[0,t]})$ is a standard Brownian motion, and the $\eufrak{h}$-valued random variables can also be interpreted as stochastic processes indexed by $\mathbb{R}_+$. It can be shown that $\tilde\delta(u)$ coincides with the It\^o integral $\int_{\mathbb{R}_+} u_t{\rm d}W_t$ for adapted integrable processes. In this case the divergence operator is also called the Hitsuda-Skorohod integral.

The derivation operator and the divergence operator satisfy the following relations
\begin{eqnarray}
\tilde{D}_h\big(\delta(u)) &=& \langle h, u\rangle + \tilde\delta(\tilde{D}_hu), \label{comm rel1}\\
\mathbb{E}\big( \tilde\delta(u)\tilde\delta(v)\big) &=& \mathbb{E}\big(\langle u,v\rangle \big) + \mathbb{E}\big( {\rm Tr}(\tilde{D}u\circ \tilde{D}v)\big), \label{skorohod cov1}\\
\tilde\delta(Fu) &=& F \tilde\delta(u) - \langle u, \tilde{D}F\rangle, \label{DFu1}
\end{eqnarray}
for $h\in\eufrak{h}$, $u,v\in\mathcal{S}_\eufrak{h}$, $F\in\mathcal{S}$.
Here $\tilde{D}$ is extended in the obvious way to $\eufrak{h}$-valued random variables, i.e.\ as $\tilde{D}\otimes{\rm id}_\eufrak{h}$. Thus $\tilde{D}u$ is an $\eufrak{h}\otimes\eufrak{h}$-valued random variable and can also be interpreted as a random variable whose values are (Hilbert-Schmidt) operators on $\eufrak{h}$. If $\{e_j;j\in\mathbb{N}\}$ is a complete orthonormal system on $\eufrak{h}$, then ${\rm Tr}(\tilde{D}u\circ \tilde{D}v)$ can be computed as
${\rm Tr}(\tilde{D}u\circ \tilde{D}v) = \sum_{i,j=1}^\infty \tilde{D}_{e_i}\langle u, e_j\rangle \tilde{D}_{e_j}\langle v,e_i\rangle$.

\section{The non-commutative Wiener space}\label{non comm Wiener}

Let again $\eufrak{h}$ be a real separable Hilbert space and let $\eufrak{h}_\mathbb{C}$ be its complexification. Then we can define a conjugation $\overline{\phantom{a}}:\eufrak{h}_\mathbb{C}\to \eufrak{h}_\mathbb{C}$ by $\overline{h_1+ih_2}=h_1-ih_2$ for $h_1,h_2\in\eufrak{h}_\mathbb{C}$. This conjugation satisfies $\left\langle \overline{h},\overline{k}\right\rangle = \overline{\langle h,k\rangle} = \langle k,h\rangle$ for all $h,k\in\eufrak{h}_\mathbb{C}$. The elements of $\eufrak{h}$ are characterized by the property $\overline{h}=h$, we will call them real. 

Let $\eufrak{H}=\Gamma_s(\eufrak{h}_\mathbb{C})$ be the symmetric Fock space over $\eufrak{h}_\mathbb{C}$, i.e.\ $\eufrak{H}= \bigoplus_{n\in\mathbb{N}} \eufrak{h}_\mathbb{C}^{\odot n}$, where `$\odot$' denotes the symmetric tensor product, and denote the vacuum vector $1+0+\cdots$ by $\Omega$. It is well-known that the symmetric Fock space is isomorphic to the complexification of the Wiener space $L^2(\Omega)$ associated to $\eufrak{h}$ in Section \ref{comm Wiener}. We will develop a calculus on the non-commutative probability space $(\mathcal{B}(\eufrak{H}),\mathbb{E})$, where $\mathbb{E}$ denotes the state defined by $\mathbb{E}(X)=\langle \Omega,X\Omega\rangle$ for $X\in \mathcal{B}(\eufrak{H})$. To emphasize the analogy with the analysis on Wiener space we call $(\mathcal{B}(\eufrak{H}),\mathbb{E})$ the non-commutative Wiener space over $\eufrak{h}$.

The exponential vectors $\{\mathcal{E}(k)= \sum_{n=0}^\infty \frac{k^{\otimes n}}{\sqrt{n!}}; k\in\eufrak{h}_\mathbb{C}\}$ are total in $\eufrak{H}$, their scalar product is given by
\[
\langle\mathcal{E}(k_1),\mathcal{E}(k_2)\rangle = e^{\langle k_1,k_2\rangle}.
\]

We can define the operators $a(h),a^+(h),Q(h),P(h)$ (annihilation, creation, position, momentum) and $U(h_1,h_2)$ with $h,h_1,h_2\in\eufrak{h}_\mathbb{C}$ on $\eufrak{H}$, see, e.g., \cite{biane93,meyer95,parthasarathy92}. The creation and annihilation operators $a^+(h)$ and $a(h)$ are closed, unbounded, mutually adjoint operators. The position and momentum operators
\[
Q(h)= \big(a(\overline{h})+a^+(h)\big), \quad \text{ and }\quad
P(h)= i\big(a(\overline{h})-a^+(h)\big)
\]
are self-adjoint, if $h$ is real. 

The commutation relations of creation, annihilation, position, and momentum are
\[
\begin{array}{rclcrcl}
[a(h),a^+(k)]&=& \langle h,k\rangle, &&
[a(h),a(k)]&=&[a^+(h),a^+(k)]=0, \\ {}
[Q(h),Q(k)]&=& [P(h),P(k)] = 0, &&
[P(h),Q(k)]&=& 2 i \langle \overline{h},k\rangle.
\end{array}
\]
The Weyl operators $U(h_1,h_2)$ can be defined
by $U(h_1,h_2)= \exp\big(iP(h_1)+iQ(h_2)\big)=\exp i \big(a(\overline{h_2}-i\overline{h_1})+ a^+(h_2-ih_1)\big)$, they satisfy
\[
U(h_1,h_2) U(k_1,k_2) = \exp i\big( \langle \overline{h}_2,k_1\rangle - \langle \overline{h}_1,k_2\rangle\big) U(h_1+h_2,k_1+k_2)
\]
Furthermore we have $U(h_1,h_2)^*=U(-\overline{h}_1,-\overline{h}_2)$ and $U(h_1,h_2)^{-1}=U(-h_1,-h_2)$. We see that $U(h_1,h_2)$ is unitary, if $h_1$ and $h_2$ are real. These operators act on the vacuum $\Omega=\mathcal{E}(0)$ as
\[
U(h_1,h_2)\Omega = \exp\left( -\frac{\langle\overline{h}_1,h_1\rangle +\langle\overline{h}_2,h_2\rangle}{2} \right)\mathcal{E}\left(h_1+ih_2\right)
\]
and on general exponential vectors $\mathcal{E}\left(f\right)= \sum_{n=0}^\infty \frac{f^{\otimes n}}{\sqrt{n!}}$ as
\[
U(h_1,h_2)\mathcal{E}\left(f\right) = \exp\left(-\langle \overline{f},h_1+ih_2\rangle -\frac{\langle\overline{h}_1,h_1\rangle +\langle\overline{h}_2,h_2\rangle}{2}\right)\mathcal{E}\left(f+h_1+ih_2\right).
\]
The operators $a(h),a^+(h),Q(h),P(h)$ and $U(h_1,h_2)$ are unbounded, but their domains contain the exponential vectors. We will want to compose them with bounded operators on $\eufrak{H}$, to do so we adopt the following convention. Let
\begin{eqnarray*}
\mathcal{L}\big(\mathcal{E}(\eufrak{h}_\mathbb{C}),\eufrak{H}\big)&=&\Big\{B\in{\rm Lin}\,\big({\rm span}(\mathcal{E}(\eufrak{h}_\mathbb{C})),\eufrak{H}\big)\Big| \exists B^* \in{\rm Lin}\,\big({\rm span}(\mathcal{E}(\eufrak{h}_\mathbb{C})),\eufrak{H}\big)\\
 &&\quad\text{ s.t. } \big\langle\mathcal{E}(f),B\mathcal{E}(g)\big\rangle = \big\langle B^*\mathcal{E}(f),\mathcal{E}(g)\big\rangle \text{ for all }f,g\in\eufrak{h}_\mathbb{C}\Big\},
\end{eqnarray*}
i.e.\ the space of linear operators that are defined on the exponential vectors and that have an ``adjoint'' that is also defined on the exponential vectors. Obviously $a(h),a^+(h),Q(h),P(h), U(h_1,h_2)\in \mathcal{L}\big(\mathcal{E}(\eufrak{h}_\mathbb{C}),\eufrak{H}\big)$. We will say that an expression of the form
$\sum_{j=1}^n X_j B_j Y_j$
with $X_1,\ldots,X_n,Y_1,\ldots,Y_n\in\mathcal{L}\big(\mathcal{E}(\eufrak{h}_\mathbb{C}),\eufrak{H})$ and $B_1,\ldots,B_n\in \mathcal{B}(\eufrak{H})$ defines a bounded operator on $\eufrak{H}$, if there exists a bounded operator $M\in \mathcal{B}(\eufrak{H})$ such that
\[
\big\langle\mathcal{E}(f),M\mathcal{E}(g)\big\rangle = \sum_{j=1}^n\big\langle X_j^*\mathcal{E}(f),B_jY_j\mathcal{E}(g)\big\rangle
\]
holds for all $f,g\in\eufrak{h}_\mathbb{C}$. If it exists, this operator is unique, because the exponential vectors are total in $\eufrak{H}$. We will then write
\[
M=\sum_{j=1}^n X_j B_j Y_j.
\]

\section{Weyl calculus}

\begin{definition}\label{def-weyl}
Let $h=(h_1,h_2)\in \eufrak{h}\otimes \mathbb{R}^2$. We set
\begin{eqnarray}
{\rm Dom}\, O_h &=& \Big\{\varphi:\mathbb{R}^2\to\mathbb{C}\Big|\exists M\in\mathcal{B}(\eufrak{H}),\forall k_1,k_2\in\eufrak{h}_\mathbb{C}:\langle\mathcal{E}(k_1),M\mathcal{E}(k_2)\rangle = \nonumber \\
&&\qquad\qquad\frac{1}{2\pi}\int\langle\mathcal{E}(k_1),U(uh_1,vh_2)\mathcal{E}(k_2)\rangle \mathcal{F}^{-1}\varphi(u,v){\rm d}u{\rm d}v\Big\} \label{defin-O_h}
\end{eqnarray}
and for $\varphi\in{\rm Dom}\,O_h$ we define $O_h(\varphi)$ to be the bounded operator $M$ appearing in Equation (\ref{defin-O_h}), it is uniquely determined due to the totality of $\{\mathcal{E}(k): k\in\eufrak{h}_\mathbb{C}\}$.
\end{definition}
We take the Fourier transform $\mathcal{F}$ as
\[
\mathcal{F}\varphi(u,v)= \frac{1}{2\pi} \int_{\mathbb{R}^2} \varphi(x,y) \exp\big(i(ux+vy)\big) {\rm d}x{\rm d}y.
\]
Its inverse is simply
\[
\mathcal{F}^{-1}\varphi(x,y)= \frac{1}{2\pi} \int_{\mathbb{R}^2} \varphi(u,v) \exp\big(-i(ux+vy)\big) {\rm d}u{\rm d}v.
\]

\begin{remark}
If $\varphi$ is a Schwartz function on $\mathbb{R}^2$, then one can check that $O_h(\varphi) = \frac{1}{2\pi}\int_{\mathbb{R}^{2}} \mathcal{F}^{-1}\varphi(u,v) \exp\big( iuP(h_1)+iv Q(h_2)\big) {\rm d}u{\rm d}v$ defines a bounded operator. It is known that the map from $\mathcal{S}(\mathbb{R}^2)$ to $B(\eufrak{H})$ defined in this way extends to a continuous map from $L^p(\mathbb{R})$ to $B(\eufrak{H})$ for all $p\in[1,2]$, but that for $p>2$ there exist functions in $L^p(\mathbb{R}^2)$ for which we can not define a bounded operator in this way, see, e.g., \cite{wong98} and the references cited therein. But it can be extended to exponential functions, since $\frac{1}{2\pi}\mathcal{F}^{-1}\exp i(x_0u+y_0v)= \delta_{(x_0,y_0)}$ and thus
\[
O_h\big(\exp i(x_0u+y_0v)\big)=U(x_0h_1,y_0h_2).
\]
\end{remark}

\begin{lemma}\label{lem-bounded}
Let $1\le p\le 2$ and $h\in\eufrak{h}\otimes\mathbb{R}^2$ such that $\langle h_1,h_2\rangle\not=0$. Then we have $L^p(\mathbb{R}^2)\subseteq{\rm Dom}\,O_h$ and there exists a constant $C_{h,p}$ such that
\[
||O_h(\varphi)||\le C_{h,p} ||\varphi||_p
\]
for all $\varphi\in L^p(\mathbb{R}^2)$.
\end{lemma}
\begin{proof}
This follows immediately from \cite[Theorem 11.1]{wong98}, where it is stated for the irreducible unitary representation with parameter $\hbar=1$ of the Heisenberg-Weyl group.
\end{proof}

As `joint density' of the pair $\big(P(h_1),Q(h_2)\big)$ we will use its Wigner distribution.
\begin{definition}
Let $\Phi$ be a state on $B(\eufrak{H})$. We will call ${\rm d}W_{h,\Phi}$ the Wigner distribution of $\big(P(h_1),Q(h_2)\big)$ in the state $\Phi$, if
\[
\int \varphi {\rm d}W_{h,\Phi} = \Phi\big(O_h(\varphi)\big)
\]
is satisfied for all Schwartz functions $\varphi$.
\end{definition}

In general, ${\rm d}W_{h,\Phi}$ is not positive, but only a signed measure, since $O_h$ does not map positive functions to positive operators. But we can show that it has a density.
\begin{proposition}
Let $h=(h_1,h_2)\in \eufrak{h}\otimes \mathbb{R}^2$ such that $\langle h_1,h_2\rangle\not=0$ and let $\Phi$ be a state on $B(\eufrak{H})$. Then there exists a function $w_{h,\Phi}\in \bigcap_{2\le p\le\infty} L^p(\mathbb{R}^2)$ such that ${\rm d}W_{h,\Phi}= w_{h,\Phi}{\rm d}x{\rm d}y$.
\end{proposition}
\begin{proof}
It is sufficient to observe that Lemma \ref{lem-bounded} implies that the map $\varphi\mapsto \Phi\big(O_h(\varphi)\big)$ defines a continuous linear functional on $L^p(\mathbb{R}^2)$ for $1\le p\le 2$.
\end{proof}

The following proposition will play the role of the Girsanov transformation in classical Malliavin calculus. If we conjugate $O_h(\varphi)$ with $U(-k_2/2,k_1/2)$ for $k\in \eufrak{h}\otimes\mathbb{R}^2$, then this amounts to a translation of the argument of $\varphi$ by $(\langle k_1,h_1\rangle,\langle k_2,h_2\rangle)$.

\begin{proposition}\label{prop girsanov}
Let $h,k\in \eufrak{h}\otimes\mathbb{R}^2$ and $\varphi\in {\rm Dom}\,O_h$. Then we have
\[
U(-k_2/2,k_1/2)O_h(\varphi)U(-k_2/2,k_1/2)^*= O_h\big(T_{(\langle k_1,h_1\rangle,\langle k_2,h_2\rangle)}\varphi\big)
\]
where $T_{(x_0,y_0)}\varphi(x,y)=\varphi(x+x_0,y+y_0)$.
\end{proposition}
\begin{proof}
For $(u,v)\in\mathbb{R}^2$, we have
\begin{eqnarray*}
&&U(-k_2/2,k_1/2)\exp\big(i(uP(h_1)+vQ(h_2))\big)U(-k_2/2,k_1/2)^* \\
&=& U(-k_2/2,k_1/2)U(uh_1,vh_2)U(-k_2/2,k_1/2)^* \\
&=&\exp -i\big(u\langle k_1, h_1\rangle+v\langle k_2, h_2\rangle\big) U(uh_1,vh_2)
\end{eqnarray*}
and therefore
\begin{eqnarray*}
&& U(-k_2/2,k_1/2)O_h(\varphi)U(-k_2/2,k_1/2)^* \\
&=&\int_{\mathbb{R}^{2}} \mathcal{F}^{-1}\varphi(u,v)\exp\left(-i\big(u\langle k_1, h_1\rangle+v\langle k_2, h_2\rangle \big)\right) \exp i\big(uP(h_1)+vQ(h_2)\big) {\rm d}u{\rm d}v \\
&=&\int_{\mathbb{C}^{d}} \mathcal{F}^{-1}T_{(\langle k_1,h_1\rangle,\langle k_2,h_2\rangle)}  \varphi(u,v) \exp i\big(uP(h_1)+vQ(h_2)\big) {\rm d}u{\rm d}v \\
&=& O_h\left( T_{(\langle k_1,h_1\rangle,\langle k_2,h_2\rangle)}\varphi\right).
\end{eqnarray*}
\end{proof}

{}From this formula we can derive a kind of integration by parts formula that can be used to get the estimates that show the differentiability of the Wigner densities.
\begin{proposition}\label{prop-int by parts}
Let $h\in \eufrak{h}\otimes \mathbb{R}^2$, $k\in \eufrak{h}_\mathbb{C}\otimes \mathbb{C}^2$, and $\varphi$ such that $\varphi,\frac{\partial\varphi}{\partial x},\frac{\partial\varphi}{\partial y}\in {\rm Dom}\,O_h$. Then $[Q(\overline{k}_1)-P(\overline{k}_2),O_h(\varphi)]$ defines a bounded operator on $\eufrak{H}$ and we have
\[
\frac{i}{2}[Q(\overline{k}_1)-P(\overline{k}_2),O_h(\varphi)] = O_h\left(\langle k_1, h_1\rangle\frac{\partial\varphi}{\partial x}+ \langle k_2, h_2\rangle\frac{\partial\varphi}{\partial y}\right)
\]
\end{proposition}
\begin{proof}
For real $k$ this is the infinitesimal version of the previous proposition, just differentiate
\[
U(\varepsilon k_2/2,\varepsilon k_1/2)O_h(\varphi)U(\varepsilon k_2/2,\varepsilon k_1/2)^*=O_h\big(T_{(\varepsilon \langle k_1,h_1\rangle,\varepsilon\langle k_2,h_2\rangle)}\varphi\big)
\]
with respect to $\varepsilon$ and set $\varepsilon =0$. For complex $k$ it follows by linearity.
\end{proof}

Like the integration by parts formula in classical Malliavin calculus, this formula follows from a Girsanov transformation. Furthermore, it can also be used to derive sufficient conditions for the existence of smooth densities.

\begin{proposition}\label{prop-regular1}
Let $\kappa\in\mathbb{N}$, $h\in\eufrak{h}\otimes\mathbb{R}^2$ with $\langle h_1,h_2\rangle\not=0$, and $\Phi$ a vector state, i.e.\ there exists a unit vector $\omega\in\eufrak{H}$ such that $\Phi(X)=\langle\omega,X\omega\rangle$ for all $X\in B(\eufrak{H})$. If there exists a $k\in\eufrak{h}_\mathbb{C}\otimes\mathbb{C}^2$ such that
\[
\omega\in\bigcap_{\kappa_1+\kappa_2\le \kappa}{\rm Dom}\,Q(k_1)^{\kappa_1}P(k_2)^{\kappa_2}\cap \bigcap_{\kappa_1+\kappa_2\le \kappa}{\rm Dom}\,Q(\overline{k}_1)^{\kappa_1}P(\overline{k}_2)^{\kappa_2}
\]
and
\[
\langle h_1,k_1\rangle \not=0\qquad\text{ and }\qquad \langle h_2,k_2\rangle \not=0,
\]
then $w_{h,\Phi}\in \bigcap_{2\le p\le \infty}H^{p,\kappa}(\mathbb{R}^2)$, i.e.\ the Wigner density $w_{h,\Phi}$ lies in the Sobolev spaces of order $\kappa$ for all $2\le p\le \infty$. 
\end{proposition}
\begin{proof}
We will show the result for $\kappa=1$, the general case can be shown similarly (see also the proof of Theorem \ref{theo regular2}). Let $\varphi$ be a Schwartz function. Let $p\in[1,2]$. Then we have
\begin{eqnarray*}
\left|\int\frac{\partial\varphi}{\partial x}{\rm d}W_{h,\Phi}\right| &=& \left|\left\langle \omega, O_h\left(\frac{\partial\varphi}{\partial x}\right)\omega\right\rangle\right| \\
&=& \left|\left\langle \omega,\frac{i}{2|\langle k_1,h_1\rangle|} \big[Q(\overline{k}_1), O_h(\varphi)\big]\omega\right\rangle\right| \\
&\le&
\frac{C_{h,p} \big(||Q(k_1)\omega||+||Q(\overline{k}_1)\omega||\big)}{2|\langle k_1,h_1\rangle|} ||\varphi||_p.
\end{eqnarray*}
Similarly, we get
\[
\left|\int\frac{\partial\varphi}{\partial y}{\rm d}W_{h,\Phi}\right| \le
\frac{C_{h,p} \big(||P(k_2)\omega||+||P(\overline{k}_2)\omega||\big)}{2|\langle k_2,h_2\rangle|} ||\varphi||_p,
\]
and together these two inequalities imply $w_{h,\Phi}\in H^{p',1}(\mathbb{R}^2)$ for $p'=\frac{p}{p-1}$.
\end{proof}
We will give a more general result of this type in Theorem \ref{theo regular2}.

\section{The derivation operator}

In this section we define a derivation operator on our non-commutative probability space and show that it satisfies similar properties as the derivation operator on Wiener space.

We want to interpret the expression in the integration by parts formula in Proposition \ref{prop-int by parts} as a directional or Fr\'echet derivative.
\begin{definition}\label{def derivation1}
Let $k\in \eufrak{h}_\mathbb{C}\otimes\mathbb{C}^2$. We set
\begin{eqnarray*}
{\rm Dom}\,D_k &=&\left\{B\in \mathcal{B}(\eufrak{H})\Big|\,\frac{i}{2}[Q(k_1)-P(k_2),B] \text{ defines a bounded operator on }\eufrak{H}\right\}
\end{eqnarray*}
and for $B\in {\rm Dom}\,D_k$, we set $D_k B=\frac{i}{2}[Q(k_1)-P(k_2),B]$.
\end{definition}

Note that $B\in{\rm Dom}\,D_k$ for some $k\in\eufrak{h}_\mathbb{C}\otimes\mathbb{C}^2$ implies $B^*\in{\rm Dom}\,D_{\overline{k}}$ and
\[
D_{\overline{k}}B^* = (D_kB)^*.
\]

\begin{example}
Let $k\in \eufrak{h}_\mathbb{C}\otimes\mathbb{C}^2$ and let $\psi\in{\rm Dom}\,P(k_2)\cap{\rm Dom}\,Q(k_1)\cap{\rm Dom}\,P(\overline{k}_2)\cap{\rm Dom}\,Q(\overline{k}_1)$ be a unit vector. We denote by $\mathbb{P}_\psi$ the orthogonal projection onto the one-dimensional subspace spanned by $\psi$. Evaluating the commutator $[Q(k_1)-P(k_2),\mathbb{P}_\psi]$ on a vector $\phi\in{\rm Dom}\,P(k_2)\cap{\rm Dom}\,Q(k_1)$, we get
\begin{eqnarray*}
[Q(k_1)-P(k_2),\mathbb{P}_\psi]\phi &=& \langle\psi,\phi\rangle \big(Q(k_1)-P(k_2)\big)(\psi) - \langle\psi,\big(Q(k_1)-P(k_2)\big)(\phi)\rangle \psi \\
&=& \langle\psi,\phi\rangle \big(Q(k_1)-P(k_2)\big)(\psi) - \langle\big(Q(\overline{k}_1)-P(\overline{k}_2)\big)\psi,\phi\rangle \psi
\end{eqnarray*}
We see that the range of $[Q(k_1)-P(k_2),\mathbb{P}_\psi]$ is two-dimensional, so it can be extended to a bounded operator on $\eufrak{H}$. Therefore $\mathbb{P}_\psi\in {\rm Dom}\, D_k$, and we get
\[
(D_k\mathbb{P}_\psi)\phi =\frac{i}{2}\Big(\langle\psi,\phi\rangle \big(Q(k_1)-P(k_2)\big)(\psi) - \langle\big(Q(\overline{k}_1)-P(\overline{k}_2)\big)\psi,\phi\rangle \psi\Big)
\]
for all $\phi\in\eufrak{H}$.
\end{example}

\begin{example}
Let $h\in \eufrak{h}\otimes\mathbb{R}^2$, $k\in \eufrak{h}_\mathbb{C}\otimes\mathbb{C}^2$. Then $\frac{i}{2}[Q(k_1)-P(k_2),U(h_1,h_2)]$ defines a bounded operator on $\eufrak{H}$, and we get
\[
D_k U(h_1,h_2) = i \big(\langle \overline{k}_1, h_1\rangle+\langle \overline{k}_2,h_2 \rangle\big) U(h_1,h_2).
\]
\end{example}

\begin{proposition}\label{prop closable1}
Let $k\in \eufrak{h}_\mathbb{C}\otimes\mathbb{C}^2$. The operator $D_k$ is a closable operator from $\mathcal{B}(\eufrak{H})$ to $\mathcal{B}(\eufrak{H})$ with respect to the weak topology.
\end{proposition}
\begin{proof}
Let $(B_n)_{n\in\mathbb{N}}\subseteq{\rm Dom}\,D_k\subseteq\mathcal{B}(\eufrak{H})$ be any sequence such that $B_n\to 0$ and $D_kB_n\to \beta$ for some $\beta\in\mathcal{B}(\eufrak{H})$ in the weak topology. To show that $D_k$ is closable, we have to show that this implies $\beta=0$. Let us evaluate $\beta$ between two exponential vectors $\mathcal{E}(h_1)$, $\mathcal{E}(h_2)$, $h_1,h_2\in\eufrak{h}_\mathbb{C}$, then we get
\begin{eqnarray*}
\langle \mathcal{E}(h_1),\beta \mathcal{E}(h_2)\rangle &=& \lim_{n\to \infty} \langle \mathcal{E}(h_1), D_kB_n\mathcal{E}(h_2)\rangle \\
&=& \lim_{n\to \infty}  \frac{i}{2} \big\langle \big(Q(\overline{k}_1)-P(\overline{k}_2)\big)\mathcal{E}(h_1), B_n\mathcal{E}(h_2)\big\rangle \\
&&-\lim_{n\to\infty}  \frac{i}{2}\big\langle\mathcal{E}(h_1), B_n\big(Q(k_1)-P(k_2)\big)\mathcal{E}(h_2)\big\rangle \\
&=& 0,
\end{eqnarray*}
and therefore $\beta=0$, as desired. 
\end{proof}

\begin{definition}\label{def derivation2}
We set
\[
\mathcal{S}= {\rm alg}\,\left\{O_h(\varphi)\Big| h\in \eufrak{h}\otimes\mathbb{R}; \varphi\in C^\infty(\mathbb{R}^2) \text{ s.t. }\frac{\partial^{\kappa_1+\kappa_2} \varphi}{\partial x^{\kappa_1}\partial y^{\kappa_2}}\in {\rm Dom}\, O_h \text{ for all } \kappa_1,\kappa_2\ge 0\right\},
\]
the elements of $\mathcal{S}$ will play the role of the smooth functionals. Note that $\mathcal{S}$ is weakly dense in $\mathcal{B}(\eufrak{H})$, i.e.\ $\mathcal{S}''=B(\eufrak{H})$, since $\mathcal{S}$ contains the Weyl operators $U(h_1,h_2)$ with $h_1,h_2\in\eufrak{h}$.

We define $D:\mathcal{S}\to\mathcal{B}(\eufrak{H})\otimes \eufrak{h}_\mathbb{C}\otimes \mathbb{C}^2$ (where the tensor product is the algebraic tensor product over $\mathbb{C}$) by setting $DO_h(\varphi)$ equal to
\[
DO_h(\varphi)= \left(\begin{array}{c} O_h\left(\displaystyle\frac{\partial \varphi}{\partial x}\right)\otimes h_1 \\[4mm] O_h\left(\displaystyle\frac{\partial \varphi}{\partial y}\right)\otimes h_2 \end{array}\right)
\]
and extending it as a derivation w.r.t.\ the $\mathcal{B}(\eufrak{H})$-bimodule structure of $\mathcal{B}(\eufrak{H})\otimes\eufrak{h}_\mathbb{C}\otimes \mathbb{C}^2$ defined by
\[
O\cdot \left(\begin{array}{c} O_1\otimes k_1 \\ O_2\otimes k_2 \end{array}\right) = \left(\begin{array}{c} OO_1\otimes k_1 \\ OO_2\otimes k_2 \end{array}\right), \quad
\left(\begin{array}{c} O_1\otimes k_1 \\ O_2\otimes k_2 \end{array}\right)\cdot O = \left(\begin{array}{c} O_1O\otimes k_1 \\ O_2O\otimes k_2 \end{array}\right) \
\]
for $O, O_1,O_2\in\mathcal{B}(\eufrak{H})$ and $k\in\eufrak{h}_\mathbb{C}\otimes\mathbb{C}^2$.
\end{definition}

\begin{example}
For $h\in \eufrak{h}\otimes\mathbb{R}^2$, we get
\begin{eqnarray*}
D U(h_1,h_2)&=& D O_h\Big(\exp i(x+y)\Big)
= i\left(\begin{array}{c}U(h_1,h_2)\otimes h_1 \\ U(h_1,h_2)\otimes h_2\end{array}\right) \\
&=& i U(h_1,h_2)\otimes h.
\end{eqnarray*}
\end{example}

\begin{definition}
We can define a $B(\eufrak{H})$-valued inner product on $\mathcal{B}(\eufrak{H})\otimes\eufrak{h}_\mathbb{C}\otimes\mathbb{C}^2$ by $\langle\cdot,\cdot\rangle:\mathcal{B}(\eufrak{H})\otimes\eufrak{h}_\mathbb{C}\otimes\mathbb{C}^2\times \mathcal{B}(\eufrak{H})\otimes\eufrak{h}_\mathbb{C}\otimes\mathbb{C}^2 \to \mathcal{B}(\eufrak{H})$ by
\[
\left\langle \left(\begin{array}{c} O_1\otimes h_1 \\ O_2\otimes h_2\end{array}\right), \left(\begin{array}{c} O'_1\otimes k_1 \\ O'_2\otimes k_2\end{array}\right)\right\rangle =  O^*_1O'_1 \langle h_1, k_1\rangle + O^*_2O'_2 \langle h_2, k_2\rangle
\]
\end{definition}
We have
\begin{eqnarray*}
\langle B,A\rangle &=& \langle A,B\rangle^* \\
O^*\langle A,B\rangle &=& \langle AO,B\rangle \\
\langle A,B\rangle O &=& \langle A,BO\rangle \\
\langle O^*A,B\rangle &=& \langle A,OB\rangle
\end{eqnarray*}
for all $A,B\in \mathcal{B}(\eufrak{H})\otimes\eufrak{h}_\mathbb{C}\otimes\mathbb{C}^2$ and all $O\in \mathcal{B}(\eufrak{H})$. This turns $\mathcal{B}(\eufrak{H})\otimes\eufrak{h}_\mathbb{C}\otimes\mathbb{C}^2$ into a pre-Hilbert module over $\mathcal{B}(\eufrak{H})$. It can be embedded in the Hilbert module $\eufrak{M}=\mathcal{B}(\eufrak{H},\eufrak{H}\otimes\eufrak{h}_\mathbb{C}\otimes\mathbb{C}^2)$ by mapping $O\otimes k\in\mathcal{B}(\eufrak{H})\otimes\eufrak{h}_\mathbb{C}\otimes\mathbb{C}^2 $ to the linear map $\eufrak{H}\ni v \mapsto Ov\otimes k\in \eufrak{H}\otimes\eufrak{h}_\mathbb{C}\otimes\mathbb{C}^2$. We will regard $\eufrak{h}_\mathbb{C}\otimes \mathbb{C}^2$ as a subspace of $\eufrak{M}$ via the embedding $\eufrak{h}_\mathbb{C}\ni k\mapsto {\rm id}_\eufrak{H}\otimes k\in \eufrak{M}$. Note that we have $O\cdot k=k\cdot O = O\otimes k$ and $\langle A,k\rangle=\langle \overline{k}, \overline{A}\rangle$ for all $k\in\eufrak{h}_\mathbb{C}\otimes\mathbb{C}^2$, $O\in\mathcal{B}(\eufrak{H})$, $A\in \eufrak{M}$, where the conjugation in $\eufrak{M}$ is defined by $\overline{O\otimes k}=O^*\otimes\overline{k}$.

\begin{proposition}\label{prop frechet}
Let $O\in \mathcal{S}$ and $k\in \eufrak{h}_\mathbb{C}\otimes\mathbb{C}^2$. Then $O\in {\rm Dom}\,D_k$ and
\[
D_k O = \langle \overline{k}, DO\rangle = \langle \overline{DO},k \rangle .
\]
\end{proposition}
\begin{proof}
For $h\in \eufrak{h}\otimes\mathbb{R}^2$ and $\varphi\in {\rm Dom}\,O_h$ s.t.\ also $\frac{\partial\varphi}{\partial x},\frac{\partial\varphi}{\partial y}\in {\rm Dom}\,O_h$, we get
\begin{eqnarray*}
\langle \overline{k}, DO_h(\varphi)\rangle &=& \left\langle \left(\begin{array}{c} \overline{k}_1 \\ \overline{k}_2\end{array}\right), \left(\begin{array}{c} O_h\left(\displaystyle\frac{\partial \varphi}{\partial x}\right)\otimes h_1 \\[4mm] O_h\left(\displaystyle\frac{\partial \varphi}{\partial y}\right)\otimes h_2\end{array}\right)\right\rangle \\
&=& O_h\left(\langle \overline{k}_1, h_1\rangle\frac{\partial\varphi}{\partial x}+ \langle \overline{k}_2, h_2\rangle\frac{\partial\varphi}{\partial y}\right)
\\
&=& \frac{i}{2}[Q(k_1)-P(k_2),O_h(\varphi)] = D_k O,
\end{eqnarray*}
where we used Proposition \ref{prop-int by parts}. The first equality of the proposition now follows, since both $O\mapsto D_kO=\frac{i}{2}[Q(k_1)-P(k_2),O]$ and $O\mapsto \langle \overline{k}, DO\rangle$ are derivations.

The second equality follows immediately.
\end{proof}

The next result is the analogue of Equation \eqref{int by parts1}.
\begin{theorem}\label{theo int by parts}
We have
\[
\mathbb{E}\big( \langle \overline{k}, DO\rangle \big) = \frac{1}{2}\mathbb{E}\big(\{P(k_1)+Q(k_2),O\}\big)
\]
for all $k\in \eufrak{h}_\mathbb{C}\otimes\mathbb{C}^2$ and all $O\in \mathcal{S}$, where $\{\cdot,\cdot\}$ denotes the anti-commutator $\{X,Y\}=XY+YX$.
\end{theorem}
\begin{proof}
This formula is a consequence of the fact that $Q(h)\Omega=h=iP(h)\Omega$ for all $h\in\eufrak{h}_\mathbb{C}$, we get
\begin{eqnarray*}
\mathbb{E}\big( \langle \overline{k}, DO\rangle \big) &=& \frac{i}{2} \Big( \langle \big(Q(\overline{k}_1)-P(\overline{k}_2)\big)\Omega,O\Omega\rangle - \langle\Omega, O\big(Q(k_1)-P(k_2)\big)\Omega\rangle\Big) \\
&=&  \frac{i}{2}\big(\langle \overline{k}_1+i\overline{k}_2 , O\Omega\rangle - \langle \Omega, O (k_1+ik_2)\rangle\big) \\
 &=& \frac{1}{2}\Big( \langle \big(P(\overline{k}_1)+Q(\overline{k}_2)\big)\Omega,O\Omega\rangle + \langle\Omega, O\big(P(k_1)+Q(k_2)\big)\Omega\rangle\Big) \\
&=&\frac{1}{2}\mathbb{E}\big(\{P(k_1)+Q(k_2),O\}\big).
\end{eqnarray*}
\end{proof}

There is also an analogue of \eqref{int by parts2}.
\begin{corollary}\label{cor int by parts}
Let $k\in \eufrak{h}_\mathbb{C}\otimes\mathbb{C}^2$, and $O_1,\ldots, O_n\in\mathcal{S}$, then
\[
\frac{1}{2}\mathbb{E}\left(\left\{P(k_1)+Q(k_2),\prod_{m=1}^n O_m\right\}\right) = \mathbb{E}\left( \sum_{m=1}^n \prod_{j=1}^{m-1} O_j\langle \overline{k}, DO_m\rangle \prod_{j=m+1}^{n} O_j\right),
\]
where the products are ordered such that the indices increase from the left to the right.
\end{corollary}
\begin{proof}
This is obvious, since $O\mapsto \langle \overline{k}, DO\rangle$ is a derivation.
\end{proof}

This formula for $n=3$ can be used to show that $D$ is a closable operator from $\mathcal{B}(\eufrak{H})$ to $\eufrak{M}$.
\begin{corollary}\label{cor closable2}
The derivation operator $D$ is a closable operator from $\mathcal{B}(\eufrak{H})$ to the $\mathcal{B}(\eufrak{H})$-Hilbert module $\eufrak{M}=\mathcal{B}(\eufrak{H},\eufrak{H}\otimes \eufrak{h}_\mathbb{C}\otimes\mathbb{C}^2)$ w.r.t.\ the weak topologies.
\end{corollary}
\begin{proof}
We have to show that for any sequence $(A_n)_{n\in\mathbb{N}}$ in $\mathcal{S}$ with $A_n\to 0$ and $DA_n\to \alpha \in \eufrak{M}$, we get $\alpha=0$. Let $f,g\in\eufrak{h}_\mathbb{C}$. Set $f_1=\frac{f+\overline{f}}{2}$, $f_2=\frac{f-\overline{f}}{2i}$, $g_1=\frac{g+\overline{g}}{2}$, and $g_2=\frac{g-\overline{g}}{2i}$, then we have $U(f_1,f_2)\Omega=e^{-||f||/2}\mathcal{E}(f)$ and $U(g_1,g_2)\Omega=e^{-||g||/2}\mathcal{E}(g)$. Thus we get
\begin{eqnarray*}
&&\exp\big((||f||^2+||g||^2)/2\big)\langle \mathcal{E}(f)\otimes \overline{h},\alpha \mathcal{E}(g)\rangle \\
&=& \exp\big((||f||^2+||g||^2)/2\big)\langle \mathcal{E}(f),\langle \overline{h},\alpha\rangle \mathcal{E}(g)\rangle \\
&=& \lim_{n\to\infty}\mathbb{E} \big(U(-f_1,-f_2) \langle \overline{h},D A_n\rangle U(g_1,g_2)\big) \\
&=& \lim_{n\to\infty}\mathbb{E}\Big( \frac{1}{2}\big\{P(h_1)+Q(h_2),U(-f_1,-f_2)A_n U(g_1,g_2)\big\} \\
&& - \big\langle \overline{h}, DU(-f_1,-f_2)\big\rangle A_n U(g_1,g_2)
- U(-f_1,-f_2)A_n\big\langle \overline{h}, DU(g_1,g_2)\big\rangle\Big) \\
&=& \lim_{n\to\infty}\big( \langle \psi_1, A_n \psi_2\rangle +\langle \psi_3, A_n \psi_4\rangle -\langle \psi_5, A_n \psi_6\rangle -\langle \psi_7, A_n \psi_8\rangle\big)  \\
&=& 0
\end{eqnarray*}
for all $h\in\eufrak{h}_\mathbb{C}\otimes\mathbb{C}^2$, where
\[
\begin{array}{rclcrcl}
\psi_1&=&\frac{1}{2}U(f_1,f_2)\big(P(\overline{h}_1)+Q(\overline{h}_2)\big)\Omega, && \psi_2&=& U(g_1,g_2)\Omega, \\
\psi_3&=&U(f_1,f_2)\Omega, && \psi_4&=& \frac{1}{2}U(g_1,g_2)\big(P(h_1)+Q(h_2)\big)\Omega, \\
\psi_5&=&\big(D_hU(-f_1,-f_2)\big)^*\Omega, && \psi_6&=&  U(g_1,g_2)\Omega, \\
\psi_7&=&U(f_1,f_2)\Omega, && \psi_8&=&D_hU(g_1,g_2)\Omega.
\end{array}
\]
But this implies $\alpha=0$, since $\{\mathcal{E}(f)\otimes \overline{h}|f\in \eufrak{h}_\mathbb{C},h\in\eufrak{h}_\mathbb{C}\otimes\mathbb{C}^2\}$ is dense in $\eufrak{H}\otimes\eufrak{h}_\mathbb{C}\otimes\mathbb{C}^2$.
\end{proof}

\begin{remark}
This implies that $D$ is also closable in stronger topologies, such as, e.g., the norm topology and the strong topology.
\end{remark}

We will denote the closure of $D$ again by the same symbol.

\begin{proposition}
Let $O\in{\rm Dom}\,D$. Then $O^*\in{\rm Dom}\,D$ and
\[
DO^* = \overline{DO}.
\]
In particular, since $D$ is a derivation, this implies that ${\rm Dom}\,D$ is a $*$-subalgebra of $\mathcal{B}(\eufrak{H})$.
\end{proposition}
\begin{proof}
It is not difficult to check this directly on the Weyl operators $U(h_1,h_2)$, $h\in\eufrak{h}\otimes\mathbb{R}^2$. We get $U(h_1,h_2)^*=U(-h_1,-h_2)$ and
\begin{eqnarray*}
D\big(U(h_1,h_2)^*\big) &=& DU(-h_1,-h_2)= -i U(-h_1,-h_2)\otimes h \\
&=& U(h_1,h_2)^* \otimes \overline{(ih)} = \overline{DU(h_1,h_2)}.
\end{eqnarray*}
By linearity and continuity it therefore extends to all of ${\rm Dom}\,D$.
\end{proof}

We will now show how $D$ can be iterated. Let $H$ be a complex Hilbert space, then we can define a derivation operator $D:\mathcal{S}\otimes H\to \mathcal{B}(\eufrak{H})\otimes\eufrak{h}_\mathbb{C}\otimes\mathbb{C}^2\otimes H$ by setting $D(O\otimes h)=DO\otimes h$ for $O\in\mathcal{S}$ and $h\in H$. Closing it, we get an unbounded derivation from the Hilbert module $\mathcal{B}(\eufrak{H}, \eufrak{H}\otimes H)$ to $\eufrak{M}(H)=\mathcal{B}(\eufrak{H}\otimes H,\eufrak{H}\otimes\eufrak{h}_\mathbb{C}\otimes \mathbb{C}^2\otimes H)$. This allows us to iterate $D$. It is easy to see that $D$ maps $\mathcal{S}\otimes H$ to $\mathcal{S}\otimes\eufrak{h}_\mathbb{C}\otimes \mathbb{C}^2\otimes H$ and so we have $D^n(\mathcal{S}\otimes H)\subseteq \mathcal{S}\otimes \left(\eufrak{h}_\mathbb{C}\otimes \mathbb{C}^2\right)^{\otimes n}\otimes H$. In particular, $\mathcal{S}\subseteq{\rm Dom}\,D^n$ for all $n\in\mathbb{N}$, and we can define Sobolev-type norms $||\cdot||_n$ and semi-norms $||\cdot||_{\psi,n}$, on $\mathcal{S}$ by
\begin{eqnarray*}
||O||^2_n &=& ||O^*O||+\sum_{j=1}^n ||\langle D^n O,D^nO\rangle||, \\
||O||^2_{\psi,n} &=& ||O\psi||^2+\sum_{j=1}^n ||\langle \psi, \langle D^n O,D^nO\rangle\psi\rangle||, \quad\text{ for } \psi\in\eufrak{H}
\end{eqnarray*}
In this way we can define Sobolev-type topologies on ${\rm Dom}\,D^n$.

We will now extend the definition of the ``Fr\'echet derivation'' $D_k$ to the case where $k$ is replaced by an element of $\eufrak{M}$. It becomes now important to distinguish between a right and a left ``derivation operator''. Furthermore, it is no longer a derivation.
\begin{definition}
Let $u\in\eufrak{M}$ and $O\in {\rm Dom}\,D$. Then we define the right gradient $\overrightarrow{D}_uO$ and the left gradient $O\overleftarrow{D}_u$ of $O$ with respect to $u$ by
\begin{eqnarray*}
\overrightarrow{D}_u O&=& \langle \overline{u},DO\rangle, \\
O\overleftarrow{D}_u &=& \langle \overline{DO},u\rangle. \\
\end{eqnarray*}
\end{definition}

We list several properties of the gradient.
\begin{proposition}
\begin{enumerate}
\item
Let $X\in\mathcal{B}(\eufrak{H})$, $O,O_1,O_2\in {\rm Dom}\,D$, and $u\in\eufrak{M}$. Then
\begin{eqnarray*}
\overrightarrow{D}_{Xu}O &=& X \overrightarrow{D}_uO, \\
\overrightarrow{D}_u(O_1O_2) &=& \left(\overrightarrow{D}_uO_1\right) O_2+\overrightarrow{D}_{uO_1}O_2 , \\
O\overleftarrow{D}_{uX} &=& (O\overleftarrow{D}_u)X, \\
(O_1O_2)\overleftarrow{D}_u &=& O_1\overleftarrow{D}_{O_2u} + O_1\left(O_2\overleftarrow{D}_u\right) , \\
\end{eqnarray*}
\item
For $k\in \eufrak{h}_\mathbb{C}\otimes \mathbb{C}^2$ and $O\in{\rm Dom}\,D$, we have
\[
D_k O = \overrightarrow{D}_{{\rm id}_\eufrak{H}\otimes k}O = O\overleftarrow{D}_{{\rm id}_\eufrak{H}\otimes k}
\]
\end{enumerate}
\end{proposition}
\begin{proof}
These properties can be deduced easily from the definition of the gradient and the properties of the derivation operator $D$ and the inner product $\langle,\rangle$.
\end{proof}

\begin{remark}
We can also define a two-sided gradient $\overleftrightarrow{D}_u:{\rm Dom}\,D\times{\rm Dom}\,D\to \mathcal{B}(\eufrak{H})$ by $\overleftrightarrow{D}_u:(O_1,O_2)\mapsto O_1\overleftrightarrow{D}_uO_2= O_1\left(\overrightarrow{D}_uO_2\right) + \left(O_1\overleftarrow{D}_u\right)O_2$. For $k\in \eufrak{h}_\mathbb{C}\otimes \mathbb{C}^2$ we have $O_1\overleftrightarrow{D}_{{\rm id}_\eufrak{H}\otimes k}  O_2= D_k(O_1O_2)$.
\end{remark}

\section{The divergence operator}

The algebra $\mathcal{B}(\eufrak{H})$ of bounded operators on the symmetric Fock space $\eufrak{H}$ and the Hilbert module $\eufrak{M}$ are not Hilbert spaces with respect to the expectation in the vacuum vector $\Omega$. Therefore we can not define the divergence operator or Skorohod integral $\delta$ as the adjoint of the derivation $D$. It might be tempting to try to define $\delta X$ as an operator such that the condition
\begin{equation}\label{duality}
\mathbb{E}\big((\delta X) B\big)\overset{?}{=} \mathbb{E}\left(\overrightarrow{D}_X B\right)
\end{equation}
is satisfied for all $B\in{\rm Dom}\,\overrightarrow{D}_X$, even though it is not sufficient to characterize $\delta X$. But the following proposition shows that this is not possible.

\begin{proposition}\label{prop no go}
Let $k\in \eufrak{h}_\mathbb{C}\otimes\mathbb{C}^2$ with $k_1+ik_2\not=0$. There exists no (possibly unbounded) operator $M$ whose domain contains the vacuum vector such that
\[
\mathbb{E}\big( M B\big) = \mathbb{E}\left(D_k B\right)
\]
holds for all $B\in{\rm Dom}\,D_k$.
\end{proposition}
\begin{proof}
We assume that such an operator $M$ exists and show that this leads to a contradiction.

Let $B\in\mathcal{B}(\eufrak{H})$ be the operator defined by $\eufrak{H}\ni\psi\mapsto \langle k_1+ik_2,\psi\rangle\Omega$, it is easy to see that $B\in{\rm Dom}\,D_k$ and that $D_kB$ is given by
\[
(D_kB)\psi =  \frac{i}{2}\langle k_1+ik_2,\psi\rangle(k_1+ik_2)-\frac{i}{2}\big\langle\big(Q(\overline{k}_1)-P(\overline{k}_2)\big)(k_1+ik_2),\psi\big\rangle \Omega,\text{ for }\psi\in\eufrak{H}.
\]
Therefore, if $M$ existed, we would have
\begin{eqnarray*}
0 &=& \langle \Omega, MB\Omega\rangle = \mathbb{E}(MB) = \mathbb{E}(D_kB) \\
&=& \langle \Omega,(D_kB)\Omega\rangle =  -\frac{i}{2}\langle k_1+ik_2,k_1+ik_2\rangle,
\end{eqnarray*}
which is clearly impossible.
\end{proof}

We introduce the analogue of smooth elementary $\eufrak{h}$-valued random variables,
\[
\mathcal{S}_\eufrak{h}=\left\{ u=\sum_{j=1}^n F_j\otimes h^{(j)}\Big| n\in\mathbb{N}, F_1,\ldots,F_n\in\mathcal{S},h^{(1)},\ldots,h^{(n)}\in\eufrak{h}_\mathbb{C}\otimes\mathbb{C}^2\right\}.
\]
If we define $A\overleftrightarrow{\delta u}B$ for $u=\sum_{j=1}^n F_j\otimes h^{(j)}\in\mathcal{S}_\eufrak{h}$ and $A,B\in B(\eufrak{H})$ by
\[
A\overleftrightarrow{\delta u}B= \frac{1}{2}\sum_{j=1}^n \left\{ P\left(h^{(j)}_1\right)+Q\left(h^{(j)}_2\right) ,AF_jB\right\} - A\sum_{j=1}^n \left(D_{h^{(j)}}F_j\right)B.
\]
then it follows from Corollary \ref{cor int by parts} that this satisfies
\begin{equation}\label{duality2}
\mathbb{E}\left(A\overleftrightarrow{\delta u}B\right) = \mathbb{E}\left(A\overleftrightarrow{D}_u B\right).
\end{equation}
But this can only be written as a product $AXB$ with some operator $X$, if $A$ and $B$ commute with $P\left(h^{(1)}_1\right)+Q\left(h^{(1)}_2\right),\ldots, P\left(h^{(n)}_1\right)+Q\left(h^{(n)}_2\right)$. We see that a condition of the form (\ref{duality}) or (\ref{duality2}) is too strong, if we require it to be satisfied for all $A,B\in{\rm Dom}\,D$. We have to impose some commutativity on $A$ and $B$ to weaken the condition, in order to be able to satisfy it. We will now give a first definition of a divergence operator that satisfies a weaker version of (\ref{duality2}), see Proposition \ref{prop duality} below. In Remark \ref{ref def delta} we will extend this definition to a bigger domain.

\begin{definition}\label{def skorohod1}
We set 
\begin{eqnarray*}
\mathcal{S}_{\eufrak{h},\delta} &=& \Big\{ u=\sum_{j=1}^n F_j\otimes h^{(j)}\in\mathcal{S}_\eufrak{h}\Big| \frac{1}{2}\sum_{j=1}^n \left\{ P\left(h^{(j)}_1\right)+Q\left(h^{(j)}_2\right) ,F_j\right\} - \sum_{j=1}^n D_{h^{(j)}}F_j\\
&& \qquad\qquad\qquad\qquad\qquad \text{ defines a bounded operator on }\eufrak{H}\Big\}
\end{eqnarray*}
and define the divergence operator $\delta:\mathcal{S}_{\eufrak{h},\delta}\to B(\eufrak{H})$ by 
\[
\delta(u) = \frac{1}{2}\sum_{j=1}^n \left\{ P\left(h^{(j)}_1\right)+Q\left(h^{(j)}_2\right) ,F_j\right\} - \sum_{j=1}^n D_{h^{(j)}}F_j.
\]
for $u=\sum_{j=1}^n F_j\otimes h^{(j)}\in\mathcal{S}_{\eufrak{h},\delta}$.
\end{definition}

It is easy to check that
\[
\delta(\overline{u}) = \big(\delta(u)\big)^*
\]
holds for all $u\in\mathcal{S}_{\eufrak{h},\delta}$.

\begin{proposition}\label{prop duality}
Let $u=\sum_{j=1}^n F_j\otimes h^{(j)}\in\mathcal{S}_{\eufrak{h},\delta}$ and
\[
A,B\in {\rm Dom}\,D\cap\left\{ P\left(h^{(1)}_1\right)+Q\left(h^{(1)}_2\right),\ldots, P\left(h^{(n)}_1\right)+Q\left(h^{(n)}_2\right)\right\}'
\]
i.e., $A$ and $B$ are in the commutant of $\left\{ P\left(h^{(1)}_1\right)+Q\left(h^{(1)}_2\right),\ldots, P\left(h^{(n)}_1\right)+Q\left(h^{(n)}_2\right)\right\}$, then we have
\[
\mathbb{E}\big(A\delta(u)B\big) = \mathbb{E}\left(A\overleftrightarrow{D}_u B\right).
\]
\end{proposition}
\begin{remark}
Note that $\delta:\mathcal{S}_\eufrak{h,\delta}\to\mathcal{B}(\eufrak{H})$ is the only linear map with this property, since for one single element $h\in\eufrak{h}_\mathbb{C}\otimes\mathbb{C}^2$, the sets
\begin{eqnarray*}
&\left\{ A^*\Omega \Big| A\in {\rm Dom}\,D\cap\left\{ P\left(h_1\right)+Q\left(h_2\right)\right\}'\right\}& \\
&\left\{ B\Omega \Big| B\in {\rm Dom}\,D\cap\left\{ P\left(h_1\right)+Q\left(h_2\right)\right\}'\right\}&
\end{eqnarray*}
are still total in $\eufrak{H}$.
\end{remark}
\begin{proof}
{}From Corollary \ref{cor int by parts} we get
\begin{eqnarray*}
\mathbb{E}\left(A\overleftrightarrow{D}_u B\right) &=&\mathbb{E}\big(A\langle \overline{u}, DB\rangle + \langle \overline{DA}, u\rangle B \big) \\
&=&\mathbb{E}\left( \sum_{j=1}^n AF_j \left( D_{h^{(j)}}B\right) +  \sum_{j=1}^n \left( D_{h^{(j)}}A\right)F_jB\right) \\
&=& \mathbb{E} \left(\frac{1}{2}\sum_{j=1}^n \left\{ P\left(h^{(j)}_1\right)+Q\left(h^{(j)}_2\right) ,AF_jB\right\} - \sum_{j=1}^n A \left(D_{h^{(j)}}F_j\right)B\right).
\end{eqnarray*}
But since $A$ and $B$ commute with $P\left(h^{(1)}_1\right)+Q\left(h^{(1)}_2\right),\ldots, P\left(h^{(n)}_1\right)+Q\left(h^{(n)}_2\right)$, we can pull them out of the anti-commutator, and we get
\begin{eqnarray*}
\mathbb{E}\left(A\overleftrightarrow{D}_u B\right) &=&\mathbb{E} \left(\frac{1}{2}\sum_{j=1}^n A \left\{ P\left(h^{(j)}_1\right)+Q\left(h^{(j)}_2\right) ,F_j\right\}B - \sum_{j=1}^n A\left(D_{h^{(j)}}F_j\right)B\right) \\
&=& \mathbb{E}(A\delta(u)B).
\end{eqnarray*}
\end{proof}

We will now give an explicit formula for the matrix elements between two exponential vectors of the divergence of a smooth elementary element $u\in \mathcal{S}_{\eufrak{h},\delta}$, this is the analogue of the first fundamental lemma in the Hudson-Parthasarathy calculus, see, e.g., \cite[Proposition 25.1]{parthasarathy92}.

\begin{theorem}\label{theo exp values}
Let $u\in\mathcal{S}_{\eufrak{h},\delta}$. Then we have the following formula
\[
\langle\mathcal{E}(k_1), \delta(u)\mathcal{E}(k_2)\rangle = \left\langle\mathcal{E}(k_1)\otimes\left(\begin{array}{c} ik_1-i\overline{k}_2 \\ k_1+\overline{k}_2\end{array}\right)   , u \mathcal{E}(k_2)\right\rangle
\]
for the evaluation of the divergence $\delta(u)$ of $u$ between two exponential vectors $\mathcal{E}(k_1)$, $\mathcal{E}(k_2)$, for $k_1,k_2\in\eufrak{h}_\mathbb{C}$.
\end{theorem}
\begin{remark}\label{ref def delta}
This suggests to extend the definition of $\delta$ in the following way: set
\begin{eqnarray}
{\rm Dom}\,\delta &=& \Big\{ u\in \eufrak{M}\Big| \exists M\in\eufrak{B}(\eufrak{H}), \forall k_1,k_2\in\eufrak{h}_\mathbb{C}: \nonumber\\
&& \qquad\qquad \langle\mathcal{E}(k_1), M \mathcal{E}(k_2)\rangle = \left\langle\mathcal{E}(k_1)\otimes\left(\begin{array}{c} ik_1-i\overline{k}_2 \\ k_1+\overline{k}_2\end{array}\right)   , u \mathcal{E}(k_2)\right\rangle\Big\} \label{eq-def-delta}
\end{eqnarray}
and define $\delta(u)$ for $u\in{\rm Dom}\,\delta$ to be the unique operator $M$ that satisfies the condition in Equation (\ref{eq-def-delta}).
\end{remark}
\begin{proof}
Let $u=\sum_{j=1}^n F_j\otimes h^{(j)}$. Recalling the definition of $D_h$ we get the following alternative expression for $\delta(u)$,
\begin{eqnarray}
\delta(u) &=&\frac{1}{2} \sum_{j=1}^n \Big(P\left(h^{(j)}_1\right)+Q\left(h^{(j)}_2\right)-iQ\left(h^{(j)}_1\right)+iP\left(h^{(j)}_2\right)\Big) F_j \nonumber \\
&& +\frac{1}{2}\sum_{j=1}^n F_j \Big(P\left(h^{(j)}_1\right) +Q\left(h^{(j)}_2\right)+iQ\left(h^{(j)}_1\right)-iP\left(h^{(j)}_2\right)\Big) \nonumber \\
&=&\sum_{j=1}^n \Big( a^+\left(h^{(j)}_2-ih^{(j)}_1\right)F_j + F_j a\left(\overline{h^{(j)}_2-ih^{(j)}_1}\right)\Big). \label{wick}
\end{eqnarray}
Evaluating this between two exponential vectors, we obtain
\begin{eqnarray*}
\langle\mathcal{E}(k_1), \delta(u)\mathcal{E}(k_2)\rangle &= & \sum_{j=1}^n\big\langle a(h_2^{(j)}-ih_1^{(j)})\mathcal{E}(k_1), F_j\mathcal{E}(k_2)\big\rangle \\
&& + \sum_{j=1}^n \big\langle \mathcal{E}(k_1), F_ja(\overline{h_2^{(j)}+ih_1^{(j)}})\mathcal{E}(k_2)\big\rangle \\
&=&  \sum_{j=1}^n \big(\overline{\langle h_2^{(j)}-ih_1^{(j)},k_1\rangle} + \langle \overline{h_2^{(j)}+ih_1^{(j)}},k_2\rangle\big)\langle \mathcal{E}(k_1), F_j\mathcal{E}(k_2)\rangle \\
&=& \sum_{j=1}^n \big(\langle k_1,h_2^{(j)}-ih_1^{(j)}\rangle + \langle \overline{k}_2, h_2^{(j)}+ih_1^{(j)}\rangle\big)\langle \mathcal{E}(k_1), F_j\mathcal{E}(k_2)\rangle \\
&=& \left\langle\mathcal{E}(k_1)\otimes\left(\begin{array}{c} ik_1-i\overline{k}_2 \\ k_1+\overline{k}_2\end{array}\right)   , u \mathcal{E}(k_2)\right\rangle
.
\end{eqnarray*}
\end{proof}

\begin{corollary}
The divergence operator $\delta$ is closable in the weak topology.
\end{corollary}
\begin{proof}
Let $(u_n)_{n\in\mathbb{N}}$ be a sequence such that $u_n\to 0$ and $\delta(u_n)\to \beta\in\mathcal{B}(\eufrak{H})$ in the weak topology. Then we get
\begin{eqnarray*}
\langle\mathcal{E}(k_1), \beta \mathcal{E}(k_2)\rangle &=& \lim_{n\to\infty}\langle\mathcal{E}(k_1), \delta(u_n) \mathcal{E}(k_2)\rangle \\
&=& \lim_{n\to\infty}
 \left\langle\mathcal{E}(k_1)\otimes\left(\begin{array}{c} ik_1-i\overline{k}_2 \\ k_1+\overline{k}_2\end{array}\right)   , u_n \mathcal{E}(k_2)\right\rangle \\
&=& 0
\end{eqnarray*}
for all $k_1,k_2\in\eufrak{h}_\mathbb{C}$, and thus $\beta=0$.
\end{proof}

We have the following analogues of Equations (\ref{comm rel1}) and (\ref{DFu1}).

\begin{proposition}
Let $u,v\in \mathcal{S}_{\eufrak{h},\delta}$, $F\in\mathcal{S}$, $h\in\eufrak{h}_\mathbb{C}\otimes\mathbb{C}^2$, then we have
\begin{eqnarray}
D_h\circ \delta(u) &=& \langle \overline{h},u\rangle + \delta \circ D_h(u) \label{comm rel}\\
\delta(Fu) &=& F\delta(u) - F\overleftarrow{D}_u  +\frac{1}{2}  \sum_{j=1}^n \left[P\left(h^{(j)}_1\right)+Q\left(h^{(j)}_2\right),F\right] F_j  \label{DFu}\\
\delta(uF) &=&  \delta(u)F - \overrightarrow{D}_uF  + \frac{1}{2} \sum_{j=1}^n \left[F,P\left(h^{(j)}_1\right)+Q\left(h^{(j)}_2\right)\right] F_j  \label{DuF}
\end{eqnarray}
\end{proposition}
\begin{proof}
\begin{enumerate}
\item
Let $u=\sum_{j=1}^n F_j\otimes h^{(j)}$. We set
\begin{eqnarray*}
X_j &=& \frac{1}{2}\left(P\left(h^{(j)}_1\right)+Q\left(h^{(j)}_2\right)\right),\\
Y_j &=& \frac{i}{2}\left(Q\left(h^{(j)}_1\right)-P\left(h^{(j)}_2\right)\right),\\
 Y &=& \frac{i}{2}\big(Q\left(h_1\right)-P\left(h_2\right)\big),
\end{eqnarray*}
then we have
$\delta(u)=\sum_{j=1}^n \big((X_j-Y_j)F_j + F_j(X_j+Y_j)\big)$,
and therefore
\[
D_h\big(\delta(u)\big) = \sum_{j=1}^n \big(Y(X_j-Y_j)F_j + Y F_j(X_j+Y_j)- (X_j-Y_j)F_jY - F_j(X_j+Y_j)Y \big).
\]
On the other hand we have $D_h(u)=\sum_{j=1}^n (YF_j-F_jY)\otimes h^{(j)}$, and
\[
\delta\big(D_h(u)\big) = \sum_{j=1}^n \big((X_j-Y_j)YF_j- (X_j-Y_j)F_jY + YF_j(X_j+Y_j)- F_jY(X_j+Y_j)\big).
\]
Taking the difference of these two expressions, we get
\begin{eqnarray*}
D_h\big(\delta(u)\big)-\delta\big(D_h(u)\big) &=&\sum_{j=1}^n \big(\big[Y,X_j-Y_j\big]F_j+ F_j[Y,X_j+Y_j]\big) \\
&=& \sum_{j=1}^n \big(\langle \overline{h}_1,h^{(j)}_1\rangle + \langle\overline{h}_2,h^{(j)}_2\rangle\big) F_j = \langle \overline{h},u\rangle.
\end{eqnarray*}
\item
A straightforward computation gives
\begin{eqnarray*}
\delta(Fu)&=& \sum_{j=1}^n \big((X_j-Y_j)FF_j + FF_j(X_j+Y_j)\big) \\
&=&F \sum_{j=1}^n \big((X_j-Y_j)F_j + F_j(X_j+Y_j)\big) - \sum_{j=1}^n [F,X_j-Y_j] F_j \\
&=& F \delta(u) - \sum_{j=1}^n [Y_j,F] F_j +  \sum_{j=1}^n [X_jF,X] F_j \\
&=&   F \delta(u) - \sum_{j=1}^n \langle F^*,h^{(j)}\rangle F_j +  \sum_{j=1}^n [X_j,F] F_j \\
&=& F\delta(u) - F\overleftarrow{D}_u  +  \sum_{j=1}^n [X_j,F] F_j
\end{eqnarray*}
where we used that  $[X_j,F]=i\left\langle\left(\begin{array}{c}-\overline{h}_2 \\ \overline{h}_1\end{array}\right) ,DF\right\rangle$ defines a bounded operator, since $F\in\mathcal{S}\subseteq{\rm Dom}\,D$. Equation (\ref{DuF}) can be shown similarly.
\end{enumerate}
\end{proof}

If we impose additional commutativity conditions, which are always satisfied in the commutative case, then we get simpler formulas that are more similar to the classical ones.

\begin{corollary}
If $u=\sum_{j=1}^n F_j\otimes h^{(j)}\in\mathcal{S}_{\eufrak{h},\delta}$ and
\[
F\in {\rm Dom}\,D\cap\left\{ P\left(h^{(1)}_1\right)+Q\left(h^{(1)}_2\right),\ldots, P\left(h^{(n)}_1\right)+Q\left(h^{(n)}_2\right)\right\}'
\]
then we have
\begin{eqnarray*}
\delta(Fu) &=& F\delta(u) - F\overleftarrow{D}_u, \\
\delta(uF) &=&\delta(u)F - \overrightarrow{D}_uF.
\end{eqnarray*}
\end{corollary}

\section{Examples and applications}

\subsection{Relation to the commutative case}

In this section we will show that the non-commutative calculus studied here contains the commutative calculus, at least if we restrict ourselves to bounded functionals.

It is well-known that the symmetric Fock space $\Gamma(\eufrak{h}_\mathbb{C})$ is isomorphic to the complexification $L^2(\Omega;\mathbb{C})$ of the Wiener space $L^2(\Omega)$ over $\eufrak{h}$, cf.\ \cite{biane93,janson97,meyer95}. Such an isomorphism $I:L^2(\Omega;\mathbb{C})\overset{\cong}{\mapsto}\Gamma(\eufrak{h}_\mathbb{C})$ can be defined by extending the map
\[
I: e^{iW(h)}\mapsto I\left(e^{iW(h)}\right) = e^{iQ(h)}\Omega= e^{-||h||^2/2} \mathcal{E}(ih), \qquad \text{ for } h\in\eufrak{h}.
\]
Using this isomorphism, a bounded functional $F\in L^\infty(\Omega;\mathbb{C})$ becomes a bounded operator $M(F)$ on $\Gamma(\eufrak{h}_\mathbb{C})$, acting simply by multiplication,
\[
M(F)\psi= I\big(FI^{-1}(\psi)\big), \qquad \text{ for } \psi\in\Gamma(\eufrak{h}_\mathbb{C}).
\]
In particular, we get $M\left(e^{iW(h)}\right)=U(0,h)$ for $h\in\eufrak{h}$.

We can show that the derivation of a bounded differentiable functional coincides with its derivation as a bounded operator.
\begin{proposition}
Let $k\in\eufrak{h}$ and $F\in L^\infty(\Omega;\mathbb{C})\cap {\rm Dom}\,\tilde{D}_k$ s.t.\ $\tilde{D}_k F\in L^\infty(\Omega;\mathbb{C})$. Then we have $M(F)\in {\rm Dom}\,D_{k_0}$, where $k_0=\left(\begin{array}{c} 0 \\ k\end{array}\right)$, and
\[
M(\tilde{D}_k F) = D_{k_0} \big(M(F)\big).
\]
\end{proposition}
\begin{proof}
It is sufficient to check this for functionals of the form $F=e^{iW(h)}$, $h\in\eufrak{h}$. We get
\begin{eqnarray*}
M(\tilde{D}_k e^{iW(h)} ) &=& M\left(i\langle k,h\rangle e^{iW(h)}\right) \\
&=& i\langle k,h\rangle U(0,h) = i\left\langle\left(\begin{array}{c} 0 \\ k\end{array}\right),\left(\begin{array}{c} 0 \\ h\end{array}\right) \right\rangle U(0,h) \\
&=& D_{k_0} U(0,h) = D_ {k_0}\big(M(e^{iW(h)})\big).
\end{eqnarray*}
\end{proof}
This implies that we also have an analogous result for the divergence.

\subsection{Sufficient conditions for the existence of smooth densities}\label{mall-smooth}

In this section we will use the operator $D$ to give sufficient conditions for the existence and smoothness of densities for operators on $\eufrak{H}$. The first result is a generalization of Proposition \ref{prop-regular1} to arbitrary states.

\begin{theorem}\label{theo regular2}
Let $\kappa\in\mathbb{N}$, $h\in\eufrak{h}\otimes\mathbb{R}^2$ with $\langle h_1,h_2\rangle\not=0$, and suppose that $\Phi$ is of the form 
\[
\Phi(X)={\rm tr}(\rho X) \text{ for all }X\in B(\eufrak{H}),
\]
for some density matrix $\rho$.
If there exist $k,\ell\in\eufrak{h}_\mathbb{C}\otimes\mathbb{C}^2$ such that
\[
\det\left(\begin{array}{cc}\langle h_1,k_1 \rangle & \langle h_2, k_2\rangle\\
\langle h_1,\ell_1 \rangle& \langle h_2, \ell_2\rangle\end{array}\right)\not=0,
\]
and $\rho\in\bigcap_{\kappa_1+\kappa_2\le\kappa} {\rm Dom}\,D_k^{\kappa_1}D_\ell^{\kappa_2}$, and
\[
{\rm tr}(|D_k^{\kappa_1}D_\ell^{\kappa_2}\rho|) < \infty \qquad\text{ for all } \kappa_1+\kappa_2\le \kappa,
\]
then $w_{h,\Phi}\in \bigcap_{2\le p\le\infty}H^{p,\kappa}(\mathbb{R}^2)$, i.e.\ the Wigner density $w_{h,\Phi}$ lies in the Sobolev spaces of order $\kappa$. 
\end{theorem}
\begin{proof}
Let
\[
A=\left(\begin{array}{cc} 
\langle h_1,k_1 \rangle & \langle h_2, k_2\rangle \\{}
\langle h_1,\ell_1 \rangle & \langle h_2, \ell_2\rangle
\end{array}\right)
\]
and set
\[
\left(\begin{array}{c} X_1 \\  X_2 \end{array}\right) = \frac{i}{2}\, A^{-1} \left(\begin{array}{c} Q(k_1)-P(k_2) \\ Q(\ell_1)-P(\ell_2) \end{array}\right), 
\]
then we have
\begin{eqnarray*}
[X_1,O_h(\varphi)]&=&\frac{\langle h_2, \ell_2\rangle D_k O_h(\varphi) - \langle h_2, k_2\rangle D_\ell O_h(\varphi)}{\det A}= O_h\left(\frac{\partial \varphi}{\partial x}\right), \\{}
[X_2,O_h(\varphi)]&=& \frac{-\langle h_1, \ell_1\rangle D_k O_h(\varphi) + \langle h_1, k_2\rangle D_\ell O_h(\varphi)}{\det A}=  O_h\left(\frac{\partial \varphi}{\partial y}\right),
\end{eqnarray*}
for all Schwartz functions $\varphi$. Therefore
\begin{eqnarray*}
\left|\int\frac{\partial^{\kappa_1+\kappa_2}\varphi}{\partial x^{\kappa_1}\partial y^{\kappa_2}}{\rm d}W_{h,\Phi}\right| &=& \left|{\rm tr}\left( \rho\, O_h\left( \frac{\partial^{\kappa_1+\kappa_2}\varphi}{\partial x^{\kappa_1}\partial y^{\kappa_2}}\right)\right)\right| \\
&=&\left|{\rm tr}\Big(\rho \underbrace{\big[ X_1,\ldots\big[X_1},\underbrace{\big[X_2,\ldots\big[X_2},O_h(\varphi)\big]\big]\big]\big]\Big)\right| \\
&&\hspace{17mm} \kappa_1 \text{ times }\hspace{7mm}\kappa_2\text{ times } \\[1mm]
&=&\left|{\rm tr}\big([ X_2,\ldots[X_2,[X_1,\ldots[X_1,\rho]]]]O_h(\varphi)\big)\right| \\[1mm]
&\le &C_{\rho,\kappa_1,\kappa_2}|| O_h(\varphi)|| \le C_{\rho,\kappa_1,\kappa_2} C_{h,p} ||\varphi||_p,
\end{eqnarray*}
for all $p\in[1,2]$, since $\rho\in\bigcap_{\kappa_1+\kappa_2\le\kappa} {\rm Dom}\,D_k^{\kappa_1}D_\ell^{\kappa_2}$ and ${\rm tr}(|D_k^{\kappa_1}D_\ell^{\kappa_2}\rho|) < \infty$ for all $\kappa_1+\kappa_2\le \kappa$, and thus
\[
C_{\rho,\kappa_1,\kappa_2} = {\rm tr}\big|[ X_2,\ldots[X_2,[X_1,\ldots[X_1,\rho]]]]\big| < \infty.
\]
But this implies that the density of ${\rm d}W_{h,\Phi}$ is contained in the Sobolev spaces $H^{p,\kappa}(\mathbb{R}^2)$ for all $2\le p\le \infty$.
\end{proof}

\begin{example}
Let $0< \lambda_1\le \lambda_2\le \cdots$ be an increasing sequence of positive numbers and $\{e_j|j\in\mathbb{N}\}$ a complete orthonormal system for $\eufrak{h}_\mathbb{C}$. Let $T_t:\eufrak{h}_\mathbb{C}\to\eufrak{h}_\mathbb{C}$ be the contraction semigroup defined by
\[
T_t e_j = e^{-t\lambda_j}e_j ,\qquad \text{ for } j\in\mathbb{N}, t\ge 0,
\]
with generator $A=\sum_{j\in\mathbb{N}} \lambda_j \mathbb{P}_j$. If the sequence increases fast enough to ensure that $\sum_{j=1}^\infty e^{-t\lambda_j}<\infty$, i.e.\ if ${\rm tr}\,T_t<\infty$ for $t>0$, then the second quantization $\rho_t=\Gamma(T_t):\eufrak{H}\to\eufrak{H}$ is a trace class operator with trace
\[
Z_t ={\rm tr}\, \rho_t = \sum_{\mathbf{n}\in\mathbb{N}_f^\infty} \langle e_{\mathbf{n}}, \rho_t e_{\mathbf{n}}\rangle
\]
where we use $\mathbb{N}_f^\infty$ to denote the finite sequences of non-negative integers and $\{e_\mathbf{n}|\mathbf{n}\in\mathbb{N}_f^\infty\}$ is the complete orthonormal system of $\eufrak{H}$ consisting of the vectors
\[
e_\mathbf{n}=e_1^{\odot n_1}\odot \cdots \odot e_r^{\odot n_r}, \qquad \text{ for } \mathbf{n}=(n_1,\ldots,n_r)\in \mathbb{N}_f^\infty,
\]
i.e.\ the symmetrization of the tensor $e_1\otimes \cdots \otimes e_1 \otimes \cdots\otimes e_r\otimes \cdots \otimes e_r$ where each vector $e_j$ appears $n_j$ times.
We get
$Z_t=\sum_{\mathbf{n}\in\mathbb{N}^\infty} \prod_{k=1}^\infty e^{-n_kt\lambda_k} = \prod_{k=1}^\infty \frac{1}{1-e^{-t\lambda_k}}$
for the trace of $\rho_t$. We shall be interested in the state defined by
\[
\Phi(X)=\frac{1}{Z_t}{\rm tr}(\rho_tX), \qquad\text{ for } X\in\mathcal{B}(\eufrak{H}).
\]
We get
\begin{eqnarray*}
\sum_{\mathbf{n}\in\mathbb{N}_f^\infty} \big|\langle e_\mathbf{n}, |\rho_{t/2} a^\ell(e_j)|^2e_\mathbf{n}\rangle\big| &=& \sum_{\mathbf{n}\in\mathbb{N}_f^\infty} || \rho_{t/2} a^\ell(e_j)||^2 \\
&=& \sum_{\mathbf{n}\in\mathbb{N}^\infty} n_j(n_j-1)\cdots (n_j-\ell+1) e^{-(n_j-\ell)t\lambda_j}\prod_{k\not=j} e^{-n_kt\lambda_k} \\
&\le & \sum_{n=0}^\infty (n+\ell)^{\ell} e^{-nt\lambda_j}   \prod_{k\not=j} \frac{1}{1-e^{-t\lambda_k}} < \infty,
\end{eqnarray*}
and therefore that $\rho_t a^\ell(e_j)$ defines a bounded operator with finite trace for all $j,\ell\in\mathbb{N}$ and $t>0$. Similarly we get
\[
{\rm tr}\, \big|a^\ell (e_j)\rho_t\big| < \infty, \quad
{\rm tr}\, \big|\rho_t \big(a^+(e_j)\big)^\ell\big| < \infty, \quad
\text{ etc. }
\]
and
\[
{\rm tr}\, \big|P^{\ell_1}(e_{j_1}) Q^{\ell_2}(e_{j_2})\rho_t\big| < \infty, \quad
{\rm tr}\, \big|P^{\ell_1}(e_{j_1}) \rho_t Q^{\ell_2}(e_{j_2})\big| < \infty, \quad \text{ etc. }
\]
for all $t>0$ and $j_1,j_2,\ell_1,\ell_2\in\mathbb{N}$. For a given $h\in\eufrak{h}\otimes\mathbb{R}^2$ with $\langle h_1,h_2\rangle\not=0$ (and thus in particular $h_1\not=0$ and $h_2\not=0$), we can always find indices $j_1$ and $j_2$ such that $\langle h_1,e_{j_1}\rangle \not=0$ and $\langle h_2,e_{j_2}\rangle \not=0$. Therefore it is not difficult to check that for all $\kappa\in\mathbb{N}$ all the conditions of Theorem \ref{theo regular2} are satisfied with $k=\left(\begin{array}{c}e_{j_1} \\ 0 \end{array}\right)$ and $\ell =\left(\begin{array}{c}0 \\ e_{j_2} \end{array}\right)$. We see that the Wigner density $w_{h,\Phi}$ of any pair $\big(P(h_1),Q(h_2)\big)$ with $\langle h_1,h_2\rangle\not=0$ in the state $\Phi(\cdot)= {\rm tr}(\rho_t\,\cdot)/Z_t$ is in $\bigcap_{\kappa\in\mathbb{N}} \bigcap_{2\le p\le\infty}H^{p,\kappa}(\mathbb{R}^2)$, in particular, its derivatives of all orders exist, and are bounded and square-integrable.
\end{example}

We will now show that this approach can also be applied to get sufficient conditions for the regularity of a single bounded self-adjoint operator, for simplicity we consider only vector states.

Given a bounded self-adjoint operator $X$, we call a measure $\mu_{X,\Phi}$ on the real line its distribution in the state $\Phi$, if all moments agree,
\[
\Phi(X^n) = \int_\mathbb{R} x^n{\rm d}\mu_{X,\Phi}, \qquad \text{ for all } n\in\mathbb{N}.
\]
Such a measure $\mu_{X,\Phi}$ always exists, it is unique and supported on the interval\\ $\big[-||X||,||X||\big]$.

\begin{proposition}\label{prop-reg-X}
Let $\omega\in\eufrak{H}$ be a unit vector and let $\Phi(\cdot)=\langle\omega,\cdot\,\omega\rangle$ be the corresponding vector state. The distribution $\mu_{X,\Phi}$ of an operator $X\in\mathcal{B}(\eufrak{H})$ in the state $\Phi$ has a bounded density, if there exists a $k\in\eufrak{h}_\mathbb{C}\otimes\mathbb{C}^2$ such that $\omega\in{\rm Dom}\big(Q(k_1)-P(k_2)\big)\cap{\rm Dom}\big(Q(\overline{k}_1)-P(\overline{k}_2)\big)$, $X\in{\rm Dom}\,D_k$, $X\cdot D_kX=D_kX\cdot X$, $D_kX$ invertible and $(D_kX)^{-1}\in {\rm Dom}\,D_k$.
\end{proposition}
\begin{proof}
Since $X\cdot D_kX=D_kX\cdot X$, we have $D_kp(X)= (D_kX)p'(X)$ for all polynomials $p$. We therefore get
\[
D_k\left( (D_kX)^{-1}p(X)\right) = p(X) D_k\left( (D_kX)^{-1}\right) + p'(X).
\]
The hypotheses of the proposition assure
\begin{eqnarray*}
\left|\langle\omega,D_k\left( (D_kX)^{-1}p(X)\right)\omega\rangle\right|&\!\le\!& \frac{||\big(\!Q(k_1)\!-\!P(k_2)\!\big)\omega||\!+\!||\big(\!Q(\overline{k}_1)\!-\!P(\overline{k}_2)\!\big)\omega||}{2}||(D_kX)^{-1}||\, ||p(X)|| \\
&\!\le\! &C_1 \sup_{x\in [ -||X||,||X||]} \big|p(x)\big|, \\[2mm]
\left|\langle\omega, p(X) D_k\left( (D_kX)^{-1}\right)\omega\rangle\right| 
&\!\le\! & \left|\left|D\left((D_kX)^{-1}\right)\right|\right|\, ||p(X)|| \\
&\!\le\! &C_2 \sup_{x\in [ -||X||,||X||]} \big|p(x)\big|,
\end{eqnarray*}
and therefore allow us to get the estimate
\begin{eqnarray*}
\left|\int_{-||X||}^{||X||} p'(x) {\rm d}\mu_{X,\Phi}(x)\right| &=&
\big|\langle\omega, p'(X)\omega\rangle\big| \\
&=& \Big|\Big\langle\omega,\left(D_k\Big( (D_kX)^{-1}p(X)\right)-p(X) D_k\left( (D_kX)^{-1}\right)\Big)\omega\Big\rangle\Big| \\
&\le & (C_1+C_2) \sup_{x\in [ -||X||,||X||]} \big|p(x)\big|
\end{eqnarray*}
for all polynomials $p$.
But this implies that $\mu_{X,\Phi}$ admits a bounded density.
\end{proof}

Let $n\in\mathbb{N}$, $n\ge 2$. We get that the density is even $n-1$ times differentiable, if in addition to the conditions of Proposition \ref{prop-reg-X} we also have
$X\in {\rm Dom}\,D_k^n$, $(D_kX)^{-1}\in {\rm Dom}\,D_k^n$, and
\[
\omega\in\bigcap_{1\le\kappa\le n}{\rm Dom}\big(Q(k_1)-P(k_2)\big)^\kappa \bigcap_{1\le\kappa\le n}{\rm Dom}\big(Q(\overline{k}_1)-P(\overline{k}_2)\big)^\kappa.
\]

The proof is similar to the proof of Proposition \ref{prop-reg-X}, we now use the formula
\[
p^{(n)}(X) = D_k^n\left( (D_kX)^{-n}p(X)\right) - \sum_{\kappa=0}^{n-1} A_\kappa p^{(\kappa)}(X),
\]
where $A_0,\ldots,A_{n-1}$ are some bounded operators, to get the necessary estimate
\[
\left|\int_{-||X||}^{||X||} p^{(n)}(x) {\rm d}\mu_{X,\Phi}(x)\right| \le
C \sup_{x\in [ -||X||,||X||]} \big|p(x)\big|
\]
by induction over $n$.

\subsection{The white noise case}\label{white noise}

Let now $\eufrak{h}=L^2(T,\mathcal{B},\mu)$, where $(T,\mathcal{B},\mu)$ is a measure space such that $\mathcal{B}$ is countably generated. In this case we can apply the divergence operator to processes indexed by $T$, i.e.\ $\mathcal{B}(\eufrak{H})$-valued measurable functions on $T$, since they can be interpreted as elements of the Hilbert module, if they are square-integrable. Let $L^2\big(T,\mathcal{B}(\eufrak{H})\big)$ denote all $\mathcal{B}(\eufrak{H})$-valued measurable functions $t\mapsto X_t$ on $T$ with $\int_T ||X_t||^2 {\rm d}t< \infty$.
Then the definition of the divergence operator becomes
\begin{eqnarray*}
{\rm Dom}\,\delta &=& \Big\{ X=(X^1,X^2)\in L^2\big(T,\mathcal{B}(\eufrak{H})\big)\oplus L^2\big(T,\mathcal{B}(\eufrak{H})\big)\Big| \\
&&\exists M\in\mathcal{B}(\eufrak{H}), \forall k_1,k_2\in \eufrak{h}_\mathbb{C}: \langle\mathcal{E}(k_1), M\mathcal{E}(k_2)\rangle \\
&&= \int_T\big(i(k_2-\overline{k}_1)\langle\mathcal{E}(k_1), X^1(t)\mathcal{E}(k_2)\rangle+(\overline{k}_1+k_2) \langle\mathcal{E}(k_1), X^2(t)\mathcal{E}(k_2)\rangle\big){\rm d}\mu(t)\Big\},
\end{eqnarray*}
and $\delta(X)$ is  equal to the unique operator satisfying the above condition. We will also use the notation
\[
\delta(X)=\int_T X^1(t){\rm d}P(t) + \int_T X^2(t){\rm d}Q(t),
\]
and call $\delta(X)$ the Hitsuda-Skorohod integral of $X$.

Belavkin \cite{belavkin91a,belavkin91b} and Lindsay \cite{lindsay93} have defined non-causal quantum stochastic integrals with respect to the creation, annihilation, and conservation processes on the boson Fock space over $L^2(\mathbb{R}_+)$ using the classical derivation and divergence operators. It turns out that our Hitsuda-Skorohod integral coincides with their creation and annihilation integrals, up to a coordinate transformation. This immediately implies that for adapted, integrable processes our integral coincides with the quantum stochastic creation and annihilation integrals defined by Hudson and Parthasarathy, cf.\ \cite{hudson+parthasarathy84,parthasarathy92}.

Let us briefly recall, how they define the creation and annihilation integral, cf.\ \cite{lindsay93}. They use the derivation operator $\tilde{D}$ and the divergence operator $\tilde{\delta}$ from the classical calculus on the Wiener space $L^2(\Omega)$, defined on the Fock space $\Gamma(L^2(\mathbb{R}_+;\mathbb{C}))$ over $L^2(\mathbb{R}_+;\mathbb{C})=L^2(\mathbb{R}_+)_\mathbb{C}$ via the isomorphism between $L^2(\Omega)$ and $\Gamma\big(L^2(\mathbb{R}_+;\mathbb{C})\big)$. On the exponential vectors $\tilde{D}$ acts as
\[
\tilde{D}\mathcal{E}(k)=\mathcal{E}(k)\otimes k, \qquad \text{ for } k\in L^2(\mathbb{R}_+,\mathbb{C}),
\]
and $\tilde{\delta}$ is its adjoint. Note that due to the isomorphism between $\Gamma\big(L^2(\mathbb{R}_+;\mathbb{C})\big)\otimes L^2(\mathbb{R}_+;\mathbb{C})$ and $L^2\Big(\mathbb{R}_+;\Gamma\big(L^2(\mathbb{R}_+;\mathbb{C})\big)\Big)$, the elements of  $\Gamma\big(L^2(\mathbb{R}_+;\mathbb{C})\big)\otimes L^2(\mathbb{R}_+;\mathbb{C})$ can be interpreted as function on $\mathbb{R}_+$. In particular, for the exponential vectors we get $\big(D\mathcal{E}(k)\big)_t=k(t)\mathcal{E}(k)$ almost surely. The action of the annihilation integral $\int F_t dA_t$ on some vector $\psi\in\Gamma\big(L^2(\mathbb{R}_+;\mathbb{C})\big) $ is then defined as the Bochner integral
\[
\int^{\rm BL}_{\mathbb{R}_+} F_t dA_r\psi = \int_{\mathbb{R}_+} F_t (D\psi)_t {\rm d}t,
\]
and that of the creation integral as
\[
\int^{\rm BL}_{\mathbb{R}_+} F_t {\rm d}A^*_t\psi = \tilde{\delta}(F_\cdot\psi).
\]
These definitions satisfy the adjoint relations
\[
\left(\int^{\rm BL}_{\mathbb{R}_+} F_t dA_r\right)^* \supset \int^{\rm BL}_{\mathbb{R}_+} F^*_t {\rm d}A^*_t, \quad \text{ and } \quad \left(\int^{\rm BL}_{\mathbb{R}_+} F_t dA^*_r\right)^* \supset \int^{\rm BL}_{\mathbb{R}_+} F^*_t {\rm d}A_t.
\]

\begin{proposition}
Let $(T,\mathcal{B},\mu)=(\mathbb{R}_+,\mathcal{B}(\mathbb{R}_+),{\rm d}x)$, i.e.\ the positive half-line with the Lebesgue measure, and let $X=(X^1,X^2)\in {\rm Dom}\,\delta$. Then we have
\[
\int_{\mathbb{R}_+} X^1(t){\rm d}P(t) + \int_{\mathbb{R}_+} X^2(t){\rm d}Q(t) = \int^{\rm BL}_{\mathbb{R}_+}(X^2-iX^1) {\rm d}A^*_t + \int^{\rm BL}_{\mathbb{R}_+}(X^2+iX^1){\rm d}A_t.
\]
\end{proposition}
\begin{proof}
To prove this, we show that the Belavkin-Lindsay integrals satisfy the same formula for the matrix elements between exponential vectors. Let $(F_t)_{t\in\mathbb{R}_+}\in L^2\big(\mathbb{R}_+,\mathcal{B}(\eufrak{H})\big)$ be such that its creation integral in the sense of Belavkin and Lindsay is defined with a domain containing the exponential vectors. Then we get
\begin{eqnarray*}
\langle\mathcal{E}(k_1),\int^{\rm BL}_{\mathbb{R}_+} F_t {\rm d}A^*_t\mathcal{E}(k_2)\rangle &=& \langle\mathcal{E}(k_1), \tilde{\delta}\big(F_\cdot\mathcal{E}(k_2)\big)\rangle \\
&=& \langle\big(\tilde{D}\mathcal{E}(k_1)\big)_\cdot, F_\cdot\mathcal{E}(k_2)\rangle \\
&=& \int\overline{k_1(t)} \langle\mathcal{E}(k_1), F_t\mathcal{E}(k_2)\rangle{\rm d}t.
\end{eqnarray*}
For the annihilation integral we deduce the formula
\begin{eqnarray*}
\langle\mathcal{E}(k_1),\int^{\rm BL}_{\mathbb{R}_+} F_t {\rm d}A_t\mathcal{E}(k_2)\rangle &=& \langle\int^{\rm BL}_{\mathbb{R}_+} F^*_t {\rm d}A^*_t\mathcal{E}(k_1),\mathcal{E}(k_2)\rangle \\
&=& \overline{\int_{\mathbb{R}_+}\overline{k_2(t)} \langle\mathcal{E}(k_2), F^*_t\mathcal{E}(k_1)\rangle{\rm d}t} \\
&=&\int_{\mathbb{R}_+} k_2(t) \langle\mathcal{E}(k_1), F_t\mathcal{E}(k_2)\rangle{\rm d}t.
\end{eqnarray*}
\end{proof}

The integrals defined by Belavkin and Lindsay are an extension of those defined by Hudson and Parthasarathy, therefore the following is now obvious.

\begin{corollary}
For adapted processes $X\in{\rm Dom}\,\delta$, the Hitsuda-Skorohod integral
\[
\delta(X)=\int_T X^1(t){\rm d}P(t) + \int_T X^2(t){\rm d}Q(t)
\]
coincides with the Hudson-Parthasarathy quantum stochastic integral defined in \cite{hudson+parthasarathy84}.
\end{corollary}

\subsection{Iterated integrals}

We give a short, informal discussion of iterated integrals of deterministic functions, showing a close relation between these iterated integrals and the so-called Wick product or normal-ordered product. Doing so, we will encounter unbounded operators, so that strictly speaking we have not defined $\delta$ for them. But everything could be made rigorous by choosing an appropriate common invariant domain for these operators, e.g., vectors with a finite chaos decomposition.

In order to be able to iterate the divergence operator, we define $\delta$ on $\mathcal{B}(\eufrak{H})\otimes \eufrak{h}_\mathbb{C}\otimes\mathbb{C}^2\otimes H$, where $H$ is some Hilbert space, as $\delta\otimes {\rm id}_{H}$.

Using Equation (\ref{wick}), on can show by induction
\[
\delta^n \left(\begin{array}{c} h^{(1)}_1 \\ h^{(1)}_2 \end{array}\right) \otimes \cdots \otimes \left(\begin{array}{c} h^{(n)}_1 \\ h^{(n)}_2 \end{array}\right)
= \sum_{I\subseteq\{1,\ldots,n\}} \prod_{j\in I} a^+(h^{(j)}_2-i h^{(j)}_1) \prod_{j\in \{1,\ldots,n\}\backslash I} a(\overline{h}^{(j)}_2+i \overline{h}^{(j)}).
\]
for $h^{(1)}=\left(\begin{array}{c} h^{(1)}_1 \\ h^{(1)}_2 \end{array}\right),\ldots,h^{(n)}=\left(\begin{array}{c} h^{(n)}_1 \\ h^{(n)}_2 \end{array}\right)\in\eufrak{h}_\mathbb{C}\otimes\mathbb{C}^2$. This is just the Wick product of $\big(P(h^{(1)}_1+Q(h^{(1)}_2)\big),\ldots,\big(P(h^{(n)}_1+Q(h^{(n)}_2)\big)$, i.e.
\[
\delta^n \left(\begin{array}{c} h^{(1)}_1 \\ h^{(1)}_2 \end{array}\right) \otimes \cdots \otimes \left(\begin{array}{c} h^{(n)}_1 \\ h^{(n)}_2 \end{array}\right)
= \big(P(h^{(1)}_1+Q(h^{(1)}_2)\big)\diamond\cdots\diamond\big(P(h^{(n)}_1+Q(h^{(n)}_2)\big),
\]
where the Wick product $\diamond$ is defined on the algebra generated by $\{P(k),Q(k)|k\in\eufrak{h}_\mathbb{C}\}$ by
\begin{eqnarray*}
P(h)\diamond X &=& X\diamond P(h) = -ia^+(h)X +iX a(\overline{h}) \\
Q(h)\diamond X &=& X\diamond Q(h) = a^+(h)X+ Xa(\overline{h})
\end{eqnarray*}
for $X\in {\rm alg}\,\{P(k),Q(k)|k\in\eufrak{h}_\mathbb{C}\}$ and $h\in\eufrak{h}_\mathbb{C}$ in terms of the momentum and position operators, or, equivalently, by
\begin{eqnarray*}
a^+(h)\diamond X &=& X\diamond a^+(h) =  a^+(h)X \\
a(h)\diamond X &=& X\diamond a(h) = Xa(h)
\end{eqnarray*}
in terms of creation and annihilation.

\section{Conclusion}

We have defined a derivation operator $D$ and a divergence operator $\delta$ on $\mathcal{B}(\eufrak{H})$ and $\mathcal{B}(\eufrak{H},\eufrak{H}\otimes\eufrak{h}_\mathbb{C}\otimes\mathbb{C}^2)$, resp., and shown that they have similar properties as the derivation operator and the divergence operator in classical Malliavin calculus. As far as we know, this is the first time that $D$ and $\delta$ are considered as operators defined on a non-commutative operator algebra, except for the free case \cite{biane+speicher99}, where the operator algebra is isomorphic to the full Fock space. To obtain close analogues of the classical relations involving the divergence operator, we needed to impose additional commutativity conditions, but Proposition \ref{prop no go} shows that this can not be avoided. Also, its domain is rather small, because we require $\delta(u)$ to be a bounded operator and so, e.g., the ``deterministic'' elements $h\in\eufrak{h}_\mathbb{C}\otimes\mathbb{C}^2$ are not integrable unless $h=0$. One of the main goals of our approach is the study of Wigner densities as joint densities of non-commutating random variables. We showed that the derivation operator can be used to obtain sufficient conditions for its regularity, see, e.g., Theorem \ref{theo regular2}. It seems likely that these results can be generalized by weakening or modifying the hypotheses. It would be interesting to apply these methods to quantum stochastic differential equations.

%% file: quasi.tex
\chapter{Quasi-Invariance Formulas for Components of Quantum L\'evy Processes}\label{ch-quasi}\markboth{QUASI-INVARIANCE FORMULAS}{QUASI-INVARIANCE FORMULAS}

\vspace*{3cm}

\begin{quote}
A general method for deriving Girsanov or quasi-invariance formulas for classical stochastic processes with independent increments obtained as components of L\'evy processes on real Lie algebras is presented. Letting a unitary operator arising from the associated factorizable current representation act on an appropriate commutative subalgebra, a second commutative subalgebra is obtained. Under certain conditions the two commutative subalgebras lead to two classical processes such that the law of the second process is absolutely continuous w.r.t.\ to the first. Examples include the Girsanov formula for Brownian motion as well as quasi-invariance formulas for the Poisson process, the Gamma process \cite{tsilevich+vershik+yor00,tsilevich+vershik+yor01}, and the Meixner process.
\end{quote}

\vspace*{3cm}

Joint work with Nicolas Privault. Published in Inf.\ Dim.\ Anal., Quantum Prob.\ and Rel.\ Topics Vol.~7, No.~1, pp.~131-145, 2004.

\newpage

\section{Introduction}

L\'evy processes, i.e.\ stochastic processes with independent and stationary increments, are used as models for random fluctuations, e.g., in physics, finance, etc. In quantum physics so-called quantum noises or quantum L\'evy processes occur, e.g., in the description of quantum systems coupled to a heat bath \cite{gardiner+zoller00} or in the theory of continuous measurement \cite{holevo01}. Motivated by a model introduced for lasers \cite{waldenfels84}, Sch\"urmann et al.\ \cite{accardi+schuermann+waldenfels88,schuermann93} have developed the theory of L\'evy processes on involutive bialgebras. This theory generalizes, in a sense, the theory of factorizable representations of current groups and current algebras as well as the theory of classical L\'evy processes with values in Euclidean space or, more generally, semigroups. For a historical survey on the theory of factorizable representations and its relation to quantum stochastic calculus, see \cite[Section 5]{streater00}.

Many interesting classical stochastic processes arise as components of these quantum L\'evy processes, cf.\ \cite{schuermann91b,biane98b,franz99,accardi+franz+skeide02}. In this Note we will demonstrate on several examples how quasi-invariance formulas can be obtained in such a situation. Our examples include the Girsanov formula for Brownian motion as well as a quasi-invariance formula for the Gamma process \cite{tsilevich+vershik+yor00,tsilevich+vershik+yor01}, which actually appeared first in the context of factorizable representations of current groups \cite{gelfand+graev+vershik83}. We also present a new quasi-invariance formula for the Meixner process.

We will consider L\'evy processes on real Lie algebras in this Note, but the general idea is more widely applicable. Furthermore, we restrict ourselves to commutative subalgebras of the current algebra that have dimension one at every point, see Subsection \ref{QUASI-class}. This allows us to use techniques familiar from the representation theory of Lie algebras and groups to get explicit expressions for the two sides of our quasi-invariance formulas.

In Section \ref{QUASI-prelim}, we present the basic facts about L\'evy processes on real Lie algebras and we recall how classical increment processes can be associated to them. The general idea of our construction is outlined in Section \ref{QUASI-idea}. Finally, in Section \ref{QUASI-examples}, we present explicit calculations for several classical increment processes related to the oscillator algebra and the Lie algebra $sl(2,\mathbb{R})$ of real $2\times2$ matrices with trace zero. In Subsections \ref{QUASI-exa-osc} and \ref{QUASI-exa-sl2}, we study representations of the Lie algebras themselves and obtain quasi-invariance formulas for random variables arising from these representations. In Subsections \ref{QUASI-exa girsanov}, \ref{QUASI-poisson girsanov}, \ref{QUASI-exa gamma}, and \ref{QUASI-exa-meixner}, we turn to L\'evy processes on these Lie algebras and derive quasi-invariance or Girsanov formulas for Brownian motion, the Poisson process, the Gamma process, and the Meixner process.

\section{L\'evy processes on real Lie algebras}\label{QUASI-prelim}

In this section we recall the definition and the main results concerning L\'evy processes on real Lie algebras, see also \cite{accardi+franz+skeide02}. This is a special case of the theory of L\'evy processes on involutive bialgebras, cf.\ \cite{schuermann93}\cite[Chapter VII]{meyer95},\cite{franz+schott99}.

\subsection{Definition and construction}
\begin{definition}\label{QUASI-def levy}
Let $\eufrak{g}$ be a real Lie algebra, $H$ a pre-Hilbert space, and $\Omega\in H$ a unit vector. We call a family $\big(j_{st}:\eufrak{g}\to\mathcal{L}(H)\big)_{0\le s\le t}$ of representations of $\eufrak{g}$ by anti-hermitian operators (i.e.\ $j_{st}(X)^*=-j_{st}(X)$ for all $0\le s\le t$, $X\in\eufrak{g}$) a {\em L\'evy process on $\eufrak{g}$} over $H$ (with respect to $\Omega$), if the following conditions are satisfied.
\begin{enumerate}
\item
(Increment property) We have
\[
j_{st}(X)+j_{tu}(X)=j_{su}(X)
\]
for all $0\le s\le t \le u$ and all $X\in\eufrak{g}$.
\item
(Independence)
We have $[j_{st}(X),j_{s't'}(Y)]=0$ for all $X,Y\in\eufrak{g}$, $0\le s\le t\le s'\le t'$ and
\[
\langle\Omega, j_{s_1t_1}(X_1)^{k_1}\cdots j_{s_nt_n}(X_n)^{k_n}\Omega\rangle = \langle\Omega, j_{s_1t_1}(X_1)^{k_1}\Omega\rangle\cdots\langle\Omega, j_{s_nt_n}(X_n)^{k_n}\Omega\rangle
\]
for all $n,k_1,\ldots,k_n\in\mathbb{N}$, $0\le s_1\le t_1\le s_2\le \cdots \le t_n$, $X_1,\ldots,X_n\in\eufrak{g}$.
\item
(Stationarity) For all $n\in\mathbb{N}$ and all $X\in\eufrak{g}$, the moments
\[
m_n(X;s,t)=\langle\Omega,j_{st}(X)^n\Omega\rangle
\]
depend only on the difference $t-s$.
\item
(Weak continuity)
We have $\lim_{t\searrow s} \langle\Omega,j_{st}(X)^n\Omega\rangle =0$ for all $n\ge 1$ and all $X\in\eufrak{g}$. 
\end{enumerate}
\end{definition}

Two L\'evy processes are called {\em equivalent}, if they have same moments.

For the classification and construction of these processes we introduce the notion of Sch\"urmann triples.

\begin{definition}
Let $\eufrak{g}$ be a real Lie algebra. A {\em Sch\"urmann triple on} $\eufrak{g}$ over some pre-Hilbert space $D$ is a triple $(\rho,\eta,L)$ consisting of
\begin{itemize}
\item
a representation $\rho$ of $\eufrak{g}$ on $D$ by anti-hermitian operators,
\item
a $\rho$-$1$-cocycle $\eta$, i.e.\ a linear map $\eta:\eufrak{g}\to D$ such that
\[
\eta\big([X,Y]\big) = \rho(X)\eta(Y)-\rho(Y)\eta(X),
\]
for all $X,Y\in\eufrak{g}$, and
\item
a linear functional $L:\eufrak{g}\to i\mathbb{R}\subseteq\mathbb{C}$ that has the map $\eufrak{g}\wedge\eufrak{g}\ni X\wedge Y\to \langle\eta(X),\eta(Y)\rangle-\langle\eta(Y),\eta(X)\rangle\in\mathbb{C}$ as a coboundary, i.e.\
\[
L\big([X,Y]\big) = \langle\eta(X),\eta(Y)\rangle-\langle\eta(Y),\eta(X)\rangle
\] 
for all $X,Y\in\eufrak{g}$.
\end{itemize}
A Sch\"urmann triple $(\rho,\eta,L)$ is called {\em surjective}, if $\eta(\eufrak{g})$ is cyclic for $\rho$.
\end{definition}

The following theorem can be traced back to the works of Araki, Streater, etc., in the form given here it is a special case of Sch\"urmann's representation theorem for L\'evy processes on involutive bialgebras, cf.\ \cite{schuermann93}.

\begin{theorem}
Let $\eufrak{g}$ be a real Lie algebra. Then there is a one-to-one correspondence between L\'evy processes on $\eufrak{g}$ (modulo equivalence) and Sch\"urmann triples on $\eufrak{g}$.

Let $(\rho,\eta,L)$ be a Sch\"urmann triple on $\eufrak{g}$ over $D$, then
\begin{equation}\label{QUASI-j-def}
j_{st}(X)=\Lambda_{st}\big(\rho(X)\big)+A^*_{st}\big(\eta(X)\big)-A_{st}\big(\eta(X)\big)+(t-s)L(X){\rm id}
\end{equation}
for $0\le s\le t$ and $X\in\eufrak{g}$ defines L\'evy process on $\eufrak{g}$ over a dense subspace $H\subseteq \Gamma\big(L^2(\mathbb{R}_+,D)\big)$ w.r.t.\ the vacuum vector $\Omega$.
\end{theorem}

\subsection{Regularity assumptions}\label{QUASI-reg assum}

In order to justify our calculations, we have to impose stronger conditions on L\'evy processes than those stated in Definition \ref{QUASI-def levy}. We will from now on assume that $(j_{st})_{0\le s\le t}$ is defined as in \eqref{QUASI-j-def} and that the representation $\rho$ in the Sch\"urmann triple can be exponentiated to a continuous unitary representation of the Lie group associated to $\eufrak{g}$. By Nelson's theorem this implies that $D$ contains a dense subspace whose elements are analytic vectors for all $\rho(X)$, $X\in\eufrak{g}$. Furthermore, we will assume that $\eta(\eufrak{g})$ consists of analytic vectors. These assumptions guarantee that $j_{st}$ can also be exponentiated to a continuous unitary group representation and therefore, by Nelson's theorem, any finite set of operators of the form $ij_{st}(X)$, $0\le s\le t$, $X\in\eufrak{g}$, is essentially selfadjoint on some common domain. Furthermore, the vacuum vector $\Omega$ is an analytic vector for all $j_{st}(X)$, $0\le s\le t$, $X\in\eufrak{g}$. These assumptions are satisfied in all the examples considered in Section \ref{QUASI-examples}.

\subsection{Classical processes}\label{QUASI-class}

Denote by $\eufrak{g}^{\mathbb{R}_+}$ the space of simple step functions with values in $\eufrak{g}$,
\[
\eufrak{g}^{\mathbb{R}_+} = \Big\{ X=\sum_{k=1}^n X_k\mathbf{1}_{[s_k,t_k[}\Big| 0\le s_1\le t_1\le s_2 \le \cdots \le t_n<\infty, X_1,\ldots,X_n\in\eufrak{g}\Big\}.
\]
Then $\eufrak{g}^{\mathbb{R}_+}$ is a real Lie algebra with the pointwise Lie bracket and any L\'evy process on $\eufrak{g}$ defines a representation $\pi$ of $\eufrak{g}^{\mathbb{R}_+}$ via
\begin{equation}\label{QUASI-fact-rep}
\pi(X) = \sum_{k=1}^n  j_{s_kt_k}(X_k)
\end{equation}
for $X=\sum_{k=1}^n X_k\mathbf{1}_{[s_k,t_k[}\in \eufrak{g}^{\mathbb{R}_+}$.

By choosing a commutative subalgebra of $\pi(\eufrak{g}^{\mathbb{R}_+})$ we can get a classical process.

Denote by $\mathcal{S}(\mathbb{R}_+)$ the algebra of real-valued simple step functions on $\mathbb{R}$,
\[
\mathcal{S}(\mathbb{R}_+) = \Big\{ f=\sum_{k=1}^n f_k\mathbf{1}_{[s_k,t_k[}\Big| 0\le s_1\le t_1\le s_2 \le \cdots \le t_n<\infty, f_1,\ldots,f_n\in\mathbb{R}\Big\},
\]
then the product $fX$ of an element $X\in\eufrak{g}^{\mathbb{R}_+}$ with a function $f\in\mathcal{S}(\mathbb{R}_+)$ is again in $\eufrak{g}^{\mathbb{R}_+}$.

\begin{theorem}
Let $(j_{st})_{0\le s\le t}$ be a L\'evy process on a real Lie algebra $\eufrak{g}$ and let $\pi$ be as in Equation \eqref{QUASI-fact-rep}. Choose $X\in\eufrak{g}^{\mathbb{R}_+}$ and define
\[
X(f)=i \pi(fX)
\]
for $f\in\mathcal{S}(\mathbb{R}_+)$.

Then there exists a classical stochastic process $(\hat{X}_t)_{t\ge 0}$ with independent increments that has the same finite distributions as $X$, i.e.
\[
\langle\Omega, g_1\big(X(f_1)\big)\cdots g_n\big(X(f_n)\big)\Omega\rangle = \mathbb{E}\Big(g_1\big(\hat{X}(f_1)\big)\cdots g_n\big(\hat{X}(f_n)\big)\Big)
\]
for all $n\in\mathbb{N}$, $f_1,\ldots,f_n\in\mathcal{S}(\mathbb{R}_+)$, $g_1,\ldots,g_n\in \mathcal{C}_0(\mathbb{R})$, where
\[
\hat{X}(f)=\int_{\mathbb{R}_+}f(t){\rm d}\hat{X}_t = \sum_{k=1}^n f_k(\hat{X}_{t_k}-\hat{X}_{s_k})
\]
for $f=\sum_{k=1}^n f_k\mathbf{1}_{[s_k,t_k[}\in\mathcal{S}(\mathbb{R}_+)$.
\end{theorem}
The existence of $(\hat{X}_t)_{t\ge 0}$ follows as in \cite[Section 4]{accardi+franz+skeide02}. Thanks to the regularity assumptions in \ref{QUASI-reg assum}, the operators $X(f_1),\ldots,X(f_n)$ are essentially self-adjoint and $g_1\big(X(f_1)\big),\ldots,g_n\big(X(f_n)\big)$ can be defined by the usual functional calculus for self-adjoint operators.

\section{Quasi-invariance formulas for components of quantum L\'evy processes}\label{QUASI-idea}

Let $(j_{st})_{0\le s\le t}$ be a L\'evy process on a real Lie algebra $\eufrak{g}$ and and fix $X\in\eufrak{g}^{\mathbb{R}_+}$ with classical version $(\hat{X}_t)_{t\ge 0}$. In order to get quasi-invariance formulas for $(\hat{X}_t)_{t\ge 0}$, we simply choose another element $Y\in\eufrak{g}^{\mathbb{R}_+}$ that does not commute with $X$ and let the unitary operator $U=e^{\pi(Y)}$ act on the algebra
\[
\mathcal{A}_X={\rm alg}\,\{X(f)|f\in\mathcal{S}(\mathbb{R}_+)\}
\]
generated by $X$. We can describe this action in two different ways.

First, we can let $U$ act on $\mathcal{A}_X$ directly. This gives another algebra
\[
\mathcal{A}_{X'}={\rm alg}\,\{UX(f)U^*|f\in\mathcal{S}(\mathbb{R}_+)\},
\]
generated by $X'(f)=UX(f)U^*$, $f\in\mathcal{S}(X)$. Since this algebra is again commutative, there exists a classical process $(\hat{X}'_t)_{t\ge 0}$ that has the same expectation values as $X'$ w.r.t.\ $\Omega$, i.e.
\[
\langle\Omega, g_1\big(X'(f_1)\big)\cdots g_n\big(X'(f_n)\big)\Omega\rangle = \mathbb{E}\Big(g_1\big(\hat{X}'(f_1)\big)\cdots g_n\big(\hat{X}'(f_n)\big)\Big)
\]
for all $n\in\mathbb{N}$, $g_1,\ldots,g_n\in \mathcal{C}_0(\mathbb{R})$, $f_1,\ldots,f_n\in\mathcal{S}(\mathbb{R}_+)$, where $\hat{X}'(f)=\int_{\mathbb{R}_+}f(t){\rm d}\hat{X}'_t$ for $f=\sum_{k=1}^n f_k\mathbf{1}_{[s_k,t_k[}\in\mathcal{S}(\mathbb{R}_+)$.

If $X'(f)$ is a function of $X(f)$, then $\mathcal{A}_X$ is invariant under the action of $U$. In this case the classical process $(\hat{X}'_t)_{t\ge 0}$ can be obtained from $(\hat{X}_t)_{t\ge 0}$ by a pathwise transformation, see \ref{QUASI-exa girsanov} and \ref{QUASI-exa gamma}. But even if this is not the case, we can still get a quasi-invariance formula that states that the law of $(\hat{X}'_t)_{t\ge 0}$ is absolutely continuous w.r.t.\ the law of $(\hat{X}_t)_{t\ge 0}$.

Second, we can let $U$ act on the vacuum state $\Omega$, this gives us a new state vector $\Omega'=U^*\Omega$. If $\Omega$ is cyclic for $\mathcal{A}_X$, then $\Omega'$ can be approximated by elements of the form $G\Omega$ with $G\in\mathcal{A}_X$. Often it is actually possible to find an element $G\in\mathcal{A}_X$ such that $G\Omega=\Omega'$. In this case the following calculation shows that the finite marginal distributions of $(\hat{X}'_t)_{t\ge 0}$ are absolutely continuous w.r.t.\ those of $(\hat{X}_t)$,
\begin{eqnarray*}
\mathbb{E}\Big(g\big(\hat{X}'(f)\big)\Big)
&=& \langle\Omega, g\big(X'(f)\big)\Omega\rangle = \langle\Omega, g\big(UX(f)U^*\big)\Omega\rangle \\
&=& \langle\Omega, Ug\big(X(f)\big)U^*\Omega\rangle = \langle U^*\Omega, g\big(X(f)\big)U^*\Omega\rangle \\
&=& \langle\Omega', g\big(X(f)\big)\Omega'\rangle = \langle G\Omega, g\big(X(f)\big)G\Omega\rangle \\
&=& \mathbb{E}\Big(g\big(X(f)\big)|\hat{G}|^2\Big).
\end{eqnarray*}
Here $G$ was a ``function'' of $X$ and $\hat{G}$ is obtained from $G$ by replacing $X$ by $\hat{X}$. This is possible, because $\mathcal{A}_X$ is commutative, and requires only standard functional calculus.

The density relating the law of $(\hat{X}'_t)_{t\ge 0}$ and that of $(\hat{X}_t)_{t\ge 0}$ is therefore given by $|\hat{G}|^2$.

The same calculation applies of course also to the finite joint distributions, i.e.\ we also have
\[
\mathbb{E}\Big(g_1\big(\hat{X}'(f_1)\big)\cdots g_n\big(\hat{X}'(f_n)\big)\Big)=\mathbb{E}\Big(g_1\big(\hat{X}(f_1)\big)\cdots g_n\big(\hat{X}(f_n)\big)|\hat{G}|^2\Big).
\]
for all $n\in\mathbb{N}$, $f_1,\ldots,f_n\in\mathcal{S}(\mathbb{R}_+)$, $g_1,\ldots,g_n\in \mathcal{C}_0(\mathbb{R})$.

\section{Examples}\label{QUASI-examples}

In this section we give the explicit calculation of the density $|G|^2$ for several examples.

We define our real Lie algebras as complex Lie algebras with an involution, because the relations can be given in a more convenient form for the complexifications. The real Lie algebra can be recovered as the real subspace of anti-hermitian elements.

\subsection{Gaussian and Poisson random variables}\label{QUASI-exa-osc}

The oscillator Lie algebra is the four dimensional Lie algebra $osc$ with basis $\{N,A^+,A^-,E\}$ and the Lie bracket given by
\[
{}
[N,A^\pm]=\pm A^\pm, \quad [A^-,A^+]=E, \quad [E,N]=[E,A^\pm]=0.
\]
We equip $osc$ furthermore with the involution $N^*=N$, $(A^+)^*=A^-$, and $E^*=E$. 

A general hermitian element of $osc$ can be written in the form $X_{\alpha,\zeta,\beta}=\alpha N+\zeta A^++\overline{\zeta}A^-+\beta E$ with $\alpha,\beta\in\mathbb{R}$, $\zeta\in\mathbb{C}$. 

Let $Y=i(wA^++\overline{w}A^-)$. We can compute the adjoint action of $g_t=e^{tY}$ on $X_{\alpha,\zeta,\beta}$ by solving the differential equation
\[
\dot{X}(t)=\frac{{\rm d}}{{\rm d}t}{\rm Ad}_{g_t}(X)=[Y,X(t)].
\]
Writing $X(t)=a(t)N+z(t)A^++\overline{z}(t)A^-+b(t)E$,
we get the system of ode's
\begin{eqnarray*}
\dot{a}(t) &=& 0, \\
\dot{z}(t) &=& -i\alpha w, \\
\dot{b}(t) &=& i (\overline{w}z-w\overline{z}),
\end{eqnarray*}
with initial conditions $a(0)=\alpha$, $z(0)=\zeta$, $b(0)=\beta$. We get 
\begin{equation}
X(t)=\alpha N+(\zeta-i\alpha w t)A^++(\overline{\zeta}+i\alpha \overline{w}t)A^-+\big(\beta +2t\Im (w\overline{\zeta})+\alpha|w|^2t^2\big)E.
\end{equation}

A representation $\rho$ of $osc$ on $\ell^2$ is defined by
\begin{eqnarray*}
\rho(N)|n\rangle &=& n |n\rangle, \\
\rho(A^+)|n\rangle &=& \sqrt{n+1}|n+1\rangle, \\
\rho(A^-)|n-1\rangle &=& \sqrt{n}|n-1\rangle, \\
\rho(E)|n\rangle &=& |n\rangle,
\end{eqnarray*}
where $|0\rangle,|1\rangle,\ldots$ is an orthonormal basis of $\ell^2$.

\begin{proposition}
The distribution of $\rho(X_{\alpha,\zeta,\beta})$ in the vacuum vector $|0\rangle$ is given by the characteristic function
\[
\langle 0|\exp\big(i\lambda \rho(X_{\alpha,\zeta,\beta})\big)|0\rangle = \left\{\begin{array}{lcl}
\exp\left(i\lambda \beta - \frac{\lambda^2}{2}|\zeta|^2\right) & \mbox{ for } & \alpha=0, \\
\exp\left(i\lambda\left(\beta-\frac{|\zeta|^2}{\alpha}\right) +\frac{|\zeta|^2}{\alpha^2}\left(e^{i\lambda\alpha} -1\right)\right) & \mbox{ for } & \alpha\not=0,
\end{array}\right.
\]
i.e.\ it is either a Gaussian random variable with variance $|\zeta|^2$ and mean $\beta$ or Poisson random variable with ``jump size'' $\alpha$, intensity $\frac{|\zeta|^2}{\alpha^2}$, and drift $\beta-\frac{|\zeta|^2}{\alpha}$.
\end{proposition}

\begin{lemma}
The vacuum vector $|0\rangle$ is cyclic for $\rho(X_{\alpha,\zeta,\beta})$, as long as $\zeta\not=0$, i.e.
\[
\overline{{\rm span\,}\{\rho(X_{\alpha,\zeta,\beta})^k|0\rangle:k=0,1,\ldots\}}=\ell^2.
\]
\end{lemma}
Therefore $v(t)=\exp\big(-t\rho(Y)\big)|0\rangle$ can be written in the form
\[
v(t)=\sum_{k=0}^\infty c_k(t) \rho(X_{\alpha,\zeta,\beta})^k|0\rangle = G(X_{\alpha,\zeta,\beta},t)|0\rangle.
\]
In order to compute the function $G$, we consider
\[
\frac{{\rm d}}{{\rm d}t} v(t) =-\exp\big(-t\rho(Y)\big)\rho(Y)|0\rangle=- i \exp\big(-t\rho(Y)\big)w|1\rangle.
\]
We can rewrite this as
\begin{eqnarray*}
\lefteqn{ 
\frac{{\rm d}}{{\rm d}t} v(t) = - 
 \frac{iw}{\tilde{\zeta}(t)}
 \exp\big(-t\rho(Y)\big)\left(X_{\alpha,\tilde{\zeta}(t),\tilde{\beta}(t)}|0\rangle-\tilde{\beta}(t)|0\rangle\right)  
} 
\\
&=& - 
 \frac{iw}{\tilde{\zeta}(t)}
 \left(X_{\alpha,\zeta,\beta}-\tilde{\beta}(t)\right)\exp\big(-t\rho(Y)\big)|0\rangle =- 
 \frac{iw}{\tilde{\zeta}(t)} 
 \left(X_{\alpha,\zeta,\beta}-\tilde{\beta}(t)\right) v(t),
\end{eqnarray*} 
 where 
\begin{eqnarray*}
\tilde{\zeta}(t) &=& \zeta - i\alpha w t, \\
\tilde{\beta}(t) &=& \beta + 2t\Im (w\overline{\zeta})+ \alpha|w|^2 t^2.
\end{eqnarray*}
This is satisfied provided $G(x,t)$ satisfies the differential equation
\[
\frac{{\rm d}}{{\rm d}t} G(x,t)=- 
 \frac{iw}{\tilde{\zeta}(t)}
 \left(x-\tilde{\beta}(t)\right)G(x,t)
\]
with initial condition $G(x,0)=1$. We get 
\[
G(x,t)=\exp \left( 
 -i w \int_0^t\left(\frac{x-\tilde{\beta}(s)}{\tilde{\zeta}(s)}\right){\rm d}s 
 \right).
\]
Evaluating the integral, this can be written as 
\[
G(x,t)= 
\left(1-i\frac{\alpha wt}{\zeta}\right)^{\frac{x-\beta}{\alpha} + \frac{\vert \zeta \vert^2}{\alpha^2} 
 } \exp\left(   i\frac{w\overline{\zeta}}{\alpha} 
 - \frac{t^2}{2} \vert w \vert^2 
 \right), 
\]
and therefore
\[
|G(x,t)|^2  = 
\left(1+2 t \alpha \Im\left(\frac{w}{\zeta}\right) 
 + t^2 \alpha^2 \frac{\vert w\vert^2}{\vert \zeta\vert^2} 
 \right)^{\frac{x-\beta}{\alpha} + \frac{\vert \zeta \vert^2}{\alpha^2} 
 } 
 e^{- \frac{\vert \zeta\vert^2}{\alpha^2} 
 \left( 
 2t\alpha \Im \frac{w}{\zeta} 
 + t^2 \alpha^2 \frac{\vert w \vert^2}{\vert \zeta\vert^2} 
 \right)  
 } 
. 
\] 
After letting $\alpha$ go to $0$ we get 
\[ 
 \vert G(x,t)\vert^2 = 
 e^{2 t 
 (x-\beta ) 
 \Im \frac{w}{\zeta} 
 - \frac{t^2}{2} \vert \zeta\vert^2 \left( 
 2 \Im \frac{w}{\zeta} 
 \right)^2 
} 
. 
\]
Note that the classical analog of this limiting procedure is 
\[
(1+\alpha )^{\alpha \lambda (N_\alpha -\lambda /\alpha^2) + \frac{\lambda^2}{\alpha}} 
 \stackrel{\alpha \to 0}{\longrightarrow}
 e^{\lambda X-\frac{1}{2}\lambda^2}, 
\] 
where $N_\alpha$ is a Poisson random variable with intensity $\lambda>0$ and $\lambda (N_\alpha -\lambda /\alpha^2)$ converges in distribution to a standard Gaussian variable $X$. No such normalization is needed in the quantum case.
\begin{proposition}
We have 
\[
\mathbb{E}\Big[ 
 g\big(X(t)\big)\Big] 
 =\mathbb{E}\left[ 
 g(X_{\alpha,\zeta,\beta})\big|G(X_{\alpha,\zeta,\beta},t)\big|^2 \right] 
\]
for all $g\in \mathcal{C}_0(\mathbb{R})$.
\end{proposition} 
For $\alpha=0$, this identity gives the relative density of two Gaussian random variables with the same variance, but different means. For $\alpha\not=0$, it gives the relative density of two Poisson random variables with different intensities.

\subsection{Brownian motion and the Girsanov formula}\label{QUASI-exa girsanov}

Let now $(j_{st})_{0\le s\le t}$ be the L\'evy process on $osc$ with the Sch\"urmann triple defined by $D=\mathbb{C}$,
\begin{gather*}
\rho(N)=1, \quad \rho(A^\pm)=\rho(E)=0, \\
\eta(A^+)=1, \quad \eta(N)=\eta(A^-)=\eta(E)=0, \\
L(N)=L(A^\pm)=0, \quad L(E)=1.
\end{gather*}
Taking for $X$ the constant function with value $-i(A^++A^-)$, we get
\[
X(f)=A^*(f)+A(f)
\]
and the associated classical process $(\hat{X}_t)_{t\ge 0}$ is Brownian motion.

We choose for $Y=h(A^+-A^-)$, with $h\in\mathcal{S}(\mathbb{R}_+)$. A similar calculation as in the previous subsection yields
\[
X'(\mathbf{1}_{[0,t]})= e^YX(\mathbf{1}_{[0,t]})e^{-Y}=X(\mathbf{1}_{[0,t]}) - 2\int_0^t h(s){\rm d}s
\]
i.e.\ $\mathcal{A}_X$ is invariant under $e^Y$ and $(\hat{X}'_t)_{t\ge 0}$ is obtained from $(\hat{X}_t)_{t\ge 0}$ by adding a drift.

$e^{\pi(Y)}$ is a Weyl operator and gives an exponential vector, if it acts on the vacuum, i.e.
\[
e^{\pi(Y)}\Omega = e^{-||h||^2/2}\mathcal{E}(h)
\]
see, e.g., \cite{parthasarathy92,meyer95}. But - up to the normalization - we can create the same exponential vector also by acting on $\Omega$ with $e^{X(h)}$,
\[
e^{X(h)}\Omega = e^{||h||^2/2}\mathcal{E}(h).
\]
Therefore we get $G=\exp\big(X(h)-||h||^2\big)$ and the well-known Girsanov formula
\[
\mathbb{E}\Big(g\big(\hat{X}'(f)\big)\Big) = \mathbb{E}\left(g\big(\hat{X}(f)\big)\exp\left(2X(h)-2\int_0^\infty h^2(s){\rm d}s\right)\right).
\]

\subsection{Poisson process and the Girsanov formula}\label{QUASI-poisson girsanov}

Taking for $X$ the constant function with value $-i(N+\nu A^++\nu A^- + \nu^2 E)$, we get
\[
X(f)=N(f)+\nu A^*(f)+ \nu A(f) + \nu^2 \int_0^\infty f(s) ds 
\]
 and the associated classical process $(\hat{X}_t)_{t\ge 0}$ is a non-compensated Poisson 
 process with intensity $\nu^2$ and jump size $1$. 
 Given $h\in {\mathcal S}(\mathbb{R}_+)$ of the form 
$$h(t) = \sum_{k=1}^n h_k 1_{[s_k,t_k[} (t),$$ 
 with $h_k> -\nu^2$, let 
$$w(t) = i(\sqrt{\nu^2+h(t)}-\nu), 
$$ 
 and $Y=w(A^+-A^-)$. 
 The calculations of Subsection \ref{QUASI-exa-osc} show that 
\[
X'(\mathbf{1}_{[0,t]})= e^YX(\mathbf{1}_{[0,t]})e^{-Y} 
\]
 is a non-compensated Poisson process with intensity $\nu^2+h(t)$. 
 We have the Girsanov formula
\begin{eqnarray*} 
\mathbb{E}\Big(g\big(\hat{X}'(f)\big)\Big) 
 & = & 
 \mathbb{E}\left(g\big(\hat{X}(f)\big) 
 \prod_{k=1}^{n} 
 \left(1+\frac{h_k}{\nu^2} \right)^{\hat{X}(1_{[s_k,t_k[})} 
 e^{- \nu^2 (t_k-s_k) h_k } 
 \right) 
\\ 
& = & 
 \mathbb{E}\left(g\big(\hat{X}(f)\big) 
 e^{- \nu^2 \int_0^\infty h(s) ds } 
 \prod_{k=1}^{n} 
 \left(1+\frac{h_k}{\nu^2} \right)^{\hat{X}(1_{[s_k,t_k[})} 
 \right) \\
&=&
 \mathbb{E}\left(g\big(\hat{X}(f)\big) 
 \exp\left(\hat{X}\left(\log\left(1+\frac{h}{\nu^2} \right)\right)
 -\nu^2\int_0^\infty h(s){\rm d}s\right)\right)
. 
\end{eqnarray*}

\subsection{$sl(2,\mathbb{R})$ and the Meixner, Gamma, and Pascal random variables}\label{QUASI-exa-sl2}

Let us now consider the three-dimensional Lie algebra $sl(2;\mathbb{R})$, with basis $B^+,B^-,M$, Lie bracket
\[
[M,B^\pm]=\pm 2 B^\pm, \qquad [B^-,B^+]=M,
\]
and the involution $(B^+)^*=B^-$, $M^*=M$.

For $\beta\in\mathbb{R}$, we set $X_\beta=B^++B^-+\beta M$. Furthermore, we choose $Y=B^--B^+$. We compute $[Y,X_\beta]=2\beta B^++2\beta B^-+2M=2\beta X_{1/\beta}$ and
\[
e^{tY/2}X_\beta e^{-tY/2}=e^{\frac{t}{2}{\rm ad}Y}X_\beta = \big(\cosh(t)+\beta\sinh(t)\big)X_{\gamma(\beta,t)},
\]
where
\[
\gamma(\beta,t)=\frac{\beta\cosh(t)+\sinh(t)}{\cosh(t) + \beta \sinh(t)}.
\]

For $c>0$ we can define representations of $sl(2;\mathbb{R})$ on $\ell^2$ by
\begin{eqnarray*}
\rho_c(B^+)|k\rangle &=& \sqrt{(k+c)(k+1)}|k+1\rangle, \\
\rho_c(M)|k\rangle &=& (2k+c)|k\rangle, \\
\rho_c(B^-)|k\rangle &=& \sqrt{k(k+c-1)}|k\rangle,
\end{eqnarray*}
where $|0\rangle,|1\rangle,\ldots$ is an orthonormal basis of $\ell^2$.

Using the Splitting Lemma for $sl(2)$, cf.\ \cite[Chapter 1, 4.3.10]{feinsilver+schott93} to write $e^{\lambda X_\beta}$ as a product $e^{\nu_+B^+}e^{\nu_0M}e^{\nu_-B^-}$, it is straight-forward to compute the distribution of $\rho_c(X_\beta)$ in the state vector $|0\rangle$.
\begin{proposition}
The Fourier-Laplace transform of the distribution of $\rho_c(X_\beta)$ w.r.t.\ $|0\rangle$ is given by
\[
\langle 0|\exp\big(\lambda \rho_c(X_\beta)\big)|0\rangle = \left(\frac{\sqrt{\beta^2-1}}{\sqrt{\beta^2-1}\cosh\big(\lambda\sqrt{\beta^2-1}\big)-\beta\sinh\big(\lambda\sqrt{\beta^2-1}\big)}\right)^c.
\]
\end{proposition}
\begin{remark}
For $|\beta|<1$, this distribution is called the Meixner distribution. It is absolutely continuous w.r.t.\ to the Lebesgue measure and the density is given by
\[
 C \exp\left(\frac{(\pi-2\arccos \beta)x}{2\sqrt{1-\beta^2}}\right)\left|\Gamma\left(\frac{c}{2}+\frac{ix}{2\sqrt{1-\beta^2}}\right)\right|^2{\rm d}x,
\]
see also \cite{accardi+franz+skeide02}. $C$ is a normalization constant.

For $\beta=\pm 1$, we get the Gamma distribution, which has the density
\[
\frac{|x|^{c-1}}{\Gamma(c)} e^{-\beta x}\mathbf{1}_{\beta\mathbb{R}_+}.
\]

Finally, for $|\beta|>1$, we get a discrete measure, the so-called Pascal distribution.
\end{remark}

\begin{lemma}
The lowest weight vector $|0\rangle$ is cyclic for $\rho_c(X_\beta)$ for all $\beta\in\mathbb{R}$, $c>0$.
\end{lemma}

Therefore we can write $v_t=e^{-t\rho_c(Y)/2}|0\rangle$ in the form
\[
v(t)=\sum_{k=0}^\infty c_k(t) \rho_c(X_{\beta})^k|0\rangle = G(X_{\beta},t)|0\rangle.
\]
In order to compute the function $G$, we consider
\[
\frac{{\rm d}}{{\rm d}t} v(t) =-\frac{1}{2}\exp\big(-t\rho_c(Y)/2\big)\rho_c(Y)|0\rangle= \frac{1}{2} \exp\big(-t\rho_c(Y)\big)\sqrt{c}|1\rangle.
\]
As in Subsection \ref{QUASI-exa-osc}, we introduce $X_\beta$ into this equation to get an ordinary differential equation for $G$,
\begin{eqnarray*}
\frac{{\rm d}}{{\rm d}t} v(t) &=&
\frac{1}{2} \exp\big(-t\rho_c(Y)/2\big)\big(\rho_c(X_{\gamma(\beta,t)})-c\gamma(\beta,t)\big)|0\rangle \\
&=& \frac{1}{2}\frac{\rho_c(X_\beta)-c\big(\beta\cosh(t)+\sinh(t)\big)}{\cosh(t)+\beta\sinh(t)}\exp\big(-t\rho_c(Y)/2\big)|0\rangle
\end{eqnarray*}
This is satisfied, if
\[
\frac{{\rm d}}{{\rm d}t}G(x,t)=\frac{1}{2}\frac{x-c\big(\beta\cosh(t)+\sinh(t)\big)}{\cosh(t)+\beta\sinh(t)}G(x,t)
\]
with initial condition $G(x,0)=1$. The solution of this ode is given by
\[
G(x,t)=\exp \frac{1}{2}\int_0^t \frac{x-c\big(\beta\cosh(s)+\sinh(s)\big)}{\cosh(s)+\beta\sinh(s)}{\rm d}s.
\]
For $\beta=1$, we get $G(x,t)=\exp\frac{1}{2}\big(x(1-e^{-t})-ct\big)$ and we recover the following identity for a Gamma distributed random variable $Z$ with parameter $c$,
\[
\mathbb{E}\big(g(e^tZ)\big)=\mathbb{E}\Big(g(Z)\exp\big(Z(1-e^{-t})-ct\big)\Big).
\]
If $|\beta|<1$, then we can write $G$ in the form
\[
G(x,t)=\exp\big(\Phi(\beta,t)x-c\Psi(\beta,t)\big),
\]
where
\begin{eqnarray}
\Phi(\beta,t) &=& \frac{1}{\sqrt{1-\beta^2}}\left(\arctan\left(e^t \sqrt{\frac{1+\beta}{1-\beta}}\right) - \arctan\left(\sqrt{\frac{1+\beta}{1-\beta}}\right)\right), \label{QUASI-phi} \\
\Psi(\beta,t) &=& \frac{1}{2}\Big(t+\log \big(1+\beta+e^{-2t}(1-\beta)\big)-\log2\Big). \label{QUASI-psi}
\end{eqnarray}

\subsection{A quasi-invariance formula for the Gamma process}\label{QUASI-exa gamma}

Let now $(j_{st})_{0\le s\le t}$ be the L\'evy process on $sl(2;\mathbb{R})$ with Sch\"urmann triple $D=\ell^2$, $\rho=\rho_2$, and
\begin{gather*}
\eta(B^+)=|0\rangle,\quad \eta(B^-)=\eta(M)=0, \\
L(M)=1,\quad L(B^\pm)=0,
\end{gather*}
cf.\ \cite[Example 3.1]{accardi+franz+skeide02}. We take for $X$ the constant function with value $-i(B^++B^++M)$, then the random variables
\[
X(\mathbf{1}_{[s,t]})=\Lambda_{st}\big(\rho(X)\big)+A^*_{st}\big(|0\rangle\big)+A_{st}(|0\rangle)+(t-s){\rm id}
\]
are Gamma distributed in the vacuum vector $\Omega$.

For $Y$ we choose $h(B^--B^+)$ with $h=\sum_{k=1}^rh_k\mathbf{1}_{[s_k,t_k[}\in\mathcal{S}(\mathbb{R}_+)$, $0\le s_1\le t_1\le s_2\le \cdots\le t_n$. As in the previous subsection, we get
\[
X'(\mathbf{1}_{[s,t]})=e^{\pi(Y)}X(\mathbf{1}_{[s,t]})e^{-\pi(Y)}=X(e^{2h}\mathbf{1}_{[s,t]}).
\]
On the other hand, using the tensor product structure of the Fock space, we can calculate
\begin{eqnarray*}
e^{-\pi(Y)}\Omega &=& e^{-\sum_{k=1}^nh_kj_{s_kt_k}(Y)}\Omega \\
&=& e^{-h_1j_{s_1t_1}(Y))}\Omega \otimes \cdots \otimes e^{-h_nj_{s_nt_n}(Y))}\Omega \\
&=& \exp \frac{1}{2}\Big( X\big((1-e^{-2h_1})\mathbf{1}_{[s_1,t_1[}\big) - (t_1-s_1) 2h {\rm d}s\Big)\Omega\otimes\cdots \\
&& \cdots \otimes \exp \frac{1}{2}\Big( X\big((1-e^{-2h_n})\mathbf{1}_{[s_n,t_n[\big)} - (t_n-s_n) 2h {\rm d}s\Big)\Omega \\
&=& \exp \frac{1}{2}\left( X(1-e^{-2h}) - \int_{\mathbb{R}_+} 2h {\rm d}s\right)\Omega,
\end{eqnarray*}
since $j_{st}$ is equivalent to $\rho_{t-s}$.

\begin{proposition}
Let $n\in\mathbb{N}$, $f_1,\ldots,f_n\in\mathcal{S}(\mathbb{R}_+)$, $g_1,\ldots,g_n\in \mathcal{C}_0(\mathbb{R})$, then we have
\begin{gather*}
\mathbb{E}\Big(g_1\big(X'(f_1)\big)\cdots g_n\big(X'(f_n)\big)\Big) \\
= \mathbb{E}\left(g_1\big(X(f_1)\big)\cdots g_n\big(X(f_n)\big)\exp \left( X(1-e^{-2h}) - \int_{\mathbb{R}_+} 2h {\rm d}s\right)\right).
\end{gather*}
\end{proposition}

\subsection{A quasi-invariance formula for the Meixner process}\label{QUASI-exa-meixner}

We consider again the same L\'evy process on $sl(2;\mathbb{R})$ as in the previous subsection. Let $\varphi,\beta\in\mathcal{S}(\mathbb{R}_+)$ with $|\beta(t)|<1$ for all $t\in\mathbb{R}_+$, and set
\[
X_{\varphi,\beta}=\varphi(B^++B^-+\beta M) \in sl(2;\mathbb{R})^{\mathbb{R}_+}.
\]
Let $Y$ again be given by $Y=h(B^--B^+)$, $h\in\mathcal{S}(\mathbb{R}_+)$. Then we get
\[
X'(t)=e^{Y(t)}X(t)e^{-Y(t)} = \varphi(t)\big(\cosh(2h)+\beta(t)\sinh(2h)\big)\Big(B^++B^-+\gamma\big(\beta(t),2h\big)M\Big)
\]
i.e.\ $X'=X_{\varphi',\beta'}$ with
\begin{eqnarray*}
\varphi'(t) &=& \varphi(t)\Big(\cosh\big(2h(t)\big)+\beta(t)\sinh\big(2h(t)\big)\Big), \\
\beta'(t) &=&\gamma\big(\beta(t),2h(t)\big).
\end{eqnarray*}

As in the previous subsection, we can also calculate the function $G$,
\[
e^{\pi(Y)}\Omega = \exp \frac{1}{2} \left(X_{\Phi(\beta,2h),\beta}-\int_{\mathbb{R}_+} \Psi\big(\beta(t),2h(t)\big)\right)\Omega,
\]
where $\Phi,\Psi$ are defined as in Equations \eqref{QUASI-phi} and \eqref{QUASI-psi}.
\begin{proposition}
The finite joint distributions of $X_{\varphi',\beta'}$ are absolutely continuous w.r.t.\ to those of $X_{\varphi,\beta}$, and the mutual density is given by
\[
\exp \left(X_{\Phi(\beta,2h),\beta}-\int_{\mathbb{R}_+} \Psi\big(\beta(t),2h(t)\big)\right).
\]
\end{proposition}

\section{Conclusion}

The results stated in Subsections \ref{QUASI-exa girsanov}, \ref{QUASI-poisson girsanov}, \ref{QUASI-exa gamma}, and \ref{QUASI-exa-meixner} have been proved with our methods for the finite joint distributions. They can be extended to the distribution of the processes either using continuity arguments for the states and endomorphisms on our operator algebras or by the use of standard tightness arguments of classical probability.

Of course the general idea also applies to classical processes obtained by a different choice of the commutative subalgebra, e.g.\ as in \cite{biane98b}, and to more general classes of quantum stochastic processes. The restrictions in this Note were motivated by the fact that they simplify the explicit calculations.

%% file: dilations.tex
\chapter{L\'evy Processes and Dilations of Completely Positive semigroups}\label{chapter-II}\markboth{L\'EVY PROCESSES AND DILATIONS}{L\'EVY PROCESSES AND DILATIONS}

In this chapter we will show how L\'evy process can be used to construct dilations of quantum dynamical semigroups on the matrix algebra $\mathcal{M}_d$. That unitary cocycles on the symmetric Fock space tensor a finite-dimensional initial space can be interpreted as a L\'evy process on a certain involutive bialgebra, was first observed in \cite{schuermann90b}. For more details on quantum dynamical semigroups and their dilations, see \cite{bhat01,bhat03} and the references therein.

\section{The non-commutative analogue of the algebra of coefficients of the unitary group}

For $d\in\mathbb{N}$ we denote by $\mathcal{U}_d$ the free non-commutative (!) $*$-algebra generated by indeterminates $u_{ij}$, $u^*_{ij}$, $i,j=1,\ldots,d$ with the relations
\begin{eqnarray*}
\sum_{j=1}^d u_{kj}u^*_{\ell j} &=& \delta_{k\ell}, \\
\sum_{j=1}^d u^*_{jk}u_{j\ell} &=& \delta_{k\ell}.
\end{eqnarray*}
The $*$-algebra $\mathcal{U}_d$ is turned into a $*$-bialgebra, if we put
\begin{eqnarray*}
\Delta(u_{k\ell}) &=& \sum_{j=1}^d u_{kj}\otimes u_{j\ell}, \\
\varepsilon(u_{k\ell}) &=& \delta_{k\ell}.
\end{eqnarray*}
This $*$-bialgebra has been investigated by Glockner and von Waldenfels, see \cite{glockner+waldenfels89}. If we assume that the generators $u_{ij}$, $u^*_{ij}$ commute, we obtain the coefficient algebra of the the unitary group $U(d)$. This is why $\mathcal{U}_d$ is often called {\em the non-commutative analogue of the algebra of coefficients of the unitary group}. It is equal to the $*$-algebra generated by the mappings
\[
\xi_{k\ell}:\mathcal{U}(\mathbb{C}^d\otimes H)\to \mathcal{B}(H)
\]
with
\[
\xi_{k\ell}(U) = P_kUP_\ell^* = U_{k\ell}
\]
for $U\in\mathcal{U}(\mathbb{C}^d\otimes H)\subseteq \mathcal{M}_d\big(\mathcal{B}(H)\big)$, where $H$ is an infinite-dimensional, separable Hilbert space and $\mathcal{U}(\mathbb{C}^d\otimes H)$ denotes the unitary group of operators on $\mathbb{C}^d\otimes H$. Moreover $\mathcal{B}(H)$ denotes the $*$-algebra of bounded operators, $\mathcal{M}_d\big(\mathcal{B}(H)\big)$ the $*$-algebra of $d\times d$-matrices with elements from $\mathcal{B}(H)$ and $P_k:\mathbb{C}^d\otimes H \to H$ the projection on the $k$-th component.

\begin{proposition}
\begin{enumerate}
\item
On $\mathcal{U}_1$ a faithful Haar measure is given by $\lambda(u^n)=\delta_{0,n}$, $n\in\mathbb{Z}$.
\item
On $\mathcal{U}_1$ an antipode is given by setting $S(x)=x^*$ and extending $S$ as a $*$-algebra homomorphism.
\item
For $d>1$ the bialgebra $\mathcal{U}_d$ does not possess an antipode.
\end{enumerate}
\end{proposition}
\begin{proof}
A Haar measure on a bialgebra $\mathcal{B}$ is a normalized linear functional $\lambda$ satisfying
\[
\lambda\star\varphi=\varphi(1)\lambda=\varphi\star\lambda
\]
for all linear functionals $\varphi$ on $\mathcal{B}$.

The verification of properties (1) and (2) is straightforward, we shall only prove (3). We suppose that an antipode exists. Then
\begin{eqnarray*}
u_{\ell k}^* &=& \sum_{n=1}^d\sum_{j=1}^d S(u_{kj})u_{j n}u^*_{\ell n} \\
&=& \sum_{j=1}^d S(u_{kj}) \sum_{n=1}^d u_{j n}u^*_{\ell n} \\
&=& \sum_{j=1}^d S(u_{kj})\delta_{j\ell} = S(u_{k\ell}).
\end{eqnarray*}
Similarly, one proves that $S(u^*_{k\ell})=u_{lk}$. Since $S$ is an antipode, it has to be an algebra anti-homomorphisms. Therefore,
\[
S\left(\sum_{j=1}^d u_{kj}u^*_{\ell j}\right) = \sum_{j=1}^d S(u^*_{\ell j})S(u_{kj}) = \sum_{j=1}^d u_{j\ell}u^*_{jk},
\]
which is not equal to $\delta_{k\ell}$, if $d>1$.
\end{proof}

\begin{remark}
Since $\mathcal{U}_d$ does not have an antipode for $d>1$, it is not a compact quantum group (for $d=1$, of course, its $C^*$-completion is the quantum group of continuous functions on the circle $S^1$). We do not know, if $\mathcal{U}_d$ has a (faithful) Haar measure for $d>1$.
\end{remark}

We have $\mathcal{U}_n=\mathbb{C}1\oplus\mathcal{U}^0_n$, where $\mathcal{U}^0_n=K_1=\text{ker }\varepsilon$ is the ideal generated by $\hat{u}_{ij}=u_{ij}-\delta_{ij}1$, $i,j=1,\ldots n$, and their adjoints. The defining relations become
\begin{equation}\label{u-hat rel}
 - \sum_{j=1}^d \hat{u}_{ij} \hat{u}^*_{kj} = \hat{u}_{ik} + \hat{u}^*_{ki} = -\sum_{j=1}^d \hat{u}^*_{ji} \hat{u}_{jk},
\end{equation}
for $i,k=1,\ldots,n$, in terms of these generators.

\section{An example of a L\'evy process on $\mathcal{U}_d$}\label{II-first example}

A one-dimensional representation $\sigma:\mathcal{U}_d\to\mathbb{C}$ is determined by the matrix $w=(w_{ij})_{1\le i,j\le d}\in\mathcal{M}_d$, $w_{ij}=\sigma(u_{ij})$. The relations in $\mathcal{U}_d$ imply that $w$ is unitary. For $\ell=(\ell_{ij})\in\mathcal{M}_d$ we can define a $\sigma$-cocycle (or a $\sigma$-derivation) as follows. We set
\begin{eqnarray*}
\eta_\ell(u_{ij}) &=& \ell_{ij}, \\
\eta_\ell(u^*_{ij}) &=& - (w^*\ell)_{ji} = -\sum_{k=1}^d \overline{w}_{kj}\ell_{ki},
\end{eqnarray*}
on the generators and require $\eta_\ell$ to satisfy
\[
\eta_\ell(ab)=\sigma(a)\eta_\ell(b)+\eta_\ell(a)\varepsilon(b)
\]
for $a,b\in\mathcal{U}_d$. The hermitian linear functional $L_{w,\ell}:\mathcal{U}_d\to\mathbb{C}$ with
\begin{eqnarray*}
L_{w,\ell}(\mathbf{1}) &=& 0, \\
L_{w,\ell}(u_{ij}) &=& \overline{L_{w,\ell}(u^*_{ij})} = -\frac{1}{2} (\ell^* \ell)_{ij} = -\frac{1}{2}\sum_{k=1}^d \overline{\ell}_{ki}\ell_{kj}, \\
L_{w,\ell}(ab) &=& \varepsilon(a)L_{w,\ell}(b) + \overline{\eta_\ell(a^*)}\eta_\ell(b) + L_{w,\ell}(a)\varepsilon(b)
\end{eqnarray*}
for $a,b\in\mathcal{U}_d$, can be shown to be a generator with Sch\"urmann triple $(\sigma,\eta_\ell,L_{w,\ell})$. The associated L\'evy process on $\mathcal{U}_d$ is determined by the quantum stochastic differential equations
\[
{\rm d}j_{st}(u_{ij}) = \sum_{k=1}^d j_{st}(u_{ik})\left(\ell_{kj}{\rm d}A^*_t + (w_{kj}-\delta_{kj}){\rm d}\Lambda_t - \sum_{n=1}^d w_{nj}\overline{\ell}_{nk}{\rm d}A_t -\frac{1}{2}\sum_{n=1}^d \overline{\ell}_{nk}\ell_{nj}{\rm d}t\right),
\]
on $\Gamma\big(L^2(\mathbb{R}_+,\mathbb{C})\big)$ with initial conditions $j_{ss}(u_{ij})=\delta_{ij}$.

We define an operator process $(U_{st})_{0\le s\le t}$ in $\mathcal{M}_d\otimes\mathcal{B}\Big(\Gamma\big(L^2(\mathbb{R}_+,\mathbb{C})\big)\Big)\cong \mathcal{B}\Big(\mathbb{C}^d\otimes \Gamma\big(L^2(\mathbb{R}_+,\mathbb{C})\big)\Big)$ by
\[
U_{st} = \big(j_{st}(u_{ij})\big)_{1\le i,j\le d},
\]
for $0\le s\le t$. Then $(U_{st})_{0\le s\le t}$ is a unitary operator process and satisfies the quantum stochastic differential equation
\[
{\rm d}U_{st} = U_{st} \left(\ell{\rm d}A_t^* + (w-\mathbf{1}){\rm d}\Lambda_t - \ell^*w{\rm d}A_t - \frac{1}{2}\ell^*\ell{\rm d}t\right)
\]
with initial condition $U_{ss}=\mathbf{1}$. The increment property of $(j_{st})_{0\le s\le t}$ implies that $(U_{st})_{0\le s\le t}$ satisfies
\begin{equation}\label{II-U-inc}
U_{0s}U_{s,s+t} = U_{0,s+t}
\end{equation}
for all $0\le s\le t$.

Let $S_t:L^2(\mathbb{R}_+,\mathcal{K})\to L^2(\mathbb{R}_+,\mathcal{K})$ be the shift operator,
\[
S_tf(s) = \left\{
\begin{array}{lcl}
f(s-t) & & \mbox{ if } s\ge t, \\
0 & & \mbox{ else,}
\end{array}
\right.
\]
for $f\in L^2(\mathbb{R}_+,\mathcal{K})$, and define $W_t:\Gamma\big(L^2(\mathbb{R}_+,\mathcal{K})\big)\otimes \Gamma\big(L^2([0,t[,\mathcal{K})\big)\to \Gamma\big(L^2(\mathbb{R}_+,\mathcal{K})\big)$ by
\[
W_t\big(\mathcal{E}(f)\otimes\mathcal{E}(g)\big) = \mathcal{E}(g+S_tf),
\]
on exponential vectors $\mathcal{E}(f)$, $\mathcal{E}(g)$ of functions $f\in L^2(\mathbb{R}_+,\mathcal{K})$, $g\in L^2([0,t[,\mathcal{K})$. Then the CCR flow $\gamma_t:\mathcal{B}\Big(\Gamma\big(L^2(\mathbb{R}_+,\mathcal{K})\big)\Big)$ is defined by
\[
\gamma_t(Z)=W_t(Z\otimes \mathbf{1})W_t^*,
\]
for $Z\in \mathcal{B}\Big(\Gamma\big(L^2(\mathbb{R}_+,\mathcal{K})\big)\Big)$. On $\mathcal{B}\Big(\mathbb{C}^d\otimes \Gamma\big(L^2(\mathbb{R}_+,\mathcal{K})\big)\Big)$ we have the $E_0$-semigroup $(\tilde{\gamma}_t)_{t\ge 0}$ with $\tilde{\gamma}_t={\rm id}\otimes\gamma_t$.

We have $U_{s,s+t}=\tilde{\gamma}_s(U_{0t})$ for all $s,t\ge 0$ and therefore increment property \eqref{II-U-inc} implies that $(U_t)_{t\ge 0}$ with $U_t=U_{0t}$, $t\ge 0$, is a left cocycle of $(\tilde{\gamma}_t)_{t\ge 0}$, i.e.\
\[
U_{s+t} = U_s \tilde{\gamma}_s(U_t),
\]
for all $s,t\ge 0$. One can check that $(U_t)_{t\ge 0}$ is also local and continuous.

Therefore we can define a new $E_0$-semigroup $(\eta_t)_{t\ge 0}$ on $\mathcal{B}\Big(\mathbb{C}^d\otimes \Gamma\big(L^2(\mathbb{R}_+,\mathcal{K})\big)\Big)$ by
\begin{equation}\label{II-e0sg-exa}
\eta_t(Z) = U_t\tilde{\gamma}_t(Z)U_t^*,
\end{equation}
for $Z\in \mathcal{B}\Big(\mathbb{C}^d\otimes \Gamma\big(L^2(\mathbb{R}_+,\mathcal{K})\big)\Big)$ and $t\ge 0$.

Let $\{e_1,\ldots,e_d\}$ be the standard basis of $\mathbb{C}^d$ and denote by $\mathbb{E}_0$ the conditional expectation from $\mathcal{B}\Big(\mathbb{C}^d\otimes \Gamma\big(L^2(\mathbb{R}_+,\mathcal{K})\big)\Big)$ to $\mathcal{B}(\mathbb{C}^d)\cong\mathcal{M}_d$ determined by
\[
\big(\mathbb{E}_0(Z)\big)_{ij} = \langle e_i\otimes \Omega, Z e_j\otimes\Omega\rangle
\]
for $Z\in \mathcal{B}\Big(\mathbb{C}^d\otimes \Gamma\big(L^2(\mathbb{R}_+,\mathcal{K})\big)\Big)$. Then
\begin{equation}\label{II-qdsg-exa}
\tau_t=\mathbb{E}_0\big(\eta_t(X\otimes \mathbf{1})\big)
\end{equation}
defines a quantum dynamical semigroup on $\mathcal{M}_d$. It acts on the matrix units $E_{ij}$ by
\begin{eqnarray*}
\tau_t(E_{ij}) & = &
\left(\begin{array}{cccc}
\langle e_1\otimes\Omega,U_t(E_{ij}\otimes\mathbf{1})U_t^*e_1\otimes\Omega\rangle & \cdots & \langle e_1\otimes\Omega,U_t(E_{ij}\otimes\mathbf{1})U_t^*e_d\otimes\Omega\rangle \\
\vdots & & \vdots \\
\langle e_d\otimes\Omega,U_t(E_{ij}\otimes\mathbf{1})U_t^*e_1\otimes\Omega\rangle & \cdots & \langle e_d\otimes\Omega,U_t(E_{ij}\otimes\mathbf{1})U_t^*e_d\otimes\Omega\rangle
\end{array}\right)
\\[1mm]
& = &
\varphi_t\left(\begin{array}{cccc}
u_{1i}u^*_{1j} & u_{1i}u^*_{2j} & \cdots & u_{1i}u^*_{dj} \\
u_{2i}u^*_{1j} & u_{2i}u^*_{2j} & \cdots&  u_{2i}u^*_{dj} \\
\vdots & \vdots & & \vdots \\
u_{di}u^*_{1j} & u_{di}u^*_{2j} & \cdots & u_{di}u^*_{dj} 
\end{array}\right),
\end{eqnarray*}
and therefore the generator $\mathcal{L}$ of $(\tau_t)_{t\ge 0}$ is given by
\[
\mathcal{L}(E_{ij}) = \big(L_{w,\ell}(u_{ki}u^*_{mj})\big)_{1\le k,m\le d},
\]
for $1\le i,j\le d$.

\begin{lemma}
The generator $\mathcal{L}$ of $(\tau_t)_{t\ge 0}$ is given by
\[
\mathcal{L}(X)=\ell^*wXw^*\ell -\frac{1}{2}\big\{X,\ell^*\ell\}
\]
for $X\in\mathcal{M}_d$.
\end{lemma}
\begin{proof}
We have, of course, $\left.\frac{{\rm d}}{{\rm d}t}\right|_{t=0}\varphi_t(u_{ki}u^*_{mj})=L_{w,\ell}(u_{ki}u^*_{mj})$. Using (\ref{I-L-coboundary}) of Chapter \ref{levy-intro} and the definition of the Sch\"urmann triple, we get
\begin{eqnarray*}
L_{w,\ell}(u_{ki}u^*_{mj}) &=& \varepsilon(u_{ki})L_{w,\ell}(u^*_{mj}) + \overline{\eta_\ell(u^*_{ki})}\eta_\ell(u^*_{mj}) + L_{w,\ell}(u_{ki})\varepsilon(u^*_{mj}) \\
&=&-\frac{1}{2}\delta_{ki}\overline{(\ell^*\ell)_{mj}} + (\ell^*w)_{ki}(w^*\ell)_{jm}-\frac{1}{2}(\ell^*\ell)_{ki}\delta_{mj}.
\end{eqnarray*}
Writing this in matrix form, we get
\[
\big(L_{w,\ell}(u_{ki}u^*_{mj})\big)_{1\le k,m\le d} = -\frac{1}{2}E_{ij}\ell^*\ell + \ell^*wE_{ij}w^*\ell -\frac{1}{2}\ell^*\ell E_{ij},
\]
and therefore the formula given in the Lemma.
\end{proof}

\section{Classification of generators on $\mathcal{U}_d$}

In this section we shall classify all L\'evy processes on $\mathcal{U}_d$, see also \cite{schuermann97} and \cite[Section 4]{franz00}.

The functionals $D_{ij}:\mathcal{U}_d\to\mathbb{C}$, $i,j=1,\ldots,d$ defined by
\begin{eqnarray*}
D_{ij}(\hat{u}_{kl}) &=& D_{ij}(u_{kl}) = i \delta_{ik}\delta_{jl}, \\
D_{ij}(\hat{u}^*_{kl}) &=& D_{ij}(u^*_{kl}) = -D_{ij}(\hat{u}_{lk}) = - i \delta_{il}\delta_{jk}, \\
D_{ij}(u) &=& 0 \quad \text{ if } u\not\in\text{span}\{\hat{u}_{ij}, \hat{u}^*_{ij}; i,j=1,\ldots,d\},
\end{eqnarray*}
for $i,j,k,l=1,\ldots,d$, are drift generators, since they are hermitian and form Sch{\"u}rmann triples together with the zero cocycle $\eta=0$ and the trivial representation $\rho=\varepsilon$.

Let $A=(a_{jk})\in\mathcal{M}_d(\mathbb{C})$ be a complex $d\times d$-matrix. It is not difficult to see that the triples $(\varepsilon, \eta_A:\mathcal{U}_d\to\mathbb{C}, G_{A})$, $i,j=1,\ldots,d$ defined by
\begin{eqnarray*}
\eta_{A}(\hat{u}_{jk}) &=& \eta_{A}(u_{jk}) = a_{jk}, \\
\eta_{A}(\hat{u}^*_{jk}) &=& \eta_{A}(u^*_{jk}) = - \eta_{A}(u_{kj}) = - a_{kj}, \\
\eta_{A}(1) &=&  \eta_{A}(uv)=0 \quad \text{ for } u,v\in\mathcal{U}^0_d,
\end{eqnarray*}
and
\begin{eqnarray*}
G_{A}(1) &=& G_{A}(\hat{u}_{jk}-\hat{u}^*_{kj}) = 0, \quad \text{ for } j,k=1,\ldots,d,\\
G_{A}(\hat{u}_{jk}+\hat{u}^*_{kj}) &=& - G_A\left(\sum_{l=1}^d \hat{u}^*_{lj}\hat{u}_{lk}\right) =  - \sum_{l=1}^d \overline{a_{lj}}a_{lk} = -(A^*A)_{jk}, \\
&& \hspace{4cm}\quad \text{ for } j,k=1,\ldots,d,  \\
G_{A}(uv) &=& \langle \eta_{A}(u^*), \eta_{A}(v)\rangle = \overline{\eta_{A}(u^*)}\eta_{A}(v),
\end{eqnarray*}
for $u,v\in\mathcal{U}_d^0$, are Sch{\"u}rmann triples. Furthermore, we have $\eta_{A}|_{K_2}=0$ and $G_{A}|_{K_3}=0$, i.e.\ the generators $G_{A}$ are Gaussian. On the elements $\hat{u}_{jk}, \hat{u}^*_{jk}$, $j,k=1,\ldots,d$, this gives
\begin{eqnarray*}
G_{A}( \hat{u}_{jk}) &=& -\frac{1}{2}(A^*A)_{jk}\\
G_{A}( \hat{u}^*_{jk}) &=& -\frac{1}{2}(A^*A)_{kj} \\
G_{A}( \hat{u}_{jk}\hat{u}_{lm}) &=& \overline{\eta_{A}(\hat{u}^*_{jk})}\eta_{A}(\hat{u}_{lm}) = - \overline{a_{kj}} a_{lm}, \\
G_{A}( \hat{u}_{jk}\hat{u}^*_{lm}) &=& \overline{a_{kj}}a_{ml}\\
G_{A}( \hat{u}^*_{jk}\hat{u}_{lm}) &=& \overline{a_{jk}}a_{lm} \\
G_{A}( \hat{u}^*_{jk}\hat{u}^*_{lm}) &=& -\overline{a_{jk}}a_{ml}
\end{eqnarray*}
for $j,k,l,m=1,\ldots,d$.

Let us denote the standard basis of $\mathcal{M}_d(\mathbb{C})$ by $E_{jk}$, $j,k=1,\ldots,d$. We define the functionals $G_{jk,lm}:\mathcal{U}_d\to\mathbb{C}$ by
\begin{eqnarray*}
G_{jk,lm}(1) &=& 0, \\
G_{jk,lm}(\hat{u}_{rs}) &=& -\frac{1}{2} \delta_{kr}\delta_{jl}\delta_{ms}=  -\frac{1}{2} (E_{jk}^*E_{lm})_{rs},  \qquad \text{ for } r,s=1,\ldots,d,\\
G_{jk,lm}(\hat{u}^*_{rs}) &=& -\frac{1}{2} \delta_{ks}\delta_{jl}\delta_{mr}=  -\frac{1}{2} (E_{jk}^*E_{lm})_{sr}, \qquad \text{ for } r,s=1,\ldots,d,\\
G_{jk,lm}(uv) &=& \langle \eta_{E_{jk}}(u^*), \eta_{E_{lm}}(v)\rangle = \overline{\eta_{E_{jk}}(u^*)}\eta_{E_{lm}}(v),
\end{eqnarray*}
for $u,v\in\mathcal{U}_n^0$, $j,k,l,m=1,\ldots,d$.

\begin{theorem}\label{class un gauss}
A generator $L:\mathcal{U}_d\to\mathbb{C}$ is Gaussian, if and only if it is of the form
\[
L=\sum_{j,k,l,m=1}^d \sigma_{jk,lm} G_{jk,lm} + \sum_{j,k=1}^d b_{jk} D_{jk},
\]
with a hermitian $d\times d$-matrix $(b_{jk})$ and a positive semi-definite $d^2\times d^2$-matrix $(\sigma_{jk,lm})$. It is a drift, if and only if $\sigma_{jk,lm}=0$ for $j,k,l,m=1,\ldots,d$.
\end{theorem}
\begin{proof}
Applying $L$ to Equation \eqref{u-hat rel}, we see that $L(\hat{u}_{jk})=-L(\hat{u}^*_{kj})$ has to hold for a drift generator. By the hermitianity we get $L(\hat{u}_{jk})=\overline{L(\hat{u}^*_{jk})}$, and thus a drift generator $L$ has to be of the form $\displaystyle\sum_{j,k=1}^n b_{ij} D_{ij}$ with a hermitian $d\times d$-matrix $(b_{ij})$.

Let $(\rho,\eta,L)$ be a Sch{\"u}rmann triple with a Gaussian generator $L$. Then we have $\rho=\varepsilon\,{\rm id}$, and $\eta(1)=0$, $\eta|_{K_2}=0$. By applying $\eta$ to Equation \eqref{u-hat rel}, we get $\eta(\hat{u}^*_{ij})=-\eta(\hat{u}_{ji})$. Therefore $\eta(\mathcal{U}_d)$ has at most dimension $d^2$ and the Sch{\"u}rmann triple $(\rho,\eta,L)$ can be realized on the Hilbert space $\mathcal{M}_d(\mathbb{C})$ (where the inner product is defined by $\langle A,B\rangle= \displaystyle\sum_{j,k=1}^d \overline{a_{jk}}b_{jk}$ for $A=(a_{jk}), B=(b_{jk})\in\mathcal{M}_d(\mathbb{C})$). We can write $\eta$ as
\begin{equation}\label{def eta}
\eta=\sum_{j,k=1}^d \eta_{A_{jk}}E_{jk} 
\end{equation}
where the matrices $A_{jk}$ are defined by $(A_{jk})_{lm}=\langle E_{lm},\eta(\hat{u}_{jk})\rangle$, for $j,k,l,m=1,\ldots,d$.

Then we get
\begin{eqnarray*}
L(\hat{u}^{\bullet_1}_{rs}\hat{u}^{\bullet_2}_{tu})&=& \langle\eta\big((\hat{u}^{\bullet_1}_{rs})^*\big),\eta(\hat{u}^{\bullet_2}_{tu}) \rangle \\
&=& \sum_{j,k,l,m=1}^d \langle\eta_{A_{jk}}\big((\hat{u}^{\bullet_1}_{rs})^*\big) E_{jk}), \eta_{A_{lm}}(\hat{u}^{\bullet_2}_{tu}) E_{lm})\rangle \\
&=&
\sum_{j,k=1}^d \overline{\eta_{A_{jk}}\big((\hat{u}^{\bullet_1}_{rs})^*\big)}\eta_{A_{jk}}(\hat{u}^{\bullet_2}_{tu}) \\
&=&
\left\{
\begin{array}{lcl}
\sum_{j,k=1}^d -\overline{(A_{jk})_{sr}}(A_{jk})_{tu} & \text{ if } & \hat{u}^{\bullet_1}=\hat{u}^{\bullet_2}=\hat{u} \\
\sum_{j,k=1}^d \overline{(A_{jk})_{sr}}(A_{jk})_{ut} & \text{ if } & \hat{u}^{\bullet_1}=\hat{u}, \hat{u}^{\bullet_2}=\hat{u}^* \\
\sum_{j,k=1}^n \overline{(A_{jk})_{rs}}(A_{jk})_{tu} & \text{ if } & \hat{u}^{\bullet_1}=\hat{u}^*,\hat{u}^{\bullet_2}=\hat{u} \\
\sum_{j,k=1}^d -\overline{(A_{jk})_{rs}}(A_{jk})_{ut} & \text{ if } & \hat{u}^{\bullet_1}=\hat{u}^{\bullet_2}=\hat{u}^*
\end{array} 
\right.
\\
&=& \sum_{l,m,p,q=1}^d \sigma_{lm,pq} G_{jk,lm}(\hat{u}^{\bullet_1}_{rs}\hat{u}^{\bullet_2}_{tu}) \\
&=& \left(\sum_{j,k,l,m=1}^d \sigma_{jk,lm} G_{jk,lm}+ \sum_{jmj=1}^n b_{jk}D_{jk}\right)(\hat{u}^{\bullet_1}_{rs}\hat{u}^{\bullet_2}_{tu})
\end{eqnarray*}
for $r,s,t,u=1,\ldots,d$, where $\sigma=(\sigma_{lm,pq})\in\mathcal{M}_{d^2}(\mathbb{C})$ is the positive semi-definite matrix defined by
\[
\sigma_{lm,pq}= \sum_{j,k=1}^d \overline{(A_{jk})_{lm}}(A_{jk})_{pq}
\]
for $l,m,p,q=1,\ldots,d$.

Setting $b_{jk}=-\frac{i}{2}L(\hat{u}_{jk}-\hat{u}^*_{kj})$, for $j,k=1,\ldots,d$, we get
\begin{eqnarray*}
L(\hat{u}_{rs}) &=& L\left(\frac{\hat{u}_{rs}+\hat{u}^*_{sr}}{2} +\frac{\hat{u}_{rs}-\hat{u}^*_{sr}}{2}\right) =  - \frac{1}{2}\sum_{p=1}^d \langle \eta(\hat{u}_{pr}),\eta(\hat{u}_{ps})\rangle+ib_{rs}  \\
&=&-\frac{1}{2} \sum_{j,k=1}^n (A_{jk}^*A_{jk})_{rs} +ib_{rs}  \\
&=&\left(\sum_{j,k,l,m=1}^d \sigma_{jk,lm} G_{jk,lm}+ \sum_{jmj=1}^d b_{jk}D_{jk}\right)(\hat{u}_{rs}) \\
L(\hat{u}^*_{sr}) &=& L\left(\frac{\hat{u}_{rs}+\hat{u}^*_{sr}}{2} - \frac{\hat{u}_{rs}-\hat{u}^*_{sr}}{2}\right) =  - \frac{1}{2}\sum_{p=1}^d \langle \eta(\hat{u}_{pr}),\eta(\hat{u}_{ps})\rangle-ib_{rs}  \\
&=&-\frac{1}{2} \sum_{j,k=1}^d (A_{jk}^*A_{jk})_{rs} -ib_{rs}  \\
&=&\left(\sum_{j,k,l,m=1}^d \sigma_{jk,lm} G_{jk,lm}+ \sum_{jmj=1}^d b_{jk}D_{jk}\right)(\hat{u}^*_{sr})
\end{eqnarray*}
where we used Equation \eqref{u-hat rel} for evaluating $L(\hat{u}_{rs}+\hat{u}^*_{sr})$. Therefore we have $\displaystyle L=\sum_{j,k,l,m=1}^d \sigma_{jk,lm} G_{jk,lm}+ \sum_{jmj=1}^n b_{jk}D_{jk}$, since both sides vanish on $K_3$ and on $\mathbf{1}$. The matrix $(b_{jk})$ is hermitian, since $L$ is hermitian,
\[
\overline{b_{jk}}= \frac{i}{2}\overline{L(\hat{u}_{jk}-\hat{u}^*_{kj})}=\frac{i}{2}L(\hat{u}^*_{jk}-\hat{u}_{kj})=b_{kj},
\]
for $j,k=1,\ldots,d$.

Conversely, let $\displaystyle L=\sum_{i,j=1}^d \sigma_{jk,lm} G_{jk,lm}+ \sum_{j,k=1}^d b_{jk}D_{jk}$ with a positive semi-definite $d^2\times d^2$-matrix $(\sigma_{jk,lm})$ and a hermitian $d\times d$-matrix $(b_{jk})$. Then we can choose a matrix $M=(m_{kl,ml})\in\mathcal{M}_{d^2}(\mathbb{C})$ such that $\displaystyle\sum_{p,q=1}^d \overline{m_{pq,jk}}m_{pq,lm}=\sigma_{jk,lm}$ for all $i,j,r,s=1,\ldots,d$. We define $\eta:\mathcal{U}_d\to\mathbb{C}^{d^2}$ by the matrices $A_{jk}$ with components $(A_{jk})_{lm}=m_{jk,lm}$ as in Equation \eqref{def eta}. It is not difficult to see that $(\varepsilon\,{\rm id}_{\mathbb{C}^{d^2}}, \eta, L)$ is a Sch{\"u}rmann triple and $L$ therefore a Gaussian generator.
\end{proof}

We can give generators of a Gaussian L\'evy process on $\mathcal{U}_n$ also in the following form, cf.\ \cite[Theorem 5.1.12]{schuermann93}
\begin{proposition}\label{II-gen-gauss}
Let $L^1,\ldots,L^n,M\in \mathcal{M}_d(\mathbb{C})$, with $M^*=M$, and let $H$ be an $n$-dimensional Hilbert space with orthonormal basis $\{e_1,\ldots,e_n\}$. Then there exists a unique Gaussian Sch\"urmann triple $(\rho,\eta,L)$ with
\begin{eqnarray*}
\rho &=& \varepsilon{\rm id}_H, \\
\eta(u_{jk}) &=& \sum_{\nu=1}^n L^\nu_{jk}e_\nu, \\
\eta(u^*_{jk}) &=& -\eta(u_{kj}), \\
L(u_{jk}) &=& \frac{1}{2} \sum_{r=1}^d \langle\eta(u^*_{jr}),\eta(u_{kr})\rangle + iM_{jk}
\end{eqnarray*}
for $1\le j,k\le d$.
\end{proposition}

The following theorem gives a classification of all L\'evy processes on $\mathcal{U}_d$.

\begin{theorem}\label{class un}
Let $H$ be a Hilbert space, $U$ a unitary operator on $H\otimes\mathbb{C}^d$, $A=(a_{jk})$ an element of the Hilbert space $H\otimes \mathcal{M}_d(\mathbb{C})$ and $\lambda=(\lambda_{jk})\in\mathcal{M}_d(\mathbb{C})$ a hermitian matrix. Then there exists a unique Sch{\"u}rmann triple $(\rho,\eta,L)$ on $H$ such that
\begin{subequations}\label{class un eq}
\begin{equation}
\rho(u_{jk}) = P_j U P^*_k,
\end{equation}
\begin{equation}
\eta(u_{jk}) = a_{jk},
\end{equation}
\begin{equation}
L(u_{jk}-u^*_{kj}) = 2i\lambda_{jk},
\end{equation}
\end{subequations}
for $j,k=1\ldots,d$, where $P_j:H\otimes\mathbb{C}^d\to H\otimes \mathbb{C}e_j\cong H$ projects a vector with entries in $H$ to its $j^{\rm th}$ component.

Furthermore, all Sch{\"u}rmann triples on $\mathcal{U}_n$ are of this form.
\end{theorem}
\begin{proof}
Let us first show that all Sch{\"u}rmann triples are of the form given in the theorem. If $(\rho,\eta,L)$ is a Sch{\"u}rmann triple, then we can use the Equations \eqref{class un eq} to define $U$, $A$, and $\lambda$. The defining relations of $\mathcal{U}_d$ imply that $U$ is unitary, since
\begin{eqnarray*}
U^*U  P^*_l &=& \sum_{j=1}^d \sum_{k=1}^d  P_j U^* P^*_k  P_k U P^*_l = \sum_{j=1}^d \sum_{k=1}^d \rho(u^*_{kj}u_{kl}) = \sum_{j=1}^d \delta_{jl} \rho(1) =  {\rm id}_{H\otimes e_l}, \\
UU^*  P^*_l &=& \sum_{j=1}^d \sum_{k=1}^d  P_j U P^*_k  P_k U^* P^*_l = \sum_{j=1}^d \sum_{k=1}^d \rho(u^*_{jk}u_{lk}) = \sum_{j=1}^d \delta_{jl} \rho(1) =  {\rm id}_{H\otimes e_l},
\end{eqnarray*}
for $l=1,\ldots,d$, where $e_1,\ldots,e_d$ denotes the standard basis of $\mathbb{C}^d$. The hermitianity of $\lambda$ is an immediate consequence of the hermitianity of $L$. 

Conversely, let $U$, $A$, and $\lambda$ be given. Then there exists a unique representation $\rho$ on $H$ such that $\rho(u_{jk}) = P_j U P^*_k$, for $j,k,=1,\ldots,d$, since the unitarity of $U$ implies that the defining relations of $\mathcal{U}_n$ are satisfied. We can set
$\eta(\hat{u}_{jk})=a_{jk}$, and extend via $\eta(u^*_{ki})=-\eta\left(\hat{u}_{ik}+\sum_{j=1}^d\hat{u}^*_{ji}\hat{u}_{jk}\right)= -a_{ik} - \sum_{j=1}^d\rho(\hat{u}_{ji})^* a_{jk}$, for $i,k=1,\ldots,d$ and $\eta(uv)=\rho(u)\eta(v)+\eta(u)\varepsilon(v)$ (i.e.\ Equation (\ref{levy-intro}.\ref{I-eta-cocycle}), for $u,v\in\mathcal{U}_d$, in this way we obtain the unique $(\rho,\varepsilon)$-cocycle with $\eta(\hat{u}_{jk})=a_{jk}$. Then we set $\displaystyle L(u_{jk})= i\lambda_{jk}-\frac{1}{2} \sum_{l=1}^d \langle a_{lj},a_{lk}\rangle$ and $\displaystyle L(u^*_{kj})= -i\lambda_{jk}-\frac{1}{2} \sum_{l=1}^d \langle a_{lj},a_{lk}\rangle$, for $j,k=1,\ldots,d$, and use Equation (\ref{levy-intro}.\ref{I-L-coboundary}) to extend it to all of $\mathcal{U}_d$. This extension is again unique, because the Relation \eqref{u-hat rel} implies $\displaystyle L(u_{jk}+u^*_{kj})=-\sum_{l=1}^d \langle a_{lj},a_{lk}\rangle$, and this together with $L(u_{jk}-u^*_{kj}) = 2i\lambda_{jk}$ determines $L$ on the generators $u_{jk}$,$u^*_{jk}$ of $\mathcal{U}_d$. But once $L$ is defined on the generators, it is determined on all of $\mathcal{U}_d$ thanks to Equation (\ref{levy-intro}.\ref{I-L-coboundary}).
\end{proof}

\section{Dilations of completely positive semigroups on $\mathcal{M}_d$}

Let $(\tau_t)_{t\ge 0}$ be a quantum dynamical semigroup on $\mathcal{M}_d$, i.e.\ a weakly continuous semigroup of completely positive maps $\tau_t:\mathcal{M}_d\to\mathcal{M}_d$.

\begin{definition}
A semigroup $(\theta_t)_{t\ge0}$ of not necessarily unital endomorphisms of $\mathcal{B}(\mathcal{H})$ with $\mathbb{C}^d\subseteq\mathcal{H}$ is called a {\em dilation} of $(\tau_t)_{t\ge 0})$, if
\[
\tau_t(X)=P\theta_t(X)P
\]
holds for all $t\ge 0$ and all $X\in \mathcal{M}_d=\mathcal{B}(\mathbb{C}^d)=P\mathcal{B}(\mathcal{H})P$. Here $P$ is the orthogonal projection from $\mathcal{H}$ to $\mathbb{C}^d$.
\end{definition}

\begin{example}\label{II-simple-dilation}
We can use the construction in Section \ref{II-first example} to get an example. Let $(\tau_t)_{t\ge 0}$ be the semigroup defined in \eqref{II-qdsg-exa}. We identify $\mathbb{C}^d$ with the subspace $\mathbb{C}^d\otimes\Omega\subseteq\mathbb{C}^d\otimes\Gamma\big(L^2(\mathbb{R}_+,\mathcal{K})\big)$. The orthogonal projection $P:\mathcal{H}\to\mathbb{C}^d$ is given by $P={\rm id}_{\mathbb{C}^d}\otimes P_\Omega$, where $P_\Omega$ denotes the projection onto the vacuum vector. Furthermore, we consider $\mathcal{M}_d$ as a subalgebra of $\mathcal{B}\Big(\mathbb{C}^d\otimes\Gamma\big(L^2(\mathbb{R}_+,\mathcal{K})\big)\Big)$ by letting a matrix $X\in\mathcal{M}_d$ act on $v\otimes w\in \mathbb{C}^d\otimes\Gamma\big(L^2(\mathbb{R}_+,\mathcal{K})\big)$ as $X\otimes P_\Omega$.

Note that we have
\[
\mathbb{E}_0(X)\otimes P_\Omega = PXP
\]
for all $X\in\mathcal{B}\Big(\mathbb{C}^d\otimes\Gamma\big(L^2(\mathbb{R}_+,\mathcal{K})\big)\Big)$.

Then the semigroup $(\eta_t)_{t\ge 0}$ defined in \eqref{II-e0sg-exa} is a dilation of $(\tau_t)_{t\ge 0}$, since
\begin{eqnarray*}
P\eta_t(X\otimes P_\Omega)P &=& PU_t\tilde{\gamma}_t(X\otimes P_\Omega)U^*_tP = PU_t(X\otimes{\rm id}_{\Gamma(L^2([0,t],\mathcal{K}))}\otimes P_\Omega)U^*_tP \\
&=& PU_t(X\otimes \mathbf{1})U^*_tP=\tau_t(X)\otimes P_\Omega
\end{eqnarray*}
for all $X\in\mathcal{M}_d$. Here we used that fact that the H-P cocycle $(U_t)_{t\ge 0}$ is adapted.
\end{example}

\begin{definition}
A dilation $(\theta_t)_{t\ge 0}$ on $\mathcal{H}$ of a quantum dynamical semigroup $(\tau_t)_{t\ge 0}$ on $\mathbb{C}^d$ is called {\em minimal}, if the subspace generated by the $\theta_t(X)$ from $\mathbb{C}^d$ is dense in $\mathcal{H}$, i.e.\ if
\[
\overline{{\rm span}\,\left\{\theta_{t_1}(X_1)\cdots\theta_{t_n}(X_n)v|t_1,\ldots,t_n\ge0,X_1,\ldots,X_n\in\mathcal{M}_d,v\in\mathbb{C}^d,n\in\mathbb{N}\right\}}=\mathcal{H}.
\]
\end{definition}

\begin{lemma}
It is sufficient to consider ordered times $t_1\ge t_2\ge\cdots\ge t_n\ge0$, since
\begin{gather*}
{\rm span}\,\left\{\theta_{t_1}(X_1)\cdots\theta_{t_n}(X_n)v|t_1\ge\ldots\ge t_n\ge0,X_1,\ldots,X_n\in\mathcal{M}_d,v\in\mathbb{C}^d,n\in\mathbb{N}\right\} \\
={\rm span}\,\left\{\theta_{t_1}(X_1)\cdots\theta_{t_n}(X_n)v|t_1,\ldots,t_n\ge0,X_1,\ldots,X_n\in\mathcal{M}_d,v\in\mathbb{C}^d,n\in\mathbb{N}\right\}
\end{gather*}
\end{lemma}

\begin{proof}
See \cite[Section 3]{bhat01}
\end{proof}

\begin{example}
We will now show that the dilation from Example \ref{II-simple-dilation} is not minimal, if $w$ and $\ell$ not linearly independent.

Due to the adaptedness of the cocycle $(U_t)_{t\ge 0}$, we can write
\[
\eta_t(X\otimes P_\Omega)= U_t\tilde{\gamma}_t(X\otimes P_\Omega)U^*_t = \big(U_t(X\otimes \mathbf{1})U^*_t\big)\otimes P_\Omega
\]
on $\Gamma\big(L^2([0,t],\mathcal{K})\big)\otimes \Gamma\big(L^2([t,\infty[,\mathcal{K})\big)$. Let
\[
\hat\eta_t(X)=\eta_t(X\otimes \mathbf{1})=U_t(X\otimes \mathbf{1})U^*_t
\]
for $X\in\mathcal{M}_d$ and $t\ge 0$, then we have
\[
\hat{\eta}_{t_1}(X_1)\cdots\hat{\eta}_{t_n}(X_n)v =\eta_{t_1}(X_1\otimes P_\Omega)\cdots\eta_{t_n}(X_n\otimes P_\Omega)v
\]
for $v\in\mathbb{C}^d\otimes\Omega$, $n\in\mathbb{N}$, $t_1\ge \cdots\ge t_n\ge0$, $X_1,\ldots,X_n\in\mathcal{M}_d$, i.e.\ time-ordered products of the $\hat{\eta}_t(X)$ generate the same subspace from $\mathbb{C}^d\otimes\Omega$ as the $\eta_{t}(X\otimes P_\Omega)$. 
Using the quantum It\^o formula, one can show that the operators $\hat{\eta}_t(X)$, $X\in\mathcal{M}_d$, satisfy the quantum stochastic differential equation.
\[
\hat{\eta}_t(X)=U_t(X\otimes \mathbf{1})U^*_t = X\otimes \mathbf{1}+\int_0^t U_s (wXw^*-X)U^*_s{\rm d}\Lambda_s, \qquad t\ge 0,
\]
if $\ell = \lambda w$ for some $\lambda\in\mathbb{C}$.

Since the quantum stochastic differential equation for $\hat{\eta}_t(X)$ has no creation part, these operators leave $\mathbb{C}^d\otimes\Omega$ invariant. More precisely, we have
\[
\left\{\hat{\eta}_{t_1}(X_1)\cdots\hat{\eta}_{t_n}(X_n)v\otimes\Omega|t_1\ge\ldots\ge t_n\ge0,X_1,\ldots,X_n\in\mathcal{M}_d,v\in\mathbb{C}^d,n\in\mathbb{N}\right\}=\mathbb{C}^d\otimes\Omega,
\]
and therefore the dilation $(\eta_t)_{t\ge 0}$ is not minimal, if $w$ and $\ell$ are not linearly independent. Note that in this case the quantum dynamical semigroup is also trivial, i.e.\ $\tau_t={\rm id}$ for all $t\ge 0$, since its generator vanishes.

One can show that the converse is also true, if $w$ and $\ell$ are linearly independent, then the dilation $(\eta_t)_{t\ge 0}$ is minimal.
\end{example}

The general form of the generator of a quantum dynamical semigroup on $\mathcal{M}_d$ was determined by \cite{gorini+kossakowski+sudarshan76,lindblad76}.
\begin{theorem}
Let $(\tau_t)_{t\ge 0}$ be a quantum dynamical semigroup on $\mathcal{M}_d$. Then there exist matrices $M,L^1,\ldots,L^n\in\mathcal{M}_d$, with $M^*=M$, such that the generator $\mathcal{L}=\frac{{\rm d}}{{\rm d}t}\tau_t$ is given by
\[
\mathcal{L}(X)=i[M,X]+\sum_{k=1}^n \left((L^k)^*XL^k-\frac{1}{2}\left\{X,(L^k)^*L^k\right\}\right)
\]
for $X\in\mathcal{M}_d$.
\end{theorem}

Note that $M,L^1,\ldots,L^n\in\mathcal{M}_d$ are not uniquely determined by $(\tau_t)_{t\ge 0}$.

Proposition \ref{II-gen-gauss} allows us to associate a L\'evy process on $\mathcal{U}_d$ to $L^1,\ldots,L^n,M$. It turns out that the cocycle constructed from this L\'evy process as in Section \ref{II-first example} dilates the quantum dynamical semigroup whose generator $\mathcal{L}$ is given by $L^1,\ldots,L^n,M$.

\begin{proposition}
Let $n\in\mathbb{N}$, $M,L^1,\cdots,L^n\in\mathcal{M}_d$, $M^*=M$, and let $(j_{st})_{0\le s\le t}$ be the L\'evy process on $\mathcal{U}_d$ over $\Gamma\big(L^2(\mathbb{R}_+,\mathbb{C}^n)\big)$, whose Sch\"urmann triple is constructed from $M,L^1,\cdots,L^n$ as in Proposition \ref{II-gen-gauss}. Then the semigroup $(\eta_t)_{t\ge 0}$ defined from the unitary cocycle
\[
U_t = \left(\begin{array}{ccc} j_{0t}(u_{11}) & \cdots & j_{0t}(u_{1d}) \\
\vdots & & \vdots \\
j_{0t}(u_{d1}) & \cdots & j_{0t}(u_{dd})
\end{array}\right)
\]
as in \eqref{II-e0sg-exa} is a dilation of the quantum dynamical semigroup $(\tau_t)_{t\ge 0}$ with generator 
\[
\mathcal{L}(X)=i[M,X]+\sum_{k=1}^n \left((L^k)^*XL^k-\frac{1}{2}\left\{X,(L^k)^*L^k\right\}\right)
\]
for $X\in\mathcal{M}_d$.
\end{proposition}
\begin{proof}
The calculation is similar to the one in Section \ref{II-first example}.
\end{proof}

We denote this dilation by $(\eta_t)_{t\ge 0}$ and define again $\hat\eta_t:\mathcal{M}_d\to\mathcal{B}\Big(C^d\otimes\Gamma\big(L^2([0,t],\mathcal{K})\big)\Big)$, $t\ge 0$ by
\[
\hat\eta_t(X)=\eta_t(X\otimes \mathbf{1})=U_t(X\otimes \mathbf{1})U^*_t
\]
for $X\in\mathcal{M}_d$.

Denote by $\mathcal{Q}_d$ the subalgebra of $\mathcal{U}_d$ generated by $u_{ij}u^*_{k\ell}$, $1\le i,j,k,\ell \le d$. This is even a subbialgebra, since
\[
\Delta(u_{ij}u^*_{k\ell}) =\sum_{r,s=1}^d u_{ir}u^*_{ks}\otimes u_{rj}u^*_{s\ell}
\]
for all $1\le i,j,k,\ell\le d$. 

\begin{lemma}
Let $\eta:\mathcal{U}_d\to H$ be the cocycle associated to $L^1,\ldots,L^n,M\in \mathcal{M}_d(\mathbb{C})$, with $M^*=M$, in Proposition \ref{II-gen-gauss}.
\begin{description}
\item[(a)]
$\eta$ is surjective, if and only if $L^1,\ldots,L^n$ are linearly independent.
\item[(b)]
$\eta|_{\mathcal{Q}_d}$ is surjective, if and only if $I,L^1,\ldots,L^n$ are linearly independent, where $I$ denotes the identity matrix.
\end{description}
\end{lemma}
\begin{proof}
\begin{description}
\item[(a)]
$\mathcal{U}_d$ is generated by $\{u_{ij}|1\le i,j\le d\}$, so, by Lemma \ref{levy-intro}.\ref{I-lemma-surj}, we have $\eta(\mathcal{U}_d)={\rm span}\{\eta(u_{ij})|1\le i,j\le d\}$.

Denote by $\Lambda_1:H\to {\rm span}\left\{(L^1)^*,\ldots,(L^n)^*\right\}\subseteq\mathcal{M}_d(\mathbb{C})$ the linear map defined by $\Lambda_1(e_\nu)=(L^\nu)^*$, $\nu=1,\ldots,n$. Then we have ${\rm ker}\,\Lambda_1=\eta(\mathcal{U}_d)^\perp$, since
\[
\langle v , \eta(u_{ij})\rangle = \sum_{\nu=1}^d \overline{v_\nu} L^\nu_{ij} = \overline{\Lambda_1(v)_{ji}},\qquad 1\le i,j\le d,
\]
for $v=\sum_{\nu=1}^n v_\nu e_\nu\in H$. 
The map $\Lambda_1$ is injective, if and only if $L^1,\ldots,L^n$ are linearly independent. Since ${\rm ker}\,\Lambda_1=\eta(\mathcal{U}_d)^\perp$, this is also equivalent to the surjectivity of $\eta|_{\mathcal{U}_d}$.
\item[(b)]
We have $\eta(\mathcal{Q}_d)={\rm span}\{\eta(u_{ij}u^*_{k\ell})|1\le i,j,k,\ell\le d\}$.

Denote by $\Lambda_2:H\to {\rm span}\{(L^1)^*\otimes I-I\otimes (L^1)^*,\ldots,(L^n)^*\otimes I-I\otimes (L^n)^*\}\subseteq\mathcal{M}_d(\mathbb{C})\otimes\mathcal{M}_d(\mathbb{C})$ the linear map defined by $\Lambda_1(e_\nu)=(L^\nu)^*\otimes I - I \otimes (L^\nu)^*$, $\nu=1,\ldots,n$. Then we have ${\rm ker}\,\Lambda_2=\eta(\mathcal{Q}_d)^\perp$, since
\begin{eqnarray*}
\langle v , \eta(u_{ij}u^*_{k\ell})\rangle &=& \langle v , \varepsilon(u_{ij})\eta(u^*_{k\ell})+\eta(u_{ij}\varepsilon(u^*_{k\ell})\rangle \\
&=& \langle v , -\delta_{ij}\eta(u_{\ell k})+\eta(u_{ij}\delta_{k\ell})\rangle \\
&=& \sum_{\nu=1}^d \overline{v_\nu} (L^\nu_{ij}\delta_{k\ell}-\delta_{ij}L^\nu_{\ell k}) \\
&=& \overline{\Lambda_2(v)_{ji,k \ell}}, \qquad 1\le i,j,k,\ell\le d,
\end{eqnarray*}
for $v=\sum_{\nu=1}^n v_\nu e_\nu\in H$. 

The map $\Lambda_2$ is injective, if and only if $L^1,\ldots,L^n$ are linearly independent and $I\not\in{\rm span}\{L^1,\ldots,L^n\}$, i.e.\ iff $I,L^1,\ldots,L^n$ are linearly independent. Since ${\rm ker}\,\Lambda_2=\eta(\mathcal{Q}_d)^\perp$, it follows that this is equivalent to the surjectivity of $\eta|_{\mathcal{Q}_d}$.
\end{description}
\qed\end{proof}

Bhat \cite{bhat01,bhat03} has given a necessary and sufficient condition for the minimality of dilations of the form we are considering.

\begin{theorem}\cite[Theorem 9.1]{bhat01}
The dilation $(\eta_t)_{t\ge 0}$ is minimal if and only if $I,L^1,\ldots,L^n$ are linearly independent.
\end{theorem}

\begin{remark}
The preceding arguments show that the condition in Bhat's theorem is necessary. Denote by $\mathcal{H}_0$ the subspace of $\Gamma\big(L^2(\mathbb{R}_+,\mathbb{C}^n)\big)$, which is generated by operators of form $j_{st}(u_{ij}u^*_{k\ell})$. By Theorem \ref{levy-intro}.\ref{I-cyclic-thm}, this subspace is dense in $\Gamma\big(L^2(\mathbb{R}_+,\eta(\mathcal{Q}_d)\big)$. Therefore the subspace generated by elements of the $\hat{\eta}_t(X)=U_t(X\otimes\mathbf{1})U^*_t$ from $\mathbb{C}^d\otimes \Omega$ is contained in $\mathbb{C}^d\otimes\mathcal{H}_0$. If $\eta$ is minimal, then this subspace in dense in $\mathbb{C}^d\otimes\Gamma\big(L^2(\mathbb{R}_+,\mathbb{C}^n)\big)$. But this can only happen if $\mathcal{H}_0$ is dense in $\Gamma\big(L^2(\mathbb{R}_+,\mathbb{C}^n)\big)$. This implies $\eta(\mathcal{Q}_d)=\mathbb{C}^d$ and therefore that $I,L^1,\ldots,L^n$ are linearly independent.

Bhat's theorem is actually more general, it also applies to dilations of quantum dynamical semigroups on the algebra of bounded operators on an infinite-dimensional separable Hilbert space, whose generator involves infitely many $L$'s, see \cite{bhat01,bhat03}.
\end{remark}